\documentclass[oneside,	a4paper]{memoir}

\usepackage[dvipsnames]{xcolor}

\usepackage[utf8]{inputenc} \usepackage{lmodern} \usepackage[ngerman, english]{babel} \usepackage{csquotes}

\usepackage{emptypage} 

\usepackage{booktabs}       \usepackage[font=footnotesize,labelfont=sc]{caption} 
\usepackage{subcaption} 

\usepackage{wrapfig} \usepackage{tikz, tikz-cd, pgfplots} \usetikzlibrary{external}
\pgfplotsset{compat=1.16}
\usepackage{tikzscale}

\usepackage[ruled,linesnumbered]{algorithm2e}


\usepackage{mathtools} \usepackage{amsthm, amssymb} \usepackage{thmtools} \usepackage{mathrsfs} \usepackage{aligned-overset} 

\usepackage{cancel} 

\usepackage[
    style=authoryear,
    backend=biber,
    seconds=true,
    date=year,
    urldate=iso,
    maxcitenames=2,
    natbib,
    autocite=footnote,
]{biblatex}
\usepackage{xurl} 

\usepackage{marginfix} \usepackage{microtype}

\usepackage{enumitem} \usepackage{appendix} \usepackage[
    unicode,
colorlinks,
    citecolor=.,
]{hyperref}

\hypersetup{hypertexnames=false}

\makeatletter

\newcommand*{\cF}{\mathcal{F}}
\newcommand*{\cE}{\mathcal{E}}
\newcommand*{\cH}{\mathcal{H}}

\newcommand*{\cX}{\mathcal{X}}
\newcommand*{\cA}{\mathcal{A}}

\newcommand*{\bb}{\mathbf{b}}

\newcommand*{\ie}{i.e.\ }
\newcommand*{\eg}{e.g.\ }

\DeclareMathOperator*{\Root}{Root}

\newcommand*{\domain}{\mathbb{X}}
\newcommand*{\range}{\mathbb{Y}}
\newcommand*{\field}{\mathbb{K}}
\newcommand*{\real}{\mathbb{R}}
\newcommand*{\complex}{\mathbb{C}}
\newcommand*{\nat}{\mathbb{N}}
\newcommand*{\integer}{\mathbb{Z}}
\newcommand*{\sphere}[1][\dims-1]{\mathbb{S}^{#1}}
\newcommand*{\ball}{\mathrm{B}}

\newcommand*{\inter}{ \mathchoice{\mskip1.5mu{:}\mskip1.5mu}{\mskip1.5mu{:}\mskip1.5mu}{{:}}{{:}}
}

\newcommand*{\transpose}{T}

\DeclareMathOperator*{\Span}{span}
\DeclareMathOperator{\diag}{diag}
\newcommand*{\ones}{\mathbf{1}}
\newcommand*{\ind}{\mathbf{1}}
\newcommand{\identity}{\mathbb{I}}

\newcommand*{\orthGroup}[1][\dims]{\ensuremath{O(#1)}}
\newcommand*{\euclGroup}[1][\dims]{\ensuremath{E(#1)}}
\newcommand*{\translGroup}[1][\dims]{\ensuremath{T(#1)}}

\newcommand*{\dims}{d}

\newcommand*{\metric}{d}
\newcommand*{\closure}[1]{\overline{#1}}
\newcommand*{\interior}[1]{\mathrm{int}\, #1}

\newcommand*{\im}{i} \newcommand*{\conj}[1]{\overline{#1}}
\newcommand*{\Ft}{\mathcal{F}} \newcommand*{\modifiedBessel}{K_\nu}
\newcommand{\scal}[2]{\langle #1,#2 \rangle }

\newcommand*{\modCont}{\omega} 
\newcommand*{\proj}[1]{P_{#1}}

\DeclareMathOperator{\mon}{mon}
\newcommand*{\uP}[1][\lambda]{C^{(#1)}} \newcommand*{\uPnorm}[1][\lambda]{\hat{C}^{(#1)}} 

\newcommand*{\fshift}{S}

\newcommand*{\info}{I}
\newcommand*{\varinfo}{J}
\newcommand*{\liminfo}{\mathscr{I}}
\newcommand*{\Info}{\mathbf{I}}
\newcommand*{\VarInfo}{\hat{\mathbf{I}}}
\newcommand*{\VarLiminfo}{\hat{\mathscr{I}}}
\newcommand*{\gradients}{G}
\newcommand*{\basis}{b}
\newcommand*{\Basis}{B}
\newcommand*{\radius}{\lambda}

\newcommand*{\rv}{\mathbf{v}}
\newcommand*{\rw}{\mathbf{w}}
\newcommand*{\rcov}{\mathbf{\Sigma}}
\newcommand*{\rmean}{\mathbf{\mu}}

  \newcommand*{\distance}{\Delta}

\newcommand*{\spectMeasure}{\sigma}
\newcommand*{\schoenbergMeas}{\nu}
\newcommand*{\PDset}{\mathcal{P}}

\newcommand*{\obj}{f} \newcommand*{\Obj}{\mathbf{f}} \newcommand*{\bump}{\varphi}
\newcommand*{\cost}{J}
\newcommand*{\Cost}{\mathbf{J}}
\newcommand*{\stepsize}{\eta}

\DeclareMathOperator*{\argmin}{arg\,min}
\DeclareMathOperator*{\argmax}{arg\,max}
\newcommand*{\bigO}{\mathcal{O}}
\newcommand*{\timestep}{n}

\newcommand*{\batchsize}{b}
\newcommand*{\Batchsize}{B}
\newcommand{\param}{w}
\newcommand{\Param}{W}
\newcommand{\Tparam}{\mathbf{m}}
\newcommand{\Loss}{\mathcal{L}}
\newcommand*{\loss}{\ell}
\newcommand*{\lr}{h}
\newcommand*{\model}{\phi}
\newcommand*{\step}{\mathbf{d}}

\newcommand*{\target}{f}

\newcommand*{\regularizer}{R}

\newcommand*{\network}{\mathfrak{N}}
\newcommand*{\nodes}{\mathbb{V}}

\newcommand*{\paramSpace}{\Theta}
\newcommand*{\polyEfficient}{\mathcal{E}_P}
\newcommand*{\efficient}{\mathcal{E}}

\newcommand*{\weight}{w}
\newcommand*{\bias}{\beta}
\newcommand*{\inSgt}{\blacklozenge}
\newcommand*{\outSgt}{*}
\newcommand*{\placeholder}{\bullet}

\newcommand*{\activation}{\psi}
\DeclareMathOperator{\sigmoid}{{\normalfont{sigmoid}}}

\newcommand*{\response}{\Psi}
\newcommand*{\tin}{\text{in}}
\newcommand*{\tout}{\text{out}}

\newcommand*{\Domain}{\mathcal{X}}

\newcommand*{\lebesgue}{\lambda}

\newcommand*{\iid}{\mathrm{iid}}
\newcommand*{\E}{\mathbb{E}}
\renewcommand{\Pr}{\mathbb{P}}
\newcommand*{\filt}{\mathcal{F}}
\newcommand*{\cG}{\mathcal{G}}
\newcommand*{\density}{\varphi}
\newcommand*{\charfct}{\psi} \DeclareMathOperator*{\support}{supp}
\newcommand{\cdf}{\Phi}
\newcommand*{\borel}{\mathcal{B}}
\newcommand*{\normal}{\mathcal{N}}
\newcommand*{\uniform}{\mathcal{U}}
\newcommand*{\rP}{\mathbf{P}}
\newcommand*{\indep}{\perp\!\!\!\perp}
\newcommand*{\as}{\text{a.s}\@ifnextchar.{}{\text{.\@} }}

\DeclareMathOperator{\LUE}{{\hypersetup{hidelinks}\hyperref[def: LUE]{LUE}}} \DeclareMathOperator{\BLUE}{{\hypersetup{hidelinks}\hyperref[def: BLUE]{BLUE}}} \DeclareMathOperator{\Cov}{Cov}
\DeclareMathOperator{\Var}{Var}
\newcommand*{\problemSpace}{\mathbb{M}}

\newcommand*{\rf}{\mathbf{f}}
\newcommand*{\rg}{\mathbf{g}}
\newcommand*{\obs}{Y} \newcommand*{\Obs}{\mathbf{Y}} \newcommand*{\hamiltonian}{\mathbf{H}}

\newcommand*{\noise}{\varsigma}
\newcommand*{\jet}[1][k]{\mathrm{J}^{#1}}

\newcommand*{\limf}{\mathfrak{f}}
\newcommand*{\limgdot}{\mathfrak{g}}
\newcommand*{\gsa}{\mathfrak{G}}

\newcommand*{\scale}{s}

\newcommand*{\C}{\mathcal{C}}
\newcommand*{\kernel}{\kappa}
\newcommand*{\ikernel}{C}
\newcommand*{\sqC}{\kernel}

\definecolor{magenta}{HTML}{D81B66}
\definecolor{blue}{HTML}{1F76C1}
\definecolor{yellow}{HTML}{FFC107}
\definecolor{teal}{HTML}{00B981}

\newcommand*{\red}[1]{{\color{red} #1}}
\newcommand*{\blue}[1]{{\color{blue} #1}}
\newcommand*{\magenta}[1]{{\color{magenta} #1}}
\newcommand*{\teal}[1]{{\color{teal} #1}}

\newcommand*{\black}[1]{{\color{black} #1}}
\newcommand*{\green}[1]{{\color{LimeGreen} #1}}

 \newcommand*{\problem}{\mathfrak{m}}
\newcommand*{\Problem}{\mathbf{M}}

\theoremstyle{plain}\newtheorem{prop}{Proposition}[section]
\newtheorem{lemma}[prop]{Lemma}
\newtheorem{corollary}[prop]{Corollary}
\newtheorem{theorem}[prop]{Theorem}

\newtheorem{extension}[prop]{Extension}

\theoremstyle{definition}
\newtheorem{definition}[prop]{Definition}
\newtheorem{example}[prop]{Example}
\newtheorem{assumption}[prop]{Assumption}
\newtheorem{motivation}[prop]{Motivation}

\theoremstyle{remark}
\newtheorem{remark}[prop]{Remark}

\setbinding{0mm}
\settrims{0pt}{0pt}
\settypeblocksize{58\baselineskip}{353bp}{*}

\setlrmargins{20mm}{*}{*}

\setulmargins{25mm}{*}{*}

\setheadfoot{2\onelineskip}{2\onelineskip}
\setheaderspaces{*}{2\baselineskip}{*}
\setmarginnotes{4mm}{48mm}{\onelineskip}
\checkandfixthelayout

\setsecnumdepth{subsection}
\maxtocdepth{subsection}

\defbibheading{secnotoc}[\bibname]{
	\section*{#1}\markboth{#1}{#1}
}

\newcommand*{\dedication}[1]{
    \clearpage
    \begin{center}
        \thispagestyle{empty}
        \vspace*{\fill}\LARGE \textit{#1}\vspace*{\fill}
    \end{center}
    \clearpage
}

\newlength{\fullwidthlen}
\newenvironment{fullwidth}{\if@twoside
        \begin{adjustwidth*}{}{-\fullwidthlen}\else
        \checkoddpage\ifoddpage
            \begin{adjustwidth*}{}{-\fullwidthlen}
        \else
            \begin{adjustwidth*}{-\fullwidthlen}{}
        \fi
    \fi
}{\end{adjustwidth*}}

\newcommand*{\fullwidthcaption}[1]{
	\sideparmargin{outer}
	\vspace{-\baselineskip}
	\sidepar{\vspace{\baselineskip}
        \caption{#1}
    }
}

\makepagestyle{tufte}
\makeoddhead{tufte}{\rightmark}{}{\thepage}

\if@twoside
    \makeevenhead{tufte}{\thepage}{}{\leftmark}
\else
    \makeevenhead{tufte}{\leftmark}{}{\thepage}
\fi
\makeevenfoot{tufte}{}{}{}
\makeoddfoot{tufte}{}{}{}

\makepsmarks{tufte}{\createmark{chapter}{left}{nonumber}{}{}
  \createmark{section}{right}{nonumber}{}{}
  \createplainmark{toc}{both}{\contentsname}
  \createplainmark{lof}{both}{\listfigurename}
  \createplainmark{lot}{both}{\listtablename}
  \createplainmark{bib}{both}{\bibname}
  \createplainmark{index}{both}{\indexname}
  \createplainmark{glossary}{both}{\glossaryname}
}
\makerunningwidth{tufte}{\headwidth}
\makeheadposition{tufte}{flushright}{flushleft}{}{}

\setsidecaps{\marginparsep}{\marginparwidth}
\sidecapmargin{outer}
\setsidecappos{t}

\makeatother
 
\title{
	A Distributional View of
	High Dimensional Optimization
}
\author{Felix Benning}

\begin{document}
{
	\setlrmargins{*}{*}{1}
\checkandfixthelayout

	\clearpage
\begin{titlingpage}
    \let\cleardoublepage\clearpage
\thispagestyle{empty}
\begin{center}
\vfill

{
    \Huge\textbf{
    A Distributional View of
	High Dimensional Optimization
    }
}
\vfill

\begin{figure}[ht]
\begin{center}
\includegraphics[scale = 0.75]{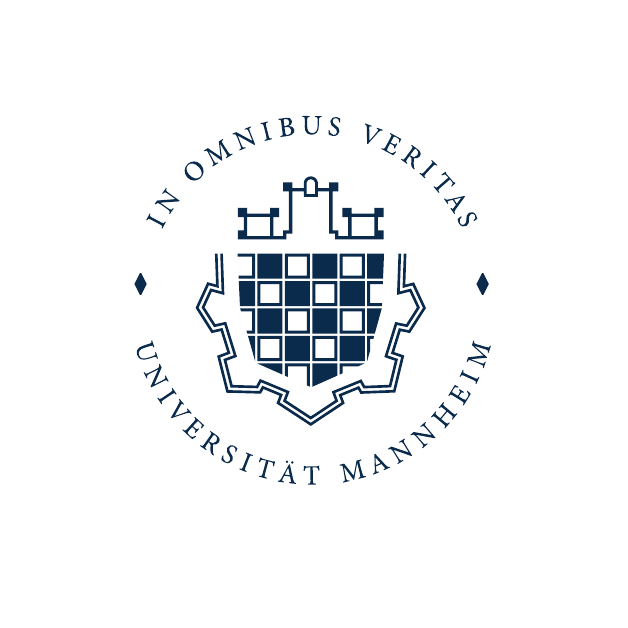}
\end{center}
\end{figure}

\vfill
{
\Large Inauguraldissertation\\
zur Erlangung des akademischen Grades\\
eines Doktors der Naturwissenschaften\\
der Universität Mannheim\par
}
\par\vspace{3\baselineskip}
{\large vorgelegt von}
\par\vspace{\baselineskip}
{\Large Felix Benning}
\par\vspace{.5\baselineskip}
{\large aus Nürtingen}
\par\vspace{3\baselineskip}
{\Large Mannheim, 2025}
\end{center}
\end{titlingpage}

\begin{titlingpage}
\let\cleardoublepage\clearpage
\vspace*{\fill} 
\begin{center}
\begin{tabular}{ll}
Dekan: & Prof. Dr.\ Claus Hertling, Universität Mannheim \\
Referent: & Prof.\ Dr.\ Leif D\"oring, Universität Mannheim \\
Korreferent: & Prof.\ Dr.\ Simon Weissmann, Universität Mannheim\\
Korreferent: & Prof. Dr. Steffen Dereich, Universität Münster\\
&\\
\multicolumn{2}{l} {Tag der mündlichen Prüfung: \, 18. Juli 2025}
\end{tabular}
\end{center}
\end{titlingpage} 	\frontmatter

	\dedication{To Veronika}
}

	{ \hypersetup{hidelinks} \tableofcontents }

	\chapter*{Acknowledgements}

Writing a dissertation in mathematics is an exercise for the privileged.
I was fortunate to be born into a stable family, where my parents -- especially
my mother -- always had the energy to discuss any rule they set. I was lucky to
have had Hans Schmiga as my mathematics teacher in school, who bothered to
introduce any topic with a derivation for those who were interested. And had the
tutor of the mathematics preparatory course for my economics bachelor not
responded to my question `What does matrix multiplication actually
\emph{mean}?'\footnote{
    I now know: A matrix is a basis representation of a linear map; matrix
    multiplication is defined to represent the composition of linear maps.
} with `I don't know, ask a
mathematician', I might never have taken the maths minor that allowed me
to gradually shift towards mathematics.

I was also fortunate that my parents had the means to ensure that I could focus
on my studies. And by chance, when I began my studies in Mannheim, my
supervisor Leif Döring had just started as a junior professor. 
He is incredibly gifted at teaching the fundamentals, and I was lucky
to be part of one of the few cohorts he taught Analysis 1 and 2 to.  His
ability to motivate abstract concepts and teach how to structure formally
correct proofs left a lasting impression on me.

Beyond these pivotal events that led to this attempt at a PhD, there were so
many people that helped me through its challenges it is impossible to
recount them all. Special thanks goes to my
friend Julie Naegelen, who supported me through the ups and downs and always
found the time to be a brutal proofreader. And I cannot put into words
how grateful I am to my girlfriend Veronika for putting up with my absentminded,
sleepless self.

\smallskip

\textit{To my family, friends, collegues -- everyone that helped me get here, on
purpose or by accident -- thank you!}

 	\footnotesinmargin

	\begin{refsection}
		\chapter{Preface}

\vspace{-\onelineskip}
\begin{quote}
    \textit{
        ``He asks this shepherd's boy: `How many seconds in
        eternity?' The shepherd's boy says, `there is this mountain of pure
        diamond. It takes an hour to climb it and an hour to go around it, and every
        hundred years a little bird comes and sharpens its beak on the diamond
        mountain. And when the entire mountain is chiseled away, the first second of
        eternity will have passed.'
        You may think that's a hell of a long time.
        Personally, I think that's a hell of a bird.''
    }\par\raggedleft--- \textup{The 12th Doctor}, ``Heaven Sent''
\end{quote}

In many areas of modern science and engineering, complex models are no longer
built manually but instead shaped through optimization. As this \emph{machine
learning} paradigm has grown in prominence, so too has the importance of
general-purpose optimization tools.

Each time a model is modified -- whether by changing its architecture, its loss
function, or its input representation -- the underlying optimization problem
shifts. In such a setting, designing a custom optimization method for every
iteration is not only impractical, but counterproductive. Consequently, there is
little appetite for crafting bespoke optimization algorithms tailored to each
new objective. Instead, practitioners rely on black-box optimization.

\textbf{Black-box optimization} provides a clean separation of concerns between the user
and the underlying optimization process. In this framework, the user only needs
to implement an oracle that can evaluate the objective function and its
derivatives up to a certain order. The black-box optimizer takes care of the rest: selecting the next query
points, exploring the objective landscape, and iterating toward an optimal
solution.
\begin{quote}
    A \(k\)-th order black box optimizers requires a \(k\)-th order oracle which
    provides all derivatives \(\obj(x), \nabla \obj(x), \dots\) up to order
    \(k\) at an arbitrary point \(x\) chosen by the optimization algorithm.
\end{quote}
For black box optimization to work, certain assumptions about the objective
function \(\obj\) are necessary. It is fairly intuitive that an optimizer cannot
optimize a discontinuous function with an infinite domain, so \emph{continuity} is the
most basic assumption.
However, as we will show in Chapter \ref{chap: worst case black box optimization}, continuity assumptions alone lead to a
curse of dimensionality: the number of oracle evaluations required to achieve a
meaningful optimization result grows exponentially with the dimension of the
problem. In practical machine learning settings, where models can have millions
of parameters, this exponential increase in complexity makes optimization seem
impossible.
In classical (convex) optimization theory this problem is solved by forcing
information into the gradient -- a lower bound on the gradient norm using the
distance of the function value to the global optimum.\footnote{
    This is known as the PL inequality. See Equation \eqref{eq: PL-inequality}.
 }
Unfortunately, this assumption also implies that all critical points
are global optima. This is unrealistic for machine learning, where
exponentially more saddle points occur than local minima
\autocite{dauphinIdentifyingAttackingSaddle2014}. Classical optimization theory
therefore suggests that non-convex, global optimization is utterly hopeless
in high dimension.

In practice, however, optimization in machine learning seems to work -- even in
high-dimensional spaces where classical optimization theory suggests failure.
This observation implies that there must be a theory that accounts for the
\emph{realistic} behavior of objective functions in machine learning, rather than
relying on the worst-case scenario that is often used in classical theory.

In Chapter \ref{chap: worst case black box optimization}, we will examine how
the worst-case objective functions that lead to the curse of dimensionality are
adversarially chosen: every oracle evaluation returns zero for the function value,
gradient and higher oder derivatives. This deprives the optimizer of useful directional
information, which would be provided by the gradient for any other value. In
contrast, a randomly selected objective function is unlikely to exhibit this
pathological behavior.

It is therefore natural to study the \textbf{optimization of random functions}
to explain the realistic behavior of optimizers. This approach is known as
`Bayesian optimization' or `Kriging' in the literature.
After the first chapter on Worst-case optimization (Chapter \ref{chap: worst
case black box optimization}), the thesis therefore structured as follows.

\paragraph{Outline of Part \ref{part: distributional optimization} (Distributional optimization).}

In Chapter \ref{chap: distributional optimization}
we properly introduce this Bayesian approach to optimization and Chapter
\ref{chap: measure theory for rf opt} addresses its measure
theoretical problems.
The perhaps most significant challenge of this approach lies in the selection of
a distribution over functions. A pragmatic solution to this problem are
uniformity assumptions, specifically isotropy, which postulates that no
coordinate system for the random function is inherently favoured. This
assumption and other approaches to uniformity are introduced in Chapter
\ref{chap: distributions over functions}.  Having established the ingredients of
random function optimization, we turn to the
challenge of optimizing \emph{high-dimensional} random functions for which
classical Bayesian optimization methods are ill suited. In Chapter \ref{chap:
random function descent} we introduce a Bayesian method, we
call `Random Function Descent', that is scalable to
high dimension. Under the isotropy assumption (Chapter \ref{chap: distributions
over functions}), this method happens to coincide with gradient descent using a
specific step size schedule -- bridging
the gap between classical worst case optimization approaches and Bayesian
optimization. Finally, in Chapter \ref{chap: predictable progres} we analyze
the behavior of classical gradient span optimization algorithms. We show that
asymptotically, as the dimension increases, the optimization progress
converges to a deterministic sequence. This implies a tight concentration
around the average case performance of optimizers in high dimension.

\paragraph*{Outline of Part \ref{part: supervised machine learning} (Supervised machine learning).}

Having motivated the random objective function assumption
as an answer to the hopelessness of worst-case optimization in Part \ref{part: distributional optimization},
we replace this argument of desperation with a more axiomatic approach in in Chapter \ref{chap: random objective from
exchangeable data}. In this approach, specific to supervised machine learning, we will derive random objective functions
from the assumption of exchangeable data. Then we outline how this results in
an elegant interpretation of the most commonly used loss functions in machine
learning as maximum a posteriori estimators. In Section \ref{sec: random linear
model} we will present evidence against stationarity of the objective function,
giving more importance to non-stationary isotropy as characterized in Section
\ref{sec:positive_definite}.

Finally, in Chapter \ref{chap: optimization landscape of shallow neural networks}
we analyze the losslandscape induced by shallow neural networks on
random regression problems. Under extremely general unversality
assumptions about the target function,\footnote{
    all continuous functions lie in the support of the distribution over target
    functions
} we find that on the domain of efficient network parameters\footnote{
    the parameters whose response functions cannot be reproduced by smaller
    networks
}
the objective function is a Morse function. This implies positive curvature in
all directions at local minima and thereby isolated local minima with strongly
convex neighborhoods.

\paragraph*{Attribution}

Chapters based on papers I contributed to are often
\emph{verbatim} copies of these papers with minor modifications.
\begin{itemize}[noitemsep,align=left]
    \item[\textbf{Chapter \ref{chap: worst case black box optimization}:}]
    Summarizes well known results on worst-case optimization theory.
    The generalization to the modulus of continuity is non-standard but 
    unlikely to be new.

    \item[\textbf{Chapter \ref{chap: distributional optimization}:}]
    Summarizes well known results about Bayesian optimization
    \autocite[e.g.][]{garnettBayesianOptimization2023}.

    \item[\textbf{Chapter \ref{chap: measure theory for rf opt}:}]
    Verbatim copy of \textcite{benningMeasureTheoryConditionally2025} except for its introduction, which
    was moved to the previous chapter.

    \item[\textbf{Chapter \ref{chap: distributions over functions}:}]
    Section \ref{sec: input invariance} is based on Section
    F in \textcite{benningRandomFunctionDescent2024},
    Section \ref{sec: characterization of covariance kernels} and \ref{sec:positive_definite} 
    are based on \textcite{benningSchoenbergCharacterizationContinuous2025}. Section \ref{sec: covariance of derivatives}
    and \ref{sec: strictly pos definite derivatives}
    are based on \textcite{benningGradientSpanAlgorithms2024}

    \item[\textbf{Chapter \ref{chap: random function descent}:}] 
    Based on \textcite{benningRandomFunctionDescent2024}.

    \item[\textbf{Chapter \ref{chap: predictable progres}:}]
    Based on \textcite{benningGradientSpanAlgorithms2024}.

    \item[\textbf{Chapter \ref{chap: random objective from exchangeable data}:}]
    To the best of our knowledge the framework of exchangable data for
    supervised learning is new, however the motivation of the MSE and
    crossentropy via the maximum likelihood is well known
    \autocite[e.g.][]{goodfellowDeepLearning2016}. The random
    linear model was already introduced in \textcite{benningRandomFunctionDescent2024}
    except for the analysis of the distribution of the noise.

    \item[\textbf{Chapter \ref{chap: optimization landscape of shallow neural networks}:}]
    Verbatim copy of \textcite{benningAlmostAllShallow2025} with shortened introduction.
\end{itemize}

 		\printbibliography[heading=secnotoc]
	\end{refsection}

	\chapter{Notation}

\begin{tabular}{p{0.15\linewidth} p{0.85\linewidth}}
    \([a,b)\) & closed-open interval \([a,b) = \{x \in \real : a\le x < b\}\)
    \\
    \([i\inter j)\) & closed-open discrete interval \([i \inter j) = [i,j) \cap \integer\) with integers \(\integer\)
    \\
    \(\fshift((x_i)_{i\in I})\) & the forward shift on sequences, \ie \(\fshift((x_i)_{i\in I}) = (x_{i+1})_{i\in I}\)
    \\
    \(\langle x,y\rangle\), \(x\cdot y\) & the inner product between \(x\) and \(y\)
    \\
    \(\|x\|\) & a norm, typically induced by the inner product \(\|x\|^2 = \langle x,x\rangle\)
    \\
    \(\sphere\) & the unit sphere in \(\real^\dims\) given by \(\sphere = \{x\in\real^\dims: \|x\|=1\}\)
    \\
    \(\nat\), \(\nat_0\) & the natural numbers starting at \(1\) and \(0\) respectively
    \\
    \(\ell^2\) & the hilbert space of sequences \(\ell^2=\{(x_i)_{i\in \nat}: \sum_{i=1}^\infty x_i^2<\infty\}\)
    \\
    \(\euclGroup\) & the group of (Euclidean) isometries, \ie \(\|\phi(x)-\phi(y)\| = \|x-y\|\) for any \(\phi\in E(\dims)\) 
    \\
    \(\orthGroup\) & the orthogonal group of linear isometries \(O(\dims) = \{\phi\in E(\dims): \phi \text{ linear}\}\) 
    \\
    \(\translGroup\) & the group of translations, \ie \(\phi(x) = x+v\) for some \(v\)
    \\
    \(\conj{z}\) & the conjugate \(\conj{z} = x - \im y\) of the complex number \(z=x + \im y\in \complex\)
    \\
    \(\identity\) & the identity (matrix)
    \\
    \(\ones\) & the function \(x\mapsto 1\) or  vector \((1,\dots, 1)\) in finite dimensional space
    \\
    \(\ind_A\) & the indicator function of the set \(A\), \ie \(\ind_A(x)\) is \(1\) if \(x\in A\)
    and \(0\) else
    \\
    \(\delta_x\) & the dirac distribution \(\delta_x(A) = \ind_A(x)\)
    \\
    \(\closure{A}\) & the closure of the set \(A\)
    \\
    \(\modCont_f(\delta)\) & the modulus of continuity of the function \(f\), see Definition \ref{def: modulus of continuity}
    \\
    \(C(\domain, \range)\) & the set of continuous functions with domain \(\domain\) and co-domain \(\range\)
    \\
    \((\Omega, \cA, \Pr)\) & the underlying probability space \(\Omega\), \(\sigma\)-algebra of measurable sets \(\cA\) and probability measure \(\Pr\)
    \\
    \(\borel(E)\) & the Borel sigma algebra induced by the open subsets of \(E\), assuming the topology is obvious
    \\
    \(\E\) & the expectation, i.e. \(\E[X] = \int X(\omega) \Pr(d\omega)\)
    \\
    \(\as\) & almost surely
    \\
    \(\iid\) & independent, identically distributed
    \\
    \(\normal(\mu, \sigma^2)\) & the normal distribution with mean \(\mu\in\real\) and variance \(\sigma^2\ge 0\)
    \\
    \(\normal(\mu, \Sigma)\) & the multivariate normal distribution with mean
    \(\mu\in \real^\dims\) and covariance matrix \(\Sigma \in \real^{\dims\times
    \dims}\)
    \\
    \(\normal(\mu, \C)\) & the distribution of a Gaussian random function with mean
    function \(\mu\) and covariance function \(\C\)
    \\
    \(\uniform(a,b)\) & the uniform distribution on the interval \((a,b)\)
    \\
    \(\mu_\rf\) & the mean function \(\mu_\rf(x) = \E[\rf(x)]\) of the random function \(\rf\)
    \\
    \(\C_\rf\) & the covariance function \(\C_\rf(x,y) = \Cov(\rf(x), \rf(y))\) of the random function \(\rf\)
\end{tabular}
 
	\mainmatter

	\begin{refsection}
		
\chapter{Worst case black box optimization}
\label{chap: worst case black box optimization}

As outlined in the preface, this chapter is concerned with the
hopelessness of worst case global optimization in high dimension. We start with
zero-order black box optimizers in Section~\ref{sec: worst case global opt is
hopeless} and proceed to explain why differentiability of the objective function
and higher order oracles change little about this fact in Section~\ref{sec:
differentiable objective function}. In Section~\ref{sec: differentiable
objective function} we further outline how classical optimization
theory solves this problem by forcing information into the gradient, at
the cost of assuming all critical points to be global minima in the process.
Since the assumption that `all critical points are global minima' is completely
implausible for many applications such as machine learning
\autocite[e.g.][]{dauphinIdentifyingAttackingSaddle2014} we will motivate
a distributional view of optimization, i.e.\ optimization on random functions
in Chapter~\ref{chap: distributional optimization}.

Before we do so -- to cover classical optimization theory `once and for
all' and to hint at a possible connection to random function optimization
theory -- there are two sections sprinkled into this chapter marked with a
star*. In these we present classical optimization results in 
terms of the modulus of continuity, which is more general than the common
Lipschitz assumptions.\footnote{
    While this generality is not commonly found in standard textbooks, we 
    do not believe that these modulus of continuity results are new.
} These result may become useful for the optimization of
random functions since there are entropy bounds on their modulus of continuity
\autocite[e.g.][]{adlerRandomFieldsGeometry2007}.

\section{High dimensional global optimization is hopeless}
\label{sec: worst case global opt is hopeless}

Without any assumptions about the objective function to optimize, it would be
necessary to try every possible input. This becomes impossible as soon as the
function has an infinite domain. A natural assumption is that the evaluation of
an objective function \(\obj\) at \(x\) is informative of \(\obj(x')\) for \(x'\) in a
neighborhood of \(x\). Intuitively this is the assumption of continuity, which classical
optimization theory always assumes in some form or another.\footnote{
    In a probabilistic framework with a random objective function \(\rf\) it is
    possible to capture the information provided by \(\rf(x)\) for \(\rf(x')\)
    with a covariance, or more generally a conditional distribution. This allows
    far more types of dependence than a simple ``similar if close-by'' relation
    that the continuity assumption represents.
}

A very intuitive continuity assumption is \(L\)-Lipschitz continuity, which limits
the rate of change of a function \(\obj\) to \(L\ge 0\), i.e.
\[
    |\obj(x) - \obj(y)| \le L \|x-y\|.
\]
A differentiable function \(\obj\) is in fact \(L\)-Lipschitz if and
only if its derivative is bounded by \(L\).

The constant \(L\) acts as a conversion factor of the distance
between \(x\) and \(y\) and their function values. Ensuring that all function
values \(\obj(y)\) of points \(y\) selected from a \(\delta\)-ball around \(x\) are
within an \(\epsilon\) tolerance of \(\obj(x)\) requires \(\delta\le
\frac{\epsilon}L\) which leads to
\[
    |\obj(x) - \obj(y)| \le L \|x-y\| \le L \delta \le \epsilon.
\]
\paragraph{Sufficient samples to find \(\epsilon\)-optimum}
To find the global optimum to \(\epsilon\)-precision it is thus sufficient to
cover the space with \(\delta = \frac{\epsilon}L\)-balls, sample from the
center of every ball and choose the largest (or respectively smallest function
value). The number of samples required is thus upper bounded by the
\(\delta\)-covering number.

\begin{definition}[Covering number]
    The \emph{\(\delta\)-covering number} of a space \(\domain\) is given by
    the smallest number of closed \(\delta\)-balls that cover \(\domain\)
    and is denoted by \(N_\domain(\delta)\).
\end{definition}

Considering a square \([0,1]^\dims\) it is intuitive that the covering number
should behave asymptotically like \(\delta^{-\dims}\). It is easiest to see this
using the sup-norm as the balls are then simply cubes with side length
\(\delta\). More generally the number \(\dims \ge 0\) which asymptotically satisfies
\[
    N_\domain(\delta) \sim C \delta^{-\dims}
\]
for some constant \(C\) and \(\delta\to 0\) is called the `box-counting dimension'
\autocite[e.g.][]{falconerFractalGeometryMathematical2003} and the
box-counting dimension of an \(m\)-dimensional submanifold of \(\real^n\) is
equal to \(m\) \autocite[Sec. 3.2]{falconerFractalGeometryMathematical2003}.

\marginskip{1em}
\begin{marginfigure}
    \includegraphics[width=\marginparwidth]{media/tikz/lipschitz_bump.tikz}
    \caption{
        The \(L\)-Lipschitz bump function \(\bump_\delta\).
    }
    \label{fig: L-lip bump}
\end{marginfigure}

\begin{marginfigure}
    \includegraphics[width=\marginparwidth]{media/tikz/zero_order_optimization.tikz}
    \caption{
        Sampling a \(4\times 4\) grid with \(15< 4\cdot 4\)
        points always leaves a tile untouched.
    }
    \label{fig: insufficient covering algo}
\end{marginfigure}
\paragraph{Necessary samples to find \(\epsilon\)-optimum} 
The inequalities used above can be made tight in the class of Lipschitz
functions by setting the slope of a bump function equal to \(L\) (cf.~Figure~\ref{fig:
L-lip bump}). Specifically, the bump function 
\[
    \bump_\delta(x) = L(\|x\|-\delta)\ind_{\|x\|\le \delta}
\]
is an \(L\)-Lipschitz function that is zero outside the \(\delta\)-ball at the
origin. To show that an algorithm cannot pick fewer than \(N_\domain(\delta)\)
evaluations we essentially place the optimum using the bump function \(g\)
in the last remaining \(\delta\)-ball (cf.~Figure~\ref{fig: insufficient
covering algo}).\footnote{
    In this picture we motivate the packing number (smallest number of disjoint
    balls, i.e. tiling), whereas in the proof we use the covering number
    again and draw balls around the points. This is because it is harder to
    visualize that there must be an uncovered point left-over than an empty
    field. While the packing number has the same asymptotic behavior \autocite{falconerFractalGeometryMathematical2003} tight
    results require the covering number.
}
To assume without loss of generality that the algorithm picks one of its
evaluation points as the decision point we gift it one additional point and
therefore only prove that there exists no algorithm which can pick \emph{fewer}
than \(N_\domain(\delta)-1\) evaluation points. We then consider what happens to
the algorithm when all evaluations are zero. This results in a sequence of
evaluation points (including the final decision) \(x_1,\dots
x_{N_\domain(\delta)-1}\) which we can interpret as centers of closed
\(\delta\)-balls. These balls cannot cover \(\domain\) as this
would be a contradiction to the definition of \(N_\domain(\delta)\). Thus
there exists a point \(x_*\) which is of at least \(\delta^+>\delta\) distance
to every evaluation point \(x_i\). Then the evaluation collected by the algorithm
are consistent with the bump function \(g_{\delta^+}(\cdot - x_*)\), which
is zero at distance \(\delta^+\) away from \(x_*\) and thus at the \(x_i\).
But the minimum of \(\bump_\delta^+\) is at \(-\epsilon^+\) with
\[
    \epsilon^+ := L\delta^+ > L\delta = \epsilon
\]
and the error of the algorithm on \(g_{\delta^+}\) is therefore greater than
\(\epsilon\).

Put together we have proven the following result
\begin{prop}
    \label{prop: optimal evaluation complexity}
    The optimal evaluation complexity \(\mathrm{EvalComp}^*(\epsilon)\) to
    guarantee an \(\epsilon\)-optimal function value in the space of
    \(L\)-Lipschitz functions with a zero order black box algorithm on the
    domain \(\domain\) is bounded by
    \[
        N_\domain(\tfrac{\epsilon}L) -1
        \le \mathrm{EvalComp}^*(\epsilon)
        \le N_\domain(\tfrac{\epsilon}L).
    \]
    In particular, we have
    \[
        \mathrm{EvalComp}^*(\epsilon)
        \sim C(\tfrac{\epsilon}L)^{-\dims}
    \]
    for some constant \(C>0\) and box-counting dimension \(\dims\) of \(\domain\).
    Moreover this complexity is achieved by a naive `grid search' algorithm,
    which covers the space by \(N_\domain(\frac{\epsilon}L)\) balls, choosing the
    centers as its evaluation points and selecting the point with the largest
    (or respectively smallest) function value.
\end{prop}

Not only does this imply that the optimal algorithm is essentially a dumb type
of grid search, but since the `Box-counting dimension' \(\dims\) of \(\domain\)
coincides with the intuitive dimension on Manifolds and in particular
\(\real^\dims\) this implies exponential complexity in the ambient dimension.
From this perspective black box optimization is essentially doomed in high
dimension \(\dims\).

\subsection{Modulus of continuity: Grid search*}

Not every continuous function is Lipschitz continuous and one may be curious,
whether or not it is possible to generalize Prop.~\ref{prop: characterization
modulus of continuity}. Indeed, with the modulus of continuity it is possible
to quantify continuity more generally.
\begin{definition}[Modulus of continuity]
    \label{def: modulus of continuity}
    The \emph{modulus of continuity of \(\obj\)} is defined by
    \[
        \modCont_\obj(\delta) := \sup\{|\obj(x)- \obj(y)| : x,y\in \domain,\; \metric(x,y) \le \delta\},
    \]
    where \(\metric\) is a metric on the domain \(\domain\) of \(\obj\).
\end{definition}

Functions which can be moduli of continuity can be characterized as follows.\footnote{
    e.g. \textcite{efimovContinuityModulus}, see Section \ref{sec: technical
    appendix} for a proof.
}.

\begin{restatable}[Characterization of moduli of continuity]{prop}{charactModCont}
    \label{prop: characterization modulus of continuity}
    For a function \(\modCont:[0,\infty)\to [0,\infty)\)
    the following are equivalent:
    \begin{enumerate}[label=(\roman*)]
        \item\label{it-mod-cont: exists func}
        There exists a uniformly continuous function \(\obj: \domain \to
        \real\), where \(\domain\) is a convex subset of a normed vectorspace,
        such that \(\modCont\) is the modulus of continuity of \(\obj\), i.e.
        \(\modCont = \modCont_\obj\)
        
        \item\label{it-mod-cont: properties}
        \begin{enumerate}
            \item \(\modCont(0) = 0\),
            \item \(\modCont\) is non-decreasing,
            \item \(\modCont\) is continuous in \(0\),
            \item \(\modCont\) is subadditive, i.e. \(\modCont(\eta+\delta) \le \modCont(\eta) + \modCont(\delta)\)
        \end{enumerate}
    \end{enumerate}
    If these statements are true, \(\modCont\) is continuous, its modulus of
    continuity is itself, i.e. for all \(\eta,\delta \ge 0\)
    \[
        |\modCont(\eta) - \modCont(\delta)| \le \modCont(|\eta-\delta|),
    \]
    and we call \(\modCont\) ``a modulus of continuity''.
\end{restatable}

With moduli of continuity characterized it is possible to define more general
continuity classes.

\begin{definition}[Continuity class]
    Functions which are dominated by some modulus of continuity \(\modCont\), i.e.
    \[
        \cF_\modCont := \{ \obj\colon \domain \to \real : \modCont_\obj(\delta) \le \modCont(\delta)\}
    \]
    are the continuity class associated to \(\modCont\).
\end{definition}

Of these classes, Lipschitz continuos functions are simply an example.

\begin{example}
    \(\cF_\modCont\) with \(\modCont(\delta) = L \delta\) is exactly the space
    of \(L\)-Lipschitz continuous functions, since for Lipschitz
    functions we have
    \[
        \modCont_\obj(\delta)
        \equiv \sup_{d(x,y)\le \delta} |\obj(x) - \obj(y)|
        \le \sup_{d(x,y)\le \delta}L d(x,y)
        = L\delta.
    \]
\end{example}
Recall that the conversion from \(\epsilon\) to \(\delta\) was previously given
by \(\delta=\frac{\epsilon}L\).  More generally the conversion from function
value scale to domain scale is given by \(\delta=\modCont^{-1}(\epsilon)\), which
necessitates the introduction of the pseudo inverse
\[
    \modCont^{-1}(\epsilon) := \sup\{\delta \ge 0: \modCont(\delta) \le \epsilon\}
\]
to cover the non-bijective case.

For completeness sake let us properly define ``grid-search''.

\begin{definition}[Generalized \((\epsilon, \modCont)\)-grid search]
    Let \(\delta=\modCont^{-1}(\epsilon)\) and let \(x_i\in \domain\) be selected
    such that the balls
    \(\ball_\delta(x_1),\dots,\ball_{\delta}(x_{N_\domain(\delta)})\) cover
    \(\domain\) (which is possible by the definition of the covering number
    \(N_\domain(\delta)\)). Then evaluate \(\obj:\domain\to\real\) at every \(x_i\) and return
    \[
        \hat{x} = \argmin_{i=1,\dots, N(\delta)} \obj(x_i).
    \]
\end{definition}

With grid search properly defined it is possible to generalize Prop.~\ref{prop: optimal evaluation complexity}.

\begin{theorem}[Grid search has optimal worst-case complexity]
    Generalized \((\epsilon, \modCont)\)-grid search is guaranteed to find an
    \(\epsilon\)-minimal function value in the class \(\cF_\modCont\) using
    \(N_\domain(\modCont^{-1}(\epsilon))\) function evaluations. The optimal
    worst-case evaluation complexity \(\mathrm{EvalComp}^*(\epsilon)\)
    is furthermore bounded by
    \[
       N_\domain(\modCont^{-1}(\epsilon)) -1  
       \le \mathrm{EvalComp}^*(\epsilon)
       \le N_\domain(\modCont^{-1}(\epsilon))
    \]
    For the box-counting dimension \(\dims\) of the domain \(\domain\) there
    thus exists a constant \(C\) such that
    \[
        \mathrm{EvalComp}^*(\epsilon) \sim C [\modCont^{-1}(\epsilon)]^\dims.
    \]
\end{theorem}

\begin{proof}
    Let \(\delta=\modCont^{-1}(\epsilon)\) and \(x_1,\dots, x_{N_\domain(\delta)}\)
    be the points chosen by grid search such that that their \(\delta\)
    neighborhoods are a cover \(\domain\).
    For any \(x\in \domain\) there exists \(x_j\) such that their distance is smaller
    than \(\delta=\modCont^{-1}(\epsilon)\) which implies
    \[
        \obj(x_j) - \obj(x) \le |\obj(x_j) - \obj(x)|  \le \modCont_\obj(\delta) \le \modCont(\delta) \le \epsilon
    \]
    Then we have for the point \(\hat{x}\) selected by grid search
    \[
        \obj(\hat{x}) = \argmin_{i} \obj(x_i) \le \obj(x_j) \le \obj(x) + \epsilon.
    \]
    As \(x\) was arbitrary we have \(\epsilon\)-optimality, i.e.
    \[
        \obj(\hat{x}) \le \inf_x \obj(x) + \epsilon.
    \]
    
    \textbf{Optimality:} We now want to show that there exists no algorithm
    which requires fewer than \(N(\modCont^{-1}(\epsilon))-1\) oracle evaluations.
    Assume that some algorithm would start with \(x_1\) and choose \(x_i\)
    if \(\obj(x_j)=0\) for all \(j<i\) and would terminate after \(n\) steps with
    \(n<N(\modCont^{-1}(\epsilon)) -1\). It then chooses point \(\hat{x}\).
    We are now going to gift this algorithm one last evaluation (setting
    \(x_{n+1} = \hat{x}\)) and can thus assume that the algorithm picks
    \[
        \hat{x} = \argmin_{i=1,\dots, n+1} \obj(x_i)
    \]
    without loss of generality. Furthermore \(n+1<N(\modCont^{-1}(\epsilon))\) so
    so the \(\ball_\delta(x_i)\) with \(\delta=\modCont^{-1}(\epsilon)\) do not cover
    \(\domain\). Therefore there exists \(x^*\in \domain\) with
    \[
       \delta^+ := \min_{i=1,\dots, n+1}\metric(x^*,x_i)> \delta,
    \]
    Since \(\delta=\modCont^{-1}(\epsilon)\) we have by definition of
    \(\modCont^{-1}\) and \(\delta^+\ge \delta\) that \(\modCont(\delta^+) > \epsilon\).
    We now define
    \[
        \obj(x) := \min\bigl\{\modCont(\metric(x,x^*)) - \modCont(\delta^+), 0\bigr\}
    \]
    Since \(\metric(x_i,x^*)\ge \delta^+\) we have by the non-decreasing property
    of \(\modCont\) that
    \[
        \modCont(\metric(x_i, x^*)) \ge \modCont(\delta^+)
    \]
    and thus \(\obj(x_i) = 0\). So the algorithm will select these \(x_i\) but we
    have
    \[
        \obj(x_*) = -\modCont(\delta^+) < - \epsilon.
    \]
    This implies
    \[
        \argmin_i \obj(x_i)
        = 0 > \inf_x \obj(x) + \epsilon.
    \]
    This would be a contradiction to \(\epsilon\) optimality if \(f\in
    \cF_\modCont\). To show \(f\in \cF_\modCont\), consider that \(x\mapsto
    \min\{x,0\}\) is \(1\)-Lipschitz resulting in
    \[
        |\obj(x) - \obj(y)|
        \le |\modCont(\metric(x,x^*)) - \modCont(\metric(y, x^*))|
        \overset{(*)}\le \modCont(\metric(x,y)).
    \]
    Where \((*)\) is the reverse triangle inequality for \(\modCont\circ\metric\),
    which is a semi-metric because it is positive, symmetric and satisfies the
    triangle inequality since \(\modCont\) is increasing and sub-additive as
    a modulus of continuity (Prop.~\ref{prop: characterization modulus of
    continuity}), so we have
    \[
        \modCont(\metric(x,y))
        \overset{\text{increasing}}\le \modCont(\metric(x, z) + \metric(z, y))
        \overset{\text{subadditive}}\le \modCont(\metric(x,z)) + \modCont(\metric(z, y)).
    \]
    Noting that \(x\) and \(y\) were arbitrary, we thus conclude
    \[
        \modCont_\obj(\delta)
        = \sup_{\metric(x,y)\le \delta} |\obj(x) - \obj(y)|
        \le \sup_{\metric(x,y)\le \delta}\modCont(\metric(x,y))
        \le \modCont(\delta).
        \qedhere
    \]
\end{proof}

\section{Differentiable objective functions}
\label{sec: differentiable objective function}

Assume the objective function \(\obj\) is differentiable with continuous
gradients to ensure they are useful. It is again instructive to start with
the \(L\)-Lipschitz continuity condition for gradients, i.e.
\begin{equation}
    \label{eq: L-smoothness, def} 
    \|\nabla \obj(x) - \nabla \obj(y) \| \le L \|x-y\|.
\end{equation}
The objective function \(\obj\) is then also called `\(L\)-smooth'.
The important question is whether this restricted function class allows for
speed-ups in optimization. In particular, does the use of gradient evaluations,
i.e. first order optimization, help?
\begin{marginfigure}
    \includegraphics[width=\marginparwidth]{media/tikz/smooth_bump.tikz}
    \includegraphics[width=\marginparwidth]{media/tikz/smooth_bump_diff.tikz}
    \caption{
        \(L\)-smooth bump function \(\bump_\delta\)
    }
    \label{fig: L-smooth bump function}
\end{marginfigure}
For intuition consider the \(L\)-smooth bump function
(cf.~Figure~\ref{fig: L-smooth bump function})\footnote{
    The proof of \(L\)-smoothness is more annoying than one might think based on
    the pictures (details in Lemma~\ref{lem: L smooth bump}).
}
\[
    \bump_\delta(x) = \begin{cases}
        \frac{L}4[2\|x\|^2 -\delta^2] & \|x\|\le \frac{\delta}2
        \\
        -\frac{L}2(\delta - \|x\|)^2 & \frac{\delta}2\le\|x\| \le \delta
        \\
        0 & \|x\| \ge \delta.
    \end{cases}
\]
While it may be no longer possible to reach the height \(L\delta\), it
is possible to reach height \(\epsilon\) within a larger \(\delta\)-ball,
specifically \(\delta = 2\sqrt{\epsilon/L}\). With the same covering argument
as before it can be shown that any algorithm which evaluates less than
\(N_\domain(\delta)\) places leaves a \(\delta^+\)-ball open where the
\(\bump_\delta\) bump function can be placed. This implies that at every
evaluation location not only the function value but also its gradients are zero.
Gradient evaluations therefore do not have any benefit at all. By the same
covering argument one can obtain 
\[
    N_\domain(2\sqrt{\epsilon/L}) -1 \le \mathrm{EvalComp}^*(\epsilon) \le N_\domain(2\sqrt{\epsilon/L})
\]
and for some constant \(C\ge 0\) and the box-counting dimension \(\dims\)
this implies
\[
    \mathrm{EvalComp}^*(\epsilon)
    \sim C (\epsilon/L)^{\dims/2}.
\]
While \(L\)-smoothness halves the dimension it therefore does not fix the issue of
exponential complexity in the dimension. Furthermore the gradient evaluations
are seemingly useless as zero-order grid search is still optimal.
More generally, global optimization has worst-case sample complexity of the order
\((\epsilon/L)^{-\dims/r}\) to get within \(\epsilon\)-distance of the minimum
of a function in \(C^r([0,1]^\dims)\) with \(r\)-th derivatives bounded by \(L\)
\autocite{wasilkowskiAverageComplexityGlobal1992}.

\subsection*{Why are gradient evaluations seemingly useless?}

Gradients are useful because they point towards a reasonable search direction,
\emph{assuming they are non-zero}. In high dimension
directional information is invaluable since there are many directions to choose
from. But the adversarially chosen bump function does not grant this
information since gradients are zero and thus directionless. Similarly, if you
had lower and higher function evaluations then, due to continuity, the
vicinities of lower function values are more likely to contain the global
minimum.  This could inform search decisions. But the adversarial bump function withholds
any information by being completely flat.

However, if information is available an algorithm ought to use it and grid
search decidedly does not. Its optimality therefore hinges on the fact that in
the worst case there is no information to use. Classical optimization theory
therefore forces information upon the function by assuming the gradient cannot
be zero outside the global minimum. One such example is the Polyak-Łojasiewicz
(PL) inequality \autocite[e.g.\label{footnote: karimi}][]{karimiLinearConvergenceGradient2016}, which
demands
\begin{equation}
    \label{eq: PL-inequality}    
    \mu (\obj(x) - \min_x \obj(x))
    \le \tfrac12\|\nabla \obj(x)\|^2
\end{equation}
for some \(\mu > 0\). And this assumption is sufficient to obtain exponential
convergence of the objective value to the global optimum with the additional
assumption of \(L\)-smoothness.\footref{footnote: karimi} Of course this
assumption also implies that all critical points are global minima, as does any
weaker assumption that forces the gradient to be non-zero outside
of a global minimum. This implies that it is sufficient to find \emph{any}
critical point, which was not the case in the adversarial example where every
evaluation point was critical.

\begin{remark}[Convexity]
    The assumption of strong convexity, i.e.\
    \[
        \obj(y) \ge \obj(x) + \langle \nabla \obj(x), y-x\rangle + \tfrac\mu2 \|y-x\|^2,
    \]
    implies a lower bound\footnote{
        Due to \(\obj(x) - \obj(x_*) \ge 0\), reordering strong convexity with
        \(y=x_*\) and the Cauchy-Schwarz inequality imply
        \begin{align*}
            \tfrac\mu2 \|x_*-x\|^2
            &\le \langle \obj(x), x_* - x\rangle
            \\
            &\le \|\nabla \obj(x)\| \|x_*-x\|.
        \end{align*} Division by \(\|x_* - x\|\) yields the claim
        \autocite[this proof is from][]{karimiLinearConvergenceGradient2016}.
    } of the gradient in terms of the distance to the unique minimizer \(x_*\)
    (that exists by convexity)
    \begin{equation}
        \label{eq: error bound}
        \mu\|x_*-x\|
        \le \|\nabla \obj(x)\|.
    \end{equation}
    For \(L\)-smooth functions \eqref{eq: L-smoothness, def}
    this is stronger than the PL-inequality, since\footnote{
        The first inequality follows from Lemma \ref{lem: taylor approximation with
        trust bound} for \(L\)-smooth \(\obj\) and \(\nabla\obj(x_*) = 0\).
    }
    \[
        \obj(x) - \obj(x_*)
        \overset{L\text{-smooth}}\le \tfrac{L}2\|x-x_*\|^2
        \overset{\eqref{eq: error bound}}\le \tfrac{L}{2\mu} \|\nabla\obj(x)\|^2.
    \]
    As minima have locally convex neighborhoods, not much generality is lost to
    assume convexity at this point, since the convex neighborhoods determine
    asymptotic behavior. And convex function theory is fully developed with
    tight lower bounds and convergence algorithms that achieve these lower
    bounds
    \autocite[e.g.][]{nesterovLecturesConvexOptimization2018,garrigosHandbookConvergenceTheorems2023,gohWhyMomentumReally2017,bottouOptimizationMethodsLargeScale2018}.
\end{remark}

\subsection{Modulus of continuity: Descent lemma and convergence*}

In this section we want to briefly mention that the ``Decent Lemma'' that
drives most convergence proofs in classical optimization theory can also
be generalized to the modulus of continuity. Accordingly, define the
modulus of continuity of the gradient to be
\[
    \modCont_{\nabla \obj}(\delta)
    := \sup\{ \|\nabla \obj(x) - \nabla\obj(y)\| : \|x-y\| \le \delta \}.
\]
Recall that \(L\)-smoothness is defined by \(L\)-Lipschitz gradients, i.e.
the special case
\begin{equation}
    \label{eq: L smoothness}
    \modCont_{\nabla\obj}(\delta) = L\delta.
\end{equation}
With the help of this modulus of continuity it is possible to generalize
the `Descent Lemma' \autocite[e.g.][Lemma~2.28]{garrigosHandbookConvergenceTheorems2023}.

\begin{lemma}[General Descent Lemma]
    \label{lem: general descent lemma}
    Steps in the direction of the gradient result in 
    \[
       \obj\Bigl(x - \stepsize \tfrac{\nabla\obj(x)}{\|\nabla\obj(x)\|}\Bigr) 
       \le \obj(x) - \int_0^{\stepsize} \|\nabla\obj(x)\| - \modCont_{\nabla\obj}(t) dt,
    \]
    where the minimizing step size is given by \(\stepsize^* =
    \modCont_{\nabla\obj}^{-1}(\|\nabla\obj(x)\|)\) and
    if \(\modCont_{\nabla\obj}\) is differentiable
    \[
       \obj\Bigl(x - \stepsize^* \tfrac{\nabla\obj(x)}{\|\nabla\obj(x)\|}\Bigr) 
       \le \obj(x) - \underbrace{
        \Bigl(\int_0^1 \frac{1-t}{\modCont_{\nabla\obj}'(t\|\nabla\obj(x)\|)} dt\Bigr)}_{
            \overset{L\text{\normalfont-smooth}}= \frac1{2L}
        } \|\nabla\obj(x)\|^2.
    \]
\end{lemma}

\begin{remark}[Step size vs. learning rate]
    \label{rem: step size vs. learning rate}
    In the case of \(L\)-smooth objective functions \eqref{eq: L smoothness}, the optimal step size is given by
    \[
        \stepsize^* = \modCont_{\nabla\obj}^{-1}(\|\nabla\obj(x)\|) =  \tfrac{\|\nabla\obj(x)\|}{L}.
    \]
    The resulting optimization step is thus characterized by
    \[
        x_{n+1}
        = x_n - \underbrace{\stepsize^*}_{\mathclap{\text{`stepsize'}}} \tfrac{\nabla\obj(x_n)}{\|\nabla \obj(x_n)\|}
        = x_n - \underbrace{\lr^*}_{\mathclap{\text{`learning rate'}}} \nabla\obj(x_n) 
    \]
    with `learning rate' \(\lr^* = \frac1L\). Classical optimization theory typically
    considers these learning rates associated to the un-normalized gradient,
    while step sizes appear to be more natural in the general theory using moduli of
    continuity. Step sizes are also more natural in the theory of
    optimization on random functions as will become apparent in the chapter on
    Random Function Descent (Chapter~\ref{chap: random function descent})
\end{remark}

We explain the significance of the Descent Lemma by proving that
convergence proofs can be easily obtained by its use. We
provide two examples, one proves convergence and the other
convergence rates.
\begin{theorem}[Convergence]
    Assume \(\obj\) is \(L\)-smooth, \(\inf_x \obj(x) > -\infty\)
    and \(x_{n+1} = x_n - \frac1L \nabla\obj(x_n)\).
    Then
    \[
        \lim_{n\to\infty}\|\nabla \obj(x_n)\| = 0.
    \]
\end{theorem}
\begin{proof}
    By Lemma~\ref{lem: general descent lemma} we have
    \[
        \obj(x_{n+1}) - \obj(x_n) \le - \tfrac1{2L}\|\nabla \obj(x_n)\|^2
    \]
    and therefore
    \[
        \sum_{k=0}^n \|\nabla \obj(x_k)\|^2
        \le 2L \sum_{k=0}^n (\obj(x_k) - \obj(x_{k+1}))
        \overset{\text{telescope}}\le 2L\bigl(\obj(x_0) - \inf_x\obj(x)\bigr).
    \]
    This bound is finite by assumption and holds uniformly over \(n\), thus the
    squared gradient norms are infinitely summable and thereby converge to
    \(0\).
\end{proof}

It is clear that if the gradients converge to zero and the
only place where the gradient may be zero is in global minima,
then the function values should also converge to the minimum.
The PL-inequality \eqref{eq: PL-inequality} ensures this is the case
and indeed, the following result provides convergence rates for the case of
\(L\)-smooth functions which satisfy the PL-inequality \eqref{eq:
PL-inequality}.

\begin{theorem}[Convergence rates]
    Let \(\obj\) be \(L\)-smooth, \(\obj_* := \inf_x \obj(x) >-\infty\),
    \(\obj\) satisfies the PL-inequality \eqref{eq:
    PL-inequality} and
    \[
        x_{n+1} = x_n - \tfrac1L \nabla \obj(x_n),
    \]
    then we have exponential convergence\footnote{typically referred to as `linear convergence' in the optimization literature due to the linear reduction in every step.}
    \[
        \obj(x_n) - \obj_* \le \bigl(1-\tfrac{\mu}{L}\bigr)^n \bigl(\obj(x_0) - \obj_*\bigr).
    \]
\end{theorem}
\begin{proof}
    Since \(\frac{\|\nabla\obj(x_n)\|}{L}\) is the minimizing step size of the
    upper bound we have by the Descent Lemma (Lemma~\ref{lem: general descent lemma})
    \begin{align*}
        \obj(x_{n+1}) - \obj_* 
        \overset{\text{Lem.~\ref{lem: general descent lemma}}}&\le \obj(x_n) - \obj_* - \tfrac1{2L}\|\nabla \obj(x_n)\|^2
        \\
        \overset{\eqref{eq: PL-inequality}}&\le
        (1- \tfrac{\mu}{L})(\obj(x_n) - \obj_*).
    \end{align*}
    Here we used the the PL-inequality \eqref{eq: PL-inequality}, i.e.
    \(\mu(\obj(x_n) - \obj_*) \le \tfrac12 \|\nabla \obj(x)\|\), in the last
    line. Induction over \(n\) finishes the proof.
\end{proof}

Clearly, weaker assumptions about continuity weaken the first inequality
provided by the Descent Lemma (Lemma~\ref{lem: general descent lemma}). 
Similarly, lower bounds on the gradient that are weaker than the
PL-inequality will weaken the second inequality. 
In essence convergence can be proven as long as the combination of these assumptions
results in a sufficient reduction factor.

In the remainder of this section we prove the Descent Lemma. The main vehicle of
the proof of Lemma~\ref{lem: general descent lemma} is the following
bound based on the Taylor approximation.

\begin{lemma}[Taylor approximation with trust bound]
    \label{lem: taylor approximation with trust bound}
    Let \(\obj\) be a continuously differentiable function whose gradient
    has modulus of continuity \(\modCont_{\nabla\obj}\). Then
    \[
        \obj(y)
        \le \obj(x)
        + \langle \nabla\obj(x), y-x\rangle
        + \underbrace{\int_0^{\|y-x\|} \modCont_{\nabla\obj}(t)dt}_{\overset{L\text{-smooth}}=\frac{L}2 \|y-x\|^2}.
    \]
\end{lemma}
\begin{proof}
    By the fundamental theorem of calculus (FTC), the Cauchy-Schwarz inequality (CS) and \(\step = y-x\)
    we have
    \begin{align*}
        \obj(y)
        &= \obj(x) + (\obj(y) - \obj(x))
        \\
        \overset{\text{FTC}}&= \obj(x) + \int_0^1 \bigl\langle \nabla\obj\bigl(x + \lambda\step\bigr),\step\bigr\rangle d\lambda
        \\
        &= \obj(x) + \langle \nabla\obj(x), y-x\rangle
        + \int_0^1 \underbrace{
            \bigl\langle \nabla\obj\bigl(x + \lambda\step\bigr) - \nabla\obj(x),\step\bigr\rangle
        }_{
            \overset{\text{CS}}\le \|\nabla\obj(x+\lambda\step) - \nabla\obj(x)\| \|\step\|
        } d\lambda
        \\
        &\le \obj(x) + \langle \nabla\obj(x), y-x\rangle
        + \int_0^1 \modCont_{\nabla\obj}(\lambda\|\step\|) \|\step\| d\lambda
        \\
        &= \obj(x) + \langle \nabla\obj(x), y-x\rangle
        + \int_0^{\|\step\|} \modCont_{\nabla\obj}(t) dt.
        \qedhere
    \end{align*}
\end{proof}

\begin{proof}[Proof of Descent Lemma (\ref{lem: general descent lemma})]
    Simply plug \(y=x - \stepsize \frac{\nabla\obj(x)}{\|\nabla\obj(x)\|}\) into
    Lemma \ref{lem: taylor approximation with trust bound} to obtain\footnote{
        Note that the gradient direction minimizes the bound of Lemma~\ref{lem:
        taylor approximation with trust bound}. This is because minimization of
        the bound
        \[
            \min_{\step : \|\step\|=\stepsize}
            \obj(x) + \langle \nabla\obj(x), \step\rangle + B(\|\step\|),
        \]
        with \(B(\stepsize) = \int_0^\stepsize \modCont_{\nabla \obj}(t)dt\),
        only depends on the direction \(\step = y-x\) in the inner product,
        which forces the gradient direction. The proof is a straightforward
        application of the Cauchy-Schwarz inequality (see Corollary \ref{cor:
        constrained minimization of scalar products}).
    }
    \begin{align}
        \nonumber
        \obj(y) &\le \obj(x)
        + \underbrace{\Bigl\langle \nabla\obj(x), -\stepsize\frac{\nabla\obj(x)}{\|\nabla\obj(x)\|}\Bigr\rangle}_{=\stepsize\|\nabla\obj(x)\|}
        + \int_0^\stepsize \modCont_{\nabla\obj}(t)dt
        \\
        \label{eq: descent bound}
        &= \obj(x) - \int_0^\stepsize \|\nabla\obj(x)\| -  \modCont_{\nabla\obj}(t)dt.
    \end{align}
    What remains is to find the optimal step size.
    Since \(\modCont_{\nabla\obj}\) is monotonously increasing, the
    integrand \eqref{eq: descent bound} is positive until
    \[
        \stepsize^*=\modCont_{\nabla\obj}^{-1}(\|\nabla\obj(x)\|)
    \]
    and negative afterwards. This is therefore the minimizing step size.
    Plugging this result into \eqref{eq: descent bound}
    we obtain
    \begin{align*}
        \obj(x - \stepsize^* \tfrac{\nabla\obj(x)}{\|\nabla\obj(x)\|})
        &\le \obj(x) - \int_0^{\stepsize^*}  \|\nabla\obj(x)\|- \modCont_{\nabla\obj}(t)dt
        \\
        &= \obj(x) - \int_0^{1} \frac{1 - t}{\modCont_{\nabla\obj}'(t\|\nabla\obj(x)\|)}ds \|\nabla\obj(x)\|^2.
    \end{align*}
    with the substitution
    \(s = \frac{\modCont_{\nabla\obj}(t)}{\|\nabla\obj(x)\|}\).
\end{proof}

\begin{subappendices}
   
\section{Technical Appendix}\label{sec: technical appendix}

\begin{lemma}[\(L\)-smooth bump function]
    \label{lem: L smooth bump}
    For any \(\delta,L>0\) we have that the function \(g\) is \(L\)-smooth,
    where \(g\) is defined by
    \[
        g:\begin{cases}
            \real^\dims \to \real
            \\
            x\mapsto h(\|x\|)
        \end{cases}
        \quad 
        h(r): \begin{cases}
            [0,\infty) \to \real
            \\
            r\mapsto 
            \begin{cases}
                \frac{L}4[2r^2 - \delta^2] & r\le \frac{\delta}2
                \\
                -\frac{L}2(\delta - r)^2 & \frac{\delta}2\le r \le \delta
                \\
                0 & r \ge \delta.
            \end{cases}
        \end{cases}
    \]
\end{lemma}
\begin{proof}
    It is straightforward to confirm (cf.~Figure~\ref{fig: L-smooth bump function})
    \[
        |\tfrac{h'(r)}{r}|\le L
        \quad \text{and}\quad |h''(r)|\le L.
    \]
    for all \(r \notin \{0, \frac{\delta}2, \delta\}\).  Furthermore, the
    Hessian of \(g\) is (for all \(\|x\| \notin \{0, \frac{\delta}2, \delta\}\))
    given by
    \[
        \nabla^2 g(x) = h''(\|x\|) \tfrac{xx^\transpose}{\|x\|^2} + \tfrac{h'(\|x\|)}{\|x\|} \Bigl( \identity - \tfrac{xx^\transpose}{\|x\|^2}\Bigr)
    \]
    We can decompose any vector \(v\) of unit length (\(\|v\| = 1\)) into the orthogonal components
    \[
        v = \alpha \tfrac{x}{\|x\|} + v^\perp,
    \]
    which implies
    \[
        |v^\transpose \nabla^2 g(x) v|
        = \bigl|h''(\|x\|) \alpha^2 + \tfrac{h'(\|x\|)}{\|x\|}(1-\alpha^2)\bigr|
        \le L.
    \]
    Thus for \(\|x\| \notin \{0, \frac{\delta}2, \delta\}\) the operator norm of
    \(\nabla^2 g(x)\) is bounded by \(L\). Since \(p(t):= \nabla g(x + t(y-x))\)
    is continuous and \emph{piecewise} differentiable except for a finite number of points,
    we can apply the fundamental theorem of calculus to obtain
    \begin{align*}
        \|\nabla g(y) - \nabla g(y)\|
        &= \|p(1) - p(0)\|
        \le \int_0^1 \underbrace{\|\nabla^2 g(x + t(y-x))\|}_{\le L} \|y-x\| dt
        \\
        &\le L \|y-x\|.
    \end{align*}
    With the understanding that the integral would formally have to be broken
    into pieces.
\end{proof}

We finally provide a proof for the well known characterization of the modulus of continuity
\autocite[e.g.][]{efimovContinuityModulus}.

\charactModCont*

\begin{proof}
    \textbf{\ref{it-mod-cont: exists func} \(\Rightarrow\) \ref{it-mod-cont: properties}:}
    \(\modCont(0)=0\) is follows immediately from the definition like the fact
    that it is non-decreasing. Continuity in \(0\) is necessitates by the fact that \(\obj\) is
    uniformly continuous. We prove sub-additivity in two steps: 
    \begin{enumerate}[wide,label={\textbf{Step \arabic*:}}]
        \item First we prove that whenever $\|x-y\|\le \eta+\delta$, there
        exists $z\in \domain$ such that $\|x-z\|\le \eta$ and $\|z-y\|\le \delta$.

        For this we define $\lambda= \frac{\delta}{\eta+\delta}$ and select
        \(z = \lambda x + (1-\lambda) y\) (which we are allowed to do by convexity). Then
        \[
        \|x-z\| = \|(1-\lambda)(x-y)\| = (1-\lambda) \|x-y\| = \frac{\eta}{\eta+\delta} \|x-y\| \le \eta.
        \]
        And similarly for the other equation.

        \item
        Now we simply consider the definition
        \begin{align*}
            \modCont_\obj(\eta+\delta)
            &= \sup \bigl\{ |\obj(x)- \obj(y)| : x,y\in \domain, \|x-y\|<\eta + \delta\bigr\}
            \\
            &\le \sup\begin{aligned}[t]
                \Bigl\{ 
                &|\obj(x)- \obj(z)| + |\obj(z)- \obj(y)|:
                \\
                &\qquad x,y\in \domain, \|x-y\|<\eta + \delta, 
                z = \tfrac{\delta}{\eta+\delta} x + \tfrac{\eta}{\eta+\delta}y \Bigr\}
            \end{aligned}
            \\
            &\le \begin{aligned}[t]
                &\sup\Bigl\{ |\obj(x)- \obj(z)| : x,z\in \domain, \|x-z\|\le \eta\Bigr\}
                \\ 
                &+ \sup\Bigl\{|\obj(z) - \obj(y)| : y,z \in \domain, \|y-z\|\le \delta\Bigr\}
            \end{aligned}
            \\
            &= \modCont_\obj(\eta) + \modCont_\obj(\delta).
        \end{align*}
        Where we use that giving more choice only increases the supremum.
    \end{enumerate}
    \textbf{\ref{it-mod-cont: properties} \(\Rightarrow\) \ref{it-mod-cont: exists func}:}
    By proving the modulus of continuity has itself as its modulus of continuity, we prove existence
    and the additional claim. By monotonicity and sub-additivity of the modulus
    of continuity for \(\eta, \delta \ge 0\)
    \[
        \modCont(\eta)
        \overset{\text{monotonicity}}\le
        \modCont(| \eta-\delta|  + |\delta|)
        \overset{\text{sub-additivity}}\le
        \modCont(|\eta - \delta|) + \modCont(\delta)
    \]
    So by symmetry of \(\delta\) and \(\eta\)
    \[
        |\modCont(\eta) - \modCont(\delta)| \le \modCont(|\eta -\delta|).
        \qedhere
    \]
\end{proof}

    \subsection{Constrained linear optimization}

Let \(U\) be a vectorspace. We define the projection of a vector \(w\) onto \(U\) by
\[
	\proj{U}(w) := \argmin_{v\in U} \|v-w\|^2
\]
\begin{lemma}[Constrained maximiziation of scalar products]
	\label{lem: constrained maximiziation of scalar products}
	For a linear subspace \(U\subseteq \real^\dims\), we have	
	\begin{align}
		\max_{\substack{v\in U \\ \|v\|=\lambda}}\langle v, w\rangle 
		&= \lambda \|\proj{U}(w)\|
		\\
		\argmax_{\substack{v\in U\\ \|v\|=\lambda}} \langle v, w\rangle
		&= \lambda \frac{\proj{U}(w)}{\|\proj{U}(w)\|}
	\end{align}
\end{lemma}
Before we get to the proof let us note that this immediately results in the following
corollary about minimization.
\begin{corollary}[Constrained minimization of scalar products]
	\label{cor: constrained minimization of scalar products}
	\begin{align}
		\min_{\substack{v\in U \\ \|v\|=\lambda}}\langle v, w\rangle 
		&= -\lambda \|\proj{U}(w)\|
		\\
		\argmin_{\substack{v\in U\\ \|v\|=\lambda}} \langle v, w\rangle
		&= -\lambda \frac{\proj{U}(w)}{\|\proj{U}(w)\|}
	\end{align}
\end{corollary}
\begin{proof}[Proof of Corollary~\ref{cor: constrained minimization of scalar products}]
	The trick is to move one `\(-\)' outside from \(w=-(-w)\)
	\[
		\min_{\substack{v\in U \\ \|v\|=\lambda}}\langle v, w\rangle 
		= - \max_{\substack{v\in U \\ \|v\|=\lambda}} \langle v, -w \rangle
		= -\lambda \|\proj{U}(w)\|
	\]
	where we have used in the last equation that the projection is linear (we can move
	the minus sign out) and the norm removes the inner minus sign. The \(\argmin\) argument
	is similar.
\end{proof}

\begin{proof}[Proof of Lemma~\ref{lem: constrained maximiziation of scalar products}]
	\begin{enumerate}[label={\textbf{Step \arabic*}:},wide]
		\item  We claim that
		\[
			v^* = \lambda \frac{\proj{U}(w)}{\|\proj{U}(w)\|}
		\]
		results in the value \(\langle v^*, w\rangle = \lambda \|\proj{U}(w)\|\).

		For this we consider
		\begin{align}
			\label{eq: projection def}
			\proj{U}(w)
			&= \argmin_{v\in U}\underbrace{\|v-w\|^2}_{= \|v\|^2 - 2\langle v,w\rangle + \|w\|^2}
			\\
			\nonumber
			&= \argmin_{v\in U} \underbrace{\|v\|^2 - 2\langle v, w\rangle}_{=: f(v)}
		\end{align}
		we know that \(t\mapsto f(t \langle w\rangle_U)\) is minimized at \(t=1\)
		by the definition of \(\langle w\rangle_U\). The first order condition implies
		\[
			0\overset!= \frac{d}{dt}	
			= 2t \|\proj{U}(w)\|^2 - 2\langle \proj{U}(w), w\rangle
		\]
		and thus
		\[
			1 = t^*  = \frac{\langle \proj{U}(w), w\rangle}{\|\proj{U}(w)\|^2}
		\]
		Multiplying both sides by \(\lambda\|\proj{U}(w)\|\) finishes this step
		\begin{equation}
			\label{eq: maximum can be achieved}
			\lambda \|\proj{U}(w)\|
			= \Bigl\langle
				\underbrace{\lambda \frac{\proj{U}(w)}{\|\proj{U}(w)\|}}_{=v^*}, w
			\Bigr\rangle.
		\end{equation}

		\item By \eqref{eq: maximum can be achieved}, we know that we can achieve
		the value we claim to be the maximum (and know the location \(v^*\) to do
		so). So if we prove that we can not exceed this value, then it is a
		maximum and \(v^*\) is the \(\argmax\). This would finish the proof. What
		remains to be shown is therefore
		\[
			\langle v, w\rangle \le \lambda \|\proj{U}(w)\| \qquad \forall v\in U: \|v\| = \lambda.
		\]
		Let \(v\in U\) with \(\|v\|=\lambda\). Then for any \(\mu\in \real\) we can plug
		\(\mu v\) into \(f\) from \eqref{eq: projection def} to get
		\begin{align*}
			\mu^2\lambda^2 - 2\mu \langle v, w\rangle
			&= f(\mu v)
			\\
			\overset{\eqref{eq: projection def}}&\ge f(\proj{U}w)
			= \|\proj{U}(w)\|^2 - 2\langle \proj{U}w , w\rangle
			\\
			&= - \langle \proj{U}w, w\rangle
		\end{align*}
		where the last equation follows from \eqref{eq: maximum can be achieved} with
		\(\lambda = \|\proj{U}w\|\). Reordering we get for all \(\mu\)
		\[
			\langle \proj{U}w, w\rangle + \mu^2\lambda^2 \ge 2\mu\langle v,w\rangle
		\]
		We now select \(\mu = \frac{\|\proj{U}w\|}{\lambda}>0\) and divide both
		sides by \(\mu\) to get
		\[
			2\langle v, w\rangle
			\le \underbrace{
				\Bigl\langle \underbrace{\frac{\proj{U}(w)}{\mu}}_{=v^*}, w\Bigr\rangle
			}_{
				 \overset{\eqref{eq: maximum can be achieved}}= \lambda\|\proj{U}w\|
			} + \lambda\|\proj{U}(w)\|
			= 2 \lambda\|\proj{U}w\|
		\]
		Dividing both sides by \(2\) yields the claim.
		\qedhere
	\end{enumerate}
\end{proof} \end{subappendices} 		\printbibliography[heading=subbibliography]
	\end{refsection}

	\part{Distributional optimization}
	\label{part: distributional optimization}

	\begin{refsection}
		\chapter{Distributional optimization}
\label{chap: distributional optimization}

While it is nice to have worst-case guarantees (when the problem exhibits
convexity, the PL-inequality or similar) the optimization problems
in machine learning are decidedly not of this type. In fact they exhibit
exponentially more saddle points than minima
\autocite{dauphinIdentifyingAttackingSaddle2014} and in particular not all
critical points are global minima.

Despite this, some optimizers from convex function theory such as gradient descent find
themselves successfully optimizing machine learning objectives. Other
convex optimizers are less successful. The geometric alternative to momentum
suggested by \textcite{bubeckGeometricAlternativeNesterov2015} for example is virtually
unknown in machine learning despite an equivalent worst case complexity to
Nesterov's momentum on convex functions. Explaining these successes and failures
requires a different mathematical framework that can tackle global optimization
questions. And the ``just find a local minimum, it is not going to be that bad,
typically'' heuristic screams for an average case analysis or more generally a
distributional analysis.

A distributional analysis has lead to algorithmic improvements in many prominent cases:
\begin{itemize}[noitemsep]
    \item \textbf{Sorting algorithms.} Quicksort
    \autocite{hoareQuicksort1962} is successful in practice with an average
    performance of \(\bigO(n\log(n))\) despite its worst case
    performance of \(\bigO(n^2)\) and the existence of Merge sort with worst
    case performance \(\bigO(n\log(n))\).
    \item \textbf{Data structures.} E.g. Hash-tables or Binary-search trees,
    which sometimes have much better performance on average than in the worst
    case.
    \item \textbf{Lossless compression algorithms.} Such algorithms have
    to be injections, so they can not map into a smaller set (a set of numbers
    with fewer bits) than the origin (pigeonhole principle). But the average
    case is bounded by entropy instead \autocite[Theorem 9, ``Shannon's source
    coding theorem'']{shannonMathematicalTheoryCommunication1948}.
    \item \textbf{Linear programming.} The simplex algorithm has polynomial
    time complexity on average despite its exponential time complexity in the
    worst case \autocite[e.g.][]{borgwardtSimplexMethodProbabilistic1986,smaleAverageNumberSteps1983,spielmanSmoothedAnalysisAlgorithms2004,toddProbabilisticModelsLinear1991}.
    \item ``Random by design'' algorithms such as \textbf{Monte Carlo integration}.
\end{itemize}

Moreover, the empirical performance of many algorithms concentrates around the
average case performance. This suggests that an analysis of the distribution of
performance should often lead to tight concentration results. But to analyze
the distribution of performance or even just the average performance
of an algorithm it is necessary to either 
\begin{enumerate}[noitemsep]
    \item randomize the algorithm itself, or
    \item make assumptions about the distribution of inputs given to
    the algorithm.
\end{enumerate}

\section{Randomized Optimizer}

The most salient randomized algorithm for global optimization
may be \textbf{simulated annealing}
\autocite[e.g.][]{vanlaarhovenSimulatedAnnealingTheory1987,salamonFactsConjecturesImprovements2002}.
If we can generate samples from the softmax (a.k.a.\@ Gibbs or Boltzmann) distribution
\[
    X \sim  \mu_{T}(dx) = \frac{\exp(\frac1T\obj(x))}{\int_\domain \exp(\frac1T \obj(y)) dy} dx
\]
then as the temperature \(T\) cools\footnote{
    the state ``anneals'', inspired by the concept of annealing in metallurgy
    where metals are heated above a certain temperature before a slow gradual
    cooling process.
}
(\(T\to 0\)), the distribution concentrates around the global maximizers of
\(\obj\) in \(\domain\). The difficulty of this method lies in sampling
from this distribution. A key insight for the case of differentiable
\(\obj\) is that \emph{Langevin Diffusion},\footnote{
    Langevin Diffusion and its discretized version are grouped under the
    umbrella term \emph{Langevin Dynamics}. Sampling from \(\mu_T\) using this
    procedure is also known as \emph{Langevin Monte Carlo}
    \autocite{chenLangevinDynamicsUnified2024}.
}
i.e. the Îto SDE given by gradient flow and disturbed by Brownian noise \(B(t)\)
\[
    \tag{Langevin Diffusion}
    dX(t) = \underbrace{\nabla \obj(X(t))dt}_{\text{gradient flow}} + \underbrace{\sqrt{2T}dB(t)}_{\text{noise}},
\]
converges to \(\mu_T\) as \(t\to \infty\). The idea is then, to reduce the
temperature \(T\) slowly enough over time \(t\) such that \(X(t)\) remains
distributed according to \(\mu_T\). For example, under suitable assumptions about
\(\obj\), the temperature schedule \(T(t) := \frac{c}{\log(t)}\) for a
sufficiently large \(c>0\) can accomplish this task and
\(X(t)\) concentrates on the global maximum
\autocite{chiangDiffusionGlobalOptimization1987}.
To turn this into an optimization algorithm it is necessary to show that these
properties can be retained for a discretized version that can be used in
practice
\autocite[further information in e.g.][]{raginskyNonconvexLearningStochastic2017,xuGlobalConvergenceLangevin2018,zouFasterConvergenceStochastic2021,chenLangevinDynamicsUnified2024}.

\begin{remark}[Stochastic gradient descent]
    Langevin dynamics suggest that noise is sometimes desirable --
    especially at the beginning of training. This may explain the
    success of Stochastic gradient descent.
\end{remark}

\begin{remark}[Randomization via preprocessing]
    \label{rem: randomization via preprocessing}
    The randomization of an algorithm also often takes the form of a pre-processing
    step of the input. For example: Random shuffling as the first step can guarantee
    that the input of a sorting algorithm is uniform. With the help of these
    pre-processing steps the analysis can often be reduced to assumptions about
    the distribution of inputs to the algorithm.
\end{remark}

\section{Random objective functions}

The input to an optimization algorithm is the objective function. A
distributional assumption over the input is thus a distribution over objective
functions. The optimization of such a \emph{random function} is a fundamental problem that has
independently emerged across multiple research domains, each developing its own
terminology for very similar methods. 

In \textbf{geostatistics}, random functions are typically referred to as \emph{random fields},
which are used to model spatial distributions, such as ore deposits at various
locations \autocite{krigeStatisticalApproachBasic1951,matheronPrinciplesGeostatistics1963,steinInterpolationSpatialData1999}. In this domain, interpolating the underlying function based
on limited sample points is known as \emph{Kriging}\footnote{
    Kriging utilizes the `Best linear unbiased estimator' (BLUE), which
    coincides with the conditional expectation in the Gaussian
    case and is otherwise the minimum variance estimator in the
    space of unbiased linear estimators.
    The BLUE is effectively used if a Gaussian prior is
    assumed despite the underlying prior not being Gaussian.
} and used to guide
subsequent evaluations -- such as where to drill the next pilot hole.

In \textbf{Bayesian optimization} (BO)
\autocite{kushnerNewMethodLocating1964,jonesEfficientGlobalOptimization1998,frazierBayesianOptimization2018,garnettBayesianOptimization2023},
concerned with general black-box function optimization, it is standard
to assume a Gaussian prior and refer to random functions as \emph{Gaussian
processes}.
In BO the subsequent evaluation point \(x_{n+1}\)
is selected based on the posterior distribution of the random function \(\rf\) conditional
on the previous evaluations \(\rf(x_0), \dots \rf(x_n)\), which coincides with
Kriging in the Gaussian case.
In recent years, BO has gained prominence in machine
learning, particularly for hyperparameter tuning, where evaluating the function
(e.g.\@ training a model) is costly.

In the field of \textbf{compressed sensing} the goal is to reconstruct
a signal \(x\) from noisy observations \(y=Ax + \noise\) with sensing
matrix \(A\) and noise \(\noise\). This task is generally achieved by minimizing
a regression objective of the form \(\rf(x) = \|Ax - y\|^2 + R(x)\) with
regularization \(R\). In the analysis of \emph{Approximate
Message Passing} algorithms, the sensing matrix \(A\) is assumed to be random
\autocite{donohoMessagepassingAlgorithmsCompressed2009,donohoMessagePassingAlgorithms2010,bayatiDynamicsMessagePassing2011}.
This turns the regression objective \(\rf\) into a random function.

Finally, in \textbf{statistical physics}, random functions appear in the study
of spin glasses
\autocite[e.g.][]{montanariOptimizationSherringtonKirkpatrickHamiltonian2019,subagFollowingGroundStates2021,huangTightLipschitzHardness2022},
where they are often referred to as \emph{Hamiltonians} but more recently also
as random functions \autocite{auffingerOptimizationRandomHighDimensional2023}. 
The optima of these energy landscapes have been extensively
studied because they correspond to stable physical states. This
research on high-dimensional random functions has also lead to insights
about the loss landscapes found in machine learning
\autocite[e.g.][]{dauphinIdentifyingAttackingSaddle2014,choromanskaLossSurfacesMultilayer2015}.

\begin{definition}[Random function]
    \label{def: random function}
    A \emph{Random function} is a collection of random variables
    \(\rf=(\rf(x))_{x\in \domain}\) indexed by the domain \(\domain\) of the random
    function \(\rf\). Since random variables are defined as functions from a
    probability space \((\Omega, \cA, \Pr)\) random functions technically have
    two inputs: Elementary events \(\omega\in \Omega\) and the `actual' input
    \(x\in \domain\). The notation
    \[
        \rf\colon \domain \to \real
    \]
    is therefore abuse of notation and stands for \(\rf\colon \Omega \times \domain \to
    \real, (\omega,x) \mapsto \rf_\omega(x)\). To emphasize the difference between
    random functions \(\rf\) and normal functions \(\obj\) we denote
    random functions in bold by convention.
\end{definition}

\paragraph*{Bayesian optimization}
\autocite[e.g.][]{kushnerNewMethodLocating1964,mockusBayesianApproachGlobal1991,frazierBayesianOptimization2018,garnettBayesianOptimization2023}
For the optimization of the random function \(\rf\) Baysian optimizers select
evaluation locations \(X_{n+1}\) in \(\domain\) based on the previously seen
evaluations \(\filt^X_n := \sigma(\rf(X_0), \dots, \rf(X_n))\). That is
\(X_{n+1}\) is selected \(\filt^X_n\) measurably using the distribution of
\(\rf(x)\) conditioned on \(\filt^X_n\).
The following definition generalizes this setting to allow for the inclusion of
derivatives of \(\rf\), accommodates for the possibility of noise on the
evaluations and permits randomized algorithms (Section \ref{sec: randomized bayesian
optimizer}).

\begin{definition}[Feasible optimization algorithm]
    \label{def: feasible optimization algorithm}
    Let \(\rf:\domain\to\real\) be a continuous random objective function, i.e.
    \(\rf\) is a random element in \(C(\domain, \range)\). Let \(\noise=(\noise_n)_{n\in \nat_0}\) be a
    sequence of (location dependent) noise in \(C(\domain,\range)\), where
    \(\noise_n\) is the noise on the \(n\)-th evaluation of \(\rf\).
    For some random element \(\Problem\) such that \(\rf\) is measurable with respect to \(\Problem\)
    (e.g. \(\Problem=\rf\)), we assume the noise sequence to be iid conditional
    on \(\Problem\) and each
    \(\noise_n\) centered conditional on \(\Problem\).\footnote{
        This formulation is both very general and motivated by
        stochastic losses in the supervised machine learning setting.
        (Section \ref{sec: stochastic loss})
    }

    If \(\rf\) and \((\noise_n)_{n\in \nat_0}\) are in \(C^k(\domain, \range)\),
    the sequence of random evaluation locations \(X=(X_n)_{n\in \nat_0}\) in \(\domain\) is
    generated by a \emph{feasible \(k\)-th order optimization algorithm}, if \(X_{n+1}\) is
    independent from \((\rf,\noise)\) conditional on\footnote{
        as introduced by \autocite{kallenbergFoundationsModernProbability2002}
        we write \(\xi \indep_{\cG} \eta\), if \(\xi\) is independent from
        \(\eta\) conditional on \(\cG\). If \(\xi\) is a random variable
        in a Polish space and \(\sigma(W)=\cG\), there exists a measurable
        function \(h\) and \(U\sim\uniform(0,1)\) independent of \(W,\eta\)
        such that
        \[
            \xi=h(W,U)
        \]
        \autocite[Prop.\@ 6.13]{kallenbergFoundationsModernProbability2002}.
        Conditional independence therefore
        generalizes the case where \(\xi=h(W)\), i.e. where \(\xi\) is
        measurable with respect to \(\cG\). In our case this allows for the
        randomized selection of \(X_{n+1}\) based on a distribution that
        is \(\filt^X_n\) measurable.
    }
    \begin{equation}
        \label{Bayes-eq: filtration}
        \filt^X_n := \sigma(\filt, \obs_0, \dots, \obs_n), \qquad n \in \{ -1,0,1,\dots\}
    \end{equation}
    where \(\filt\) is the initial information and \(\obs_n = \Obs_n(X_n)\)
    is the \(n\)-th observation with
    \[
        \Obs_n(x) := \jet \rf(x) + \jet\noise_n(x),
    \]
    where \(\jet f(x)\) is the jet of differentials of \(f\) in
    \(x\) up to \(k\)-th order\footnote{
        e.g.\ for \(k=1\) the jet is of the form \(\jet\rf(x) = (\rf(x),
        \nabla\rf(x))\)
    }.
\end{definition}

Bayesian optimization algorithms select the next evaluation point based on this
conditional distribution with the help of \emph{acquisition functions}.
In the following we present prominent examples with maximization as the
objective.
\begin{example}[Acquisition functions]
    \label{ex: acquisition functions}
    \leavevmode
    \begin{itemize}[leftmargin=*]
        \item \emph{Probability of improvement} (PI) chooses
        \[
            X_{n+1} = \argmax_{x\in \domain} \Pr\Bigl(\rf(x) > \tau_n \mid \filt^X_n \Bigr)
        \]
        for some \(\tau_n\) adapted to the filtration \(\filt^X_n\).
        In the noiseless case where the values \((\rf(X_k))_{k<n}\) are contained in
        \(\filt^X_n\) the canonical selection is essentially the running maximum
        \begin{equation}
            \label{eq: tau canonical choice}    
            \tau_n := \max_{k=0,\dots, n-1} \rf(X_k)+\epsilon
        \end{equation}
        where \(\epsilon\ge 0\) is a minimal improvement.
        
        \item \emph{Expected improvement} (EI) chooses
        \[
            X_{n+1} = \argmax_{x\in \domain} \E\Bigl[ \bigl(\rf(x) - \tau_n\bigr)_+ \mid \filt^X_n \Bigr]
        \]
        with \((a)_+ = \max\{a, 0\}\) and \(\tau_n\) adapted to \(\filt^X_n\).
        The canonical choice for \(\tau_n\) in the noiseless setting is again
        given by \eqref{eq: tau canonical choice}.
        
        \item \emph{Upper confidence bound} (UCB)\footnote{
            The name of this acquisition function clashes with the
            presentation of Bayesian optimization with minimization as
            the objective.
        } maximization selects
        \[
            X_{n+1} = \argmax_{x \in \domain} \mu_n(x) + \sqrt{\beta_n} \sigma_n(x)
        \]
        where \(\mu_n(x) = \E[\rf(x) \mid \filt^X_n]\) represents the conditional expectation,
        \(\sigma_n(x) = \Var(\rf(x) \mid \filt^X_n)\) the conditional variance and \(\beta_n\) is a
        parameter that weighs the value of exploration at time step \(n\).
    \end{itemize}
\end{example}

With \(\epsilon=0\) in \eqref{eq: tau canonical choice} the PI acquisition function highly values safe, infinitesimal
improvements. To encourage exploration, non-zero \(\epsilon\) have to be selected
in practice to get reasonable convergence. However EI can work `hyperparameter
free'\footnote{
    ignoring the selection of the distributional model.
} with \(\epsilon=0\) as it
inherently values larger improvements.

\begin{remark}[Convergence proof]
    For the zero order case (without gradients) there are convergence proofs
    for the UCB algorithm \autocite{srinivasGaussianProcessOptimization2010} and for
    the EI algorithm
    \autocite{vazquezConvergencePropertiesExpected2010}. Even if the prior distribution
    of \(\rf\) is not known and has to be estimated some convergence results are
    available \autocite{berkenkampNoRegretBayesianOptimization2019}.
\end{remark}

\subsection{Randomized Bayesian optimizer}
\label{sec: randomized bayesian optimizer}

Before we take a closer look at Bayesian optimization algorithms let
us quickly note that it is also possible to combine the assumption about
the random objective function with randomized algorithms.

\emph{Thomson sampling} \autocite{thompsonLikelihoodThatOne1933} randomly
samples the subsequent evaluation point \(X_{n+1}\) according to the posterior
distribution of the \(\argmax\) of the random objective \(\rf\)
\[
    X_{n+1} \sim \Pr\bigl(\argmax_{x} \rf(x) \in \cdot \mid \filt^X_n\bigr)
\]
However, due to the complex nature of this posterior distribution,\footnote{
    except for special priors on the objective or small domains \(\domain\)
}
Thomson sampling is considered impractical for implementation
\autocite[Sec.~8.7]{garnettBayesianOptimization2023}.
\subsection{Optimal Bayesian algorithms}
\label{sec: optimal bayesian algorithms}

Intuitively, the best performance an optimizer can achieve with zero further
evaluations and available information \(\filt\) is given by the value
\[
    V_0\{\filt\} := \max_{x\in\domain} \E[\rf(x) \mid \filt].
\]
It is achieved with the greedy choice of
\[
    Q_0\{\filt\} := \argmax_{x\in \domain} \E[\rf(x) \mid \filt].
\]
To consider the case of \(n\) evaluations we want to proceed by backwards
induction. After a single evaluation, one has \(n-1\) evaluations left and those
evaluations are affected by the shifted noise sequence \(\fshift(\noise)= (\noise_{n+1})_{n\in
\nat_0}\).

Thus, if \(V_{n-1}^{\fshift(\noise)}\{\filt\}\) is the best performance that is
achievable with information \(\filt\) and \(n-1\) further evaluations,
then with \(n\) evaluations available, the best one can do with information
\(\filt\) is
\begin{align*}
    \tag{\(n\)-step lookahead value}
    V_n^\noise\{\filt\}
    &:= \max_{x\in \domain} \E\bigl[ V_{n-1}^{\fshift(\noise)}\{\filt, \Obs_{0}(x)\} \mid \filt\bigr]
    \\
    \tag{\(n\)-step lookahead choice}
    Q_n^\noise\{\filt\}
    &:= \argmax_{x\in \domain} \E\bigl[ V_{n-1}^{\fshift(\noise)}\{\filt, \Obs_{0}(x)\} \mid \filt\bigr].
\end{align*}
Clearly, with \(V_0^\noise:=V_0\) and \(Q_0^\noise:= Q_0\) this results in
well defined \(V_n^\noise\) and \(Q_n^\noise\).

\begin{remark}[Knowledge Gradient]
    \label{rem: knowledge gradient}
    The one step lookahead, which selects
    \[
        X_n := Q_1^{\fshift^n(\noise)}\{\filt^X_{n-1}\}
        = \argmax_{x\in \domain} \E\Bigl[\max_y \E\bigl[\rf(y) \mid \filt^X_{n-1}, \Obs_n(x)\bigr] \Bigm| \filt^X_{n-1}\Bigr]
    \]
    is known as the \emph{knowledge gradient}\footnote{
        see e.g.\ \textcite{frazierKnowledgeGradientPolicyCorrelated2009},
        therein a convergence proof can also be found.
    }
     (KG) acquisition
    function in Bayesian optimization.
\end{remark}

The following theorem shows that backwards induction using the \(Q\) functions results in the
algorithm with optimal average case performance
after \(N\) evaluations. 

\begin{theorem}[\(N\)-evaluation optimal algorithm]
    \label{thm: N-eval optimal}
    In the setting of feasible optimization algorithms (Definition \ref{def:
    feasible optimization algorithm}), the algorithm that chooses
    \begin{equation}
        \label{eq: optimal N-eval algo}
        X_n := Q_{N-n}^{\fshift^n(\noise)}\{\filt^X_{n-1}\}
        \qquad\text{for} \quad n\in \{0,\dots, N\},
    \end{equation}
    and \(X_n=X_N\) for \(n> N\), is a feasible \(k\)-th order optimization
    algorithm that achieves the best average performance in the set of all
    feasible algorithms after \(N\) evaluations of \(k\)-th order information.
    That is, if \(\tilde{X}_0, \dots, \tilde{X}_N\) is generated by a feasible
    \(k\)-th order algorithm,
    then almost surely
    \[
        \E[\rf(X_N) \mid \filt]
        = V_N^\noise\{\filt\}
        \ge \E[\rf(\tilde{X}_N)\mid \filt],
    \]
    where we recall that \(\filt=\filt_{-1}\) is the initial information.
\end{theorem}

\begin{remark}[Randomization does not provide an advantage]
    Observe that \(X_n\) is selected as a function of \(\filt^X_{n-1}\) and
    randomization as in Thomson sampling \ref{sec: randomized bayesian
    optimizer} is unnecessary to obtain optimal performance.
\end{remark}

\begin{proof}
    \begin{enumerate}[label={\textbf{Step \arabic*:}},wide,labelindent=0pt]
        \item More than the first equality, we will prove for all \(n\in \{0,\dots, N\}\)
        \[
            \E[\rf(X_N)\mid \filt]
            = \E[V_{N-n}^{\fshift^n(\noise)}\{\filt^X_{n-1}\}\mid \filt]
        \]
        by backwards induction. 
        Since \(V_N^\noise\{\filt^X_{-1}\}\) is measurable with
        respect to \(\filt^X_{-1} = \filt\) by definition,
        the case \(n=0\) implies the first part of the claim
        \[
            \E[\rf(X_N) \mid \filt]
            = \E[V_{N}^{\fshift^0(\noise)}\{\filt^X_{-1}\}\mid \filt]
            = V_N^\noise(\filt).
        \]
        The induction start with \(n=N\) is given by
        \begin{align*}
            \E\bigl[V_0^{\fshift^N(\noise)}\{\filt^X_{N-1}\}\mid \filt\bigr]
            &= \E\bigl[\max_{x\in \domain}\E[\rf(x) \mid \filt^X_{N-1}]\mid \filt\bigr]
            \\
            &= \E\Bigl[\E\bigl[\rf(Q_0\{\filt^X_{N-1}\}) \mid \filt^X_{N-1}\bigr]\mid \filt\Bigr]
            \\
            &= \E\bigl[\E[\rf(X_N) \mid \filt^X_{N-1}]\mid \filt\bigr]
            \\
            \overset{\text{tower}}&= \E[\rf(X_n) \mid \filt].
        \end{align*}
        For the induction step \((n+1) \to n\), recall that by definition of
        \(X_n\) and \(Q_n^\noise\) we have
        \begin{align*}
            X_n
            &= Q_{N-n}^{\fshift^n(\noise)}\{\filt^X_{n-1}\}
            \\
            &= \argmax_{x\in\domain} \E\bigl[
                V_{N-n-1}^{\fshift^{n+1}(\noise)}\{\filt^X_{n-1}, \obs_0(x)\}
            \mid \filt^X_{n-1}\bigr].
        \end{align*}
        And thus by definition of \(V_n^\noise\) and \(\filt^X_n\) and Example \ref{ex: conditional minimization}
        \begin{align*}
            V_{N-n}^{\fshift^n(\noise)}\{\filt^X_{n-1}\}
            &= \E\bigl[
                V_{N-n-1}^{\fshift^{n+1}(\noise)}\{\filt^X_{n-1}, \obs_0(X_n)\}
            \mid \filt^X_{n-1}\bigr].
            \\
            &= \E\bigl[
                V_{N-(n+1)}^{\fshift^{n+1}(\noise)}\{\filt^X_{(n+1)-1}\}
            \mid \filt^X_{n-1}\bigr].
        \end{align*}
        Since we have the claim for \(n+1\) we thus obtain for \(n\) by induction
        \[
            \E[V_{N-n}^{\fshift^n(\noise)}\{\filt^X_{n-1}\} \mid \filt]
            \overset{\text{tower}}= 
            \E\bigl[
                    V_{N-(n+1)}^{\fshift^{n+1}(\noise)}\{\filt^X_{(n+1)-1}\}
                \mid \filt
            \bigr]
            \overset{\text{ind.}}= \E[\rf(X_n) \mid \filt].
        \]

        \item To avoid notational clutter we drop the tilde and show
        for any initial information \(\filt\) and any \(X=(X_n)_{n\in \nat_0}\)
        generated by a feasible \(k\)-th order optimization algorithm 
        by induction over \(n\)
        \[
            \E[\rf(X_n) \mid \filt] \le V_n^{\fshift^m(\noise)}\{\filt\}.
        \]
        for all \(m\in \nat_0\). The selection \(n=N\) and \(m=0\) then yields the claim.

        Observe that \(X_0\) is independent from \(\rf,\noise\) conditional on
        \(\filt\) and therefore there exists \(U\sim\uniform(0,1)\) independent
        from \((\rf,\noise, \filt)\) such that
        \begin{equation}
            \label{eq: conditional indep. representation of X_0}    
            X_0=h(W, U)
        \end{equation}
        for some random element \(W\) that generates \(\filt\)\footnote{
            the identity map from the underlying
            probability space \((\Omega, \cA)\) into \((\Omega, \filt)\)
            does the job.
        } and a measurable
        function \(h\).
        For the induction start with \(n=0\) this implies for
        all \(m\in
        \nat_0\)
        \begin{align*}
            \E[\rf(X_0) \mid \filt]
            \overset{\text{tower}}&=
            \E\bigl[ \E[\rf(X_0) \mid \rf, W] \mid \filt\bigr]
            \\
            \overset{\eqref{eq: conditional indep. representation of X_0}}&= \E\Bigl[ \int_0^1 \rf(h(W, u)) du \mid \filt\Bigr]
            \\
            &= \int_0^1 \underbrace{\E\bigl[\rf(h(W,u)) \mid \filt\bigr]}_{\le \max_{x} \E[\rf(x) \mid \filt]} du
            \le V_0\{\filt\} = V_0^{\fshift^m(\noise)}\{\filt\},
        \end{align*}
        where we use that \(h(W,u)\) is an \(\filt\) measurable random variable
        and Example \ref{ex: conditional minimization}.

        Let us now consider the induction step \(n\to n+1\). Since the shifted
        process \(\fshift(X) := (X_{n+1})_{n\in \nat_0}\) retains the
        conditionally independence property with respect to the shifted
        filtration \(\fshift(\filt^X) := (\filt^X_{n+1})_{n\in \{-1,0,\dots\}}\) we can
        apply the induction assumption for \(n\) and \(\filt = \filt^X_0\) to obtain for all \(m\in \nat_0\)
        \begin{align*}
            \E[\rf(X_{n+1}) \mid \filt^X_0]
            &= \E\bigl[\rf(\fshift(X)_n) \mid \fshift(\filt^X)_{-1}\bigr]
            \\
            \overset{\text{ind.}}&\le
            V_n^{\fshift^{m+1}(\noise)}\{\fshift(\filt^X)_{-1}\}
            \\
            &= V_n^{\fshift^{m+1}(\noise)}\{\filt^X_0\}
        \end{align*}
        We therefore have
        \begin{align*}
            \E[\rf(X_{n+1}) \mid \filt]
            \overset{\text{tower}}&= \E\bigl[\E[\rf(X_{n+1}) \mid \filt^X_0]\mid \filt\bigr]
            \\
            \overset{\text{monot.}}&\le
            \E\bigl[V_n^{\fshift^{m+1}(\noise)}\{\filt^X_0\} \mid \filt\bigr]
            \\
            &= \E\bigl[V_n^{\fshift^{m+1}(\noise)}\{ \filt, \Obs_0(X_0)\} \mid \filt\bigr]
            \\
            \overset{\text{tower}}&= \E\Bigl[\E\bigl[V_n^{\fshift^{m+1}(\noise)}\{ \filt, \Obs_0(X_0)\} \mid \rf, \noise, W\bigr] \mid \filt\Bigr]
            \\
            \overset{\eqref{eq: conditional indep. representation of X_0}}&=
            \int_0^1 \E\Bigl[V_n^{\fshift^{m+1}(\noise)}\{ \filt, \Obs_0(h(W, u))\} \mid \filt\Bigr]du
            \\
            \overset{(*)}&\le \max_{x} \E[V_n^{\fshift^{m+1}(\noise)}\{ \filt, \Obs_0(x)\} \mid \filt]
            \\
            &= V_{n+1}^{\fshift^m(\noise)}(\filt).
        \end{align*}
        For \((*)\) we use that \(h(W,u)\) is \(\filt\) measurable
        and Example \ref{ex: conditional minimization}.
        \qedhere
    \end{enumerate}    
\end{proof}

\begin{remark}[Utility]
    One may assign a utility \(u\) to reaching a certain function value \(\rf\).
    But in this case maximizing \(\rf(X_n)\) turns into maximizing
    \(u(\rf(X_n))\) and one can just consider \(\tilde{\rf} = u \circ \rf\)
    instead.
\end{remark}

\begin{remark}[Regret]
    The performance of the intermediate optimization steps
    \(\rf(X_n)\) with \(n\le N\) is completely disregarded for the optimality of
    Theorem~\ref{thm: N-eval optimal}. Only the final performance matters.
    This allows algorithms to prioritize information gain and therefore
    exploration up until timestep \(N\).
    One could instead consider a weighted average of \(\rf(X_n)\) as a
    performance target. To avoid the requirement that the weights \(\gamma =
    (\gamma_n)_{n\in \nat_0}\) have to be be summable one can instead consider
    the weighted regret
    \[
        R_\gamma := 
            \sum_{n=0}^\infty \gamma_n \bigl( \max_{y\in \domain} \rf(y) - \rf(X_n)\Bigr)
    \]
    which is equivalent but allows \(\gamma_n=1\) for all \(n\) as long as \(\rf(X_n)\) converges quickly enough
    against the global maximum.\footnote{
        Theorem~\ref{thm: N-eval optimal} treats optimality for the case
        \(\gamma_n=\ind_{N=n}\) where the stated algorithm achieves minimal
        expected regret.
    }

    The iteration we considered translates to
    \[
        V_\gamma^{\noise}\{\filt\} := \max_{x\in \domain} \E\Bigl[
            \underbrace{\gamma_0 \bigl(\rf(x) - \max_{y\in\domain}\rf(y)\bigr)}_{
                \substack{\text{immeditate}\\\text{value}}
            }
            + \underbrace{V_{\fshift(\gamma)}^{\fshift(\noise)}\{\filt, \Obs_0(x)\}}_{\text{future value}} \mid \filt
        \Bigr]
    \]
    using the shift operator \(\fshift(\gamma) := (\gamma_{n+1})_{n\in \nat_0}\). 
    If the sequence of \(\gamma\) is eventually zero, i.e. \(\fshift^k(\gamma) = \mathbf 0\)
    for some \(k\), then termination with \(V_{\mathbf{0}}\{\filt\} :=0\) ensures \(V_\gamma\)
    is well defined.\footnote{
        For \(\gamma_n=\ind_{N=n}\) it follows
        \[V_\gamma\{\filt\} = V_N\{\filt\} - \E\bigl[\max_{y\in \domain} \rf(y)\mid \filt\bigr].\]
    }
    In the general case one has to define the value function \(V_\gamma^\noise\)
    as the attained expected negative regret of the algorithm that maximizes
    the expected negative regret and then prove this recursion as a theorem.
    This is similar to value functions in reinforcement learning.
\end{remark}

\subsection{Practical Bayesian optimization}
\label{sec: practical bayesian optimization}

Section~\ref{sec: optimal bayesian algorithms} demonstrates
that optimal average performance algorithms are relatively easy to obtain in an
abstract form. The problem is that abstract maximizers of acquisition functions
are not necessarily computable. To make the acquisition functions calculable
it is typically assumed that \(\rf\) is Gaussian.

More generally, assume \(\rf\) and the random variables generating the
filtration \(\filt^X_n\) to be jointly Gaussian.\footnote{
    Treating the, typically random, evaluation points \(X_n\) as deterministic
    (cf.\ Chapter \ref{chap: measure theory for rf opt}, in particular Corollary \ref{cor: gaussian previsible sampling})
} The posterior distribution of
\(\rf(x)\) given \(\filt^X_n\) is then Gaussian by Theorem~\ref{thm: conditional
gaussian distribution}, i.e.
\[
    \rf(x) \mid \filt^X_n \sim \normal(\mu_n(x), \sigma_n^2(x)),
\]
where we denote by \(\mu_n(x) := \E[\rf(x) \mid \filt^X_n]\) the conditional mean
and by \(\sigma_n^2(x) := \Var(\rf(x) \mid \filt^X_n)\) the conditional variance.

Recall the UCB acquisition function in Example~\ref{ex: acquisition functions}
was already defined using only this conditional mean and variance. In the following
theorem we present results that show that PI and EI can also be expressed using
\(\mu_n(x)\) and \(\sigma_n^2(x)\) assuming a Gaussian distribution.

\begin{theorem}[More explicit acquisition functions]
    \label{thm: more explicit acquisition functions}
    Let \(\rf\) and  be jointly Gaussian. Then
    \begin{enumerate}
        \item Probability of improvement (PI) 
        \[
            \argmax_{x\in \domain}\Pr\bigl(\rf(x) > \tau_n \mid \filt^X_n\bigr)
            = \argmax_{x\in \domain}\underbrace{\frac{\mu_n(x) - \tau_n}{\sigma_n(x)}}_{=:\pi_n(x)}
        \]
        with \(\tau_n\) adapted to the filtration \(\filt^X_n\).
        \item Expected improvement (EI)
        \[
            \E\Bigl[ \bigl(\rf(x) - \tau_n\bigr)_+ \mid \filt^X_n \Bigr]
            = \sigma_n(x) \Bigl(\density\bigl(\pi_n(x)\bigr) + \pi_n(x) \cdf\bigl(\pi_n(x)\bigr)\Bigr),
        \]
        where \(\tau_n\) and \(\pi_n\) are as for PI, and \(\density\) is the density and
        \(\cdf\) the cumulative distribution function of the standard normal
        distribution.
    \end{enumerate}
\end{theorem}
\begin{proof}[Proof {(e.g. \cite[chap.~8]{garnettBayesianOptimization2023})}]
    Since the posterior distribution of \(\rf(x)\) given  \(\filt^X_n\) is Gaussian with mean \(\mu_n\)
    and variance \(\sigma_n^2\) we have that
    \[
        Z_n(x) := \frac{\rf(x) -\mu_n(x)}{\sigma_n(x)}
    \]
    is standard normal conditional on \(\filt^X_n\).  Reordering the inequality in
    the PI acquisition function therefore implies it is equal to 
    \[
        \Pr\bigl(\rf(x) > \tau_n \mid \filt^X_n\bigr)
        = \cdf\Bigl(\frac{\mu_n(x) - \tau_n}{\sigma_n(x)}\Bigr)
        = \cdf(\pi_n(x)).
    \]
    Since the cumulative distribution function \(\cdf\) of the standard
    normal distribution is monotonously increasing, maximization of this
    term is equivalent to maximizing its input \(\pi_n(x)\).

    Similarly transforming \(\rf(x)\) and \(\tau_n\) in the EI acquisition
    we obtain
    \[
        \E[(\rf(x) - \tau_n)_+ \mid \filt^X_n]
        = \sigma_n(x) \E\bigl[(Z_n(x) - \pi_n(x))_+ \mid \filt^X_n\bigr]
    \]
    since \(\sigma_n(x)\) is non-negative and \(\filt^X_n\) measurable. Since
    \(\pi_n(x)\) is \(\filt^X_n\) measurable and \(Z_n(x)\) standard normal
    conditional on \(\filt^X_n\) it is thus sufficient to calculate for \(Y\sim \normal(0,1)\)
    and \(x\in \real\)
    \begin{align*}
        \E[(Y + x)_+]
        &= \frac1{\sqrt{2\pi}} \int_{-x}^\infty (\magenta{y}+\teal{x}) \exp(-\tfrac{y^2}2)dy
        \\
        &= \int_{-x}^\infty \magenta{\frac{d}{dy}}\Bigl(- \tfrac1{\sqrt{2\pi}} \exp(-\tfrac{y^2}2)\Bigr)dy
        + \teal{x}(1-\cdf(-x))
        \\
        &= \density(x) + x\cdf(x).
        \qedhere
    \end{align*}
\end{proof}

While the knowledge gradient acquisition function (Remark~\ref{rem: knowledge
gradient}) may be evaluated in \(\bigO(n \log(n))\)
time\autocite{frazierKnowledgeGradientPolicyCorrelated2009} for domains
\(\domain\) of size \(n=\#\domain\) it becomes
intractable for infinite domains in all but a few special cases
\autocite[Sec.~8.3]{garnettBayesianOptimization2023}.

Observe that the acquisition functions of Theorem~\ref{thm:
more explicit acquisition functions} contain an \(\argmax\) that may be
difficult to obtain depending on the nature of the posteriors \(\mu_n\) and
\(\sigma_n\). Roughly speaking: If the posterior mean \(\mu_n\) is
as difficult to optimize as the function \(\rf\) then one has not gained
anything with this approach.
This is the motivation for the `Random Function Descent' acquisition function
(Chapter \ref{chap: random function descent}).

\paragraph{Higher evaluation costs justify more sophisticated acquisition
functions}
If the evaluation of \(\rf\) is significantly more expensive 
than evaluations of \(\mu_n\) and \(\sigma_n\), then the conversion of the
problem from \(\argmax_x \rf(x)\) to the surrogate maximization problem
\(\argmax_x F_n(\mu_n(x), \sigma_n(x))\), for some function \(F_n\) induced
by the acquisition function, can be worth it. An example is `hyperparameter
tuning' in machine learning, where the evaluation of one `hyperparameter' is the
final error after training the entire model using this hyperparameter. Compared
to the training of a machine learning model the evaluation of most posteriors is
very cheap. 

 \begin{subappendices}
    \section{Conditional Gaussian distribution}
\label{sec: conditional Gaussian distribution}

In Section~\ref{sec: practical bayesian optimization} we claimed that
the posterior distribution of \(\rf(x)\) given \(\filt^X_n\) is Gaussian.
For this we assumed that the noisy evaluations \(\obs_0,\dots, \obs_n\)
generating the filtration\footnote{
    Here we omit the initial information \(\filt\)
    which may be added as another Gaussian random variable \(\obs_{-1}\)
    such that all variables are joint Gaussian.
} \(\filt^X_n\) are joint Gaussian together with \(\rf\), treating
the evaluation locations as deterministic.\footnote{
    see Chapter \ref{chap: measure theory for rf opt} especially Example
    \ref{ex: gaussian conditional} and Corollary \ref{cor: gaussian previsible sampling}
}
This implies \((\rf(x),\obs_0,\dots, \obs_n)\) is a Gaussian vector for any
\(x\) and to calculate the posterior distribution of \(\rf(x)\) given \(\obs_0,
\dots, \obs_n\) we need to be able to obtain the posterior distribution
of multivariate Gaussian random vectors.

Theorem~\ref{thm: conditional gaussian distribution}
presents the well known explicit conditional distribution
of multivariate Gaussian random vectors 
\autocite[e.g.][Prop.~3.13.]{eatonMultivariateStatisticsVector2007}.

\begin{theorem}[Conditional Gaussian distribution]
    \label{thm: conditional gaussian distribution}
    Let \(X\sim\normal(\mu,\Sigma)\) be a multivariate Gaussian vector where
    the mean is a block vector and the covariance matrix a block matrix of the form
    \[
        \mu = \begin{bmatrix}
            \mu_1\\ \mu_2
        \end{bmatrix}
        \quad \text{and} \quad
        \Sigma = \begin{bmatrix}
            \Sigma_{11} & \Sigma_{12}
            \\
            \Sigma_{21} & \Sigma_{22}
        \end{bmatrix},
    \]
    then the conditional distribution of \(X_2\) given \(X_1\) is
    \[
        X_2\mid X_1 \sim \normal(\mu_{2\mid 1}, \Sigma_{2\mid 1}),
    \]
    with conditional mean and variance
    \begin{align*}
        \mu_{2\mid 1} &:= \mu_2 + \Sigma_{21}\Sigma_{11}^{-1}(X_1-\mu_1)
        \\
        \Sigma_{2\mid 1} &:= \Sigma_{22} -  \Sigma_{21}\Sigma_{11}^{-1}\Sigma_{12}.
    \end{align*}
    where \(\Sigma_{11}^{-1}\) may be a generalized inverse
    \autocite[cf.][]{eatonMultivariateStatisticsVector2007}.
\end{theorem}
\begin{proof}
    For the general statement we point to the standard literature
    \autocite[e.g.][]{eatonMultivariateStatisticsVector2007}.
    For this proof we will assume for simplicity that \(\Sigma\) is invertible.

    Let \(\bar{X} := X-\mu\) be the centered version of \(X\). Then there exists
    a unique lower triangular matrix \(L\) such that \(\Sigma = LL^T\) (i.e.\ the Cholesky
    Decomposition). This results in the following representation\footnote{
        It is straightforward to check that \(Y := L^{-1}(X-\mu)\)
        is a multivariate Gaussian vector of centered, uncorrelated entries with variance \(1\).
        They are therefore \(\iid\) standard normal since uncorrelated multivariate Gaussian
        random variables are independent. Sometimes this representation is even
        taken as the definition of a multivariate normal distribution.
        This is the only real place we use the invertibility of \(\Sigma\) in
        its entirety, the invertibility of \(\Sigma_{11}^{-1}\) is used later on
        because we do not want to get into the business of generalized inverses
        here.
    }
    \[
        X - \mu =:
        \begin{bmatrix}
            \bar{X}_1
            \\
            \bar{X}_2
        \end{bmatrix}
        = \begin{bmatrix}
            L_{11} & 0
            \\
            L_{21} & L_{22}
        \end{bmatrix}
        \begin{bmatrix}
            Y_1 \\ Y_2
        \end{bmatrix}
        = LY
    \]
    with independent standard normal \(Y_i\), i.e. \(Y\sim\normal(0, \identity)\).
    \(L_{11}\) is invertible since \(\Sigma_{11}\) and therefore the map from
    \(Y_1\) to \(X_1\) is bijective. Conditioning on \(X_1\) is therefore
    equivalent to conditioning on \(Y_1\). But we have
    \begin{equation}
        \label{eq: underlying explicit decomposition}
        X_2 = \mu_2 + \bar{X}_2
        = \underbrace{\mu_2 + L_{21}Y_1}_{\text{conditional expectation}} + \underbrace{L_{22} Y_2}_{\text{conditional distribution}}
    \end{equation}
    So it follows that
    \[
        X_2\mid X_1 \sim  \normal(\mu_{2\mid 1}, \Sigma_{2\mid 1})
    \]
    with
    \begin{alignat*}{3}
        \mu_{2\mid 1} &:= \E[X_2 \mid X_1]
        &=& \mu_2 + L_{21}Y_1
        \\
        \Sigma_{2\mid 1} &:= \Cov[X_2 \mid X_1] = \E\Bigl[ \bigl(X_2 - \E[X_2 \mid X_1]\bigr)^2 \mid X_1\Bigr]
        &=& L_{22}L_{22}^T.
    \end{alignat*}
    What is left to do, is to find a representation for the \(L_{ij}\)
    using the block matrices of \(\Sigma\). For this we observe
    \begin{equation}
        \label{eq: sigma to L conversion}
        \setlength\arraycolsep{4pt}
        \begin{bmatrix}
            \Sigma_{11} & \Sigma_{12}\\
            \Sigma_{21} & \Sigma_{22}
        \end{bmatrix}
        = \Sigma
        = LL^T = \begin{bmatrix}
            L_{11}L_{11}^T & L_{11}L_{21}^T
            \\
            L_{21}L_{11}^T & L_{22}L_{22}^T + L_{21}L_{21}^T
        \end{bmatrix}.
    \end{equation}
    Using \(Y_1 = L_{11}^{-1}\bar{X}_1\) and the insertion of an identity matrix this implies
    \[
        L_{21}Y_1 = L_{21}(L_{11}^T L_{11}^{-T})(L_{11}^{-1}\bar{X}_1)
        \overset{\eqref{eq: sigma to L conversion}}= \Sigma_{21}\Sigma_{11}^{-1}(X_1-\mu_1).
    \]
    resulting in the desired conditional expectation \(\mu_{2\mid 1}\). The conditional variance follows alike
    \begin{align*}
        \Sigma_{2\mid 1} = L_{22}L_{22}^T
        \overset{\eqref{eq: sigma to L conversion}}&= \Sigma_{22} - L_{21}L_{21}^T
        \\
        &= \Sigma_{22} - \underbrace{L_{21}(L_{11}^T}_{\overset{\eqref{eq: sigma to L conversion}}=\Sigma_{21}}\underbrace{L_{11}^{-T}) (L_{11}^{-1}}_{\overset{\eqref{eq: sigma to L conversion}}=\Sigma_{11}^{-1}}\underbrace{L_{11})L_{21}^T}_{\overset{\eqref{eq: sigma to L conversion}}=\Sigma_{12}}.
        \qedhere
    \end{align*}
\end{proof}

\begin{remark}[Decomposition]
    \label{rem: decomposition}
    \(X_2\mid X_1 \sim \normal(\mu_{2\mid 1}, \Sigma_{2\mid 1})\) always implies
    that there exists there exists \(Y_2\sim\normal(0,\identity)\) independent
    of \(X_1\) such that the following equality is true in distribution
    \begin{equation}
        \label{eq: explicit decomposition}
        X_2 = \mu_{2\mid 1} + \sqrt{\Sigma_{2\mid 1}} Y_2,
    \end{equation}
    where \(\sqrt{\Sigma}\) denotes the cholesky decomposition of \(\Sigma\).
    If the covariance matrix is moreover invertible, then \(Y_2\) can
    be determined constructively as in the proof of Theorem~\ref{thm:
    conditional gaussian distribution} such that equation holds always and not
    just in distribution (cf.~\eqref{eq: underlying explicit decomposition}).
    This might be true in the non-invertible case as well
    but would require deeper insights about singular matrices.
\end{remark}

     \section{Sample regularity of random functions}
\label{sec: sample regularity of random functions}

In the previous section we made the assumption that \((\rf, \obs_0,\dots,
\obs_n)\) was joint Gaussian, treating the evaluation locations as deterministic.\footnote{see Chapter \ref{chap: measure theory for rf opt}.}
Since the observations are of the form
\[
	\obs_i = \Obs_i(X_i) = \jet\rf(X_i) + \jet\noise(X_i),
\]
we need \(\jet\rf\) and \(\jet\noise\) to be joint Gaussian random functions.
For the case \(k=0\) it is sufficient if \(\rf\) and \(\noise\) are joint
Gaussian random functions. But for \(k>0\) the question is what requirements 
\(\rf\) and \(\noise\) need to satisfy such that \(\jet\rf\) and \(\jet\noise\)
are joint Gaussian. But the question is not really whether
or not it is joint Gaussian, since derivatives are Gaussian
as limits of Gaussian random variables
\begin{equation}
	\label{eq: derivative as limit}	
	D_v\rf(x) = \lim_{t\to 0} \tfrac1t (\rf(x+tv) - \rf(x)).
\end{equation}
The real question is, do derivatives exist? Specifically, since Gaussian random
functions are understood in terms of mean and covariance --
what requirements do mean and covariance function need to satisfy
to ensure the samples \(\rf\) are regular enough for \(\rf\) to
be \(k\) times continuously differentiable.

Clearly, if the mean is in \(C^k(\domain, \range)\), then we may add
it to a centred Gaussian random function in \(C^k(\domain, \range)\)
to obtain a Gaussian random function in \(C^k(\domain, \range)\) with
the desired mean and covariance. We can thus assume centered Gaussian
random functions without loss of generality.

Gaussian processes, like most infinite dimensional random
objects, are typically constructed as a collection of random variables
\(\rf=(\rf(x))_{x\in \domain}\) using Kolmogorov's extension theorem
\autocite[e.g.][Thm.\ 14.36]{klenkeProbabilityTheoryComprehensive2014}.
For such a collection of random variables \(\rf\) one first needs to prove there
exists a version\footnote{
	if \(\rf\) is in \(C^k(\domain, \range)\) up to a probability zero set
	\(N\), it can be set to zero on \(N\) resulting in the version \(\tilde\rf_\omega(x) =
	\rf_\omega(x) \ind_{N^\complement}(\omega)\), which is in \(C^k(\domain, \range)\) for all
	elementary events \(\omega \in \Omega\) such that \(\tilde\rf\) becomes a
	random variable in \(C^k(\domain,\range)\), assuming the Borel
	\(\sigma\)-algebra on \(C^k(\domain, \range)\) is compatible with the
	product \(\sigma\)-algebra, see Remark \ref{rem: construction of a cont.
	random function}
} that lies in \(C^k(\domain, \range)\).

The requirements on the covariance kernel for sample path regularity
are well known.\footnote{
	see \eg \textcite{costaSamplePathRegularity2024}, \textcite[chap. 5]{scheurerComparisonModelsMethods2009},
	\textcite[Ch.~IV, \S 3, Def.~3 and
	following]{gihmanTheoryStochasticProcesses1974},
	\textcite{talagrandRegularityGaussianProcesses1987} and
	references therein.
} 
To understand these requirements consider a centered Gaussian random
function \(\rf\). If the derivative \(\partial_{x_i}\rf(x)\) exists, then
by \eqref{eq: derivative as limit} and and an
argument why the expectation and limit may be exchanged \autocite[Sec.\ 5.3.2]{scheurerComparisonModelsMethods2009}
\[
	\E[\partial_{x_i}\rf(x)] = 0
\]
and
\begin{equation}
	\label{eq: covariance of derivatives of rf, general case}	
	\Cov(\partial_{x_i}\rf(x), \partial_{y_j}\rf(y))
	= \E[\partial_{x_i}\rf(x) \partial_{y_j}\rf(y)]
	= \partial_{x_i} \partial_{y_j} \C_\rf(x,y).
\end{equation}
So clearly, for \(\partial_{x_i}\rf(x)\) to exist, \(\partial_{x_i} \partial_{y_j} \C_\rf(x,y)\)
must exist. And in fact the existence of \(\partial_{x_i} \partial_{y_i}
\C_\rf(x,y)\) is equivalent to the existence of \(\partial_{x_i}\rf(x)\) as an
\(L^2\) limit of \eqref{eq: derivative as limit} \autocite[Thm.\
5.3.10]{scheurerComparisonModelsMethods2009}. In this case
derivative and expectation may be freely exchanged, \ie we also have
for the mixed derivatives
\begin{equation}
	\label{eq: mixed derivatives of rf}	
	\Cov(\partial_{x_i}\rf(x), \rf(y)) = \partial_{x_i} \C_\rf(x,y).
\end{equation}
However for \(\partial_{x_i}\rf(x)\) to exist as a continuous random function,
slightly more regularity on the covariance function is necessary. It turns out
that Hölder continuity of \(\partial_{x_i} \partial_{y_i} \C_\rf(x,y)\) in
both entries leads to Hölder continuity of the sample function \(\partial_{x_i}\rf(x)\)
\autocite[Thm.\ 7]{costaSamplePathRegularity2024}.

In summary, if \(\rf\) is a Gaussian random function in \(C(\domain, \range)\),
then \(\jet\rf\) is a Gaussian random function and the covariance function
can be calculated using the rules in \eqref{eq: mixed derivatives of rf} and
\eqref{eq: covariance of derivatives of rf, general case}.

 \end{subappendices}
 		\printbibliography[heading=subbibliography]
	\end{refsection}
	\begin{refsection}
		\chapter{Measure theory for random function optimization}
\label{chap: measure theory for rf opt}

In this chapter we highlight and address measure theoretical challenges
that must be overcome for a rigorous mathematical treatment of random function
optimization. In particular, we show how non-rigorous but intuitive claims
about the conditional distributions of the objective function \(\rf\)
can be made rigorous. Before that though, we need to address measurability itself.

\paragraph*{Measurability of random evaluations.}

The evaluation \(\rf(X)\) of a random function \(\rf=(\rf(x))_{x\in \domain}\)
at a random location \(X\) is a complicated measure theoretical object.
Ex ante not even the \textbf{measurability of \(\rf(X)\)}, i.e. its existence as a
random variable, is certain. In the study of stopping times this problem is
sometimes defined away by the requirement that the stochastic process
\(\rf=(\rf(s))_{s\in \real_+}\) is
\emph{progressive} \autocite[e.g.][Lem.
7.5]{kallenbergFoundationsModernProbability2002}. This leaves the user
with the burden to confirm that their process is in fact progressive. However, the main reason
we do not take this approach is that it is tailored for stochastic processes with
one dimensional input and stopping times.

We prefer the assumption that the evaluation function \(e(f,x):= f(x)\) is
\emph{measurable}, which immediately implies that \(\rf(X) = e(\rf, X)\) is a random
variable for any random variable \(X\). This assumption holds with
great generality:
If \(\rf\) is a \emph{continuous} random function with locally compact,
separable, metrizable domain \(\domain\) and polish co-domain \(\range\),
then the evaluation function \(e\) is continuous and thereby measurable (cf.
Theorem~\ref{thm: continuous function}). This accounts for almost all
continuous applications, in particular the case \(\domain \subseteq \real^\dims\) and \(\range =
\real^n\).

\paragraph*{Conditional distributions.}
During the optimization of a random function \(\rf=(\rf(x))_{x\in \domain}\) one
typically selects evaluation locations \(X_{n+1}\) in \(\domain\) based on the
previously seen evaluations \(\filt_n:=\sigma(\rf(X_0), \dots, \rf(X_n))\).
Since \(X_{n+1}\) is thereby measurable with respect to \(\filt_n\), 
the sequence \((X_n)_{n\in \nat}\) is called \emph{previsible}\footnote{see also Definition \ref{def: previsible}} with respect to
the filtration \((\filt_n)_{n\in\nat}\). Our main result is a formalization
of the intuitive notion that previsible evaluation locations may be
treated as if they are deterministic during the calculation of the conditional
distribution
\[
    \Pr(\rf(X_n)\in \cdot \mid \rf(X_0), \dots, \rf(X_{n-1})).
\]
In the case of Gaussian random functions \(\rf\) for example,
\((\rf(x_0), \dots, \rf(x_n))\) is a multivariate Gaussian vector
with well known conditional distribution \(\rf(x_n)\) given \((\rf(x_0),\dots,
\rf(x_{n-1}))\) when the evaluation locations are deterministic.  But
\(\rf(X)\) is not necessarily Gaussian if \(X\) is random\footnote{
    consider \(X= \argmin_{x\in K}\rf(x)\) for some compact set \(K\subseteq
    \domain\).
} and the calculation of conditional distributions becomes much more
difficult. Treating previsible inputs as deterministic ensures the
calculation is feasible but it lacks theoretical foundation.

We formalize the hope, that previsible evaluation locations may be treated as
deterministic as follows.
Let \((\kernel_{x_{[0\inter n]}})_{x_{[0\inter n]}\in \domain^{n+1}}\) be a
collection of regular conditional distributions for \(\rf(x_n)\) given
\(\rf(x_{[0\inter n)})=(\rf(x_0), \dots, \rf(x_{n-1}))\) \autocite[e.g.][Def.\@
8.28]{klenkeProbabilityTheoryComprehensive2014} indexed by the evaluation
locations \(x_{[0\inter n]} = (x_0,\dots, x_n)\), where we use the following notation
for discrete intervals
\[
    \tag{discrete intervals}
    [i\inter j] := [i,j] \cap \integer,
    \qquad
    [i\inter j) := [i,j)\cap \integer,
    \qquad
    \text{etc.}
\]
For all locations \(x_{[0\inter n]}\) we thus
have for all measurable sets \(A\)
\[
    \Pr\bigl(\rf(x_n) \in A \mid \rf(x_{[0\inter n)})\bigr)
    \overset{\as}= \kernel_{x_{[0\inter n]}}\bigl(\rf(x_{[0\inter n)}); A\bigr).
\]
Recall that this collection is easy to come by in the Gaussian case.
For previsible locations \(X_{[0\inter n]}\)
the \emph{hope} is therefore that for all measurable sets \(A\)
\begin{equation}
    \label{eq: hope equation}    
    \Pr\bigl(\rf(X_n) \in A \mid \rf(X_{[0\inter n)})\bigr)
    \overset\as= \kernel_{X_{[0\inter n]}}\bigl(\rf(X_{[0\inter n)}); A\bigr).
\end{equation}
In practice, this hope is often naively treated as self-evident \autocite[e.g.][Lemma 5.1, p.
3258]{srinivasInformationTheoreticRegretBounds2012}.  But while the
collection of probability kernels \((\kernel_{x_{[0\inter
n]}})_{x_{[0\inter n]} \in \domain^{n+1}}\) may be treated as a function in
\(x_{[0\inter n]}\), there is no guarantee this function is even measurable.
This means that the term on the right in \eqref{eq: hope equation} might not
even be a well defined random variable, let alone satisfy the equation (cf.\@
Example \ref{ex: inconsistent joint conditional distribution}).

\paragraph*{Outline} In Section \ref{sec: main results}
we introduce the concepts needed to formalize the treatment of previsible
random variables as deterministic and state our main results for
continuous random functions. In Section \ref{sec: previsible sampling}
we prove the building blocks for our main results, which are then used
in Section \ref{sec: conditionally independent sampling} to state
and prove our main result with the additional generalization to conditionally
independent evaluation locations and noisy evaluations. Section \ref{sec:
topological foundation} is concerned with the topological foundations that ensure
the evaluation function \(e(f,x)=f(x)\) is measurable on the space of continuous
functions; and the limitations of this approach.

 \section{Main results}
\label{sec: main results}

To ensure that we may plug random variables into the index of a collection of
regular conditional distributions, we introduce the concept of a `joint'
probability kernel.

\begin{definition}[Joint probability kernels]
    \label{def: joint and consinstent}
    A probability kernel \(\kernel\) is a \emph{joint probability kernel}
    for the collection \((\kernel_x)_{x \in I}\) of probability
    kernels with index set \(I\), if for all \(x\in I\)
    \[
        \kernel_x(\omega;A) = \kernel(\omega,x; A)  \qquad \forall \omega, A.
    \]
    We call \(\kernel\) a \emph{joint conditional distribution},
    if the collection of probability kernels \((\kernel_x)_{x \in I}\) is a
    collection of regular conditional distributions.
\end{definition}

Due to the measurability of a probability kernel for fixed sets \(A\), the
existence of a joint conditional distribution ensures that the object in
\eqref{eq: hope equation} is a well defined random variable.

In the following, we provide sufficient conditions for a joint regular
conditional distribution to be \emph{consistent} -- a term we use informally to
describe the setting in which random evaluation locations can be treated as if
they were deterministic, as conjectured in \eqref{eq: hope equation}.

The meaning of `consistent' is often clear from context and should therefore help
with the interpretation of our results. Nevertheless there are slightly
different requirements (e.g.\@ measurable/previsible/conditionally independent) on the random evaluation locations in every case,
so we accompany any occurrence of the term with a clarification of its precise meaning.

\begin{example}[A joint conditional distribution that is not consistent]
    \label{ex: inconsistent joint conditional distribution}
    Consider a standard uniform random variable \(U\sim \uniform(0,1)\), an
    independent standard normal random variable \(Y\sim \normal(0,1)\) and
    define \(\rf(x) = Y\) for all \(x\in \domain:=[0,1]\). As a constant
    function \(\rf\) is clearly a continuous Gaussian random function. Let
    \[
        \kernel(u, x; B) := \Pr_{Y}(B)
        \qquad
        \tilde\kernel(u,x; B):= \Pr_Y(B) \ind_{U\neq x} + \delta_0(B) \ind_{U=x}.
    \]
    Then clearly for all \(x\in \domain\)
    \[
        \kernel(U, x; B)
        \overset{\as}= \Pr(\rf(x)\in B \mid U)
        \overset{\as}= \tilde\kernel(U,x; B)
    \]
    and thereby both \(\kernel\) and \(\tilde\kernel\) are joint regular probability
    kernels for \(\rf(x)\) conditioned on \(U\). But \(\tilde\kernel\)
    is not consistent, because for most measurable sets \(B\)
    \[
        \Pr(\rf(U) \in B \mid U)
        = \Pr_Y(B) \neq \delta_0(B)
        = \tilde\kernel(U,U; B),
    \]
    even though \(U\) is clearly measurable with respect to \(U\).
\end{example}

Observe that for any fixed \(x\), the kernels in the example above coincide
almost surely. That is, they coincide up to a null set \(N_x\), specifically
\(N_x=\{U=x\}\). But the union of these null sets over all possible values of
\(x\) is not a null set. Joint null set can ensure that the consistency of
one kernel implies the other. A sufficient criterion for such a joint null set is a notion
of continuity.

The following result implies that a continuous, joint conditional distribution
is consistent and such a conditional distribution exists.

\begin{theorem}[Consistency for dependent evaluations \(\rf(x)\)]
    \label{thm: dependent consistent}
	Let \(\filt\) be a sub \(\sigma\)-algebra of the underlying
    probability space \((\Omega, \cA, \Pr)\) and \(\rf\) a random variable in
    the space of continuous function \(C(\domain,\range)\), with
    locally compact, separable metrizable domain \(\domain\) and Polish co-domain \(\range\).
	Then there exists a consistent joint conditional distribution \(\kernel\)
	for \(\rf(x)\) given \(\filt\), which means
	\[
		\Pr(\rf(X) \in B \mid \filt)(\omega)
		= \kernel(\omega, X(\omega); B)
	\]
	for all \(\filt\)-measurable \(X\).

    Furthermore \(x\mapsto \kernel(\omega, x; \cdot)\) is continuous with
	respect to the weak topology on the space of measures for all \(\omega\).
	Let \(\tilde\kernel\) be another joint conditional distribution for
	\(\rf(x)\) given \(\filt\), which is continuous in this sense. Then there
	exists a joint null set \(N\) such that for all \(\omega\in N^\complement\),
	all \(x\in \domain\) and all borel sets \(B\in \borel(\range)\)
	\[
		\kernel(\omega,x; B) = \tilde\kernel(\omega,x; B).
	\]
	In particular, \(\tilde\kernel\) is also a consistent joint conditional distribution.
\end{theorem}

This Theorem will be proven as a special case of Theorem \ref{thm: consistent cond.
distribution for dependent}. 
The following example demonstrates its use in the optimization of random
functions.

\begin{example}[Conditional minimization]
    \label{ex: conditional minimization}
    Assume the setting of Theorem \ref{thm: dependent consistent} and that
    \(\rf(x)\) is integrable. Using the consistent joint conditional
    distribution \(\kernel\) there exists a measurable function
    \(
        H(\omega, x) := \int y\, \kernel(\omega, x; dy)
    \)
    such that by disintegration \autocite[e.g.][Thm.\@
    6.4]{kallenbergFoundationsModernProbability2002} we have for
    \(\Pr\)-almost all \(\omega\) and all \(x\in \domain\)
    \begin{align*}
        \E[\rf(x) \mid \filt](\omega)
        = H(\omega, x)
        \quad\text{and}\quad
        \E[\rf(X) \mid \filt](\omega)
        = H(\omega, X(\omega))
    \end{align*}
    for any \(\filt\)-measurable \(X\). That is, we can treat \(\filt\)-measurable
    random variables \(X\) in \(\domain\) as if they were deterministic inputs
    to \(\E[\rf(x)\mid \filt]\). This implies for any such \(X\)
    \begin{equation}
        \inf_{x\in \domain} \E[\rf(x) \mid \filt]
\le  \E[\rf(X) \mid \filt].
    \end{equation}
    And if \(X_* := \argmin_{x\in \domain} \E[\rf(x) \mid \filt]\) is a
    \(\filt\)-measurable random variable,\footnote{
        This is a non-trivial problem by itself \autocite[see e.g.][Thm.\@
        18.19]{aliprantisInfiniteDimensionalAnalysis2006}.
    } we have
    \begin{equation}
        \inf_{x\in \domain} \E[\rf(x) \mid \filt]
        = \inf_{\substack{\text{\(X\) \(\filt\)-meas.}\\ \text{r.v. in \(\domain\)}}} \E[\rf(X) \mid \filt]
        = \E[\rf(X_*) \mid \filt].
    \end{equation}
\end{example}

So far, the random function evaluation \(\rf(x)\) only occurred as a dependent
variable. In the following we analyze the case where \(\rf(x)\) is conditioned
on. While consistency was an issue when \(\rf(x)\) is a dependent variable, it
turns out that every joint conditional distribution is consistent when
\(\rf(x)\) is conditioned on.

\begin{definition}[Previsible]
    \label{def: previsible}
    The \emph{previsible setting} is
    \begin{itemize}[noitemsep]
        \item an underlying probability space \((\Omega, \cA, \Pr)\),
        \item a sub-\(\sigma\)-algebra \(\filt\) (the `initial
        information') with \(W\) a random element such that \(\filt =
        \sigma(W)\),\footnote{
            There always exists such a random element \(W\) since the identity
            map from the measurable space \((\Omega, \cA)\) into \((\Omega,
            \filt)\) is measurable and clearly generates \(\filt\).
        }
        \item \(\rf\) a random function in the space of continuous functions
        \(C(\domain, \range)\), where \(\domain\) is a locally compact,
        separable metrizable space and \(\range\) a polish space.
    \end{itemize}
    A sequence \(X=(X_n)_{n\in \nat_0}\) of random evaluation locations in
    \(\domain\) is called \emph{previsible}, if \(X_{n+1}\) is measurable with
    respect to
    \[
        \filt_n^X := \sigma(\filt, \rf(X_0), \dots, \rf(X_n))
        \qquad \text{for } n\ge -1.
    \]
\end{definition}

\begin{theorem}[Previsible sampling]
    \label{thm: previsible sampling}
    Assume the previsible setting (Def.\@ \ref{def: previsible}).
    \begin{enumerate}[label={(\roman*)}]
        \item 
        Let \(Z\) be a random variable in a standard Borel space
        \((E,\borel(E))\) and \(\kernel\) a joint conditional distribution for
        \(Z\) given \(\filt, \rf(x_0), \dots, \rf(x_n)\), i.e.
        \[
            \Pr(Z\in A \mid \filt, \rf(x_0), \dots, \rf(x_n))
            \overset{\as}= \kernel(W, \rf(x_0), \dots, \rf(x_n), x_{[0\inter n]}; B)
        \]
        for all \(x_{[0\inter n]}\in \domain^{n+1}\) and \(A\in\borel(E)\). Then
        \(\kernel\) is consistent, i.e. for all previsible sequences
        \((X_k)_{k\in \nat_0}\) and all \(A\in\borel(E)\)
        \[
            \Pr(Z\in A \mid \filt_n^X)
            \overset\as=
            \kernel(W, \rf(X_0), \dots, \rf(X_n), X_{[0\inter n]}; A).
        \]
        \item Let \(\kernel\) be a joint conditional distribution for
        \(\rf(x_n)\) given \(\filt, \rf(x_{[0\inter n)})\) such that
        \[
            x_n \mapsto \kernel(y_{[0\inter n)}, x_{[0\inter n]}; \,\cdot\,)
        \]
        is continuous with respect to the weak topology on the space of measures
        for all \(x_{[0\inter n)}\in \domain^n\) and \(y_{[0\inter n)}\in
        \domain^n\).

        Then \(\kernel\) is consistent, i.e.\@ for all previsible \((X_k)_{k\in
        \nat_0}\) and \(B\in \borel(\range)\)
        \[
            \Pr(\rf(X_n) \in B \mid \filt_{n-1}^X)
            \overset\as=
            \kernel(W, \rf(X_{[0\inter n)}), X_{[0\inter n]}; B).
        \]
    \end{enumerate}
\end{theorem}

Theorem \ref{thm: previsible sampling} is proven as a special case of
Theorem \ref{thm: conditionally independent sampling}. There we allow \(X_{n+1}\)
to be random, conditional on \(\filt_n\) and only require it to be independent from
\(\rf\) conditional on \(\filt_n\). This result also covers noisy evaluations of \(\rf\).

We want to highlight that continuity of the kernel is only required for the case
where function evaluations are dependent variables. While consistency
is therefore never an issue, the existence of such a joint
conditional distribution is uncertain in general. However in the Gaussian case,
the joint conditional distribution is known explicitly.

\begin{example}[Gaussian case]
    \label{ex: gaussian conditional}
    Let \(\rf=(\rf(x))_{x\in \domain}\)
    be a Gaussian random function with mean and covariance functions 
    \[
        \mu_0(x) = \E[\rf(x)]
        \qquad\text{and}\qquad
        \C_\rf(x,y) = \Cov(\rf(x), \rf(y)).
    \]
    The conditional distribution of \(\rf(x_n)\) given
    \(\rf(x_{[0\inter n)})\) is the conditional distribution of a multivariate
    Gaussian random vector \(\rf(x_{[0\inter n]})\), which is well known to be
    \(\normal\bigl(\mu_n(x_{[0\inter n]}, \rf(x_{[0\inter n)})),
    \Sigma_n(x_{[0\inter n]})\bigr)\),\footnote{
        e.g.\@ Theorem \ref{thm: conditional gaussian distribution}
        or \textcite[e.g.][Prop.\@ 3.13]{eatonMultivariateStatisticsVector2007}
    } with
    \begin{align*}
        \mu_n(x_{[0\inter n]}, y_{[0\inter n)})
        &:= \mu_0(x_n) + \sum_{i,j=0}^{n-1} \C_\rf(x_n,x_i) \bigl[\Sigma_0(x_{[0\inter n)})^{-1}\bigr]_{ij}\bigl(y_j - \mu_0(x_j)\bigr)
        \\
        \Sigma_n(x_{[0\inter n]})
        &:= \Sigma_0(x_n) - \sum_{i,j=0}^{n-1} \C_\rf(x_n, x_i) \bigl[\Sigma_0(x_{[0\inter n)})^{-1}\bigr]_{ij} \C_\rf(x_j, x_n),
    \end{align*}
    where \(\Sigma_0(x_{[0\inter n)})_{ij} := \Cov(\rf(x_i), \rf(x_j))\).
    This induces the joint conditional distribution
    \begin{equation}
        \label{eq: consistent gaussian joint conditional distribution}
        \kernel(y_{[0\inter n)}, x_{[0\inter n]}; B)
        \propto \int_B \exp\Bigl(- \tfrac12(t-\mu_n)^\transpose \Sigma_n(x_{[0\inter n]})^{-1}(t-\mu_n)\Bigr)
        dt
    \end{equation}
    with \(\mu_n := \mu_n(x_{[0\inter n]}, y_{[0\inter n)})\). Now if \(\rf\) is a
    \emph{continuous} Gaussian random function, Theorem \ref{thm: previsible sampling}
    \ref{it: consistent without dependent variable} is immediately applicable.
    For the applicability of \ref{it: consistent with dependent variable} observe 
    that a continuous Gaussian random function must have continuous mean
    \(\mu_0\) and covariance \(\C_\rf\) \autocite[e.g.][Thm.\@
    3]{costaSamplePathRegularity2024,talagrandRegularityGaussianProcesses1987}.
    And the continuity of \(\mu_0\) and \(\C_\rf\) is sufficient for \(\mu_n\)
    and \(\Sigma_n\) to be continuous in \(x_n\). This implies the
    characteristic function of the joint conditional distribution
    \[
        \hat\kernel(y_{[0\inter n)}, x_{[0\inter n]}; t)
        = \exp\Bigl(it^\transpose \mu_n(x_{[0\inter n]}, y_{[0\inter n)}) - \tfrac12 t^\transpose\Sigma_n(x_{[0\inter n]}) t\Bigr)
    \]
    is continuous in \(x_n\), which implies continuity of \(\kernel\) in the
    weak topology by Lévy's continuity theorem \autocite[e.g.][Thm.\@
    5.3]{kallenbergFoundationsModernProbability2002}.
\end{example}

\begin{corollary}[Gaussian previsible sampling]
    \label{cor: gaussian previsible sampling}
    For a continuous Gaussian random function \(\rf\), \eqref{eq: consistent
    gaussian joint conditional distribution} is a consistent joint probability kernel
    for \(\rf\) in the sense of Theorem \ref{thm: previsible sampling}.
\end{corollary}

 \section{Previsible sampling}
\label{sec: previsible sampling}

In this section we establish the building blocks for our main results.
As the analysis of multiple function evaluations will ultimately proceed via
induction, we restrict our attention here to the case of a single function
evaluation. This section thereby lays the groundwork for the most general
form of our results, presented in Section \ref{sec: conditionally independent
sampling}.

Throughout this section any \(\sigma\)-algebra is implicitly
assumed to be a sub \(\sigma\)-algebra of the underlying probability space
\((\Omega, \cA, \Pr)\). \(\domain\) always denotes a locally
compact, separable metrizable space and \(\range\) a
polish space.

\begin{theorem}[Consistency for dependent \(\rf(x)\)]
	\label{thm: consistent cond. distribution for dependent}
	Let \(\filt\) be a \(\sigma\)-algebra, \(Z\) a random variable in the standard
	borel space \((E, \borel(E))\) and \(\rf\) a random variable in \(C(\domain,
	\range)\).
	Then there exists a consistent joint conditional distribution \(\kernel\)
	for \(Z,\rf(x)\) given \(\filt\). That is
	\[
		\Pr( Z, \rf(X) \in B \mid \filt)(\omega)
		= \kernel(\omega, X(\omega); B)
	\]
	for all \(\filt\)-measurable \(X\).

	Furthermore \(x\mapsto \kernel(\omega, x; \,\cdot\,)\) is continuous with
	respect to the weak topology on the space of measures for all \(\omega\in \Omega\). If \(\tilde\kernel\)
	is another joint conditional distribution for \(Z, \rf(x)\) given \(\filt\),
	that is continuous in this sense, then there exists a joint null set \(N\)
	such that for all \(\omega\in N^\complement\), all \(x\in \domain\) and all
	borel sets \(B\)
	\[
		\kernel(\omega,x; B) = \tilde\kernel(\omega,x; B).
	\]
	In particular, \(\tilde\kernel\) is also a consistent joint conditional distribution.
\end{theorem}

\begin{remark}[Existence of consistent joint conditional distribution]
	\label{rem: existence of consistent joint cond. dist}
	Note that for the existence of a consistent joint conditional distribution
	we only require a regular conditional distribution for \(Z,\rf\) given
	\(\filt\) to exist and measurability of the evaluation map \(e\). This part
	of the result can therefore be made to hold with greater generality.
\end{remark}

\begin{proof}
	Observe that \(E\times C(\domain, \range)\) is a standard borel space since \(C(\domain,
	\range)\) is Polish (Theorem~\ref{thm: continuous function}). There
	therefore exists a regular conditional probability distribution
	\(\kernel_{Z,\rf \mid \filt}\) \autocite[e.g.][Thm.\@ 6.3]{kallenbergFoundationsModernProbability2002}.
	Using this probability kernel, we define the kernel
	\[
		\kernel(\omega, x; B)
		:= \int \ind_B(z, e(f,x)) \kernel_{Z, \rf\mid \filt}(\omega; dz \otimes df)
	\]
	which is a measure in \(B\in \borel(E)\otimes \borel(\range)\) by linearity of the integral,
	so we only need to prove it is measurable in \((\omega, x)\in
	\Omega\times\domain\) to prove it is a probability kernel.
	This follows from measurability of the evaluation function \(e\)
	(Theorem~\ref{thm: continuous function}) and the application of Lemma 14.20
	by \citet{klenkeProbabilityTheoryComprehensive2014}
	to the probability kernel \(\tilde\kernel(\omega, x; A) := \kernel_{Z,\rf\mid \filt}(\omega; A)\)
	in the equation above. By `disintegration'
	\autocite[e.g.][Thm~6.4]{kallenbergFoundationsModernProbability2002}
	this probability kernel is moreover a regular conditional version of
	\(\Pr(Z,\rf(X)\in B\mid \filt)\) for all \(\filt\)-measurable \(X\), i.e. for all \(B\in \borel(E)\otimes
	\borel(\range)\) and for \(\Pr\)-almost all \(\omega\)
	\begin{align*}
		\Pr\bigl(Z, \rf(X) \in B \mid \filt\bigr)(\omega)
		\overset{\text{disint.}}&= \int \ind_B(z, e(f,X(\omega))) \kernel_{Z,\rf\mid \filt}(\omega; dz \otimes df)
		\\
		\overset{\text{def.}}&= \kernel(\omega, X(\omega); B).
	\end{align*}
	The kernel is thereby consistent (and a joint kernel since the constant map
	\(X\equiv y\) is \(\filt\) measurable).
	
	For continuity observe that we have \(\lim_{x\to y} (z,f(x)) = (z,f(y))\)
	for any \(f\in C(\domain, \range)\). For open \(U\) this implies
	\[
		\liminf_{x\to y} \ind_U(z,f(x)) \ge \ind_U(z,f(y)),
	\]
	because if \((z,f(y)) \in U\), then eventually \((z,f(x))\) in \(U\) due to
	openness of \(U\). An application of Fatou's lemma \autocite[e.g.][Thm.\@ 4.21]{klenkeProbabilityTheoryComprehensive2014}
	yields for all open \(U\)
	\[
		\liminf_{x\to y}\kernel(\omega, x; U)
		\ge \int \liminf_{x\to y} \ind_U(z,f(x)) \kernel_{Z,\rf\mid \filt}(\omega; dz\otimes df)
		\ge \kernel(\omega, y; U).
	\]
	And we can conclude weak convergency by the Portemanteau theorem
	\autocite[Thm.\@ 13.16]{klenkeProbabilityTheoryComprehensive2014} since \(E\times \range\)
	is metrizable.

	Let \(\tilde\kernel\) be another continuous joint probability kernel.
	Since \(E\times \range\) is second countable, there is a countable base
	\(\{U_n\}_{n\in \nat}\) of its topology, which generates the Borel
	\(\sigma\)-algebra \(\borel(E) \otimes \borel(\range)\). And since \(\domain\)
	is separable, it has a countable dense subset \(Q\). 
	There must therefore exist a zero set \(N\) such that
	\[
		\kernel(\omega, q; U_n)
		= \tilde\kernel(\omega, q; U_n),
		\qquad \forall \omega \in N^\complement,\ n\in \nat,\ q\in Q,
	\]
	because both kernels are regular conditional version of \(\Pr(Z,\rf(q)\in U_n; \cG)\)
	and the union over \(\nat\times Q\) is a countable union. Since
	\(\{U_n\}_{n\in \nat}\) generates the \(\sigma\)-algebra, we deduce for all
	\(\omega \in N^\complement\) and all \(q\in Q\) that \(\kernel(\omega,q;
	\cdot) = \tilde\kernel(\omega, q; \cdot)\).
	As \(Q\) is dense in \(\domain\) we have by continuity of the joint kernels
	for all \(\omega\in N^\complement\) and all \(x\in\domain\)
	\[
		\kernel(\omega,x; \cdot) 
		= \tilde\kernel(\omega, x; \cdot).
		\qedhere
	\]
\end{proof}

In Example \ref{ex: inconsistent joint conditional distribution} we showed a
joint probability distribution to exist, which is not consistent. In this example,
the random function was a dependent variable. In Theorem \ref{thm: consistent
cond. distribution for dependent} we gave a sufficient condition for a 
unique continuous and consistent conditional distributions to exist in this
case where the function value is the dependent variable. It is now time to
consider the case where functions evaluated at random points are
conditional variables. The following result shows that
we never have to worry about the consistency of probability kernels where the
function value is a conditional variable, if a consistent joint probability
kernel exists for function values as dependent variables.

\begin{prop}[Consistency shuffle]\label{prop: consistency shuffle}
	Let \((\Omega, \cA, \Pr)\) be a probability space and let \(\xi_1,
	(\xi_2^{y})_{y\in D},\xi_3\) be random variables in the measurable spaces
	\((E_i, \cE_i)\) with measurable domain \((D, \mathcal D)\).
	\begin{quote}
	If there exists a consistent joint probability kernel \(\magenta{\kappa_{3,2\mid 1}}\)
	for \(\xi_3, \xi_2^y\) given \(\xi_1\), then \textbf{any} joint probability kernel
	\(\teal{\kappa_{3\mid 2,1}}\) for \(\xi_3\) given \(\xi_2^y, \xi_1\) is
	consistent.
	\end{quote}
	Formally, assume there exists a consistent probability kernel
	\(\magenta{\kappa_{3,2\mid 1}}\) for \(\xi_3, \xi_2^y\) given \(\xi_1\), i.e.\@
	for all \(A\in \cE_3\otimes \cE_2\) and all measurable functions \(g\colon E_1 \to D\)
	\begin{equation}
		\label{eq: consistent outer kernel}
		\Pr\bigl(\xi_3, \xi_2^{g(\xi_1)} \in A \mid \xi_1\bigr)
		\overset{\as}= \magenta{\kappa_{3,2\mid 1}}(\xi_1, g(\xi_1); A),
	\end{equation}
	where we assume \(\xi_2^{g(\xi_1)}\) is a random variable, i.e. measurable.
	Then if there exists a joint conditional probability kernel
	\(\teal{\kappa_{3\mid1,2}}\) for \(\xi_3\) given \(\xi_1 ,\xi_2^{y}\) such
	that for all \(y\in D\) and \(A_3\in \cE_3\)
	\begin{equation}
		\label{eq: joint inner kernel}
		\Pr(\xi_3\in A_3 \mid \xi_1, \xi_2^{y})
		\overset{\as}= \teal{\kappa_{3\mid 2, 1}}(\xi_1, \xi_2^{y}, y; A_3),
	\end{equation}
	then \(\teal{\kappa_{3\mid 2, 1}}\) is consistent, i.e. we have for all \(A_3\in \cE_3\)
	and measurable \(g\colon E_1 \to D\)
	\[
		\Pr(\xi_3 \in A_3 \mid \xi_1, \xi_2^{g(\xi_1)})
		\overset{\as}= \teal{\kappa_{3\mid2,1}}\bigl(\xi_1, \xi_2^{g(\xi_1)}, g(\xi_1); A_3\bigr).
	\]
\end{prop}
\begin{remark}[Possible generalization]
	Note that we keep \(g\) fixed throughout the proof. So if consistency of \(\magenta{\kappa_{3,2\mid 1}}\)
	only holds for a specific \(g\), then we also obtain consistency of \(\teal{\kappa_{3\mid2,1}}\)
	only for this specific function \(g\). For consistency of \(\teal{\kappa_{3\mid2,1}}\) it is therefore
	sufficient to find a \(\magenta{\kappa_{3,2\mid 1}^g}\) that is only
	consistent w.r.t. \(g\) for each \(g\).
\end{remark}
\begin{proof}
	Let \(g\colon E_1\to D\) be a measurable function.
	By definition of the conditional expectation
	we need to show for all \(A_3\in \cE_3\) and all \(A_{1,2}\in \cE_1 \otimes
	\cE_2\)
	\[
		\E\Bigl[
			\blue{\ind_{A_{1,2}}(\xi_1, \xi_2^{(g(\xi_1))})}
			\teal{\kappa_{3\mid2,1}}\bigl(\xi_1, \xi_2^{(g(\xi_1))}, g(\xi_1); A_3\bigr)
		\Bigr]
		=\E\bigl[\blue{\ind_{A_{1,2}}(\xi_1, \xi_2^{(g(\xi_1))})}\ind_{A_3}(\xi_3)\bigr] 
	\]
	Without loss of generality we may only consider \(A_{1,2}= A_1\times A_2\in \cE_1\times \cE_2\)
	since the product sigma algebra \(\cE_1\otimes \cE_2\) is generated by these
	rectangles. 
	Since
	\[
		\magenta{\kappa_{2\mid 1}^g}(x_1; A_2):= \magenta{\kappa_{3,2\mid 1}}(x_1, g(x_1); E_3\times A_2)
	\]
	is a regular conditional version of \(\Pr(\xi_2^{(g(\xi_1))}\in \cdot \mid \xi_1)\)
	by assumption \eqref{eq: consistent outer kernel}
	we may apply disintegration \autocite[e.g.][Thm. 6.4]{kallenbergFoundationsModernProbability2002} to the 
	measurable function
	\[
		\varphi(x_1, x_2) \mapsto \blue{\ind_{A_2}(x_2)}\teal{\kappa_{3\mid2,1}}(x_1, x_2, g(x_1); A_3)
	\]
	to obtain
	\begin{equation}
		\label{eq: disintegration of h}\begin{aligned}
		\E\bigl[\varphi\bigl(\xi_1, \xi_2^{g(\xi_1)}\bigr)\mid \xi_1\bigr]
		\overset{\as}&= \int \varphi(\xi_1, x_2)\magenta{\kappa_{2\mid 1}^g}(\xi_1; dx_2)
		\\
		\overset{\text{def.}}&= \int \varphi(\xi_1, x_2)\magenta{\kappa_{3,2\mid 1}}(\xi_1, g(\xi_1); E_3 \times dx_2).
	\end{aligned}	
	\end{equation}
	We thereby have
	\begin{align*}
		&\E\Bigl[
			\blue{\ind_{A_1}(\xi_1)\ind_{A_2}(\xi_2^{g(\xi_1)})}
			\teal{\kappa_{3\mid2,1}}\bigl(\xi_1, \xi_2^{g(\xi_1)}, g(\xi_1); A_3\bigr)
		\Bigr]
		\\
		&=\E\Bigl[\blue{\ind_{A_1}(\xi_1)}\varphi\bigl(\xi_1, \blue{\xi_2^{g(\xi_1)}}\bigr)\Bigr]
		\\
		\overset{\eqref{eq: disintegration of h}}&=
		\E\Bigl[
			\blue{\ind_{A_1}(\xi_1)}
			\int \varphi(\xi_1, \blue{x_2})
			\magenta{\kappa_{3,2\mid 1}}(\xi_1, g(\xi_1); E_3 \times dx_2)\Bigr]
		\\
		\overset{\text{Lemma}~\ref{lem: coupled kernel 3|2,1}}&=
		\E\Bigl[
			\blue{\ind_{A_1}(\xi_1)}
			\magenta{\kappa_{3,2\mid 1}}(\xi_1, g(\xi_1); A_3 \times \blue{A_2})
		\Bigr]
		\\
		\overset{\eqref{eq: consistent outer kernel}}&=
		\E\Bigl[
			\blue{\ind_{A_1}(\xi_1)}
			\ind_{A_3 \times \blue{A_2}}(\xi_3, \xi_2^{g(\xi_1)})
		\Bigr]
		\\
		&= \E\Bigl[
			\blue{\ind_{A_1}(\xi_1)\ind_{A_2}(\xi_2^{g(\xi_1)})}
			\ind_{A_3}(\xi_3)
		\Bigr].
	\end{align*}
	The crucial step is the application of Lemma~\ref{lem: coupled kernel
	3|2,1}, which provides an integral representation of a regular conditional
	distribution of \(\xi_3, \xi_2^y \mid \xi_1\) that \emph{couples} the two
	conditional kernels.

	\begin{lemma}\label{lem: coupled kernel 3|2,1}
		For all \(A_2\in \cE_2\), \(A_3\in \cE_3\), all \(y\in \domain\) and \(\Pr_{\xi_1}\)-almost all \(x_1\)
		\begin{align}
			\label{eq: coupled kernel 3|2,1}
			\magenta{\kappa_{3,2\mid 1}}(x_1, y;A_3 \times A_2)
			&=\int \varphi(x_1,x_2) \magenta{\kappa_{3,2\mid 1}}(x_1, y; E_3 \times dx_2)
			\\
			\nonumber
			&= \int \ind_{A_2}(x_2) \teal{\kappa_{3\mid 2,1}}(x_1, x_2, y; A_3) \magenta{\kappa_{3,2\mid 1}}(x_1, y; E_3 \times dx_2).
		\end{align}
	\end{lemma}
	In the remainder of the proof we will show this Lemma.
	To this end pick any \(A_1 \in \cE_1\). Then by definition of the
	conditional expectation \autocite[e.g.\label{footnote: test}][chap.~8]{klenkeProbabilityTheoryComprehensive2014}
	\begin{align*}
		&\E[\blue{\ind_{A_1}(\xi_1)}\magenta{\kappa_{3,2\mid 1}}(x_1, y; A_3 \times A_2)]
		\\
		\overset{\eqref{eq: consistent outer kernel}}&= \E[\blue{\ind_{A_1}(\xi_1)} \ind_{A_2}(\xi_2^y) \ind_{A_3}(\xi_3)]
		\\
		\overset{\eqref{eq: joint inner kernel}}&=
		\E[\blue{\ind_{A_1}(\xi_1)} \ind_{A_2}(\xi_2^y) \teal{\kappa_{3\mid 2,1}}(\xi_1, \xi_2^y; A_3)]
		\\
		\overset{(*)}&= 
		\E\Bigl[\blue{\ind_{A_1}(\xi_1)} \int \ind_{A_2}(x_2) \teal{\kappa_{3\mid 2,1}}(\xi_1, x_2, y; A_3) \magenta{\kappa_{3,2\mid 1}}(\xi_1, y; E_3 \times dx_2)\Bigr]
	\end{align*}
	Note that the constant function \(g\equiv y\) is always measurable for the application of \eqref{eq: consistent outer kernel}.
	The last step \((*)\) is implied by disintegration \autocite[e.g.][Thm. 6.4]{kallenbergFoundationsModernProbability2002}
	\[
		\E[f(\xi_1, \xi_2^y) \mid \xi_1] \overset{\as}= \int f(\xi_1, x_2) \magenta{\kappa_{2\mid 1}^y}(\xi_1; dx_2)
	\]
	of the measurable function
	\(
		f(x_1, x_2)
		:= \ind_{A_2}(x_2)\teal{\kappa_{3\mid 2,1}}(x_1, x_2; A_3)
	\)
	using the probability kernel
	\[
		\magenta{\kappa_{2\mid 1}^y}(x_1; A_2):= \magenta{\kappa_{3,2\mid 1}}(x_1, y; E_3\times A_2)
	\]
	which is a regular conditional version of \(\Pr(\xi_2^y \in \cdot \mid \xi_1)\) by
	assumption \eqref{eq: consistent outer kernel}, since the constant function
	\(g\equiv y\) is measurable.
\end{proof}

\begin{corollary}[Automatic consistency]
	\label{cor: automatic consistency}
	Let \(Z\) be a random variable in a standard borel space \((E, \borel(E))\),
	\(\rf\) a random variable in \(C(\domain, \range)\).
	Let \(W\) be a random element in an arbitrary measurable space \((\Omega,
	\cF)\). If there exists a joint conditional distribution \(\kernel\) for
	\(Z\) given \(W, \rf(x)\) then
	\(\kernel\) is automatically consistent.
	That is, for all \(B\in \borel(E)\) all \(\sigma(W)\) measurable \(X\)
	\[
		\Pr(Z\in B \mid W, \rf(X)) \overset\as= \kernel(W, \rf(X), X; B),
	\]
\end{corollary}
\begin{proof}
		Since \(X=g(W)\) for some measurable function \(g\), Proposition
		\ref{prop: consistency shuffle} with \((\xi_3, \xi_2^y,
		\xi_1)=(Z,\rf(y), W)\) yields the claim, since a consistent joint
		probability kernel for \(\xi_3, \xi_2^y\) given \(\xi_1\) exists by
		Theorem \ref{thm: consistent cond. distribution for dependent}.
\end{proof}

 \section{Conditionally independent sampling}
\label{sec: conditionally independent sampling}

In this section we will first introduce generalizations to
of our main result (Theorem \ref{thm: previsible sampling}) and then proceed to
prove this more general result (Theorem \ref{thm: conditionally independent
sampling}).

\paragraph*{Conditional independence}
Sometimes \(X_{n+1}\) is not previsible itself, but sampled from a previsible
distribution.  That is, a distribution constructed from previously seen
evaluations (e.g.\@ Thompson sampling \autocite{thompsonLikelihoodThatOne1933}). In
this case,
\(X_{n+1}\) is not previsible, but independent from \(\rf\) conditional on
\(\filt_n\) denoted by \(X_{n+1} \indep_{\filt_n} \rf\).\footnote{
	as introduced in \citet[p.\@
	109]{kallenbergFoundationsModernProbability2002}.
}
Note that conditional independence is almost always equivalent to \(X_{n+1} =
h(\xi, U)\) for a measurable function \(h\), a random (previsible) element
\(\xi\) that generates \(\filt_n\) and a standard uniform random variable \(U\)
independent from \((\rf, \filt_n)\) \autocite[Prop.\@
6.13]{kallenbergFoundationsModernProbability2002}.

\paragraph*{Noisy evaluations} In many optimization applications only
noisy evaluations of the random objective function \(\rf\) at \(x\) may be
obtained. We associate the noise \(\noise_n\) to the \(n\)-th evaluation
\(x_n\), such that the function \(\rf_n = \rf + \noise_n\) returns the \(n\)-th
observation \(Y_n = \rf_n(x_n)\). While the noise may simply be independent, 
identically distributed constants, observe that this framework allows for much
more general location-dependent noise. The only requirement is that
\(\noise_n\) is a continuous function, such that \(\rf_n\) is a continuous
function.

\begin{definition}[Conditionally independent evolution]
    \label{def: cond. indep. evolution}
	The \emph{general conditional independence setting} is given by
    \begin{itemize}[noitemsep]
		\item an underlying probability space \((\Omega, \cA, \Pr)\),
        \item a sub-\(\sigma\)-algebra \(\filt\) (the `initial
        information') with \(W\) a random element such that \(\filt =
        \sigma(W)\),
        \item A sequence \((\rf_n)_{n\in \nat_0}\) of continuous random functions
        in \(C(\domain, \range)\), where \(\domain\) is a locally compact, separable
		metrizable space and \(\range\) a polish space.
		\item a random variable \(Z\) in a standard borel space \((E,\borel(E))\)
		(representing an additional quantity of interest).
    \end{itemize}
	A sequence \(X=(X_n)_{n\in \nat_0}\) of random evaluation locations in \(\domain\) 
	is called a \emph{conditionally independent evolution}, if \(X_{n+1} \indep_{\filt_n^X} (Z, (\rf_n)_{n\in \nat_0})\)
	for the filtration
	\[
		\filt_n^X := \sigma(\filt, \rf_0(X_0),\dots, \rf_n(X_n), X_{[0\inter n]})
		\qquad \text{for \(n\ge -1\).}
	\]
\end{definition}

\begin{theorem}[Conditionally independent sampling]
	\label{thm: conditionally independent sampling}
	Assume the general conditional independence setting (Definition \ref{def: cond. indep. evolution}).
    \begin{enumerate}[label={(\roman*)}]
        \item\label{it: consistent without dependent variable}
        Let \(\kernel\) be a joint conditional distribution for \(Z\) given
        \(\filt, \rf_0(x_0), \dots, \rf_n(x_n)\), i.e. for all \(x_{[0\inter n]}\in
        \domain^{n+1}\) and \(A\in\borel(E)\)
        \[
            \Pr\bigl(Z\in A \mid \filt, \rf_0(x_0), \dots, \rf_n(x_n)\bigr)
            \overset{\as}= \kernel\bigl(W, \rf_0(x_0), \dots, \rf_n(x_n), x_{[0\inter n]}; B\bigr).
        \]
        Then for all conditionally independent evolutions
        \((X_k)_{k\in \nat_0}\) and \(B\in\borel(E)\)
        \[
            \Pr\bigl(Z\in A \mid \filt_n^X\bigr)
            \overset\as=
            \kernel\bigl(W, \rf_0(X_0), \dots, \rf_n(X_n), X_{[0\inter n]}; B\bigr).
        \]
        \item\label{it: consistent with dependent variable}
		Let \(\kernel\) be a joint conditional distribution for
        \(Z,\rf_n(x_n) \mid \filt, (\rf_k(x_k))_{k\in[0\inter n)}\)
        such that
        \[
            x_n \mapsto \kernel(y_{[0\inter n)}, x_{[0\inter n]}; \,\cdot\,)
        \]
        is continuous in the weak topology for all \(x_{[0\inter
        n)}\in \domain^n\) and all \(y_{[0\inter n)}\in \domain^n\).

		Then for all conditionally independent evolutions \((X_k)_{k\in \nat_0}\)
		and \(B\in\borel(\range)\)
        \[
            \Pr\bigl(Z,\rf_n(X_n) \in B \mid \filt_{n-1}^X, X_n\bigr)
            \overset\as=
            \kernel\bigl(W, (\rf_k(X_k))_{k\in [0\inter n)}, X_{[0\inter n]}; B\bigr).
        \]
    \end{enumerate}
\end{theorem}

\begin{remark}[Gaussian case]
    If \((\rf_n)_{n\in \nat_0}\) is a sequence of continuous, joint Gaussian
    random functions, then the same procedure as in Example \ref{ex: gaussian
    conditional} leads to a consistent joint conditional distribution
    in the sense of Theorem \ref{thm: conditionally independent sampling}.
\end{remark}

The proof of \ref{it: consistent without dependent variable} will follow from
repeated applications of the following lemma. \ref{it: consistent with dependent variable}
will follow from \ref{it: consistent without dependent variable} and an application
of Theorem \ref{thm: consistent cond. distribution for dependent}.

\begin{lemma}[Consistency allows conditional independence]
	\label{lem: consistency allows conditional independence}
	Let \(Z\) be a random variable in a standard borel space \((E, \borel(E))\),
	\(\rf\) a continuous random function in \(C(\domain, \range)\). Let \(W\) be
	a random element in an arbitrary measurable space \((\Omega, \cF)\). If
	there exists a joint conditional distribution \(\kernel\) for \(Z\) given
	\(W, \rf(x)\) then for any random variable \(X \indep_W (Z,\rf)\) in
	\(\domain\)
	\[
		\Pr(Z\in B \mid W, X, \rf(X)) = \kernel(W, \rf(X), X; B).
	\]
\end{lemma}

\begin{proof}
	Observe that \(X\) is clearly measurable with respect to \(W^+ := (W,X)\). Our
	proof strategy therefore relies on constructing a joint conditional distribution
	for \(Z\) given \((W^+, \rf(x))\) using \(\kernel\) and apply Corollary \ref{cor: automatic
	consistency}.
	
	Since \(X\) independent from \((Z, \rf)\) conditional on \(W\)
	there exists a standard uniform \(U\sim \uniform(0,1)\) independent from \((W,Z,\rf)\) such
	that \(X = h(W,U)\) for some measurable function \(h\) \autocite[Prop.\@
	6.13]{kallenbergFoundationsModernProbability2002}.
	Since \(U\) is independent from \(W,Z,\rf\) we have by \autocite[Prop.\@ 6.6]{kallenbergFoundationsModernProbability2002}
	\[
		\Pr(Z\in B \mid W, U, \rf(x))
		= \Pr(Z\in B \mid W, \rf(x))
		= \kernel(W, \rf(x), x; B)
	\]
	for all \(x\in \domain\). Since \(\sigma(W,X,\rf(x)) \subseteq \sigma(W, U, \rf(x))\) and
	\((W,\rf(x))\) is measurable with respect to \(\sigma(W,X,\rf(x))\) we therefore have
	\[
		\Pr(Z\in B \mid \underbrace{W, X}_{=W^+}, \rf(x))
		= \kernel(W, \rf(x),x;B)
		=: \kernel^+(\underbrace{W, X}_{=W^+}, \rf(x),x;B),
	\]
	where \(\kernel^+\) is defined as constant in the second input.
	An application of Corollary \ref{cor: automatic consistency} to
	\(\kernel^+\) yields the claim.
\end{proof}

\begin{proof}[Proof of Theorem \ref{thm: conditionally independent sampling}]
	We will prove  \ref{it: consistent without dependent variable} by induction over \(k\in\{0,\dots,n+1\}\).
	For any conditionally independent evolution \((X_n)_{n\in \nat_0}\), the induction claim is
	\begin{equation}
		\label{eq: induction claim, conditionally indep. sampling}
		\begin{aligned}
		&\Pr\bigl(
			Z \in A
			\mid W, (\rf_i(X_i))_{i\in[0\inter k)},
			(\rf_i(x_i))_{i\in[k\inter n]}, X_{[0\inter k)}
		\bigr)
		\\
		&\qquad\quad= \kernel\bigl(W, (\rf_i(X_i))_{i\in [0\inter k)}, (\rf_i(x_i))_{i\in [k\inter n]}, X_{[0\inter k)}, x_{[k\inter n]}; A\bigr).
	\end{aligned}\end{equation}
	The induction start with \(k=0\) is given by assumption and \(k=n+1\) is the claim,
	so we only need to show the induction step \(k\to k+1\). For this purpose we
	want to define \(\tilde W
	= \bigl(W, (\rf_i(X_i))_{i\in [0\inter k)}, (\rf_i(x_i))_{i\in(k\inter n]}, X_{[0\inter k)}\bigr)\)
	and the kernel
	\[
		\tilde\kernel_{x_{(k\inter n]}}\bigl(\tilde W, \rf_k(x_k), x_k; A\bigr)
		:= \kernel\bigl(W, (\rf_i(X_i))_{i\in [0\inter k)}, (\rf_i(x_i))_{i\in [k\inter n]}, X_{[0\inter k)}, x_{[k\inter n]}; A\bigr),
	\]
	which is formally defined for any fixed \(x_{(k\inter n]}\) by mapping the elements of \(\tilde W\) into the right 
	position. By induction \eqref{eq: induction claim, conditionally indep. sampling} we thereby have
	\[
		\Pr(Z \in A \mid \tilde W, \rf(x_k))
		= \tilde\kernel_{x_{(k\inter n]}}(\tilde W, \rf(x_k), x_k; A).
	\]
	We can thereby finish the induction using Lemma \ref{lem: consistency allows conditional independence}
	if we can prove \(X_k\) is independent from \((Z, \rf)\) conditional on \(\tilde W\).
	For this we will use the characterization of conditional independence in
	Proposition~6.13 by \citet{kallenbergFoundationsModernProbability2002}.

	Since \(X_k\) is independent from \((Z,\rf)\) conditionally on \(\filt_{k-1}\) there exists,
	by this Proposition, a uniform random variable \(U\sim \uniform(0,1)\)
	independent from \((Z,\rf, \filt_{k-1})\)
	such that \(X_k = h(\xi, U)\) for some measurable function \(h\) and a random
	element \(\xi\) that generates \(\filt_{k-1}\). Due to \(\filt_{k-1}\subseteq \sigma(\tilde W)\)
	the element \(\xi\) is a measurable function of \(\tilde W\) and therefore \(X_k=
	\tilde h(\tilde W, U)\) for some measurable function \(\tilde h\).
	Since \(U\) is independent from \((Z,\rf,\filt_{k-1})\), it is independent from
	\(\tilde W\) as \(\sigma(\tilde W)\subseteq \sigma(\rf, \filt_{k-1})\). Using Prop.~6.13 from \citet{kallenbergFoundationsModernProbability2002}
	again, \(X_k\) is thereby independent from \((Z, \rf)\) conditional on \(\tilde W\).

	What remains is the proof of \ref{it: consistent with dependent variable}.
    Let \(x_n\) be fixed and define \(\tilde Z=(Z,\rf_n(x_n))\). Since \(\kernel\) is a joint conditional
    distribution for \(Z,\rf_n(x_n)\) given \(\filt, (\rf_k(x_k))_{k\in [0\inter
    n)}\) the kernel
    \[
        \kernel_{x_n}(y_{[0\inter n)}, x_{[0\inter n)}; B) :=  \kernel(y_{[0\inter n)}, x_{[0\inter n]}; B)
    \]
    clearly satisfies the requirements of \ref{it: consistent without dependent variable} and thereby
    \[
        \Pr\bigl(Z,\rf(x_n) \in B \mid \filt_{n-1}^X\bigr)
        = \kernel\bigl(\underbrace{W, (\rf_k(X_k))_{k\in [0\inter n)}, X_{[0\inter n)}}_{=:\tilde W}, x_n; B\bigr)
        =: \tilde\kernel(\tilde W, x_n; B).
    \]
    Since \(\tilde W\) generates \(\filt_{n-1}^X\) we are almost in the setting of
    Theorem \ref{thm: consistent cond. distribution for dependent}, as we have
    continuity in \(x_n\). However, since \(X_n\) is not previsible we have to
    repeat the same trick we used in the proof of Lemma \ref{lem: consistency
    allows conditional independence}. Namely, \(X_n\) is measurable with
    respect to \(W^+ := (\tilde W, X)\) and we will have to construct a
    joint conditional distribution for \(Z,\rf(x_n)\) given \(W^+\).

    Since \(X_n\) is independent from \(Z,\rf\) conditional on \(\tilde W\),
    there exists a standard uniform \(U\sim \uniform(0,1)\) independent from
    \((\tilde W, Z,\rf)\) such that \(X_n = h(\tilde W, U)\) for some measurable
    function \(h\) \autocite[Prop.\@
    6.13]{kallenbergFoundationsModernProbability2002}.  Since \(U\) is
    independent from \((\tilde W, Z, \rf)\), we have by
    \autocite[Prop.\@ 6.6]{kallenbergFoundationsModernProbability2002}
    \[
        \Pr(Z,\rf(x_n) \in B \mid \tilde W, U)
        \overset\as= \Pr(Z,\rf(x_n) \in B \mid \tilde W)
        \overset\as= \tilde\kernel(\tilde W, x_n; B)
    \]
    for all \(x_n\in \domain\). Since \(\sigma(\tilde W,X_n) \subseteq \sigma(\tilde W, U)\) and
	\(\tilde W\) is measurable with respect to \(\sigma(\tilde W,X)\) we therefore have
	\[
		\Pr(Z,\rf(x_n)\in B \mid \underbrace{\tilde W, X_n}_{=W^+})
		\overset\as= \tilde\kernel(\tilde W,x;B)
		=: \kernel^+(\underbrace{\tilde W, X_n}_{=W^+}, x_n;B),
	\]
	where \(\kernel^+\) is defined as constant in the second input. Clearly, by
    definition of \(\kernel^+\) via \(\kernel\), \(\kernel^+\) is continuous in \(x_n\)
    and as a continuous joint conditional distribution it is consistent by
    Theorem \ref{thm: consistent cond. distribution for dependent}.
    This finally implies the claim
    \begin{align*}
        \Pr(Z,\rf(X_n) \in B \mid \smash{\underbrace{\filt_{n-1}, X_n}_{=W^+}})
        \overset\as&= \kernel^+(W^+, X_n; B)
        \\
        &= \kernel(W, (\rf_k(X_k))_{k\in [0\inter n)}, X_{[0\inter n)}, X_n; B).
        \qedhere
    \end{align*}
\end{proof}

 \section{Topological foundation}
\label{sec: topological foundation}

In this section we show the evaluation function to be continuous and therefore
measurable for continuous random functions.
For compact \(\domain\) this result can be collected from various
sources.\footnote{e.g.\@ \textcite[Thm.\@ 4.2.17]{engelkingGeneralTopologyRevised1989}
and \textcite[Thm.\@ 4.19]{kechrisClassicalDescriptiveSet1995}} But we could not
find a reference for the result in this generality, so we provide a proof.

\begin{theorem}[Continuous functions]
    \label{thm: continuous function}
    Let \(\domain\) be a locally compact, separable and metrizable space\footnote{
        \label{footnote: regular second countable}
        technically, we do not need \(\domain\) to be metrizable but only regular
        and second countable, which is equivalent to separability in metrizable
        spaces \autocite[Cor.\@ 4.1.16]{engelkingGeneralTopologyRevised1989}.
        With this definition it is more obvious that a locally compact polish
        space satisfies the requirements, but we will assume the more general
        setting in the proof.
    },
    \(\range\) a polish space and \(C(\domain, \range)\) the space of continuous
    functions equipped with the \emph{compact-open}\footnote{
        For \(K\subseteq \domain\) compact and \(U\subseteq \range\) open, the
        sets
        \[\begin{aligned}[t]
            &M(K,U)
            \\
            &:= \{f\in C(\domain, \range): f(K) \subseteq U\}
        \end{aligned}
        \]
        form a sub-base of the \emph{compact-open} topology \autocite[e.g.][Sec.
        3.4]{engelkingGeneralTopologyRevised1989}. I.e.\@ the compact-open
        topology it is the smallest topology such that all \(M(K,U)\) are open.
        Recall that the set of finite intersections of a sub-base form a base of the topology
        and elements from the topology can be expressed as unions of base elements. 
    } topology. Then
    \begin{enumerate}[label=\normalfont{(\roman*)}]
        \item\label{it: eval function is continuous}
        the evaluation function
        \[
            e\colon \begin{cases}
                C(\domain, \range) \times \domain \to \range
                \\
                (f,x) \mapsto f(x)
            \end{cases}
        \]
        is continuous and therefore measurable.
        \item\label{it: continuous function polish with metric}
        \(C(\domain,\range)\) is a \textbf{polish space}, whose topology is
        generated by the metric
        \[
            \metric(f,g) := \sum_{n=1}^\infty 2^{-n} \frac{\metric_n(f,g)}{1+\metric_n(f,g)}
            \quad \text{with}\quad
            \metric_n(f,g) := \sup_{x\in K_n} \metric_\range(f(x), g(x))
        \]
        for any metric \(\metric_\range\) that generates the topology of \(\range\)
        and any \emph{compact exhaustion}\footnote{
            The set \(\domain\) is \emph{hemicompact} if it can be
            \emph{exhausted by the compact sets} \((K_n)_{n\in \nat}\), which
            means that the compact set \(K_n\) is contained in the interior of
            \(K_{n+1}\) for any \(n\) and \[\textstyle\domain =\bigcup_{n\in \nat} K_n.\]
        } \((K_n)_{n\in \nat}\) of \(\domain\),
        that always exists because \(\domain\) is hemicompact!
        \item\label{it: borel coincides} The Borel \(\sigma\)-algebra of \(C(\domain, \range)\) is equal to
        the restriction of the product sigma algebra of \(\range^\domain\) to
        \(C(\domain, \range)\), i.e. \(\borel(C(\domain, \range)) =
        \borel(\range)^{\otimes \domain}\bigr|_{C(\domain, \range)}\).
    \end{enumerate}
\end{theorem}

\begin{remark}[Topology of pointwise convergence]
    \label{rem: top of pointwise convergence}
    The topology of point-wise convergence ensures that all
    projection mappings \(\pi_x(f) = f(x)\) are continuous. It coincides with
    the product topology \autocite[Prop.
    2.6.3]{engelkingGeneralTopologyRevised1989}.
    Thm.~\ref{thm: continuous function} \ref{it: borel coincides} ensures
    that the Borel-\(\sigma\)-algebra generated by the topology of
    point-wise convergence coincides with the Borel
    \(\sigma\)-algebra generated by the compact-open
    topology.
\end{remark}
\begin{remark}[Construction]
    \label{rem: construction of a cont. random function}
    The main tool for the construction of probability measures, Kolmogorov's
    extension theorem \autocite[e.g.][Sec.\@ 14.3]{klenkeProbabilityTheoryComprehensive2014}, allows for the construction of random measures
    on product spaces. This is only compatible with the product topology, i.e.
    the topology of point-wise convergence. But the evaluation map is generally
    not continuous with respect to this topology \autocite[Prop.\@
    2.6.11]{engelkingGeneralTopologyRevised1989}. \ref{it: borel coincides}
    ensures that this does not pose a problem as long as \(\domain\) and
    \(\range\) satisfy the requirements of Theorem~\ref{thm: continuous
    function} and the constructed random process has a continuous
    version\footnote{cf.\@ \textcite[Thm.
    3]{talagrandRegularityGaussianProcesses1987},
    \textcite{costaSamplePathRegularity2024} and references therein}.
\end{remark}

\begin{remark}[Limitations]
    \label{rem: limits of the result}
    While the compact-open topology can be defined for general topological
    spaces, the continuity of the evaluation map crucially depends on \(\domain\)
    being locally compact \autocite[Thm.\@ 3.4.3 and comments
    below]{engelkingGeneralTopologyRevised1989}. For \(\domain\) and \(\range\)
    polish spaces, this implies \(C(\domain,\range)\) is generally only 
    well behaved if \(\domain\) is locally compact.
\end{remark}

\begin{remark}[Discontinuous case]
    \label{rem: discontinuous case}
    Without continuity it is already difficult to obtain a random
    function \(\rf\) that is almost surely measurable and can be evaluated point-wise. The construction of
    Lévy processes in càdlàg\footnote{
        french: continue à droite, limite à gauche, ``right-continuous with left-limits''
    } space only works on ordered domains such as \(\real\), where
    `right-continuous' has meaning. Typically, discontinuous random
    functions are therefore only constructed as generalized functions in the
    sense of distributions\footnote{
        The set of distributions is defined as the topological dual to a
        set of test functions. In particular, distributions are \emph{continuous}
        linear functionals acting on the test functions. Thereby one may
        hope that Theorem~\ref{thm: continuous function} is applicable, but
        the set of test functions is typically not locally compact (cf.\@ Remark \ref{rem: limits of the result}).
    } that cannot be evaluated point-wise
    \autocite[e.g.][]{schafflerGeneralizedStochasticProcesses2018}. In particular, we
    cannot hope to evaluate generalized random functions at random locations.
\end{remark}

\begin{proof}
    Since \(\domain\) is locally compact, \ref{it: eval function is continuous}
    follows from Proposition 2.6.11 and Theorem 3.4.3. by
    \citet{engelkingGeneralTopologyRevised1989}.
    
    For \ref{it: continuous function polish with metric} let us begin to show that
    \(\domain\) is \textbf{hemicompact/exhaustible by compact sets}. Since the space \(\domain\) is
    locally compact, pick a compact neighborhood for every point. The interiors
    of these compact neighborhoods obviously cover \(\domain\). Since every
    regular, second countable space\footref{footnote: regular second countable} is Lindelöf \autocite[Thm.\@ 3.8.1]{engelkingGeneralTopologyRevised1989}, we can pick a countable
    subcover, such that the interiors of the sequence \((C_i)_{i\in \nat}\) of
    compact sets cover the domain \(\domain\). We inductively
    define a compact exhaustion \((K_n)_{n\in \nat}\) with \(K_1 := C_1\).
    Observe that the set \(K_n\) is covered by the interiors \((\interior{C_i})_{i\in
    \nat}\). Since \(K_n\) is compact, we can choose a finite sub-cover
    \((\interior{C_i})_{i\in I}\) and define \(K_{n+1} := \bigcup_{i \in I} C_i\cup C_{n+1}\).
    Then by definition \(K_n\) is contained in the interior of the compact set
    \(K_{n+1}\) and due to \(C_n \subseteq K_n\) this sequence also covers the space \(\domain\)
    and is thereby a compact exhaustion.

    It is straightforward to check that the metric defined in \ref{it:
    continuous function polish with metric} is a metric, so we will only prove
    this \textbf{metric induces the compact-open} topology.
    \begin{enumerate}[wide,label={(\Roman*)}]
        \item \textbf{The compact-open topology is a subset of the metric topology.}
        \label{it: compact-open is subset of metric top}
        We need to show that the sets \(M(K,U)\) are open with respect to the metric. 
        This requires for any \(f\in M(K,U)\) an \(\epsilon>0\) such that the
        epsilon ball \(B_\epsilon(f)\) is contained in \(M(K,U)\).

        We start by constructing a finite cover of \(f(K)\).
        For any \(x\in K\) there exists \(\delta_x>0\) with 
        \(B_{2\delta_x}(f(x))\subseteq U\) for balls induced by the metric
        \(\metric_\range\) as \(U\) is open. Since \(K\) is compact,
        \(f(K)\subseteq U\) is a compact set covered by the balls
        \(B_{\delta_x}(f(x))\). This yields a finite subcover
        \(B_{\delta_1}(f(x_1)), \dots, B_{\delta_m}(f(x_m))\) of \(f(K)\).

        Using this cover we will prove the following criterion: Any \(g\in C(\domain, \range)\) is in \(M(K, U)\) if
        \begin{equation}
            \label{eq: sup bound}
            \sup_{x\in K} \metric_\range(f(x), g(x))
            < \delta
            := \min\{\delta_1,\dots, \delta_m\}.
        \end{equation}
        For this criterion note that for any \(x\in K\) there exists \(i\in \{1,\dots, m\}\) such that \(f(x) \in B_{\delta_i}(f(x_i))\).
        This implies 
        \[
            \metric(g(x), f(x_i)) \le \metric(g(x), f(x)) + \metric(f(x), f(x_i)) \le 2\delta_i,
        \]
        which implies \(g(K)\subseteq \bigcup_{i=1}^m B_{2\delta_i}(f(x_i))\subseteq
        U\) and therefore \(g\in M(K,U)\).

        Consequently, if there exists \(\epsilon>0\) such that \(g\in B_\epsilon(f)\) implies
        criterion \eqref{eq: sup bound}, then we have \(B_\epsilon(f)
        \subseteq M(K,U)\) which finishes the proof. And this is what we will show.
        Since \(K\) is compact and the interiors of \(K_n\) cover the space, there exists
        a finite sub-cover \(K\subseteq \bigcup_{i\in I}K_i\) and therefore some \(m=\max I\) such that \(K\) is in the
        interior of \(K_m\). By definition of \(\metric_m\) it is thus clearly sufficient
        to ensure \(\metric_m(f,g) < \delta\). And since \(\varphi(x)= \frac{x}{1+x}\) is a strict
        monotonous function \(\epsilon := 2^{-m} \varphi(\delta)\) does the job, since
        \(
            2^{-m}\varphi(\metric_m(f,g))\le \metric(f,g) \le \epsilon
        \)
        implies \(\metric_m(f,g) \le \delta\).

        \item \textbf{The metric topology is a subset of the compact-open topology.}
        Since the balls \(B_\epsilon(f)\) form a base of the metric topology it is sufficient
        to prove them open in the compact-open topology. If 
        for any \(g\in B_\epsilon(f)\) there exists a compact \(C_1,\dots, C_m\subseteq \domain\)
        and open \(U_1,\dots, U_m\subseteq \range\) such that \(g\in \bigcap_{j=1}^m M(C_j, U_j) \subseteq B_\epsilon(f)\),
        then the ball is open since these finite intersections are open sets in
        the compact-open topology and their union over \(g\) remains open. But
        since there exists \(r>0\) such that \(B_r(g)\subseteq B_\epsilon(f)\),
        it is sufficient to prove for any \(r>0\) that there exist compact
        \(C_j\) and open \(U_j\) such that
        \begin{equation}
            \label{eq: compact open subset}
            \textstyle
            g\in V:=\bigcap_{j=1}^m M(C_j, U_j) \subseteq B_r(g).
        \end{equation}
        For this purpose pick \(K_N\) from the compact exhaustion with sufficiently large \(N\) such that
        \(2^{-N} < \frac{r}2\).
        Pick a finite cover \(O_{x_1}, \dots O_{x_m}\) of \(K_N\) from the cover
        \(\{O_x\}_{x\in K_N}\) with \(O_x := g^{-1}(B_{r/5}(g(x)))\)
        and define the sets
        \[
            C_j := \closure{O_{x_j}} \cap K_N\qquad U_j := B_{r/4}(f(x_j)).
        \]
        Clearly the \(U_j\) are open and the \(C_j\) are compact and we will now
        prove they satisfy \eqref{eq: compact open subset}. Observe that \(g\in V\)
        since for all \(j\)
        \[
            g(C_j)\subseteq g(\closure{O_{x_j}}) \subseteq \closure{B_{r/5}(f(x_j))}\subseteq U_j.
        \]
        Pick any other \(h\in V\). Then for all \(x\in K_N\) there exists \(i\)
        such that \(x\in O_{x_i}\subseteq C_i\) and by definition of \(V\) this implies
        \(h(x)\in U_i\) and also \(g(x)\in U_i\) and thereby
        \(\metric_\range(h(x), g(x)) \le r/2\). This uniform bound implies \(\metric_N(h,g) \le r/2\) and therefore
        \[
            \metric(g,h) \le \Bigl(\sum_{n=1}^N 2^{-n} \metric_n(g,h)\Bigr) + 
            \Bigl(\sum_{n=N+1}^\infty 2^{-n}\Bigr)
            \le \metric_N(g,h) + 2^{-N} < r,
        \]
        since \(\metric_n(f,g)\le \metric_N(f,g)\) for \(n\le N\).
        Thus \(h\in B_r(g)\) which proves \eqref{eq: compact open subset}.
    \end{enumerate}
    As \(C(\domain, \range)\) is clearly metrizable, what is left to prove are its
    separability and completeness. Separability could be proven directly similarly to the
    proof of Theorem 4.19 in \citet{kechrisClassicalDescriptiveSet1995} but for
    the sake of brevity this result follows from Theorem 3.4.16 and Theorem
    4.1.15 (vii) by \citet{engelkingGeneralTopologyRevised1989} and the fact that
    \(\domain\) and \(\range\) are second countable. Completeness follows from the fact that any
    Cauchy sequence \(f_n\) induces a Cauchy sequence \(f_n(x)\) for any \(x\) by
    definition of the metric. And by completeness of \(\range\) there must exist a
    limiting value \(f(x)\) for any \(x\). The continuity of \(f\) follows from
    the uniform convergence on compact sets, since every compact set is contained
    in some \(K_n\) from the compact exhaustion (cf. last paragraph in \ref{it:
    compact-open is subset of metric top}).
    
    What is left to prove is \ref{it: borel coincides}. Since the projections are continuous with respect to the compact open
    topology, they are measurable with respect to the Borel-\(\sigma\)-algebra.
    The product sigma algebra, which is the smallest sigma algebra to ensure all
    projections are measuralbe, restricted to the continuous functions is
    therefore a subset of the Borel \(\sigma\)-algebra. To prove the opposite inclusion,
    we need to show that the open sets are contained in the product \(\sigma\)-algebra.
    Since the space is second countable \autocite[Cor.\@
    4.1.16]{engelkingGeneralTopologyRevised1989} and every open set thereby a
    countable union of its base, it is sufficient to check that the open ball
    \(B_\epsilon(f_0)\) for \(\epsilon>0\) and \(f_0\in C(\domain, \range)\) is in
    the product sigma algebra restricted to \(C(\domain, \range)\). But since
    \(B_\epsilon(f_0) = H^{-1}([0,\epsilon))\) with \(H(f):=\metric(f,f_0)\), it is
    sufficient to prove \(H\) is \(\sigma(\pi_x : x\in \domain)\)-\(\borel(\real)\)-measurable,
    where \(\pi_x\) are the projections. \(H\) is measurable if \(H_n(f) =
    \metric_n(f,f_0)\) is measurable, as a limit, sum, etc.\@ \autocite[Thm.\@ 1.88-1.92]{klenkeProbabilityTheoryComprehensive2014} of measurable functions.
    But since \(\domain\) is separable \autocite[Cor.\@ 1.3.8]{engelkingGeneralTopologyRevised1989}, i.e.\@ has a countable dense subset \(Q\),
    we have by continuity of \(f\) and \(f_0\)
    \[
        \metric_n(f, f_0) = \sup_{x\in K_n} \metric_\range(f(x), f_0(x))
        = \sup_{x\in K_n\cap Q} \metric_\range(\pi_x(f), \pi_x(f_0)).
    \]
    Since \(\metric_\range\) is continuous and thereby measurable \autocite[Thm.\@
    1.88]{klenkeProbabilityTheoryComprehensive2014},
    \(H_n\) is measurable as a countable supremum of measurable
    functions \autocite[Thm.\@ 1.92]{klenkeProbabilityTheoryComprehensive2014}.
\end{proof}

  		\printbibliography[heading=subbibliography]
	\end{refsection}

	\begin{refsection}
		\chapter{Distributions over functions}
\label{chap: distributions over functions}

To perform Bayesian optimization it is necessary to choose a prior distribution
over objective functions or at least a family of prior distributions.

A minimal assumption is that every continuous function lies
in the support of the objective function distribution. A slightly
weaker form of such a \emph{universality}\footnote{
    Universality is typically defined as `The Reproducing kernel
    hilbertspace (RKHS) lies dense in the continuous functions'
    \autocite[e.g.][]{micchelliUniversalKernels2006,carmeliVectorValuedReproducing2010},
    but since the RKHS generally lies dense in the support of a random
    function \autocite[Thm.\ 3.6.1]{bogachevGaussianMeasures1998} this definition implies this more intuitive concept.
} assumption will be used in
Chapter \ref{chap: optimization landscape of shallow neural networks}.

But for most applications, such a minimal assumption leaves a
family of distributions that is way too big to make useful statements. 
In this section we therefore introduce uniformity assumptions about the
distribution of the objective function as an effective way to shrink the set of
prior distributions to a manageable level.

A uniformity assumption can be justified in two ways: First, if we do not have
any information about the distribution it seems like a reasonable zero-knowledge
prior to weigh every possible input equally. And second, it may sometimes be
possible to enforce this uniformity, for example via a preprocessing step (cf.
Remark~\ref{rem: randomization via preprocessing}).

But what \emph{is} a uniform prior in the case of functions? We will
present two different approaches: The first will motivate the Gaussian 
assumption by infinite divisibility (Section~\ref{sec: infinite divisibility}).
The second motivates the stationary, isotropy assumption with general input
invariance (Section~\ref{sec: input invariance}). In other words, we will
motivate \emph{stationary, isotropic Gaussian random functions} as a zero-knowledge
prior. Unfortunately this assumptions does not result in a single distribution
but in a family of distributions. In Section \ref{sec: characterization of
covariance kernels} we thus present characterizations of various
input-invariant families. The remaining sections are dedicated to the proof of
the characterization of the non-stationary covariance kernels. These have
not been characterized so far and are motivated by the fact that objective
functions in machine learning are typically not stationary but may be isotropic
cf.\ Section~\ref{sec: random linear model}.

\section{Infinite divisibility}
\label{sec: infinite divisibility}

A naive approach to build a uniform prior over functions is to sample
every value \(N(x)\) of a function \(N\) independent and identically
distributed. Since optimization requires some continuity one may want
to smooth this white noise, perhaps with a Gaussian blur.
Convolution with a mollifier such as the Gaussian blur is integration
and integration requires measurability; but unfortunately, the naive white
noise is almost surely \emph{not} measurable. A more fruitful approach is
to construct white noise as a generalized function with measurability in
mind. Accordingly, consider what properties the integrals
\[
    W(A) = \int \ind_A(x) N(x)dx
\]
\emph{should} have if they were well defined. We can then take these
properties as definition for the measure \(W\), which is a generalized
function.\footnote{ Observe that if we define \(W\), then the noise at a certain
point \(N(x)\)
    is heuristically equal to the limit \(\lim_{\epsilon\to 0}
    \frac{W(\ball_\epsilon(x))}{|\ball_\epsilon(x)|}\), where \(\ball_\epsilon(x)\) is an
    \(\epsilon\) ball around \(x\) and \(|\ball_\epsilon(x)|\) its Lebesgue
    measure (volume). Of course this limit does not actually exist, since \(N\) is
    discontinuous. But this conveys the intuition 
    why a measure can generalize a function (see `generalized function' or
    `Schwarz distribution').
} Since the  values \(N(x)\) are iid we would expect:
\begin{enumerate}[label={(\alph*)}]
    \item\label{it: noise independence} \(W(A)\) and \(W(B)\) to be independent for disjoint \(A\) and \(B\)
    (by independence),
    \item\label{it: noise summability} \(W(A \cup B) = W(A) + W(B)\) for disjoint \(A\) and \(B\) (by
    linearity of the integral),
    \item\label{it: noise identically distributed} \(W(A)\) to be equal in distribution to \(W(B)\) if the volume of \(A\)
    is equal to that of \(B\) (since the \(N(x)\) are identically distributed).
\end{enumerate}
Observe that any set \(A\) of non-zero volume can be subdivided into \(n\)
equally sized disjoint parts, which are then iid distributed by \ref{it: noise independence}
and \ref{it: noise identically distributed}. Hence, the summability
in \ref{it: noise summability} implies the distribution of \(W(A)\) is
infinitely divisible. The Lévy-Khintchine representation
Theorem \autocite[e.g.][Theorem~8.1]{ken-itiLevyProcessesInfinitely1999}
then characterizes the possible distributions of \(W(A)\).

Gaussian white noise \autocite[e.g.][Section 1.4.3 and Section 5.3]{adlerRandomFieldsGeometry2007}, characterized by \(W(A) \sim \normal(0, |A|)\)
in place of \ref{it: noise identically distributed}, may still be
an arbitrary selection from the general set of Lévy white
noise \autocite[e.g.][]{fageotDomainDefinitionLevy2021} (see also `Lévy random
field'); but -- while the Gaussian assumption was a forgone conclusion for
practical reasons (cf.~Section \ref{sec: practical bayesian optimization}) --
the fact that it is picked from a reduced, plausible set of infinitely divisible
distributions may alleviate some concern.

\paragraph*{Smoothing noise} Before we move on, let us quickly touch on
two prominent smoothing approaches of white noise to obtain distributions
over continuous functions.

\begin{example}[Brownian motion \& Lévy process]
    Let \(W\) be Gaussian white noise on \([0,\infty)\), then
    the integrated white noise \(B(t) := W([0,t))\) is the Brownian motion.\footnote{
        see e.g.\@ \textcite[Sec.~1.4.3]{adlerRandomFieldsGeometry2007} including
        a generalization to `Brownian sheets' on \([0,\infty)^d\).
    }
    For general Lévy noise \(W\) defined on \([0,\infty)\), \(W([0,t))\) is a Lévy process which
    the the Lévy-Itô decomposition breaks up into, drift, Brownian motion and a
    jump process. \autocite[e.g.][Theorem 2.4.16]{applebaumLevyProcessesStochastic2009}
    In particular the assumption of continuity eliminates all distributions except for the Gaussian one.
\end{example}

\begin{example}[Convolution with mollifier/Blur]
    \label{ex: convolution with mollifier}
    Let \(W\) be Lévy white noise and assume
    \(\E[W(A)]=0\) and \(\E[W(A)^2] = |A|\), where \(|A|\) is the
    Lebesgue measure of \(A\). Then integrals and thereby convolutions with
    mollifiers are well defined \autocite[Sec.\@
    5.2]{adlerRandomFieldsGeometry2007}
    \[
        \rf(t) = k * W(t) = \int k(t-x)W(dx).
    \]
    The resulting random function \(\rf\) is Gaussian if
    \(W\) is Gaussian noise. \(\rf\) is further always centered with covariance
    \autocite[Lemma 5.3.1]{adlerRandomFieldsGeometry2007}
    \[
        \E[\rf(t)\rf(s)] = \int k(t-x)k(s-x)dx = \int k((t-s) -y)k(y)dy =: C(t-s).
    \]
    If we choose the Gaussian blur \(k(x) = \exp(-\frac{\|x\|^2}{\scale^2})\)
    as the mollifier, then the covariance is given by
    \begin{align*}
        \E[\rf(t)\rf(s)]
        &= \int \exp\Bigl(-\frac{\|t-x\|^2 + \|s-x\|^2}{\scale^2}\Bigr) dx
        \\
        \overset{y=x-\frac{s+t}2}&=
        \int \exp\Bigl(-\frac{\|\frac{t-s}2-y\|^2 + \|\frac{s-t}2-y\|^2}{\scale^2}\Bigr) dy
        \\
        &= \int \exp\Bigl(
            -\frac{
                2\|\frac{t-s}2\|^2 - 2\langle (t-s) + (s-t), y\rangle + 2\|y\|^2
            }{\scale^2}
        \Bigr) dy
        \\
        &= e^{-\frac{\|t-s\|^2}{2\scale^2}} \int e^{-\frac{\|y\|^2}{2(\scale^2/4)}}dy
        \\
        &=\sigma^2 \exp\Bigl(-\frac{\|t-s\|^2}{2\scale^2}\Bigr),
    \end{align*}
    for an appropriate choice of \(\sigma^2\) that may be changed by scaling of
    the mollifier \(k\). This covariance, known as
    the `Squared exponential', `Radial Basis Function' or `Gaussian' kernel,
    plays a crucial role in the characterization of stationary isotropic
    covariance kernels (cf.\ Section \ref{sec: characterization of covariance
    kernels}).

    \begin{remark}[Representation]
        The `spectral representation theorem' \autocite[Theorem 5.4.2]{adlerRandomFieldsGeometry2007}
        gives a characterizing result in the opposite direction: Any continuous stationary
        complex random function \(\rf\) may be obtained by application of the
        mollifier \(k(t,x) = e^{i\langle t,x\rangle}\) to complex \(\nu\)-noise,
        where \(\nu\) is a intensity measure replacing the Lebesgue measure in the definition of white noise.
    \end{remark}
\end{example}

\section{Input invariance}
\label{sec: input invariance}

A different approach to get something akin to a uniform distribution
on the space of measures, is to consider axiomatic invariance properties.
On subsets of \(\real^\dims\) with non-zero volume, the uniform distribution is
defined as the restriction of the Lebesgue measure to these sets. The Lebesgue
measure in turn is characterized by its translation invariance which
encapsulates the notion that any two points are alike.\footnote{
    Translation invariance circumvents the issue that all continuous measures 
    assign zero measure to points. It is therefore not sufficient to require that
    the probability of any two points are alike to define a uniform
    distribution in \(\real^\dims\).
}
Translations in function space is however the wrong operations to demand invariance
from:
\begin{enumerate}[noitemsep]
    \item There does not exist a non-zero translation invariant measure in
    infinite dimensional banach spaces \autocite[Sec.~5.3,
    Thm.~4]{gelfandGeneralizedFunctionsVolume1964}.

    \item Translation in
    function space is an extremely noisy operation which completely changes the
    landscape of the function, because translation entails the addition of an arbitrarily
    rough function.
\end{enumerate}

Invariance with respect to transformations of the input to the function are much
more natural however. If these transformations are bijections, the
minimizer found in the transformed space can be translated back into a minimizer
of the original function and these transformations are thus also amendable to 
preprocessing.\footnote{
    Preprocessing akin to random shuffling requires a uniform distribution over
    the set of transformations. This exists for the orthogonal group
    \(\orthGroup\) but not the group of translations \(\translGroup\) and
    foreshadows that the stationarity assumption may be problematic,
    see Section \ref{sec: random linear model}.
}

\begin{definition}[Distributional \(\Phi\)-invariance]
    \label{def: distributional input invariance}
    We say a random function \(\rf=(\rf(x))_{x\in \domain}\) is (distributionally) \(\Phi\)-invariant to
    a set of transformations \(\phi\colon \domain \to \domain\), typically a group, if\footnote{ Since
    random functions technically have two inputs (cf.~Definition~\ref{def:
    random function}) we should clarify that the transformation \(\phi\)
    only applies to \(x\), i.e. \(\rf \circ \phi(x)  = \rf(\phi(x))\).
    }
    \[
        \Pr_{\rf \circ \phi} = \Pr_\rf \qquad \forall \phi \in \Phi.
    \]
    There are a few prominent examples for \(\domain = \real^\dims\):
    \begin{itemize}
        \item If \(\Phi\) is the set of \emph{translations} \translGroup, we call \(\rf\)
        \textbf{stationary}.

        \item If \(\Phi\) is the set of \emph{linear isometries}, i.e. the
        orthogonal group \orthGroup, we call \(\rf\) (non-stationary) \textbf{isotropic}.

        \item If \(\Phi\) is the set of (Euclidean) \emph{isometries} \euclGroup, we call \(\rf\) \textbf{stationary isotropic}.
    \end{itemize}
    We further say a random function \(\rf\) is \(n\)-weakly \(\Phi\)-invariant, if
    for all \(\phi\in \Phi\), all \(k\le n\) and all \(x_i\)
    \[
        \E[\rf(\phi(x_1))\cdot \dots \cdot \rf(\phi(x_k))] = \E[\rf(x_1)\cdot\dots\cdot\rf(x_k)].
    \]
    Since second moments fully determine Gaussian distributions, \(2\)-weakly invariance
    is special, because it is equivalent to full input invariance in the
    Gaussian case. So an omitted \(n\) equals \(2\).
\end{definition}

\begin{remark}[Almost sure invariance]
    We defined a distributional invariance. The function \(\rf\) itself is
    typically not input invariant to \(\Phi\). Almost sure \(\Phi\)-invariance, i.e.
    \(\rf(\phi(x)) = \rf(x)\) almost surely, results in
    \(n\)-weak invariance of the form
    \[
        \E[\rf(\phi_1(x_1))\cdot \dots \cdot \rf(\phi_2(x_k))] = \E[\rf(x_1)\cdot\dots\cdot\rf(x_k)].
    \]
    for all selections of \(\phi_i\in \Phi\) and \(x_i\in \domain\) and \(k\le n\).
    Since the distributional invariance we introduced only allows for the selection of a
    single transformation \(\phi\) in its weak form, some authors refer to
    distributional invariance as `simultaneous invariance' or `total invariance'
    \autocite[e.g.][]{haasdonkInvariantKernelFunctions2007,brownSampleefficientBayesianOptimisation2025}.
\end{remark}

Different notions of weak invariance have simple
characterizations in terms of the mean and covariance functions.\footnote{
    The mean and covariance function of a random function \(\rf\) are given by
    \(\mu_\rf(x) := \E[\rf(x)]\) and \(\C_\rf(x,y) := \Cov(\rf(x), \rf(y))\)
    respectively
}
We present these in Theorem~\ref{thm: characterization of weak input
invariances} of which the stationary isotropic and stationary case are already
well known.

\begin{theorem}[Functional characterization of weak input invariance]
    \label{thm: characterization of weak input invariances}
    The random function \(\rf\colon\real^\dims \to \real\) is
    \begin{enumerate}[topsep=0pt,itemsep=0pt,partopsep=0pt]
        \item weakly stationary, if and only if there exists \(\mu\in\real\) and
        covariance function \(\ikernel:\real^\dims \to \real\) such that for all \(x,y\)
        \[
            \mu_{\rf}(x) = \mu,
            \qquad \C_{\rf}(x,y) =
            \ikernel(x-y),
        \]

        \item weakly non-stationary isotropic, if and only if there exists a
        mean function \(\mu\colon\real_{\ge 0}\to\real\) and  covariance
        function \(\kernel\colon D\to\real\) with
        \(
            D = \{ \lambda\in \real_{\ge 0}^2 \times \real : |\lambda_3| \le 2\sqrt{\lambda_1\lambda_2}\}
            \subseteq\real^3
        \)  
        such that for all \(x,y\)
        \begin{align*}
            \mu_{\rf}(x) &= \mu\bigl(\tfrac{\|x\|^2}2\bigr)
            \\
            \C_\rf(x,y) &= \kernel\bigl(\tfrac{\|x\|^2}2, \tfrac{\|y\|^2}2, \langle x, y\rangle\bigr),
        \end{align*}
        
        \item\label{item: weak isotropy characterization} 
        weakly stationary isotropic, if and only if there exists \(\mu\in \real\) and a function
        \(\ikernel:\real_{\ge 0}\to \real\) such that for all \(x,y\)
        \[
            \mu_{\rf}(x) = \mu,
            \qquad
            \C_{\rf}(x,y) = \ikernel\bigl(\tfrac{\|x-y\|^2}2\bigr).
        \]
    \end{enumerate}
\end{theorem}

\begin{proof}
	The proof essentially follow as a corollary from a characterization
	of isometries (Proposition~\ref{prop: characterization of isometries}). The
	form of the covariance is proven in the non-stationary isotropic case as
	part of Theorem \ref{thm: isotropy characterization}.\footnote{
        There we prove the covariance is of the form
        \[
            \C_\rf(x,y) = \kernel(\|x\|, \|y\|, \langle x, y\rangle).
        \]
        Through \(\tilde\kernel(t,s, \gamma):= \kernel(\sqrt{2t}, \sqrt{2s}, \gamma)\)
        this is clearly equivalent.
        However, this parametrization is much better suited for derivatives (cf.\
        Section \ref{sec: covariance of derivatives}).
    }
    We leave the other characterizations as an exercise.
\end{proof}

\begin{remark}
    While straightforward to prove, these characterizations are profound: They
    provide information about the covariance of \(\rf(x)\) and \(\rf(y)\)
    \emph{within} the same function \(\rf\) from a uniformity assumption
    \emph{over} possible realizations of \(\rf\). 
\end{remark}

Since weak invariance is equivalent to invariance in the Gaussian case,
the conditions on the mean and covariance above fully characterize the
invariant Gaussian distributions. 
 {
    \renewcommand*{\kernel}{K}
\section{Characterization of covariance kernels}
\label{sec: characterization of covariance kernels}

While any value \(\mu\in \real\)
in the stationary case, or function \(\mu\colon\real_{\ge 0}\to \real\)
in the isotropic case,
is admissible as a mean function,\footnote{
    simply take a centered random function and add this mean
} not every kernel function \(\kappa\) results in a
valid covariance kernel. Specifically, since the variance of random variables
is necessarily positive, the covariance function \(\C_\rf\) of a random function
\(\rf\) necessarily has to satisfy
\[
    0 \le \Var\Bigl(\sum_{i=1}^n a_i \rf(x)\Bigr) = \sum_{i,j=1}^m a_i \conj{a_j}\C_\rf(x_i, x_j)
\]
for all evaluation locations \(x_i\) and (complex valued) scalars \(a_i\).
The covariance \(\C_\rf\) is thus necessarily `positive definite'.
A function \(\kernel\colon \domain \times \domain \to \complex\) is
called a \emph{positive definite kernel}, if for all \(m\in \nat\), distinct
locations \(x_1,\dots, x_m\in \domain\) and \(c=(c_1,\dots,c_m) \in \complex^m\)
\begin{align}
	    \sum_{i,j=1}^m c_i \conj{c_j}\kernel(x_i, x_j) \ge 0.
	\label{eq:positive_definite_definition}
\end{align}
Note that any positive definite kernel is also hermitian \autocite[e.g.][Problem 2.13]{scholkopfLearningKernelsSupport2002}
and it is well known that for any positive definite kernel \(\kernel\) a centered
Gaussian random function \(\rf\) exists with \(\C_\rf = \kernel\).\footnote{
    The proof of this fact is essentially an application of Kolmogorov's
    extension theorem \autocite[e.g.][Thm.\ 14.36]{klenkeProbabilityTheoryComprehensive2014}, since positive definite functions induce consistent
    finite dimensional multivariate Gaussian distributions.
} The positive definite kernels thereby coincide with the covariance kernels.
If the inequality in \eqref{eq:positive_definite_definition} is strict whenever \(c \neq 0\), we call the kernel
\emph{strict positive definite}.\footnote{
    this implies that no linear combination of finitely many random function
    values is deterministic and is a weak (finite) form of `universality'
    \autocite{micchelliUniversalKernels2006}.
}
To make sense of the vast space of positive definite kernels,
characterizations have been identified for classes of kernels invariant to
certain input transformations.  
Corresponding to the notion of distributional
invariance for random functions in Definition \ref{def: distributional input invariance}
we define invariance for kernels in Definition \ref{def: input invariance kernels} and we summarize a list
of such characterizations in Table~\ref{table: characterizations}.

\clearmargin \blockmargin 
\begin{definition}[\(\Phi\)-invariance, kernels]
    \label{def: input invariance kernels}
	Let $\Phi$ be a set of transformations $\phi\colon\domain\to\domain$, typically
    a group.  A kernel \(\kernel\) is called \(\Phi\)-invariant if
    \[
        \kernel(\phi(x), \phi(y)) = \kernel(x,y) \qquad \forall \phi \in \Phi, x,y\in \domain,
    \]
    and we denote the set of \(\Phi\)-invariant kernels by \(\PDset_\Phi:=\PDset_\Phi(\domain)\). And
    for \(\domain=\real^\dims\) the kernel \(\kernel\) is called
    \begin{itemize}[noitemsep]
        \item \emph{stationary}, if \(\Phi\) is the group of translations \(\translGroup\)
        \item \emph{isotropic}, if \(\Phi\) is the orthogonal group \(\orthGroup\) of
        linear isometries\footnote{
            a linear isometry \(U\) is norm preserving, and thus
            for all \(x\)
            \[
                \|x\|^2= \|U x\|^2 = x^\transpose U^\transpose U x
            \]
            This forces \(UU^\transpose =\identity\) such that \(U\) is
            orthogonal. Conversely, any orthogonal matrix
            \(U\) induces a linear isometry.
        }
        \item \emph{stationary isotropic}, if \(\Phi\) is the group of (Euclidean) isometries \(\euclGroup\).
    \end{itemize}
\end{definition}
\unblockmargin[2.7cm]

\begin{table*}
\begin{fullwidth}
    \caption[]{
        \label{table: characterizations}
        Characterizations of continuous pos. def. kernels \(\kernel\colon \domain\times \domain \to \complex\).
        Here \(\Omega_\dims(t):=\E[e^itX]\) for \(X\sim \uniform(\sphere)\),
        \(\uPnorm_n\) are the normalized Gegenbauer polynomial of degree \(n\) and index
        \(\lambda=\frac{\dims-2}2\),
        \(
            \sphere[\infty] := \{x\in \ell^2 : \|x\|=1\}
        \)
        and \(\real^\infty := \ell^2\).
    }
    \centering
    \(\def\arraystretch{1.2}
    \begin{array}{l l l l l}
        \domain & \Phi & \kernel(x,y) & \text{Characterization}
        \\
        \toprule
        \real^\dims
        & \translGroup
        & \kappa(x-y)
        & \kappa(u) = \int e^{i \langle v, u\rangle} d\mu(v)
        & \mu \text{ finite measure}
        \\
        \real^\dims
        & \euclGroup
        & \kappa(\|x-y\|)
        & \kappa(r) = \int_0^\infty \Omega_\dims(rs) \mu(ds)
        & \mu \text{ finite measure}
        \\
        \real^\infty
        & \euclGroup[\infty]
        & \kappa(\|x-y\|)
        & \kappa(r) = \int_0^\infty \exp(-s^2\frac{r^2}2) \mu(ds)
        & \mu \text{ finite measure}
        \\
        \sphere
        & \orthGroup
        & \kappa(x\cdot y)
        & \kappa(t) = \sum_{n=0}^\infty \gamma_n \uPnorm_n(t)
        & \gamma_n\ge 0
        \\
        \sphere[\infty]
        & \orthGroup[\infty]
        & \kappa(x\cdot y)
        & \kappa(t) = \sum_{n=0}^\infty \gamma_n t^n
        & \gamma_n\ge 0
        \\
        X \times \sphere
        & \identity \otimes \orthGroup
        & \kappa(x_X, y_X, x_{\sphere[]} \cdot y_{\sphere[]})
        & \kappa(r,s,t) = \sum_{n=0}^\infty \alpha_n(r,s) \uPnorm_n(t)
        & \alpha_n \text{ pos. definite}
        \\
        X \times \sphere[\infty]
        & \identity \otimes \orthGroup[\infty]
        & \kappa(x_X, y_X, x_{\sphere[]} \cdot y_{\sphere[]})
        & \kappa(r,s,t) = \sum_{n=0}^\infty \alpha_n(r,s) t^n
        & \alpha_n \text{ pos. definite}
        \\
        \real^\dims \times \sphere
        & \translGroup \otimes \orthGroup
        & \kappa(x_\real - y_\real, x_{\sphere[]} \cdot y_{\sphere[]})
        & \kappa(r,s,t) = \sum_{n=0}^\infty \alpha_n(r-s) \uPnorm_n(t)
        & \alpha_n \text{ pos. definite}
        \\
        \real^\dims \times \sphere[\infty]
        & \translGroup \otimes \orthGroup[\infty]
        & \kappa(x_\real - y_\real, x_{\sphere[]} \cdot y_{\sphere[]})
        & \kappa(r,s,t) = \sum_{n=0}^\infty \alpha_n(r-s) t^n 
        & \alpha_n \text{ pos. definite}
        \\
        \midrule
        \real^\dims
        & \orthGroup
        & \kappa(\|x\|, \|y\|, x\cdot y)
        & \kappa(r,s, t) = \sum_{n=0}^\infty \alpha_n(r,s) \uPnorm_n(\tfrac{t}{rs})
        & \substack{
            \text{\(\alpha_n\) pos. definite \&}\\
            \text{zero at origin}
        }
        \\
        \real^\infty
        & \orthGroup[\infty]
        & \kappa(\|x\|, \|y\|, x\cdot y)
        & \kappa(r,s, t) = \sum_{n=0}^\infty \alpha_n(r,s)\bigl(\tfrac{t}{rs}\bigr)^n 
        & \substack{
            \text{\(\alpha_n\) pos. definite \&}\\
            \text{zero at origin}
        }
    \end{array}
    \)
    \end{fullwidth}
	\sideparmargin{outer}
	\vspace{-\baselineskip}
	\sidepar{\vspace{\baselineskip}
        \footnotesize
        Our contribution (cf.~Thm.\@ \ref{thm: isotropy characterization}) is
        summarized in the last two lines of Table \ref{table: characterizations}. 
        The first three results can be found in the book by
        \citet{sasvariMultivariateCharacteristicCorrelation2013}. The very first
        result goes back to \citet{bochnerMonotoneFunktionenStieltjessche1933}
        and can be generalized to groups. Results 2-5 go back to
        \citet{schoenbergMetricSpacesPositive1938,schoenbergPositiveDefiniteFunctions1942}.
        In his honor results 6-9 have been called Schoenberg characterizations \parencite{guellaSchoenbergsTheoremPositive2019, estradeCovarianceFunctionsSpheres2019,bergSchoenbergCoefficientsSchoenberg2017}.
	}
\end{table*}

Since any kernel invariant to a set \(\Phi\) of transformations is also invariant to any subset of
transformations, we have \(\PDset_\Psi \supseteq \PDset_\Phi\) whenever \(\Psi
\subseteq \Phi\). Consequently, the set of stationary isotropic kernels
\(\PDset_{\euclGroup}\) is a subset of \(\PDset_{\translGroup}\) and \(\PDset_{\orthGroup}\).
While the stationary kernels \(\PDset_{\translGroup}\) have long been
characterized \autocite{bochnerMonotoneFunktionenStieltjessche1933}, translation
invariance is typically assumed before invariance with respect to the orthogonal
group is even considered \autocite[e.g.][Sec.
4.2]{rasmussenGaussianProcessesMachine2006}.
This leaves the non-stationary isotropic positive definite kernels without
characterization, which is surprising: 
Virtually all kernel families used in practice are isotropic, up to `geometric
anisotropies' \autocite[p.\@ 50]{steinInterpolationSpatialData1999} that reduce to
isotropy after a change of coordinate system.
In contrast, the generalization of
\(\PDset_{\euclGroup}\) to \(\PDset_{\orthGroup}\) includes many practical
non-stationary kernel families, such as the dot-product kernels
\autocite[e.g.][Sec.\@ 4.2.2]{rasmussenGaussianProcessesMachine2006} the kernels
arising from infinite width neural networks, which are all non-stationary isotropic (cf.
Example \ref{ex: multilayer neural networks are isotropic at initialization}) and the covariance kernel of an objective
function induced by a random linear model (Section \ref{sec: random linear
model}).

Our characterization of non-stationary isotropic kernels, based on the work of
\citet{guellaSchoenbergsTheoremPositive2019}, thereby fills a fundamental gap,
capturing virtually all kernel families used in practice.

We further provide necessary and sufficient conditions for isotropic kernels to be
strictly positive definite (Section~\ref{sec: strict
positive definite}) and provide examples in \ref{sec: applications}.

 \section{Characterization of isotropic kernels on \texorpdfstring{\(\real^{\lowercase{\dims}}\)}{Rd} }\label{sec:positive_definite}

To prove Theorem \ref{thm: isotropy characterization}, we decompose
isotropic kernels $\kernel$ in two scalar and one spherical component.
The spherical component translates to the Gegenbauer polynomials \autocite[see
e.g.][]{reimerGegenbauerPolynomials2003}, which were already used by
\citet{schoenbergPositiveDefiniteFunctions1942} to characterize the
$\orthGroup$-invariant kernels on $\sphere$.
We denote the Gegenbauer  polynomial  of index \(\lambda =
\frac{\dims-2}2\) by \(\uP_n\). 
For the subsequent characterization we use the \emph{normalized
Gegenbauer polynomials} defined by\footnote{
    note that \(|\uP_n(t)|\le \uP_n(1)\) and \(0<\uP_n(1)<\infty\) \autocite[e.g.][Eq.
    (2.16), (2.13)]{reimerGegenbauerPolynomials2003}.
}
\[
    \uPnorm_n(t) := \begin{cases}
        \uP_n(t)/\uP_n(1) & \lambda < \infty
        \\
        t^n & \lambda = \infty.
    \end{cases}
\]
The limiting case follows from Lemma 1 by \citet{schoenbergPositiveDefiniteFunctions1942},
which shows for any $t\in [-1,1]$ that \(\uPnorm_n(t)\) converges uniformly in \(n\) to \(t^n\) as
\(\lambda\to\infty\). In our setting, the index \(\lambda=\frac{\dims-2}2\)
is determined by the dimension \(\dims\) and the high dimensional case can
therefore be well approximated by the case \(\dims=\infty\) which corresponds to
the sequence space
\(
    \real^\infty := \ell^2 = \{(x_i)_{i\in \nat}: \|x\|^2<\infty\}.
\)
The limiting case \(\uPnorm[\infty]_n(t) = t^n\) is not only an approximation of
the finite dimensional cases, it also corresponds to kernels that are valid in
all dimensions (Lemma \ref{lem: valid in all dimensions}).
\begin{theorem}[Characterization of positive definite isotropic kernels]
    \label{thm: isotropy characterization}
    Let \(d \in \{2,3,\dots, \infty\}\) and \(\lambda = \tfrac{\dims-2}2\).
    For a continuous kernel \(\kernel:
    \real^\dims \times \real^\dims \to \complex\), the following assertions are
    equivalent:
    \begin{enumerate}[label={(\roman*)}]
        \item\label{it: positive, definite isotropic}
        The kernel is positive definite and isotropic, i.e.
        \[
            \kernel(x,y) = \kernel(Ux, Uy)
        \]
        for every linear isometry \(U\).
        \item\label{it: pos definite, weak representation}
        The kernel is positive definite and has representation
        \[
            \kernel(x,y) = \kappa\bigl(\|x\|, \|y\|, \langle x, y\rangle\bigr)
        \]
        for a continuous \(\kappa: M \to \complex\) with \(M=\{ (r,s, \gamma) \in [0,\infty)^2 \times \real : \gamma \le rs\}\).
        \item\label{it: representation}
        With the convention \(\uPnorm_n\bigl(\bigl\langle \tfrac{x}{\|x\|}, \tfrac{y}{\|y\|}\bigr\rangle\bigr) := 1\) if $x=0$ or $y=0$, the kernel has
        representation
        \[
            \kernel(x,y)
            = \sum_{n=0}^\infty \alpha^{(\dims)}_n \bigl(\|x\|, \|y\|\bigr)
            \uPnorm_n\bigl(\bigl\langle \tfrac{x}{\|x\|}, \tfrac{y}{\|y\|}\bigr\rangle\bigr),
        \]
        where \(\alpha^{(\dims)}_n:[0,\infty)^2\to\complex\) are continuous,
        positive definite kernels. These kernels are summable, i.e.
        \(\sum_{l=0}^\infty \alpha^{(\dims)}_n(r,r) < \infty\) for any \(r\in
        [0,\infty)\), and zero at the origin, i.e. 
        \begin{equation}
            \label{eq: zero in origin}
            \alpha^{(\dims)}_n(r,0) = \alpha^{(\dims)}_n(0,r) = 0
            \qquad \forall r \in [0,\infty),\; n>0.
        \end{equation}
    \end{enumerate}
    Moreover for \(\dims <\infty\) the kernels \(\alpha^{(\dims)}_n\) can be
    represented as
    \begin{align}
    	\label{eq:alpha_formula}
        \alpha^{(\dims)}_n(r,s) = Z \int_{-1}^1 \kappa(r,s, rs\rho) \uPnorm_n(\rho)(1-\rho^2)^{\frac{\dims-3}2}d\rho,
    \end{align}
    where the normalization constant \(Z\) is defined by 
    \[
        Z = Z(n, \dims) = \Bigl[\int_{-1}^1 [\uPnorm_n(\rho)]^2(1-\rho^2)^{\frac{\dims-3}2}d\rho\Bigr]^{-1} .
    \]
\end{theorem}

\begin{remark}[Embedding]
    \label{rem: embedding}
    Embedding the finite dimensional spaces \(\real^\dims\) into the space of
    sequences \(\ell^2=\{(x_i)_{i\in \nat}: \|x\|^2<\infty\}\), we naturally obtain
    \[
        \real^1 \subseteq \real^{2} \subseteq \dots \subseteq \real^\infty := \ell^2~.
    \]
    Since any positive definite kernel remains positive definite on a subspace it
    follows that
    \begin{equation}
        \label{eq: pd subset relation}    
        \PDset_{O(1)}(\real) \supseteq \PDset_{O(2)}(\real^2) \supseteq \dots \supseteq \PDset_{\orthGroup[\infty]}(\real^\infty).
    \end{equation}
    
    While the representation \(\kernel(x,y) = \kappa(\|x\|, \|y\|, \scal{x}{y})\)
    can be used to lift any isotropic kernel into \(\ell^2\), the lifted kernel is not necessarily
    positive definite.
\end{remark}

In light of the remark above, one may be interested in the set of isotropic kernels that are positive definite
in all finite dimensional spaces \(\real^\dims\), that is on \(\bigcup_{\dims\in
\nat}\real^\dims\), which corresponds to the space of eventually zero sequences.
But these kernels coincide with the positive definite kernels on the entire
sequence space \(\ell^2\) as shown in the following lemma.

\begin{lemma}[Valid in all dimensions]
    \label{lem: valid in all dimensions}
    The set of isotropic positive definite kernels that are valid in all dimensions coincides
    with the set of isotropic positive definite kernels on \(\ell^2\), \ie
    \[
        \textstyle
        \bigcap_{\dims \in \nat} \PDset_{O(d)} = \PDset_{\orthGroup[\infty]} (\ell^2).
    \]
\end{lemma}
\begin{proof}
    `\(\supseteq\)' follows from \eqref{eq: pd subset relation}. For
    `\(\subseteq\)' we need to prove that the representation \(\kappa\) induces a positive definite
    kernel on \(\ell^2\). For this consider vectors \(x_1,\dots,x_n \in \ell_2\)
    and map their \(n\) dimensional subspace via a linear isometry \(U\in
    \orthGroup[\infty]\) into \(\real^n\). By positive definiteness in
    \(\PDset_{O(n)}(\real^n)\) Equation \eqref{eq:positive_definite_definition}
    follows.
\end{proof}

In particular the isotropic positive definite kernels on \(\ell^2\) coincide with the
isotropic positive definite kernels on the space of eventually zero sequences.

The core idea of the proof of Theorem~\ref{thm: isotropy characterization} is to translate results about
positive definite kernels on \((0,\infty) \times \sphere\), independently
discoverd by \citet{guellaSchoenbergsTheoremPositive2019} and
\citet{estradeCovarianceFunctionsSpheres2019}, to a characterization of
isotropic kernels on \(\real^\dims \setminus \{0\}\), and continuously extend
these results to \(\real^\dims\). 
For your convenience we restate Theorem 2.3 by
\citet{guellaSchoenbergsTheoremPositive2019} as Theorem~\ref{thm: guella XxS
characterization} in our notation.

\begin{lemma}[Continuous extension]
	\label{lem: continuous extension}
	Let \(\kernel: \domain \times \domain \to \complex\) be a continuous kernel
	on some topological space \(\domain\) and let $A\subset \domain$ be a
	set such that the closure of \(\domain \setminus A\) is $\domain$.
	If $\kernel$ is positive definite on $\domain\setminus A$, then \(\kernel\) is positive definite on \(\domain\).
\end{lemma}
\begin{proof}
	Let \(x_0,\dots, x_m\in \domain\) and \(c_0,\dots, c_m \in \complex\)
	and choose for $i=1,\dots, m$ a sequence $x_i^{(n)} \in \domain\setminus A$ such that $x_i^{(n)}\overset{n\to \infty}{\longrightarrow} x_i$.
	By continuity and positive definiteness of $\kernel$ we have 
	\[
	\sum_{i,j=0}^m c_i \conj{c_j} \kernel(x_i, x_j)= \lim_{n\to \infty} \sum_{i,j=0}^m  c_i\conj{c_j} \kernel(x_i^{(n)}, x_j^{(n)})\ge 0.
	\qedhere
	\]
\end{proof}

\begin{proof}[Proof of Theorem \ref{thm: isotropy characterization}]
    \textbf{\ref{it: positive, definite isotropic} \(\Leftrightarrow\) \ref{it: pos definite, weak representation}:} 
    ``\(\Leftarrow\)'': This follows
    directly from the fact that \(U\) is a linear isometry that preserves norms
    and inner products. \\
    ``\(\Rightarrow\)'':
    Let $M_\ge := \{(r,s,\gamma) \in M : r \ge s\}$ .  Using two orthonormal
    vectors \(e_1, e_2\), which exist due to \(\dims\ge 2\), we define the 
    \emph{continuous} maps
    $\varphi,\psi: M_\ge \to \real^\dims$ by
    \[
    	\varphi(r,s,\gamma) := r e_1,\qquad
    	\psi(r,s,\gamma) := \begin{cases}
    		\frac{\gamma}{r} e_1 + e_2 \sqrt{s^2 - \frac{\gamma^2}{r^2}} & r>0,\\
    		0 & r=0.
    	\end{cases}
    \]
    By the construction of $M_{\ge}$ all \((r,s,\gamma)\in M_\ge\) fulfill $r\ge s$ and $rs\ge \vert\gamma\vert$ and hence $\psi$ is continuous away from $(0,0,0)$. For the continuity of $\psi$ in \((0,0,0)\) 
    consider
    \begin{align*}
    	\|\psi(r,s,\gamma)\| \le \tfrac{rs}{r} + s = 2s \to 0
        \quad \text{for}\quad (r,s,\gamma) \to (0,0,0).
    \end{align*}
    The goal of $\varphi$ and $\psi$ is to create two vectors of
    specified length and angle. Indeed, for \(\theta=(r,s,\gamma)\) we have
    \begin{equation}
    	\label{eq:preserved_directions}
    	\|\varphi(\theta)\| = r, \quad \|\psi(\theta)\| = s, \quad \langle \varphi(\theta), \psi(\theta)\rangle = \gamma.
    \end{equation}
    The continuity requirement prevented the definition of \(\varphi\) and \(\psi\)
    on the entirety of $M$. To deal with $r<s$ we define the swap function
    \(\tau(r,s,\gamma) := \bigl(s,r,\gamma\bigr)\).
    With its help we finally define $\kappa\colon M\to \complex$ for $\theta = (r,s,\gamma)$ as
    \begin{equation}
        \label{eq: definition of kappa}
        \kappa(r,s,\gamma) := 
        \begin{cases}
        	\kernel\bigl(\varphi(\theta), \psi(\theta)\bigr) &  \text{ if } r\ge s \iff \theta \in M_\ge  \\
        	\kernel\bigl(\psi(\tau(\theta)), \varphi(\tau(\theta)) \bigr)	& \text{ if } r\le s.
        \end{cases}
    \end{equation}
    Clearly, \(\kappa\) is continuous if it is well-defined in \(r=s\) since \(\varphi\), \(\psi\) and \(\tau\) are
    continuous. For the case \(r=s\),
    Proposition \ref{prop: characterization of isometries} and Equation
    \eqref{eq:preserved_directions} yield a linear isometry $U\in
    \orthGroup$ such that $U\varphi (\theta) = \psi(\tau(\theta))$
    and $U \psi(\theta) = \varphi(\tau(\theta))$ and thus
    \[
    	\kernel\bigl(\varphi(\theta) , \psi(\theta)\bigr)
        \overset{\text{isotropy}}= \kernel\bigl( U \varphi(\theta), U\psi(\theta) \bigr)
        = \kernel\bigl(\psi(\tau(\theta)), \varphi(\tau(\theta))\bigr).
    \]
	It remains to show that $\kappa(\| x\|,\|y\|,\langle x,y\rangle)= \kernel
	(x,y)$ holds for all $x,y\in \real^\dims$. Since the other case works analogous
    we assume w.l.o.g.\@ $\|x\|\ge \|y\|$ and denote by $\theta :=(\|x\|,\|y\|,\langle x,y
	\rangle)$ the input to \(\kappa\). Then by Prop.~\ref{prop: characterization of isometries}
    and \eqref{eq:preserved_directions} there again exists
    $U\in  \orthGroup$  such that \(U x = \varphi(\theta)\) and \(Uy = \psi(\theta)\)
    and we conclude
    \begin{align*}
    	\kernel(x,y)
        = \kernel(Ux, Uy)
        =  \kernel(\varphi(\theta), \psi(\theta))
        = \kappa(\theta).
    \end{align*}

    \textbf{\ref{it: positive, definite isotropic} \(\Rightarrow\) \ref{it: representation}:}
    Over the homeomorphism \(T\colon (0,\infty) \times \sphere\to
    \real^d\setminus\{0\}\) with \(T(r,v) = rv\)
    we define the continuous kernel $\kernel_{\sphere[]}\colon
    ((0,\infty) \times \sphere)^2 \to \complex $ by
    \[
    	\kernel_{\sphere[]} ((r,v),(s,w))
        := \kernel(rv,sw)
        \overset{\ref{it: pos definite, weak representation}}=
        \kappa\bigl(r, s, rs \langle v, w\rangle \bigr).
    \]
    With \(f_{\sphere[]}(r,s,\rho):= \kappa(r,s,rs\rho)\)
    it is clear that Theorem~\ref{thm: guella XxS characterization} is
    applicable and we obtain by its representation of \(f_{\sphere[]}\) for any
    \(x,y\in \real^\dims\setminus\{0\}\) 
    \[
        \kernel(x,y)
        = f_{\sphere[]}\Bigl(\|x\|,\|y\|,\bigl\langle \tfrac{x}{\|x\|},\tfrac{y}{\|y\|}\bigr\rangle\Bigr)
    	= \sum_{n=0}^\infty \alpha_n^{(\dims)} (\|x\|,\|y\|) \uPnorm_n\Bigl(\bigl\langle \tfrac{x}{\|x\|},\tfrac{y}{\|y\|}\bigr\rangle\Bigr),
    \]
    where $\alpha_n^{(\dims)}$ are positive definite kernels on
    $(0,\infty)$ that are summable on the diagonal. What remains to be shown is that
    these kernels are continuous and can be extended to $[0,\infty)$, that this
    extension indeed resembles $\kernel$ and that \eqref{eq: zero in
    origin} and \eqref{eq:alpha_formula} hold. To this end we investigate the cases
    $\dims<\infty$ and $\dims=\infty$ separately.
    
    \begin{enumerate}[wide]
    	\item \textbf{Case \(\dims<\infty\)}:
    	For all $n\in\nat$ and all $r,s\in (0,\infty)$ we have by Theorem~\ref{thm: guella XxS characterization} 
\[
            \alpha^{(\dims)}_n(r,s)
            = Z \int_{-1}^1
            \overbrace{\kappa(r,s,rs\rho)}^{=f_{\sphere[]}(r,s, \rho)} \uPnorm_n(\rho)(1-\rho^2)^{\frac{\dims-3}2}d\rho.
    	\]
        If we extend this representation to \([0,\infty)\), we obtain
        the representation \eqref{eq:alpha_formula} by definition.
        Since \(\kappa\) is continuous by \ref{it: pos definite, weak
        representation}, dominated convergence using \(|\uPnorm_n(\rho)(1-\rho^2)^{\frac{\dims-3}2}|\le 1\)
        \autocite[e.g.][Eq. (2.16)]{reimerGegenbauerPolynomials2003} 
            implies we can extend \(\alpha^{(\dims)}_n\) to a continuous function
            on $[0,\infty)^2$.
This continuous extension remains positive definite by Lemma~\ref{lem:
        continuous extension}. To see that this extension results in a representation
        of $\kernel$, observe that, by the orthogonality of
        the Gegenbauer polynomials with respect to the weight function $w(\rho)
        = (1-\rho^2)^{\lambda - \frac12}$
        \autocite[Theorem~2.3]{reimerGegenbauerPolynomials2003}, $\uPnorm_0 \equiv 1$
        and the definition of \(Z\) as a normalization factor, we have for $r\ge
        0$
        \[
            \alpha^{(\dims)}_n(r,0)
            = Z \int_{-1}^1
            \kappa(r,0,0) \uPnorm_0(\rho)\uPnorm_n(\rho) (1-\rho^2)^{\frac{\dims-3}2}d\rho
            = \kappa (r,0,0) \delta_{n0},
        \]
        where $\delta_{n0}$ denotes the Kronecker delta. This implies \eqref{eq: zero in origin}
        and we also obtain
        \[
            \kernel(x,0)
            = \kappa(\|x\|,0, 0 ) = \sum_{n=0}^\infty\underbrace{\alpha^{(\dims)}_n(\|x\|, \|0\|)}_{=\kappa (\|x\|,0,0) \delta_{n0}}
            \underbrace{\uPnorm_n\bigl(\langle \tfrac{x}{\|x\|}, \tfrac{0}{\|0\|}\rangle \bigr)}_{\overset{\text{def.}}=1}.
        \]
        By an analogous argument the representation also holds for \(\kernel(0,x)\).

    	\item \textbf{Case \(\dims=\infty\):}
    	First note that the \(\alpha^{(\dims)}_n\) are \emph{continuous} kernels on
    	\((0,\infty)\) by
    	\citet[Example~2.4]{guellaSchoenbergsTheoremPositive2019}.  For the
    	extension we select two orthonormal vectors \(e_1,e_2\) and observe
    	\begin{align*}
    		\kernel(r e_1, r e_1) - \kernel(r e_1, r e_2)
    		&= \sum_{n=0}^\infty \alpha^{(\dims)}_n(r,r)\langle e_1, e_1\rangle^n
    		- \sum_{n=0}^\infty \alpha^{(\dims)}_n(r,r)\langle e_1, e_2\rangle^n
    		\\
    		&= \sum_{n=\red{1}}^\infty \alpha^{(\dims)}_n(r,r).
    	\end{align*}
    	By continuity of \(\kernel\) we obtain
        \[
    		0 = \lim_{r\to 0}\kernel(r e_1, r e_1) - \kernel(r e_1, r e_2)
    		= \lim_{r\to 0}\sum_{n=1}^\infty \alpha_n^{(d)}(r,r).
        \]
    	Since \(\alpha^{(\dims)}_n(r,r)\ge 0\) holds for all $n\in\nat$ since the $\alpha^{(\dims)}_n$ are kernels, this implies
    	\(\alpha^{(\dims)}_n(r,r) \to 0\) for all \(n\ge 1\).
    	Assume \(r=0\) or \(s=0\) and 
    	let \((r_l, s_l)\overset{n\to \infty}{\longrightarrow}(r,s)\). Then the Cauchy-Schwarz inequality yields
    	\[
    	|\alpha^{(\dims)}_n(r_l, s_l) - \underbrace{\alpha^{(\dims)}_n(r,s)}_{=0}|
    	\overset{\text{C.S.}}\le \sqrt{\alpha^{(\dims)}(r_l, r_l) \alpha^{(\dims)}_n(s_l, s_l)} \overset{l\to\infty} \longrightarrow 0.
    	\]
    	Hence, we can continuously extend $\alpha^{(d)}_n$ by zero as claimed in \eqref{eq: zero in origin}.
    	Lemma~\ref{lem: continuous extension} implies for \(n>0\) that this continuous extension of positive
    	definite kernels \(\alpha^{(\dims)}_n\)  is positive definite. 
    	We continue with the extension in the case \(n=0\). For this consider
    	\[
    	\kernel(r e_1, s e_2) = \sum_{n=0}^\infty \alpha^{(\dims)}_n(r,s) \langle e_1, e_2\rangle^n = \alpha^{(\dims)}_0(r,s).
    	\]
    	This yields the continuous extension $\alpha^{(\dims)}_0(r,s)
    	:= \kernel(re_1,se_2)$. Lemma~\ref{lem: continuous extension} shows
    	that this extension remains positive definite.  What is left to prove is
    	that our continuous extension of the
    	\(\alpha^{(\dims)}_n\) kernels does in fact lead to a representation 
    	of \(\kernel\). To this end, let $x,y\in \real^\dims$ and w.l.o.g.\@ assume $x=0$ to
        obtain
        \[
    		K(x,y) = K(0,\|y \| e_2 ) = \alpha_0^{(\dims)}(0,\|y\| ) = \sum_{n=0}^\infty \alpha_n^{(\dims)}(0,\|y\|) \langle \tfrac{0}{\|0\|}, y\rangle^n.
        \]
    \end{enumerate}
    
    \textbf{\ref{it: representation} \(\Rightarrow\) \ref{it: positive, definite isotropic}:} 
    Isotropy is obvious and it remains to show
    positive definiteness. We split the series representation into two parts
    \begin{equation}
        \label{eq: kernel decomposition} 
        \kernel(x,y)
        = \underbrace{\alpha^{(\dims)}_0(\|x\|, \|y\|)}_{\kernel_0(x,y)}
        + \underbrace{\sum_{n=1}^\infty \alpha^{(\dims)}_n(\|x\|, \|y\|)\uPnorm_n\bigl(\langle \tfrac{x}{\|x\|}, \tfrac{y}{\|y\|}\rangle\bigr)}_{=\kernel_{>0}(x,y)},
    \end{equation}
    where we use \(\uPnorm_0 \equiv 1\).  
    Let $m\in \nat , x_1,\dots,x_m\in \real^\dims, c_1,\dots,c_m\in \complex$
    and let without loss of generality $x_i=0$ if and only if $i=1$. Since
    $\alpha_n^{(\dims)}(0,r)=0$ for $n\ge 1$ we obtain
    \[
    	\sum_{i,j=1}^m c_i \conj{c_j} \kernel(x_i,x_j)
        = \sum_{i,j=1}^m c_i\conj{c_j} \kernel_0(x_i,x_j) + \sum_{i,j=2}^m c_i\conj{c_j} \kernel_{>0}(x_i,x_j).
    \]
    Since $\alpha_0^{(d)}$ is a kernel on $[0,\infty)$ we get that the first
    summand is non-negative. The second summand is non-negative, since it only
    contains $x_i\in \real^\dims \setminus\{0\}$ and the restriction
    $\kernel_{>0}\vert_{(\real^\dims\setminus \{0\})^2}$ is positive definite
    using the homeomorphism $T\colon (0,\infty) \times \sphere\to \real^\dims
    \setminus \{0\},\; (r,v)\mapsto rv$ and Theorem~\ref{thm: guella XxS
    characterization} characterizing kernels on \((0,\infty)\times \sphere\).
\end{proof}

 \subsection{Strict positive definiteness}
\label{sec: strict positive definite}

In this section we characterize the requirements on a continuous, isotropic
positive definite kernel to be \emph{strictly} positive definite.\footnote{
    For readers interested in \emph{universal} kernels we want to point out
    their characterization on the sphere by \citet[Theorem
    10]{micchelliUniversalKernels2006} that may also be helpful via the
    connection to \((0,\infty)\times \sphere\) for the characterization of
    universal kernels in \(\PDset_{\orthGroup}\).
}
We will again utilize the connection between \((0,\infty)\times \sphere\) and
\(\real^\dims\setminus\{0\}\) and base our results on these by
\citet{guellaSchoenbergsTheoremPositive2019}. Since they had to single out
the case \(\dims=2\) case, we have to do so too.\footnote{
    Warning: We consider \(\sphere\) whereas they considered \(\sphere[\dims]\) and
    therefore singled out \(\dims=1\).
}

\begin{theorem}[Characterization of strictly positive definite isotropic kernels]
    Let \(\dims \in \{3,\dots, \infty\}\).  For a continuous positive definite isotropic
    kernel \(\kernel\colon \real^\dims
    \times \real^\dims \to \complex\)  with representation as in Theorem~\ref{thm: isotropy
    characterization}, the following assertions are equivalent:
    \begin{enumerate}[label={(\roman*)}]
        \item\label{it: strict positive definite}
        The kernel \(\kernel\) is strictly positive definite;
        \item\label{it: infinitely even and odd elements in set}
        \(\alpha_0^{(\dims)}(0,0)>0\), and
        for any \(m\in \nat\), distinct \(r_1,\dots, r_m\in (0,\infty)\) and \(c \in \complex^m\setminus\{0\}\)
        the set
        \[
            \Bigl\{
                n \in \nat_0 :
                c^\transpose \Bigl[\alpha^{(\dims)}_n(r_i, r_j)\Bigr]_{i,j=1}^m \conj{c} > 0
            \Bigr\}
        \]
        contains infinitely many even and infinitely many odd integers.
        \item\label{it: even and odd kernels strict pos. def.}
        \(\alpha_0^{(\dims)}(0,0)>0\), and for each \(\gamma \ge 0\) the even and odd kernels
        \[
            \alpha_\gamma^e(r,s) := \sum_{2n \ge \gamma} \alpha_{2n}^{(\dims)}(r,s)
            \quad\text{and}\quad
            \alpha_\gamma^o(r,s) := \sum_{2n+1 \ge \gamma} \alpha_{2n+1}^{(\dims)}(r,s)
        \]
        are both strictly positive definite on \((0,\infty)\).
    \end{enumerate}
    Furthermore, the following statement is sufficient for strict positive definiteness:
    \begin{enumerate}[label={(\roman*)},resume]
        \item\label{it: infinitely many even and odd n strict pos. def.}
        \(\alpha_0^{(\dims)}(0,0)>0\), and \(\alpha_n^{(\dims)}\) is strictly positive
        definite on \((0,\infty)\) for infinitely many even and odd \(n\in \nat_0\).
    \end{enumerate}
    If \(\alpha_0^{(\dims)}(0,0) =0\) but all other requirements in \ref{it:
    infinitely even and odd elements in set}, \ref{it: even and odd kernels
    strict pos. def.}, \ref{it: infinitely many even and odd n strict pos. def.}
    are satisfied, then \(\kernel\) is strictly positive definite on
    \(\real^\dims\setminus\{0\}\). \ref{it: infinitely even and odd elements in
    set} and \ref{it: even and odd kernels strict pos. def.} without \(\alpha_0^{(\dims)}(0,0) > 0\)
    are also necessary for strict positive definiteness on \(\real^\dims\setminus\{0\}\).
\end{theorem}
\begin{proof}
Recall the decomposition \eqref{eq: kernel decomposition}, i.e.
\[
    \kernel(x,y)
    = \underbrace{\alpha^{(\dims)}_0(\|x\|, \|y\|)}_{=:\kernel_0(x,y)}
    + \underbrace{\sum_{n=1}^\infty \alpha^{(\dims)}_n(\|x\|, \|y\|) \uPnorm_n(\langle \tfrac{x}{\|x\|}, \tfrac{y}{\|y\|}\rangle)}_{=:\kernel_{>0}(x,y)},
\]
with positive definite kernels $\kernel_0$ and $\kernel_{>0}$.  Let \(m\in
\nat\), \(c\in (\complex\setminus\{0\})^m\) and let  \(x_1,\dots, x_m\in
\real^\dims\) be distinct locations where \(x_1=0\) w.l.o.g.\@, then we have
    \begin{equation}
        \label{eq: kernel decomposition 2}
    	\sum_{i,j=1}^m c_i \conj{c_j} \kernel(x_i,x_j)
        = \sum_{i,j=1}^m c_i\conj{c_j} \kernel_0(x_i,x_j)
        + \sum_{i,j=2}^m c_i\conj{c_j} \kernel_{> 0}(x_i,x_j)
\end{equation}
    using that \(\kernel_{>0}(x,0) = \kernel_{>0}(0,x) = 0\) for all \(x\in \real^\dims\)
    due to \eqref{eq: zero in origin}
    in Theorem~\ref{thm: isotropy characterization}.

    ``\ref{it: infinitely even and odd elements in set}, \ref{it: even and odd kernels strict pos. def.}, \ref{it: infinitely many even and odd n strict pos. def.} \(\Rightarrow\) \ref{it: strict positive definite}'':
    For strict positive definiteness of \(\kernel\) we need to show \eqref{eq: kernel decomposition 2}
    is non-zero. Since both summands are non-negative it is sufficient if either one is non-zero.
    With the usual homeomorphism the conditions in \ref{it: infinitely even and odd elements in set},
    \ref{it: even and odd kernels strict pos. def.} and \ref{it: infinitely many
    even and odd n strict pos. def.} 
    imply the conditions in Theorem 3.5, 3.7 and 3.9 of
    \citet{guellaSchoenbergsTheoremPositive2019} on \((0,\infty)\times
    \sphere\). Since these conditions are satisfied for \(\kernel_{>0}\) if
    they are satisfied for \(\kernel\) we thereby see that \(\kernel_{>0}\) is strictly positive definite on
    \(\real^\dims \setminus\{0\}\).
    If \(m>1\) the second term in \eqref{eq: kernel decomposition 2} is therefore
    non-zero.  
    If $m=1$ then the first term is
    \[
        \sum_{i,j=1}^m c_i\conj{c_j} \kernel_0(x_i,x_j)
        = |c_1|^2 \alpha_0^{(\dims)}(0,0) >0.
    \]

    ``\ref{it: strict positive definite} \(\Rightarrow\) \ref{it: infinitely even and odd elements in set}, \ref{it: even and odd kernels strict pos. def.}'':
    If \(\kernel\) is strictly positive definite on \(\real^\dims\), then
    \(\alpha^{(\dims)}_0(0,0)>0\) by the selection of \(m=1\), \(c=1\),
    \(x_1=0\) in \eqref{eq: kernel decomposition}. 
    Mapping strict positive definiteness of \(\kernel\) on \(\real^\dims\setminus\{0\}\)
    into strict positive definiteness on \((0,\infty)\times \sphere\), the
    assertions \ref{it: infinitely even and odd elements in set} and \ref{it:
    even and odd kernels strict pos. def.} follow from Theorem 3.5 and 3.7 of
    \citet{guellaSchoenbergsTheoremPositive2019}.

    The final remark about strict positive definiteness on \(\real^\dims\setminus\{0\}\)
    follows similarly from the homeomorphism to \((0,\infty)\times \sphere\) and
    the results of \citet{guellaSchoenbergsTheoremPositive2019}.
\end{proof}
By the same arguments Theorem 3.6, 3.8 and 3.9 of
\citet{guellaSchoenbergsTheoremPositive2019} for the case \(\dims=2\) yield to the
following theorem.
\begin{theorem}[Strict positive definiteness for \(\dims=2\)] 
    Let \(\kernel\colon\real^2 \times\real^2 \to \complex\) be a continuous positive
    definite isotropic kernel with decomposition as in Theorem~\ref{thm:
    isotropy characterization}. Then \emph{necessary} conditions for
    strict positive definiteness include
    \begin{enumerate}[label={(\roman*)}]
        \item \(\alpha_0^{(2)}(0,0)>0\), and
        for any \(m\in \nat\), \(c \in \complex^m\setminus\{0\}\) and \(r_1,\dots, r_m\in (0,\infty)\)
        the set
        \[
            \Bigl\{
                n \in \integer :
                c^\transpose \Bigl[\alpha^{(2)}_{|n|}(r_i, r_j)\Bigr]_{i,j=1}^m \conj{c} > 0
            \Bigr\}
        \]
        intersects every full arithmetic progression in \(\integer\).

        \item \(\alpha_0^{(2)}(0,0)>0\), and for each full arithmetic progression \(S\) in \(\integer\)
        the kernel
        \[
            (r,s) \mapsto \sum_{n\in S} \alpha^{(2)}_{|n|}(r,s) \qquad r,s\in (0,\infty).
        \]
        is strictly positive definite.
    \end{enumerate}
    and a \textbf{sufficient} condition for strict positive definiteness of \(\kernel\)
    is
    \begin{enumerate}[label={(\roman*)},resume]
        \item \(\alpha_0^{(2)}(0,0)>0\), and the set \(\{n \in \integer : \alpha_{|n|}^{(2)} \text{ is strictly positive definite}\}\)
        intersects every full arithmetic progression in \(\integer\).
    \end{enumerate}
\end{theorem} \subsection{Applications and examples}
\label{sec: applications}

In applications such as support vector machines and
Bayesian modelling one has to choose a family of covariance kernels
and fit its parameters to the data. The family of non-stationary isotropic
covariance kernels represents a joint generalization of stationary
isotropic covariance kernels and `dot product' kernels, which have so far been
treated as separate cases
\autocite[e.g.][]{roosHighDimensionalGaussianProcess2021,amentScalableFirstOrderBayesian2022}.

Popular stationary isotropic covariance kernels include the squared exponential\footnote{
	cf.\ Gaussian smoothing in Example \ref{ex: convolution with mollifier}
	and Schoenbergs characterization of stationary isotropic
	kernels in \(\ell^2\) in Table \ref{table: characterizations}.
}
\[
\kernel(x,y) = \sigma^2\exp\bigl(-\tfrac{\|x-y\|^2}{2s^2}\bigr)
= \sum_{n=0}^\infty \underbrace{
	\frac{\sigma^2\|x\|^n \|y\|^n}{n!\,s^{2n}} e^{-\tfrac{\|x\|^2}{2s^2}} e^{-\tfrac{\|y\|^2}{2s^2}} 
}_{=\alpha_n^{(\infty)}(\|x\|, \|y\|)}\scal{\tfrac{x}{\|x\|}}{\tfrac{y}{\|y\|}}^n,
\]
as well as the Matérn and rational quadratic kernels. Dot product kernels
include the linear kernel
\[
\kernel(x,y) = \sigma_b^2 + \sigma_w^2\langle x,y\rangle
\]
used for Bayesian linear regression as well as the polynomial kernels which are
powers of the linear kernel. The linear kernel is also the covariance function
of the mean squared error applied to a random linear model (Section \ref{sec:
random linear model}). The union of stationary isotropic kernels and dot product
kernels captures essentially all kernel families
used in practice \autocite[e.g.][Sec. 4.2]{rasmussenGaussianProcessesMachine2006} --
up to `geometric anisotropies' \autocite[p.\ 50]{steinInterpolationSpatialData1999}.

Strictly non-stationary isotropic covariance kernels have also naturally come
up as ``Neural network kernels'', i.e. kernels that arise from taking
the layer widths of neural networks to infinity and initializing them randomly.
This approach was pioneered by \citet{williamsComputingInfiniteNetworks1996}
where the limit of a one layer neural network with the `error function' as
activation function resulted in the kernel
\[
\kernel(x,y) = \arcsin\Bigl(\frac{x\cdot y}{\sqrt{(1+\|x\|^2)(1+\|y\|^2)}}\Bigr)
\]
Another such example is the arc-cosine kernel \autocite{choKernelMethodsDeep2009} of
the form 
\[
\kernel(x,y) = \frac1\pi\|x\|^n \|y\|^n J_n(\cos^{-1}(x\cdot y))
\]
for certain \(J_n\) that arise when the activation function
\(\phi(x) = \max\{0,x\}^n\) is used instead. This includes the popular ReLU
activation function with \(n=1\).

\begin{example}
	\label{ex: multilayer neural networks are isotropic at initialization}
	Multilayer neural networks with activation function \(\phi\) are
	(non-stationary) isotropic, centered Gaussian random functions at initialization.
\end{example}
\begin{proof}
	\citet{leeDeepNeuralNetworks2018} show
	that multilayer neural networks are centered Gaussian processes at initialization and that
	the kernel of the \(l+1\)-th layer is recursively given by
	\[
	\kernel_{l+1}(x,y)
	= \sigma_b^2 + \sigma_w^2F_\phi(\kernel_l(x,y), \kernel_l(x,x), \kernel_l(y,y))
	\]
	for some function \(F_\phi\) that depends on the activation function \(\phi\)
	with \(\kernel_0(x,y) = \sigma_b^2 + \sigma_w^2(x\cdot y)\).
	With this formula one can show recursively that the kernels are of the form
	\(\kernel_l(x,y) = \kappa_l(\|x\|, \|y\|, x\cdot y)\) and therefore
	non-stationary isotropic. In particular the output layer
	and thus the entire network is an isotropic Gaussian random function.
\end{proof}

 }
\section{Covariance of derivatives}
\label{sec: covariance of derivatives}

In Section \ref{sec: sample regularity of random functions} we explained that the
covariance of derivatives of \(\rf\) is given by the derivatives of covariances
\[
    \Cov(\partial_{x_i}\rf(x), \partial_{y_j} D_w\rf(y))
    = \partial_{x_i}\partial_{y_j} \C_\rf(x,y)
\]
In this section we want to calculate these covariances explicitly for
the case of stationary and non-stationary isotropic random functions.
To calculate the covariances of more general directional derivatives we
introduce some notation. For functions with one input, \(D_v f(x)\)
unmistakeably denotes the directional derivative
of \(f\) in the direction of \(v\). But for functions with multiple inputs
\(D_vf(x,y)\) is unclear. To better mark the place the directional
derivative is applied to, we introduce the notation
\[
    f(D_v x,y) := \tfrac{d}{dt} f(x + tv, y).
\]
such that
\begin{equation}
    \label{eq: covariance of directional derivatives}
    \Cov(D_v\rf(x), D_w\rf(y))
    = \C_\rf(D_v x, D_w y).
\end{equation}

First, we calculate the covariance of derivatives for stationary
isotropic random functions \(\rf\). Recall these have covariance functions of
the form \(\C_{\rf}(x,y) = \ikernel\bigl(\tfrac{\|x-y\|^2}2\bigr)\) by
Theorem \ref{thm: characterization of weak input invariances}.
\begin{lemma}[Covariance of derivatives under stationary isotropy]
    \label{lem: cov of derivatives, isotropy}
    Let \(\rf\sim\normal(\mu, \ikernel)\) be an isotropic random function and let \(\distance = x-y\), then
    \begin{equation*}
    \def\arraystretch{1.5}
    \begin{tabular}{c | c c}
        \(\Cov\) & \(\rf(y)\) & \(D_w\rf(y)\) \\	
        \hline
        \(\rf(x)\)
        &  \(\ikernel\bigl(\frac{\|\distance\|^2}{2}\bigr)\)
        & \(-\ikernel'\bigl(\frac{\|\distance\|^2}2\bigr)\langle \distance, w\rangle\)
        \\
        \(D_v\rf(x)\)
        &  \(\ikernel'\bigl(\frac{\|\distance\|^2}{2}\bigr)\langle \distance, v\rangle\)
        & \(
            - \Bigl[
                \ikernel''\bigl(\tfrac{\|\distance\|^2}2\bigr)\langle \distance, w\rangle\langle \distance, v\rangle
                + \ikernel'\bigl(\tfrac{\|\distance\|^2}2\bigr)\langle w, v\rangle
            \Bigr]
        \)
    \end{tabular}
    \end{equation*}
\end{lemma}
\begin{proof}
    Straightforward application of chain and product rules applied to \eqref{eq: covariance of directional derivatives}. 
\end{proof}

Second we will consider the non-stationary isotropic kernels. Recall that these
are of the form \(\C_\rf(x,y) = \kernel\bigl(\tfrac{\|x\|^2}2, \tfrac{\|y\|^2}2,
\langle x, y\rangle\bigr)\) by Theorem \ref{thm: characterization of weak input
invariances}. Since the (non-stationary) isotropic kernel \(\kernel\) have
multiple inputs, we further introduce the notation
\[
    \kernel_i(\lambda_1,\lambda_2,\lambda_3)
    := \tfrac{d}{d\lambda_i} \kernel(\lambda_1, \lambda_2, \lambda_3)
    \qquad i\in \{1,2,3\}.
\]

\begin{restatable}[Covariance of derivatives under (non-stationary) isotropy]{lemma}{covOfDerivatives}
    \label{lem: cov of derivatives, non-stationary isotropic case}
    Let \(\rf\sim\normal(\mu, \kernel)\) be (non-stationary) isotropic
    (Definition~\ref{def: distributional input invariance}). Then the
    expectation of the directional derivative is given by
    \begin{equation}
        \label{eq: expect df}
        \E[D_v\rf(x)] = \mu'\bigl(\tfrac{\|x\|^2}2\bigr) \langle x, v\rangle,
    \end{equation}
    and, using the abbreviation \(\kernel :=
    \kernel\bigl(\tfrac{\|x\|^2}2,\tfrac{\|y\|^2}2, \langle x,y\rangle \bigr)\)
    to omit inputs, the covariance of directional derivatives is given by
    \begin{align}
        \label{eq: cov df, f}
        \Cov(D_v \rf(x), \rf(y))   
        &= \frac1\dims\Bigl[\kernel_1 \langle x, v\rangle + \kernel_3 \langle y, v\rangle\Bigr]
        \\
        \label{eq: cov df, df}
        \Cov(D_v \rf(x), D_w\rf(y))   
        &= \frac1\dims\Bigl[\underbrace{\begin{aligned}[t]
            &\kernel_{12} \langle x, v\rangle\langle y, w\rangle + \kernel_{13} \langle x, v\rangle\langle x, w\rangle
            \\
            & + \kernel_{32} \langle y, v\rangle\langle y, w\rangle + \kernel_{33}\langle y, v\rangle\langle y, w\rangle
            \\
        \end{aligned}}_{\text{(I)}}
        + \underbrace{\kernel_3 \langle v, w\rangle}_{\text{(II)}}\Bigr].
    \end{align}
    In particular, if the directions \(v,w\) are orthogonal to the points
    \(x,y\) then (I) of \eqref{eq: cov df, f} and \eqref{eq: cov df, df} are zero.
    (II) in turn is zero, if the directions \(v,w\) are orthogonal.
\end{restatable}
\begin{proof}
    Straightforward application of chain and product rules applied to \eqref{eq:
    covariance of directional derivatives}.
\end{proof} \section{Strict positive definite derivatives}
\label{sec: strictly pos definite derivatives}

While covariance matrices are always positive definite, they are not always
strict positive definite (i.e. invertible).
Recall, that a random function \(\rf\) and its covariance function \(\C_{\rf}\)
are \emph{strict} positive definite if the matrix \((\C_{\rf}(x_k,
x_l))_{k,l=1,\dots, n}\) is strict positive definite for any finite
\(x_1,\dots, x_n\) (cf. Section \ref{sec: characterization of covariance
kernels}). It is already known that, whenever \(\rf\) is stationary isotropic
and non-constant, \(\rf\) is strict positive definite if the
dimension satisfies \(\dims \ge 2\)
\autocite[Theorem~3.1.5]{sasvariMultivariateCharacteristicCorrelation2013}.
But since we require the jet \(\jet[1]\rf = (\rf, \nabla\rf)\)
to be strict positive definite in Chapter \ref{chap: predictable
progres}, we prove a generalization of Theorem~3.1.6 in
\citet{sasvariMultivariateCharacteristicCorrelation2013}. This
generalization will be sufficient to prove that all stationary isotropic
covariance functions which are valid in all dimensions (cf. Lemma \ref{lem:
valid in all dimensions}) are covered as we will see in Corollary~\ref{cor:
strict positive definite}. This result will be used to remove
the strict positive definiteness assumption in Corollary~\ref{cor:
asymptotically deterministic behavior}.

If \(\rg\) is a random function with multivariate output, then \(\C_\rg(x_k, x_j)\) is already a
matrix for fixed \(k,j\) and the collection over \(k,j\) is really a tensor. To
avoid introducing this machinery and explaining positive definiteness for
tensors, we will take the following equivalent statement for positive
definiteness of random functions as definition.

\begin{definition}[Strict positive definite random function]
	\label{def: strict positive definite random function}
    The covariance \(\C_\rg\) of a random function \(\rg: \real^\dims\to\real^m\)
    and the random function itself is called \emph{strict positive definite}
    if for all \(w_k\in \real^m\) and distinct \(x_1,\dots, x_n\in \real^\dims\) the
    equality
    \[
        0 = \Var\Bigl[\sum_{k=1}^n w_k^T\rg(x_k)\Bigr]
        \overset{\eqref{eq: equivalent formulation}}= \sum_{k,l=1}^n w_k^T\C_\rg(x_k, x_l) w_l
    \]
    implies \(w_k = 0\) for all \(k\).
\end{definition}
Note, that the second equality marked with \eqref{eq: equivalent formulation} is
always true and only represents an equivalent formulation, because after
centering \(\rg\) (without loss of generality, using \(\tilde{\rg} = \rg -
\E[\rg]\)), we have
\begin{equation}
    \label{eq: equivalent formulation}
    \sum_{k,l=1}^n w_k^T \C_\rg(x_k,x_l) w_l
    = \E\Bigl[\sum_{k,l=1}^n w_k^T \rg(x_k) \rg(x_l)^T w_l\Bigr]
    = \Var\Bigl[ \sum_{k=1}^n w_k^T \rg(x_k)\Bigr].
\end{equation}

\begin{restatable}[Strict positive definite derivatives]{theorem}{strictlyPosDefiniteDerivatives}
    \label{thm: strict positive definite derivatives}
    Let \(\rf: \real^\dims \to \real\) be a stationary random function
    (cf. Definition~\ref{def: distributional input invariance} and Theorem~\ref{thm:
    characterization of weak input invariances}) with
    continuous covariance such that the support of its spectral measure contains
    a non-empty open set. Then up to any order \(k\) up to which \(\rf\) is
    almost surely differentiable, the jet
    \[
        \jet\rf
        = (
			\rf, \nabla\rf,
			(\partial_{ij}\rf)_{i\le j},
			\dots,
			(\partial_{i_1\cdots i_k} \rf)_{i_1\le \dots\le i_k}
		)
    \]
    is strict positive definite.
\end{restatable}
Theorem~\ref{thm: strict positive definite derivatives}
is not more general than
\citet[Theorem~3.1.6]{sasvariMultivariateCharacteristicCorrelation2013} in the
sense that conditions are weakend, but it is more general in its implication. That is, 
\(\rg\) is strict positive definite and not just \(\rf\).

The following corollary shows that this fully covers the stationary isotropic
covariance functions valid in all dimensions (cf.\ Lemma~\ref{lem: valid in all dimensions}). We use this result about stationary isotropic
random functions in Corollary~\ref{cor: asymptotically deterministic behavior}
to omit the assumption that \((\rf, \nabla\rf)\) has to be strict
positive definite, which is needed for our general result (Theorem~\ref{THM: ASYMPTOTICALLY DETERMINISTIC BEHAVIOR}).

\begin{restatable}{corollary}{corStrictPosDef}
    \label{cor: strict positive definite}
	Assume that \(\ikernel\) is a continuous stationary isotropic covariance
	kernel valid in all dimensions (\ie defined on \(\ell^2\) by Lemma \ref{lem:
	valid in all dimensions}).  Let \(\rf\sim\normal(\mu, \ikernel)\) be a
	Gaussian random function, which is not almost surely constant. Then the jet
    \[
        \jet \rf 
        = (\rf, \nabla\rf, (\partial_{ij}\rf)_{i\le j}, \dots, (\partial_{i_1\cdots i_k} \rf)_{i_1\le \dots\le i_k})
    \]
    is strict positive definite for any \(k\in \nat\).
\end{restatable}

\subsection{Proofs}

\begin{proof}[Proof of Theorem \ref{thm: strict positive definite derivatives}]
    For any \(n\) finite and distinct \(x_1,\dots,x_n\in\real^\dims\) we need that
    \[
        0=\Var \Bigl[\sum_{j=1}^n w_j^T \jet\rf(x_j)\Bigr]
    \]
    implies \(w_i=0\) for all \(i=1,\dots, n\). Without loss of generality we
    assume \(\rf\) (and thus \(\jet\rf\)) to be centered. We can then rewrite this as
    \[
        \Var \Bigl[\sum_{j=1}^n w_j^T \jet\rf(x_j)\Bigr]
        = \Var\Bigl[\sum_{j=1}^n \sum_{l=0}^k \sum_{i_1\le \dots \le i_l} w_j^{(l, i_1,\dots, i_l)}\partial_{i_1\dots i_l}\rf(x_j)\Bigr],
    \]
    with appropriate indexing of the \(w_j\). Note that \(l\) ranges over the
    order of differentiation contained in \(\rg\) before we sum over all the
    partial derivatives in the inner sum.
    For \(k=1\) this is for example
    \[
        \Var\Bigl[\sum_{j=1}^n w_j^{(0)}\rf(x_j) + \sum_{i=1}^\dims w_j^{(1,i)}\partial_i\rf(x_j)\Bigr]
        = \Var\Bigl[\sum_{j=1}^n w_j^{(0)}\rf(x_j) + D_{w_j^{(1,\cdot)}}\rf(x_j)\Bigr] = 0.
    \]
    We now consider the linear differential operator
    \[
        T := \sum_{j=1}^n \sum_{l=0}^k \sum_{i_1\le \dots \le i_l}
        (-1)^lw_j^{(l, i_1,\dots, i_l)}\partial_{i_1\dots i_l}\delta_{x_j},
    \]
    where \(\delta_x(f) = f(x)\) is the dirac delta function and \(\partial_i
    \delta_x(f) = - \partial_i f(x)\) is its derivative in the sense of distributions. 
    Recall this means for test functions \(\phi\), that we have
    \[
        \langle \partial_i\delta_x, \phi\rangle
        = -\langle\delta_x, \partial_i\phi\rangle
        = -\partial_i\phi(x).
    \]
    The higher order derivatives in the sense of distributions are similarly
    defined. Using this operator we now obtain more succinctly
    \[
        0 = \Var\Bigl[\sum_{j=1}^n \sum_{l=0}^k \sum_{i_1\le \dots \le i_l} w_j^{(l, i_1,\dots, i_l)}\partial_{i_1\dots i_l}\rf(x_j)\Bigr]
        = \Var[T\rf].
    \]
    While \(T\rf\) is well defined, we now want to move this operator outside. For this
    we want to use the bilinearity of the covariance. But the covariance has two
    inputs and \(T \C_{\rf}\) is then not well defined because the differential
    operator \(T\) might be applied to the first or second input of \(\C_{\rf}\).
    To avoid this issue, we define with some abuse of notation \(T^t f(t) := Tf\)
    such that we can write \(T^t\C_{\rf}(t, s)\) when we mean \(T \C_{\rf}(\cdot, s)\).
    Note that \(T\) is a linear combination of basis elements
    \(\partial_{i_1\dots i_l}\delta_x\). So by bilinearity of the covariance it
    is sufficient to check, that we can move these basis elements out of
    the covariance, i.e.
    \begin{align*}
        \Cov(\partial_{i_1\dots i_l}\delta_x\rf, \partial_{j_1\dots j_{l'}}\delta_y\rf)
        &= \Cov(\partial_{i_1\dots i_l}\rf(x), \partial_{j_1\dots j_{l'}}\rf(y))
        \\
        &= (\partial_{i_1\dots i_l}\delta_x)^t\; (\partial_{j_1\dots j_{l'}}\delta_y)^s \C_{\rf}(t,s).
    \end{align*}
    Since \(T\) is a linear combinations of these, we get by the bilinearity of
    the covariance
    \[
        0
        = \Var[T \rf]
        = T^x T^y \C_{\rf}(x,y).
    \]
    In the remainder of the proof we will essentially show, that this
    variance can be represented as an integral over the absolute value
    of the fourier transform of \(T\) with respect to the spectral measure
    of \(\C_{\rf}\). This forces the fourier transform of \(T\) and therefore
    \(T\) to be zero.

    Using the spectral representation of \(\C_{\rf}\) \autocite[e.g.][Theorem
    1.7.4]{sasvariMultivariateCharacteristicCorrelation2013} given by
    \[
        \C_{\rf}(x,y)
        = \int e^{\im \langle x-y, t\rangle} \spectMeasure(dt),
    \]
    we move the operator into the integral (by moving sums and derivatives
    into the integral)
    \[
        0 = T^x T^y \C_{\rf}(x,y)
        = \int T^x e^{\im\langle x,t\rangle} T^y e^{-\im\langle y,t\rangle}\spectMeasure(dt)
        = \int |T e^{\im\langle \cdot, t\rangle}| \spectMeasure(dt),
    \]
    where we use \(T\bar{f} = \overline{T f}\) for the last equation, which follows from
    \begin{itemize}
        \item \(\delta_x(\overline{f}) =\overline{\delta_x(f)}\)
        \item \(D_v\delta_x(\overline{f}) = D_v\overline{f(x)}= \overline{D_v f(x)} =
        \overline{D_v\delta_x(f)}\) because from \(|z| = |\overline{z}|\) follows
        \[
            \lim_{t\to 0}\frac{|\overline{f(x+tv)} - \overline{f(x)} - \overline{D_vf(x)} t|}{t}
            = 0,
        \]
        \item induction for higher order derivatives.
    \end{itemize}
    So we have that \(P(t):=Te^{\im\langle \cdot, t\rangle}\) must be zero \(\spectMeasure\)-almost everywhere.
    Since
    \[
        P(t)
        = \sum_{j=1}^n \sum_{l=0}^k \sum_{i_1\le \dots \le i_l} w_j^{(l, i_1,\dots, i_l)}\partial_{i_1\dots i_l}e^{\im\langle x_j, t\rangle}
    \]
    is continuous in \(t\), it can only be zero \(\spectMeasure\)-almost everywhere if it is zero on
    the support of \(\spectMeasure\). Since the support of the spectral measure
    \(\spectMeasure\) contains an open subset by assumption of the theorem, \(P\)
    must be zero on this open subset. But then \(P\) must be zero everywhere as an
    analytic function. As \(P(t)\) is the fourier transform of \(T\) in the
    sense of distributions, i.e.
    \[
        P(t)
        = \int \Bigl(\sum_{j=1}^n \sum_{l=0}^k \sum_{i_1\le \dots \le i_l} w_j^{(l, i_1,\dots, i_l)}\partial_{i_1\dots i_l}\delta_{x_j}\Bigr)(x)
        e^{\im\langle x, t\rangle}dx
        = \Ft[T](t),
    \]
    this requires \(T\) to be zero by linearity and invertibility of the Fourier transform. But since
    the \(\partial_{i_1\dots i_l}\delta_{x_j}\) are linear independent\footnote{
        To see the linear independence of \(\partial_{i_1\dots
        i_l}\delta_{x_j}\), consider that the finite \(x_k\) are distinct, so
        they have a minimal distance. Rescale a bump
        function with zero slope at the top, e.g.
        \[
            \phi(x) = \begin{cases}
                \exp\Bigl(-\frac{1}{1-\|x\|^2}\Bigr) & \|x\| < 1
                \\
                0 & \|x\| \ge 1
            \end{cases}
        \]
        such that it is centered on some \(x_j\) and it is zero at all other \(x_k\)
        (and zero at all derivatives). This implies
        \[
            0 = \langle T, \phi\rangle = w_j^{(0)} \phi(x_j)
        \]
        and thus \(w_j^{(0)}=0\).
        Then construct similar test functions to ensure the other prefactors have to
        be zero, by placing a non-zero slope at \(x_j\) while ensuring it is zero at
        all other \(x_k\).
    }
    for distinct \(x_j\) the only way \(T\) can be zero is, if all the \(w_j\)
    are zero. This finishes the proof.
\end{proof}

The general requirement, the existence of a non-empty open sets in the support
of the spectral measure, is satisfied for the stationary isotropic random
functions valid in all dimensions (Corollary~\ref{cor: strict positive
definite}). To prove this result we require the following lemma.
It proves that stationary isotropic random function in \(\ell^2\)
is either almost surely constant or the `Schoenberg measure' has positive measure
on \((0,\infty)\). Where the Schoenberg measure \(\schoenbergMeas\) refers to 
the measure in Schoenberg's characterization of stationary isotropic covariance
kernels on \(\ell^2\) given by\footnote{
    cf.\ Section \ref{sec: characterization of covariance kernels}, Table
    \ref{table: characterizations}
}
\begin{equation}
    \label{eq: schoenberg charact}    
    \C_\rf(x,y) = \ikernel(\tfrac{\|x-y\|^2}2) = \int_{[0,\infty)} \exp(-t^2 \tfrac{\|x-y\|^2}2)\schoenbergMeas(dt)
\end{equation}

\begin{lemma}[Constant random functions]
    \label{lem: constant random functions}
    Let \(\ikernel\) be a stationary isotropic covariance kernel valid in all
    dimensions and let \(\rf\sim \normal(\mu, \ikernel)\) be a
    continuous stationary isotropic random function.

    If the Schoenberg measure \(\schoenbergMeas\) of \(\ikernel\) has no mass on
    \((0,\infty)\), i.e.  \(\schoenbergMeas((0,\infty))=0\), then \(\rf\)
    is almost surely constant.
\end{lemma}
\begin{proof}
    By Schoenberg's characterization of stationary isotropic positive definite
    kernels in \(\ell^2\) \eqref{eq: schoenberg charact} we have
    \[
        \ikernel(r)
        = \int_{[0,\infty)} \exp(-t^2r)\schoenbergMeas(dt)
        \overset{\schoenbergMeas((0,\infty)=0)}= \schoenbergMeas(\{0\}).
    \]
    This implies a constant covariance
    \[
        \Cov(\rf(x), \rf(y)) = \schoenbergMeas(\{0\}).
    \]
    Assuming \(\rf\) to be centered (without loss of generality)
    by switching to \(\tilde{\rf}_\dims = \rf - \mu\), this implies
    \[
        \E[(\rf(x) - \rf(y))^2] = \E[\rf(x)^2] - 2\E[\rf(x)\rf(y)] + \E[\rf(y)^2] = 0.
    \]
    Thus \(\rf(x) = \rf(y)\) almost surely for all \(x,y\). Via the union over a
    dense countable subset of \(\real^\dims\) and \as continuity of \(\rf\) we get
    \[
        \Pr(\rf \text{ \as constant})
        = \Pr\Bigl(\rf(x) = \rf(y) \quad \forall x,y \in \real^\dims\Bigr)
        = 1.
        \qedhere
    \]
\end{proof}

Using this lemma, we can prove that all stationary isotropic covariance kernels
on \(\ell^2\) have strictly positive definite derivatives. 

\corStrictPosDef*
\begin{proof}[Proof of Corollary \ref{cor: strict positive definite}]
    Since \(\charfct(x) = e^{-\frac{\|x\|^2t^2}2}\) is the characteristic function of
    \(Y_t\sim\normal(0, t\identity_{\dims\times\dims})\), Schoenbergs
    characterization \eqref{eq: schoenberg charact} of the covariance  implies up to scaling
    \begin{align*}
        \C_{\rf}(x,\tilde{x})
        &= \ikernel\bigl(\tfrac{\|x-\tilde{x}\|^2}2\bigr)
        = \int_{[0,\infty)} e^{-\frac{t^2\|x-\tilde{x}\|^2}2} \schoenbergMeas(dt)
        \\
        &= \int_{[0,\infty)} \E[e^{\im\langle x-\tilde{x}, Y_t\rangle}] \schoenbergMeas(dt)
        \\
        &= \int e^{\im\langle x-\tilde{x}, y\rangle}\underbrace{
            \int_{(0,\infty)}  e^{-\frac{\|y\|^2}{2t}} \schoenbergMeas(dt)
        }_{=:\density_\spectMeasure(y)} dy
        + \schoenbergMeas(\{0\})
        \\
        &= \int e^{\im\langle x-\tilde{x}, y\rangle} \spectMeasure(dy)
    \end{align*}
    with spectral measure
    \[
        \spectMeasure(A) :=\int_A \density_\spectMeasure(y) dy + \tfrac1\dims\schoenbergMeas\{0\} \delta_0(A).
    \]
    Since \(\schoenbergMeas((0,\infty))\neq 0\) because \(\rf\) is not almost surely
    constant (Corollary~\ref{lem: constant random functions}),
    its density \(\density_\spectMeasure\) is continuous and positive
    \(\density_\spectMeasure(y) >0\) for all \(y\in \real^\dims\) and thus
    \(\support(\spectMeasure)= \real^\dims\). In particular Theorem~\ref{thm:
    strict positive definite derivatives} is applicable and finishes the proof.
\end{proof}
 
\printbibliography[heading=subbibliography]

\begin{subappendices}
{
    \renewcommand*{\kernel}{K}
    \section{Appendix}

Since we are interested in \(\real^\dims\) we consider \(\sphere\), whereas the
source we build our characterization of isotropic kernels on \autocite{guellaSchoenbergsTheoremPositive2019}
considers \(\sphere[\dims]\). We also work with the normalized Gegenbauer
polynomials while they stated their result with the unnormalized ones. For
your convenience we therefore restate their main theorem in our setting and notation.

\begin{theorem}[\(\domain \times \sphere\)-characterization, \cite{guellaSchoenbergsTheoremPositive2019}]
	\label{thm: guella XxS characterization}
	Let \(\dims \in \{2,3,\dots, \infty\}\) and define \(\lambda:=\frac{\dims-2}2\). For a
	kernel \(\kernel: (\domain \times \sphere)^2 \to \complex\) which is
	isotropic on the sphere, i.e.  of
	the form 
	\[
        \kernel((r,v), (s,w)) = f_{\sphere[]}(r,s, \langle v ,w\rangle),
	\]
	the following are equivalent 
	\begin{enumerate}[label=(\roman*)]
		\item The kernel \(\kernel \) is positive definite and for all $r,s\in \domain$ the function
		\(f_{\sphere[]}(r,s, \rho)\) is continuous in \([-1,1]\).
		
		\item The kernel has series representation of the form
		\[
                f_{\sphere[]}(r,s, \rho)
		= \sum_{n=0}^\infty \alpha^{(\dims)}_n(r,s) \uPnorm_n(\rho), \qquad \rho\in [-1,1],
		\]
		where \(\alpha^{(\dims)}_n: \domain^2 \to \complex\) are positive definite kernels for all \(n\in \nat_0\)
		that are summable, i.e.  \(\sum_{n=0}^\infty \alpha^{(\dims)}_n(r,r) <
		\infty\) for all \(r\in \domain\).
	\end{enumerate}
	If the kernel satisfies one of these equivalent conditions and \(\dims<\infty\)
	the kernels \(\alpha^{(\dims)}_n\) moreover have the representation
	\[
        \alpha^{(\dims)}_n(r,s) = Z \int_{-1}^1 f_{\sphere[]}(r,s, \rho) \uPnorm_n(\rho)(1-\rho^2)^{(\dims-3)/2} d\rho
	\]
	with normalizing constant
	\[
	Z = Z(n, \dims) = \Bigl[\int_{-1}^{1} [\uPnorm_n(\rho)]^2(1-\rho)^{(\dims-3)/2}d\rho\Bigr]^{-1}.
	\]
\end{theorem}
\begin{proof}
    This is essentially Theorem~2.3 by \citet{guellaSchoenbergsTheoremPositive2019} except
    that we translate their \(a_n^{\dims}(r,s)\) into \(\alpha^{(\dims)}_n(r,s) := a_n^{(\dims)}(r,s) \uP_n(1)\)
    for the use with \emph{normalized} ultraspherical polynomials.
\end{proof}

A key component for the characterization of invariant kernels and distributionally
invariant random functions is a characterization of isometries that is most likely
well known.

\begin{prop}[Characterizing isometries]\label{prop: characterization of isometries}
    Let \(\cX\) be a vector space over \(\field\in \{\real, \complex\}\) and \(x_i,y_i \in \cX\) for \(i\in 1,\dots, n\),
    then the following pairs of statements are equivalent
	\begin{enumerate}
        \item \label{item: translation characterization}
        \begin{enumerate}
            \item\label{itm: translation charac. b}
            there exists a \textbf{translation}\footnote{
                A translation is a map on \(\cX\) of the form \(x\mapsto x +v\) for some \(v\in \cX\)
            } \(\varphi\) with \(\varphi(x_i) = y_i\) for all
                \(i\).

            \item\label{itm: translation charac. a}
            \(x_i - x_j = y_i - y_j\) for all \(i,j\)
        \end{enumerate}
    \end{enumerate}
    If \(\cX\) is furthermore a hilbertspace, then
    \begin{enumerate}[resume*]
		 \item \label{item: linear isometry characterization}
		 \begin{enumerate}
            \item \label{itm: lin. isometry charac. lin isom}
            there exists a \textbf{linear isometry}\footnote{
                An isometry \(\phi\colon \cX\to\cX\) is a map that retains distances,
                \ie for all \(x,y\in\cX\)
                \[
                    \|\phi(x) - \phi(y)\| = \|x-y\|.
                \]
                For linear maps the isometry property is equivalent to norm
                preservation, \ie \(\|\phi(x)\| = \|x\|\) for all \(x\in \cX\).
            } \(\phi\) with \(\phi(x_i) = y_i\) for all \(i\)
            
            \item \label{itm: lin. isometry charac. norm}
            \(\|x_i - x_j\| = \|y_i-y_j\|\) and \(\|x_i\|=\|y_i\|\) for all \(i,j\) 
            
            \item \label{itm: lin. isometry charac. inner prod}
            \(\langle x_i, y_i\rangle\) for all \(i,j\)
		 \end{enumerate}

		 \item \label{item: affine isometry characterization}
		 \begin{enumerate}
			  \item\label{itm: affine isometry charac. a}
			  there exists an \textbf{isometry} \(\phi\) with \(\phi(x_i) = y_i\) for all \(i\).

			  \item\label{itm: affine isometry charac. b}
			  \(\|x_i - x_j\| = \|y_i - y_j\|\) for all \(i,j\)

              \item\label{itm: affine isometry charac. affine map}
              there exists a linear isometry \(\phi\) and \(v\in \cX\)
              such that \(y_i = \phi(x_i) + v\) for all \(i\)
		 \end{enumerate}
	\end{enumerate}
\end{prop}

\begin{proof}
	\begin{description}
        \item[\eqref{itm: translation charac. a} \(\Rightarrow\) \eqref{itm: translation charac. b}:] 
        we define \(\varphi(x) := x + (y_0-x_0)\), which implies
        \[
            \varphi(x_i) = x_i - x_0 + y_0 = (y_i-y_0) + y_0 = y_i.
        \]
        
        \item[\eqref{itm: translation charac. b} \(\Rightarrow\) \eqref{itm: translation charac. a}:] 
        Let \(\varphi(x) = x + c\) for some \(c\). Then we immediately have
        \[
            y_i - y_j = \varphi(x_i) - \varphi(x_j) = x_i + c - (x_j + c) = x_i - x_j.
        \]

        \item[\eqref{itm: lin. isometry charac. lin isom} \(\Rightarrow\) \eqref{itm: lin. isometry charac. norm}:] 
        Isometries preserve distances by definition. This implies \(\|x_i - x_j\|
        = \|y_i - y_j\|\). And linear functions map
        \(0\) to \(0\), so we have 
        \[
            \|x_i\| = \|x_i - 0\| = \|\phi(x_i) - \phi(0)\| = \|y_i\|.
        \]

        \item[\eqref{itm: lin. isometry charac. norm} \(\Rightarrow\) \eqref{itm: lin. isometry charac. inner prod}:] 
        The polarization formula \(\langle x,y\rangle = \frac12(\|x\|^2 +
        \|y\|^2 - \|x-y\|^2)\) yields the claim.

        \item[\eqref{itm: lin. isometry charac. inner prod} \(\Rightarrow\) \eqref{itm: lin. isometry charac. lin isom}:] 
        With the application of the Gram-Schmidt orthonormalization procedure to
        both \(x_i\) and \(y_i\) we obtain
        \begin{align*}
            U_m &:= \Span\{u_1, \dots, u_{k_m}\} = \Span\{x_1, \dots, x_m\}
            \\
            V_m &:= \Span\{v_1, \dots, v_{\tilde k_m}\} = \Span\{y_1, \dots, y_m\}
        \end{align*}
        for orthonormal \(u_i\) and \(v_i\) where we skip \(x_m\) (or \(y_m\) respectively)
        if it is already contained in \(U_{k_{m-1}}\) (or \(V_{k_{m-1}}\) resp.) resulting in \(k_m=k_{m-1}\)
        (or \(\tilde k_m=\tilde k_{m-1}\)) respectively. We will prove inductively over
        \(m\le n\) that \(k_m = \tilde k_m\) and for all \(i \le n\) and \(j\le k_m\)
        \begin{equation}
            \label{eq: inner products remain the same}    
            \langle x_i, u_j\rangle = \langle y_i, v_j\rangle.
        \end{equation}
        The induction start with \(m=0\) is trivial. For the induction step
        observe that \(x_m\) being contained in \(U_{k_{m-1}}\) is equivalent to
        \[
            \langle x_m, x_m\rangle
            = \Bigl\langle x_m, \sum_{j=1}^{k_{m-1}} \langle x_m, u_j\rangle u_j \Bigr\rangle
            = \sum_{j=1}^{k_{m-1}} \langle x_m, u_j\rangle^2
        \]
        And since all the inner products are the same by induction \eqref{eq: inner products remain the same}, \(x_m \in U_{k_{m-1}}\) is
        satisfied if and only if \(y_m\in V_{k_{m-1}} = V_{\tilde k_{m-1}}\). We therefore have
        \(k_m = \tilde k_m\). In the case where \(k_m\) increments we have by the induction
        assumption \eqref{eq: inner products remain the same} and \eqref{itm: lin. isometry charac. inner prod}
        \begin{align*}
            \langle x_i, u_{k_m}\rangle
            &= \Bigl\langle x_i, x_m - \sum_{j=1}^{k_{m-1}}\langle x_m, u_j\rangle u_j\Bigr\rangle
            \\
            \overset{\text{\eqref{eq: inner products remain the same}, \eqref{itm: lin. isometry charac. inner prod}}}&= \langle y_i, y_m\rangle - \sum_{j=1}^{k_{m-1}} \langle y_i, v_j\rangle \langle y_m, v_j\rangle
            = \langle y_i, v_{k_m}\rangle
        \end{align*}
        for all \(x_i\) with \(i\le n\). This finishes our induction.

        Since \(W = \Span(U_n \cup V_n)\) has finite dimension \(\dims\le 2n\),
        we can extend the orthonormal basis \(u_1,\dots, u_{k_n}\) of \(U_n = \Span\{x_1,\dots,x_n\}\)
        to an orthonormal basis \(u_1,\dots, u_\dims\) of \(W\)
        and similarly the orthonormal basis \(v_1,\dots, v_{k_n}\) of \(V_n =
        \Span\{y_1,\dots, y_n\}\) to an orthonormal basis \(v_1,\dots, v_\dims\) of \(W\).

        We finally define the linear map \(\phi\colon \cX \to \cX\), which acts as the
        identity on \(W^\perp\) and maps the basis elements \(u_j\) to \(v_j\).
        This is an isometry since for all \(x\in \cX\) there exists 
        a basis representation \(x=\sum_{i=1}^n \lambda_i u_i + w\) with
        \(\lambda_i\in \field^n\) and \(w\in W\), which implies
        \[
            \|\phi(x)\|^2 
            = \bigl\|\sum_{i=1}^n \lambda_i \phi(v_i) + w\Bigr\|^2
            = \|\lambda\|^2 + \|w\|^2
            = \|x\|^2.
        \]
        This isometry does the job, because
        \[
            \phi(x_m)
            = \phi\Bigl(\sum_{j=1}^n \langle x_m, u_j\rangle u_j\Bigr)
            = \sum_{j=1}^n \langle x_m, u_j\rangle \phi(u_j)
            \overset{\eqref{eq: inner products remain the same}}= \sum_{j=1}^n \langle y_m, v_j\rangle v_j 
            = y_m.
        \]

        \item[\eqref{itm: affine isometry charac. a} \(\Rightarrow\) \eqref{itm: affine isometry charac. b}] 
        This is precisely the distance preserving property of isometries.

        \item[\eqref{itm: affine isometry charac. b} \(\Rightarrow\) \eqref{itm: affine isometry charac. affine map}:] 
        We define \(\tilde{x}_i = x_i - x_0\)
        and similarly for \(y\). In particular, \(\tilde{x}_1 = \tilde{y}_1=0\). Since
        \(\tilde{x}_i\) and \(\tilde{y}_i\) satisfy the requirements of
        \eqref{itm: lin. isometry charac. norm}, there exists a linear
        isometry matrix \(\phi\) with \(\phi(\tilde{x}_i) = \tilde{y}_i\).
        Define \(v = y_1 - \phi(x_1)\)
        to obtain
        \[
            \phi(x_i) + v = \phi(\tilde{x}_i) + y_1 = \tilde{y}_i + y_1 = y_i
        \]

        \item[\eqref{itm: affine isometry charac. affine map} \(\Rightarrow\) \eqref{itm: affine isometry charac. a}:] 
        It is straightforward to show \(\varphi(x) = \phi(x) + v\) is an isometry.
        \qedhere
	\end{description}
\end{proof}

 }
\end{subappendices}

 	\end{refsection}

	\begin{refsection}
		{
			\renewcommand*{\param}{x}
			\renewcommand*{\Param}{X}
			\chapter{Random Function Descent}
\label{chap: random function descent}

In Section~\ref{sec: practical bayesian optimization} we have shown that
the most popular acquisition functions in Bayesian optimization can be reduced
to functions of the posterior mean and covariance in the Gaussian case.
While this certainly makes these algorithms more explicit, it has still
left us with the task of computing the conditional mean and standard deviation
and the optimization of the acquisition function. In Section~\ref{sec: conditional
Gaussian distribution} we have learned how such a conditional Gaussian distribution
can be computed. Recall that to compute the conditional distribution \(Z_{0}
\mid Z_1, \dots Z_N\) of the multivariate Gaussian random vector \((Z_0,\dots, Z_N)\)
the covariance matrix of \(Z_1,\dots, Z_N\) has to be inverted. Further recall
that gradients of a Gaussian random function \(\obj:\real^\dims \to \real\)
contain \(\dims\) Gaussian entries. In the case of first order Bayesian
optimization the evaluation filtration
\(\filt_n\) \eqref{Bayes-eq: filtration} therefore contains \(\bigO(n\dims)\)
Gaussian random variables.  Inverting their covariance matrix thus incurs a cost
of \(\bigO(n^3\dims^3)\) in general, which renders this approach infeasible in
high-dimension.

In contrast to these computationally complex Bayesian optimization algorithms
classical algorithms, such as gradient descent, are very explicit and cheap
to compute with \(\bigO(n\dims)\) complexity. Taking inspiration
from the classical approach we introduce the acquisition function 
`Random Function Descent' in this chapter. It turns out that the minimizer
of this acquisition function can be explicitly computed and is in fact
equal to gradient descent with a specific step size.
While the step size generally has to be tuned in the classical setting
this Bayesian view provides a theoretical foundation for the choice of step
size.

\section{The random function descent algorithm}\label{sec: rfd}

Classically, gradient descent is motivated as the minimizer of
the trust bound based on the first Taylor approximation and a bound
on the error (Lemma~\ref{lem: taylor approximation with trust bound}).
In the case of \(L\)-smoothness it is of the form
\[
	\obj(y) \le \underbrace{
		\obj(\param) + \langle \nabla\obj(\param), y-\param\rangle
	}_{
		=: T[\obj(y) \mid \obj(\param), \nabla\obj(\param)]
	} + \tfrac{L}2\|y-\param\|.
\]
While the function \(\obj\) may be hard to minimize, this trust bound
is a quadratic function and thereby easy to minimize with
minimizer \(\param - \tfrac1L \nabla \obj(\param)\) (cf. Lemma \ref{lem: general
descent lemma}).

The eccentric notation for the Taylor approximation,\footnote{
	which is a linear approximation of \(\obj\) in \(\param\)
} \(T[ \obj(y)\mid \obj(\param),\nabla \obj(\param)]\), is meant to
highlight the connection to the stochastic Taylor approximation we define below.

\begin{figure}
	\begin{sidecaption}[outer]{
		The stochastic Taylor approximation naturally contains a trust bound
		in contrast to the classical one. Here \(\Obj\) is a Gaussian
		random function (with covariance as in Equation~\eqref{eq: sqExp
		covariance model}, with length scale \(\scale=2\) and variance \(\sigma^2=1\)).
		The ribbon represents two conditional standard
		deviations around the conditional expectation.
	\label{fig: visualize conditional expectation}
	}
	\centering
	\def\svgwidth{1\linewidth}
	\vspace*{-2cm}
	\begingroup \makeatletter \providecommand\color[2][]{\errmessage{(Inkscape) Color is used for the text in Inkscape, but the package 'color.sty' is not loaded}\renewcommand\color[2][]{}}\providecommand\transparent[1]{\errmessage{(Inkscape) Transparency is used (non-zero) for the text in Inkscape, but the package 'transparent.sty' is not loaded}\renewcommand\transparent[1]{}}\providecommand\rotatebox[2]{#2}\newcommand*\fsize{\dimexpr\f@size pt\relax}\newcommand*\lineheight[1]{\fontsize{\fsize}{#1\fsize}\selectfont}\ifx\svgwidth\undefined \setlength{\unitlength}{300bp}\ifx\svgscale\undefined \relax \else \setlength{\unitlength}{\unitlength * \real{\svgscale}}\fi \else \setlength{\unitlength}{\svgwidth}\fi \global\let\svgwidth\undefined \global\let\svgscale\undefined \makeatother \begin{picture}(1,0.75)\lineheight{1}\setlength\tabcolsep{0pt}\put(0,0){\includegraphics[width=\unitlength,page=1]{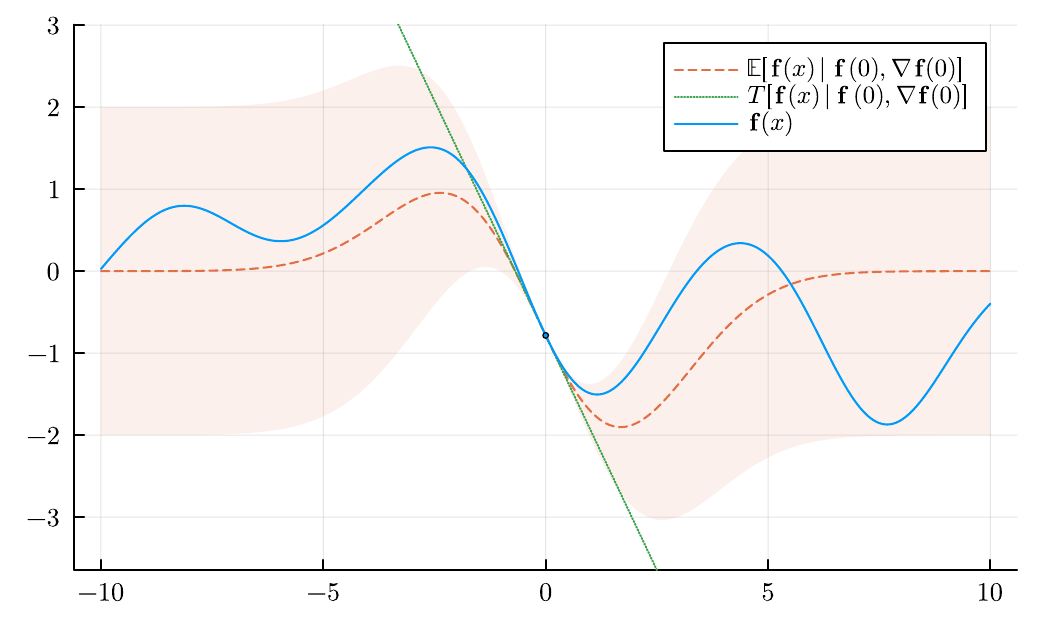}}\end{picture}\endgroup  	\end{sidecaption}	
\end{figure}

\begin{definition}[Stochastic Taylor approximation]
	We define the first order stochastic Taylor approximation of a random
	(cost) function \(\Obj\) around \(\param\) by the
	conditional expectation
	\[
		\E[\Obj(y) \mid \Obj(\param), \nabla\Obj(\param)].
	\]
	This is the best \(L^2\) approximation
	\autocite[Cor.~8.17]{klenkeProbabilityTheoryComprehensive2014} of
	\(\Obj(y)\) provided first order knowledge of \(\Obj\) around
	\(\param\).
\end{definition}

We call this the `stochastic Taylor approximation' because this approximation
only makes use of derivatives in a single point. While the standard Taylor
approximation is a polynomial approximation, the `stochastic Taylor
approximation' is the best approximation in an \(L^2\) sense.
Since the stochastic Taylor approximation is already
mean-reverting by itself, it naturally incorporates covariance-based trust
(cf.~Figure~\ref{fig: visualize conditional expectation})
so we omit a trust bound. A variance based trust may be
added (cf. Section \ref{sec: conservative rfd}) but
to obtain step sizes that are optimal on \emph{average}
we start without.\footnote{
	As we will see in Chapter \ref{chap: predictable progres}, the average case
	is highly representative of the actual behavior in high-dimension.
}\footnote{
	Observe that the UCB algorithm from Example~\ref{ex: acquisition functions}
	(for maximization!) does the opposite of a trust bound. Instead of
	a penalty for uncertainty, it encourages exploration
	with a variance-based bonus.
}

While \(L\)-smoothness-based trust \emph{guarantees} that the gradient still
points in the direction we are going (for learning rates smaller \(1/L\)), covariance
based trust tells us whether the derivative is still negative \emph{on average}.
Minimizing the stochastic Taylor approximation is therefore optimized for the
average case. Since convergence proofs for gradient descent typically rely on an
improvement \emph{guarantee}\footnote{the `Descent Lemma'
(Lemma \ref{lem: general descent lemma})}, proving convergence is significantly
harder in the average case and we answer this question only partially in
Corollary~\ref{prop: convergence}.

\begin{definition}[Random Function Descent -- RFD]
	\label{def: rfd}
	Select \(\param_{n+1}\) as the minimizer\footnote{
		we ignore throughout the main body that \(\argmin\) could be set-valued and
		that the \(\param_n\) would be random variables (cf.~Section~\ref{sec: formal rfd}
		for a formal approach).
	}
	of the first order stochastic Taylor
	approximation
	\begin{equation*}
\param_{n+1}
		:= \argmin_{\param}\E[\Obj(\param)\mid \Obj(\param_n), \nabla\Obj(\param_n)].
	\end{equation*}
\end{definition}

\paragraph*{Properties of RFD}
Before we make RFD more explicit in Section~\ref{subsec: explicit rfd}, we
discuss some properties which are easier to see in the abstract Bayesian form.

First, observe that RFD is greedy and forgetful in the same way gradient descent is greedy
and forgetful when derived as the minimizer of the regularized first Taylor
approximation, or the Newton method as the minimizer of the second Taylor
approximation. This is because the Taylor approximation only uses derivatives
from the last point \(\param_n\) (forgetful), and we minimize this approximation
(greedy). Since momentum methods retain some information about past gradients,
they are not as forgetful. We therefore expect a similar improvement
could be made for RFD in the future.

Second, it is well known that classical gradient descent with exogenous step sizes (and
most other first order methods) lack the scale invariance property of the Newton
method \autocite[e.g.][]{hanssonOptimizationLearningControl2023,deuflhardAffineInvariantConvergence1979}. Scale invariance means that scaling the input parameters \(\param\) or the
cost itself (e.g. by switching from the mean squared error to the sum squared
error) does not change the points selected by the optimization method.

\begin{restatable}[Scale invariance and equivariance]{advantage}{scaleInvariance}
	\label{advant: scale invariance}
	RFD is invariant to additive shifts and positive scaling of the cost
	\(\Obj\). RFD is also equivariant\footnote{
		a function \(\phi\) is \emph{invariant} with respect to a group \(G\),
		if \(\phi(\rho(g) x) = \phi(x)\), where \(\rho(g)\) is a
		representation of the element \(g\in G\). \(\phi\) is \emph{equivariant}
		with respect to the group \(G\), if \(\phi(\rho(g) x) = \rho(g)\phi(x)\)
		for all \(g\in G\).
	} with respect to transformations of the
	parameter input of \(\Obj\) by differentiable bijections whose Jacobian is
	invertible everywhere (e.g. invertible linear maps).
\end{restatable}
While equivariance with respect to bijections of inputs is much stronger than
the affine equivariance offered by the Newton method, non-linear bijections will
typically break the isotropy assumption (Chapter \ref{chap: distributions over functions}) about the random function \(\Obj\) which makes RFD explicit.
This equivariance should therefore be viewed as an
opportunity to look for the bijection of inputs which ensures isotropy (e.g. a
whitening transformation). The discussion of geometric anisotropy in
Section~\ref{sec: geometric anisotropy} is conducive to build an understanding
of this.
 \section{Relation to gradient descent}
\label{subsec: explicit rfd}

While we were able to define RFD abstractly without any assumptions on the
distribution \(\Pr_{\Obj}\) of the random cost \(\Obj\), an explicit
calculation of this acquisition function requires distributional assumptions.
As we motivated in Section~\ref{sec: practical bayesian optimization}
and Chapter~\ref{chap: distributions over functions}, a sensible
assumption are stationary isotropic Gaussian random functions.
This assumption allows for an explicit version of the stochastic Taylor
approximation which then immediately leads to an explicit version of RFD.

\begin{restatable}[Explicit first order stochastic Taylor approximation]{lemma}{firstStochTaylor}
	\label{lem: first stoch Taylor}
	For \(\Obj\sim\normal(\mu, \ikernel)\), the first order stochastic Taylor
	approximation is given by
	\[
		\E[\Obj(\param-\step)\mid \Obj(\param),\nabla\Obj(\param)]
		= \mu + \frac{\ikernel\bigl(\frac{\|\step\|^2}2\bigr)}{\ikernel(0)}
		(\Obj(\param)-\mu) - \frac{\ikernel'\bigl(\frac{\|\step\|^2}2\bigr)}{\ikernel'(0)}
		\langle \step, \nabla\Obj(\param)\rangle.
	\]
\end{restatable}

The explicit version of RFD follows by fixing the step size
\(\stepsize = \|\step\|\) and optimizing over the direction first.

\begin{restatable}[Explicit RFD]{theorem}{explicitRFD}
	\label{thm: explicit rfd}
	Let \(\Obj\sim\normal(\mu, \ikernel)\), then RFD coincides with gradient
	descent
	\[
		\param_{n+1}
		= \param_n  -\stepsize_n^*\tfrac{\nabla\Obj(\param_n)}{\|\nabla\Obj(\param_n)\|},
	\]
	where the RFD step sizes are given by
	\begin{equation}
		\label{eq: rfd step size}
		\stepsize_n^*
		:= \argmin_{\stepsize\in\real} \frac{\ikernel\bigl(\frac{\stepsize^2}2\bigr)}{\ikernel(0)}(\Obj(\param_n)-\mu)
		- \stepsize \frac{\ikernel'\bigl(\frac{\stepsize^2}2\bigr)}{\ikernel'(0)}\|\nabla\Obj(\param_n)\|.
	\end{equation}
\end{restatable}

\begin{table*}[t]
\begin{fullwidth}
	\centering
	\begin{tabular}{p{1.8cm} llll l}
	  	\toprule

	  	\multicolumn{2}{c}{Model} & \multicolumn{2}{c}{RFD step size \(\stepsize^*\) for \(\Obj(\param) \le \mu\)} 
		& A-RFD
		\\
		\cmidrule(r){3-5}
		& & General case \(\bigl(\text{with }\Theta=\tfrac{\|\nabla\Obj(\param)\|}{\mu-\Obj(\param)}\bigr)\)
		& \(\Obj(\param) = \mu\)
		& \(\Theta\to 0\)
		\\
	  	\midrule
		Matérn & \(\nu\) \\
	  	\cmidrule(r){1-2}
		& \(3/2\) 
		& \(
			\frac{\scale}{\sqrt{3}}\frac{1}{\left(
				1 + \frac{\sqrt{3}}{s\Theta}
			\right)}
		\)
		& \(\approx 0.58\scale\) 
		& \(\frac{1}{3}\scale^2 \Theta\)
		\\
		& \(5/2\) 
		& \(
			\frac{\scale}{\sqrt{5}}\frac{
				(1-\zeta)+\sqrt{4 + (1+\zeta)^2}
			}{2(1+\zeta)}\)
			with \(\zeta := \frac{\sqrt{5}}{3\scale\Theta}.
		\)
		& \(\approx 0.72\scale\)
		& \(\frac{3}{5} \scale^2\Theta\)
		\\
		Squared-exponential
		& \(\infty\)
		& \(
			\frac{
				\scale^2
			}{
				\sqrt{\bigl(\frac{\mu-\Obj(\param)}{2}\bigr)^2+\scale^2\|\nabla\Obj(\param)\|^2}
				+ \frac{\mu-\Obj(\param)}{2}
			}\|\nabla\Obj(\param)\|
		\)
		& \(\scale\)
		& \(\scale^2\Theta\)
		\\
	  	\cmidrule(r){1-2}
		Rational quadratic & \(\beta\)
		& \(
		\scale \sqrt{\beta} \Root\limits_\stepsize\Bigl\{
			- 1 + \frac{\sqrt{\beta}}{\scale\Theta}\stepsize
			+ (1+\beta)\stepsize^2 + \frac{\sqrt{\beta}}{\scale\Theta}\stepsize^3
		\Bigr\}
		\) &
		\(
			\scale \sqrt{\frac{\beta}{1+\beta}}
		\)
		& \(\scale^2\Theta\)
		\\
		\bottomrule
	\end{tabular}
    \end{fullwidth}
	\sideparmargin{outer}
	\vspace{-\baselineskip}
	\sidepar{\vspace{\baselineskip}
	\caption{
		RFD step size (cf.~Figure~\ref{fig: rfd step sizes} and
		Eq.~\eqref{eq: sqExp covariance model}, \eqref{eq: rational
		quadratic}, \eqref{eq: matern model} for the formal definitions of the models).
		In particular, \(\scale\) is the length scale in all covariance models.
	} \label{table: optimal step size}
	}
\end{table*} 
While the descent direction is a universal property for all isotropic
Gaussian random functions, it follows from \eqref{eq: rfd step size} that the
step sizes depend much more on the specific covariance structure. In particular
it depends on the decay rate of the covariance acting as the trust bound.

\begin{remark}[Scalable complexity]
	While Bayesian optimization typically has computational complexity \(\bigO(n^3\dims^3)\)
	in number of steps \(n\) and dimensions \(\dims\)
	\autocite{wuBayesianOptimizationGradients2017,roosHighDimensionalGaussianProcess2021},
	RFD under the isotropy assumption has the same computational complexity as
	gradient descent (i.e. \(\bigO(n\dims)\)).
\end{remark}

\begin{remark}[Step until the given information is no longer informative]
While \(L\)-smoothness-based trust prescribes step sizes that \emph{guarantee}
the slope to point downwards over the entire step, RFD prescribes steps which
are exactly large enough that the gradient is no longer correlated to the
previously observed evaluation. This is because the first order condition
demands
\[
	0 \overset{!}= \nabla \E[\Obj(\param)\mid \Obj(\param_n), \nabla\Obj(\param_n)]
	= \E[\nabla \Obj(\param)\mid \Obj(\param_n), \nabla\Obj(\param_n)].
\]
And for measurable functions \(\phi:\real^{\dims+1}\to\real\) such that \(\Phi=\phi(\Obj(\param_n),
\nabla\Obj(\param_n))\) is sufficiently integrable, \(\Phi\) is then uncorrelated
from \(\partial_{i}\Obj(\param)\) by the first order condition
\[
	\Cov(\partial_{i}\Obj(\param), \Phi)
	= \E\Bigl[
		\underbrace{\E[\partial_{i}\Obj(\param) \mid \Obj(\param_n), \nabla\Obj(\param_n)]}_{=0}
		(\Phi-\E[\Phi]) \Bigr] = 0.
\]
\end{remark}

\section{The RFD step size schedule}
\label{subsec: rfd step sizes}

While classical \(L\)-smooth theory leads to `learning rates', RFD suggests
`step sizes' applied to normalized gradients\footnote{
	which we also encountered in the worst case theory with the modulus of
	continuity (Remark~\ref{rem: step size vs. learning rate})
}
representing the actual length of the step size. In the following we thus make the distinction
\[
	\param_{n+1}
	= \param_n - \underbrace{\lr_n}_{\mathclap{\text{`learning rate'}}} \nabla\Obj(\param_n)
	= \param_n - \underbrace{\stepsize_n}_{\text{\text{`step size'}}}
	\tfrac{\nabla\Obj(\param_n)}{\|\nabla\Obj(\param_n)\|}.
\]
To get a better feel for the step sizes suggested by RFD, it is enlightening to
divide \eqref{eq: rfd step size} by \(\mu-\Obj(\param_n)\) which results in the
equivalent minimization problem
\begin{equation}
	\label{eq: simplified step size opt}
	\stepsize^* : =\stepsize^*(\Theta) := \argmin_\stepsize q_\Theta(\stepsize)
	\qquad\text{for} \qquad
	q_\Theta(\stepsize) := -\frac{\ikernel\bigl(\frac{\stepsize^2}2\bigr)}{\ikernel(0)}
		- \stepsize \frac{\ikernel'\bigl(\frac{\stepsize^2}2\bigr)}{\ikernel'(0)}\Theta,
\end{equation}
which is only parametrized by the ``gradient cost quotient''
\[
	\Theta_n = \frac{\|\nabla\Obj(\param_n)\|}{\mu-\Obj(\param_n)},
\]
\begin{figure}[b]
	\begin{sidecaption}{
		RFD step sizes as a function of \(\Theta=\frac{\|\nabla\Obj(\param)\|}{\mu -
		\Obj(\param)}\) assuming scale \(\scale=1\) (cf.~Table~\ref{table: optimal step
		size}). A-RFD (Definition~\ref{def: a-rfd}) is plotted as dashed lines.
		A-RFD of the rational quadratic coincides with A-RFD of the squared
		exponential covariance.
	\label{fig: rfd step sizes}
	}
	\centering	
	\def\svgwidth{0.8\linewidth}
	\begingroup \makeatletter \providecommand\color[2][]{\errmessage{(Inkscape) Color is used for the text in Inkscape, but the package 'color.sty' is not loaded}\renewcommand\color[2][]{}}\providecommand\transparent[1]{\errmessage{(Inkscape) Transparency is used (non-zero) for the text in Inkscape, but the package 'transparent.sty' is not loaded}\renewcommand\transparent[1]{}}\providecommand\rotatebox[2]{#2}\newcommand*\fsize{\dimexpr\f@size pt\relax}\newcommand*\lineheight[1]{\fontsize{\fsize}{#1\fsize}\selectfont}\ifx\svgwidth\undefined \setlength{\unitlength}{450bp}\ifx\svgscale\undefined \relax \else \setlength{\unitlength}{\unitlength * \real{\svgscale}}\fi \else \setlength{\unitlength}{\svgwidth}\fi \global\let\svgwidth\undefined \global\let\svgscale\undefined \makeatother \begin{picture}(1,0.66666667)\lineheight{1}\setlength\tabcolsep{0pt}\put(0,0){\includegraphics[width=\unitlength,page=1]{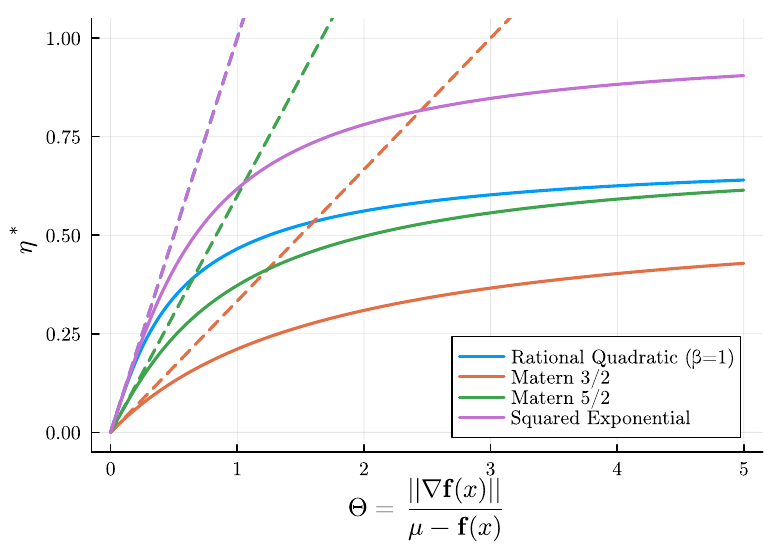}}\end{picture}\endgroup  	\end{sidecaption}
\end{figure}
i.e. \(\stepsize_n^* = \stepsize^*(\Theta_n)\).
This minimization problem can be solved explicitly for the most common
\autocite[ch.~4]{rasmussenGaussianProcessesMachine2006} differentiable isotropic
covariance models, see Table~\ref{table:
optimal step size}, Figure~\ref{fig: rfd step sizes} and Appendix~\ref{appendix:
covariance models} for details.

Figure~\ref{fig: rfd step sizes} can be interpreted as follows: At the start of
optimization, the cost should be roughly equal to the average cost \(\mu \approx
\Obj(\param)\), so the gradient cost quotient \(\Theta\) is infinite and the
step sizes are therefore given by \(\stepsize^*(\infty)\) (also listed
in its own column in Table~\ref{table: optimal step size}). As we start
minimizing, the difference \(\mu-\Obj(\param)\) becomes positive. Towards the
end of minimization this difference no longer changes as the cost no longer
decreases. I.e. towards the end the gradient cost quotient \(\Theta\) is roughly
linear in the gradient \(\|\nabla\Obj(\param)\|\). The derivative
\(\frac{d}{d\Theta}\stepsize^*(0)\) of \(\stepsize^*(\Theta)\) at zero then effectively
results in a constant asymptotic learning rate.

\subsection{Asymptotic learning rate}
\label{subsec: asymptotic RFD step sizes}

To explain the claim above, note that the gradient cost quotient \(\Theta\)
converges to zero towards the end of optimization, because the gradient norm
converges to zero. A first order Taylor expansion of \(\stepsize^*\) would
therefore imply
\[
	\stepsize^*(\Theta)
	\approx \stepsize^*(0) + \tfrac{d}{d\Theta}\stepsize^*(0) \Theta
	= \underbrace{\tfrac{\tfrac{d}{d\Theta}\stepsize^*(0)}{\mu-\Obj(\param)}}_{\text{asymptotic learning rate}} \|\nabla\Obj(\param)\|
\]
assuming \(\stepsize^*(0) = 0\) and differentiability of \(\stepsize^*\), which
is a reasonable educated guess based on the the examples in Figure~\ref{fig: rfd
step sizes}.
But since the RFD step sizes \(\stepsize^*\) are abstractly defined as an
\(\argmin\), it is necessary to formalize this intuition for general covariance
models. First, we define asymptotic step sizes as an object towards which we can
prove convergence. Then we prove convergence, proving they are well defined. In
addition, we obtain a more explicit formula for the asymptotic learning rate.

\begin{definition}[A-RFD]
	\label{def: a-rfd}
	We define the step sizes of ``asymptotic RFD'' (A-RFD) to be the minimizer of the
	second order Taylor approximation \(T_2q_\Theta\) of \(q_\Theta\)
	around zero
	\[
		\hat{\stepsize}(\Theta)
		:= \argmin_\stepsize T_2q_\Theta(\stepsize)
		= \tfrac{\ikernel(0)}{-\ikernel'(0)}\Theta
		= \underbrace{\tfrac{\ikernel(0)}{\ikernel'(0)(\Obj(\param)-\mu)}}_{\text{asymptotic learning rate}}\|\nabla\Obj(\param)\|.
	\]
\end{definition}
In the following we prove that these are truly asymptotically equal to the step
sizes \(\stepsize^*\) of RFD.
\begin{prop}[A-RFD is well defined]
	\label{prop: a-rfd well defined}
	Let \(\Obj\sim\normal(\mu, \ikernel)\) and assume there exists \(\stepsize_0>0\) such
	that the correlation for larger distances \(\stepsize\ge \stepsize_0\) are
	bounded smaller than \(1\), i.e. \(\frac{\ikernel(\stepsize^2/2)}{\ikernel(0)} < \rho \in (0,1)\).
	Then the step sizes of RFD are asymptotically equal to the step sizes of
	A-RFD, i.e.
	\[
		\hat{\stepsize}(\Theta)\sim \stepsize^*(\Theta)
		\quad \text{as}\quad \Theta\to 0.
	\]
\end{prop}

Note that the assumption is essentially always satisfied, since the Cauchy-Schwarz inequality implies
\[
	\ikernel\bigl(\tfrac{\|\param-\tilde{\param}\|^2}2\bigr) = \Cov(\Obj(\param),\Obj(\tilde{\param}))
	\le \sqrt{\Var(\Obj(\param))\Var(\Obj(\tilde{\param}))} = \ikernel(0),
\]
where equality requires the random variables to be almost surely equal
\autocite{klenkeProbabilityTheoryComprehensive2014}. If the random function is
not periodic or constant, this will generally be strict.
In the proof, this requirement is only needed to ensure that \(\stepsize^*\) is not
very large. The smallest local minimum of \(q_\Theta\) is always close to
\(\hat{\stepsize}\) even without this assumption (which ensures it is a global minimum).

Figure~\ref{fig: rfd step sizes} illustrates that \(\stepsize^*\to 0\) should
imply \(\Theta\to 0\), resulting in a weak convergence guarantee.

\begin{restatable}{corollary}{convergence}\label{prop: convergence}
	Assume \(\stepsize^*\to 0\) implies \(\Theta\to 0\), the cost \(\Obj\) is
	bounded, has continuous gradients and RFD converges to some point
	\(\param_\infty\). Then \(\param_\infty\) is a critical point and
	the RFD step sizes \(\stepsize^*\) are asymptotically equal to
	\(\hat{\stepsize}\).
\end{restatable}

For the squared exponential covariance model we formally prove that \(\stepsize^*\)
is strictly monotonously increasing in \(\Theta\) and thus \(\stepsize^*\to 0\)
implies \(\Theta\to 0\) (Prop.\ \ref{prop: sq exp is
strictly monotonous in xi}). The `bounded' and `continuous gradients'
assumptions are almost surely satisfied for all sufficiently smooth covariance
functions \autocite[cf.][]{adlerRandomFieldsGeometry2007}, where three times
differentiable is more than enough smoothness.

\subsection{RFD step sizes explain common step size heuristics}
\label{subsec: step size heuristics}

Asymptotically, RFD suggests constant learning rates, similar to the classical
\(L\)-smooth setting. We thus define these asymptotic learning rates (as the limit
of the learning rates \(\lr_n\) of iteration \(n\)) to be
\begin{equation}
	\label{eq: asymptotic learning rate}
	\lr_\infty := \frac{\ikernel(0)}{\ikernel'(0)(\Obj(\param_\infty) - \mu)},
\end{equation}
where \(\Obj(\param_\infty)\) is the cost we reach in the limit. If we used
these asymptotic learning rates from the start, step sizes would become too
large for large gradients, as RFD step sizes exhibit a plateau
(cf.~Figure~\ref{fig: rfd step sizes}). To emulate the behavior of RFD with a
piecewise linear function, we could introduce a cutoff whenever our step size
exceeds the initial step size
\(\stepsize^*(\infty)\), i.e.
\[
	\tag{gradient clipping}
	\param_{n+1}
	= \param_n - \min\Bigl\{
		\lr_\infty, \frac{\stepsize^*(\infty)}{\|\nabla\Obj(\param_n)\|}
	\Bigr\}\nabla\Obj(\param_n).
\]
At this point we have rediscovered `gradient clipping'
\autocite{pascanuDifficultyTrainingRecurrent2013}. As the rational quadratic
covariance has the same asymptotic learning rate \(\lr_\infty\) for every
\(\beta\), its parameter \(\beta\) controls the step size bound
\(\stepsize^*(\infty)\) of gradient clipping (cf.~Table~\ref{table: optimal step
size}, Figure~\ref{fig: rfd step sizes}).

\citet{pascanuDifficultyTrainingRecurrent2013} motivated gradient clipping with the
geometric interpretation of movement towards a `wall' placed behind the minimum.
This suggests that clipping should happen towards the end of training. This
stands in contrast to a more recent
step size heuristic, ``(linear) warmup''
\autocite{goyalAccurateLargeMinibatch2018}, which suggests smaller learning
rates at the start (i.e.
\(\lr_0 = \frac{\stepsize^*(\infty)}{\|\nabla\Obj(\param_0)\|}\)) and gradual ramp-up
to the asymptotic learning rate \(\lr_\infty\). In other words, gradients are not clipped
due to some wall next to the minimum, but because the step sizes would be too
large at the start otherwise. \citet{goyalAccurateLargeMinibatch2018} further observe that
`constant warmup' (i.e. a step change of learning rates akin to gradient clipping) performs
worse than gradual warmup. Since RFD step sizes suggest this gradual increase,
we argue that they may have discovered RFD step sizes empirically (also cf.~Figure~\ref{fig:
summary figure}). \section{Mini-batch loss and covariance estimation}
\label{sec: stochastic loss and covariance estimation}

In practice we typically do not have access to evaluations of the cost \(\Obj\)
but rather access to noisy evaluations 
\(\obs_n = \jet\Obs_n(X_n)\) with  \(\Obs_n(\param) =
\Obj(\param)+\noise_n(\param)\), where we typically assume
the noise \(\noise_n\) is centered and \(\iid\) conditional on \(\Obj\).\footnote{
    More precisely iid conditional on some \(\Problem\) such hat \(\Obj\) is
    measurable w.r.t.\ \(\Problem\), see Definition \ref{def: feasible
    optimization algorithm}
}
But at the cost of multiple evaluations the variance of the noise can be
reduced. Specifically, let \(\batchsize_n\) be the number of
evaluations (the `mini-batch size') used to produce the observation \(\Obs_n\),
then
\begin{equation}
    \label{eq: mini batch loss}
    \Obs_n(\param)
    = \frac1{\batchsize_n}\sum_{i=1}^{\batchsize_n} \Obj(\param)
    + \noise_n^{(i)}(\param)
    = \Obj(\param) + \underbrace{\frac1{\batchsize_n}\sum_{i=1}^{\batchsize_n}\noise_n^{(i)}(\param)}_{=:\noise_n(x)},
\end{equation}
where we assume the noise \(\noise_n^{(i)}\) is \(\iid\) and centered conditional
on \(\Obj\). This implies
\begin{equation}
    \label{eq: variance of mini batch}
    \Var(\Obs_n(\param))
    = \Var(\Obj(\param)) + \tfrac1{\batchsize_n}\Var(\noise_n^{(1)}(\param))
    \overset{\text{isotropy}}=
    \ikernel(0) + \tfrac1{\batchsize_n}\ikernel_\noise(0),
\end{equation}
where we assume \(\Obj\sim\normal(\mu, \ikernel)\) and
\(\noise_n^{(i)}\sim\normal(0, \ikernel_\noise)\) in the last equation for
simplicity.\footnote{Although this step did not yet require the distributional Gaussian assumption
beyond the mean and variance.} We therefore have control over the variance of
\(\noise_n\). This will help us estimate both \(\ikernel\) and
\(\ikernel_\noise\) from evaluations.

\subsection{Variance estimation}
\label{subsec: non-parametric covariance estimation}

Recall that the asymptotic learning rate \(\lr_\infty\) in Equation~\eqref{eq:
asymptotic learning rate} only depends on \(\ikernel(0)\)
and \(\ikernel'(0)\). So if we estimate these values, we are certain to get
the right RFD step sizes asymptotically without knowing the entire covariance
kernel \(\ikernel\).

Equation~\eqref{eq: variance of mini batch} reveals that
for \(Z_n := (\Obs_n(\param) - \mu)^2\) we have
\[
    \E[Z_n] = \beta_0 + \tfrac1{\batchsize_n} \beta_1
    \qquad\text{i.e.}\qquad
    Z_\batchsize = \beta_0 + \tfrac1{\batchsize_n} \beta_1 + \text{noise}
\]
with bias \(\beta_0=\ikernel(0)\) and slope
\(\beta_1=\ikernel_\epsilon(0)\). So a linear regression on samples
\((\tfrac1{\batchsize_k}, Z_k)_{k\le n}\) allows for the estimation of
\(\beta_0\) and \(\beta_1\).

As the noise functions \(\noise_k^{(i)}\) are all conditionally independent and
therefore uncorrelated, we have
\[
	\Cov(\Obs_k(\param_k), \Obs_l(\param_l))
	= \Cov(\Obj(\param_k), \Obj(\param_l))
	= \ikernel\bigl(\tfrac{\|\param_k - \param_l\|^2}2\bigr)
\]
Since the covariance kernel \(\ikernel\) is typically monotonously falling in the
distance of parameters \(\|\param_k - \param_l\|^2\), we want to select them as spaced
out as possible to minimize the covariance of \(\Obs_k(\param_k)\) to minimize
the bias. Randomly selecting \(\param_i\) with Glorot initialization
\autocite{glorotUnderstandingDifficultyTraining2010} ensures a good
spread.\footnote{
    Note that Glorot initialization places all parameters approximately on the same
    sphere. This is because Glorot initialization initializes all parameters independently,
    therefore their norm is the sum of independent squares, which converges by the
    law of large numbers due to the normalization Glorot uses.
    Since stationary isotropy and non-stationary isotropy coincides on the sphere,
    this is an important effect to consider (cf.~Section~\ref{sec: input invariance}).
}

By the Gaussian assumption the variance of \(Z_n\) is the (centered)
fourth moment of \(\Obs_n\), which is
given by
\[
    \sigma_{\batchsize_n}^2 := \Var(Z_n)
    = \E[Z_n^4] - \E[Z_n^2]^2
    =  2 \Var(\Obs_n(\param))^2
    = 2(\beta_0 + \tfrac1{\batchsize_n} \beta_1)^2.
\]
In particular the variance of \(Z_n\) depends on the batch size \(\batchsize_n\).
The linear regression is therefore heteroskedastic. Weighted least squares (WLS) 
\autocite[e.g.][Theorem~4.2]{kayFundamentalsStatisticalSignal1993} is designed to handle
this case, but for its application the variance of \(Z_n\) is needed. Since
\(\beta_0,\beta_1\) are the parameters we wish to estimate, we find ourselves in
the paradoxical situation that we need \(\beta\) to obtain \(\beta\). Our solution to
this problem is to start with a guess of \(\beta_0, \beta_1\), apply WLS to
obtain a better estimate and repeat this bootstrapping procedure until
convergence.

The same procedure can be applied to obtain \(\ikernel'(0)\), where the counterpart of
Equation~\eqref{eq: variance of mini batch} is given by
\[
    \Var(\partial_i \Obs_n(\param))
    = \Var(\partial_i\Obj(\param)) + \tfrac1{\batchsize_n}\Var(\partial_i\noise_n^{(1)}(\param))
    \overset{\text{isotropy}}=
    -(\ikernel'(0) + \tfrac1{\batchsize_n}\ikernel_\noise'(0)).
\]
\begin{remark}
Under the isotropy assumption the partial derivatives are iid, so the expectation of
\(\|\nabla\Obs_n(\param)\|^2
= \sum_{i=1}^\dims (\partial_i\Obs_n(\param))^2\)
is this variance scaled by \(\dims\). In particular the variance needs to
scale with \(\frac1\dims\) to keep the gradient norms (and thus the Lipschitz
constant of \(\Obj\)) stable. This observation is closely related to
``isoperimetry'' \autocite[e.g.][]{bubeckUniversalLawRobustness2021}, for details
see \cite{benningGradientSpanAlgorithms2024}. Removing the isotropy
assumption and estimating the variance component-wise is most likely how
``adaptive'' step sizes
\autocite[e.g.][]{duchiAdaptiveSubgradientMethods2011,kingmaAdamMethodStochastic2015},
like the ones used by Adam, work (cf.~Sec.~\ref{sec: geometric anisotropy}).
\end{remark}

\paragraph*{Batch size distribution}

Before we can apply linear regression to the samples  \((\tfrac1{\batchsize_k},
Z_k)_{k\le n}\), it is necessary to choose the batch sizes
\(\batchsize_k\). As this choice is left to us, we
calculate the variance of our estimator \(\hat{\beta}_0\) of \(\beta_0\)
explicitly (Lemma~\ref{lem: variance of beta_0}), in order to minimize this
variance subject to a sample budget \(\alpha\) over the selection of batch sizes
\begin{equation}
    \label{eq: bsize optimization problem}
    \min_{n, \batchsize_1,\dots, \batchsize_n}\Var(\hat{\beta_0}) 
    \quad\text{s.t.}\quad \underbrace{\sum_{k=1}^n \batchsize_k}_{\text{samples used}}\le \alpha.
\end{equation}
Since this optimization problem is very difficult to solve, we rephrase it in terms of
the empirical distribution of batch sizes \(\nu_n = \frac1n\sum_{i=1}^n
\delta_{\batchsize_i}\). Optimizing over distributions is still difficult, but we
explain in Section~\ref{sec: batch size distribution} how to heuristically arrive at the parametrization
\[
    \nu(b) \propto
    \exp\bigl(\lambda_1 \tfrac1{\sigma_\batchsize^2} - \lambda_2 \batchsize\bigr),
    \qquad b\in \nat
\]
where the parameters \(\lambda_1,\lambda_2\ge 0\) can then be used to optimize
\eqref{eq: bsize optimization problem}. Due to our usage of
\(\sigma_\batchsize^2\) this has to be bootstrapped.

\paragraph*{Covariance estimation} While the variance estimates above ensure correct
asymptotic learning rates, we motivated in Section~\ref{subsec: step size heuristics} that
asymptotic learning rates alone would result in too large step sizes at the
beginning. We therefore use the estimates of \(\ikernel(0)\) and \(\ikernel'(0)\)
to fit a covariance model, effectively acting as a gradient clipper while
retaining the asymptotic guarantees. Note that covariance models with less than
two parameters are generally fully determined by these values.

\subsection{Stochastic RFD (S-RFD)}

It is reasonable to ask whether there is a `stochastic gradient descent'-like
counterpart to the `gradient descent'-like RFD. The answer is yes, and we already
have all the required machinery. 

\begin{restatable}[S-RFD]{extension}{srfd}
    \label{ext: s-rfd}
    For loss \(\Obj\sim\normal(\mu, \ikernel)\) and stochastic errors
    \(\noise_i\overset{\iid}\sim\normal(0, \ikernel_\noise)\) we have 
    \[
        \argmin_{\step}\E[\Obj(\param - \step) \mid \Obs_n(\param), \nabla \Obs_n(\param)]
        = \stepsize^*(\Theta)\tfrac{\nabla\Obs_n(\param)}{\|\nabla\Obs_n(\param)\|}
    \]
    with the same step size function \(\stepsize^*\) as for RFD, but modified \(\Theta\)
    \[
        \Theta =
        \frac{
            \ikernel'(0)
        }{
            \ikernel'(0)+ \tfrac1{\batchsize_n}\ikernel_\noise'(0)
        }
        \frac{
            \ikernel(0)+\tfrac1{\batchsize_n}\ikernel_\noise(0)
        }{\ikernel(0)}
        \frac{\|\nabla\Obs_n(\param)\|}{\mu-\Obs_n(\param)}.
    \]
\end{restatable}
Note, that our non-parametric covariance estimation already provides us with estimates
of \(\ikernel_\noise(0)\) and \(\ikernel_\noise'(0)\), so no further adaptions are
needed. The resulting asymptotic learning rate is given by
\begin{equation}
    \label{eq: asymptotic lr s-rfd}
    \lr_\infty = \frac{\ikernel(0)+\tfrac1{\batchsize_n}\ikernel_\noise(0)}{
        (\ikernel'(0)+ \tfrac1{\batchsize_n}\ikernel_\noise'(0))(\Obs_n(\param_\infty) - \mu)
    }.
\end{equation}
 \section{MNIST case study}\label{sec: mnist case study}

For our case study we use the negative log likelihood loss to train a neural network
{
	\DeclareFieldFormat{postnote}{#1}
	\autocite[M7]{anEnsembleSimpleConvolutional2020}
}
on the MNIST dataset
\autocite{lecunMNISTDATABASEHandwritten2010}. We choose this model as one of the
simplest state-of-the-art models at the time of selection, consisting only of
convolutional layers with ReLU activation interspersed by batch
normalization layers and a single dense layers at the end with softmax
activation.
Assuming isotropy, we estimate \(\mu\), \(\ikernel(0)\) and \(\ikernel'(0)\) as
described in Section~\ref{subsec: non-parametric covariance estimation} and
deduce the parameters \(\sigma^2\) and \(\scale\) of the respective covariance
model. We then use the step sizes listed in Table~\ref{table: optimal step size}
for the `squared exponential' and `rational quadratic' covariance in our RFD
algorithm.

In Figure~\ref{fig: summary figure}, RFD is benchmarked against step size tuned Adam
\autocite{kingmaAdamMethodStochastic2015} and stochastic gradient descent (SGD). 
Even with early stopping, their tuning would typically require more than 1 epoch
worth of samples, \emph{in contrast to RFD} (Section~\ref{sec: sampling
efficiency and stability}). We
highlight that A-RFD performs significantly worse than either of the RFD versions
which effectively implement some form of learning rate warmup. This is despite the
RFD learning rates converging to the asymptotic one within one epoch (ca. \(30\) out
of \(60\) steps per epoch). The step sizes on the other hand are (up to noise)
monotonously decreasing. This stands in contrast to the ``wall next to the
minimum'' motivation of gradient clipping.

\begin{figure}
\begin{sidecaption}{
	Training on the MNIST dataset (batch size \(1024\)). Ribbons describe the
	range between the \(10\%\) and \(90\%\) quantile of \(20\) repeated
	experiments while lines represent their mean. SE stands for the squared exponential \eqref{eq: sqExp
	covariance model} and RQ for the rational quadratic \eqref{eq:
	rational quadratic} covariance. The validation loss uses the test data set, 
	which provides a small advantage to Adam and SGD, as we also use it for tuning.
	\label{fig: summary figure}
}
	\centering
	\includegraphics[width=\linewidth]{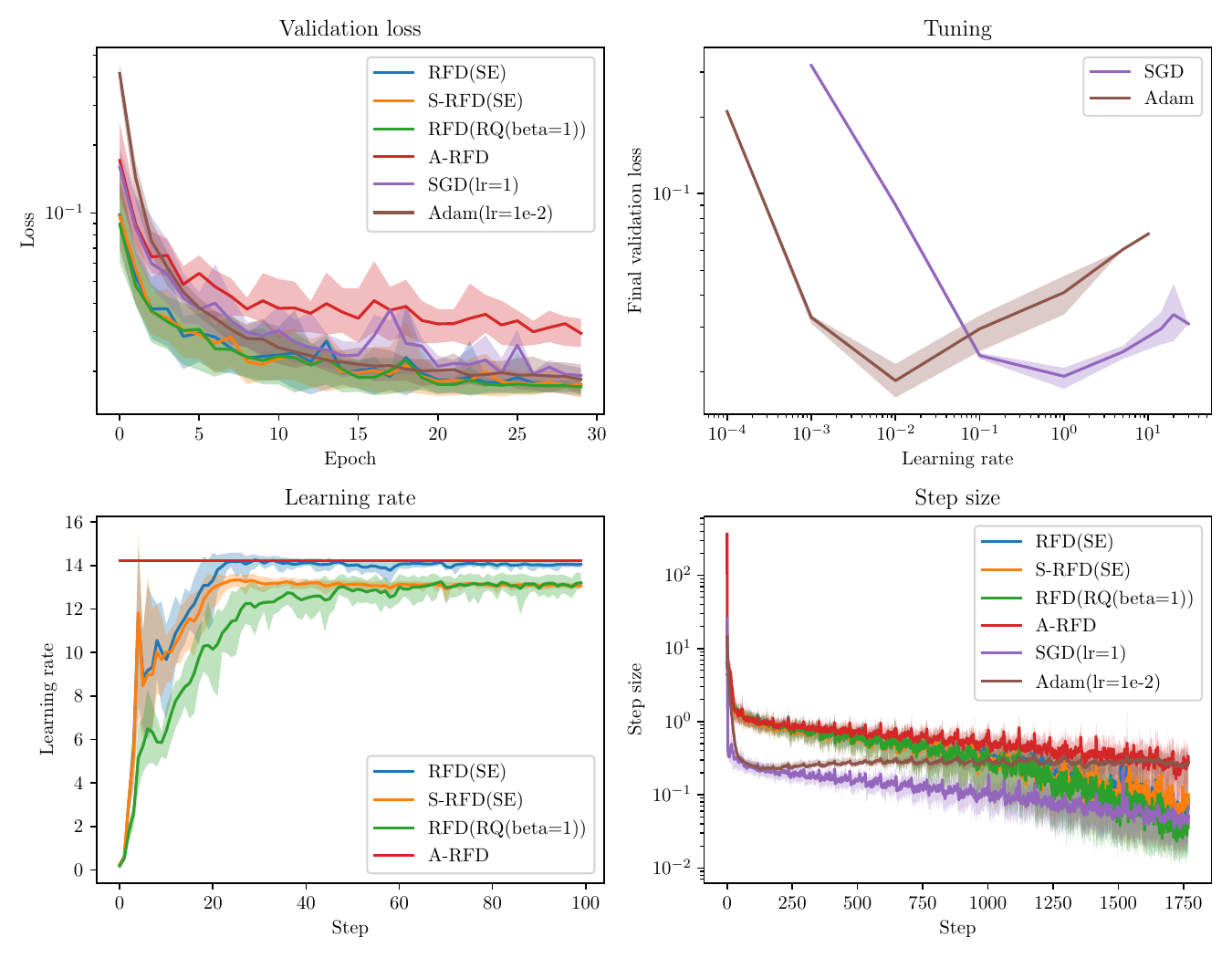}
\end{sidecaption}
\end{figure}

\textbf{Code availability:}
Our implementation of RFD can be found at \url{https://github.com/FelixBenning/pyrfd} and
the package can also be installed from \href{https://pypi.org/project/pyrfd/}{PyPI} via
`\verb|pip install pyrfd|'. \section{Limitations and extensions}
\label{sec: extensions and limitations}

\clearmargin
To cover the vast amount of ground that lays between the `formulation of a
general average case optimization problem' and the `prototype of a working
optimizer with theoretical backing',
\begin{enumerate}[leftmargin=1cm]
	\item we used the common 
	\autocite[][]{frazierBayesianOptimization2018,wangBayesianOptimizationBillion2016,wuBayesianOptimizationGradients2017,dauphinIdentifyingAttackingSaddle2014}
	\emph{isotropic} and \emph{Gaussian} distributional assumption for \(\Obj\),
	\item we used very \emph{simple covariance models} for the actual implementation,
	\item we used WLS in our variance estimation procedure despite the \emph{violation of independence}.
\end{enumerate}
Since RFD is defined as the minimizer of an average instead of an upper bound -- 
making it more risk affine --
it naturally loses the improvement guarantee driving classical convergence
proofs. It is therefore impossible to extend classical optimization proofs and new
mathematical theory must be developed. This risk-affinity can also be observed
in its comparatively large step sizes (cf.~Fig.~\ref{fig: summary
figure}
and Sec.~\ref{sec: experiments}). On CIFAR-100
\autocite{krizhevskyLearningMultipleLayers2009}, the step sizes were \emph{too}
large and it is an open question whether assumptions were violated or whether
RFD is simply too risk-affine. But since the variance of random functions
vanishes asymptotic with high dimension\footnote{see Chapter \ref{chap: predictable progres}}
we highly suspect the former, see Remark~\ref{rem: high dimension} for further details.

Future work will therefore have to target these assumptions. Some of the assumptions
were already simplifications for the sake of exposition, and we deferred their
relaxation to the appendix. The Gaussian assumption can be relaxed with a
generalization to the `BLUE' (Sec.~\ref{sec: BlUE}),
isotropy can be generalized to `geometric anisotropies' (Sec.~\ref{sec:
geometric anisotropy}) and the risk-affinity of RFD can be reduced with
confidence intervals (Sec.~\ref{sec: conservative rfd}). Since simple random
linear models already violate stationary isotropy (Sec.~\ref{sec: random linear model}),
we believe that stationarity is the most important assumption to attack in future
work.

 \section{Conclusion}

In this paper we have demonstrated the \textbf{viability} (computability and
scalable complexity) and \textbf{advantages} (scale invariance, explainable step
size schedule which does not require expensive tuning) of replacing the
classical ``convex function'' framework with the ``random function'' framework.
Along the way we
bridged the gap between Bayesian optimization (not scalable so far) and
classical optimization methods (scalable). This theoretical framework not only
sheds light on existing step size heuristics, but can also be used to develop
future heuristics.

We envision the following improvements to RFD in the future:

\begin{enumerate}
	\item The \emph{reliability} of RFD can be improved by generalizing the
	distributional assumptions to cover more real world scenarios. In particular
	we are interested in the generalization to non-stationary isotropy because we
	suspect that regularization such as weight and batch normalization
	\autocite{salimansWeightNormalizationSimple2016,ioffeBatchNormalizationAccelerating2015}
	are used to patch violations of stationarity (cf. Section~\ref{sec: random linear model}).

	\item The \emph{performance} of RFD can also be improved. Since RFD is forgetful
	while momentum methods retains some information it is likely fruitful to
	relax the full forgetfulness. Furthermore, we suspect that adaptive learning
	rates \autocite[e.g.][]{duchiAdaptiveSubgradientMethods2011,kingmaAdamMethodStochastic2015}, such as those used by Adam, can be incorporated with
	geometric anisotropies (cf. Sec. E.1). Performance could also be further
	improved by estimating the covariance (locally) online instead of globally at
	the start. Finally, the implementation itself can be made more performant.
\end{enumerate}

\section*{Acknowledgement}
We extend our sincere gratitude to our colleagues at the University of
Mannheim, with special thanks to Rainer Gemulla and Julie Naegelen for
insightful discussions and invaluable feedback. 
The Experiments in this work were partially carried out on the compute cluster of
the state of Baden-Würtemberg (bwHPC).

\printbibliography[heading=subbibliography]

\clearpage

\begin{subappendices}
	\begin{center}
		\textbf{\Large Appendix: Random Function Descent}
	\end{center}

\section{Experiments}
\label{sec: experiments}

\subsection{Covariance estimation}

\begin{figure}
	\begin{fullwidth}
	\includegraphics[width=\linewidth]{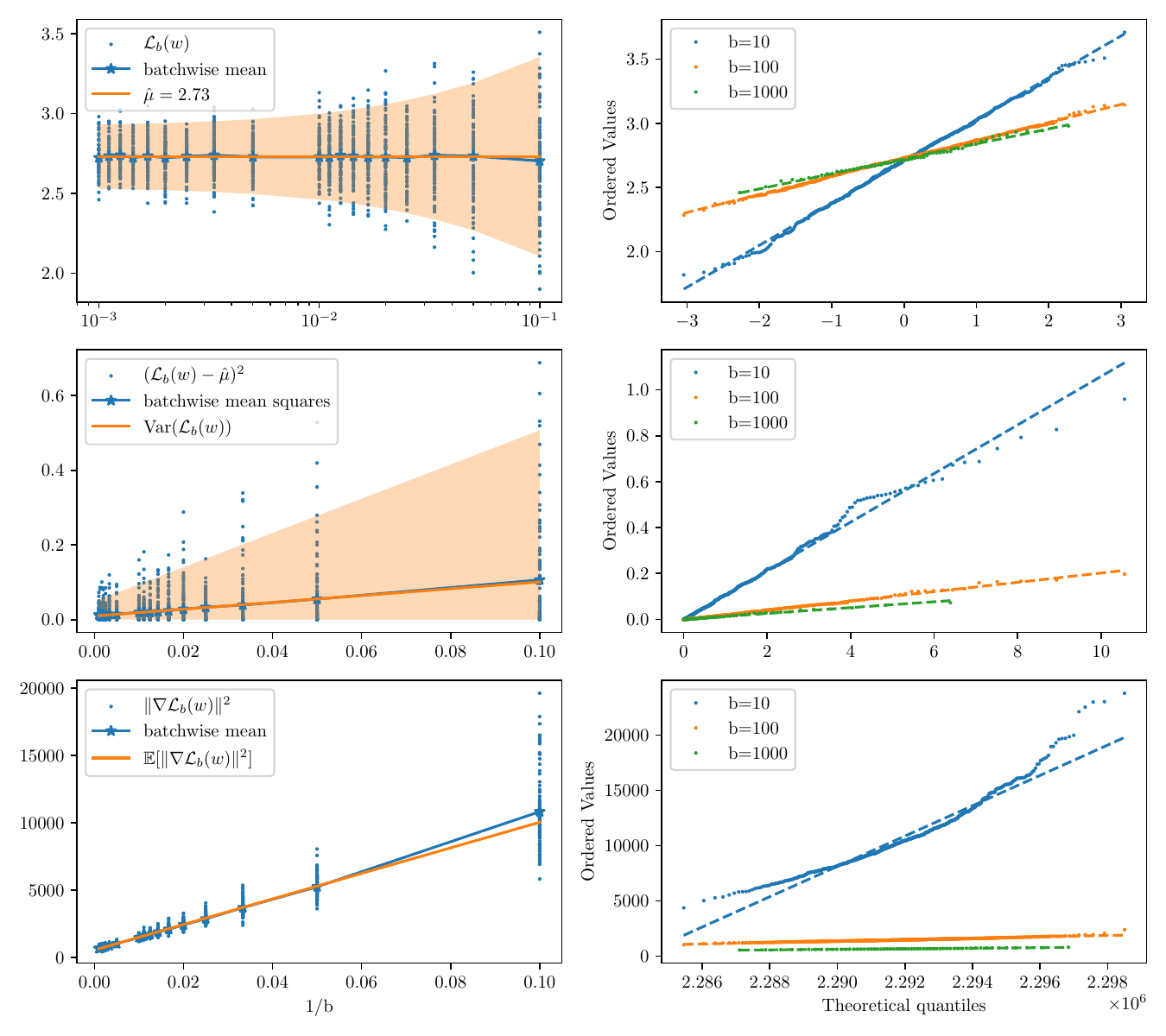}
	\end{fullwidth}
	\fullwidthcaption{
		Visualization of the variance estimation (Section \ref{subsec:
		non-parametric covariance estimation}) with \(95\%\)-confidence intervals
		based on the assumed distribution. Here, \(\Loss_\batchsize\) are the
		mini-batch losses of batch-size \(\batchsize\) that we
		more generally denote as the observations
		\(\Obs_n\) with batch-size \(\batchsize_n\).
		Quantile-quantile (QQ) plots of the
		losses (against a normal distribution), squared losses (against a
		\(\chi^2(1)\) distribution) and squared gradient norms (against a
		\(\chi^2(\dims)\)-distribution) are displayed on the right for a selection
		of batch sizes.
		\label{fig: covariance estimation}
	}
\end{figure}

In Figure~\ref{fig: covariance estimation} we visualize weighted
least squares (WLS) regression of the covariance estimation from
Section~\ref{subsec: non-parametric covariance estimation}.
Note, that we sampled much more samples per batch size for these plots than RFD
would typically require by itself in order to be able to plot batch-wise means and
batch-wise QQ-plots. The batch size distribution we described in
Section~\ref{sec: batch size distribution} would avoid sampling the
same batch size multiple times to ensure better stability of the regression and
generally requires much fewer samples than were used for this visualization
(cf.~\ref{sec: sampling efficiency and stability})

We can observe from the QQ-plots on the right, that the Gaussian assumption is
essentially justified for the losses, resulting in a \(\chi^2(1)\) distribution
for the squared losses and a \(\chi^2(\dims)\) distribution for the gradient
norms squared. The confidence interval estimate for the squared norms appears to
be much too small (it is plotted, but too small to be visible). Perhaps this
signifies a violation of the isotropy assumption as the variance of
\[
	\|\nabla\Obs_n(\param)\|^2
	= \sum_{i=1}^\dims (\partial_i \Obs_n(\param))^2
\]
does not appear to be the variance of independent \(\chi^2(\dims)\) Gaussian
random variables, and the independence only follows from the isotropy
assumption.

\subsubsection{Sampling efficiency and stability}
\label{sec: sampling efficiency and stability}

To evaluate the sampling efficiency and stability of our variance estimation process,
we repeated the covariance estimation of the model model M7
\parencite{anEnsembleSimpleConvolutional2020} applied to the
MNIST dataset \(20\) times (Figure~\ref{fig: sampling efficiency and stability}).
We used a tolerance of \(\mathrm{tol}=0.3\) as a stopping criterion for
the estimated relative standard deviation \eqref{eq: relative std}.

At this tolerance, the asymptotic learning rate already seems relatively stable (in the
same order of magnitude) and the sample cost is quite cheap. The majority of
runs (\(16/20\) runs or \(80\%\)) required less than \(60\,000\) samples (1
epoch). There was one large outlier which used \(500\,589\) samples. A closer
inspection revealed, that after the initial sample to estimate the optimal batch size
distribution, it sampled almost exclusively at batch sizes \(20\) (which was the minimal
cutoff to avoid instabilities caused by batch normalization) and batch sizes between
\(1700\) and \(1900\). It therefore seems like the initial batch of samples caused
a very unfavorable batch size distribution which then required a lot of samples to
recover from. Our selection of an initial sample size of \(6000\) might
therefore have been too small.

A more extensive empirical study is needed to tune this estimation process, but the
process promises to be very sample efficient. Classical step size tuning would
train models for a short duration in order to evaluate the performance of a particular
learning rate \parencite[e.g.][]{smithDisciplinedApproachNeural2018}, but a single
epoch worth of samples is very hard to beat.

Our implementation of this process on the other hand is very inefficient as of writing.
Piping data of differing batch sizes into a model is not a standard use case. We
implement this by repeatedly initializing data loaders, which is anything but
performance friendly.

\begin{figure}
	\begin{sidecaption}{
		\(20\) repeated covariance estimations of model M7
		\parencite{anEnsembleSimpleConvolutional2020} applied to the
		MNIST dataset. On the left are the resulting asymptotic learning rates
		(assuming a final loss of zero) and on the right are the
		samples used until the stopping criterion interrupted sampling. 
	\label{fig: sampling efficiency and stability}
	}
	\includegraphics[width=0.49\linewidth]{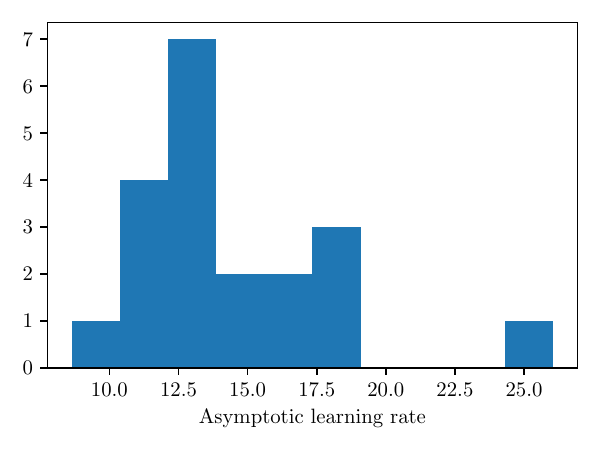}
	\includegraphics[width=0.49\linewidth]{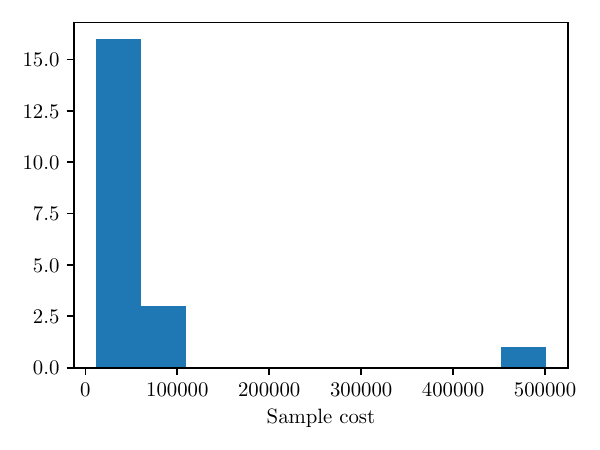}
	\end{sidecaption}
\end{figure}

\subsection{Other models and datasets}

To estimate the effect of the batch size on RFD, we trained the same model (M7
\parencite{anEnsembleSimpleConvolutional2020}) on MNIST with batch size \(128\)
(Figure~\ref{fig: mnist cnn7 bsize=128}). We can see that the asymptotic learning
rate of S-RFD is reduced at a smaller batch size (cf.~Equation~\ref{eq:
asymptotic lr s-rfd}) but the performance is barely different. Overall, RFD
seems to be slightly too risk-affine, selecting larger step sizes than
the tuned SGD models.

\begin{figure}
	\begin{fullwidth}
	\includegraphics[width=\linewidth]{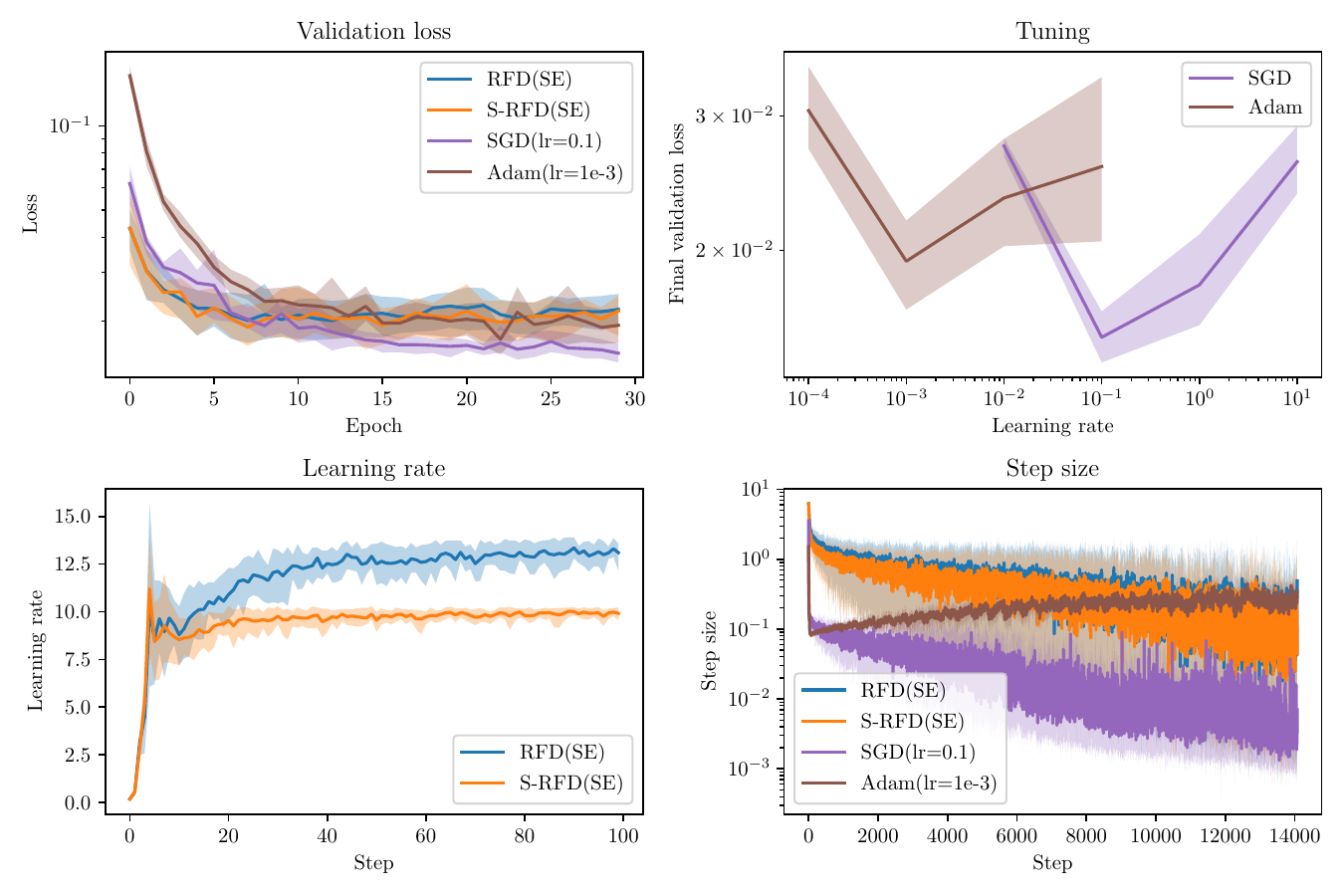}
	\end{fullwidth}
	\fullwidthcaption{
		Training model M7 \parencite{anEnsembleSimpleConvolutional2020} with batch
		size \(128\) on MNIST \parencite{lecunMNISTDATABASEHandwritten2010}.
		\label{fig: mnist cnn7 bsize=128}
	}
\end{figure}

We also trained a different model (M5 \parencite{anEnsembleSimpleConvolutional2020})
on the Fashion MNIST dataset \parencite{xiaoFashionMNISTNovelImage2017}
with batch size \(128\) (Figure~\ref{fig: fashion mnist}). Since the validation loss increases after
epoch \(5\), early stopping would have been appropriate. We therefore include
Adam with learning rate \(10^{-3}\), despite Adam with learning rate \(10^{-4}\)
technically performing better at the end of training. We can generally see, that
RFD comes very close to tuned performance at the time early stopping would have been
appropriate. Again, learning rates seem to be slightly too large (risk-affine)
in comparison to tuned SGD.

\begin{figure}
	\begin{fullwidth}
		\includegraphics[width=\linewidth]{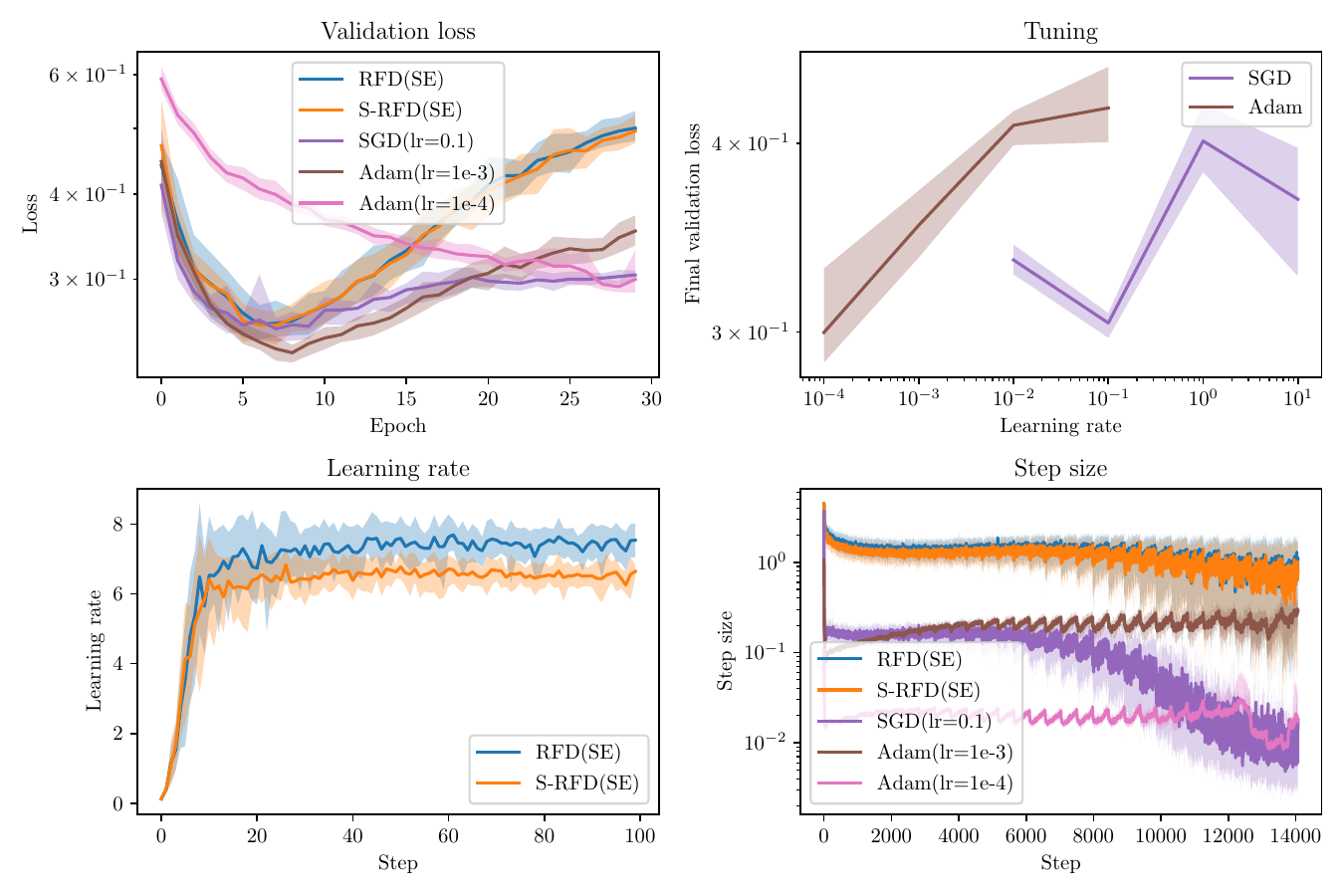}
	\end{fullwidth}
	\fullwidthcaption{
		Model M5 \parencite{anEnsembleSimpleConvolutional2020} trained on Fashion
		MNIST \parencite{xiaoFashionMNISTNovelImage2017} with batch size \(128\).
		\label{fig: fashion mnist}
	}	
\end{figure}

\section{Batch size distribution}
\label{sec: batch size distribution}

Since we plan to use the data set \((\frac1{\batchsize_k}, \Obs_k(\param_k))_{k\le n}\)
for weighted least squares (WLS) regression and do not have a selection process for
the batch sizes \(\batchsize_k\) yet, it might be appropriate to select the batch
sizes \(\batchsize_k\) in such a way, that the variance of our estimator \(\hat{\beta}_0\)
of \(\beta_0\) is minimized. Here we choose \(\Var(\hat{\beta}_0)\) and not
\(\Var(\hat{\beta}_1)\) as our optimization target, since \(\beta_0 =
\ikernel(0)\) is used to fit the covariance model, while \(\beta_1 =
\ikernel_\epsilon(0)\) is only required for S-RFD. Without deeper analysis
\(\beta_0\) therefore seems to be more important.

Optimization over \(n\) parameters \(\batchsize_k\) is quite difficult, but we can
simplify this optimization problem by considering the empirical batch size distribution
\[
	\nu_n = \frac1n\sum_{k=1}^n \delta_{\batchsize_k}.
\]
Using a random variable \(\Batchsize\) distributed according to \(\nu_n\), the
total number of sample losses can then be expressed as 
\[
	\sum_{k=1}^n \batchsize_k	= n \E[B] = \text{samples used}.
\]
Under an (unrealistic) independence assumption, the variance
\(\Var(\hat{\beta}_0)\) also has a simple representation in terms of \(\nu_n\)
(Lemma~\ref{lem: variance of beta_0}).  We now want to minimize this variance
subject to compute constraint \(\alpha\)
limiting the number of sample losses we can use resulting in the optimization problem
\begin{equation}
    \label{eq: batch size distribution opt problem}
    \Var(\hat{\beta_0}) = \frac1n 
    \underbrace{\frac{1}{\E[\frac1{\sigma_\Batchsize^2}]}}_{\text{`variance of \(Z_\Batchsize\)'}}
    \underbrace{
        \frac{\E[\frac1{\sigma_\Batchsize^2 \Batchsize^2}]}{\E\Bigl[\frac1{\sigma_\Batchsize^2}
        \bigl(\frac1\Batchsize - \E[\frac1{\Batchsize\sigma_\Batchsize^2 \E[1/\sigma_\Batchsize^2]}]\bigr)^2\Bigr]}
    }_{\text{inverse of the `spread'  of \(\frac1\Batchsize\)}}
    \quad\text{s.t.}\quad n\E[\Batchsize]\le \alpha.
\end{equation}
where we recall that \(\sigma_\Batchsize^2\) is the variance of
\(Z_\Batchsize=Z_n = (\Obs_n(x) - \mu)^2\) if \(\Batchsize=\batchsize_n\). 
Intuitively, the first therm is therefore the average variance of \(Z_\Batchsize\). The second half is the fraction of a weighted
second moment divided by the weighted
variance. Unless the mean is at zero, the former will be larger. In particular
we want a spread of data otherwise the variance would be zero. This is in some
conflict with the variance of \(Z_\Batchsize\).

But first, let us get rid of \(n\). Note that we would always increase \(n\)
until our compute budget is used up, since this always reduces variance. So we
approximately have \(n\E[\Batchsize] = \alpha\). Thus
\[
    \Var(\hat{\beta_0})
	 = \frac{\E[\Batchsize]}{\alpha}
    \frac{1}{\E[\frac1{\sigma_\Batchsize^2}]}
        \frac{\E[\frac1{\sigma_\Batchsize^2 \Batchsize^2}]}{\E\Bigl[\frac1{\sigma_\Batchsize^2}
        \bigl(\frac1\Batchsize - \E[\frac1{\Batchsize\sigma_\Batchsize^2 \E[1/\sigma_\Batchsize^2]}]\bigr)^2\Bigr]}
\]
Since \(\alpha\) is now just resulting in a constant factor, it can be assumed
to be \(1\) without loss of generality. Over batch size distributions \(\nu\) we
therefore want to solve the minimization problem
\begin{equation}
	\label{eq: the optimization target}
	\min_{\nu}	
	\underbrace{
   \frac{\E[\Batchsize]}{\E[\frac1{\sigma_\Batchsize^2}]}
	}_{\text{moments}}
	\underbrace{
		\frac{\E[\frac1{\sigma_\Batchsize^2 \Batchsize^2}]}{\E\Bigl[\frac1{\sigma_\Batchsize^2}
		\bigl(\frac1\Batchsize - \E[\frac1{\Batchsize\sigma_\Batchsize^2 \E[1/\sigma_\Batchsize^2]}]\bigr)^2\Bigr]}
	}_{\text{spread}}
\end{equation}

\begin{example}[If we did not require spread]
	If we were not concerned with the variance of batch sizes, we could select a
	constant \(\Batchsize=\batchsize\).
	Then it is straightforward to minimize the moments factor manually
	\[
		\min_{\batchsize} \frac{\E[\batchsize]}{\E[\frac1{\sigma_{\batchsize}^2}]}
		= \batchsize \sigma_\batchsize^2
		= 2\batchsize\bigl(\beta_0 + \tfrac1\batchsize \beta_1\bigr)^2,
	\]
	resulting in \(\frac1\batchsize = \frac{\beta_0}{\beta_1}\).
	In other words: If we did not have to be concerned with the spread of \(\Batchsize\)
	there is one optimal selection to minimize the first factor. But in reality we
	have to trade-off this target with the spread of \(\Batchsize\).
\end{example}

To ensure a good spread of data, we use the maximum entropy distribution for \(B\),
with the moment constraints
\begin{align*}
    \E[-B] &\ge -\frac{\alpha}{n} \quad &&\text{average sample usage}\\
    \E\bigl[\tfrac1{\sigma_\Batchsize^2}\bigr] &\ge \theta \quad &&\text{\(Z_\Batchsize\) variance}
\end{align*}
which capture the first factor. Maximizing entropy under moment
constraints is known \autocite{jaynesInformationTheoryStatistical1957}
to result in the Boltzmann (a.k.a. Gibbs) distribution
\[
    \nu(\batchsize)
	 = \Pr(\Batchsize = \batchsize) \propto
    \exp\Bigl(\lambda_1 \frac1{\sigma_\batchsize^2} - \lambda_2 \batchsize\Bigr),
\]
where \(\lambda_1,\lambda_2\) depend on the momentum constraints. We can
now forget the origin of this distribution and use \(\lambda_1,\lambda_2\) as
parameters for the distribution \(\nu\) in Equation~\eqref{eq: the
optimization target} to get close to its minimum. In practice we use a zero
order black box optimizer (Nelder-Mead
\autocite{gaoImplementingNelderMeadSimplex2012}). One could calculate the expectations
of \eqref{eq: the optimization target} under this distribution explicitly and take
manual derivatives with respect to \(\lambda_i\) to investigate this further,
but we wanted to avoid getting too distracted by this tangent. 

We also use the estimated relative standard deviation
\begin{equation}
	\label{eq: relative std}
	\mathrm{rel\_std} = \frac{\sqrt{\widehat{\Var}(\hat{\beta}_0)}}{\hat{\beta}_0}
\end{equation}
as a stopping criterion for sampling. Without extensive testing we found a
tolerance of \(\mathrm{rel\_std}<\mathrm{tol}=0.3\) to be reasonable, cf.~Section~\ref{sec:
sampling efficiency and stability}.

\begin{lemma}[Variance of \(\hat{\beta}_0\) in terms of the empirical batch size distribution]
	\label{lem: variance of beta_0}
	Assuming independence of the samples \(((\frac1{\batchsize_k}),
	Z_{\batchsize_k})_{k\le n}\), the variance of \(\hat{\beta}_0\) is given by
	\[
		\Var(\hat{\beta_0}) = \frac1n 
		\frac{1}{\E[\frac1{\sigma_\Batchsize^2}]}
			\frac{\E[\frac1{\sigma_\Batchsize^2 B^2}]}{\E\Bigl[\frac1{\sigma_\Batchsize^2}
			\bigl(\frac1B - \E[\frac1{B\sigma_\Batchsize^2 \E[1/\sigma_\Batchsize^2]}]\bigr)^2\Bigr]}
	\]
	where \(B\) is distributed according to the empirical batch size distribution
	\(\nu_n = \frac1n\sum_{k=1}^n \delta_{\batchsize_k}\).
\end{lemma}
\begin{proof}
	With the notation \(\sigma_{k}^2 = \sigma^2_{\batchsize_k}\) to describe the
	variance of \(Z_{\batchsize_k}\) it follows from
	\autocite[cf.][Thm.~4.2]{kayFundamentalsStatisticalSignal1993}
	that the variance of the estimator \(\hat{\beta}\) of \(\beta\) using \(n\) samples
	is given by
	\begin{align*}
		\Var(\hat{\beta})
		&= (H^T C^{-1} H)^{-1}
		\\
		&= \frac{1}{
			\bigl(\sum_k \frac1{\sigma_k^2}\bigr)\bigl(\sum_k \frac1{(\sigma_k\batchsize_k)^2}\bigr)
			- (\sum_k \frac1{\sigma_k^2\batchsize_k})^2
		}
		\begin{pmatrix}
			\sum_k \frac1{(\sigma_k\batchsize_k)^2}
			& - \sum_k \frac1{\sigma_k^2\batchsize_k}
			\\
			- \sum_k \frac1{\sigma_k^2\batchsize_k}
			&  \sum_k \frac1{\sigma_k^2}
		\end{pmatrix}
	\end{align*}
	where
	\[
		C := \begin{pmatrix}
			\sigma_{1}^2
			\\
			& \ddots
			\\
			&& \sigma_{n}^2
		\end{pmatrix}
		\qquad
		H := \begin{pmatrix}
			1 & \frac1{\batchsize_1}
			\\
			\vdots
			\\
			1 & \frac1{\batchsize_n}
		\end{pmatrix}.
	\]
	In particular we have
	\[
		\Var(\hat{\beta}_0)= \frac{
			\sum_k \frac1{\sigma_k^2\batchsize_k^2}
		}{
			\bigl(\sum_k \frac1{\sigma_k^2}\bigr)\bigl(\sum_k \frac1{\sigma_k^2\batchsize_k^2}\bigr)
			- (\sum_k \frac1{\sigma_k^2\batchsize_k})^2
		}.
	\]
	With the help of \(\theta :=  \sum_j \frac1{\sigma_j^2}\) and \(\lambda_k :=
	\frac1{\sigma_k^2\theta}\),
	we can reorder the divisor. For this note that since the \(\lambda_k\) sum to \(1\)
	we have
	\begin{align*}
			\sum_k \lambda_k \Bigl(\frac1{\batchsize_k}
			- \sum_j \lambda_j\frac{1}{\batchsize_j}\Bigr)^2
			&=
			\sum_k \lambda_k \Bigl(\frac1{\batchsize_k^2}
			- 2 \frac1{\batchsize_k}\sum_j \lambda_j\frac{1}{\batchsize_j}
			+ \Bigl(\sum_j \lambda_j\frac{1}{\batchsize_j}\Bigr)^2
			\Bigr)
			\\
			&=
			\sum_k \lambda_k \frac1{\batchsize_k^2}
			- 2  \Bigl(
				\sum_k \lambda_k\frac1{\batchsize_k}
			\Bigr)
			+ \Bigl(\sum_k \lambda_k\frac{1}{\batchsize_j}\Bigr)^2
			\\
			&=
			\sum_k \lambda_k \frac1{\batchsize_k^2}
			-  \Bigl(
				\sum_k \lambda_k\frac1{\batchsize_k}
			\Bigr)^2
	\end{align*}
	Where the above is essentially the well known statement \(\E[(Y-\E[Y])^2] =
	\E[Y^2] - \E[Y]^2\) for an appropriate selection of \(Y\). This implies that
	our divisor is given by a weighted variance
	\[
			\theta^2\sum_k \lambda_k \Bigl(\frac1{\batchsize_k}
			- \sum_j \lambda_j\frac{1}{\batchsize_j}\Bigr)^2
			= 
			\theta \sum_k \frac1{\sigma_k^2\batchsize_k^2}
			-  \Bigl(
				\sum_k \frac1{\sigma_k^2\batchsize_k}
			\Bigr)^2,
	\]
	where it is only necessary to plug in the definition of \(\theta\) to see the
	right term is exactly our divisor. Expanding both the enumerator as well as the divisor
	by \(\frac1n\), we obtain
	\[
		\Var(\hat{\beta}_0)
		= \frac1{\theta} \frac{
			\frac1n \sum_{k} \frac1{\sigma_k^2 \batchsize_k^2}
		}{
			\frac1n \sum_{k} \frac1{\sigma_k^2} \bigl(
				\frac1{\batchsize_k} - \sum_{j} \lambda_j \frac1{\batchsize_j}
			\bigr)^2
		}
	\]
	Since \(\theta = n \E[1/\sigma_\Batchsize^2]\) for \(\Batchsize\sim \frac1n\sum_{k=1}^n \delta_{\batchsize_k}\)
	and \(\lambda_k = \frac{1}{n\sigma_k^2 \E[1/\sigma_\Batchsize^2]}\), the above can thus
	be written as
	\[
		\Var(\hat{\beta_0}) = \frac1n 
		\frac{1}{\E[1/\sigma_\Batchsize^2]}
			\frac{\E[\frac1{\sigma_\Batchsize^2 \Batchsize^2}]}{\E\Bigl[\frac1{\sigma_\Batchsize^2}
			\bigl(\frac1\Batchsize - \E[\frac1{\Batchsize\sigma_\Batchsize^2 \E[1/\sigma_\Batchsize^2]}]\bigr)^2\Bigr]},
	\]
	which proves our claim.
\end{proof}
 	\section{Covariance models}\label{appendix: covariance models}

In this section we calculate the step sizes of the covariance models listed
in Table~\ref{table: optimal step size} and plotted in Figure~\ref{fig: rfd step sizes}.
Additionally we calculate the asymptotic learning rate of A-RFD and prove an
Assumption of Corollary~\ref{prop: convergence} for the squared exponential covariance
(Prop.~\ref{prop: sq exp is strictly monotonous in xi}).

\subsection{Squared exponential}

The squared exponential covariance function is given by
\begin{equation}
	\label{eq: sqExp covariance model}
	\ikernel\bigl(\tfrac{\|x-y\|^2}2\bigr)
	= \sigma^2 \exp\bigl(-\tfrac{\|x-y\|^2}{2\scale^2}\bigr).
\end{equation}
Note that \(\sigma^2\) will play no role in the step sizes of RFD due to its
scale invariance (cf. Advantage~\ref{advant: scale invariance}).

\begin{theorem}
	Let \(\Obj\sim \normal(\mu, \ikernel)\) where \(\ikernel\) is the
	\textbf{squared exponential} covariance function \eqref{eq: sqExp covariance
	model}, then we have
	\begin{equation*}
		\stepsize^* \frac{\nabla\Obj(\param)}{\|\nabla\Obj(\param)\|}
		= \argmin_{\step} \E[\Obj(\param - \step)\mid \Obj(\param),\nabla\Obj(\param)]
	\end{equation*}	
	with RFD step size
	\[
		\stepsize^*
		= \frac{
			\scale^2\|\nabla\Obj(\param)\|
		}{
			\sqrt{\bigl(\frac{\mu-\Obj(\param)}{2}\bigr)^2+\scale^2\|\nabla\Obj(\param)\|^2}
			+ \frac{\mu-\Obj(\param)}{2}
		}.
	\]
\end{theorem}
\begin{proof}
	The covariance function \(\ikernel\) is of the form
	\[
		\ikernel(h)
		= \sigma^2 e^{-\frac{h}{\scale^2}}.
	\]
	By Equation~\eqref{eq: simplified step size opt}
	\[
		\stepsize^*
		= -\argmin_{\stepsize} \tfrac{\ikernel(\frac{\stepsize^2}2)}{\ikernel(0)}
		- \stepsize\tfrac{\ikernel'(\frac{\stepsize^2}2)}{\ikernel'(0)}\Theta.
	\]
	where \(\Theta = \frac{\|\nabla\Obj(\param)\|}{\mu-\Obj(\param)}\). We calculate
	\[
		 -\frac{\ikernel\bigl(\frac{\stepsize^2}2\bigr)}{\sqC(0)}
		- \stepsize \frac{\ikernel'\bigl(\frac{\stepsize^2}2\bigr)}{\ikernel'(0)}\Theta
		= -e^{-\frac{\stepsize^2}{2\scale^2}}
		(1 + \stepsize\Theta).
	\]
	This results in the first order condition
	\[
		0 \overset{!}{=} \frac{\stepsize}{\scale^2} e^{-\frac{\stepsize^2}{2\scale^2}}
		(1 + \stepsize\Theta)
		- e^{-\frac{\stepsize^2}{2\scale^2}}\Theta
		= \frac{e^{-\frac{\stepsize^2}{2\scale^2}}}{\scale^2}
		(\stepsize^2\Theta + \stepsize - \scale^2\Theta).
	\]
	Since the exponential can never be zero, we have to solve a quadratic equation.
	Its solution results in
	\begin{equation}
		\label{eq: unstable formula}
		\stepsize^*(\Theta)
		= \sqrt{
			\bigl(\tfrac1{2\Theta}\bigr)^2 + \scale^2 
		}
		- \tfrac1{2\Theta}.
	\end{equation}
	At this point we could stop, but the result is numerically unstable as it suffers
	from catastrophic cancellation. To solve this issue we set \(x=\frac1{2\Theta}\)
	and reorder
	\[
		\stepsize^*
		= \sqrt{x^2 + \scale^2} - x
		= (\sqrt{x^2 + \scale^2} - x) \frac{\sqrt{x^2 + \scale^2} + x}{\sqrt{x^2 + \scale^2} + x}
		= \frac{\cancel{x^2} + \scale^2 - \cancel{x^2}}{\sqrt{x^2+\scale^2} + x}.
	\]
	Re-substituting \(x=\frac1{2\Theta} = \frac{\mu-\Obj(\param)}{2\|\nabla\Obj(\param)\|}\),
	we finally get
	\[
		\stepsize^*
		= \frac{\scale^2}{\sqrt{\bigl(\frac{\mu-\Obj(\param)}{2\|\nabla\Obj(\param)\|}\bigr)^2+\scale^2} + \frac{\mu-\Obj(\param)}{2\|\nabla\Obj(\param)\|}}
		= \frac{\scale^2\|\nabla\Obj(\param)\|}{\sqrt{\bigl(\frac{\mu-\Obj(\param)}{2}\bigr)^2+\scale^2\|\nabla\Obj(\param)\|^2} + \frac{\mu-\Obj(\param)}{2}}.
		\qedhere
	\]
\end{proof}

\begin{prop}[A-RFD for the Squared Exponential Covariance]
	\label{prop: A-RFD for squared exponential covariance}
	If \(\Obj\) is isotropic with squared exponential covariance \eqref{eq:
	sqExp covariance model}, then the step size of A-RFD is given by
	\[
		\hat{\stepsize}
		= \frac{\scale^2}{\mu-\Obj(\param)} \|\nabla\Obj(\param)\|,
	\]
\end{prop}
\begin{proof}
	By Definition~\ref{def: a-rfd} of A-RFD and \(\Theta=\frac{\|\nabla\Obj(\param)\|}{\mu-\Obj(\param)}\)
	we have
	\[
		\hat{\stepsize}(\Theta)
		= \frac{\ikernel(0)}{-\ikernel'(0)}\Theta
		= \frac{\sigma^2\exp(0)}{\frac{\sigma^2}{\scale^2}\exp(0)} \frac{\|\nabla\Obj(\param)\|}{\mu-\Obj(\param)}
		= \scale^2\frac{\|\nabla\Obj(\param)\|}{\mu-\Obj(\param)}.
		\qedhere
	\]
\end{proof}
\begin{prop}
	\label{prop: sq exp is strictly monotonous in xi}
	If \(\Obj\) is isotropic with squared exponential covariance \eqref{eq:
	sqExp covariance model}, then the RFD step sizes are strictly monotonously increasing
	in \(\Theta\).
\end{prop}
\begin{proof}
	Since we know that \(\Theta \to 0\) implies \(\stepsize^* \sim \hat{\stepsize}\to 0\)
	strict monotonicity of \(\stepsize^*\) in \(\Theta\) is sufficient to show that
	\(\stepsize^*\to 0\) also implies \(\Theta\to 0\). So we take the derivative
	of \eqref{eq: unstable formula} resulting in
	\[
		\frac{d}{d\Theta}\stepsize^*
		= 	\frac{1-\frac{1}{\sqrt{1+\scale^2(2\Theta)^2}}}{2\Theta^2},
	\]
	which is greater zero for all \(\Theta>0\).
\end{proof}
 \subsection{Rational quadratic}

The \textbf{rational quadratic} covariance function is given by
\begin{equation}
	\label{eq: rational quadratic}
	\ikernel\bigl(\tfrac{\|x-y\|}{2}\bigr)
	= \sigma^2\left(1+\frac{\|x-y\|^2}{\beta\scale^2}\right)^{-\beta/2}
	\quad \beta > 0.
\end{equation}
It can be viewed as a scale mixture of the squared exponential and converges
to the squared exponential in the limit \(\beta\to\infty\)
\autocite[87]{rasmussenGaussianProcessesMachine2006}.

\begin{theorem}[Rational Quadratic]
	For \(\Obj\sim\normal(\mu, \ikernel)\) where \(\ikernel\) is the \textbf{rational
	quadratic covariance} we have for \(\Theta=\frac{\|\nabla\Obj(\param)\|}{\mu-\Obj(\param)}\ge 0\)
	that the RFD step size is given by
	\begin{align*}
		\stepsize^*
		&= \scale \sqrt{\beta} \Root_\stepsize\left(
			-1 + \tfrac{\sqrt{\beta}}{\scale\Theta}\stepsize
			+ (1+\beta)\stepsize^2 + \tfrac{\sqrt{\beta}}{\scale\Theta}\stepsize^3
		\right).
	\end{align*}
	The unique root of the polynomial in \(\stepsize\) can be found either directly 
	with a formula for polynomials of third degree (e.g. using Cardano's method)
	or by bisection as it is contained in \([0, 1/\sqrt{1+\beta}]\). 
\end{theorem}

\begin{proof}
	By Theorem~\ref{thm: explicit rfd}	we have
	\begin{equation*}
		\stepsize^*
		= \argmin_{\stepsize}-\frac{\ikernel\bigl(\frac{\stepsize^2}2\bigr)}{\ikernel(0)}
		-  \stepsize\frac{\ikernel'\bigl(\frac{\stepsize^2}2\bigr)}{\ikernel'(0)}\Theta
	\end{equation*}
	for \(\ikernel(x) = \sigma^2(1+\frac{2x}{\beta\scale^2})^{-\beta/2}\).
	We therefore
	need to minimize
	\begin{equation*}
		f\bigl(\frac\stepsize{\sqrt{\beta}\scale}\bigr) := -\left(1+\frac{\stepsize^2}{\beta\scale^2}\right)^{-\beta/2}
		- \stepsize \left(1+\frac{\stepsize^2}{\beta\scale^2}\right)^{-\beta/2-1}\Theta.
	\end{equation*}
	Substitute in \(\tilde{\stepsize}:=\frac\stepsize{\sqrt{\beta}\scale}\), then the first
	order condition is
	\begin{equation*}
		0\overset!=f'(\tilde{\stepsize})
		= -\frac{d}{d\tilde{\stepsize}}
		\left(1+\tilde{\stepsize}^2\right)^{-\beta/2}
		+ \sqrt{\beta}\scale\tilde{\stepsize}\left(1+\tilde{\stepsize}^2\right)^{-\beta/2-1}\Theta
	\end{equation*}
	Dividing both sides by \(\sqrt{\beta}\scale\Theta\) we get
	\begin{align*}
		0 = \frac{f'(\tilde{\stepsize})}{\sqrt{\beta}\scale\Theta}&= \tfrac\beta2(1+\tilde{\stepsize}^2)^{-\frac\beta2-1}
			2\tilde{\stepsize}
			\tfrac1{\sqrt{\beta}\scale\Theta}
			+ (1+\tilde{\stepsize}^2)^{-\frac\beta2 -2}
			\left[1 + \tilde{\stepsize}^2 - (\tfrac\beta2+1)2\tilde{\stepsize}^2 \right]
		\\
		&= (1+\tilde{\stepsize}^2)^{-\frac\beta2-2}\underbrace{\left[
			\beta\tilde{\stepsize}\tfrac{1}{\sqrt{\beta}\scale\Theta}
			(1+\tilde{\stepsize}^2)
			- [1-\tilde{\stepsize}^2(1+\beta)]
		\right]}_{
			= -1 + \frac{\sqrt{\beta}}{\scale\Theta}\tilde{\stepsize}
			+ (1+\beta)\tilde{\stepsize}^2
			+ \frac{\sqrt{\beta}}{\scale\Theta}\tilde{\stepsize}^3
		}
	\end{align*}
	Since \(\Theta\ge 0\) and \(\beta > 0\) all coefficients of the polynomial are
	positive except for the shift. The polynomial thus starts out at \(-1\) in
	zero and only increases from there. Therefore there exists a unique positive
	critical point which is a minimum.

	At the point \(\tilde{\stepsize} = \sqrt{1+\beta}\) the quadratic
	term is already larger than \(1\) so the polynomial is positive and we have
	passed the root. The minimum is
	therefore contained in the interval \([0, \sqrt{1+\beta}]\).
	
	After finding the minimum in \(\tilde{\stepsize}\) we return to \(\stepsize\)
	by multiplication with \(\sqrt{\beta}s\).
\end{proof}

\begin{prop}[A-RFD for the Rational Quadratic Covariance]
	If \(\Obj\) is isotropic with rational quadratic covariance \eqref{eq:
	rational quadratic}, then the step size of A-RFD is given by
	\[
		\hat{\stepsize} = \frac{\scale^2}{\mu-\Obj(\param)} \|\nabla\Obj(\param)\|.
	\]
\end{prop}
\begin{proof}
	\(\ikernel(x) = \sigma^2(1+\frac{2x}{\beta\scale^2})^{-\beta/2}\)
	implies by Definition~\ref{def: a-rfd} of A-RFD and \(\Theta=\frac{\|\nabla\Obj(\param)\|}{\mu-\Obj(\param)}\)
	\[
		\hat{\stepsize}(\Theta)
		= \frac{\ikernel_\Obj(0)}{-\ikernel_\Obj'(0)}\Theta 
		= \frac{\sigma^2(1+ 0)^{-\beta/2}}{\frac{\sigma^2}{\scale^2}(1+0)^{-\beta/2-1}} \frac{\|\nabla\Obj(x)\|}{\mu-\Obj(x)}
		= \scale^2\frac{\|\nabla\Obj(x)\|}{\mu-\Obj(x)}.
		\qedhere
	\]
\end{proof}
 \subsection{Matérn}

\begin{definition}
	The Matérn model parametrized by \(\scale >0, \nu\ge 0, \sigma^2\ge 0\) is given by
	\begin{equation}
		\label{eq: matern model}
		\ikernel\bigl(\tfrac{\|x-y\|^2}{2}\bigr)
		= \sigma^2 \frac{2^{1-\nu}}{\Gamma(\nu)}
		\left(\tfrac{\sqrt{2\nu}\|x-y\|}{\scale}\right)^\nu
		\modifiedBessel\left(\tfrac{\sqrt{2\nu}\|x-y\|}{\scale}\right)
	\end{equation}
	where \(\modifiedBessel\) is the modified Bessel function. 
	
	For \(\nu=p+\frac12\) with \(p\in\nat_0\), it can be simplified
	\autocite[cf.][sec.~4.2.1]{rasmussenGaussianProcessesMachine2006} to
	\[
		\ikernel\bigl(\tfrac{\|x-y\|^2}{2}\bigr)
		= \sigma^2 e^{- \tfrac{\sqrt{2\nu} \|x-y\|}{s}}\tfrac{p!}{(2p)!}
		\sum_{k=0}^p \tfrac{(2p-k)!}{(p-k)!k!}\left(\tfrac{2\sqrt{2\nu}}{s}\|x-y\|\right)^k
	\]
\end{definition}

The Matérn model encompasses \citet{rasmussenGaussianProcessesMachine2006}
\begin{itemize}
	\item the \textbf{nugget effect} for \(\nu=0\) (independent randomness)
	\item the \textbf{exponential model} for \(\nu=\tfrac12\) (Ornstein-Uhlenbeck process)
	\item the \textbf{squared exponential model} for \(\nu\to \infty\) with the same
	scale \(\scale\) and variance \(\sigma^2\).
\end{itemize}
The random functions induced by the Matérn model are a.s.
\(\lfloor \nu\rfloor\)-times differentiable \citet{rasmussenGaussianProcessesMachine2006}, i.e. the smoothness of the model
increases with increasing \(\nu\). While the exponential covariance model with
\(\nu=\tfrac12\) results in a random function
which is not yet differentiable, larger \(\nu\) result in increasing
differentiability. As differentiability starts with \(\nu=\frac32\) and we have
a more explicit formula for \(\nu=p+\tfrac12\) the cases \(\nu=\frac32\) and
\(\nu=\frac52\) are
of particular interest.
\begin{quote}
	``[F]or \(\nu \ge 7/2\), in the absence of explicit prior knowledge about the existence
	of higher order derivatives, it is probably very hard from finite noisy
	training examples to distinguish between values of \(\nu \ge 7/2\) (or even to
	distinguish between finite values of \(\nu\) and \(\nu \to\infty\), the smooth squared
	exponential, in this case)'' \autocite[85]{rasmussenGaussianProcessesMachine2006}.
\end{quote}

\begin{theorem}
	Assuming \(\Obj\sim\normal(\mu, \ikernel)\) is a random function where
	\(\ikernel\) is the Matérn covariance such that \(\nu = p+\tfrac12\) with
	\(p\in \{1,2\}\). Then the RFD step
	is given for \(\Theta := \frac{\|\nabla\Obj(\param)\|}{\mu-\Obj(\param)}\ge 0\) by
	\begin{itemize}
		\item \(p=1\)
		\begin{equation*}
			\stepsize^*
			= \frac{s}{\sqrt{3}}\frac{1}{\left(
				1 + \frac{\sqrt{3}}{s\Theta}
			\right)}
		\end{equation*}

		\item \(p=2\)
		\begin{equation*}
			\stepsize^*
			= \frac{\scale}{\sqrt{5}}\frac{
				(1-\zeta)+\sqrt{4 + (1+\zeta)^2}
			}{2(1+\zeta)} \qquad \zeta := \frac{\sqrt{5}}{3\scale\Theta}.
		\end{equation*}
	\end{itemize}
\end{theorem}

\begin{proof}
	We define \(\C(\stepsize) := \ikernel(\frac{\stepsize^2}2)\), which implies
	\[
		\C'(\stepsize) = \ikernel'\bigl(\tfrac{\stepsize^2}2\bigr)\stepsize 
	\]
	or conversely
	\begin{equation}
		\label{eq: quadratic representation from absolute}
		\ikernel'\bigl(\tfrac{\stepsize^2}2\bigr) = \frac{1}\stepsize \C'(\stepsize).
	\end{equation}
	By Theorem~\ref{thm: explicit rfd}, we need to calculate
	\begin{equation}
		\label{eq: reminder minimization theorem}	
        \stepsize^*
        = \argmin_{\stepsize} -\tfrac{\ikernel(\frac{\stepsize^2}2)}{\ikernel(0)}
        - \stepsize\tfrac{\ikernel'(\frac{\stepsize^2}2)}{\ikernel'(0)}\Theta.
	\end{equation} Discarding \(\sigma\) w.l.o.g. due to scale invariance
	(Advantage~\ref{advant: scale invariance}), we have in the case \(p=1\) 
	\[
		\C(\stepsize) = \Bigl(1+\frac{\sqrt{3}}{\scale}\stepsize\Bigr)
		\exp\Bigl(-\frac{\sqrt{3}}{\scale}\stepsize\Bigr).
	\]
	The derivative is then given by
	\[
		\C'(\stepsize) = -\bigl(\tfrac{\sqrt{3}}\scale\bigr)^2\stepsize \exp\bigl(-\tfrac{\sqrt{3}}{\scale}\stepsize\bigr)
	\]
	which implies using \eqref{eq: quadratic representation from absolute}
	\begin{equation}
		\label{eq: derivative matern p=1}
		\ikernel'\bigl(\tfrac{\stepsize^2}2\bigr) = -\bigl(\tfrac{\sqrt{3}}\scale\bigr)^2 \exp\bigl(-\tfrac{\sqrt{3}}{\scale}\stepsize\bigr)
	\end{equation}
	We therefore need to minimize \eqref{eq: reminder minimization theorem} which is given by
	\[
	\begin{aligned}[t]
		&\argmin_\stepsize
		- \bigl(1+\tfrac{\sqrt{3}}{\scale}\stepsize\bigr)
		\exp\bigl(-\tfrac{\sqrt{3}}{\scale}\stepsize\bigr)
		- \stepsize \exp\bigl(-\tfrac{\sqrt{3}}{\scale}\stepsize\bigr)\Theta
		\\
		&= \argmin_\stepsize
		- \bigl(1+(\tfrac{\sqrt{3}}{\scale} + \Theta)\stepsize\bigr)
		\exp\bigl(-\tfrac{\sqrt{3}}{\scale}\stepsize\bigr).
	\end{aligned}
	\]
	The first order condition is
	\[
		0\overset!=
		\Bigl(
			\tfrac{\sqrt{3}}{\scale}\bigl(\cancel{1}+(\tfrac{\sqrt{3}}{\scale} + \Theta)\stepsize\bigr)
			- (\cancel{\tfrac{\sqrt{3}}{\scale}} + \Theta)
		\Bigr)
\]
	which (divided by \(\Theta\) and noting that the exponential can never be zero) is equivalent to
	\[
			0\overset!=\tfrac{\sqrt{3}}{\scale}(\tfrac{\sqrt{3}}{\scale\Theta} + 1)\stepsize - 1
	\]
	reordering for \(\stepsize\) implies
	\[
		\stepsize \overset! = \frac{\scale}{\sqrt{3}} \frac1{\bigl(1+\tfrac{\sqrt{3}}{\scale\Theta}\bigr)}.
	\]
	It is also not difficult to see that this is the point where the derivative switches
	from negative to positive (i.e. a minimum).

	Let us now consider the case \(p=2\), i.e.
	\[
		\C(\stepsize)
		= \bigl(1 + \tfrac{\sqrt{5}}{\scale}\stepsize + \tfrac{5}{3\scale^2}\stepsize^2\bigr)
		\exp\bigl(-\tfrac{\sqrt{5}}{\scale}\stepsize\bigr),
	\]
	which results in
	\[
		\C'(\stepsize)
		= -\tfrac{5}{3\scale^2}\bigl(\stepsize+\tfrac{\sqrt{5}}{\scale}\stepsize^2\bigr)
		\exp\bigl(-\tfrac{\sqrt{5}}{\scale}\stepsize\bigr),
	\]
	i.e. by \eqref{eq: quadratic representation from absolute}
	\begin{equation}
		\label{eq: derivative matern p=2}
		\ikernel'\bigl(\tfrac{\stepsize^2}2\bigr)
		= -\tfrac{5}{3\scale^2}\bigl(1+\tfrac{\sqrt{5}}{\scale}\stepsize\bigr)
		\exp\bigl(-\tfrac{\sqrt{5}}{\scale}\stepsize\bigr).
	\end{equation}
	We therefore need to minimize \eqref{eq: reminder minimization theorem} which is given by
	\[
		\underbrace{
			\Bigl(-\bigl(1 + \tfrac{\sqrt{5}}{\scale}\stepsize + \tfrac{5}{3\scale^2}\stepsize^2\bigr)
			- \stepsize\bigl(1+\tfrac{\sqrt{5}}{\scale}\stepsize\bigr)\Theta\Bigr)
		}_{
			= -\bigl(1 + \bigl(\tfrac{\sqrt{5}}{\scale}+\Theta\bigr)\stepsize
			+ \bigl(\tfrac{5}{3\scale^2} +\tfrac{\sqrt{5}}{\scale}\Theta\bigr)\stepsize^2
			\bigr)
		}
		\exp\bigl(-\tfrac{\sqrt{5}}{\scale}\stepsize\bigr).
	\]
	The first order condition results in
	\begin{align*}
		0 \overset{!}&{=}	
		\tfrac{\sqrt{5}}\scale\Bigl(
			\cancel{1}
			+ \bigl(\tfrac{\sqrt{5}}{\scale}+\Theta\bigr)\stepsize
			+ \bigl(\tfrac{5}{3\scale^2} + \tfrac{\sqrt{5}}{\scale}\Theta\bigr)\stepsize^2
		\Bigr)
		- \Bigl(
			\bigl(\cancel{\tfrac{\sqrt{5}}{\scale}}+\Theta\bigr)
			+ 2\bigl(\tfrac{5}{3\scale^2} +\tfrac{\sqrt{5}}{\scale}\Theta\bigr)\stepsize
		\Bigr)
		\\
		&= -\Theta + \bigl(\tfrac{5}{3\scale^2} - \tfrac{\sqrt{5}}\scale\Theta\bigr)\stepsize
		+\tfrac{\sqrt{5}}\scale\bigl(\tfrac{5}{3\scale^2}+\tfrac{\sqrt{5}}{\scale}\Theta\bigr)\stepsize^2
	\end{align*}
	Dividing everything by \(\Theta\) and using \(\zeta := \frac{\sqrt{5}}{3\scale\Theta}\) we get
	\[
		0\overset!= -1 -\bigl(\zeta - 1\bigr)\bigl(\tfrac{\sqrt{5}}\scale\stepsize\bigr)
		+ \bigl(\zeta + 1\bigr)\bigl(\tfrac{\sqrt{5}}\scale\stepsize\bigr)^2
	\]
	Taking a closer look at the sign changes of the derivative it becomes
	obvious, that the positive root is the
	minimum, i.e.
	\[
		\frac{\sqrt{5}}{\scale} \stepsize
		\overset{!}= \frac{
			(1-\zeta)+\sqrt{(1-\zeta)^2 +4(1+\zeta)}
		}{2(1+\zeta)}
		= \frac{
			(1-\zeta)+\sqrt{4 + (1+\zeta)^2}
		}{2(1+\zeta)}.
		\qedhere
	\]
\end{proof}

\begin{prop}[A-RFD for the Matérn Covariance]
	If \(\Obj\) is isotropic with Matérn covariance \eqref{eq: matern model}
	such that \(\nu = p+\tfrac12\), then the step size of A-RFD for \(p\in
	\{1,2\}\) is given by
	\begin{itemize}
		\item \(p=1\)
		\[
			\hat{\stepsize}
			= \frac{\scale^2}{3} \frac{\|\nabla\Obj(x)\|}{\mu - \Obj(x)}
		\]
		\item \(p=2\)
		\[
			\hat{\stepsize}
			= \frac{3\scale^2}{5} \frac{\|\nabla\Obj(x)\|}{\mu - \Obj(x)}
		\]
	\end{itemize}
\end{prop}
\begin{proof}
	Noting \(\Theta = \frac{\|\nabla\Obj(x)\|}{\mu - \Obj(x)}\), we have by
	Definition~\ref{def: a-rfd} of A-RFD for \(p=1\)
	\[
		\hat{\stepsize}
		= \frac{\ikernel(0)}{-\ikernel'(0)} \Theta
		\overset{\eqref{eq: derivative matern p=1}}=
		\frac{\scale^2}{3}\Theta,
	\]
	and in the case \(p=2\)	
	\[
		\hat{\stepsize}
		= \frac{\ikernel(0)}{-\ikernel'(0)} \Theta
		\overset{\eqref{eq: derivative matern p=2}}=
		\frac{3\scale^2}{5} \Theta.
		\qedhere
	\]
\end{proof}  	\section{Proofs}
\label{sec: proofs}

In this section we prove all the claims made in the main body.

\subsection{Section~\ref{sec: rfd}: Random function descent}

\subsubsection{Formal RFD}\label{sec: formal rfd}

As we mentioned in a footnote at the definition of RFD, the fact that the parameters
become random variables as they are selected by random gradients poses some mathematical
challenges which would have been distracting to address in the main body.
In following paragraphs leading up to Definition~\ref{def: formal rfd} we introduce
and discuss the probability theory required to provide a mathematically
sound definition.

For a fixed cost distribution \(\Pr_\Obj\) and any weight vectors \(\param\)
and \(\tilde{\param}\) the conditional distribution
\[
	\E[\Obj(\tilde{\param}) \mid \Obj(\param), \nabla\Obj(\param)]
\]
is by its axiomatic definition a \((\Obj(\param),
\nabla\Obj(\param))\)-measurable random variable. By the factorization lemma
\autocite[Cor.~1.9.7]{klenkeProbabilityTheoryComprehensive2014},
there therefore exists a measurable function \((j,g)\mapsto \varphi_{\param, \tilde{\param}}(j, g)\) such that
the following equation holds almost surely
\begin{equation}
	\label{eq: factorization of conditional expec}
	\varphi_{\param, \tilde\param}(\Obj(\param), \nabla\Obj(\param))
	= \E[\Obj(\tilde{\param}) \mid \Obj(\param), \nabla\Obj(\param)].
\end{equation}
Since it is possible to calculate \(\varphi_{\param, \tilde{\param}}\)
explicitly in the Gaussian case (cf.~\ref{thm: conditional gaussian
distribution}), the function  
\[
		\Phi_{\Pr_\Obj}:
		\begin{cases}
			(\real^\dims \times \real \times \real^\dims) \to \real\\
			(\param, j, g) \mapsto \argmin_{\tilde{\param}} \varphi_{\param, \tilde{\param}}(j, g),
		\end{cases}
\]
which implements some tie-breaker rules for set valued \(\argmin\) is measurable
when \(\Obj\) is Gaussian and its covariance function is sufficiently smooth.
To prove measurability in the general case is a difficult problem of its own,
which we do not attempt to solve here, since we would not utilize the conditional
expectation outside of the Gaussian case anyway (cf.~Section~\ref{sec: BlUE}). For deterministic \(\param\), we
therefore have
\begin{align*}
	\Phi_{\Pr_\Obj}(\param, \Obj(\param), \nabla\Obj(\param))
	&= \argmin_{\tilde{\param}} \varphi_{\param, \tilde{\param}}(
		\Obj(\param), \nabla\Obj(\param)
	)
	\\
	\overset{\eqref{eq: factorization of conditional expec}}&=
	\argmin_{\tilde{\param}} \E[
		\Obj(\tilde{\param}) \mid \Obj(\param), \nabla\Obj(\param)
	].
\end{align*}
So if the parameter vectors \(\param_n\) were deterministic, our formal definition of
RFD and our initial definition would coincide.
But for random weights \(\Param\) \eqref{eq: factorization of conditional
expec} stops to hold in general, i.e.
\[
	\varphi_{\Param, \tilde{\param}}(\Obj(\Param), \nabla\Obj(\Param))
	\neq \E[\Obj(\tilde{\param})\mid \Obj(\Param), \nabla\Obj(\Param)].
\]
If this equation does not need to hold, we similarly have in general
\[
	\Phi_{\Pr_\Obj}(\Param, \Obj(\Param), \nabla\Obj(\Param))
	\neq \argmin_{\tilde{\param}} \E[\Obj(\tilde{\param}) \mid \Obj(\Param), \nabla\Obj(\Param)].
\]
So the following definition is not just a restatement of the original definition of RFD.

\begin{definition}[Formal RFD]
	\label{def: formal rfd}
	For a Gaussian random cost function \(\Obj\), we define the RFD algorithm
	with starting point \(\Param_0 = \param_0 \in \real^\dims\) by
	\[
		\Param_{n+1} := \Phi_{\Pr_\Obj}(\Param_n, \Obj(\Param_n), \nabla\Obj(\Param_n))
	\]
\end{definition}
This is what we effectively do in Theorem~\ref{thm: explicit rfd} under the
additional isotropy assumption, where we calculate the \(\argmin\) under the
assumption that \(\param\) is deterministic (i.e.  we determine
\(\Phi_{\Pr_\Obj}\)), before we plug-in the random variables
\(\Param_n\) to obtain \(\Param_{n+1}\). Similarly this is how the step size
prescriptions of RFD actually work. We first assume deterministic weights and later
plug the random variables into our formulas. For this reason, we avoided large letters
indicating random variables for parameters \(\param\) in the main body.

\subsubsection{Scale invariance}

\scaleInvariance*

Before we get to the proof, let us quickly formulate the statement in mathematical terms.
Let \(\param_n\) be the parameters selected optimizing \(\Obj\) starting in
\(\param_0\) and \(\tilde{\param}_n\) the parameters selected by the same optimizer
optimizing \(\tilde{\Obj}\) starting in \(\tilde{\param}_0\).

If we apply affine linear scaling to cost \(\Obj\) such that
\(\tilde{\Obj}(\param)= a\Obj(\param) + b\) and start optimization in the same
point, i.e. \(\param_0 = \tilde{\param}_0\), then we expect a scale invariant optimizer
to select
\[
	\param_n = \tilde{\param}_n.
\]
If we scale inputs on the other hand (or more generally map them with a bijection \(\phi\)),
then we expect for \(\tilde{\Obj} := \Obj \circ \phi\) and starting point
\(
	\tilde{\param}_0 = \phi^{-1}(\param_0)
\),
that this relationship is retained by an equivariant optimizer, i.e.
\[
	\tilde{\param}_n = \phi^{-1}(\param_n).
\]
Why do we use a different starting point? As an illustrating example, assume that
\(\phi\) maps miles into kilometers. Then \(\tilde{\Obj}\) accepts miles, while
\(\Obj\) accepts kilometers. Then we have to map the initial starting point
\(\param_0\) of \(\Obj\) measured in kilometers into miles
\(\tilde{\param}_0\). \(\phi^{-1}\) is precisely this transformation from
kilometers into miles. An equivariant optimizer should retain this relation, i.e.
no matter if the input is measured in miles or kilometers the same points are
selected.

\begin{proof}
	The following proof will be split into three parts. The first two parts of
	the proof will address a more general audience and ignore the mathematical
	subtleties we discussed in Section~\ref{sec: formal rfd}. In the third part
	we explain to the interested probabilists how to resolve these issues.

\begin{enumerate}[wide]
	\item \textbf{Invariance with regard to affine linear scaling}

	Let \(\tilde{\Obj}(\param):= a\Obj(\param) + b\)
	where \(a>0\) and \(b\in \real\) and assume
	\(\tilde{\param}_0 = \param_0\). With the induction start given, we
	only require the induction step to prove \(\tilde{\param}_n = \param_n\).

	For the induction step, we assume this equation holds up to \(n\).
	Since \(\phi(x) = ax +b\) is a measurable bijection, the sigma algebra\footnote{
		if you are unfamiliar with sigma algebras read them as ``information''.
	}
	generated by
	\[
		(\tilde{\Obj}(\param_n),\nabla\tilde{\Obj}(\param_n))
		= (\phi\circ\Obj(\param_n),a\nabla\Obj(\param_n))
	\]
	is therefore equal to the sigma algebra generated by
	\((\Obj(\param_n),\nabla\Obj(\param_n))\). This implies
	\begin{equation}
	\label{eq: output scale invariance}	
	\begin{aligned}
		\tilde{\param}_{n+1}
		&= \argmin_{\param}\E[\tilde{\Obj}(\param) \mid \tilde{\Obj}(\tilde{\param}_n), \nabla\tilde{\Obj}(\tilde{\param}_n)]
		\\
		\overset{\text{induction}}&= \argmin_\param \E[\tilde{\Obj}(\param) \mid \tilde{\Obj}(\param_n), \nabla\tilde{\Obj}(\param_n)]
		\\
		\overset{\text{sigma alg.}}&= \argmin_\param \E[\tilde{\Obj}(\param) \mid \Obj(\param_n), \nabla\Obj(\param_n)]
		\\
		\overset{\text{linearity}}&= \argmin_\param a\E[\Obj(\param) \mid \Obj(\param_n), \nabla\Obj(\param_n)] + b
		\\
		\overset{\text{monotonicity}}&= \argmin_\param \E[\Obj(\param) \mid \Obj(\param_n), \nabla\Obj(\param_n)]
		\\
		\overset{\text{def.}}&= \param_{n+1}
	\end{aligned}
	\end{equation}
	
	Where we have used the linearity of the conditional expectation and the
	strict monotonicity of \(\phi(x) = ax + b\).

	\item \textbf{Equivariance with regard to certain input bijections}
		
	Let \(\phi\) be a differentiable bijection whose jacobian is invertible everywhere
	and assume \(\tilde{\Obj} := \Obj \circ \phi\). Since \(\phi\) is a
	bijection, \(\phi(M)\) is the domain of \(\Obj\) whenever \(M\) is the domain
	of \(\tilde{\Obj}\).

	For a starting point \(\param_0\in \phi(M)\) we now assume \(\tilde{\param}_0
	= \phi^{-1}(\param_0) \in M\) and are again going to prove the claim 
	\[
		\tilde{\param}_n = \phi^{-1}(\param_n).
	\]
	by induction. Assume that we have this claim up to \(n\). Then we have by
	induction
	\begin{equation}
		\tilde{\Obj}(\tilde{\param}_n)
		= \Obj \circ \phi(\phi^{-1}(\param_n))
		= \Obj(\param_n)
	\end{equation}
	and
	\[
		\nabla\tilde{\Obj}(\tilde{\param}_n)
		= \nabla_{\tilde{\param}_n} (\Obj \circ \phi(\tilde{\param}_n))
		= \phi'(\tilde{\param}_n)(\nabla\Obj)(\phi(\tilde{\param}_n))
		= \phi'(\tilde{\param}_n) \nabla\Obj(\param_n).
	\]
	Since \(\phi'(\tilde{\param}_n)\) is invertible by assumption, the
	\(\sigma\)-algebras generated by the vectors \((\tilde{\Obj}(\tilde{\param}_n),
	\nabla\tilde{\Obj}(\tilde{\param}_n))\) and \(\Obj(\param_n), \nabla\Obj(\param_n)\)
	are identical. But this results in the induction step
	\begin{equation}
	\label{eq: input scale invariance}	
	\begin{aligned}
		\tilde{\param}_{n+1}
		&= \argmin_{\param\in M}\E[\tilde{\Obj}(\param) \mid \tilde{\Obj}(\tilde{\param}_n), \nabla\tilde{\Obj}(\tilde{\param}_n)]
		\\
		\overset{\text{sigma alg.}}&=
		\argmin_{\param\in M}\E[\tilde{\Obj}(\param) \mid \Obj(\param_n), \nabla\Obj(\param_n)]
		\\
		\overset{\text{def.}}&= \argmin_{\param\in M}\E[\Obj\circ \phi(\param) \mid \Obj(\param_n), \nabla\Obj(\param_n)]
		\\
		&= \phi^{-1}\Bigl(
			\underbrace{
				\argmin_{\theta\in \phi(M)}
				\E[\Obj(\theta) \mid \Obj(\param_n), \nabla\Obj(\param_n)]
			}_{\overset{\text{def.}}= \param_{n+1}}
		\Bigr).
	\end{aligned}
	\end{equation}
	where we simply optimize over \(\theta = \phi(\param)\) instead of \(\param\) and
	correct the \(\argmin\) at the end.

	\item \textbf{Addressing the subtleties}

	In equation \eqref{eq: output scale invariance} we have really proven for deterministic \(\param\)
	\[
		\Phi_{\Pr_{\tilde{\Obj}}}(\param, \tilde{\Obj}(\param), \nabla\tilde{\Obj}(\param))
		= \Phi_{\Pr_\Obj}(\param, \Obj(\param), \nabla\Obj(\param)).
	\]
	But this implies with the induction assumption \(\Param_n = \tilde{\Param}_n\)
	\[
		\tilde{\Param}_{n+1} = 
		\Phi_{\Pr_{\tilde{\Obj}}}(
			\tilde{\Param}_n,
			\tilde{\Obj}(\tilde{\Param}_n),
			\nabla\tilde{\Obj}(\tilde{\Param}_n)
		)
		\overset{\text{ind.}}= \Phi_{\Pr_\Obj}(\Param_n, \Obj(\Param_n), \nabla\Obj(\Param_n))
		= \Param_{n+1}.
	\]
	Similarly we have proven in \eqref{eq: input scale invariance} that
	\[
		\Phi_{\Pr_{\tilde{\Obj}}}\bigl(
			\phi^{-1}(\param),
			\tilde{\Obj}(\phi^{-1}(\param)),
			\nabla\tilde{\Obj}(\phi^{-1}(\param))
		\bigr)
		= \phi^{-1}\bigl(
			\Phi_{\Pr_\Obj}(\param, \Obj(\param), \nabla\Obj(\param))
		\bigr).
	\]
	By the induction assumption \(\tilde{\Param}= \phi^{-1}(\Param_n)\), this
	implies
	\begin{align*}
		\tilde{\Param}_{n+1}
		&= \Phi_{\Pr_{\tilde{\Obj}}}(
				\tilde{\Param}_n,
				\tilde{\Obj}(\tilde{\Param}_n),
				\nabla\tilde{\Obj}(\tilde{\Param}_n))
		\\
		\overset{\text{ind.}}&=\Phi_{\Pr_{\tilde{\Obj}}}\bigl(
			\phi^{-1}(\Param_n),
			\tilde{\Obj}(\phi^{-1}(\Param_n)),
			\nabla\tilde{\Obj}(\phi^{-1}(\Param_n))
		\bigr)
		\\
		&= \phi^{-1}\bigl(
			\Phi_{\Pr_\Obj}(\Param_n, \Obj(\Param_n), \nabla\Obj(\Param_n))
		\bigr)
		\\
		&= \phi^{-1}(\Param_{n+1}).
		\qedhere
	\end{align*}
\end{enumerate}

\end{proof}

\subsection{Section~\ref{subsec: explicit rfd}: Relation to gradient descent}

\firstStochTaylor*
\begin{proof}
	\((\Obj(\param),\nabla\Obj(\param),\Obj(\param-\step))\) is a Gaussian vector for
	which the conditional distribution is well known. It is only necessary to calculate
	the covariance matrix. The key ingredient here is to observe that \(\Obj(\param), \partial_1\Obj(\param), \dots, \partial_\dims\Obj(\param)\)
	are all independent, trivializing matrix inversion.

	More formally, by Lemma~\ref{lem: cov of derivatives, isotropy} we have
	\[
		\Cov\Bigl(\begin{pmatrix}
			\Obj(\param)\\
			\nabla\Obj(\param)
		\end{pmatrix}\Bigr)
		= \begin{pmatrix}
			\ikernel(0) & \\
			& -\ikernel'(0)\identity_{\dims\times\dims}
		\end{pmatrix}
	\]
	and
	\[
		\Cov\Bigl(\Obj(\param-\step), \begin{pmatrix}
			\Obj(\param)\\
			\nabla\Obj(\param)
		\end{pmatrix}\Bigr)
		= \begin{pmatrix}
				\ikernel(\frac{\|\step\|^2}{2})\\
				\ikernel'(\frac{\|\step\|^2}{2})\step
		\end{pmatrix}.
	\]
	By Theorem~\ref{thm: conditional gaussian distribution} we therefore know that
	\[
	\begin{aligned}[t]
		&\E[\Obj(\param-\step) \mid \Obj(\param),\nabla\Obj(\param)]	
		\\
		&=  \mu +
		\begin{pmatrix}
				\ikernel(\frac{\|\step\|^2}{2})\\
				\ikernel'(\frac{\|\step\|^2}{2})\step
		\end{pmatrix}^T
		\begin{pmatrix}
			\ikernel(0) & \\
			& -\ikernel'(0)\identity_{\dims\times\dims}
		\end{pmatrix}^{-1}
		\begin{pmatrix}
			\Obj(\param) - \mu\\
			\nabla\Obj(\param)
		\end{pmatrix},
	\end{aligned}
	\]
	which immediately yields the claim.
\end{proof}

\explicitRFD*
\begin{proof}
	The explicit version of RFD follows essentially by fixing the step size
	\(\stepsize = \|\step\|\) and optimizing over the direction first.
	With Lemma~\ref{lem: first stoch Taylor} we have
	\begin{align*}
		&\min_\step 
		\E[\Obj(\param-\step)\mid \Obj(\param),\nabla\Obj(\param)]
		\\
		&= \min_{\stepsize \ge 0} \min_{\step: \|\step\|=\stepsize}
		\mu + \frac{\ikernel\bigl(\frac{\stepsize^2}2\bigr)}{\ikernel(0)}
		(\Obj(\param)-\mu) - \frac{\ikernel'\bigl(\frac{\stepsize^2}2\bigr)}{\ikernel'(0)}
		\langle \step, \nabla\Obj(\param)\rangle
		\\
		&= \min_{\stepsize \ge 0}
		\mu + \frac{\ikernel\bigl(\frac{\stepsize^2}2\bigr)}{\ikernel(0)}
		(\Obj(\param)-\mu) - \frac{\ikernel'\bigl(\frac{\stepsize^2}2\bigr)}{\ikernel'(0)}
		\begin{cases}
			\displaystyle \max_{\step: \|\step\|=\stepsize}
			\langle \step, \nabla\Obj(\param)\rangle
			& \frac{\ikernel'(\frac{\stepsize^2}2)}{\ikernel'(0)} \ge 0
			\\
			\displaystyle \min_{\step: \|\step\|=\stepsize}
			\langle \step, \nabla\Obj(\param)\rangle
			& \frac{\ikernel'(\frac{\stepsize^2}2)}{\ikernel'(0)} < 0.
		\end{cases}
	\end{align*}
	By Lemma~\ref{lem: constrained maximiziation of scalar products} and
	Corollary~\ref{cor: constrained minimization of scalar products} the maximizing or
	minimizing step direction is then given by
	\[
		\step(\stepsize) = \pm \stepsize \frac{\nabla\Obj(\param)}{\|\nabla\Obj(\param)\|}.
	\]
	Where it is typically to be expected, that we have a positive sign. Since that depends
	on the covariance though, we avoid this problem with the following argument:
	Since \(\stepsize\) only appears as \(\stepsize^2\) in the remaining equation, we can
	optimize over \(\stepsize\in \real\) in the outer minimization instead of
	over \(\stepsize \ge 0\) to move the sign into the step size \(\stepsize\) and set
	without loss of generality
	\[
		\step(\stepsize) = \stepsize \frac{\nabla\Obj(\param)}{\|\nabla\Obj(\param)\|}.
	\]
	Since \(\langle \step(\stepsize), \nabla\Obj(\param)\rangle = \stepsize \|\nabla\Obj(\param)\|\)
	the remaining outer minimization problem over the step size is then given by
	\[
		\min_{\stepsize\in\real} \frac{\ikernel\bigl(\frac{\stepsize^2}2\bigr)}{\ikernel(0)}(\Obj(\param)-\mu)
		- \stepsize \frac{\ikernel'\bigl(\frac{\stepsize^2}2\bigr)}{\ikernel'(0)}\|\nabla\Obj(\param)\|,
	\]
	Its minimizer is by definition the RFD step size as given in the Theorem.
\end{proof} 
\subsection{Section~\ref{subsec: rfd step sizes}: RFD-step sizes}

\begin{prop}[Tayloring the step size optimization problem]
	The second order Taylor approximation of the step size optimization problem
	\[
		q_\Theta(\stepsize)
= -\frac{\ikernel\bigl(\tfrac{\stepsize^2}2\bigr)}{\ikernel(0)} 
		- \stepsize\frac{\ikernel'\bigl(\tfrac{\stepsize^2}2\bigr)}{\ikernel'(0)}\Theta
	\]
	around zero is given by
	\[
		T_2q_\Theta(\stepsize)
		= -1 
		- \stepsize \Theta
		+ \stepsize^2\frac{-\ikernel'(0)}{2\ikernel(0)},
	\]
	and minimized by	
	\[
		\hat{\stepsize}
		:= \argmin_\stepsize T_2q_\Theta(\stepsize)
		= \tfrac{\ikernel(0)}{-\ikernel'(0)}\Theta.
	\]
	Furthermore, the Taylor residual is bounded by
	\[
		\bigl|
			q(\stepsize)
			- T_2q(\stepsize)
		\bigr|
		\le
		\stepsize^3 c_0
		\bigl(\tfrac{\stepsize}{4} + \Theta\bigr)
	\]
	with  \(c_0=\frac12\max\{\sup_{\theta\in[0,1]}|\ikernel''(\theta)|, |\ikernel'(0)|\}(\tfrac1{\ikernel(0)} + \tfrac1{|\ikernel'(0)|})<\infty\).
\end{prop}

\begin{proof}
	Using the Taylor approximation with the mean value reminder for \(\ikernel\), we get
	\[\begin{aligned}
		\ikernel\bigl(\tfrac{\stepsize^2}2\bigr)
		&= \ikernel(0) + \ikernel'(0)\tfrac{\stepsize^2}2
		+ \ikernel''(\theta_2)\frac{\bigl(\frac{\stepsize^2}2\bigr)^2}{2!}
		\\
		\ikernel'\bigl(\tfrac{\stepsize^2}2\bigr)
		&= \ikernel'(0) + \ikernel''(\theta_1)\tfrac{\stepsize^2}2
	\end{aligned}
	\]
	for some \(\theta_1,\theta_2\in [0, \frac{\stepsize^2}2]\). This implies
	\[
		q(\stepsize)
		- \underbrace{\Bigl(
			-\bigl(1+\tfrac{\ikernel'(0)}{\ikernel(0)}\tfrac{\stepsize^2}2\bigr)
			- \stepsize \Theta 
		\Bigr)}_{=:T_2q_\Theta(\stepsize)}
		= -\frac{\ikernel''(\theta_2)}{\ikernel(0)}\frac{\stepsize^4}{2^3}
		- \frac{\ikernel''(\theta_1)}{\ikernel'(0)}\frac{\stepsize^3}2\Theta
	\]
	By the following error the optimistically defined \(T_2q_\Theta(\stepsize)\) is
	really the second Taylor approximation (which can be confirmed manually, but we
	deduce it by arguing that its residual is in \(\bigO(\stepsize^3)\)). More
	specifically,
	\begin{align*}
		\bigl|
			q(\stepsize)
			- T_2q(\stepsize)
		\bigr|
		&\le
		\stepsize^3
		\Bigl(
			\tfrac{\sup_{\theta\in [0, \frac{\stepsize^2}2]}|\ikernel''(\theta)|}{2\ikernel(0)}
			\frac{\stepsize}{4}
			+ \tfrac{\sup_{\theta\in [0, \frac{\stepsize^2}2]}|\ikernel''(\theta)|}{2|\ikernel'(0)|}
			\Theta
		\Bigr)
		\\
		\overset{\text{Lem.~\ref{lem: bound on the second derivative of the covariance}}}&\le
		\stepsize^3 c_0
		\bigl(\tfrac{\stepsize}{4} + \Theta\bigr)
	\end{align*}
	It is easy to see for \(\Obj(\param)<\mu\) that \(T_2q(\stepsize)\) is a convex
	parabola due to \(\ikernel'(0)<0\). We thus have
	\[
		\hat{\stepsize}
		:= \argmin_\stepsize T_2q_\Theta(\stepsize)
		= \tfrac{\ikernel(0)}{-\ikernel'(0)}\Theta.
		\qedhere
	\]
\end{proof}

\begin{theorem}[Details of Proposition~\ref{prop: a-rfd well defined}]
	Let \(\Obj\sim\normal(\mu, \ikernel)\) and assume there exists \(\stepsize_0>0\) such
	that the correlation for larger distances \(\stepsize\ge \stepsize_0\) are
	bounded smaller than \(1\), i.e. \(\frac{\ikernel(\stepsize^2/2)}{\ikernel(0)} < \rho \in (0,1)\).
	Then there exists \(K, \Theta_0>0\) such that for all \(\Theta < \Theta_0\)
	\[
		1-K\Theta \le \frac{\stepsize^*(\Theta)}{\hat{\stepsize}(\Theta)}\le 1+ K\Theta.
	\]
	In particular we have \(\stepsize^*(\Theta)\sim \hat{\stepsize}(\Theta)\) as \(\Theta\to 0\) or
	equivalently as \(\hat{\stepsize}\to 0\).
	
\end{theorem}
\begin{proof}
	This follows immediately from Lemma \ref{lem: small step sizes},
	\ref{lem: medium step sizes} and \ref{lem: large step sizes}.
\end{proof}

\convergence*
\begin{proof}
	Assuming RFD converges, its step sizes \(\stepsize^*\) converge to zero. But this implies
	\(\Theta\to 0\) by assumption, i.e.	
	\[
		\Theta = \frac{\|\nabla\Obj(\param)\|}{\mu - \Obj(\param)} \to 0
	\]
	Since \(\Obj(\param)\) is bounded, this implies \(\|\nabla\Obj(\param)\|\to 0\)
	and by continuity the of the gradient, it is zero in its limit. Thus we converge
	to a stationary point. The asymptotic equality follows by Lemma~\ref{lem: small step sizes}
	and \ref{lem: medium step sizes}, as we know \(\stepsize^*\) converges so we do
	not require the assumptions of Lemma~\ref{lem: large step sizes}.
\end{proof}

\subsubsection{Locating the Minimizer}

In the following we want to rule out locations for the RFD step size \(\stepsize^*\) by
proving \(q_\Theta(\stepsize) > q_\Theta(\hat{\stepsize})\) for a wide range of \(\stepsize\).
For this endeavour the relative position of the step size \(\stepsize\) relative to \(\hat{\stepsize}\)
is a useful re-parametrization \[
	\stepsize := \stepsize(\lambda) = \lambda \hat{\stepsize}.
\]
Due to \(\hat{\stepsize}=\frac{\ikernel(0)}{-\ikernel'(0)}\Theta\)
we obtain
\[
	T_2q_\Theta(\stepsize)
	= -1 - \stepsize \Theta + \tfrac{\stepsize^2}2\tfrac{-\ikernel'(0)}{\ikernel(0)}
	= -1 + \lambda(\tfrac{\lambda}2 - 1)\hat{\stepsize}\Theta
\]
On the other hand we have for the bound
\begin{align*}
	|q_\Theta(\stepsize)- T_2q_\Theta(\stepsize)|
	\le 
	\lambda^3\hat{\stepsize}^3 c_0\Bigl(\lambda\tfrac{\ikernel(0)}{4|\ikernel'(0)|} + 1\Bigr)\Theta
\end{align*}
Since \(\hat{\stepsize} = \stepsize(1)\) we thus obtain
\begin{align}
	\nonumber
	\frac{q_\Theta(\stepsize) - q_\Theta(\hat{\stepsize})}{\hat{\stepsize}\Theta}
	&\ge \frac{
		\magenta{T_2q_\Theta(\stepsize)}
		- |q_\Theta(\stepsize) - T_2q_\Theta(\stepsize)| 
		- \teal{T_2q_\Theta(\hat{\stepsize})}
		- \blue{|q_\Theta(\hat{\stepsize}) - T_2q_\Theta(\hat{\stepsize})|}
	}{\hat{\stepsize}\Theta}
	\\
	\nonumber
	&\ge
	\underbrace{\bigl(
		\magenta{\lambda(\tfrac{\lambda}2 - 1)} - (\teal{-\tfrac12})
	\bigr)}_{
		= \tfrac12 - \lambda + \tfrac{\lambda^2}2
	}
	- \hat{\stepsize}^2 c_0\Bigl[
		\lambda^3\Bigl(\lambda\tfrac{\ikernel(0)}{4|\ikernel'(0)|} + 1\Bigr)
		+ \blue{\Bigl(\tfrac{\ikernel(0)}{4|\ikernel'(0)|} + 1\Bigr)}
	\Bigr]
	\\
	\label{eq: minimizer rule-out equation}
	&=
	\tfrac12 (1-\lambda)^2
	- \hat{\stepsize}^2 c_0\Bigl[
		\lambda^3\Bigl(\lambda\tfrac{\ikernel(0)}{4|\ikernel'(0)|} + 1\Bigr)
		+ \Bigl(\tfrac{\ikernel(0)}{4|\ikernel'(0)|} + 1\Bigr)
	\Bigr].
\end{align}
This equation will be the basis of a number of lemmas ruling out various step sizes as minimizers.

\begin{lemma}[Ruling out small step sizes]
	\label{lem: small step sizes}
	If the step size is (much) smaller than the asymptotic step size
	\(\hat{\stepsize}=\hat{\stepsize}(\Theta)\), then it can not be a minimizer. More specifically
	\[
		\frac{\stepsize}{\hat{\stepsize}} \in [0, 1-c_1\Theta) \implies q_\Theta(\stepsize) > q_\Theta(\hat{\stepsize})
	\]	
	where \(c_1 := 2 \tfrac{\ikernel(0)}{|\ikernel'(0)|} \sqrt{c_0\bigl(\tfrac{\ikernel(0)}{4|\ikernel'(0)|} + 1\bigr)}<\infty\).
\end{lemma}
\begin{proof}
	Here we consider the case \(\stepsize\le \hat{\stepsize}\), i.e. \(\lambda \in [0,1]\).
	By \eqref{eq: minimizer rule-out equation} we have
	\begin{align*}
		\frac{q_\Theta(\stepsize) - q_\Theta(\hat{\stepsize})}{\hat{\stepsize}\Theta}
		&\ge
		\tfrac12 (1-\lambda)^2
		- \hat{\stepsize}^2 c_0\Bigl[
			\lambda^3\Bigl(\lambda\tfrac{\ikernel(0)}{4|\ikernel'(0)|} + 1\Bigr)
			+ \Bigl(\tfrac{\ikernel(0)}{4|\ikernel'(0)|} + 1\Bigr)
		\Bigr]
		\\
		&\ge
		\tfrac12 (1-\lambda)^2
		- 2\hat{\stepsize}^2 c_0 
			\Bigl(\tfrac{\ikernel(0)}{4|\ikernel'(0)|} + 1\Bigr)
		\\
		\overset{!}&> 0
	\end{align*}
	for which
	\[
		(1-\lambda)^2
		> 4 \hat{\stepsize}^2 c_0\Bigl(\tfrac{\ikernel(0)}{4|\ikernel'(0)|} + 1\Bigr)
	\]
	is sufficient or equivalently
	\[
		\lambda < 1-2 \hat{\stepsize} \sqrt{c_0\Bigl(\tfrac{\ikernel(0)}{4|\ikernel'(0)|} + 1\Bigr)}
		= 1- \Theta \underbrace{
			2 \tfrac{\ikernel(0)}{|\ikernel'(0)|} \sqrt{c_0\Bigl(\tfrac{\ikernel(0)}{4|\ikernel'(0)|} + 1\Bigr)}
		}_{=:c_1}
	\]
	So for \(\lambda \in [0, 1-\Theta c_1)\) we have \(q_\Theta(\stepsize) > q_\Theta(\hat{\stepsize})\).
\end{proof}

\begin{lemma}[Ruling out medium sized step sizes as minimizer]
	\label{lem: medium step sizes}
	For \(c_2 = 2c_1\) and \(\Theta \le \Theta_0 := \frac1{5c_1}\), we have
	\[
		\frac{\stepsize}{\hat{\stepsize}} \in \bigl(1+ c_2\Theta, \tfrac1{c_2\Theta}\bigr)
		\implies q_\Theta(\stepsize) > q_\Theta(\hat{\stepsize})
	\]
\end{lemma}
\begin{proof}
	Here we consider the case \(\lambda \ge 1\), i.e. \(\stepsize>\hat{\stepsize}\).
	Again starting with \eqref{eq: minimizer rule-out equation} we get
	\begin{align*}
		\frac{q(\stepsize) - q(\hat{\stepsize})}{\hat{\stepsize}\Theta}
		&\ge
		\tfrac12 (1-\lambda)^2
		- \hat{\stepsize}^2 c_0\Bigl[
			\lambda^3\Bigl(\lambda\tfrac{\ikernel(0)}{4|\ikernel'(0)|} + 1\Bigr)
			+ \Bigl(\tfrac{\ikernel(0)}{4|\ikernel'(0)|} + 1\Bigr)
		\Bigr]
		\\
		&\ge
		\tfrac12 (\lambda-1)^2
		- 2\lambda^4\hat{\stepsize}^2 c_0\Bigl(\tfrac{\ikernel(0)}{4|\ikernel'(0)|} + 1\Bigr)
		\\
		\overset{!}&> 0,
	\end{align*}
	for which
	\[
		\lambda -1 > 2 \lambda^2\hat{\stepsize}\sqrt{c_0\Bigl(\tfrac{\ikernel(0)}{4|\ikernel'(0)|} + 1\Bigr)}
		= c_1\Theta \lambda^2
	\]
	or equivalently
	\[
		\lambda -1 - c_1\Theta\lambda^2
		> 0
	\]
	is sufficient. Note that this is a concave parabola in \(\lambda\). So it is positive
	between its zeros which are characterized by
	\[
		c_1\Theta\lambda^2 - \lambda + 1 = 0.
	\]
	They are thus given by
	\[
		\lambda_{1/2} = \frac{1 \pm \sqrt{1 - 4c_1\Theta}}{2c_1\Theta}.
	\]
	So whenever \(\lambda \in (\lambda_1, \lambda_2)\) we have that
	\(q_\Theta(\stepsize) > q_\Theta(\hat{\stepsize})\). In particular for \(4c_1\Theta\le 1\)
	or equivalently \(\Theta \le \frac1{4c_1}\) we have
	\[
		\lambda_2 \ge \frac1{2c_1\Theta}	 = \frac1{c_2\Theta}
	\]
	To get a bound on \(\lambda_1\) note that the original equation was essentially
	\[
		\lambda \ge 1+c_1\Theta\lambda^2
	\]
	with equality for \(\lambda=\lambda_1\), if \(\Theta\) is reduced, the inequality remains,
	which implies that \(\lambda_1\) is decreasing with \(\Theta\). So
	assuming the inequality is satisfied for a particular \(\lambda\) e.g.
	\(\lambda=\sqrt{2}\) which requires
	\[
		\sqrt{2} \ge 1+ 2c_1\Theta  \iff \Theta \le \tfrac{\sqrt{2} -1}{2c_1},
	\]
	then we know that \(\lambda_1\le \sqrt{2}\) for all smaller \(\Theta\). This
	implies for \(\Theta \le \Theta_0 = \frac1{5c_1} \le \frac{\sqrt{2} -1}{2c_1}\)
	\[
		\lambda_1 = 1+ c_1\Theta\lambda_1^2 \le 1+ \underbrace{2c_1}_{c_2}\Theta.
		\qedhere
	\]
\end{proof}

\begin{lemma}[Ruling out large step sizes as minimizer]
	\label{lem: large step sizes}
	If there exists step size \(\stepsize_0>0\) such that the correlation is
	bounded by some \(\rho<1\), i.e.
	\[
		\frac{\ikernel\bigl(\frac{\stepsize^2}2\bigr)}{\ikernel(0)} \le \rho\in (0,1),
	\]
	for larger step sizes \(\stepsize\ge \stepsize_0\), then there exist
	\(\Theta_0>0\) such that for all \(\Theta<\Theta_0\)
	\[
		\frac{\stepsize}{\hat{\stepsize}} \in
		\bigl(1+c_2\Theta,\infty\bigr)
		\implies q(\stepsize) > q(\hat{\stepsize}),
	\]
	where \(c_2\) is the constant from Lemma~\ref{lem: medium step sizes}.
\end{lemma}
\begin{proof}
	The upper bound \(\frac{1}{c_2\Theta}\) in Lemma~\ref{lem: medium step sizes} is
	only due to the loss of precision of the Taylor approximation. To remove it,
	we take a closer look at the actual \(q_\Theta\) itself. We have the following
	bound for our asymptotic minimum
	\begin{align*}
		\frac{q_\Theta(\hat{\stepsize})}{\Theta}
		&\le \frac{T_2 q_\Theta(\hat{\stepsize}) + |q_\Theta(\hat{\stepsize}) - T_2q_\Theta(\hat{\stepsize})|}{\Theta}
		= -\frac1\Theta - \frac12 \hat{\stepsize}
		+ \hat{\stepsize}^3 \underbrace{c_0\Bigl(\tfrac{\ikernel(0)}{4|\ikernel'(0)|} + 1\Bigr)}_{=:c_3}
		\\
		&\le - \frac1\Theta + \hat{\stepsize}^3 c_3
	\end{align*}
	Which means we have for
	\begin{align*}
		\frac{q_\Theta(\stepsize) - q_\Theta(\hat{\stepsize})}{\Theta}
		&\ge \Bigl(1- \frac{\ikernel\bigl(\tfrac{\stepsize^2}2\bigr)}{\ikernel(0)}\Bigr)\frac1\Theta
		- \stepsize\frac{\ikernel'\bigl(\tfrac{\stepsize^2}2\bigr)}{\ikernel'(0)}
		- \hat{\stepsize}^3 c_3
		\\
		\overset{\text{Lemma~\ref{lem: bound on first derivative of covariance}}}&\ge
		\Bigl(1- \frac{\ikernel\bigl(\tfrac{\stepsize^2}2\bigr)}{\ikernel(0)}\Bigr)\frac1\Theta
		- \tfrac{\sqrt{\ikernel(0)}}{\sqrt{-\ikernel'(0)}}
		- \hat{\stepsize}^3 c_3
		\\
		&\ge
		\bigl(1- M\bigr)\frac1\Theta
		- \tfrac{\sqrt{\ikernel(0)}}{\sqrt{-\ikernel'(0)}}
		- \hat{\stepsize}^3 c_3
		\\
		\overset{!}&> 0,
	\end{align*}
	where we use the assumption that there exists \(\rho\in (0,1)\) such that
	\(\rho \ge \frac{\ikernel\bigl(\tfrac{\stepsize^2}2\bigr)}{\ikernel(0)}\) for
	all \(\stepsize\ge \stepsize_0\) and the fact that we only need
	to consider \(\stepsize \ge \frac{1}{c_2\Theta}\) (due to Lemma~\ref{lem: medium
	step sizes}) which allows a translation of \(\stepsize_0\) into some maximal
	\(\Theta_0\).  Note that \(\hat{\stepsize}\sim \Theta\) vanishes as \(\Theta\to 0\),
	so eventually the term \((1-M)\frac1\Theta\) dominates. Selecting
	\(\Theta_0\) small small enough is thus sufficient to cover everything that is not
	already covered by Lemma~\ref{lem: medium step sizes}.
\end{proof}

\subsubsection{Technical bounds}

\begin{lemma}[Bound on the first derivative of the covariance]
	\label{lem: bound on first derivative of covariance}	
	\[
		\sup_{\stepsize\ge 0}|\ikernel'\bigl(\tfrac{\stepsize^2}2\bigr)\stepsize|
		\le \sqrt{-\ikernel'(0)\ikernel(0)}
	\]
\end{lemma}
\begin{proof}
	Since we have	
	\[
		\Cov(D_v \Obj(x), \Obj(y))
		= \ikernel'\bigl(\tfrac{\|x-y\|^2}2\bigr)\langle x-y, v\rangle
	\]
	we have for a standardized vector \(\|v\|=1\) and \(x-y = \stepsize v\)
	by Cauchy-Schwarz
	\[
		|\ikernel'\bigl(\tfrac{\stepsize^2}2\bigr)\stepsize|
		= |\Cov(D_v \Obj(x), \Obj(y))|
		\overset{\text{C.S.}}\le \sqrt{\Var(D_v \Obj(x)) \Var(\Obj(y))}
		= \sqrt{-\ikernel'(0)\ikernel(0)}.
	\]
	As the bound is independent of \(\stepsize\) this yields the claim.
\end{proof}

\begin{lemma}[Bound on the second derivative of the covariance]
	\label{lem: bound on the second derivative of the covariance}
	\[
		\sup_{\theta\ge 0}|\ikernel''(\theta)|
		\le \max\Bigl\{\sup_{\theta\in[0,1]}|\ikernel''(\theta)|, |\ikernel'(0)|\Bigr\}
	\]
\end{lemma}
\begin{proof}
	Note that
	\[
		\Cov(D_v \Obj(x), D_w\Obj(y))
		= -\ikernel''\bigl(\tfrac{\|x-y\|^2}{2}\bigr)\langle x-y, v\rangle\langle x-y, w\rangle
		- \ikernel'\bigl(\tfrac{\|x-y\|^2}{2}\bigr)\langle v, w\rangle
	\]
	Selecting \(v,w\) as orthonormal vectors (e.g. \(v=e_1, w=e_2\)) and \(x-y := \stepsize (v+w)\)
	for some \(\stepsize > 0\) results in \(\|x-y\|^2 = 2\stepsize^2\) and thus
	by the Cauchy-Schwarz inequality
	\begin{align*}
		\bigl|-\ikernel''(\stepsize^2)\stepsize^2\bigr|
		&= \bigl| \Cov(D_v \Obj(x), D_w\Obj(y))\bigr|
		\\
		\overset{\text{C.S.}}&\le\sqrt{\Var(D_v\Obj(x))\Var(D_w\Obj(y))}
		\\
		&= \sqrt{(-\ikernel'(0))^2}
	\end{align*}
	This implies the claim.
\end{proof}
 
\subsection{Section~\ref{sec: stochastic loss and covariance estimation}: Stochastic loss}

\srfd*
\begin{proof}
	Since \(\epsilon_i\) are conditionally independent between each other and to \(\Obj\),
	as entire functions, the same holds true for \(\nabla\epsilon_i\). As all the mixed
	covariances disappear we have
	\begin{align*}
		\Cov\Bigl(\begin{pmatrix}
			\Obs_n(\param)\\
			\nabla\Obs_n(\param)
		\end{pmatrix}\Bigr)	
		&= 
		\Cov\Bigl(\begin{pmatrix}
			\Obj(\param)\\
			\nabla\Obj(\param)
		\end{pmatrix}\Bigr)	
		+ \frac1{\batchsize_n^2}
		\sum_{i=1}^{\batchsize_n}
		\Cov\Bigl(\begin{pmatrix}
			\noise_n^{(i)}(\param)\\
			\nabla\noise_n^{(i)}(\param)
		\end{pmatrix}\Bigr)	
		\\
		&= 
		\begin{pmatrix}
			\ikernel(0) &
			\\
			& -\ikernel'(0)\identity_{\dims\times\dims}
		\end{pmatrix}
		+ \frac1{\batchsize_n^2}
		\sum_{i=1}^{\batchsize_n}
		\begin{pmatrix}
			\ikernel_\noise(0) &
			\\
			& -\ikernel_\noise'(0)\identity_{\dims\times\dims}
		\end{pmatrix}
		\\
		&= 
		\begin{pmatrix}
			\ikernel(0) + \frac1{\batchsize_n}\ikernel_\noise(0) &
			\\
			& -\Bigl(\ikernel'(0) + \frac1{\batchsize_n} \ikernel_\noise'(0)\Bigr)
			\identity_{\dims\times\dims}.
		\end{pmatrix}
	\end{align*}
	by Lemma~\ref{lem: cov of derivatives, isotropy}. If you want to break up the first
	step we recommend considering individual entries of the covariance matrix to
	convince yourself that all the mixed covariances disappear. Together with the
	fact
	\begin{align*}
		&\Cov\Bigl(
			\Obj(\param-\step),
			\begin{pmatrix}
				\Obs_n(\param)\\
				\nabla\Obs_n(\param)
			\end{pmatrix}
		\Bigr)	
		\\
		&= 
		\Cov\Bigl(
			\Obj(\param-\step),
			\begin{pmatrix}
				\Obj(\param)\\
				\nabla\Obj(\param)
			\end{pmatrix}
		\Bigr)
		+ \frac1{\batchsize_n^2}\sum_{i=1}^{\batchsize_n}
		\underbrace{
		\Cov\Bigl(
			\Obj(\param-\step),
			\begin{pmatrix}
				\noise_n^{(i)}(\param)\\
				\nabla\noise_n^{(i)}(\param)
			\end{pmatrix}
		\Bigr)
		}_{=0}
		\\
		&= \begin{pmatrix}
				\ikernel(\frac{\|\step\|^2}{2})\\
				\ikernel'(\frac{\|\step\|^2}{2})\step
		\end{pmatrix}.
	\end{align*}
	The rest is analogous to Lemma~\ref{lem: first stoch Taylor} and
	Theorem~\ref{thm: explicit rfd}, so we only sketch the remaining steps.
	
	Applying Theorem~\ref{thm: conditional gaussian distribution} as in Lemma~\ref{lem:
	first stoch Taylor} we obtain a
	stochastic version of the stochastic Taylor approximation
	(``stochastic\({}^2\) Taylor approximation'' perhaps?)
	\[\begin{aligned}[t]
		&\E[\Obj(\param-\step) \mid \Obs_n(\param), \nabla\Obs_n(\param)]	
		\\
		&= \mu
		+ \frac{\ikernel\bigl(\frac{\|\step\|^2}2\bigr)}{\ikernel(0) + \frac1{\batchsize_n}\ikernel_\noise(0)}(\Obs_n(\param) - \mu)
		-  \frac{\ikernel'\bigl(\frac{\|\step\|^2}2\bigr)}{\ikernel'(0) + \frac1{\batchsize_n}\ikernel_\noise'(0)}\langle \step, \Obs_n(\param)\rangle.
	\end{aligned}\]
	Minimizing this subject to a constant step size as in Theorem~\ref{thm: explicit rfd}
	results in 
	\begin{align*}
		\stepsize^*
		&= \argmin_{\stepsize\in \real}
		\frac{\ikernel\bigl(\frac{\|\step\|^2}2\bigr)}{\ikernel(0) + \frac1{\batchsize_n}\ikernel_\noise(0)}
		(\Obs_n(\param) - \mu)
		-  \stepsize\frac{\ikernel'\bigl(\frac{\|\step\|^2}2\bigr)}{\ikernel'(0) + \frac1{\batchsize_n}\ikernel_\noise'(0)}
		\|\Obs_n(\param)\|
		\\
		&= \argmin_{\stepsize\in \real}
		- \frac{\ikernel\bigl(\frac{\|\step\|^2}2\bigr)}{\ikernel(0)}
		-  \stepsize\frac{\ikernel'\bigl(\frac{\|\step\|^2}2\bigr)}{\ikernel'(0) + \frac1{\batchsize_n}\ikernel_\noise'(0)}
		\frac{\ikernel(0)}{\ikernel(0) + \frac1{\batchsize_n}\ikernel_\noise(0)}
		\frac{\|\Obs_n(\param)\|}{\mu -\Obs_n(\param)},
	\end{align*}
	where we divided the term by \(\frac{\ikernel(0)}{\ikernel(0) +
	\frac1{\batchsize_n}\ikernel_\noise(0)}\frac{1}{\mu
	-\Obs_n(\param)} \ge 0\)
	to obtain the last equation. The claim follows by definition of \(\stepsize^*(\Theta)\) and
	our redefinition of \(\Theta\).
\end{proof}
  	\section{Extensions}
\label{sec: extensions}

In this section we present a few possible extensions to Theorem~\ref{thm:
explicit rfd}, which are all composable, i.e. it is possible to combine these
extensions without any major problems
(including S-RFD, i.e. Extension~\ref{ext: s-rfd}).

\subsection{Geometric anisotropy/Adaptive step sizes}
\label{sec: geometric anisotropy}

In this section, we generalize isotropy to ``geometric anisotropies''
\autocite[17]{steinInterpolationSpatialData1999}, which provide good insights into
the inner workings of adaptive learning rates (e.g. AdaGrad \autocite{duchiAdaptiveSubgradientMethods2011} and
Adam \autocite{kingmaAdamMethodStochastic2015}).
\begin{definition}[Geometric Anisotropy]
	We say a random function \(\Obj\) exhibits a ``geometric anisotropy'', if there exists
	an invertible matrix \(A\) such that \(\Obj(x) = \rg(Ax)\) for some isotropic random
	function \(\rg\).
\end{definition}

This implies that the expectation of \(\Obj\) is still constant (\(\E[\Obj(x)] =
\E[\rg(Ax)] = \mu\)) and the covariance function of \(\Obj\) is given by
\begin{equation}
	\label{eq: geometric anisotropy characterization}
	\Cov(\Obj(x), \Obj(y))
	= \Cov(\rg(Ax), \rg(Ay))
	= \ikernel\Bigl(\frac{\|A(x-y)\|^2}2\Bigr)
	= \ikernel\Bigl(\frac{\|x-y\|_{A^TA}^2}2\Bigr)
\end{equation}
where \(\|\cdot\|_\Sigma\) is the norm induced by \(\langle x,y\rangle_\Sigma :=
\langle x, \Sigma y\rangle\) for some strictly positive definite matrix \(\Sigma= A^T A\).
Here \eqref{eq: geometric anisotropy characterization} characterizes the set of
random functions with a geometric anisotropy in the Gaussian case, because for
an \(\Obj\) with such a covariance we can always obtain an isotropic \(\rg\) by
\(\rg(x):= \Obj(A^{-1}x)\). This is the whitening transformation we suggest looking
for in order to ensure isotropy in the context of scale invariance (Section~\ref{sec: rfd}).

An important observation is, that Theorem~\ref{thm: characterization of weak input
invariances} implies that \(\Obj\) is still stationary, so the distribution of
\(\Obj\) is still invariant to translations. If stationarity is a problem, this is
therefore not the solution. But geometric anisotropies are a beautiful model to
explain preconditioning and adaptive step sizes. For this, we first determine the RFD steps.

\begin{restatable}[RFD steps under geometric anisotropy]{extension}{geometricAnisotropyRFD}
	\label{ext: geometric anisotropy rfd steps}
	Let \(\Obj\) be a Gaussian random function which exhibits a ``geometric
	anisotropy'' \(A\) and is based on an isotropic random function
	\(\rg\sim\normal(\mu, \ikernel)\). Then the RFD steps are given by
	\[
		\stepsize^* \frac{\Sigma^{-1}\nabla\Obj(\param)}{\|\Sigma^{-1}\nabla\Obj(\param)\|_\Sigma}
		= \argmin_{\step} \E[\Obj(\param- \step) \mid \Obj(\param), \nabla\Obj(\param)]
	\]
	with
	\[
		\stepsize^* = \argmin_\stepsize q_\Theta(\stepsize)	
		\quad\text{where}\quad
		\Theta = \frac{\|\Sigma^{-1}\nabla\Obj(\param)\|_\Sigma}{\mu - \Obj(\param)}.
	\]	
\end{restatable}
\begin{proof}[Proofsketch]
	There are two ways to see this. Either we apply scale invariance
	(Advantage~\ref{advant: scale invariance}) directly to translate the steps on \(\rg\) into
	steps on \(\Obj\). Alternatively one can manually retrace the steps of the proof.
	Details in Subsection~\ref{sec: proof of rfd steps under geometric anisotropy}
\end{proof}

The step direction is therefore
\[
	\Sigma^{-1}	\nabla\Obj(x)
\]
and \(\Sigma^{-1}\) acts as a preconditioner. So how would one obtain \(\Sigma\)? As it
turns out the following holds true (by Lemma~\ref{lem: cov of derivatives, isotropy})
\[
	\E[\nabla\Obj(\param)\nabla\Obj(\param)^T]
	= A^T \E[\nabla\rg(\param)\nabla\rg(\param)^T] A
	= A^T (-\ikernel'(0)\identity)  A
	= -\ikernel'(0)\Sigma
\]
In their proposal of the first ``adaptive'' method, AdaGrad,
\textcite{duchiAdaptiveSubgradientMethods2011} suggest to use the matrix
\[
	G_t = \sum_{k=1}^t \nabla\Obj(\param_k)\nabla\Obj(\param_k)^T,
\]
which is basically already looking like an estimation method of \(\Sigma\).
They then restrict themselves to \(\diag(G_t)\) due to the computational costs of a
full matrix inversion. This results in entry-wise (``adaptive'') learning rates.
Later adaptive methods like RMSProp \autocite{hintonNeuralNetworksMachine2012},
AdaDelta
\autocite{zeilerADADELTAAdaptiveLearning2012} and in particular Adam
\autocite{kingmaAdamMethodStochastic2015} replace this sum with an exponential mean
estimate, i.e. in the case of Adam the decay rate \(\beta_2\) is used to get an
exponential moving average
\[
	v_t
	= \beta_2 v_{t-1} + (1-\beta_2)\diag(\nabla\Obj(\param_t)\nabla\Obj(\param_t)^T)
	= \beta_2 v_{t-1} + (1-\beta_2)(\nabla\Obj(\param_t))^2.
\]
They then take the expectation
\begin{align*}
	\E[v_t]
	&= \E\Bigl[(1-\beta_2)\sum_{k=1}^t \beta_2^{t-k}\nabla\Obj(\param_k)^2\Bigr]
	\\
	&= \E\Bigl[(1-\beta_2)\sum_{k=1}^t \beta_2^{t-k}\nabla\Obj(\param_k)^2\Bigr]
	= (1-\beta_2^t)\underbrace{\E[\nabla\Obj(x_t)^2]}_{\propto \diag(\Sigma)}
\end{align*}
So \(\hat{v}_t = v_t/(1-\beta_2^t)\) in the Adam optimizer is essentially an
estimator for \(\diag(\Sigma)\). It is noteworthy, that
\textcite{kingmaAdamMethodStochastic2015} already used the expectation symbol.
This is despite the fact, that they did not yet model the optimization objective
\(\Obj\) as a random function.

We can not yet explain why they then use the square root of their estimate
\(\diag(\Sigma)^{-1/2}\) instead of \(\diag(\Sigma)^{-1}\) itself. This might have something
to do with the fact that the estimation of \(G_t\) happens online and the
\(\Obj(\param_k)\) are therefore highly correlated. Another reason might be
that the inverse of an estimator has different properties than the estimator
itself. Finally, the fact that only the diagonal is used might also be the reason, if
the preconditioner \(\diag(\Sigma)^{-1/2}\) is simply better when we restrict ourselves
to diagonal matrices.

\subsubsection{Proof of Extension~\ref{ext: geometric anisotropy rfd steps}}
\label{sec: proof of rfd steps under geometric anisotropy}
	Since the application of scale invariance provides no intuition, we provide a
	proof which retraces some of the steps of the original proof.

	Recall, that for an isotropic random function \(\rg\) we have the stochastic
	Taylor approximation
	\[
		\E[\rg(\param - \step) \mid \rg(x), \nabla\rg(x)]
		= \mu + \frac{\ikernel\bigl(\frac{\|\step\|^2}2\bigr)}{\ikernel(0)}(\rg(\param)-\mu)
		+ \frac{\ikernel'\bigl(\frac{\|\step\|^2}2\bigr)}{\ikernel'(0)}\langle \step, \nabla \rg(\param)\rangle
	\]
	This implies for a random function with geometric anisotropy \(\Obj(\param) =
	\rg(A\param)\) that
	\begin{align*}
		&\E[\Obj(\param - \step) \mid \Obj(\param), \nabla\Obj(\param)]
		\\
		&=\E[\rg(A(\param - \step)) \mid \rg(A\param), \nabla\rg(A\param)]
		\\
		&= \mu + \frac{\ikernel\bigl(\frac{\|A\step\|^2}2\bigr)}{\ikernel(0)}(\rg(A\param)-\mu)
		- \frac{\ikernel'\bigl(\frac{\|A\step\|^2}2\bigr)}{\ikernel'(0)}\langle A\step, \nabla \rg(A\param)\rangle
		\\
		&= \mu + \frac{\ikernel\bigl(\frac{\|\step\|_{\Sigma}^2}2\bigr)}{\ikernel(0)}(\Obj(\param)-\mu)
		- \frac{\ikernel'\bigl(\frac{\|\step\|_{\Sigma}^2}2\bigr)}{\ikernel'(0)}\langle \step, \underbrace{A^T\nabla \rg(A\param)}_{=\nabla\Obj(\param)}\rangle
	\end{align*}
	with \(\Sigma:= A^T A\). As in the original proof, we now optimize over the
	direction first, while keeping the step size constant, although we now fix the
	step size with regard to the norm \(\|\cdot\|_\Sigma\) (which basically means that
	we still do the optimization in the isotropic space). Note that
	\[
		\max_{\step} \langle \step, \nabla \Obj(x)\rangle
		\quad \text{s.t.} \quad
		\|\step\|_\Sigma = \stepsize
	\]
	is equivalent to
	\[
		\max_{\step} \langle \step, \Sigma^{-1}\nabla \Obj(x)\rangle_\Sigma
		\quad \text{s.t.} \quad
		\|\step\|_\Sigma = \stepsize
	\]
	which is solved by
	\[
		\pm\stepsize\frac{\Sigma^{-1}\nabla\Obj(x)}{\|\Sigma^{-1}\nabla\Obj(x)\|_\Sigma}
	\]
	The remainder of the proof is exactly the same as in the original.
 \subsection{Conservative RFD}\label{sec: conservative rfd}

In the first paragraph of Section~\ref{sec: rfd} we motivated the relation
between RFD and classical optimization with the observation, that gradient
descent is the minimizer of a regularized first order Taylor approximation
\[
	\tfrac1L \nabla\Obj(\param)
	= \argmin_\step T[\Obj(\param-\step)\mid \Obj(\param),\nabla\Obj(\param)] + \tfrac{L}2 \|\param\|^2.
\]
This regularized Taylor approximation is in fact an upper bound on our function
under the \(L\)-smoothness assumption \autocite{nesterovLecturesConvexOptimization2018},
i.e.
\[
	\Obj(\param-\step)
	\le T[\Obj(\param-\step)\mid \Obj(\param),\nabla\Obj(\param)]
	+ \tfrac{L}2 \|\step\|^2
\]
An improvement on of this upper bound compared to \(\Obj(x)\) therefore
guarantees an improvement of the loss. This guarantee was lost with the
conditional expectation (on purpose, as we wanted to consider the average case).
Losing this guarantee also makes convergence proofs more difficult as they
typically make use of this improvement.
In view of the confidence intervals of Figure~\ref{fig: visualize conditional
expectation}, it is natural to ask for a similar upper bound in the random
setting, where this can only be the top of an confidence
interval. This is provided in the following theorem
\begin{lemma}[An \(\gamma\)-upper bound]
	\label{lem: an 1-eps upper bound}
	We have
	\[
		\Pr\Bigl(
			\Obj(\param-\step)
			\le \E[\Obj(\param - \step)\mid \Obj(\param), \nabla\Obj(\param)]
			+ \rho_\gamma(\|\step\|^2)
		\Bigr) \ge \gamma 
	\]
	for
	\(
		\rho_\gamma(\stepsize^2)	
		:= \cdf^{-1}(\gamma)\sigma(\stepsize^2)
	\)
	with 
	\[
		\sigma^2(\stepsize^2):=
			\ikernel(0) - \frac{\ikernel\bigl(\frac{\stepsize^2}2\bigr)^2}{\ikernel(0)}
			- \frac{\ikernel'\bigl(\frac{\stepsize^2}2\bigr)^2}{-\ikernel'(0)}\stepsize^2
	\]
	where \(\cdf\) is the cumulative distribution function (cdf) of
	the standard normal distribution.
\end{lemma}

\begin{proof}
	Note that the conditional variance is with the usual argument about the
	covariance matrices (cf. the proof of Thoerem~\ref{thm: explicit rfd}) using
	Lemma~\ref{lem: cov of derivatives, isotropy} and an application of
	Theorem~\ref{thm: conditional gaussian distribution} given by
	\[
		\sigma^2(\|\param\|^2)
		:=\Var[\Obj(\param-\step) \mid \Obj(\param), \nabla\Obj(\param)]
		= \ikernel(0) - \frac{\ikernel\bigl(\frac{\|\step\|^2}2\bigr)^2}{\ikernel(0)}
		- \frac{\ikernel'\bigl(\frac{\|\step\|^2}2\bigr)^2}{-\ikernel'(0)}\|\step\|^2.
	\]
	Since the conditional distribution is normal (by~Theorem~\ref{thm:
	conditional gaussian distribution}), we have
	\[
		\frac{
			\Obj(\param- \step) - \E[\Obj(\param - \step)\mid \Obj(\param), \nabla\Obj(\param)]
		}{
			\sigma(\|\param\|^2)
		}
		\sim \normal(0,1).
	\]
	But this implies the claim
	\begin{align*}
		&\Pr\Bigl(
			\Obj(\param-\step)
			\le \E[\Obj(\param - \step)\mid \Obj(\param), \nabla\Obj(\param)]
			+ \rho_\gamma(\|\step\|^2)
		\Bigr)
		\\
		&=\Pr\Bigl(
			\frac{
				\Obj(\param-\step)
				- \E[\Obj(\param - \step)\mid \Obj(\param), \nabla\Obj(\param)]
			}{\sigma(\|\param\|^2)}
			\le \cdf^{-1}(\gamma)
		\Bigr)
		\\
		&= \cdf(\cdf^{-1}(\gamma)) = \gamma.
	\end{align*}
	To avoid the Gaussian assumption, one could apply the Markov inequality instead, or
	another applicable concentration inequality.
\end{proof}

Using this upper bound, we obtain a natural conservative extension of RFD
\begin{extension}[\(\gamma\)-conservative RFD]
	Let \(\Obj\sim\normal(\mu, \ikernel)\) be a Gaussian random function and
	\(\rho_\gamma(\stepsize^2) := \cdf^{-1}(\gamma) \sigma(\stepsize^2)\), where
	\(\sigma\) is the conditional
	standard deviation as defined in Lemma~\ref{lem: an 1-eps upper bound}.
	Then the conservative RFD step direction is given by
	\[
		\stepsize^* \frac{\nabla\Obj(\param)}{\|\nabla\Obj(\param)\|}
		= \argmin_\step \E[\Obj(\param-\step) \mid \Obj(\param), \nabla\Obj(\param)]
		+ \rho_\gamma(\|\step\|^2)
	\]
	and the \(\gamma\)-conservative RFD step size is given by
	\[
		\stepsize^* = 
		\argmin_{\stepsize}\frac{\ikernel\bigl(\frac{\stepsize^2}2\bigr)}{\ikernel(0)}
		(\Obj(\param)-\mu)
		-  \stepsize\frac{\ikernel'\bigl(\frac{\stepsize^2}2\bigr)}{\ikernel'(0)} \|\nabla\Obj(\param)\|
		+ \rho_\gamma(\stepsize^2).
	\]
\end{extension}
\begin{proof}
	The proof is the same as in Theorem~\ref{thm: explicit rfd} with
	Lemma~\ref{lem: first stoch Taylor} replaced by Lemma~\ref{lem: an 1-eps
	upper bound}.
\end{proof}

Taking multiple steps should generally have an averaging effect, so we expect
faster convergence for almost risk neutral minimization of the conditional
expectation (i.e. \(\gamma\approx \frac12\)). Here \(\gamma\) is a natural
parameter to vary conservatism. In a software implementation it might be a
good idea to call this parameter `conservatism' and rescale it to be in
\([0,1]\) instead of \([\tfrac12, 1]\). But formulas look cleaner with \(\gamma\). 

In Bayesian optimization it is much more common to reverse this approach and
minimize a lower confidence bound (`conservatism' \(< 0\) or \(\gamma <
\tfrac12\)) in order to encourage exploration. But since
RFD is forgetful, this is not a good idea for RFD.

\begin{remark}[Conservative RFD coincides asymptotically with RFD in high dimension]
	\label{rem: high dimension}
	While conservative RFD might seem like a good approach to fix the instability
	of RFD under the isotropy assumption on some optimization problems, the
	variance generally vanishes in high dimension\footnote{see Chapter \ref{chap: predictable progres}} and conservative RFD
	coincides asymptotically with RFD. We therefore believe that the underlying
	issue is not an overly risk-affine algorithm but rather that distributional
	assumptions, in particular the stationarity assumption, are violated when
	instabilities occur (cf. Section \ref{sec: random linear model} and
	\ref{sec: characterization of covariance kernels}).
\end{remark}

Nevertheless, conservative RFD might be a good approach for lower dimensional,
risk-sensitive applications.
\clearmargin
 \subsection{Beyond the Gaussian assumption}\label{sec: BlUE}

In this section we sketch how the extension beyond the Gaussian case using the
``best linear unbiased estimator'' BLUE
\autocite[e.g.][ch.~7]{johnsonAppliedMultivariateStatistical2007} works.\footnote{
	this is the approach is common in geostatistics
}

For this we recapitulate what a BLUE is. A \textbf{linear estimator}
\(\hat{Y}\) of \(Y\) using \(X_1,\dots,X_n\) is of the form 
\begin{equation*}
	\hat{Y}\in \Span\{X_1,\dots,X_n\} + \real.
\end{equation*}
The set of \textbf{linear unbiased estimators} is defined as
\begin{align}\label{def: LUE}
	\LUE
	&:= \LUE{}[Y\mid X_1,\dots, X_n]
	\nonumber
	\\
	&:= \bigl\{ \hat{Y} \in \Span\{X_1,\dots,X_n\} + \real : \E[\hat{Y}] = \E[Y]\bigr\}
	\\
	\nonumber
	&= \{ \hat{Y} + \E[Y] : \hat{Y} \in \Span\{X_1-\E[X_1],\dots,X_n-\E[X_n]\}\}.
\end{align}
And the BLUE is the \textbf{best linear unbiased estimator}, i.e.
\begin{equation}\label{def: BLUE}
	\BLUE[Y\mid X_1,\dots, X_n] := \argmin_{\hat{Y}\in\LUE} \E[\|\hat{Y} - Y\|^2].
\end{equation}
Other risk functions to minimize are possible, but this is the usual one.

\begin{lemma}\label{lem: blue is cond. expectation}
	If \(X,Y_1,\dots, Y_n\) are multivariate normal distributed, then we
	have
	\begin{align*}
		\BLUE[Y\mid X_1,\dots,X_n]
		&= \E[Y\mid X_1,\dots, X_n]\\
		\bigg(&=
		\argmin_{\hat{Y}\in\{f(X_1,\dots, X_n) : f\text{ meas.}\}} \E[\|Y-\hat{Y}\|^2]
		\bigg).
	\end{align*}
\end{lemma}
\begin{proof}
	We observe that the conditional expectation of Gaussian
	random variables is linear (Theorem~\ref{thm: conditional gaussian
	distribution}). So as a linear function its \(L^2\) risk must be larger or
	equal to that of the BLUE. And as an \(L^2\) projection
	\autocite[Cor.~8.17]{klenkeProbabilityTheoryComprehensive2014} the
	conditional expectation was already optimal.
\end{proof}

If we now replace the conditional expectation with the BLUE, then all our theory
remains the same because the result in Theorem~\ref{thm: conditional gaussian distribution} 
remains the BLUE for general distributions
\autocite{johnsonAppliedMultivariateStatistical2007}.  Instead of minimizing
\[
	\min_\step\E[\Obj(\param-\step) \mid \Obj(\param), \nabla\Obj(\param)]	
\]
we can therefore always minimize
\[
	\min_\step \BLUE[\Obj(\param-\step) \mid \Obj(\param), \nabla\Obj(\param)]
\]
without the Gaussian assumption and all our results can be translated to this case.
The reader only needs to replace all mentions of Theorem~\ref{thm: conditional
gaussian distribution} with the BLUE equivalent and replace all ``idependence''
claims with ``uncorrelated''.

\end{subappendices}

 		}
	\end{refsection}

	\begin{refsection}
		{
			\newcommand*{\dimV}{d}
			\renewcommand*{\dims}{N}
			\renewcommand*{\param}{x}
			\renewcommand*{\Param}{X}
			\chapter{Predictable progress}
\label{chap: predictable progres}

In this chapter we prove that general first order optimization algorithms
applied to realizations of high dimensional (Gaussian random functions) GRFs
essentially make deterministic progress,
regardless of the initialization and the realization. Since all quantities that
appear are computable we call this property \textbf{predictable progress} of the
optimizer.
\begin{motivation}[Machine learning]
    \label{mot: ml motivation}
    Our research is motivated by an astonishing phenomenon from the training of
    (large) machine learning models which this chapter tackles in a rigorous way.
    `Training' is the process of minimizing a cost (also known as error or risk)
    function \(f:\real^\dims\to\real\) over parameter vectors \(\param\in
    \real^\dims\) by selecting successively parameter vectors
    \(\param_n\in\real^\dims\) based on the evaluations
    \(f(\param_k)\) and gradients \(\nabla f(\param_k)\) at the previous parameter
    vectors \(\param_0,\dots \param_{n-1}\). The parameters are typically the weights
    of a neural network and $N$ is very large. The phenomenon is that over multiple
    initializations \(\param_0^{(1)},\param_0^{(2)},\dots \in \real^\dims\) the
    progress a particular optimizer makes during the optimization is approximately
    equal over different optimization runs:
    \[
        (f(\param_0^{(i)}), f(\param_1^{(i)}), \dots )
        \approx (f(\param_0^{(j)}), f(\param_1^{(j)}), \dots ).
    \]
    Figure~\ref{fig: loss plots} shows predictable progress in practice.
    We demonstrate predictable progress ourselves for the training on the MNIST
    dataset using a standard convolutional neural network (cf. Figure~\ref{fig:
    self loss plots}). The key to predictable progress is high dimensionality of
    parameters (e.g.  the neural network of Figure~\ref{fig: self loss plots}
    has roughly \(\dims = 2.3\) million parameters, which is small in comparison
    to neural networks used for tasks beyond the ``toy-problem'' MNIST).
    Figure~\ref{fig: other val plots} by \autocite{kleinLearningCurvePrediction2017}
    demonstrates that this phenomenon is well known and moreover relied upon for
    training heuristics.
    \begin{figure}[h]
        \begin{fullwidth}
        \begin{subfigure}[t]{0.48\linewidth}
            \centering
            \def\svgwidth{0.86\linewidth}
            \begingroup \makeatletter \providecommand\color[2][]{\errmessage{(Inkscape) Color is used for the text in Inkscape, but the package 'color.sty' is not loaded}\renewcommand\color[2][]{}}\providecommand\transparent[1]{\errmessage{(Inkscape) Transparency is used (non-zero) for the text in Inkscape, but the package 'transparent.sty' is not loaded}\renewcommand\transparent[1]{}}\providecommand\rotatebox[2]{#2}\newcommand*\fsize{\dimexpr\f@size pt\relax}\newcommand*\lineheight[1]{\fontsize{\fsize}{#1\fsize}\selectfont}\ifx\svgwidth\undefined \setlength{\unitlength}{262.5bp}\ifx\svgscale\undefined \relax \else \setlength{\unitlength}{\unitlength * \real{\svgscale}}\fi \else \setlength{\unitlength}{\svgwidth}\fi \global\let\svgwidth\undefined \global\let\svgscale\undefined \makeatother \begin{picture}(1,0.85714286)\lineheight{1}\setlength\tabcolsep{0pt}\put(0,0){\includegraphics[width=\unitlength,page=1]{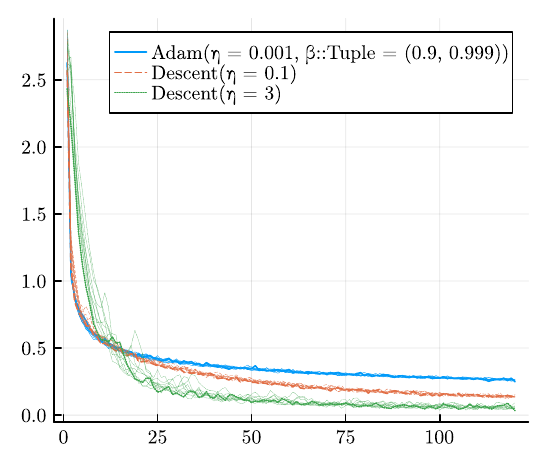}}\end{picture}\endgroup              \caption{
            The plot shows an empirical approximation of the cost sequence resulting
            from the training of a standard convolutional neural network on the MNIST dataset
            \autocite{lecunMNISTDATABASEHandwritten2010}.
            We plot the objective \(f(\param_0^{(i)}),\dots, f(\param_{120}^{(i)})\) against the
            steps \(0,\dots,120\) on the \(x\)-axis. The minimization is performed with three optimization
            algorithms: Adam \autocite{kingmaAdamMethodStochastic2015} (with learning rate \(\eta\) and
            momentum \(\beta\)) in blue and two version of gradient descent
            (learning rate \(\eta=0.1\) and \(\eta=3\)) in red and
            green. Each optimizer was run \(10\) times from randomly selected initializations 
            \(\param_0^{(i)}\) using the (random) default initialization procedure known
            as Glorot initialization \autocite{glorotUnderstandingDifficultyTraining2010}.
            }
            \label{fig: self loss plots}
        \end{subfigure}
        \hspace{\marginparsep}
        \begin{subfigure}[t]{0.48\linewidth}
            \includegraphics[width=\linewidth]{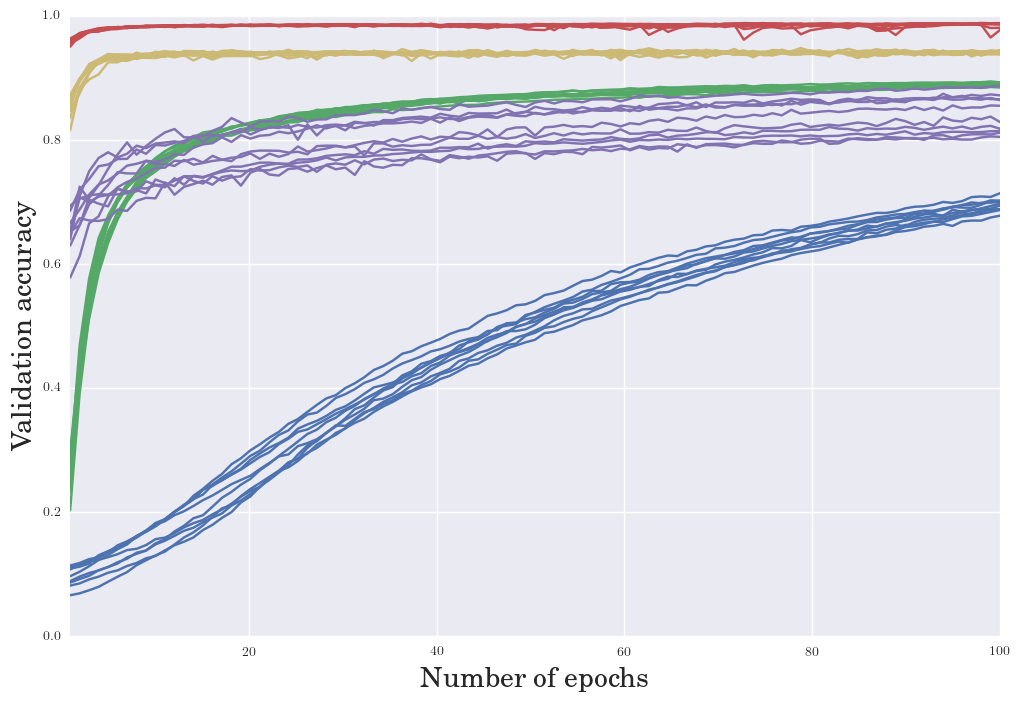}
            \caption{
                The plot, taken from \citet[Figure~3]{kleinLearningCurvePrediction2017},
                shows 5 different optimization hyperparameter configurations
                each evaluated 10 times (grouped by color).
                The predictable progress of the validation accuracy per
                configuration is used in \citet{kleinLearningCurvePrediction2017} to argue that it is
                sufficient to try a configuration once. The overall
                goal of \citet{kleinLearningCurvePrediction2017} is furthermore to fit a parametric models
                to the `learning curves' in order to stop training early and switch to a
                different configuration if the progression does not seem
                promising.
                \\
                The training heuristics that aspire to be `AutoML' (i.e.
                automatically fit data without human intervention) are therefore
                built on this phenomenon. And the fact that these training
                heuristics are so successful demonstrates how ubiquitous
                predictable progress is in practice.
            }
            \label{fig: other val plots}
        \end{subfigure}
\end{fullwidth}
        \sideparmargin{outer}
        \vspace{-\baselineskip}
        \sidepar{
            \makeatletter
            \ifthenelse{\boolean{@twoside}}{
                \checkoddpage\ifoddpage \else \vspace{\baselineskip}\fi
            }{}
            \makeatother
        \caption{Examples of predictable progress in machine learning practice}
        \label{fig: loss plots}
        }
    \end{figure} 
\end{motivation}

The predictable progress phenomenon is certainly surprising since there is no
reason to believe that function values along an optimizer  should be approximatively equal for different initializations. And, in view of handcrafted deterministic counterexamples,
predictable progress is obviously not provable for arbitrary high-dimensional functions.
Instead, we will prove this behavior under the assumption that the cost function \(f\) is
the realization of a Gaussian random function \(\rf\) (GRF), because GRFs have
turned out to be a fruitful model for the statistical properties of cost functions encountered in machine learning 
\autocite[e.g.][]{pascanuDifficultyTrainingRecurrent2013,dauphinIdentifyingAttackingSaddle2014,choromanskaLossSurfacesMultilayer2015}.
Modelling the cost function of Motivation~\ref{mot: ml motivation} as a realization
of a random function could for instance be motivated by the assumption that the
available sequence of data is exchangeable. Since the cost is typically the infinite
data limit over the average sample loss, where the
loss measures the prediction error on a single sample of the data,
de Finetti's representation theorem
\autocite[e.g.][Theorem~1.49]{schervishTheoryStatistics1995} implies the
existence of an underlying random measure driving the randomness of the cost.

While we have motivated predictable progress only for different initializations
\(\param_0\) in Motivation~\ref{mot: ml motivation}, the mathematical results are
much deeper. We show that predictable progress for Gaussian random functions is
also independent of the particular realization, that
is, for a given covariance model the values seen in an optimization run are also
independent of the realization of the GRF \(\rf\).\smallskip

Recently, predictable progress has been proven for \textbf{random quadratic
functions} generated from the mean squared error applied to linear
models
\autocite{paquetteHaltingTimePredictable2022,paquetteUniversalityConjugateGradient2022,deiftConjugateGradientAlgorithm2021,pedregosaAccelerationSpectralDensity2020,scieurUniversalAverageCaseOptimality2020}.
Under the term `state evolution' the asymptotic behavior of `approximate message passing' (AMP)
algorithms has also been analyzed on random quadratic functions
\autocite{bayatiDynamicsMessagePassing2011,javanmardStateEvolutionGeneral2013}. AMP algorithms can be related to
gradient span algorithms in this quadratic setting induced by linear regression
problems \autocite{celentanoEstimationErrorGeneral2020}.
While quadratic functions are a simplified convex setting, these contributions
offered a first explanation for the observed phenomenon of predictable behavior
in high dimensions. We will extend these results to the 
much more general setting of Gaussian random functions. This setting was already used in the machine learning literature for instance by  
\citet{pascanuSaddlePointProblem2014} and
\citet{dauphinIdentifyingAttackingSaddle2014} to explain why
the overwhelming share of critical points are saddle points in high dimension,
whereas the critical points of low dimensional GRFs are dominated by minima and
maxima \autocite{rasmussenGaussianProcessesMachine2006}.
Similarly, our results will show once again that high-dimensional behavior often
eludes low dimensional intuitions.
The specific modelling assumption of \citet{pascanuSaddlePointProblem2014} and
\citet{dauphinIdentifyingAttackingSaddle2014} were stationary isotropic GRFs.
\smallskip

Similar to the convention of capital letters for random variables, we use
bold font to denote random functions \(\rf\). Since our
results are limit theorems in the dimension \(\dims\) we mark this important
dependency in the index. While the parameter sequences would also have to
be indexed by the dimension \(\dims\) we omit this index for notational
clarity. In simplified terms we will prove the following theorem for predictable
optimization in high dimension, precise theorems are given in Section~\ref{sec:setting}:\smallskip

\tikz[baseline]{\draw[thick] (0,1) -- (0,-1); \node[anchor=west] at (0,0) {\begin{minipage}[t]{\columnwidth-2\parindent} \textbf{
        If \(\rf_\dims\) is a (non-stationary) isotropic GRF on a 
        high-dimensional domain \(\real^\dims\), then the (random)  sequence 
        \(\rf_\dims(\Param_0),    \rf_\dims(\Param_1),...\) along the (random) parameter sequence \(\Param_0,
        \Param_1,...\) selected by a standard first order
        optimization algorithm are close to a deterministic sequence
        \(\limf_0,\limf_1,...\) with high probability. 
        }
    \end{minipage}
};
}\smallskip

To be a bit more
precise, we will prove the following in Theorem~\ref{THM: ASYMPTOTICALLY
DETERMINISTIC BEHAVIOR}: Suppose
\(\Param_0,\Param_1,...\) is the (random) sequence of parameter points
obtained by running an optimizer on the (random) function \(\rf_\dims:
\real^\dims \to \real\) initialized at an independent, possibly random point
\(\Param_0\). For `gradient span algorithms' \(\gsa\) (e.g. gradient
descent) and a sequence of (non-stationary) isotropic GRFs
\((\rf_\dims)_{\dims\in\nat}\) we construct a sequence of deterministic real numbers 
\(\limf_0, \limf_1,...\) such that, for $n\in\mathbb N$, 
\[
\lim_{\dims\to\infty}\Pr\Bigl(|\rf_\dims(\Param_\timestep) - \limf_\timestep| > \epsilon\Bigr) = 0,\quad \forall\epsilon>0.
\]
The proof is completely constructive and, in principle, allows for the numerical
computation of the limiting values \(\limf_\timestep\) with complexity
\(\bigO(\timestep^6)\) given an algorithm \(\gsa\), the mean and covariance kernel of
the random functions \(\rf_\dims\) and the length of the initialization vector
\(\Param_0\). In Corollary~\ref{cor: asymptotically identical cost curves} we
will show that an application of the stochastic
triangle inequality yields approximately equal
optimization progress given different initialization points (as can be seen in
Figure~\ref{fig: loss plots}).

The general modelling assumption of (non-stationary) isotropy can be introduced in
two different ways. The first is through the covariance function
(Section~\ref{sec: class of distributions}), the second is axiomatically
through certain invariances. To motivate the latter in the context of
predictability let us highlight a case that must be avoided. Think about some
one-dimensional function
\(f_1:\real\to\real\) where the values along the
optimization path (say of gradient descent) strongly depend on the
initialization. Lift this one dimensional function into \(\real^\dims\) by
\(f_\dims(x) := f_1(\langle v,x\rangle)\) for some direction vector \(v\). Then
\(f_\dims\) certainly does not share the features of high dimensional cost
functions encountered in machine learning, since the cost along the optimization
path depends heavily on the initialization again. To explain the phenomenon
of predictable cost sequences  (cf.~Figure~\ref{fig: loss plots})
the model for cost functions therefore should not be reliant on particular directions
\(v\). Thus, we assume that directions should be exchangeable, which suggests at least 
rotation and reflection invariant random functions as a model. This assumption 
was introduced in Chapter \ref{sec: input invariance} as `(non-stationary)
isotropy' to generalize the well known setting of stationary isotropy.
\smallskip

In spirit, this chapter is close to the \textbf{spin glass} community.
(Non-stationary) isotropy captures the setting of spin-glasses, we use the same
type of dimensional scaling and prove similar limiting results. Finally, there
has also been recent progress in understanding optimization of spin glasses
(see Section~\ref{sec: spin glasses} for further details and e.g.
\citet{auffingerOptimizationRandomHighDimensional2023} for a review of
optimization of spin glasses).

\paragraph*{Outline}

In Sections \ref{sec:algos}-\ref{sec:mainresults} we formalize the setting
(gradient span optimization algorithms and (non-stationary) isotropic GRFs) and state the
main result.
Our mathematical formulation of the original problem using
limit formulations requires a dimensional scaling of the random functions that
is also standard in spin glass theory. This scaling, crucial to our approach, is
discussed in Section \ref{sec: discuss scaling}.
The proof of our main result is given in Section \ref{sec:proof}, first as a
sketch assuming stationary isotropy
and then in detail.

 \section{Setting and main results}\label{sec:setting}

To enable our asymptotic analysis, a careful setting of the scene is required.
Both, the optimization algorithms considered and the distribution of \(\rf_\dims\)
need a representation independent of the dimension \(\dims\) to allow for an analysis of
\(\dims\to\infty\).

\subsection{The class of optimization algorithms \texorpdfstring{\(\gsa\)}{}}\label{sec:algos}

Given a sufficiently smooth function \(f:\mathbb R^\dims\to \mathbb R\),
a naive optimization algorithm is given by gradient descent, whose evaluation points are
recursively defined as
\begin{align*}
	\param_{\timestep}=\param_{\timestep-1}- \alpha \nabla f (\param_{\timestep-1})
\end{align*}
with some initialization \(\param_0\in\real^\dims\) and a learning rate \(\alpha\).
The reader might want to keep gradient descent in mind as a toy example, but
everything we prove holds for a much larger class of first order algorithms that
contains many standard optimizers. Gradient span algorithms (GSA) are this very
general class of first order algorithms which pick the \(\timestep\)-th point
\(\param_\timestep\) from the previous span of gradients
\begin{equation}
	\label{eq: gsa property}
	\param_\timestep
	\in \Span\{
		\param_0, \nabla f(\param_0), \dots, \nabla f(\param_{\timestep-1})
	\}.
\end{equation}
GSAs contain classic gradient descent, momentum methods
such as heavy-ball momentum \autocite{polyakMethodsSpeedingConvergence1964}
and Nesterov's momentum \autocite[e.g.][]{nesterovLecturesConvexOptimization2018},
and also the conjugate
gradient method \autocite{hestenesMethodsConjugateGradients1952}. GSAs do not only
encompass minimizers but also maximizers and all sorts of other algorithms. While the initial point \(\param_0\) is typically not included in the
span, we admit this generalization to allow for concepts such as
`weight normalization' \autocite{salimansWeightNormalizationSimple2016}, which
project points back to the sphere (see Remark~\ref{rem: projection} for further
details). Since the defining property \eqref{eq: gsa property} of gradient span algorithms 
is not sufficient to identify a particular GSA, we define a very general
parametric family of GSAs in Definition~\ref{def: general gsa}. This family
includes all the algorithms mentioned above. Importantly, the parametrization chosen 
does not use dimension specific information and a fixed gradient span algorithm (such as gradient descent) can therefore be used for all dimensions.

\begin{definition}[General gradient span algorithm]
	\label{def: general gsa}
    For a starting point \(\param_0\in \real^\dims\) and a function
    \(f:\real^\dims\to\real\), a \emph{gradient span algorithm} \(\gsa\)
    selects
    \begin{equation}
		\label{eq: definition gsa}
        \param_\timestep := \gsa(f, x_0, \timestep)
		= \lr^{(\param)}_{\timestep} \param_0
		+ \sum_{k=0}^{\timestep-1}\lr_{\timestep,k}^{(g)}\nabla f(\param_k),
	\end{equation}
	where we assume the prefactors \(\lr_\timestep = (\lr^{(\param)}_\timestep) \cup
	(\lr^{(g)}_{\timestep,k})_{k=0,\dots,\timestep-1}\), using the
	union \(\cup\) to indicate a concatenation of tuples,
	to be functions \(\lr_\timestep = \lr_\timestep(\info_{\timestep-1})\)
	of the previous dimensionless information \(\info_{\timestep-1}\) available at time \(\timestep\), that is
	\[
	\begin{aligned}
		\info_\timestep
		&:= \Bigl(f(\param_k): k \le \timestep\Bigr)
		\cup \Bigl(
			\langle v, w\rangle :
			v,w \in (\param_0) \cup \gradients_\timestep
		\Bigr)
		\\
		&\text{with}\quad
		\gradients_\timestep := \bigl(\nabla f(\param_k): k \le \timestep\bigr).
	\end{aligned}
	\]
    The algorithm \(\gsa\) is called \emph{\(\param_0\)-agnostic}, if 
    \begin{itemize}
        \item it remains in the gradient span shifted by \(\param_0\), i.e.
        \[
            \lr_\timestep^{(x)} = 1
            \qquad\forall \timestep\in \nat.
        \]
        \item The prefactors \(\lr_\timestep\) do not use the inner products
        with \(\param_0\), that is they are functions of the reduced information
        \[
            \info^{\setminus \param_0}_\timestep
            = \Bigl(f(\param_k), \langle \nabla f(\param_k), \nabla f(\param_l)\rangle : k, l \le \timestep\Bigr).
        \]
    \end{itemize}
\end{definition}
As mentioned above the reader might just think about ordinary gradient descent,
which is also \(\param_0\)-agnostic.
We introduce the `\(\param_0\)-agnostic' property, to allow for arbitrary initialization
distributions in the stationary isotropic case. Without this property, the algorithm
could, for example, use the length \(\|\Param_0\|\) as an indicator to switch between optimizers and
therefore cause non-deterministic behavior.

Since the main objective of this article is a dimensional limit statement for
fixed algorithms (i.e. fixed choice of prefactors \(\lr\)) it is important to
note that every fixed gradient span algorithm \(\gsa\) is well-defined in all
dimensions because the scalar product \(\langle\cdot,\cdot\rangle\) is well
defined for every dimension \(\dims \in \nat\).\smallskip

It should perhaps be noted, that the initial point \(\param_0\) is not the only
point that may be used by gradient span algorithms, since we have
\[
	\Span\{ (\param_0) \cup G_{n-1}\}
	= \Span\{ \param_0, \dots, \param_{\timestep-1}, \nabla f(\param_0), \dots, \nabla f(\param_{\timestep-1})\}
\]
by induction over \(\timestep\) as \(\param_\timestep\) is selected from this span. We claimed that the inclusion of the initial point \(\param_0\) into the span,
or equivalently the  prefactor \(\lr^{(\param)}_\timestep\), would allow for
algorithms projecting back to the sphere or ball. In the following,
we show that our  general gradient span algorithms also contain gradient span
algorithms with spherical projections.

\begin{remark}[Projection]
	\label{rem: projection}
	Assume that the point \(\param_\timestep\) is defined by
	\[
		\param_\timestep := P \tilde{\param}_\timestep
		\quad\text{with}\quad
		\tilde{\param}_\timestep
		= \tilde{\lr}^{(\param)}_\timestep \param_0
		+ \sum_{k = 0}^{\timestep-1} \tilde{\lr}^{(g)}_{\timestep, k}\nabla f(\param_k),
	\]
	where \(P\) is either a projection to the sphere or the ball. Note that this implies
	either a division by \(\|\tilde{\param}_\timestep\|\) in case of the sphere,
	or a division by \(\max\{\|\tilde{\param_\timestep}\|,1\}\) in case of the
	ball. But the norm
	\[
		\|\tilde{\param}_\timestep\|^2
		=\begin{aligned}[t]
			(\tilde{\lr}^{(\param)}_\timestep)^2
			\underbrace{\|\param_0\|^2}_{\in \info_\timestep}
			&+ 2\tilde{\lr}^{(\param)}_\timestep \sum_{k=0}^{\timestep-1}
			\tilde{\lr}^{(g)}_{\timestep, k}\underbrace{
				\langle \nabla f(\param_k), \param_0\rangle
			}_{\in \info_\timestep}
			\\
			&+ \sum_{k,l=0}^{\timestep-1}
			\tilde{\lr}^{(g)}_{\timestep, k}\tilde{\lr}^{(g)}_{\timestep, l}
			\underbrace{
				\langle \nabla f(\param_k),\nabla f(\param_l)\rangle
			}_{\in \info_{\timestep}}
		\end{aligned}
	\]
	is a function of the information in \(\info_{\timestep-1}\) and can therefore be
	used to define
	\[
		\lr^{(\param)}_\timestep := \tfrac{\tilde{\lr}^{(\param)}_\timestep}{\|\tilde{\param}_\timestep\|}
		\qquad
		\lr^{(g)}_{\timestep,k} := \tfrac{\tilde{\lr}^{(g)}_{\timestep,k}}{\|\tilde{\param}_\timestep\|}
	\]
	in the case of the projection to the sphere and similarly for the projection to the ball.
\end{remark}

\subsection{The class of random function distributions}
\label{sec: class of distributions}

For a random function\footnote{
	Recall that `random process', `stochastic process' and `random field' are
	used synonymously for `random function' in the literature, that their law is
	characterized by all finite dimensional marginals
	\autocite[e.g.][Thm.~14.36]{klenkeProbabilityTheoryComprehensive2014} and that
	Gaussian random functions are fully determined by their mean and covariance.
} 
\(\rf_\dims:\real^\dims\to\real\) on some probability space $(\Omega, \mathcal A, \mathbb P)$ 
recall we denote the mean and covariance functions by
\[
	\mu_{\rf_\dims}(x) := \E[\rf_\dims(x)]
	\quad\text{and}\quad 
	\C_{\rf_\dims}(x,y):= \Cov(\rf_\dims(x), \rf_\dims(y)),
	\qquad x,y\in \real^\dims.
\]
As we want to prove a limit theorem for \(\dims\to\infty\) to make predictions
about the high dimensional cost functions found in practice, the mean \(\mu_{\rf_\dims}\)
and covariance \(\C_{\rf_\dims}\) needs to be defined for every dimension. They
should also remain constant in some sense to avoid arbitrary sequences.
As the domain \(\real^\dims\) changes over \(\dims\) this requires a
parametric form provided by covariance kernels.
To obtain non-trivial results, some scaling is also required (Section~\ref{sec:
discuss scaling}). The scaling used is consistent with the typical scaling
for spin glasses as we discuss in Section
\ref{sec: spin glasses}.
\begin{definition}[Scaled sequence of isotropic Gaussian random functions]
	\label{def: isotropic gaussian random function}	
	\leavevmode
	\begin{enumerate}[label=(\roman*)]
		\item Suppose \(\mu:\real\to\real\) and 	
		\(\kernel:D \to \real\) with 
		\(
			D = \{ \lambda\in \real_{\ge 0}^2 \times \real : |\lambda_3| \le 2\sqrt{\lambda_1 \lambda_2}\}
			\subseteq\real^3.
		\)
		A scaled sequence \((\rf_\dims)_{\dims\in\nat}\) of GRFs
		\(\rf_\dims:\mathbb R^\dims \to \mathbb R\) is called \emph{(non-stationary)
		isotropic} if
		\begin{equation}
			\label{eq: canonical parametrization}
			\mu_{\rf_\dims}(x) = \mu\bigl(\tfrac{\|x\|^2}2\bigr)
			\quad\text{and}\quad
			\C_{\rf_\dims}(x,y)
			= \red{\tfrac1{\dims}} \kernel\bigl(\tfrac{\|x\|^2}2,\tfrac{\|y\|^2}2, \langle x,y\rangle\bigr)
			,\quad x,y\in \real^{\red{\dims}}.
		\end{equation}
		In that case we write \(\rf_\dims\sim\normal(\mu, \kernel)\).
	\item
	Suppose \(\ikernel: \real_{\ge 0}\to \real\). A scaled sequence \((\rf_\dims)\) of
	GRFs \(\rf_\dims:\mathbb R^\dims \to \mathbb R\) is
	called \emph{stationary isotropic} if
	\[
			\mu_{\rf_\dims}(x) = \mu \in \real
			\qquad\text{and}\qquad
			\C_{\rf_\dims}(x,y)
			= \red{\tfrac1{\dims}}\ikernel\bigl(\tfrac{\|x-y\|^2}2\bigr),
			\quad x,y\in \real^{\red\dims}.
		\]
		In that case we write \(\rf_\dims\sim\normal(\mu, \ikernel)\).
	\end{enumerate}	
\end{definition}

\begin{remark}[``Function'']
	In the generality of this definition, \(\rf_\dims\) might only exist in the
	sense of distribution but for our results we will add continuity assumptions
	(cf.~Assumption~\ref{assmpt: smoothness}) that ensure \(\rf_\dims\) to exist as
	a continuous function.
\end{remark}

\begin{remark}[Alternative definition]
	It is also possible to define a scaled sequence of GRFs by scaling one big
	(non-stationary) isotropic centered Gaussian random function \(\rf\) defined on
	an infinite domain as
	\[
		\rf_\dims(\param)
		:= \mu\bigl(\tfrac{\|\param\|^2}2\bigr) + \tfrac1{\sqrt{\dims}}\rf(\param),
		\quad \param\in\real^\dims.
	\]
	We refer the reader who prefers this point of view to sequences of random
	functions with varying domains to Remark \ref{remark:scaling} below.\smallskip
\end{remark}

\begin{remark}[Scaling]
	The importance of the particular dimensional scaling
	\(\red{\tfrac1\dims}\) in \eqref{eq: canonical parametrization} will be discussed in 
	Section~\ref{sec: discuss scaling}. We only want to note that if we
	were to define a Hamiltonian \(\hamiltonian_\dims\) with the canonical scaling used for spin glasses
	\[
		\mu_{\hamiltonian_\dims}(x) = \red{\dims}\mu\bigl(\tfrac{\|x\|^2}2\bigr)
		\quad\text{and}\quad
		\C_{\hamiltonian_\dims}(x,y)
		= \red{\dims} \kernel\bigl(\tfrac{\|x\|^2}2,\tfrac{\|y\|^2}2, \langle x,y\rangle\bigr),
		\quad x,y\in \real^{\red{\dims}},
	\]
	then \(\rf_\dims := \frac{\hamiltonian_\dims}{\dims}\) would have the same distribution
	as in our definition. More details on the relation to spin glasses are given in Section \ref{sec: spin glasses}.
\end{remark}

(Non-stationary) isotropy as a distributional assumption for cost functions was
introduced in Chapter \ref{chap: distributions over functions} as a generalization of the common
stationary isotropy assumption
\autocite[e.g.][]{dauphinIdentifyingAttackingSaddle2014}.
The axiomatic approach to isotropy in Definition \ref{def: distributional input invariance}
explains our characterization of (non-stationary) isotropy as a rotation and
reflection invariant distribution resulting in exchangeable directions.
Non-stationary isotropy is further motivated in Section \ref{sec: random linear model}
with the fact that simple linear models already break the stationary isotropy
assumptions. 

Before getting to the main results, let us recall one more concept. For
a given mean function \(\mu\) and covariance kernel \(\kernel\) we are interested in
a sequence of (non-stationary) isotropic random functions \(\rf_\dims\sim\normal(\mu,
\kernel)\) over the dimension \(\dims\). But it is not clear that a covariance
function \(\kernel\) (resp. \(\ikernel\)) corresponds to a random function for every
dimension. If it does we speak of validity in all dimensions.
Since one can always restrict the domain of a random function
\(\rf_\dims\) to a lower dimensional subspace it should be clear that the only
possible type of restriction is an upper bound on the dimension (cf. Remark
\ref{rem: embedding}).
Further recall, that the isotropic kernels valid in all dimensions correspond to
the isotropic kernels defined on \(\ell^2\) (Lemma \ref{lem: valid in all
dimensions}), which are characterized in Theorem \ref{thm: isotropy characterization}.
Further recall that the stationary isotropic covariance functions which are valid
in all dimensions are given exactly by the kernels of the form\footnote{
	\eg \textcite[Theorem~3.8.5]{sasvariMultivariateCharacteristicCorrelation2013}
	or \textcite{schoenbergMetricSpacesPositive1938}, 
	also listed in Table \ref{table: characterizations}
}
\[
    \ikernel(r) = \int_{[0,\infty)} \exp(-t^2 r) \schoenbergMeas(dt),\qquad r\geq 0,
\]
for probability measure \(\schoenbergMeas\) on \([0,\infty)\), which we
refer to as the Schoenberg measure.
An interesting property of the stationary isotropic covariance kernels
valid in all dimensions is, is that they are all in \(C^\infty\) and the
realization of \(\rf\) thereby smooth (Section \ref{sec: sample regularity of
random functions}).

\subsection{Main result: predictable progress in high dimensions}\label{sec:mainresults}

Our main result proves that the optimization path
\(\rf_\dims(\Param_0),\rf_\dims(\Param_1),...\) is asymptotically deterministic for large dimension \(\dims\). In addition, we
prove a similar statement about the gradient norms
\(\|\nabla\rf_\dims(\Param_\timestep)\|^2\) and more
generally about their scalar products. As we assume continuity of the
prefactors \(\lr_\timestep\) for the gradient span algorithm \(\gsa\), they are
also asymptotically deterministic due to continuous mapping. To reduce the initial complexity,
we first state the simplified corollary for the stationary isotropic case.
In this case some of the technical assumptions can be removed which should allow
the reader to focus on the core message.

\begin{corollary}[Predictable progress of gradient span algorithms, isotropic case]
	\label{cor: asymptotically deterministic behavior}
	Let \(\ikernel\) be a stationary isotropic kernel valid in all dimensions
	(Lemma \ref{lem: valid in all dimensions}) and
	\(\rf_\dims \sim \normal(\mu, \ikernel)\) be a sequence of scaled isotropic
	Gaussian random functions (Definition~\ref{def: isotropic gaussian random
	function}) in \(\dims\).
	Let \(\gsa\) be a general gradient span algorithm (Definition~\ref{def:
	general gsa}), where we assume that its prefactors
	\(\lr_\timestep=\lr_\timestep(\info_{\timestep-1})\) are continuous in the
	information \(\info_{\timestep-1}\), and utilize the most recent gradients, i.e.
	\(\lr^{(g)}_{\timestep,\timestep-1} \neq 0\) for all \(\timestep\in \nat\).
	Let the algorithm \(\gsa\) be applied to \(\rf_\dims\) with independent,
	possibly random starting point \(\Param_0\in \real^\dims\) and denote by
	\(\Param_\timestep\) the points resulting from \(\gsa\). 
	Finally, assume that \(\|\Param_0\|=\radius\) almost surely.
	Then there exist characteristic real numbers
	\begin{align*}
		\limf_\timestep &= \limf_\timestep(\gsa, \mu, \kernel, \radius) \in \real
		\quad \text{and}\quad
		\limgdot_{\timestep,k} = \limgdot_{\timestep,k}(\gsa, \mu, \kernel, \radius) \in \real,
	\end{align*}
	such that, for all \(\epsilon>0\), \(\timestep, k\in\nat\),
	\begin{align*}
		&\lim_{\dims\to \infty}\Pr\bigl(|\rf_\dims(X_\timestep) - \limf_\timestep| > \epsilon\bigr) = 0
		\quad\text{and}\quad
		\\
		&\lim_{\dims\to \infty}\Pr\Bigl(
			\bigl|\langle \nabla\rf_\dims(X_\timestep), \nabla\rf_\dims(X_k)\rangle - \limgdot_{\timestep,k}\bigr| > \epsilon
		\Bigr) = 0.
	\end{align*}
	If furthermore the algorithm is \(\param_0\)-agnostic (e.g. gradient descent), then the limiting values do not depend
	on \(\radius\) and the starting point \(\Param_0\) can have arbitrary distribution
	independent of \(\rf_\dims\).
\end{corollary}
	
Since most algorithms are \(\param_0\)-agnostic, this result explains the
approximately deterministic behavior in high dimension for stationary
isotropic GRFs. The general case is covered by the following remark. 

\begin{remark}[Initialization on the Sphere]
	Initialization procedures used in machine learning like Glorot
	initialization \autocite{glorotUnderstandingDifficultyTraining2010} select the
	entries of
	\(\Param_0\) independent, essentially identically distributed, and scaled in
	such a way that the norm does not diverge.
	The norm therefore obeys a law of large numbers and is thus plausibly
	deterministic in high dimension. The assumption \(\|\Param_0\|=\radius\)
	almost surely is therefore realistic.
	This explains the phenomenon of (approximately) predictable progress of
	optimizers started in a randomly selected initial point using
	Glorot initialization, as observed in Figure~\ref{fig: loss plots}.
\end{remark}

The following proof of Corollary~\ref{cor: asymptotically deterministic behavior} should
also serve as a reading aid for Theorem~\ref{THM: ASYMPTOTICALLY DETERMINISTIC BEHAVIOR},
which states the same result in greater generality.
The proof shows how the general statements of Theorem~\ref{THM: ASYMPTOTICALLY DETERMINISTIC BEHAVIOR}
simplify to the more digestible statement of stationary isotropic GRFs.

\begin{proof}
	Observe that \(\rf_\dims(\Param_k)\) and \(\langle \nabla\rf_\dims(\Param_k), \nabla\rf_\dims(\Param_l)\rangle\)
	for \(k,l\le \timestep\) are members of the random information vector
	\(\Info_\timestep\) defined in Theorem~\ref{THM: ASYMPTOTICALLY
	DETERMINISTIC BEHAVIOR}. Their convergence therefore follows immediately
	from the convergence of the information vector \(\Info_\timestep\) proven in
	Theorem~\ref{THM: ASYMPTOTICALLY DETERMINISTIC BEHAVIOR}.
	Since all stationary isotropic random functions are (non-stationary) isotropic
	(cf.\ Section~\ref{sec: class of distributions}), this Corollary therefore follows immediately from
	Theorem~\ref{THM: ASYMPTOTICALLY DETERMINISTIC BEHAVIOR} once we verified the additional
	required assumptions.\smallskip

	The first assumption is the smoothness assumption on the covariance function  (Assumption~\ref{assmpt: smoothness}). It follows from the fact that
	\emph{all} stationary isotropic random function, which are valid in all dimensions,
	are infinitely differentiable by the Schoenberg characterization (cf.
	Table~\ref{table: characterizations} or \eqref{eq: schoenberg charact} and Section \ref{sec: sample
	regularity of random functions}).
	The second assumption is the strict positive definiteness of
	\((\rf_\dims, \nabla\rf_\dims)\), which also holds for all stationary isotropic random
	functions valid in all dimensions by Corollary~\ref{cor: strict positive
	definite} below. Note that Corollary~\ref{cor: strict positive
	definite} requires that the random function \(\rf_\dims\) is not almost
	surely constant. But if \(\rf_\dims\) were almost surely constant, then we
	get asymptotically deterministic behavior from the fact that
	\(\rf_\dims(\param_0) \sim \normal(\mu, \frac1\dims\ikernel(0))\), i.e.
	\(\rf_\dims(\param_0)\to \mu\). Since \(\rf_\dims\) is almost surely constant,
	we thus obtain \(\rf_\dims \to \mu\) uniformly in probability.
	The \(\param_0\)-agnostic case follows from the stationary case of
	Proposition~\ref{prop: wlog deterministic starting points}.
\end{proof}
We now come to the main theorem of the paper, predictable progress of gradient span algorithms on (non-stationary) isotropic random functions in high dimensions. 
While kernels of stationary isotropic GRFs that are valid in all dimensions are always smooth, and the mean function is always constant,
the same may not be the case for the more general situation of (non-stationary) isotropy. We will therefore assume the following smoothness properties:
\begin{assumption}[Sufficiently smooth]
	\label{assmpt: smoothness}
	Let \(\kernel_{i}(\lambda_1,\lambda_2,\lambda_3)
	:= \frac{d}{d\lambda_i}\kernel(\lambda_1,\lambda_2,\lambda_3)\),
	we assume the partial derivatives
	\begin{equation}
		\label{eq: necessary derivatives}
		\kernel_{12},\quad \kernel_{13},\quad \kernel_{23} \quad\text{and}\quad \kernel_{33}
	\end{equation}
	of the kernel \(\kernel\) exist and are continuous. Furthermore, we assume 
	the derivative of the mean \(\mu\) exists and is continuous.
\end{assumption}
In Section \ref{sec: sample regularity of random functions} we have seen
that the covariance of derivatives is directly related to the derivatives
of the covariance function. For \(\nabla\rf_\dims\) to exist, it is therefore natural
that \(\C_{\rf_\dims}\) has to be differentiable. The covariance has to be two times
differentiable in a sense, as we require the following to be well defined
\[
	\Cov(\partial_i \rf_\dims(x), \partial_j\rf_\dims(y)) = \partial_{x_i}\partial_{x_j}\C_{\rf_\dims}(x,y).
\]
In the case of (non-stationary) isotropic functions, this requires the existence
of \eqref{eq: necessary derivatives} by Lemma~\ref{lem: cov of derivatives,
non-stationary isotropic case}. And the existence of the terms in \eqref{eq:
necessary derivatives} 
is in fact necessary and sufficient for \(\nabla\rf_\dims\)
to exist in an \(L^2\)-sense (see Section \ref{sec: sample regularity of random functions}
for further details).

Here is our main result:
\begin{theorem}[Predictable progress of gradient span algorithms, general case]
	\label{THM: ASYMPTOTICALLY DETERMINISTIC BEHAVIOR}
Let \(\kernel\) be a kernel valid in all dimensions
	(Lemma~\ref{lem: valid in all dimensions}) and  \(\rf_\dims \sim \normal(\mu, \kernel)\) be a sequence of scaled (non-stationary) isotropic Gaussian random
	functions (Definition~\ref{def: isotropic gaussian random
	function}) in $\dims$. Assume that \(\mu\) and \(\kernel\)  are sufficiently smooth
	(Assumption~\ref{assmpt: smoothness}) and that the
	covariance of \((\rf_\dims,\nabla\rf_\dims)\) is strictly positive definite
	(Definition~\ref{def: strict positive definite random function}). Let \(\gsa\) be a general gradient span algorithm (Definition~\ref{def:
	general gsa}), where we assume that its prefactors
	\(\lr_\timestep=\lr_\timestep(\info_{\timestep-1})\) are continuous in the
	information \(\info_{\timestep-1}\), and utilize the most recent gradients, i.e.
	\(\lr^{(g)}_{\timestep,\timestep-1} \neq 0\) for all \(\timestep\in \nat\).
	Let the algorithm \(\gsa\) be applied to \(\rf_\dims\) with independent, possibly
	random starting point \(\Param_0\in \real^\dims\) and denote by
	\(\Param_\timestep\) the points resulting from \(\gsa\). Finally, assume that
	\(\|\Param_0\|=\radius\) almost surely.  Then, for any \(\timestep\in
	\nat\), the random information vector
	\[\begin{aligned}
		\Info_\timestep
		&:= \Bigl(\rf_\dims(\Param_k): k \le \timestep\Bigr)
		\cup \Bigl(
			\langle v, w\rangle :
			v,w \in (\param_0) \cup \gradients_\timestep
		\Bigr)
		\qquad\text{with}
		\\
		\gradients_\timestep &:= \bigl(\nabla \rf_\dims(\Param_k): k \le \timestep\bigr)
	\end{aligned}\]
	converges, as the dimension increases \((\dims\to\infty)\), in probability against a deterministic information vector
	\(\liminfo_\timestep=\liminfo_\timestep(\gsa, \mu, \kernel, \lambda)\).
\end{theorem}

An application of the stochastic triangle inequality yields the following
corollary, which perfectly describes Figure~\ref{fig: loss plots}.

\begin{corollary}[Asymptotically identical progress over multiple initializations]
	\label{cor: asymptotically identical cost curves}
	Assume the setting of Theorem~\ref{THM: ASYMPTOTICALLY DETERMINISTIC
	BEHAVIOR} and let \(\Param_0^{(1)}\), \(\Param_0^{(2)}\) be random initialization
	points selected independently from \(\rf_\dims\). 
	In the non-stationary case, we additionally assume that we we have almost
	surely \(\|\Param_0^{(1)}\| = \|\Param_0^{(2)}\| = \radius\in \real\). Let the sequences
	\((\Param_\timestep^{(i)})_{\timestep\in\nat}\) be generated from
	the initialization \(\Param_0^{(i)}\) by the same general gradient
	span algorithm \(\gsa\).
	Then we have, for all \(\timestep\in \nat\) and all \(\epsilon>0\),
	\[
		\lim_{\dims\to\infty}\Pr\Bigl(
			\max_{k\le \timestep}\Bigl|
				\rf_\dims(\Param_k^{(1)})-\rf_\dims(\Param_k^{(2)})
			\Bigr| > \epsilon
		\Bigr) = 0.
	\]
\end{corollary}
\begin{proof}
	The proof is essentially an application of the stochastic triangle inequality
	\[
		\Bigl\{ |X - Z| > \epsilon\Bigr\}
		\subseteq 
		\Bigl\{ |X - Y| > \frac{\epsilon}2\Bigr\}
		\cup \Bigl\{ |Y - Z| > \frac{\epsilon}2\Bigr\},
	\]
	which yields by Theorem~\ref{THM: ASYMPTOTICALLY DETERMINISTIC BEHAVIOR}
	\begin{align*}
		&\lim_{\dims\to\infty}\Pr\Bigl(
			\max_{k\le \timestep}\Bigl|
				\rf_\dims(\Param_k^{(1)})-\rf_\dims(\Param_k^{(2)})
			\Bigr| > \epsilon
		\Bigr)
		\\
		&\le \lim_{\dims\to\infty}\sum_{k\le \timestep}
		\Bigl[\Pr\Bigl(
			\Bigl|
				\rf_\dims(\Param_k^{(1)})-\limf_k
			\Bigr| > \frac{\epsilon}2
		\Bigr)
		+ \Pr\Bigl(
			\Bigl|
				\rf_\dims(\Param_k^{(2)})-\limf_k
			\Bigr| > \frac{\epsilon}2
		\Bigr)
		\Bigr]
		\\
		&=0.
		\qedhere
	\end{align*}
\end{proof}

Similar to \citet{paquetteHaltingTimePredictable2022}, we obtain so-called
asymptotically deterministic halting times as a corollary. 
The \(\epsilon\)-halting is defined as
\[
	T_\epsilon := \inf\{\timestep > 0: \|\nabla\rf_\dims(X_\timestep)\|^2 \le \epsilon\}.
\]
We are going to prove it is asymptotically equal to the
asymptotic \(\epsilon\)-halting time
\[
	\tau_\epsilon := \inf\{\timestep > 0: \limgdot_\timestep \le \epsilon\} \in \nat,
\]
where  \(\limgdot_\timestep := \limgdot_{\timestep,\timestep}\) is the stochastic limit of
\(\|\nabla\rf_\dims(X_\timestep)\|^2\) in the dimension.
In practical terms this means that optimization always stops at roughly
the same time in high dimension.

\begin{corollary}[Asymptotically deterministic halting times]
	If \(\epsilon\notin (\limgdot_\timestep)_{\timestep\in\nat}\), then
	the halting time is asymptotically deterministic
	\[
		\lim_{\dims\to\infty} \Pr(T_\epsilon = \tau_\epsilon) = 1.
	\]
	If \(\epsilon= \limgdot_\timestep\) for some \(\timestep\), then
	\[
		\lim_{\dims\to\infty} \Pr(T_\epsilon \in [\tau_\epsilon, \tau^+_\epsilon]) = 1.
		\quad \text{with}\quad
		\tau^+_\epsilon := \inf\{\timestep > 0 : \limgdot_\timestep < \epsilon\}.
	\]
\end{corollary}
\begin{proof}
	The proof is identical to the proof of Theorem~4 of \citet{paquetteHaltingTimePredictable2022}.
\end{proof}

In the random quadratic setting random matrix theory already allowed
for the derivation of bounds on the limiting values \(\limf_\timestep\) and \(\limgdot_\timestep\).
This in turn provides bounds on the asymptotic halting
times in the random quadratic setting. 
In our much more general setting we cannot yet give such explicit convergence
guarantees but want to emphasize again that all quantities appearing in our approach
are in principle computable up to finite steps \(\timestep\).

 \section{A discussion of the dimensional scaling}
\label{sec: discuss scaling}

Readers unfamiliar with the \(\frac1\dims\) scaling of the covariance in
\eqref{eq: canonical parametrization} might hypothesize that this reduction of
the variance simply
collapses the random function \(\rf_\dims\) to
the mean \(\mu\) in the asymptotic limit. Since asymptotically deterministic
behavior would then follow trivially, it is important to understand why this is
not the case. The reader familiar with spin glasses  will likely be satisfied by the explanation that
we essentially consider \(\rf_\dims = \frac{H_\dims}\dims\) for a
Hamiltonian \(H_\dims\) whose variance scales with \(\dims\). Therefore the
variance of \(\rf_\dims\) should scale with \(\frac1\dims\). In the following we thus assume the reader to be unfamiliar with spin glasses.\smallskip

To simplify the argument consider a stationary isotropic GRF, i.e.
\(\rf_\dims\sim\normal(\mu,\ikernel)\), the arguments work similarly for any (non-stationary) isotropic GRF. For any \emph{fixed} parameter \(x\) we have \(\rf_\dims(\param)\sim\normal(\mu, \frac1\dims\ikernel(0))\),
so the function value \(\rf_\dims(\param)\) indeed collapses exponentially fast to the mean
by a standard Chernoff-bound
\begin{equation}
    \label{eq: Chernoff bound}
    \Pr\bigl(|\rf_\dims(x) - \mu| \ge t\bigr)
    \le 2\exp\bigl(-\dims\tfrac{t}{2\ikernel(0)}\bigr).
\end{equation}
While the function value \(\rf_\dims(\param)\) collapses to the mean, we will
proceed to show that the gradient \(\nabla\rf_\dims(\param)\) does not.
And if the gradient does not collapse, it is sufficient to follow the gradient
to find points where \(\rf_\dims\) stays away from \(\mu\). So the function \(\rf_\dims\)
does not \emph{uniformly} collapse to the mean.

In Section \ref{sec: sample regularity of random functions} we explained that
the covariance of derivatives of a random function
is of the form
\begin{equation}
	\label{eq: how to calculate the covariance of derivatives}	
    \Cov(\partial_{x_i} \rf_\dims(x), \partial_{y_j}\rf_\dims(y))
= \partial_{x_i}\partial_{y_j} \C_{\rf_\dims}(x,y).
\end{equation}
In the case of a stationary isotropic covariance \(\C_{\rf_\dims}(x,y) =
\frac1\dims \ikernel\bigl(\frac{\|x-y\|^2}2\bigr)\) this implies 
\[
    \partial_{x_i}\partial_{y_j} \C_{\rf_\dims}(x,y)
    = \tfrac{1}{\dims}\Bigl[
        \underbrace{\ikernel''\bigl(\tfrac{\|x-y\|^2}2\bigr)(x_i-y_i)(y_j-x_j)}_{\text{(I)}}
        - \underbrace{\ikernel'\bigl(\tfrac{\|x-y\|^2}2\bigr) \delta_{ij}}_{\text{(II)}}
    \Bigr],
\]
where \(\delta_{ij}\) is the Kronecker delta. Part (I) is zero if \(x=y\).
The Kronecker delta in part (II) then implies that the entries of \(\nabla\rf(x)\)
are uncorrelated and therefore independent by the
Gaussian assumption. Specifically, we have
\[
    \partial_{x_i} \rf_\dims(x)
    \overset{\iid}\sim \normal\bigl(0, -\tfrac1\dims\ikernel'(0)\bigr).
\]
The gradient norm therefore experiences a law of large numbers
\begin{equation}
	\label{eq: slope}
    \|\nabla \rf_\dims(x)\|^2
    = \sum_{i=1}^\dims (\partial_{\param_i}\rf_\dims(x))^2
    \underset{\dims\to\infty}{\overset{p}\to} -\ikernel'(0).
\end{equation}
Thus, while the function value \(\rf_\dims(x)\) collapses to the mean exponentially
fast \eqref{eq: Chernoff bound}, the slope of the function
(counterintuitively) does not. Any other type of scaling would result in vanishing or exploding gradients. Since the supremum over gradient norms is exactly the Lipschitz constant of
\(\rf_\dims\), this scaling is therefore necessary to stay in a class of Lipschitz
functions, as is often assumed in optimization theory. And the Parisi formula 
\eqref{eq: parisi formula} shows in the case of spin glasses, that this
scaling also stabilizes the global maximum.

\begin{remark}[Isoperimetry]
    In view of the counterintuitive observations above, one might ask 
    whether the assumption of isotropy results in very peculiar Lipschitz
    functions.  To understand why that is not the case, consider the concept of
    isoperimetry \autocite[e.g.][]{bubeckUniversalLawRobustness2021}. A random
    variable \(X\) on \(\real^\dims\) is said to satisfy c-isoperimetry, if
    for any \(L\)-Lipschitz function \(f:\real^\dims\to\real\) and any \(t\ge
    0\)	an exponential concentration bound holds
    \[
        \Pr\Bigl(|f(X) - \E[f(X)]| \ge t\Bigr)
        \le 2 \exp\bigl(-\dims\tfrac{ t^2}{2cL}\bigr).
    \]
    Observe that this concentration bound is the mirror image of \eqref{eq: Chernoff bound}.
    While we consider random functions, isoperimetry is concerned with random input \(X\) to
    deterministic \(L\)-Lipschitz functions. In particular Gaussian or uniform input
    satisfies isoperimetry \autocite{bubeckUniversalLawRobustness2021}. Intuitively, this paints the following picture of very high dimensional
    Lipschitz functions:
    Most of the function is equal to the mean except for a few peculiar points.
    In the case of a deterministic functions this requires the exclusion of specific
    deterministic points. In the case of random functions this implies any deterministic point is allowed
    but not the use gradient information. In either case the point we pick must be `independent' 
    of the function. And we show in Theorem~\ref{THM: ASYMPTOTICALLY DETERMINISTIC
    BEHAVIOR} that if both function and input is random, independence is
    essentially a sufficient condition.
\end{remark}

The remark above shows that any Lipschitz function in high dimension appears to
`collapse to the mean'.

Our definition of scaled sequences of random functions can be also be reframed
in terms of a single random function defined on \(\ell^2\), which is
externally scaled.

\begin{remark}[External scaling]\label{remark:scaling}

One could define a a (non-stationary)
isotropic GRF \(\rf\) with covariance
\[
	\C_\rf(x,y) = \kernel\bigl(\tfrac{\|x\|^2}2,\tfrac{\|y\|^2}2, \langle x,y\rangle\bigr),
	\qquad x,y\in \ell^2
\]
on the hilbertspace of sequences \(\ell^2\) which contains the eventually-zero
sequences \(\ell^2_0\). These can also be viewed as the union of
\(\real^\dims\)
\[
	\ell^2_0 := \bigcup_{\dims\in \nat}\real^\dims
	\quad\text{with}\quad
	\real^\dims := \{ (\param_i)_{i\in \nat}\in \ell^2 : \param_i = 0 \quad \forall i > N\}.
\]
Since the scaling \(\frac1{\sqrt{\dims}}\) would scale away any mean
the function \(\rf\) might possess, we assume \(\rf\) to be \emph{centered} and
define
\[
	\rf_\dims(\param)
	:= \mu\bigl(\tfrac{\|\param\|^2}2\bigr) + \tfrac1{\sqrt{\dims}}\rf(\param),
	\quad \param\in\real^\dims.
\]
The distribution of \(\rf_\dims\) is then equivalent to the one in
Definition~\ref{def: isotropic gaussian random function}.
\end{remark}

 \subsection{Relation to spin glasses}\label{sec: spin glasses}

In our asymptotic viewpoint of scaled (non-stationary) isotropic GRFs is very close in spirit to the
dimensional scaling of spin glasses. Spin glasses are modelled with a `Hamiltonian' defined on the domain
\(\{-1, +1\}^\dims\) which represents \{spin up, spin down\}-configurations. This
domain is embedded in the sphere \(\sphere(\sqrt{\dims})\) of radius
\(\sqrt{\dims}\). It turned out that some mathematical question were easier
to answer in this extended domain of spherical spin glasses. The notation
further simplifies if spin glasses are defined on the unit sphere\footnote{
	The definition on the unit sphere removes the need for the definition of the
	`overlap', a scaled scalar product which effectively turns the sphere of
	radius \(\sqrt{\dims}\) into a unit sphere.
} and the domain can also be extended to \(\real^\dims\).
Specifically, define the pure \(p\)-spin Hamiltonian to be
\[
	\hamiltonian_{p, \dims}(x)
	:= \sqrt{\dims}\sum_{i_1,\dots, i_p=0}^{\dims-1} Z_{i_1\dots i_p} x^{(i_1)}\cdots x^{(i_p)}
\]
for \(x=(x^{(0)},\dots,x^{(\dims-1)})\in \real^\dims\),
where \((Z_{i_1\dots i_p})_{i_1\dots i_p}\) are identically distributed
standard normal random variables. The covariance of a \(p\)-spin Hamiltonian
\(\hamiltonian_{p, \dims}\) is then straightforward to calculate and given by
\[
	\C_{\hamiltonian_{p,\dims}}(x,y)
	= \Cov(\hamiltonian_{p,\dims}(x), \hamiltonian_{p,\dims}(y))
	= \dims\langle x,y \rangle^p.
\]
The general mixed spin Hamiltonian \(\hamiltonian_\dims := \sum_{p=0}^\infty c_p
\hamiltonian_{p,\dims}\) then has the covariance
\begin{equation}
	\label{eq: mixed spin Hamiltonian}
	\C_{\hamiltonian_\dims}(x,y) = \dims \xi\bigl(\langle x, y\rangle\bigr),
	\quad\text{with}\quad
	\xi(s) = \sum_{p=0}^\infty c_p^2 s^p.
\end{equation}
While the covariance of the Hamiltonian \(\hamiltonian_\dims\) scales with
\(\dims\), the object under consideration is typically
\(\rf_\dims = \frac{\hamiltonian_\dims}\dims\). The covariance of \(\rf_\dims\) therefore
scales with \(\frac1\dims\), which is the scaling we chose in Definition \ref{def: isotropic gaussian random function}.
Setting
\[
	\mu \equiv 0 \quad\text{and}\quad
	\kernel(\lambda_1,\lambda_2, \lambda_3) := \xi(\lambda_3)
\]
reveals that spin glasses are a subset of (non-stationary) isotropic Gaussian random
functions.
While the definition of spin glasses on \(\real^\dims\) causes no issues, they are
typically restricted to the sphere. Since the first two inputs to \(\kernel\)
are the lengths of \(x\in \sphere\), these are constant while the
third input, the scalar product, can be reconstructed from the distance.
The (non-stationary) isotropic and stationary isotropic functions therefore coincide on the
sphere. By \citet{schoenbergPositiveDefiniteFunctions1942} it follows that spin
glasses actually exhaust all possible isotropic covariance functions on the
sphere \(\sphere\) which are valid in all dimensions \(\dims\) (cf.
Remark \ref{rem: embedding}).
\smallskip

Being the end goal of optimization, the analysis of the global optimum is
complementary to our analysis of the optimization path.
The Parisi formula \autocite{parisiSequenceApproximatedSolutions1980}
shows that the global optimum is also asymptotically deterministic
\begin{equation}
    \label{eq: parisi formula}
    \lim_{\dims\to \infty}\; \max_{\param\in \sphere} \rf_\dims(x)
    = \lim_{\dims\to\infty} \E\Bigl[\; \max_{\param\in\sphere} \rf_\dims(x)\Bigr]
    = \inf_\gamma\mathcal{P}(\gamma),
\end{equation}
where \(\mathcal{P}\) is the Parisi functional. 
The first equation in \eqref{eq: parisi formula}
follows directly from the Borel-TIS inequality
\autocite[e.g.][Theorem~2.1.1]{adlerRandomFieldsGeometry2007}. 
A formal proof of the Parisi formula
can be found in the book by \citet{panchenkoSherringtonKirkpatrickModel2013}.
In subsequent papers the proof was modified to be more
constructive
\autocite{subagFreeEnergyLandscapes2018,huangConstructiveProofSpherical2024}. This
allowed for the identification of an achievable algorithmic barrier
\autocite{subagFollowingGroundStates2021,sellkeOptimizingMeanField2024},
which no Lipschitz continuous algorithms can pass \autocite{huangTightLipschitzHardness2022}
\[
	\text{ALG} = \int_0^1 \sqrt{\xi''(s)}ds.
\]
Recall that \(\xi\) is the function from \eqref{eq: mixed spin Hamiltonian} and 
note that the algorithmic barrier assumes a restriction to the sphere. A recent
review can be found in
\citet{auffingerOptimizationRandomHighDimensional2023}, see also the lecture
notes of \citet{sellkeSTAT291Random2024}.

So far, the algorithms used to show this barrier to be achievable have not
been practical for optimization. But as we approach the topic from the viewpoint of machine learning
with plausible first order optimizers,
it might be possible to obtain sufficiently tight bounds on the limiting 
function values \(\limf_\timestep\) in order to prove convergence 
towards this algorithmic barrier. It is perhaps even possible to identify
convergence rates and optimize these rates over the chosen optimizer.\smallskip

Finally, building a bridge from the spin glass viewpoint towards machine learning,
\citet{choromanskaLossSurfacesMultilayer2015} show, under strong
simplifying assumptions, that the loss function of \(p\)-layer neural networks
can be related to \(p\)-spin Hamiltonians. While they restricted the
parameters of the networks to the sphere, it is likely that much of the spin glass
theory could be extended from the sphere to (non-stationary) isotropic random functions on
\(\real^\dims\).

 \section{Proof of Theorem~\ref{THM: ASYMPTOTICALLY DETERMINISTIC BEHAVIOR}}\label{sec:proof}
We now present the proofs of our main result. First, we provide a sketch
to motivate the overall picture (Section~\ref{sec:sketch}). After we explain how the covariance of derivatives are determined
(Section~\ref{sec: covariance of derivatives}), we reforge
Theorem~\ref{THM: ASYMPTOTICALLY DETERMINISTIC BEHAVIOR} into an even more general
Theorem~\ref{THM: ASYMPTOTICALLY DETERMINISTIC BEHAVIOR VARIANT}.
This theorem is more natural to prove but requires the definition of a special
orthonormal coordinate system (Definition~\ref{def: previsible orthonormal
coordinate system}).  Finally, we give the main proof of 
Theorem~\ref{THM: ASYMPTOTICALLY DETERMINISTIC BEHAVIOR VARIANT}. Proofs for some instrumental results are deferred to the
appendix.

\subsection{Sketch of the proof of Theorem~\ref{THM: ASYMPTOTICALLY DETERMINISTIC BEHAVIOR}}
\label{sec:sketch}

To better guide the reader, we first sketch the idea behind the proof. For
simplicity, we will assume the random function to be centered and stationary isotropic
\(\rf_\dims \sim \normal(0, \ikernel)\).

\subsubsection*{Idea~1: Induction over Gaussian Conditionals}

We interpret the objects of interest
\(\rg = (\rf_\dims, \nabla \rf_\dims)\) (the heights joined with the gradients) as an
auxiliary discrete-time stochastic process \((\rg(\Param_\timestep))_{\timestep\in\nat_0}\) and
prove the claimed convergence in \(\dims\) by induction over the steps \(\timestep\).
Please pretend for now that the inputs \(\Param_\timestep\) were deterministic \(\param_\timestep\).
We will address this problem in Idea~3 of the sketch. The process
\((\rg(\param_\timestep))_{\timestep\in \nat_0}\) is then a Gaussian process and
\[
    \rg(\param_{[0:\timestep]}) := (\rg(\param_0), \dots, \rg(\param_\timestep))
\]
is therefore a multivariate Gaussian vector. 
This might be a good time to recall our notation for discrete ranges
\begin{equation}
    \label{eq: discrete interval} 
    [n\inter m] := [n,m] \cap \integer,
    \qquad
    [n\inter m) := [n,m) \cap \integer,
    \qquad
    \text{etc.}
\end{equation}
as we will make use of them throughout the following sections. It is well known
that the conditional distribution \(\rg(\param_\timestep)\mid
\rg(\param_{[0:\timestep)})\) is then also normal distributed (cf.\
Theorem~\ref{thm: conditional gaussian distribution}).
By the induction hypothesis \(\rg(\param_{[0:\timestep)})\) is already converging in
probability to something deterministic, so it is perhaps
natural to decompose \(\rg(\param_\timestep)\) into
\begin{equation}
    \label{eq: the induction engine}
    \rg(\param_\timestep) = \E[\rg(\param_\timestep) \mid \rg(\param_{[0:\timestep)})]
    + \sqrt{\Cov[\rg(\param_\timestep)\mid \rg(\param_{[0:\timestep)})]}\begin{pmatrix}
        Y_0
        \\
        \vdots
        \\
        Y_{\dims}
    \end{pmatrix}
\end{equation}
for independent standard normal distributed \(Y_i\) independent of \(\rg(x_{[0:\timestep)})\).
The conditional expectation is then given (cf.\ Theorem~\ref{thm: conditional gaussian distribution}) by
\[
	\E[\rg(\param_\timestep) \mid \rg(\param_{[0:\timestep)})]
	= \blue{\Cov(\rg(\param_{[0:\timestep)}), \rg(\param_\timestep))}^T[\magenta{\Cov(\rg(\param_{[0:\timestep)}))}]^{-1} \rg(\param_{[0:\timestep)})
\]
and conditional variance by
\begin{align*}
	&\Cov[\rg(\param_\timestep) \mid \rg(\param_{[0:\timestep)})]
	\\
	&= \green{\Cov[\rg(\param_\timestep)]} - \blue{\Cov(\rg(\param_{[0:\timestep)}), \rg(\param_\timestep))}^T[\magenta{\Cov(\rg(\param_{[0:\timestep)}))}]^{-1}\blue{\Cov(\rg(\param_{[0:\timestep)}), \rg(\param_\timestep))}.
\end{align*}
The vector \(\rg(x_{[0:\timestep)})\) already converges by induction as the dimension
\(\dims\) increases. What is therefore left to study are the covariance matrices.

\subsubsection*{Idea~2: Finding sparsity in the covariance matrices with a custom coordinate system}

The covariance matrices which make up the conditional expectation and variance
are increasing in size as the dimension \(\dims\) grows, because the gradient
\(\nabla\rf_\dims\) contained in \(\rg\) increases in size. Additionally we somehow need
to group the independent \(Y_i\) into something which converges.
It turns out that the standard coordinate system is ill suited to achieve these
goals and get sufficiently sparse (and therefore manageable) matrices.
For this reason we will now phrase everything in terms of directional derivatives.

In Lemma~\ref{lem: cov of derivatives, isotropy} we show how the covariances of the
directional derivatives of an isotropic \(\rf_\dims\sim\normal(0, \ikernel)\) can be
calculated explicitly. In particular we have for vectors \(v,w\in \real^\dims\) and points
\(x,y\in \real^\dims\) with distance \(\distance = x-y\) 
\begin{equation}
    \label{eq: the two sparsity requirements}
    \Cov(D_v\rf_\dims(x), D_w\rf_\dims(y))
    = - \tfrac1\dims\Bigl[
        \underbrace{\ikernel''\bigl(\tfrac{\|\distance\|^2}2\bigr)\langle \distance, w\rangle\langle \distance, v\rangle}_{\text{(I)}}
        + \underbrace{\ikernel'\bigl(\tfrac{\|\distance\|^2}2\bigr)\langle w, v\rangle}_{\text{(II)}}
    \Bigr].
\end{equation}
We are now going to first explain how an orthogonal basis keeps (II) sparse
before we explain where the use of the standard basis breaks down for (I).
For simplicity, let us forget the height $\rf_\dims$ for now and pretend \(\rg\) consists only of the
gradient, i.e. \(\rg=\nabla\rf_\dims\), and first consider the covariance matrix at a fixed point 
\(\nabla\rf_\dims(x_0)\) only. Since we only consider the derivatives
at a single point we always have \(x= y = x_0\) so the distance \(\distance\)
vanishes. This completely removes the first part (I) (which will later require a
special basis). The second part (II) is zero if \(v\) and \(w\) are orthogonal.
Assuming we use an orthonormal coordinate system (e.g. the standard basis), we
therefore get
\[
    \green{\Cov[\rg(x_0)]}= \Cov[ \nabla\rf_\dims(x_0) ] = \frac1\dims\begin{pmatrix}
        -\ikernel'(0) &
        \\
        & \ddots
        \\
        & & -\ikernel'(0)
    \end{pmatrix}
    = \tfrac{-\ikernel'(0)}\dims \identity.
\]
As \(\rg(\param_{[0:0)})\) is the empty set we therefore have in the first step of the induction in $n$ 
(cf.~Equation~\eqref{eq: the induction engine})
\begin{align}
    \nonumber
    \rg(\param_0)
    &= \E[\rg(\param_0) \mid \rg(\param_{[0:0)})] 
    + \sqrt{\Cov[\rg(\param_0)\mid \rg(\param_{[0:0)})]}\begin{pmatrix}
        Y_0
        \\
        \vdots
        \\
        Y_{\dims}
    \end{pmatrix}
    \\
    \label{eq: first law of large numbers}
    &= \underbrace{\E[\rg(\param_0)]}_{=0} + \sqrt{\green{\Cov[\rg(x_0)]}} \begin{pmatrix}
        Y_0
        \\
        \vdots
        \\
        Y_{\dims}
    \end{pmatrix}
    \\
    \nonumber
    &= \sqrt{\tfrac{-\ikernel'(0)}\dims}\begin{pmatrix}
        Y_0
        \\
        \vdots
        \\
        Y_{\dims}
    \end{pmatrix}.
\end{align}
In particular, we have by the law of large numbers,
\[
    \langle \nabla \rf_\dims(x_0), \nabla\rf_\dims(x_0) \rangle
    = \frac{-\ikernel'(0)}{\dims} \sum_{i=1}^\dims Y_i^2 \overset{p}{\underset{\dims\to \infty}\to} -\ikernel'(0).
\]
For the induction start and (II) there is therefore no need to change the coordinate system.
But for the induction step with \(\timestep>1\) the situation in (I) becomes more delicate.\smallskip

We now want to understand where the issues in (I) with the standard coordinate
system come form. Recall that we also want to determine
limiting values of \(\langle \nabla \rf_\dims(x_0), \nabla\rf_\dims(x_\red{1}) \rangle\).

\textbf{Problem:} Since \(x_1\) will be \(x_0\) plus a gradient step,
\(\distance=x_1 - x_0\) is therefore a multiple of the gradient. If we continue to use
the standard coordinate system, (I) causes the covariance matrix to be dense. This is because the
entries of \(\distance\) are proportional to \(\partial_i \rf_\dims(x_0)\), which are
multiples of the \(Y_i\) and therefore almost surely never zero.

\textbf{Solution:} We use an adapted coordinate system. For this we select
\(\tilde{\rv}_0 := \nabla \rf_\dims(\param_0)\) and normalize this vector \(\rv_0 =
\frac{\tilde{\rv}_0}{\|\tilde{\rv}_0\|}\). Again, we are going to pretend that this vector
was deterministic and mark this fact with the notation \(v_0\).
We then extend it by \(w_{[1:\dims)}\) to an orthonormal basis. Since the \(w_i\)
are orthogonal to the gradient span, they are orthogonal to the distances
\(\param_0 - \param_1\), i.e. \(\langle \distance, w_i\rangle = 0\) and therefore
(I) is zero for almost all directions (except for \(v_0\)). As the coordinate
system is orthonormal, (II) is sparse again. More generally, we chose a coordinate system
\(v_{[0:\timestep)}\) capturing the span of the first \(\timestep\) gradients \(\nabla\rf_\dims(\param_{[0:\timestep)})\) and
therefore also the span of the first \(\timestep\) parameters and extend this coordinate
system by \(w_{[\timestep:\dims)}\) to an orthonormal coordinate system such that
we have essentially
\begin{align*}
    \rg(\param_{[0:\timestep)})
    &= (\nabla\rf_\dims(\param_0),\dots, \nabla\rf_\dims(\param_{\timestep-1}))
    \\
    &= U\begin{pmatrix}
        D_{v_0}\rf_\dims(\param_0) & \cdots & D_{v_0}\rf_\dims(\param_{\timestep-1})
        \\
        \vdots & & \vdots
        \\
        D_{v_{\timestep-1}}\rf_\dims(\param_0) & \cdots & D_{v_{\timestep-1}}\rf_\dims(\param_{\timestep-1})
        \\
        D_{w_\timestep}\rf_\dims(\param_0) & \cdots & D_{w_\timestep}\rf_\dims(\param_{\timestep-1})
        \\
        \vdots & & \vdots
        \\
        D_{w_{\dims-1}}\rf_\dims(\param_0) & \cdots & D_{w_{\dims-1}}\rf_\dims(\param_{\timestep-1})
    \end{pmatrix}U^T.
\end{align*}
Here \(U\) is a basis change matrix which we are going to suppress for simplicity in the following.
Note that the covariance matrix is dense for the vectors \(v_{[0:\timestep)}\) because they span the subspace
of the evaluation points \(x_k\) and therefore cause (I) to be non-zero. But (I)
remains zero for the vectors \(w_i\) and the covariance matrix is therefore sparse there.
We want to group the chaos together (i.e. group the directional
derivatives of \(v_{[0:\timestep)}\)), so we view \(\rg(\param_{[0:\timestep)})\)
as a row-major matrix, i.e.  whenever we treat \(\rg(\param_{[0:\timestep)})\) as a vector we
concatenate the rows. It then turns out that the covariance matrix is of the
form
\[
    \magenta{\Cov[\rg(\param_{[0:\timestep)})]}
    = \left(\begin{array}{c c c c c}
        \cline{1-1}
        \multicolumn{1}{|c|}{}
        \\
        \multicolumn{1}{|c|}{\Cov[D_{v_{[0:\timestep)}\rf_\dims(\param_{[0:\timestep)})}]} 
        \\
        \multicolumn{1}{|c|}{}
        \\
        \cline{1-2}
        &\multicolumn{1}{|c|}{\Sigma_\timestep}
        \\
        \cline{2-3}
        & &\multicolumn{1}{|c|}{\quad}
        \\
        \cline{3-3}
        & & & \ddots
        \\
        \cline{5-5}
        & & & &\multicolumn{1}{|c|}{\Sigma_{\dims-1}}
        \\
        \cline{5-5}
    \end{array}\right),
\]
where \(\Sigma_i = \Cov[D_{w_i}\rf_\dims(\param_{[0:\timestep)})]\) are \(\timestep\times \timestep\) matrices
and \(\Cov[D_{[0:\timestep)}\rf_\dims(\param_{[0:\timestep)})]\) is a dense \(\timestep^2\times \timestep^2\) block matrix.
We can assume \(\dims-1 > \timestep\) without loss of generality as we let \(\dims \to
\infty\).
In the proof we will make the observation that all the
\(\Sigma_i\) are in fact equal.
Note that the mixed covariance matrix \(\blue{\Cov(\rg(\param_{[0:\timestep)}),
\rg(\param_\timestep))}\) used in the conditional expectation and covariance have
similar block form and \(\green{\Cov[\rg(\param_\timestep)]}\)
is always diagonal.
The key is now to understand that
as \(\dims\) grows to infinity, the scaling \(\frac1\dims\) keeps the \(\timestep^2\times
\timestep^2\) chaos in check while we get a similar law of large numbers as in the first
step (cf.~\eqref{eq: first law of large numbers}) by grouping the identical
\(\Sigma_i\) to get the length of
the new gradient \(\nabla\rf_\dims(x_\timestep)\) projected to
the subspace of \(w_{[\timestep:\dims)}\) which is the orthogonal space of \(v_{[0:\timestep)}\).
This projection makes up the new coordinate \(v_\timestep\) and the induction continues.

\subsubsection*{Idea~3: Treating the random input with care}

At the beginning of the sketch we asked the reader to ignore the fact that the
evaluation points \((\Param_\timestep)_{\timestep\in \nat_0}\) are random themselves. Moreover,
we ignored that the adapted orthogonal basis \(\rv_{[0:\timestep)}\) was constructed from
random gradients \(\nabla\rf_\dims(\Param_k)\) which implies it is random as well.
Instead we denoted it by \(v_{[0:\timestep)}\) and pretended it was deterministic.
In the following, we sketch how to get rid of the randomness of
\(\Param_\timestep\) and the randomness of the  coordinate system \(\rv_{[0:\timestep)}\). It is
perhaps most obvious for the random coordinate system, that
\[
    D_{\rv_0}\rf_\dims(\param_0) = \langle \rv_0, \nabla\rf_\dims(x_0)\rangle = \|\nabla\rf_\dims(x_0)\|^2
\]
is certainly not Gaussian by construction (recall \(\rv_0 = \nabla\rf_\dims(x_0)\)).
Similarly, the random points \(\Param_\timestep\) will typically change the distribution
of \(\rf_\dims(\Param_\timestep)\) and taking (conditional)
expectations and covariances becomes quite delicate. The key to solve this issue
is Corollary \ref{cor: gaussian previsible sampling} which (assuming the choice
of a consistent joint conditional distribution) states that random inputs can be
treated as deterministic via
\[
    \E[h(\rg(\Param_\timestep))\mid \filt_{\timestep-1}]
    = \Bigl(
        (x_{[0:\timestep]})\mapsto \E[h(\rg(x_\timestep)) \mid \rg(x_{[0:\timestep)})]
    \Bigr)(X_{[0:\timestep]}), 
\]
if \(\Param_\timestep\) is previsible, i.e. \(\Param_{\timestep+1}\) is measurable with
regard to the filtration \(\filt_\timestep = \sigma(\rg(X_k): k \le \timestep)\). Applying this
result to the evaluations points
\(\Param_\timestep\) is relatively straightforward. But it is not obvious how to treat the
random coordinate system. The key is to turn the coordinate system into an
input. The definition
\[
    Z(v; \param) := D_v \rf_\dims(\param)
\]
allows us to apply Corollary \ref{cor: gaussian previsible sampling} to both coordinate
system and points. The delicate part is that the input must be previsible with
respect to the filtration \(\filt_\timestep\). So now the coordinate system must also be
previsible. We therefore cannot use future gradients for our coordinate
system. This forces us to define a custom coordinate system for every step. That
means after \(\timestep\) steps we have \(\timestep\) different coordinate
systems. We will carefully reduce this to just one coordinate system per
iteration in the formal proof. Turning this strategy into a rigorous proof is a
bit tedious and requires a careful induction including a few additional
technical induction hypothesis.
 
\subsection{Proof of Theorem~\ref{THM: ASYMPTOTICALLY DETERMINISTIC BEHAVIOR}}\label{sec:main_proof}

After the preparations of the previous sections we can now give the details of
the proof sketched in Section \ref{sec:sketch}. In the sketch of the proof,
we highlighted the importance of an adapted coordinate system in order to
keep the covariance matrices sparse. In particular, part (I) of \eqref{eq: the
two sparsity requirements} requires orthogonality to the walking directions.
We also needed the coordinate system to be orthogonal to keep (II) as sparse as
possible. In Idea 3, we sketched why the coordinate system had to be previsible.
That is, why we needed a different coordinate system in every step.
We therefore begin the proof by building these coordinate systems.\smallskip

Let us start by defining the natural filtration
\[
    \filt_\timestep
    := \sigma\bigl( \rf_\dims(\Param_k), \nabla\rf_\dims(\Param_k) : 0\le k \le \timestep\bigr).
\]
The following previsible \emph{vector space of evaluation points}
\begin{equation}
    \label{eq: vector space of evaluation points}    
    V_\timestep
    := \Span\Bigl(
        \{\param_0\} \cup \{\nabla\rf_\dims(\Param_k) : 0\le k \mathrel{\red{<}}\timestep\}
    \Bigr),
    \qquad \dimV_\timestep := \dim(V_\timestep),
\end{equation}
contains \(\Param_{[0:\timestep]}\) by definition of the generalized gradient
span algorithm (Definition~\ref{def: general gsa}) and is similarly measurable
with respect to \(\filt_{\timestep-1}\), i.e. previsible.
Also recall that \([0\inter\timestep]=[0, \timestep]\cap \integer\) is notation
for a discrete intervals as (re-)introduced in Equation \ref{eq: discrete interval}.
We will now use this chain of vector spaces to define previsible coordinate
systems for every step \(\timestep\).

\begin{definition}[Previsible orthonormal coordinate systems]
    \label{def: previsible orthonormal coordinate system}
    We define an \(\filt_{\timestep-1}\)-measurable orthonormal basis
    \(\rv_{[0:\dimV_\timestep)}\) of \(V_\timestep\) inductively such that
    \(\rv_{[0:\dimV_k)}\) is a basis of the subspace \(V_k\).
    
    \paragraph*{Induction start} Assuming \(\param_0\neq 0\) and thus \(\dimV_0 =
    1\), we define
    \[
        \rv_0 :=\frac{\param_0}{\|\param_0\|}.
    \]
    In any case this results in a basis \(\rv_{[0:\dimV_0)}\) of \(V_0\)
    since \([0:\dimV_0)=\emptyset\) if \(\dimV_0 = 0\).

    \paragraph*{Induction step (Gram-Schmidt)}
    Assuming we have a basis \(\rv_{[0:\dimV_\timestep)}\) of \(V_\timestep\), let
    us construct a basis for \(V_{\timestep+1}\). For this we define the
    following Gram-Schmidt procedure candidate
    \begin{equation}
        \label{eq: definition of v candidate}
        \tilde{\rv}_{\timestep}
        := \nabla\rf_\dims(\Param_\timestep) - P_{V_\timestep}\nabla\rf_\dims(\Param_\timestep)
        = \nabla\rf_\dims(\Param_\timestep) - \sum_{i < \dimV_\timestep} \langle \nabla \rf_\dims(\Param_\timestep), \rv_i\rangle \rv_i,
    \end{equation}
    where \(P_{V_\timestep}\) is the projection to \(V_\timestep\).
    If \(\tilde{\rv}_{\timestep} = 0\), then \(\nabla\rf_\dims(\Param_\timestep) \in
    V_\timestep\) and we thus have \(V_{\timestep+1} = V_\timestep\) by definition
    \eqref{eq: vector space of evaluation points}. In that
    case we already have a basis for \(V_{\timestep+1}\).

    For any \(\timestep\) where the dimension increases such that
    \(\dimV_{\timestep+1} = \dimV_\timestep + 1\), we define
    \begin{equation}
        \label{eq: definition of v}
        \rv_{\dimV_{\timestep+1}-1}
        =\rv_{\dimV_\timestep}
        := \frac{\tilde{\rv}_{\timestep}}{\|\tilde{\rv}_{\timestep}\|}.
    \end{equation}
    We thus have obtained an orthonormal basis \(\rv_{[0:\dimV_{\timestep+1})}\)
    of \(V_{\timestep+1}\) which is \(\filt_\timestep\) measurable.
    
    \paragraph*{Basis extensions} With this construction of basis elements of
    \(V_\timestep\) done, we \(\filt_{\timestep-1}\)-measurably select an
    arbitrary orthonormal basis \(\rw^{(\timestep)}_{[\dimV_\timestep:\dims)}\)
    of \(V_\timestep^\perp\) to obtain an \(\filt_{\timestep-1}\)-measurable
    orthonormal basis
    \[
        B_\timestep := (\rv_{[0:\dimV_\timestep)}, \rw^{(\timestep)}_{[\dimV_\timestep:\dims)})
    \]
    for every \(\timestep\in\nat\). The coordinate systems \(B_\timestep\) are
    thus previsible.
\end{definition}

Using this specialized coordinate system we state an extension of
Theorem \ref{THM: ASYMPTOTICALLY DETERMINISTIC BEHAVIOR}. This extension looks less friendly but is more
natural to prove. It implies Theorem~\ref{THM: ASYMPTOTICALLY DETERMINISTIC
BEHAVIOR}, proves three additional claims and relaxes the assumptions on
the prefactors of the gradient span algorithm. The additional claims cannot
be stated separately as they are all proved in one laborious induction. And to
prove the claim we are most interested in \ref{ind: information convergence}
we also require the other claims in the induction step.
Finally, with the help of Proposition~\ref{prop: wlog deterministic starting
points}, we can (and will) assume deterministic starting points in
Theorem~\ref{THM: ASYMPTOTICALLY DETERMINISTIC BEHAVIOR VARIANT}
without loss of generality.

\begin{theorem}[Predictable progress of gradient span algorithms {[Extension of Theorem~\ref{THM: ASYMPTOTICALLY DETERMINISTIC BEHAVIOR}]}]
	\label{THM: ASYMPTOTICALLY DETERMINISTIC BEHAVIOR VARIANT}
    Let \(\kernel\) be a kernel valid in all dimensions (\ie defined on \(\ell^2\)
    by Lemma \ref{lem: valid in all dimensions}) and  \(\rf_\dims \sim \normal(\mu,
    \kernel)\) be a sequence of scaled (non-stationary) isotropic Gaussian random
	functions (Definition~\ref{def: isotropic gaussian random
	function}) in $\dims$. 
    Assume that \(\mu\) and \(\kernel\)  are
	sufficiently smooth (Assumption~\ref{assmpt: smoothness}) and assume the
	covariance of \((\rf_\dims,\nabla\rf_\dims)\) is strictly positive definite.
	Let \(\gsa\) be a general gradient span algorithm (Definition~\ref{def:
	general gsa}), which is \emph{asymptotically} continuous and uses
    the most recent gradient \emph{asymptotically} (cf.~Assumption~\ref{assmpt: generalized assumptions}).

	Let \(\gsa\) be applied to \(\rf_\dims\) with starting points
	\(\param_0\in \real^\dims\) such that \(\Param_\timestep = \gsa(\rf_\dims,
	\param_0, \timestep)\). We assume that the deterministic initialization point
	\(\param_0\in \real^\dims\) is of constant length \(\|\param_0\|=\radius\)
    over the dimension \(\dims\).
	Using \(\gradients_\timestep := \bigl(\nabla \rf_\dims(\Param_k): k \le
	\timestep\bigr)\), we further define the modified random information
	vector
    \[
		\VarInfo_\timestep
		:= \Bigl(\rf_\dims(\Param_k): k \le \timestep\Bigr)
		\cup \Bigl(
			\langle v, w\rangle :
			v \in \rv_{[0:\dimV_{\timestep+1})}, w \in \gradients_\timestep
		\Bigr).
    \]
    The following inductive claims hold for all \(\timestep\in \nat\), where
    \ref{ind: full dimension} implies that the vector \(\VarInfo_\timestep\) is
    almost surely of constant length.
    \begin{enumerate}[label=(Ind-\Roman*), wide=0pt,leftmargin=\parindent]
        \item\label{ind: information convergence}
        \textbf{Information convergence:}
        There exists some deterministic limiting information vector
        \(\VarLiminfo_\timestep=\VarLiminfo_\timestep(\gsa, \mu, \kernel)\),
        such that
        \[
            \VarInfo_\timestep
            \underset{\dims\to\infty}{\overset{p}\to}
            \VarLiminfo_\timestep.
        \]
        The limiting information is split into \(\VarLiminfo_\timestep =
        (\limf_k, \gamma_k^{(i)})_{k,i}\),
        where the limiting elements \(\limf_k\) of \(\rf_\dims(\Param_k)\) were already
        defined in Theorem~\ref{THM: ASYMPTOTICALLY DETERMINISTIC BEHAVIOR}.
        We further define the limiting inner products of gradients by
        \begin{equation}
            \label{eq: definition gamma}
            \gamma^{(i)}_k := \lim_{\dims\to\infty}
            \langle \nabla\rf_\dims(X_k), \rv_i\rangle.
        \end{equation}

        \item\label{ind: full dimension} \textbf{(Asymptotic) Full rank}:
        If the most recent gradient is always used (and not just asymptotically),
        then the previsible vector space of evaluation points \(V_{\timestep+1}\) has almost
        surely full rank (assuming \(\param_0\neq 0\)). Specifically, for all \(m\le \timestep+1\)
        we have almost surely
        \[
            \dimV_m = m + \ind_{\param_0 \neq 0}.
        \]
        This always holds in the limit (even when the most recent gradient is only
        used asymptotically), which means that for all \(k\le \timestep\) the
        Gram-Schmidt candidate \(\tilde{\rv}_{k}\) defined in \eqref{eq:
        definition of v} is not zero in the limit
        \begin{equation}
            \label{eq: limiting full rank}
            \gamma_k^{(\dimV_k)} = \lim_{\dims\to\infty}\langle \nabla\rf_\dims(\Param_k), \rv_{\dimV_k}\rangle = \lim_{\dims\to\infty}\|\tilde{\rv}_{k}\| \neq 0.
        \end{equation}

        \item\label{ind: representation}
        \textbf{Representation:} For all \(k\le m\le \timestep +1\) there exist
        limiting representation vectors
        \[
            y_k = y_k(\gsa, \mu, \kernel)
            \quad\text{with}\quad 
            y_k
            = (y_k^{(0)}, \dots y^{(\dimV_m-1)})\in \real^{\dimV_m}
        \]
        of \(\Param_k\), such that for all \(i < \dimV_m\) we have
        \begin{equation}
            \label{eq: <X, v> converges to gamma}
            \langle \Param_k, \rv_i\rangle
            \underset{\dims\to\infty}{\overset{p}\to} y_k^{(i)}
            \quad\text{and}\quad
            \|\Param_k\|^2
            = \sum_{i=0}^{\dimV_m-1} \langle \Param_k, \rv_i\rangle^2
            \underset{\dims\to\infty}{\overset{p}\to}
            \|y_k\|^2.
        \end{equation}
        The asymptotic distances are then given for all \(k,l\le m\) by
        \begin{equation}
            \label{eq: limiting distances}
            \|\Param_k - \Param_l\|^2
            = \sum_{i=0}^{\dimV_m-1} \langle \Param_k - \Param_l, \rv_i\rangle^2
            \underset{\dims\to\infty}{\overset{p}\to} 
            \; \rho_{kl}^2
            \;:= \|y_k - y_l\|^2.
        \end{equation}

        \item\label{ind: evaluation points asymptotically different}
        \textbf{The evaluation points are asymptotically different}, i.e.
        \[
            \rho_{kl}^2>0 \quad \text{for all}\quad k,l\leq \timestep+1 \quad \text{with}\quad k\neq l.
        \]
    \end{enumerate}
\end{theorem}

Before we prove that Theorem~\ref{THM: ASYMPTOTICALLY DETERMINISTIC BEHAVIOR} follows from 
Theorem~\ref{THM: ASYMPTOTICALLY DETERMINISTIC BEHAVIOR VARIANT} let
us state the `asymptotic continuity' and `use of the most recent gradient' assumptions on the prefactors.

\begin{assumption}[Assumptions on prefactors]
    \label{assmpt: generalized assumptions}
    Recall, that we defined the prefactors \(\lr_{\timestep}\) of a GSA
    to be functions of the information \(\Info_{\timestep-1}\) (Definition~\ref{def: general gsa}).
    We assume that
    \begin{itemize}
		\item \(\lr_\timestep\) is continuous in the point \(\liminfo_{\timestep-1}\) for all \(\timestep\), and  
		\item the most recent gradients are used at least in the asymptotic
		limit, i.e.
		\[		
			\hat{\lr}^{(g)}_{\timestep,\timestep-1}
			:=\lr^{(g)}_{\timestep,\timestep-1}(\liminfo_{\timestep-1}) \Bigl(= \lim_{\dims\to\infty}\lr^{(g)}_{\timestep,\timestep-1}(\Info_{\timestep-1})\Bigr) \neq 0.
		\]
	\end{itemize}
    By Lemma~\ref{lem: info is continuous function} it is apparent, 
    that we can equivalently assume the prefactors \(\lr_\timestep\) to be functions in
    \(\VarInfo_{\timestep-1}\) continuous in \(\VarLiminfo_{\timestep-1}\),
    which use the most recent gradient in the asymptotic limit \(\hat{\lr}_{\timestep,\timestep-1} = \lr_{\timestep,\timestep-1}(\VarLiminfo_{\timestep-1})\).
\end{assumption}

At first it might seem circular to make use of limiting elements in an assumption
which is necessary to prove the limiting elements exist.
But a closer inspection reveals that \(\lr_\timestep\) is only used
for the convergence of \(\Info_{k}\) to \(\liminfo_{k}\) for times \(k\ge \timestep\).
We can therefore interleave this assumption on \(\lr_\timestep\) with the
inductive existence proof of \(\liminfo_{\timestep-1}\).

The proof of Theorem~\ref{THM: ASYMPTOTICALLY DETERMINISTIC BEHAVIOR} now
follows readily from Theorem~\ref{THM: ASYMPTOTICALLY DETERMINISTIC BEHAVIOR
VARIANT} via a couple of preliminary results:

\begin{prop}\label{prop: wlog deterministic starting points}
    We can assume without loss of generality that the independent initialization
    \(\Param_0\) is a deterministic \(\param_0\in \real^\dims\) of length
    \(\|\param_0\|=\radius\) over all dimensions \(\dims\).

    \textbf{Stationary case:} Assume the random functions \(\rf_\dims\) are
    furthermore stationary isotropic and the algorithm \(\param_0\)-agnostic
    (Definition~\ref{def: general gsa}). Then the limiting information is
    independent of \(\radius\) and the random
    initialization \(\Param_0\) does not need to satisfy
    \(\|\Param_0\|=\radius\) almost surely.
\end{prop}

The proof of this Proposition is accomplished via two lemmas:

\begin{lemma}[The `information' is invariant to linear isometries]
    \label{lem: information invariant}
    Let \(f\) be some function and \((x_0, \dots, x_n)\) evaluation points.
    We further define a change in coordinates
    \[
        g(y) := f \circ \phi
        \quad \text{and}\quad
        y_k := \phi^{-1}(x_k) \quad\text{for}\quad k\in \{0,\dots, n\}
    \]
    via a linear isometry \(\phi(x) = Ux\) for some orthonormal matrix \(U\). Then the information
    \[\begin{aligned}
        \info_\timestep
        &:= \Bigl(f(\param_k): k \le \timestep\Bigr)
        \cup \Bigl(
            \langle v, w\rangle :
            v,w \in (\param_0) \cup \gradients_\timestep
        \Bigr)
        \quad\text{with}
        \\
        \gradients_\timestep &:= \bigl(\nabla f(\param_k): k \le \timestep\bigr)
    \end{aligned}\]
    is exactly equal to the information
    \[\begin{aligned}
        \varinfo_\timestep
        &:= \Bigl(g(y_k): k \le \timestep\Bigr)
        \cup \Bigl(
            \langle v, w\rangle :
            v,w \in (y_0) \cup \tilde{\gradients}_\timestep
        \Bigr)
        \quad\text{with}
        \\
        \tilde{\gradients}_\timestep &:= \bigl(\nabla g(y_k): k \le \timestep\bigr).
    \end{aligned}\]
    If \(\phi\) is a general isometry of the form \(\phi(x) = Ux + b\), then we
    still have equality for the the reduced information, i.e.
    \[
        \info^{\setminus x_0}_\timestep = \varinfo^{\setminus y_0}_\timestep.
    \]
\end{lemma}

\begin{proof}
    We have by definition
    \[
        f(\param_k)
        = f\circ \phi(\phi^{-1}(\param_k))
        = g(y_k).
    \]
    Let us therefore turn to the inner products. Since \(\phi(x) = Ux\) or
    \(\phi(x) = Ux + b\), we have for the gradient
    \begin{equation}
        \label{eq: basis change gradient}
        \nabla g(y) = U^\transpose \nabla f(\phi(y)).
    \end{equation}
     this implies all the inner products are equal 
    \begin{alignat*}{3}
        \langle \nabla g(y_k), \nabla g(y_l)\rangle
        &= \langle U^\transpose\nabla f(x_k), U^\transpose \nabla f(x_l)\rangle
        &&= \langle \nabla f(x_k), \nabla f(x_l)\rangle.
        \\
    \intertext{
        In the linear isometry case, where \(\phi(x) = Ux\), we have \(y_n = U^\transpose x_n\) and thus
    }
        \langle y_0, \nabla g(y_l)\rangle
        &= \langle U^\transpose x_0, U^\transpose \nabla f(x_l)\rangle
        &&= \langle x_0, \nabla f(x_l)\rangle
        \\
        \langle y_0, y_0\rangle
        &= \langle U^\transpose x_0, U^\transpose x_0\rangle
        &&= \langle x_0, x_0\rangle.
    \end{alignat*}
    The information is therefore exactly the same.
\end{proof}

\begin{lemma}[Linear isometry equivariance of general gradient span algorithms]
    \label{lem: gsa linear isometry invariant}
    Let \(\gsa\) be a general gradient span algorithm (Definition~\ref{def:
    general gsa}), let \(f\) be a function, \(x_0\) a starting point.
    We define a change of basis
    \[
        g(y) := f \circ \phi
        \quad \text{and}\quad
        y_0 := \phi^{-1}(x_0).
    \]
    via a linear isometry \(\phi(x) = Ux\) with orthonormal matrix \(U\).
    Then the optimization paths
    \[
        x_n := \gsa(f, \param_0, \timestep)
        \quad \text{and}\quad
        y_n := \gsa(g, y_0, \timestep)
    \]
    are equivariant, i.e. the simple basis change
    \[
        y_n = \phi^{-1}(x_n)
    \]
    is retained for all \(n\in \nat\).
    The same holds true for all isometries \(\phi\), if the algorithm is
    \(x_0\)-agnostic.
\end{lemma}
\begin{proof}
    We proceed by induction over \(n\), where the induction start is obvious.

    For the induction step \((\timestep-1)\to \timestep\) we use the induction
    claim \(y_k = \phi^{-1}(x_k)\) for all \(k\le \timestep-1\) to obtain that
    the information \(\info_{\timestep-1}\) is invariant (by Lemma~\ref{lem:
    information invariant}). This implies by definition of the gradient
    span algorithm and \eqref{eq: basis change gradient}
    \[
        y_\timestep
		= \lr^{(\param)}_{\timestep} y_0
		+ \sum_{k=0}^{\timestep-1}\lr_{\timestep,k}^{(g)}\nabla g(y_k)
		= U^\transpose\Bigl(\lr^{(\param)}_{\timestep} x_0
            + \sum_{k=0}^{\timestep-1}\lr_{\timestep,k}^{(g)}\nabla f(x_k)
        \Bigr)
        = \phi^{-1}(x_\timestep),
    \]
    where the invariance of the information was used implicitly as the \(\lr_\timestep\)
    are functions of \(\info_{\timestep-1}\).

    In the \(\param_0\)-agnostic case, the induction step is almost the same
    except we have for isometries \(\phi(x) = Ux + b\) that \(y_k = \phi^{-1}(x_k) = U^T (x_k-b)\)
    and thus
    \[
        y_\timestep
		= y_0
		+ \sum_{k=0}^{\timestep-1}\lr_{\timestep,k}^{(g)}\nabla g(y_k)
		= U^\transpose\Bigl( x_0 - b
            + \sum_{k=0}^{\timestep-1}\lr_{\timestep,k}^{(g)}\nabla f(x_k)
        \Bigr)
        = \phi^{-1}(x_\timestep),
    \]
    using the fact that \(\lr^{(x)}_\timestep = 1\). We similarly use that the
    reduced information \(\info^{\setminus x_0}_{\timestep-1}\) is retained for
    the prefactors.
\end{proof}

We are now ready to prove that we can assume without loss of generality that
the initialization is deterministic.

\begin{proof}[Proof of Proposition~\ref{prop: wlog deterministic starting points}]
    Since we have
    \[
        \Pr\bigl(\VarInfo_\timestep \in A\bigr)
        = \E\Bigl[\Pr\bigl(\VarInfo_\timestep \in A \mid \Param_0\bigr)\Bigr],
    \]
    it is sufficient to show that \(\Pr(\VarInfo_\timestep \in A \mid
    \Param_0=\param_0)\) is only dependent on \(\|\param_0\|=\radius\),
    which is constant over the distribution of \(\Param_0\).

    For any \(\param_0\), \(y_0\) with \(\|\param_0\|=\|y_0\|=\radius\), there
    exists a linear isometry \(\phi\) such that \(\param_0 = \phi^{-1}(y_0)\).

    Let \(\Info_\timestep = \Info_\timestep(\rf_\dims, y_0)\) be
    the information vector generated from running the gradient
    span algorithm on \(\rf_\dims\) with starting point \(y_0\) for \(\timestep\)
    steps. Since \(\rf_\dims \overset{(d)}= \rf_\dims \circ \phi\) in
    distribution, due to (non-stationary) isotropy (cf.~Definition~\ref{def: distributional input invariance}), we have
    \[
        \Info_\timestep(\rf_\dims, y_0)
        \overset{(d)}=
        \Info_\timestep(\rf_\dims\circ \phi, y_0)
        = \Info_\timestep(\rf_\dims, x_0)
    \]
    where we have used Lemma~\ref{lem: information invariant} and
    Lemma~\ref{lem: gsa linear isometry invariant} for the last equation.
    With the note that \(\VarInfo_\timestep\) is a deterministic map of
    \(\Info_\timestep\) as it only requires Gram-Schmidt orthogonalization as
    outlined in Definition~\ref{def: previsible orthonormal
    coordinate system} we can conclude this proof
    \begin{align*}
        \Pr\Bigl( \VarInfo_\timestep(\rf_\dims, \Param_0) \in A \mid \Param_0= y_0\Bigr)
        &= \Pr\Bigl( \VarInfo_\timestep(\rf_\dims, y_0) \in A\Bigr)
        \\
        &= \Pr\Bigl( \VarInfo_\timestep(\rf_\dims, x_0) \in A\Bigr)
        \\
        &= \Pr\Bigl( \VarInfo_\timestep(\rf_\dims, X_0) \in A \mid \Param_0= x_0\Bigr).
    \end{align*}

    In the stationary case with \(\param_0\)-agnostic algorithm, we can make use
    of arbitrary isometries \(\phi\). In particular we can chose \(\phi(x) = x-x_0\)
    that maps \(x_0\) to zero. With the same arguments as above (using the
    \(x_0\)-agnostic version of Lemma~\ref{lem: information invariant} and
    Lemma~\ref{lem: gsa linear isometry invariant}) we get
    \[
        \Pr\Bigl( \VarInfo_\timestep(\rf_\dims, \Param_0) \in A \mid \Param_0= x_0\Bigr)
        = \Pr\Bigl( \VarInfo_\timestep(\rf_\dims, \Param_0) \in A \mid \Param_0= 0\Bigr).
    \]
    The distribution is thus completely independent of \(\param_0\). In particular,
    this forces the limiting values to be independent of \(\param_0\) and therefore
    also independent of \(\lambda\).
\end{proof}

\begin{lemma}
    \label{lem: info is continuous function}
    For any fixed \(\param_0\),
    \(\Info_\timestep\) is a continuous function of \(\VarInfo_\timestep\) for all \(\timestep\in \nat\).
\end{lemma}
\begin{proof}
    Let us first recall the definition of \(\Info_\timestep\) and
    \(\VarInfo_\timestep\). With \(\gradients_\timestep := \bigl(\nabla \rf_\dims(\Param_k): k \le \timestep\bigr)\)
    we have in a direct comparison:
    \begin{align*}
		\Info_\timestep
		&:= \Bigl(\rf_\dims(\Param_k): k \le \timestep\Bigr)
		\cup \Bigl(
			\langle v, w\rangle :
			v,w \in (\param_0) \cup \gradients_\timestep
		\Bigr)
        \\
		\VarInfo_\timestep
		&:= \Bigl(\rf_\dims(\Param_k): k \le \timestep\Bigr)
		\cup \Bigl(
			\langle v, w\rangle :
			v \in \rv_{[0:\dimV_{\timestep+1})}, w \in \gradients_\timestep
		\Bigr).
    \end{align*}
    As the identity is a continuous map, we simply map the function values
    \(\rf_\dims(\Param_k)\) from \(\VarInfo_\timestep\) to itself in
    \(\Info_\timestep\). We therefore only need to find a way to continuously
    construct the inner products of \(\Info_\timestep\) from
    \(\VarInfo_\timestep\). Since \(\nabla\rf_\dims(\Param_k)\) is contained in
    \(V_{\timestep+1}\) by its definition \eqref{eq: vector
    space of evaluation points} and \(\rv_{[0:\dimV_{\timestep+1})}\) is a basis
    of \(V_{\timestep+1}\) by construction (Definition~\ref{def: previsible
    orthonormal coordinate system}), we have for all \(k,l\le \timestep\)
    \begin{align}
        \nonumber
        \langle \nabla \rf_\dims(\Param_k), \nabla\rf_\dims(\Param_l)\rangle
        &= \Bigl\langle
            \sum_{i=1}^{\dimV_{\timestep+1}-1}\langle \nabla \rf_\dims(\Param_k), \rv_i\rangle \rv_i,
            \sum_{j=1}^{\dimV_{\timestep+1}-1}\langle \nabla \rf_\dims(\Param_l), \rv_j\rangle\rv_j
        \Bigr\rangle
        \\
        \label{eq: representation of gradient inner prod}
        &= \sum_{i=0}^{\dimV_{\timestep+1}-1}
        \underbrace{\langle \nabla\rf_\dims(\Param_k), \rv_i\rangle}_{\in \VarInfo_\timestep}
        \underbrace{\langle\nabla\rf_\dims(\Param_l), \rv_i\rangle}_{\in \VarInfo_\timestep}
    \end{align}
    Since \(\dimV_{\timestep +1}\) is almost surely constant by \ref{ind: full
    dimension}\footnote{
        We sketch a path to get rid of the strict positive definiteness assumption
        in the outlook (Section~\ref{sec: outlook}). In this case
        \ref{ind: full dimension} and \ref{ind: evaluation
        points asymptotically different} are lost. This argument then has
        to be replaced using an upper bound on the \(\dimV_{\timestep+1}\) and
        Lemma~\ref{lem: gamma is triangular}.
    }, this covers most of the inner products of \(\Info_\timestep\).
    What is left are the inner products using \(\param_0\). If \(\param_0 = 0\),
    then all those inner products are zero and the zero map does the job. 

    In the following we therefore assume \(\param_0 \neq 0\). Because \(\param_0\)
    is deterministic, we do not need to construct \(\|\param_0\|^2 = \radius^2\)
    from \(\VarInfo_\timestep\) as a constant map does the job. Since
    \(\param_0\neq 0\) implies \(\rv_0 = \frac{\param_0}{\|\param_0\|}\) by
    construction (Definition~\ref{def: previsible orthonormal coordinate
    system}), we have for all \(k\le \timestep\)
    \[
        \langle \param_0, \nabla\rf_\dims(\Param_k)\rangle
        = \|\param_0\| \underbrace{\langle \rv_0, \nabla\rf_\dims(\Param_k)\rangle}_{\in \VarInfo_\timestep}.
    \]
    We have thus continuously constructed all inner products in \(\Info_\timestep\)
    form \(\VarInfo_\timestep\).
\end{proof}

Here is how the main theorem of the paper follows from Theorem \ref{THM: ASYMPTOTICALLY DETERMINISTIC BEHAVIOR VARIANT}.
\begin{proof}[Proof of Theorem~\ref{THM: ASYMPTOTICALLY DETERMINISTIC BEHAVIOR}]
    By Proposition~\ref{prop: wlog deterministic starting points}, we can assume
    without loss of generality that the initial point is deterministic.
    Since the other assumptions of Theorem~\ref{THM: ASYMPTOTICALLY DETERMINISTIC
    BEHAVIOR VARIANT} are the same (or even more general), we only need to prove that
    for every \(\timestep\in\nat\) there exists some \(\liminfo_\timestep\) such
    that
    \[
        \Info_\timestep\underset{\dims\to\infty}{\overset{p}\to} \liminfo_\timestep.
    \]
    As \(\Info_\timestep\) is a continuous function of \(\VarInfo_\timestep\)
    by Lemma~\ref{lem: info is continuous function} this follows immediately
    from continuous mapping by the inductive claim \ref{ind: information convergence}
    of Theorem~\ref{THM: ASYMPTOTICALLY DETERMINISTIC BEHAVIOR VARIANT}.
\end{proof}
The rest of the section is an inductive proof of Theorem~\ref{THM: ASYMPTOTICALLY DETERMINISTIC BEHAVIOR VARIANT}.

\subsection{Proof of Theorem~\ref{THM: ASYMPTOTICALLY DETERMINISTIC BEHAVIOR VARIANT}}

The heart of the proof will
be a lengthy induction. Before, we want to address the additional assumptions
of
\begin{enumerate}
    \item representable limit points \ref{ind: representation} and 
    \item the claim that these points are different \ref{ind: evaluation points asymptotically different},
\end{enumerate}
which we introduced in Theorem~\ref{THM: ASYMPTOTICALLY DETERMINISTIC BEHAVIOR
VARIANT} but did not use in the proof of Theorem~\ref{THM: ASYMPTOTICALLY
DETERMINISTIC BEHAVIOR}. Together with the strictly positive definite covariance of
\((\rf_\dims,\nabla\rf_\dims)\) these assumptions will allow us to argue for converging entries of
covariance matrices and invertible limiting covariance matrices, which are used
in the conditional expectation and conditional variance.

\subsubsection{Complexity reduction}

Before we start the induction, we will perform some complexity reductions.

\begin{enumerate}
    \item We show that \ref{ind: representation} follows from \ref{ind:
    information convergence} in Lemma~\ref{lem: conv inf -> representation}.
    
    \item We show that \ref{ind: evaluation points asymptotically different} follows from
    \ref{ind: full dimension} and \ref{ind: information convergence} in
    Lemma~\ref{lem: conv info, full rank -> different eval pts}.
    
    \item We reduce the work necessary to prove \ref{ind: information convergence}.
\end{enumerate}

In the actual induction we will therefore be able to focus on the claims
\ref{ind: information convergence} and \ref{ind: full dimension}.

\begin{lemma}\label{lem: conv inf -> representation}
    For all fixed \(\timestep\in\nat\), the convergence of information \ref{ind: information convergence} implies the
    representation \ref{ind: representation}.
\end{lemma}

\begin{proof}
    Assuming \ref{ind: information convergence} we have that \(\VarInfo_\timestep \to
    \VarLiminfo_\timestep\). By Lemma~\ref{lem: info is continuous function} this
    also implies \(\Info_\timestep \to \liminfo_\timestep\). Since the prefactors
    \(\lr_m\) are functions of \(\Info_{m-1}\) continuous in \(\liminfo_{m-1}\),
    this implies for all \(m\le \timestep+1\) that the prefactors converge
    \[
        \lr_m = \lr_m(\Info_{m-1})
        \underset{\dims\to\infty}{\overset{p}\to} \lr_m(\liminfo_{m-1})
        =: \hat{\lr}_m.
    \]
    Recall that by \eqref{eq: definition gamma} of \ref{ind: information
    convergence} we have for all \(k\le \timestep\) and all \(i <
    \dimV_{\timestep+1}\) 
    \[
        \langle \nabla\rf_\dims(X_k), \rv_i\rangle
        \underset{\dims\to\infty}{\overset{p}\to}
        \gamma^{(i)}_k
    \]
    for some \(\gamma^{(i)}_k\in \real\). For all \(m\le \timestep+1\) we therefore
    have by definition of \(\Param_m\) (Definition~\ref{def: general gsa})
    \begin{align}
        \nonumber
        \langle  \Param_m, \rv_i\rangle
        &= \lr^{(\param)}_m\langle \param_0, \rv_i\rangle
        + \sum_{k=0}^{m-1} \lr^{(g)}_{m,k} \langle \nabla\rf_\dims(\Param_k), \rv_i\rangle
        \\
        \label{eq: inner product convergence}
        \overset{p}&{\underset{\dims\to\infty}\to}
        \; y^{(i)}_m \quad
        := \hat{\lr}^{(\param)}_m\|\param_0\|\delta_{0i}
        + \sum_{k=0}^{m-1} \hat{\lr}^{(g)}_{m,k} \gamma^{(i)}_k,
    \end{align}
    where \(\delta_{ij}\) denotes the Kronecker delta.
    Since \(\Param_k\) is contained in the vector space of evaluation points \(V_m\)
    for all \(k\le m\) by definition \eqref{eq: vector space of evaluation points},
    its norms converges
    \[
        \|\Param_k\|^2
        = \sum_{i=0}^{\dimV_m-1} \langle \Param_k, \rv_i\rangle^2 
        \overset{p}{\underset{\dims\to\infty}\to}
        \sum_{i=0}^{\dimV_m-1} (y_k^{(i)})^2 = \|y_k\|^2,
    \]
    and likewise their distances for all \(k,l\le m\)
    \[
        \|\Param_k - \Param_l\|^2
        = \sum_{i=0}^{\dimV_m-1} \langle \Param_k - \Param_l, \rv_i\rangle^2
        \underset{\dims\to\infty}{\overset{p}\to} 
        \rho_{kl}^2
        := \|y_k - y_l\|^2.
    \]
    This proves the limiting representation \ref{ind: representation}.
\end{proof}

In the following we will prove, assuming \ref{ind: information convergence}
and \ref{ind: full dimension}, that the limiting distances \(\rho_{kl}\) are
greater zero, i.e. \ref{ind: evaluation points asymptotically different}.
The main ingredients are \ref{ind: full dimension} and the assumed asymptotic
use of the last gradient.

\begin{lemma}\label{lem: conv info, full rank -> different eval pts}
    For all fixed \(\timestep\in\nat\), the convergence of information \ref{ind:
    information convergence} together with the asymptotic full rank \ref{ind:
    full dimension} imply asymptotically different evaluation points \ref{ind:
    evaluation points asymptotically different}.
\end{lemma}

\begin{proof}
    We will proceed by induction over \(m\), where we assume the claim to
    be shown for all \(k,l\le m\le
    \timestep+1\). The induction start \(m=0\) is trivial since a single point
    is always distinct. For the induction step we assume to have the statement
    for \(m\), that is
    \[
        \rho_{kl} > 0
        \quad\text{for all}\quad
        k, l\le m
        \quad\text{with}\quad
        k\neq l.
    \]
    To show the statement for \(m+1\le \timestep+1\), we only need to check the
    distances to the point \(y_{m+1}\) as we have the others by induction.

    To show that \(y_{m+1}\) is distinct from any \(y_k\) with \(k\le m\), we use the fact
    that the most recent gradient is used asymptotically by assumption of the
    Theorem~\ref{THM: ASYMPTOTICALLY DETERMINISTIC BEHAVIOR VARIANT}
    (cf.~Assumption~\ref{assmpt: generalized assumptions}), i.e.
    \begin{equation}
        \label{eq: using the most recent gradient}
        \hat{\lr}^{(g)}_{m+1, m} = \lr^{(g)}_{m+1, m}(\liminfo_m) \neq 0.
    \end{equation}
    The claim \ref{ind: full dimension} ensures that \(V_{m+1}\) has a larger
    dimension than \(V_m\), that is \(\dimV_{m+1}=\dimV_m+1\).  The last basis
    element of \(V_{m+1}\) is thus given by \(\rv_{\dimV_{m+1}-1}=\rv_{\dimV_m}\).
    By \eqref{eq: definition of v candidate} this element is produced from the last
    gradient \(\nabla\rf_\dims(\Param_m)\), which cannot be used for \(\Param_k\) with
    \(k\le m\) but must be used for \(\Param_{m+1}\) asymptotically by \eqref{eq: using
    the most recent gradient}. Therefore \(\Param_{m+1}\) asymptotically contains a
    component of \(\rv_{\dimV_m}\), that no other \(\Param_k\) has, which translates
    to the inequality of the asymptotic representations \(y_{m+1}\) and \(y_k\).

    Formally, since \(\Param_0,\dots\Param_m\) is contained in \(V_m\) spanned by
    \(\rv_{[0:\dimV_m)}\), we have for all \(k\le m\)
    \begin{equation}
        \label{eq: limiting components <= m}
        y_k^{(\dimV_m)}
        \overset{\eqref{eq: inner product convergence}}=
        \lim_{\dims\to\infty} \langle \Param_k, \rv_{\dimV_m}\rangle = 0.
    \end{equation}
    For the asymptotic representation of the last evaluation point \(\Param_{m+1}\)
    on the other hand, we have
    \begin{equation}
        \label{eq: limiting component m+1}
        y_{m+1}^{(\dimV_m)}
        \overset{\eqref{eq: inner product convergence}}= 
        \hat{\lr}^{(x)}_{m+1}\|\param_0\|\underbrace{\delta_{0 \dimV_m}}_{=0} + \sum_{k=0}^{m}\hat{\lr}^{(g)}_{m+1,k}\gamma_k^{(\dimV_m)}
        = \hat{\lr}^{(g)}_{m+1,m} \gamma_m^{(\dimV_m)}.
    \end{equation}
    The last equation is due to \eqref{eq: gamma later are zero}. This follows from
    the fact that for all \(k < m\) the gradient \(\nabla\rf_\dims(\Param_k)\) is
    contained in \(V_m\) spanned by \(\rv_{[0:\dimV_m)}\). By definition of
    \(\gamma_k^{(i)}\) \eqref{eq: definition gamma} we therefore have
    \begin{equation}
        \label{eq: gamma later are zero}
        \gamma_k^{(\dimV_m)}
        \overset{\eqref{eq: definition gamma}}= \lim_{\dims\to\infty} \langle \rf_\dims(\Param_k), \rv_{\dimV_m}\rangle
        = 0.
    \end{equation}
    Putting \eqref{eq: limiting components <= m} and \eqref{eq: limiting component m+1} together
    we have
    \begin{align*}
        \rho_{(m+1)k}^2
        &= \|y_{m+1} - y_k\|^2
        \ge (y_{m+1}^{(\dimV_m)} - y_k^{(\dimV_m)})^2
        = (\hat{\lr}^{(g)}_{m+1,m}\gamma_m^{(\dimV_m)})^2
        > 0,
    \end{align*}
    where the last inequality is due to the asymptotic use of the most recent
    gradient \eqref{eq: using the most recent gradient} and the limiting full rank
    claim \eqref{eq: limiting full rank} of \ref{ind: full dimension}.
\end{proof}

In \ref{ind: information convergence} we claim convergence of \(\VarInfo_\timestep\),
where \(\VarInfo_\timestep\) contains inner products of the form
\[
    \langle \rv_i, \nabla\rf_\dims(\Param_k)\rangle
\]
for \(k\le \timestep\) and \(i< \dimV_{\timestep+1}\). The following lemma
essentially implies that the restriction on \(i\) was unnecessary. So
when we increase \(\timestep-1\) to \(\timestep\) in the induction step, we do
not have to revisit the inner products of old gradients and can focus solely on
\(\nabla\rf_\dims(\Param_\timestep)\).

\begin{lemma}
    \label{lem: gamma is triangular}

    For all \(k\in \nat\) and \(i\ge \dimV_{k+1}\) we have
    \[
        \langle \nabla\rf_\dims(\Param_k), \rv_i\rangle \underset{\dims\to\infty}{\overset{p}\to} 0 = \gamma_k^{(i)}.
    \]
\end{lemma}
\begin{proof}
    For all \(i\ge \dimV_{k+1}\) the basis vector \(\rv_i\) is constructed to be orthogonal to
    \(\nabla\rf_\dims(\Param_k)\) contained in \(V_{k+1}\) spanned by
    \(\rv_{[0:\dimV_{k+1})}\). This implies
    \[
        \lim_{\dims\to\infty}\langle \nabla\rf_\dims(\Param_k), \rv_i\rangle = 0 = \gamma_k^{(i)}
    \]
    where the last equation is simply the definition of \(\gamma_k^{(i)}\) of \eqref{eq: definition gamma}.
\end{proof}

\begin{remark}[Gamma are triangular]\label{rem: gamma is triangular}
    Lemma~\ref{lem: gamma is triangular} can be visualized using a triangular
    matrix as follows. With the definition
    \[
        \real^n := \{ (\param_i)_{i\in [0:\infty)} : \param_i = 0 \quad \forall i\ge n\},
    \]
    we can view \(\real^m\) as a subspace of \(\real^n\) for \(m\le n\). The vector
    \[
        \gamma_k := (\gamma_k^{(i)})_{i\in [0:\dimV_{k+1})} \in \real^{\dimV_{k+1}},
    \]
    is then (by Lemma~\ref{lem: gamma is triangular}) also a member of
    \(\real^m\) for \(m\ge \dimV_{k+1}\). Concatenating the vectors
    \[
        \gamma_{[0:\timestep]}
        = \begin{pmatrix}
            \gamma_0^{(0)} & \dots & \gamma_\timestep^{(0)}\\
            \vdots & & \vdots\\
            \gamma_0^{(\dimV_{\timestep+1}-1)}
            & \dots
            & \gamma_\timestep^{(\dimV_{\timestep+1}-1)}
        \end{pmatrix}
        = \begin{pmatrix}
            \gamma_0^{(0)} & \dots & \gamma_{\timestep-1}^{(0)} & \gamma_\timestep^{(0)}\\
            \gamma_0^{(1)} & \dots & \gamma_{\timestep-1}^{(1)} & \gamma_\timestep^{(1)}
            \\
            0 & 
            \\
            \vdots  & \ddots & & \vdots
            \\
            0 & \dots & 0 & \gamma_\timestep^{(\dimV_{\timestep+1}-1)}
        \end{pmatrix}
    \]
    therefore results in a upper triangular matrix above an offset
    diagonal, since we always have \(\dimV_k \le k+1\) due to the definition of
    \(V_k\) \eqref{eq: vector space of evaluation points}.

    Note that, for visualization purposes, we assumed \(\dimV_{\timestep+1}
    =\dimV_{\timestep} +1\) in the second representation. If the dimension stays
    constant at some times \(k\), then the zeros encroach above this diagonal.
\end{remark}

\subsubsection{Induction start with \texorpdfstring{\(\timestep=0\)}{\timestep=0}}

We start the induction by proving the claim for \(\timestep=0\). We have
structured the induction start similar to the induction step. We therefore
suggest the reader to familiarize themselves with this strategy here, as it is
easier to get lost in the details of the induction step.

By Lemma~\ref{lem: conv inf -> representation} and Lemma~\ref{lem: conv
info, full rank -> different eval pts} we only need to prove \ref{ind:
information convergence} and \ref{ind: full dimension}. For \ref{ind:
information convergence} we need to prove that
\[
    \VarInfo_0
    = \bigl(\rf_\dims(\param_0)\bigr)
    \cup
    \bigl(
        \langle \nabla\rf_\dims(\param_0), v\rangle:
        v \in \rv_{[0:\dimV_1)}
    \bigr)
    = \bigl(
        \rf_\dims(\param_0),
        \langle \nabla\rf_\dims(\param_0), \rv_{[0:\dimV_1)}\rangle
    \bigr)
\]
converges to a limiting \(\VarLiminfo_0 = (\limf_0, \gamma_0^{(0)}, \dots, \gamma_0^{(\dimV_1-1)})\)
in probability. Note that \(\dimV_1 \le 2\) as the vector space of
evaluation points is defined in \eqref{eq: vector space of evaluation points} to be previsible,
that is
\[
    V_1 = \Span\{\param_0, \nabla\rf_\dims(\param_0)\}.
\]
So we have \(\dimV_1 -1\le 1\) depending on whether the starting point
\(\param_0\) is zero. The vector
\[
    \gamma_0 = (\gamma_0^{(i)})_{i< \dimV_1} = (\gamma_0^{(0)}, \dots, \gamma_0^{(\dimV_1-1)})
\]
might therefore collapse to a single entry \(\gamma_0^{(0)}\). The approach we
will now take mirrors the approach in the induction step. Beyond the pedagogical
benefit, this order is also quite natural for the induction start when using the
notation we introduced.
\begin{enumerate}[label=\emph{Step \arabic*},wide=0pt, leftmargin=\parindent]
    \item\label{it: step 1 induction start}
    First we prove that
    \[
        \bigl(\rf_\dims(\param_0)\bigr) \cup
        \bigl( \langle \nabla\rf_\dims(\param_0), v\rangle : v\in \rv_{[0:\red{\dimV_0})}\bigr)
    \]
    converges. We highlight that the limit \red{\(\dimV_0\)} is purposefully different from the
    limit in \(\VarInfo_0\).

    \item\label{it: step 2 induction start}
    We will then squeeze in the proof of \ref{ind: full dimension}, which ensures the dimension
    actually increases
    \[
        \dimV_1 = \dimV_0 + 1.
    \]
    
    \item\label{it: step 3 induction start}
    Finally we prove that \(\langle \nabla\rf_\dims(\param_0), \rv_{\dimV_0}\rangle = \langle \nabla\rf_\dims(\param_0), \rv_{\dimV_1-1}\rangle\) converges.
\end{enumerate}
While we could prove \ref{ind: full dimension} before \ref{ind: information convergence} in
the induction start, we will require the results of \ref{it: step 1 induction start} for
\ref{ind: full dimension} in the induction step. The element \(\rv_{\dimV_0}\) is not
only different from \(\rv_{[0:\dimV_0)}\) in the sense that the dimension increase needs
to be shown, it is also the first truly random basis element here. In the induction
step the difference will be between the previsible basis elements and the last
non-previsible element, which needs to be treated differently.

\paragraph*{\ref{it: step 1 induction start}}
Since our random function is (non-stationary) isotropic with \(\rf_\dims\sim \normal(\mu, \kernel)\), we have
\[
    \rf_\dims(\param_0)
    \sim \normal\Bigl(
        \mu\bigl(\tfrac{\|\param_0\|^2}2\bigr),
        \frac1\dims \kernel\bigl(\tfrac{\|\param_0\|^2}2,\tfrac{\|\param_0\|^2}2, \|\param_0\|^2 \bigr)
    \Bigr).
\]
this immediately implies convergence of the first component
\[
    \rf_\dims(\param_0)
    \overset{p}{\underset{\dims\to\infty}\to}
    \mu\bigl(\tfrac{\|\param_0\|^2}2\bigr)
    = \mu\bigl(\tfrac{\radius^2}2\bigr) =: \limf_0.
\]
Let us now turn to the convergence of the inner products. We consider two
cases.

\textbf{Case (\(\param_0\neq 0\)):} In this case we have \(\rv_0 =
\frac{\param_0}{\|\param_0\|}\) by Definition~\ref{def: previsible orthonormal
coordinate system}. Therefore \(\rv_0\) is deterministic and by an application
of Lemma~\ref{lem: cov of derivatives, non-stationary isotropic case} we have
\[
    D_{\rv_0}\rf_\dims(\param_0) \sim \normal\Bigl(
        \mu'\bigl(\tfrac{\|\param_0\|^2}2\bigr)\|\param_0\|,
        \frac1\dims\bigl[(\kernel_{12} + \kernel_{13} + \kernel_{32} + \kernel_{33})\|\param_0\|^2
        + \kernel_3\bigr]
    \Bigr),
\]
using the notation
\[
    \kernel
    := \kernel\bigl(\tfrac{\|\param_0\|^2}2, \tfrac{\|\param_0\|^2}2, \|\param_0\|^2\bigr)
    = \kernel\bigl(\tfrac{\radius^2}2, \tfrac{\radius^2}2, \radius^2\bigr)
\]
to omit inputs to the kernel \(\kernel\).
But this immediately implies
\[
    \langle \nabla\rf_\dims(\param_0), \rv_0\rangle 
    = D_{\rv_0}\rf_\dims(\param_0)
    \overset{p}{\underset{\dims\to\infty}\to}
    \mu'\bigl(\tfrac{\|\param_0\|^2}2\bigr)\|\param_0\|
    =: \gamma_0^{(0)}.
\]
As \(\dimV_0 = \dim(V_0) = \dim(\Span(\param_0)) = 1\), we have covered all elements
in \(\rv_{[0:\dimV_0)}\).

\textbf{Case (\(\param_0 = 0\)):} In this case, \(\dimV_0 = \dim(V_0) = 0\)
and \(\rv_{[0:\dimV_0)}\) is empty. We therefore do not have to do anything.

\paragraph*{\ref{it: step 2 induction start}}

To prove the dimension increases, we consider the Gram-Schmidt candidate
\begin{align*}
    \tilde{\rv}_0
    &:= \nabla\rf_\dims(\param_0) - P_{V_0}\rf_\dims(\param_0)
    = \nabla\rf_\dims(\param_0) - \sum_{i<\dimV_0} \langle\nabla\rf_\dims(\Param_0), \rv_i\rangle \rv_i
    \\
    &= \sum_{i= \dimV_0}^{\dims-1} \langle\nabla\rf_\dims(\Param_0), \rw^{(0)}_i\rangle \rw^{(0)}_i,
\end{align*}
where \(\rw^{(0)}_{[\dimV_0:\dims)}\) is the basis defined to be orthogonal to
\(V_0\) (Definition~\ref{def: previsible orthonormal coordinate system}). Since
\(V_0\) is deterministic, this orthogonal basis is also deterministic. In the
induction step, both will only be previsible. Note that in the case \(\param_0 = 0\),
\(\tilde{\rv}\) is simply the entire gradient by definition of \(\dimV_0\).

Since the \(\rw^{(0)}_i\) are deterministic, and orthogonal to \(\param_0\) in
either case, we can apply Lemma~\ref{lem: cov of derivatives, non-stationary isotropic case} to obtain
\[
    D_{\rw^{(0)}_i}\rf_\dims(\param_0) \overset{\iid}\sim \normal\Bigl(0, \frac{\kernel_3}\dims\Bigr).
\]
This implies that there exist \(Y_i\overset{\iid}{\sim}\normal(0,1)\) such that
\[
    D_{\rw^{(0)}_i}\rf_\dims(\param_0) = \sqrt{\frac{\kernel_3}\dims} Y_i.
\]
But since \(\dimV_0 \le 1\) and in particular \(\dimV_0\) is finite, this implies
with the law of large numbers
\begin{equation}
    \label{eq: LLN induction start}
    \|\tilde{\rv}_0\|^2
    = \sum_{i=\dimV_0}^{\dims-1} \langle \nabla\rf_\dims(\param_0), w_i\rangle^2
    = \kernel_3 \cdot\frac1\dims \sum_{i=\dimV_0}^{\dims-1} Y_i^2
    \overset{p}{\underset{\dims\to\infty}\to} \kernel_3.
\end{equation}
Since we have \(\kernel_3>0\) by Lemma~\ref{lem: kappa_3 positive},
\(\|\tilde{\rv}_0\|^2\) is almost surely strictly greater than zero, which
implies \(\dim(V_1) > \dim(V_0)\) almost surely. Additionally,
the dimension also increases in the limit, i.e.
\(\lim_{\dims\to\infty}\|\tilde{\rv}_0\|=\kernel_3>0\).
We have therefore shown all of \ref{ind: full dimension}.

\paragraph*{\ref{it: step 3 induction start}}

We have by definition of \(\rv_{\dimV_0}\) in \eqref{eq: definition of v} and
the definition of \(\tilde{\rv}_0\)
\[
    \langle \nabla\rf_\dims(\param_0), \rv_{\dimV_0}\rangle 
    = \bigl\langle \nabla\rf_\dims(\param_0), \tfrac{\tilde{\rv}_0}{\|\tilde{\rv}_0\|}\bigr\rangle
    = \|\tilde{\rv}_{0}\|
    \overset{p}{\underset{\dims\to\infty}\to} \sqrt{\kernel_3}
    =: \gamma_0^{(\dimV_0)}.
\]
This finishes the last step and therefore the induction start.

\subsubsection{Induction step \texorpdfstring{\((\timestep-1)\to \timestep\)}{}}

We now get to the main body of the proof. Before we start, let us recapitulate
the lemmas we proved for complexity reduction. By Lemma~\ref{lem: conv inf ->
representation} and Lemma~\ref{lem: conv info, full rank -> different eval pts}
it is sufficient to prove the statements \ref{ind: information convergence} and
\ref{ind: full dimension}, as the other follow.

With \(\gradients_\timestep
:= \bigl(\nabla \rf_\dims(\Param_k): k \le \timestep\bigr)\) the modified information \(\VarInfo_\timestep\)
of the information convergence \ref{ind: information convergence} was given by
\[
    \VarInfo_\timestep
    := \Bigl(\rf_\dims(\Param_k): k \le \timestep\Bigr)
    \cup \Bigl(
        \langle v, w\rangle :
        v \in \rv_{[0:\dimV_{\timestep+1})}, w \in \gradients_\timestep
    \Bigr).
\]
By induction we already have \(\VarInfo_{\timestep-1}\to \VarLiminfo_{\timestep-1}\).
And due to our discussion of \(\langle \rv_i, \nabla\rf_\dims(\Param_k)\rangle\) for
\(i\ge \dimV_{\timestep+1} \ge \dimV_{k+1}\) in Lemma~\ref{lem: gamma is triangular}
we therefore only need to prove
\begin{equation}
    \label{eq: new column is convergence claim} 
    \begin{pmatrix}
        \rf_\dims(\Param_\timestep)
        \\
        \bigl\langle \rv_{[0:\dimV_{\timestep+1})}, \nabla\rf_\dims(\Param_\timestep)\bigr\rangle
    \end{pmatrix}
    = \begin{pmatrix}
        \rf_\dims(\Param_\timestep)
        \\
        D_{\rv_0}\rf_\dims(\Param_\timestep)
        \\
        \vdots
        \\
        D_{\rv_{\dimV_{\timestep+1}-1}}\rf_\dims(\Param_\timestep)
    \end{pmatrix}
    \underset{\dims\to\infty}{\overset{p}\to} \begin{pmatrix}
        \limf_\timestep
        \\
        \gamma_\timestep
    \end{pmatrix},
\end{equation}
where we used the definition \(\gamma_\timestep = (\gamma_\timestep^{(i)})_{i\in
[0:\dimV_{\timestep+1})}\) from Remark~\ref{rem: gamma is triangular}. In this
remark we have also explained how \(\gamma_k \in \real^{\dimV_{k+1}}\) can
be viewed as a member of \(\real^{\dimV_{\timestep+1}}\) and laid out the
concatenated vectors \(\gamma_{[0:\timestep]}\) in triangular matrix form.
In the following we will arrange the random information into a matching
matrix form by the definition of an auxiliary function \(Z\).

This function has multiple purposes: First, it turns the direction vectors of the
directional derivatives into inputs for an application of previsible sampling
(Corollary \ref{cor: gaussian previsible sampling}), which we explained in Idea 3 of the
proof sketch (Section~\ref{sec:sketch}). Second, it structures the information
vector \(\VarInfo_\timestep\) into a better readable matrix form for human
digestion. Third, it allows us to rearrange the elements of the column vectors
we concatenate into row vectors. This effectively rearranges the covariance
matrix of \(Z\) into a more sparse block form.

With some abuse of notation we define \(Z\) for a varying number\footnote{
    formally, one can view this as slices of a function \(\real^{\dims^2} \times
    \real^{(\timestep+1)\dims} \to \real^{(\dims+1)(\timestep+1)}\).
} 
of direction vectors \(w_1,\dots,w_m\in \real^\dims\) and evaluation points
\(\param_1,\dots, \param_k\in\real^\dims\) for \(m\le \dims\) and \(k\le \timestep\)
\[
    Z(w_1,\dots, w_m; \param_0,\dots, \param_k) := \begin{pmatrix}
        \rf_\dims(\param_0) & \dots & \rf_\dims(\param_k)\\
        D_{w_1}\rf_\dims(\param_0) & \cdots & D_{w_1}\rf_\dims(\param_k)
        \\
        \vdots & & \vdots
        \\
        D_{w_m}\rf_\dims(\param_0) & \cdots & D_{w_m}\rf_\dims(\param_k)\\
    \end{pmatrix}.
\]
The number of inputs is therefore not variable and we separate their
type by a semicolon. Note that we are not interested in this matrix as an
operator. The two dimensional layout is only used for better readability and we
will generally view it as a vector. For this purpose we treat the matrix as
`row-major', i.e. whenever we treat it like a vector, we concatenate the rows.
This grouping is purposefully different from the column layout of the \(\gamma_k\)
as it will enable a block matrix sparsity in the covariance matrices later on.

By induction we have information convergence \ref{ind: information convergence}
for \(\timestep-1\), which can be expressed using \(Z\) as
\begin{equation}
    \label{eq: main induction assumption}
    Z(\rv_{[0:\dimV_\timestep)}; \Param_{[0:\timestep)})
    = \begin{pmatrix}
        \rf_\dims(\Param_0) & \cdots & \rf_\dims(\Param_{\timestep-1}) 
        \\
        D_{\rv_0}\rf_\dims(\Param_0) & \cdots & D_{\rv_0}\rf_\dims(\Param_{\timestep-1})
        \\
        \vdots & & \vdots
        \\
        D_{\rv_{\dimV_\timestep-1}}\rf_\dims(\Param_0) & \cdots & D_{\rv_{\dimV_\timestep-1}}\rf_\dims(\Param_{\timestep-1})
    \end{pmatrix}
    \overset{p}{\underset{\dims\to\infty}\to}
    \begin{pmatrix}
        \limf_{[0:\timestep)}\\
        \gamma_{[0:\timestep)}
    \end{pmatrix}.
\end{equation}
Observe that the induction step simply adds the additional column given in
\eqref{eq: new column is convergence claim} to the right of the matrix,
where we use Lemma~\ref{lem: gamma is triangular} to argue that we can arbitrarily
extend the number of rows.

We will now follow the same strategy as in the induction start:

\begin{enumerate}[label=\emph{Step \arabic*}]
    \item\label{it: step 1 induction step}
    First, we prove the convergence of the `new column'
    \[
        Z(\rv_{[0:\red{\dimV_{\timestep}})}; \Param_\timestep)
        = \begin{pmatrix}
            \rf_\dims(\Param_\timestep)\\
            \langle \rv_{[0:\red{\dimV_{\timestep}})}, \nabla\rf_\dims(\Param_\timestep)\rangle
        \end{pmatrix}
        \underset{\dims\to\infty}{\overset{p}\to}
        \begin{pmatrix}
            \limf_\timestep\\
            \gamma_\timestep^{([0:\dimV_\timestep))}
        \end{pmatrix}
    \]
    \item\label{it: step 2 induction step}
    We squeeze in the proof of \ref{ind: full dimension}, which ensures asymptotically
    \[
        \dimV_{\timestep+1} = \dimV_\timestep + 1.
    \]
    \item\label{it: step 3 induction step}
    We finally prove convergence of the `new corner element'
    \[
        Z(\rv_{\dimV_\timestep}; \Param_\timestep)
        = \langle \rv_{\dimV_\timestep}, \nabla\rf_\dims(\Param_\timestep)\rangle
        \underset{\dims\to\infty}{\overset{p}\to}
        \gamma_\timestep^{(\dimV_\timestep)}.
    \]
\end{enumerate}
The reason for this split are twofold. First, we have that
\(\rv_{[0:\dimV_\timestep)}\) is previsible, that is
\(\filt_{\timestep-1}\)-measurable, while \(\rv_\timestep\) is not.
\(\rv_\timestep\) must therefore be treated differently. Second,
we need to construct \(\rv_{\dimV_\timestep}\) from \(\tilde{\rv}_\timestep\),
which will naturally prove \ref{ind: full dimension} before we get
to the convergence of the inner product in \ref{it: step 3 induction step}.

This strategy can get a bit lost, as we will spend the majority of our time
with \ref{it: step 1 induction step} before wrapping up \ref{it: step 2
induction step} and \ref{it: step 3 induction step} fairly quickly.
An additional reason for this strategy to get lost is, that we will not
just consider
\((
    \rf_\dims(\param_\timestep),
    \langle \rv_{[0:\dimV_\timestep)}, \nabla \rf_\dims(\Param_\timestep)\rangle
)\), but instead consider the conditional distribution of
\((
    \rf_\dims(\param_\timestep),
    \langle \Basis_\timestep, \nabla \rf_\dims(\Param_\timestep)\rangle
)\)
for the entire previsible basis \(\Basis_\timestep = (\rv_{[0:\dimV_\timestep)},
\rw^{(\timestep)}_{[\dimV_\timestep:\dims)})\). We will aggregate the
directional derivatives in the directions \(\rw^{\timestep}_i\) later into
\(\tilde{\rv}_\timestep\), which means that we work towards \ref{it: step 1
induction step} and \ref{it: step 2 induction step} simultaneously.
The first objective will be, to apply the previsible sampling lemma
such that we can treat the previsible inputs as deterministic.

\paragraph*{Getting rid of the random input}

In the proof sketch we outlined how we needed to treat the random input with
care (cf.~\ref{sec:sketch},  Idea~3). In particular we need to consider the
(basis; point) pairs
\[
    (\Basis_0; \Param_0),\cdots, (\Basis_\timestep;\Param_\timestep)
    := \Bigl(\rv_{[0:\dimV_0)}, \rw_{[\dimV_0:\dims)}^{(0)}; \Param_0\Bigr),
    \cdots,
    \Bigl(\rv_{[0:\dimV_\timestep)}, \rw_{[\dimV_\timestep:\dims)}^{(\timestep)}; \Param_\timestep\Bigr).
\]
Recall that we defined the basis \(\Basis_k\) in Definition~\ref{def: previsible
orthonormal coordinate system} to be previsible just like \(\Param_k\), i.e.\
measurable with regard \(\filt_{k-1}\). Moreover, since all basis
changes are previsible and invertible it is straightforward
to check inductively that there exists a measurable bijective map between
different basis representations of the derivative information up to \(\timestep\)
\[
    \bigl(\rf_\dims(\Param_k),\nabla\rf_\dims(\Param_k)\bigr)_{k<\timestep}
    \leftrightsquigarrow \bigl(Z(\Basis_k; \Param_k)\bigr)_{k<\timestep}.
\]
This implies that their generated sigma algebras are the same
\[
    \filt_{\timestep-1}
    \overset{\text{def.}}=
    \sigma\Bigl(\bigl(\rf_\dims(\Param_k),\nabla\rf_\dims(\Param_k)\bigr), k\le\timestep-1\Bigr)
    = \sigma\bigl(Z(\Basis_k; \Param_k), k\le\timestep-1\bigr).
\]
By Application of Corollary \ref{cor: gaussian previsible sampling} we then have, for all
bounded measurable \(h\),
\[
    \E\Bigl[
        h\bigl(Z(\Basis_\timestep;\Param_\timestep)\bigr) \mid \filt_{\timestep-1}
    \Bigr]
    = F(\Param_{[0:\timestep]}; \Basis_{[0:\timestep]}),
\]
for the function \(F\) defined with deterministic evaluation points \(\param_k\)
and basis \(\basis_k\) as
\[
    F(\param_{[0:\timestep]}; \basis_{[0,\timestep]})
    := \int h(z) \Pr\Bigl[
        Z(\basis_\timestep;\param_\timestep) \in dz \mid 
        Z(\basis_0;\param_0), \dots, Z(\basis_{\timestep-1};\param_{\timestep-1})
    \Bigr],
\]
using the standard consistent joint conditional distribution for Gaussian distributions
(cf.\ Example \ref{ex: gaussian conditional}).
In other words, we are allowed to treat the inputs \((\Basis_k; \Param_k)\)
as deterministic when calculating the conditional distribution. But keeping
track of \(\timestep+1\) different basis \(\Basis_k\) for every is evaluation
point \(\Param_k\) is very inconvenient. So our next goal is to reduce the
number of basis to one.
For this we optimistically define for a single basis \(\basis\)
\[
    G(\param_{[0:\timestep]}; \basis)
    := \E\Bigl[
        h\bigl(Z(\basis;\param_\timestep)\bigr) \mid 
        Z(\basis;\param_0), \dots, Z(\basis;\param_{\timestep-1})
    \Bigr],
\]
where we assume here and in the following, that this conditional distribution is
defined via the standard consistent joint conditional distribution of Example
\ref{ex: gaussian conditional}.
We are now going to observe that \(F\) is constant in all but the most recent
basis \(\basis_\timestep\), i.e.
\begin{equation}
    \label{eq: basis reduction}
    F(\param_{[0:\timestep]}; \basis_{[0:\timestep]})
    = G(\param_{[0:\timestep]}; \basis_\timestep).
\end{equation}
That is because the sigma algebras generated by different basis representations
are identical as they can be bijectively translated into each other, i.e.
for any \(\basis\)
\[
    \sigma(Z(\basis_0;\param_0), \dots, Z(\basis_{\timestep-1};\param_{\timestep-1}))
    = \sigma(Z(\basis;\param_0), \dots, Z(\basis;\param_{\timestep-1})).
\]
In particular this is true for \(\basis = \basis_\timestep\) and thus we have
\eqref{eq: basis reduction}, i.e.
\begin{align*}
    &\E\Bigl[
        h\bigl(Z(\basis_\timestep;\param_\timestep)\bigr) \mid 
        Z(\basis_0;\param_0), \dots, Z(\basis_{\timestep-1};\param_{\timestep-1})
    \Bigr]
    \\
    &= \E\Bigl[
        h\bigl(Z(\basis_\timestep;\param_\timestep)\bigr) \mid 
        Z(\basis_\timestep;\param_0), \dots, Z(\basis_\timestep;\param_{\timestep-1})
    \Bigr].
\end{align*}
Since we only need to consider the most recent basis from now on, we
drop the index and write \(\rw_k := \rw_k^{(\timestep)}\) and summarize our
result as
\begin{equation}
    \label{eq: inputs as deterministic}
    \E\Bigl[
        h\bigl(Z(\rv_{[0:\dimV_\timestep)}, \rw_{[\dimV_\timestep:\dims)}; \Param_\timestep)\bigr) \mid \filt_{\timestep-1}
    \Bigr]
    = G(\Param_0,\dots, \Param_\timestep; \rv_{[0:\dimV_\timestep)}, \rw_{[\dimV_\timestep:\dims)}).
\end{equation}
Recall that we use the semicolon to separate the basis elements from the evaluation points,
which should be treated as concatenated vectors respectively.

In essence, by definition of \(G\) we can now treat our evaluation points
\(\Param_{[0:\timestep]}\) and coordinate system \(\Basis_\timestep =
(\rv_{[0:\dimV_\timestep)}, \rw_{[\dimV_\timestep:\dims)})\) as deterministic when
calculating the conditional distribution.

\paragraph*{The conditional distribution is known in the Gaussian case}

Since \(Z\) is a Gaussian random function, we have for non-random input
\((\basis, \param)\) that its conditional distribution
is also Gaussian (cf.~Theorem~\ref{thm: conditional gaussian distribution}).
Equation~\eqref{eq: inputs as deterministic} translates this result to the
random input and we can therefore conclude that
\[
    Z(\Basis_\timestep; \Param_\timestep) \mid \filt_{\timestep-1}
    \sim \normal\Bigl(
        \E[Z(\Basis_\timestep; \Param_\timestep) \mid \filt_{\timestep-1}],
        \Cov[Z(\Basis_\timestep; \Param_\timestep) \mid \filt_{\timestep-1}]
    \Bigr).
\]
In particular there exist independent \(Y_1,\dots, Y_\dims \sim \normal(0,1)\)
independent of \(\filt_{\timestep-1}\) such that we have in distribution
\begin{equation}
    \label{eq: decomposition in distribution}
    Z(\Basis_\timestep; \Param_\timestep)
    = \E[Z(\Basis_\timestep; \Param_\timestep) \mid \filt_{\timestep-1}]
    + \sqrt{\Cov[Z(\Basis_\timestep; \Param_\timestep) \mid \filt_{\timestep-1}]} \begin{pmatrix}
        Y_0
        \\
        \vdots
        \\
        Y_{\dims}
    \end{pmatrix},
\end{equation}
where we denote the cholesky decomposition of a matrix \(A\) by the squareroot \(\sqrt{A}\).
\begin{remark}
    As we only wish to prove convergence in probability to deterministic numbers,
    this distributional equality will be enough for our purposes, but if the
    covariance was invertible one could get a almost sure  equality (cf.~Remark~\ref{rem:
    decomposition}). Later, we will also see that the covariance is at least
    asymptotically invertible due to strict positive definiteness of \(Z\)
    and asymptotically different evaluation points \ref{ind: evaluation points asymptotically different}.
\end{remark}

\paragraph*{Calculating the conditional expectation and covariance}

Since the first two moments can be calculated by treating the inputs as deterministic
(cf.~Equation~\eqref{eq: inputs as deterministic}), we can calculate the
conditional expectation and variance using the well known formulas for
Gaussian conditionals (cf.~Theorem~\ref{thm:
conditional gaussian distribution})
\begin{align}
    \label{eq: cond expect}
    \E[Z(\Basis_\timestep; \Param_\timestep) \mid \filt_{\timestep-1}]
    &= \rmean_{\timestep}^\dims + \blue{\rcov^\dims_{[0:\timestep),\timestep}}^T [\magenta{\rcov^\dims_{[0:\timestep)}}]^{-1} (Z(\Basis_\timestep; \Param_{[0:\timestep)})- \rmean_{[0:\timestep)}^\dims)
    \\
    \label{eq: cond variance}
    \Cov[Z(\Basis_\timestep; \Param_\timestep) \mid \filt_{\timestep-1}]
    &= \tfrac1\dims\Bigl[\green{\rcov^\dims_\timestep}
    - \blue{\rcov^\dims_{[0:\timestep),\timestep}}^T
    [\magenta{\rcov^\dims_{[0:\timestep)}}]^{-1}
    \blue{\rcov^\dims_{[0:\timestep),\timestep}}\Bigr],
\end{align}
where we define the mixed covariance by
\begin{align*}
    \tfrac1\dims\blue{\rcov^\dims_{[0:\timestep),\timestep}}
    &:= \C_Z(\Basis_\timestep; \Param_{[0:\timestep)}, \Param_\timestep)
    &
    \C_Z(\basis; \param_{[0:\timestep)}, \param_\timestep)
    &:= \Cov(Z(\basis; \param_{[0:\timestep)}), Z(\basis; \param_\timestep)),
\intertext{
    the autocovariance matrices by
}
    \tfrac1\dims\magenta{\rcov^\dims_{[0:\timestep)}}
    &:= \C_Z(\Basis_\timestep; \Param_{[0:\timestep)})
    &
    \C_Z(\basis; \param_{[0:\timestep)})
    &:= \Cov[Z(b; \param_{[0:\timestep)})]
    \\
    \tfrac1\dims\green{\rcov^\dims_\timestep}
    &:= \C_Z(\Basis_\timestep; \Param_\timestep)
    &
    \C_Z(\basis; \param_\timestep)
    &:=\Cov[Z(\basis; \param_\timestep)]
\intertext{
    and the expectations by
}
    \rmean_\timestep^\dims
    &:= \mu_Z(\Basis_\timestep; \Param_\timestep)
    &
    \mu_Z(\basis; \param_\timestep)
    &:= \E[Z(\basis; \param_\timestep)]
    \\
    \rmean_{[0:\timestep)}^\dims
    &:= \mu_Z(\Basis_\timestep; \Param_\timestep)
    &
    \mu_Z(\basis; \param_{[0:\timestep)})
    &:= \E[Z(\basis; \param_{[0:\timestep)})]
\end{align*}
The detour over the functions \(\C_Z\) and \(\mu_Z\) was necessary
to define the unconditional covariance matrices properly, because
we may only treat previsible random variables as deterministic
in the calculation of conditional distributions (using Corollary \ref{cor: gaussian previsible sampling}). So in general, since \(\Basis_\timestep\) and
\(\Param_\timestep\) are random variables,
we have for the unconditional covariance
\[
    \C_Z(\Basis_\timestep; \Param_\timestep) \neq \Cov[Z(\Basis_\timestep; \Param_\timestep)].
\]
To ensure the inputs are treated as deterministic as Equation~\eqref{eq: inputs
as deterministic} demands, we therefore had to make sure the covariance was
already calculated before we plugged in our random input.

Note that we have moved the dimensional scaling \(\frac1\dims\) of the covariances
(cf. Lemma \ref{lem: cov of derivatives, non-stationary isotropic case}) outside of our definition of
\(\magenta{\rcov^\dims_{[0:\timestep)}}\), \(\blue{\rcov^\dims_{[0:\timestep),\timestep}}\)
and \(\green{\rcov^\dims_\timestep}\), as we are now going to prove their entries, and
the entries of \(\rmean_{[0:\timestep)}^\dims\) and \(\rmean_\timestep^\dims\),
converge. This will eventually allow us to prove that the conditional expectation and
covariance will converge, which leads to the stochastic convergence of \(Z\) we
want to obtain for \ref{ind: information convergence}. Moving the dimensional scaling out also made the
dimensional scaling of the conditional covariance in \eqref{eq: cond variance}
much more visible.

\paragraph*{Covariance matrix entries converge}

Recall that \(Z\) is made up of evaluations of \(\rf_\dims\) and
\(D_v\rf_\dims\) for directions \(v\). Thus the entries of its covariance matrices are
given by Lemma~\ref{lem: cov of derivatives, non-stationary isotropic case}, which we restate here for
your convenience.

\covOfDerivatives*

Recall that we have \eqref{eq: <X, v> converges to gamma} of \ref{ind: representation}
for \(\timestep-1\) by induction assumption, i.e. we have for all \(k\le \timestep\)
and all \(i < \dimV_\timestep\)
\[
    \langle \Param_k, \rv_i\rangle
    \underset{\dims\to\infty}{\overset{p}\to} y_k^{(i)}
    \quad\text{and}\quad
    \|\Param_k\|^2
    = \sum_{i=0}^{\dimV_\timestep-1} \langle \Param_k, \rv_i\rangle^2
    \underset{\dims\to\infty}{\overset{p}\to}
    \|y_k\|^2.
\]
Since the \(\Param_k\) for \(k\le \timestep\) are contained in \(V_\timestep\)
orthogonal to \(\rw_{[\dimV_\timestep: \dims)}\) we also have for all \(i\ge
\dimV_\timestep\)
\begin{equation}
    \label{eq: param and w is orthogonal}
    \langle \Param_k, \rw_i\rangle = 0.
\end{equation}
Put together, we have that all the inner products \(\langle \Param_k, v\rangle\)
for \(k\le \timestep\) and \(v\in \Basis_\timestep\) converge.

The entries of 
\(\magenta{\rcov^\dims_{[0:\timestep)}}\),
\(\blue{\rcov^\dims_{[0:\timestep),\timestep}}\),
\(\green{\rcov^\dims_\timestep}\), \(\rmean_{[0:\timestep)}^\dims\) and
\(\rmean_\timestep^\dims\) calculated with
Lemma~\ref{lem: cov of derivatives, non-stationary isotropic case} therefore all converge in probability
by the continuous mapping theorem, since \(\kernel\) and \(\mu\) are
sufficiently smooth by Assumption~\ref{assmpt: smoothness} used in
Theorem~\ref{THM: ASYMPTOTICALLY DETERMINISTIC BEHAVIOR VARIANT}.
Observe that it was very important to remove the dimensional scaling \(\frac1\dims\)
of \eqref{eq: cov df, f} and \eqref{eq: cov df, df}
from the covariance matrices \(\magenta{\rcov^\dims_{[0:\timestep)}}\),
\(\blue{\rcov^\dims_{[0:\timestep),\timestep}}\) and
\(\green{\rcov^\dims_\timestep}\), as their entries would otherwise all
converge to zero. As we want to invert \(\magenta{\rcov^\dims_{[0:\timestep)}}\)
this would have been very inconvenient.

We further note that the sizes of these covariance matrices change with the dimension,
because the number of directional derivatives increases with \(\dims\) which
increases the size of \(Z(\rv_{[0:\dimV_\timestep)}, \rw_{[\dimV_\timestep:\dims)}; \Param_{[0:\timestep)})\)
and therefore the size of its covariance matrix. The convergence of their
entries is therefore not yet sufficient for a limiting object to be well defined.

\paragraph*{Splitting the increasing matrices into block matrices of constant size}

This is the heart of the proof, which relies heavily on our custom
coordinate system \(\Basis_\timestep=(\rv_{[0:\dimV_\timestep)}, \rw_{[\dimV_\timestep:\dims)})\). To understand this,
let us focus on the covariance of derivatives given in \eqref{eq: cov df, df} of
Lemma~\ref{lem: cov of derivatives, non-stationary isotropic case}, i.e.
\begin{equation}
    \tag{\ref{eq: cov df, df}}
    \Cov(D_v \rf_\dims(x), D_w\rf_\dims(y))   
    = \frac1\dims\Bigl[\underbrace{\begin{aligned}[t]
        &\kernel_{12} \langle x, v\rangle\langle y, w\rangle + \kernel_{13} \langle x, v\rangle\langle x, w\rangle
        \\
        & + \kernel_{32} \langle y, v\rangle\langle y, w\rangle + \kernel_{33}\langle y, v\rangle\langle y, w\rangle
        \\
    \end{aligned}}_{\text{(I)}}
    + \underbrace{\kernel_3 \langle v, w\rangle}_{\text{(II)}}\Bigr].
\end{equation}
Notice that for the \(\rw_i\) the part (I) is always zero by \eqref{eq: param
and w is orthogonal}. This is why we defined \(\rw_{[\dimV_\timestep:\timestep)}\)
to be orthogonal to \(V_\timestep\) in Definition~\ref{def: previsible
orthonormal coordinate system}. Since we also defined our basis to be an
orthonormal basis, (II) is only non-zero when covariances of the directional
derivatives in the same direction are taken.

As we treat our \(Z\) matrix \eqref{eq: main induction assumption} as row-major,
the directional derivatives are grouped by direction \(\rw_i\) in \(Z\), which therefore
results in the following block matrix structure.
\begin{align}
    \nonumber
    \tfrac1\dims\magenta{\rcov_{[0:\timestep)}^{\dims}}
    &=\C_Z(\rv_{[0:\dimV_\timestep)}, \rw_{[\dimV_\timestep:\dims)}; \Param_{[0:\timestep)})
    \\
    \nonumber
    &= \begin{pmatrix}
        \C_Z(\rv_{[0:\dimV_\timestep)}; \Param_{[0:\timestep)})
        \\
        & \C_{D_{\rw_{\dimV_\timestep}}\rf_\dims}(\Param_{[0:\timestep)})
        \\
        & & \ddots 
        \\
        & & & \C_{D_{\rw_{\dims-1}}\rf_\dims}(\Param_{[0:\timestep)})
    \end{pmatrix}
    \\
    &=: \frac1\dims\left(\begin{array}{c c c c}
        \cline{1-1}
        \multicolumn{1}{|c|}{}
        \\
        \multicolumn{1}{|c|}{\quad\magenta{\rcov^{\black{v},\dims}_{[0:\timestep)}}\quad} 
        \\
        \multicolumn{1}{|c|}{}
        \\
        \cline{1-2}
        &\multicolumn{1}{|c|}{\magenta{\rcov^{\black{w},\dims}_{[0:\timestep)}}}
        \\
        \cline{2-2}
        & & \ddots
        \\
        \cline{4-4}
        & & &\multicolumn{1}{|c|}{\magenta{\rcov^{\black{w},\dims}_{[0:\timestep)}}}
        \\
        \cline{4-4}
    \end{array}\right),
    \label{eq: block matrix autocovariance}
\end{align}
For \(\magenta{\rcov^{\black{w},\dims}_{[0:\timestep)}}\) to be well defined, we
require the blocks \(\C_{D_{\rw_{i}}\rf_\dims}(\Param_{[0:\timestep)})\) to not depend
on \(i\). But since (I) is zero and \(\langle \rw_i, \rw_i\rangle =1\) we have
by \eqref{eq: cov df, df} of Lemma~\ref{lem: cov of derivatives, non-stationary isotropic case}
\begin{align*}
    &\C_{D_{\rw_i}\rf_\dims}(\Param_{[0:\timestep)})
    \\
    &= \!\frac1\dims\!\! \begin{pmatrix}
        \kernel_3\bigl(\tfrac{\|\Param_0\|^2}2, \tfrac{\|\Param_0\|^2}2, \langle \Param_0, \Param_0\rangle\bigr)
        &\!\!\!\!\cdots\!\!\!\!\!
        & \kernel_3\bigl(\tfrac{\|\Param_0\|^2}2, \tfrac{\|\Param_{\timestep-1}\|^2}2, \langle \Param_0, \Param_{\timestep-1}\rangle\bigr)
        \\
        \vdots & & \vdots
        \\
        \kernel_3\bigl(\tfrac{\|\Param_{\timestep-1}\|^2}2, \tfrac{\|\Param_0\|^2}2, \langle \Param_{\timestep-1}, \Param_0\rangle\bigr)
        &\!\!\!\!\cdots\!\!\!\!\!
        & \kernel_3\bigl(\tfrac{\|\Param_{\timestep-1}\|^2}2, \tfrac{\|\Param_{\timestep-1}\|^2}2, \langle \Param_{\timestep-1}, \Param_{\timestep-1}\rangle\bigr)
    \end{pmatrix}
    \\
    &=: \tfrac1\dims\magenta{\rcov^{\black{w},\dims}_{[0:\timestep)}}.
\end{align*}
In particular there is no dependence on \(\rw_i\) so
\(\magenta{\rcov^{\black{w},\dims}_{[0:\timestep)}}\) is well defined. Since
these block matrices are of constant size, they converge if
all their finitely many entries converge. But we already argued that the entries
converge (in the previous paragraph) and we therefore have\footnote{
    The increasing number of identical
    \(\magenta{\rcov^{\black{w},\dims}_{[0:\timestep)}}\) will drive the law of large numbers of
    \[
        \|P_{V_\timestep^\perp}\!\nabla\rf_\dims(\param_0)\|^2 =\!\! \sum_{i=\dimV_{\timestep}}^{\dims-1}\! (D_{\rw_i}\rf_\dims(\param_0))^2
    \]
    similar to the first step (cf.~Equation~\eqref{eq: LLN induction start}). At the moment they are still matrices but this will change when
    combined with the mixed covariances \(\blue{\rcov^{\black{w}, \dims}_{[0:\timestep),\timestep}}\)
    finally resulting in \eqref{eq: limiting residual variance}.
}
\[
    \magenta{\rcov^{\black{v},\dims}_{[0:\timestep)}}
    \underset{\dims\to\infty}{\overset{p}\to}
    \magenta{\Sigma^{\black{v},\infty}_{[0:\timestep)}}
    \in \real^{\timestep (\dimV_\timestep+1) \times \timestep (\dimV_\timestep+1)}
    \quad \text{and} \quad
    \magenta{\rcov^{\black{w},\dims}_{[0:\timestep)}}
    \underset{\dims\to\infty}{\overset{p}\to} \magenta{\Sigma^{\black{w},\infty}_{[0:\timestep)}}
    \in \real^{\timestep \times \timestep}.
\]
For the mixed covariance we similarly have
\begin{align}
    \nonumber
    &\tfrac1\dims\blue{\rcov^\dims_{[0:\timestep),\timestep}}
    \\
    \nonumber
    &=\C_Z(\rv_{[0:\dimV_\timestep)},\rw_{[\dimV_\timestep:\dims)};\Param_{[0:\timestep)}, \Param_\timestep)
    \\
    \nonumber
    &= \begin{pmatrix}
        \C_Z(\rv_{[0:\timestep)};\Param_{[0:\timestep)},\Param_\timestep)
        \\
        & \C_{D_{w_\timestep}\rf_\dims}(\Param_{[0:\timestep)}, \Param_\timestep)
        \\
        & & \ddots
        \\
        & & & \C_{D_{w_{\dims-1}}\rf_\dims}(\Param_{[0:\timestep)},\Param_\timestep)
    \end{pmatrix}
    \\
    &=: \frac1\dims\left(\begin{array}{c c c c}
        \cline{1-1}
        \multicolumn{1}{|c|}{}
        \\
        \multicolumn{1}{|c|}{}
        \\
        \multicolumn{1}{|c|}{\;\blue{\rcov^{\black{v},\dims}_{[0:\timestep),\timestep}}\;} 
        \\
        \multicolumn{1}{|c|}{}
        \\
        \multicolumn{1}{|c|}{}
        \\
        \cline{1-2}
        &\multicolumn{1}{|c|}{}
        \\
        &\multicolumn{1}{|c|}{\blue{\rcov^{\black{w},\dims}_{[0:\timestep),\timestep}}}
        \\
        &\multicolumn{1}{|c|}{}
        \\
        \cline{2-2}
        & & \ddots
        \\
        \cline{4-4}
        & & &\multicolumn{1}{|c|}{}
        \\
        & & &\multicolumn{1}{|c|}{\blue{\rcov^{\black{w},\dims}_{[0:\timestep),\timestep}}}
        \\
        & & &\multicolumn{1}{|c|}{}
        \\
        \cline{4-4}
    \end{array}\right),
    \label{eq: block matrix mixed cov}
\end{align}
with a similar argument why \(\blue{\rcov^{\black{w},\dims}_{[0:\timestep),\timestep}}\) is well defined as for
\(\magenta{\rcov^{\black{w},\dims}_{[0:\timestep)}}\).
And again by the discussion of the previous paragraph establishing convergence
of the entries, we have that these block matrices of constant size converge
\[
    \blue{\rcov^{\black{v},\dims}_{[0:\timestep),\timestep}}
    \underset{\dims\to\infty}{\overset{p}\to}
    \blue{\Sigma^{\black{v},\infty}_{[0:\timestep),\timestep}}
    \in \real^{\timestep(\dimV_\timestep+1) \times (\dimV_\timestep+1)}
    \quad\text{and}\quad
    \blue{\rcov^{\black{w},\dims}_{[0:\timestep),\timestep}}
    \underset{\dims\to\infty}{\overset{p}\to}
    \blue{\Sigma^{\black{w},\infty}_{[0:\timestep),\timestep}}
    \in \real^{\timestep \times 1}.
\]
We can also split up the autocovariance \(\green{\rcov^{\dims}_\timestep}\) in a similar fashion with
\[
    \green{\rcov^{\black{v},\dims}_\timestep}
    \underset{\dims\to\infty}{\overset{p}\to}
    \green{\Sigma^{\black{v},\infty}_\timestep}
    \in \real^{(\dimV_\timestep +1) \times (\dimV_\timestep+1)}
    \quad\text{and}\quad
    \green{\rcov^{\black{w},\dims}_\timestep}
    \underset{\dims\to\infty}{\overset{p}\to}
    \green{\Sigma^{\black{w},\infty}_\timestep}
    \in \real^{1 \times 1}.
\]
Finally, the expectation functions can also be split into a block containing
the directional derivatives \(\rv_{[0:\dimV_\timestep)}\) spanning \(V_\timestep\)
and the expectations of directional derivatives in the \(\rw_i\) directions.
More specifically we have
\begin{equation}
    \label{eq: expectation v} 
    \begin{aligned}
    \mu_\timestep^{v,\dims} = \mu_Z(\rv_{[0:\dimV_\timestep)}; \Param_\timestep)
    \underset{\dims\to\infty}&{\overset{p}\to}
    \mu_\timestep^{v, \infty}\in\real^{\dimV_\timestep +1}
    \qquad\text{and}\quad
    \\
    \mu_{[0:\timestep)}^{v,\dims}
    \underset{\dims\to\infty}&{\overset{p}\to}
    \mu_{[0:\timestep)}^{v, \infty}\in\real^{\timestep(\dimV_\timestep +1)},
    \end{aligned}
\end{equation}
which converge since the inner products and norms used in \eqref{eq: expect df}
converge.
Since the directions \(\rw_i\) are selected orthogonal to all evaluation
points \(\Param_\timestep\in V_\timestep\), the inner product in the expectation
of the directional derivative \eqref{eq: expect df} is always zero
and we therefore have
\begin{equation}
    \label{eq: expectation w}
    \mu_\timestep^{w,\dims} = \mu_Z(\rw_{[\dimV_\timestep:\dims)}; \Param_\timestep)
    = 0 \in\real
    \quad\text{and}\quad
    \mu_{[0:\timestep)}^{w,\dims}
    = 0 \in\real^{\timestep}.
\end{equation}

\paragraph*{Convergence of the conditional expectation}

Let us take a step back and review what we have done and want to do. We have
argued that \(Z(\Basis_\timestep; \Param_{[0:\timestep)})\) is conditionally Gaussian and that it
is therefore enough to understand the conditional expectation and covariance
(cf.~\eqref{eq: decomposition in distribution}). We have then applied a
well known result about the conditional distribution of Gaussian random
variables to obtain explicit formulas for the conditional expectation \eqref{eq: cond expect}
and covariance \eqref{eq: cond variance} made up of unconditional covariance
matrices. We proved that their entries converged and split these covariance
matrices into block form, such that these blocks of constant size converge in
probability. What is left to do, is to put these results together to prove
convergence results about \(Z(\Basis_\timestep; \Param_{[0:\timestep)})\). We start with the
conditional expectation and then get to \(Z(\Basis_\timestep; \Param_{[0:\timestep)})\) itself
by a consideration of the conditional covariance.

So let us take a look at the conditional expectation.
Since the \(\rw_i\) are by definition orthonormal to the previsible running span
of evaluation points \(V_\timestep\) \eqref{eq: vector space of evaluation
points}, we have \(D_{\rw_k}\rf_\dims(\Param_{[0:\timestep)})=0\). By applying the
block matrix structure \eqref{eq: block matrix autocovariance}, \eqref{eq: block
matrix mixed cov} to the the representation of the conditional expectation
\eqref{eq: cond expect} we therefore get using \eqref{eq: expectation w}
\begin{align}
\nonumber
    \E[D_{\rw_i}\rf_\dims(X_\timestep)\mid \filt_{\timestep-1}]
    &= \mu_\timestep^{w,\dims}
    + \blue{\rcov^{\black{w},\dims}_{[0:\timestep),\timestep}}^\transpose
    [\magenta{\rcov^{\black{w},\dims}_{[0:\timestep)}}]^{-1}
    (D_{\rw_k}\rf_\dims(\Param_{[0:\timestep)}) - \mu_{[0:\timestep)}^{w,\dims})
    \\
\label{eq: conditional expectation in direction w}
    &= 0.
\end{align}
The only interesting part of the vector \(\E[Z(\rv_{[0:\dimV_\timestep)},
\rw_{[\dimV_\timestep:\dims)}; \Param_\timestep) \mid \filt_{\timestep-1}]\) are therefore the first
\(\dimV_\timestep\) entries \(\E[Z(\rv_{[\dimV_\timestep:\timestep)}; \Param_\timestep) \mid
\filt_{\timestep-1}]\). For this we use the formula for the conditional expectation
\eqref{eq: cond expect} and our block matrix decompositions \eqref{eq: block
matrix autocovariance}, \eqref{eq: block matrix mixed cov} again to obtain 
\begin{equation}
    \label{eq: relevant part of cond expec}
    \E[Z(\rv_{[0:\dimV_\timestep)}; \Param_\timestep) \mid \filt_{\timestep-1}]
    = \mu_\timestep^{(v,\dims)}
    +\blue{\rcov^{\black{v},\dims}_{[0:\timestep),\timestep}}^\transpose
    [\magenta{\rcov^{\black{v},\dims}_{[0:\timestep)}}]^{-1}
    (Z(\rv_{[0:\dimV_\timestep)}; \Param_{[0:\timestep)})-\mu_{[0:\timestep)}^{(v,\dims)}).
\end{equation}
Since \(Z(\rv_{[0:\dimV_\timestep)}; \Param_{[0:\timestep)})\) converges in probability by the
induction assumptions \ref{ind: information convergence}  summarized in
\eqref{eq: main induction assumption}, and since we have spent the previous
paragraphs proving that the covariance matrices and expectations converge, it
almost seems like we are done. But while the conditional expectation \eqref{eq:
cond expect} can handle non-invertible matrices with a generalized inverse
(Theorem~\ref{thm: conditional gaussian distribution}), an application of
continuous mapping requires continuity and the matrix inverse is only continuous
at invertible matrices.

\begin{lemma}
    \label{lem: strict positive cov matrix}
    \(\magenta{\Sigma^{\black{v},\infty}_{[0:\timestep)}}\) is strictly positive definite and
    therefore invertible.
\end{lemma}
\begin{proof}
    The essential ingredients are that \((\rf_\dims, \nabla\rf_\dims)\) is strictly positive
    definite by assumption of Theorem~\ref{THM: ASYMPTOTICALLY DETERMINISTIC
    BEHAVIOR VARIANT} and that the evaluation points are asymptotically
    different \(\rho_{ij} >0\) by
    \ref{ind: evaluation points asymptotically different}. The technical
    complications are that
    \(\magenta{\Sigma^{\black{v},\infty}_{[0:\timestep)}}\) is only defined as
    the limit of covariance matrices (which are themselves not necessarily
    strictly positive definite themselves) and that this limit involves changing
    the domain of \((\rf_\dims, \nabla\rf_\dims)\) as we increase its dimension.

    Let us first address the problem of the changing domain.
    Since we are only interested in the block matrix of the derivatives
    in the directions contained in \(V_\timestep\) and the evaluation points
    \(\Param_{[0:\timestep)}\) are also contained in \(V_\timestep\) we can map
    them via an isometry to \(\real^{\dimV_\timestep}\) which retains all distances and inner
    products and therefore retains the covariance matrices
    \(\magenta{\rcov^{\black{v},\dims}_{[0:\timestep)}}\) (cf.\ Lemma~\ref{lem: cov of derivatives, non-stationary isotropic case}).
    Note that the dimensional scaling \(\frac1\dims\) is not an issue here, as
    we have already removed the scaling from
    \(\magenta{\rcov^{\black{v},\dims}_{[0:\timestep)}}\).
    Therefore \(\magenta{\rcov^{\black{v},\dims}_{[0:\timestep)}}\) can be viewed as a sequence
    of covariance matrices of \((\rf_\dims,\nabla\rf_\dims)\) with \(\rf_\dims:\real^{\dimV_\timestep}\to \real\).
    This solves the problem of the changing domain.

    Finally, we need that the limiting matrix
    \(\magenta{\Sigma^{\black{v},\infty}_{[0:\timestep)}}\) is in fact a covariance matrix of
    \((\rf_\dims, \nabla\rf_\dims)\) and not just a limit of covariance
    matrices. For this we apply \ref{ind: representation}, which
    ensures that the limiting scalar products \eqref{eq: <X, v> converges to gamma}
    and distances \eqref{eq: limiting distances} which make up
    \(\magenta{\Sigma^{\black{v},\infty}_{[0:\timestep)}}\) (cf.~Lemma~\ref{lem:
    cov of derivatives, non-stationary isotropic case}) are realizable by points
    \(y_0,\dots, y_{\timestep-1}\in\real^{\dimV_\timestep}\).
    \(\magenta{\Sigma^{\black{v},\infty}_{[0:\timestep)}}\) is therefore a
    covariance matrix of \((\rf_\dims, \nabla\rf_\dims)\) viewed as a random function with
    domain \(\real^{\dimV_\timestep}\). 
    Since the distances \(\rho_{ij}\) are positive by \ref{ind: evaluation
    points asymptotically different}, the asymptotic representations \(y_k\) are
    distinct and their covariance matrix
    \(\magenta{\Sigma^{\black{v},\infty}_{[0:\timestep)}}\) is therefore strictly
    positive definite by the strict positive definiteness of \((\rf_\dims, \nabla\rf_\dims)\).
\end{proof}
The matrix inverse is therefore continuous in
\(\magenta{\Sigma^{\black{v},\infty}_{[0:\timestep)}}\) and by the convergence
of the matrices in probability and the convergence of \(Z(\rv_{[0:\dimV_\timestep)};
\Param_{[0:\timestep)})\) by induction assumptions \ref{ind: information convergence} for
\(\timestep-1\) restated in \eqref{eq: main induction assumption}, the
conditional expectation \eqref{eq: relevant part of cond expec} converges in
probability
\begin{align}
    \nonumber
    &\E[Z(\rv_{[0:\dimV_\timestep)}; \Param_\timestep) \mid \filt_{\timestep-1}]
    \\
    \label{eq: converging conditional expectation}
    \underset{\dims\to\infty}&{\overset{p}\to}
    \mu_\timestep^{v, \infty} + \blue{\Sigma^{v,\infty}_{[0:\timestep),\timestep-1}}^T
    [\magenta{\Sigma^{v,\infty}_{[0:\timestep)}}]^{-1}
    \Bigl(
        \begin{pmatrix}
            \limf_{[0:\timestep)}\\
            \gamma_{[0:\timestep)}
        \end{pmatrix}
        -\mu_{[0:\timestep)}^{v, \infty}
    \Bigr)
    =: \begin{pmatrix}
        \limf_\timestep
        \\
        \gamma_\timestep^{([0:\dimV_\timestep))}
    \end{pmatrix}.
\end{align}
Together with the conditional expectation of the directional derivatives in the
directions \(\rw_i\) in \eqref{eq: conditional expectation in direction w}, we
have now fully determined the conditional expectation.

An analysis of the conditional variance will immediately give us an understanding
of the distribution and we can therefore put these considerations together
to obtain convergence of \(Z\) itself.

\paragraph*{\ref{it: step 1 induction step}: Convergence of \(Z(\rv_{[0:\dimV_\timestep)}; \Param_\timestep)\)}

Recall, we outlined at the start of the induction that we would first prove
convergence of \(Z(\rv_{[0:\red{\dimV_\timestep})}; \Param_\timestep)\) (i.e.\
\ref{it: step 1 induction step}) before proving \(V_{\timestep+1}\) to
asymptotically have full rank (\ref{it: step 2 induction step}) and the
convergence of the new corner element
\(D_{\rv_{\dimV_\timestep}}\rf_\dims(\Param_\timestep)\) (\ref{it: step 3 induction
step}).

In \eqref{eq: converging conditional expectation} we obtained the limit of the
conditional expectation of \(Z(\rv_{[0:\dimV_\timestep)}; \Param_\timestep)\).
This limit is also going to be the limit of the new column
\(Z(\rv_{[0:\dimV_\timestep)}; \Param_\timestep)\) itself. 
To prove this we only need to control the variance. For this purpose we recall that the
unconditional variance is given by 
\[
    \green{\rcov^\dims_\timestep}
    = \begin{pmatrix}
        \green{\rcov^{\black{v},\dims}_\timestep}
        \\
        & &  \kernel_3\bigl(\tfrac{\|\Param_\timestep\|^2}2,\tfrac{\|\Param_\timestep\|^2}2, \|\Param_\timestep\|^2\bigr) \identity_{(\dims-\dimV_\timestep)\times (\dims-\dimV_\timestep)}
    \end{pmatrix}.
\]
The conditional covariance matrix is therefore given by
\begin{align}
    \nonumber
    &\Cov[
        Z(\rv_{[0:\dimV_\timestep)},\rw_{[\dimV_\timestep:\dims)}; \Param_\timestep)
        \mid \filt_{\timestep-1}
    ]
    \\
    \label{eq: conditional covariance}
    &= \frac1\dims \begin{bmatrix}
        \green{\rcov^{\black{v},\dims}_\timestep}
        - \blue{\rcov^{\black{v},\dims}_{[0:\timestep),\timestep}}^T
        [\magenta{\rcov^{\black{v},\dims}_{[0:\timestep)}}]^{-1}
        \blue{\rcov^{\black{v},\dims}_{[0:\timestep),\timestep}}
        \\
        & \sigma_{w,d}^2
        \identity_{(\dims-\dimV_\timestep)\times (\dims-\dimV_\timestep)}
    \end{bmatrix},
\end{align}
where the conditional variance in the direction of the \(\rw_i\) is given by
\begin{equation}
    \label{eq: residual variance}
    \sigma_{w,\dims}^2
    := \bigl(
        \green{\kernel_3\bigl(\tfrac{\|\Param_\timestep\|^2}2,\tfrac{\|\Param_\timestep\|^2}2, \|\Param_\timestep\|^2\bigr)}
        - \blue{\rcov^{\black{w},\dims}_{[0:\timestep),\timestep}}^T
        [\magenta{\rcov^{\black{w},\dims}_{[0:\timestep)}}]^{-1}
        \blue{\rcov^{\black{w},\dims}_{[0:\timestep),\timestep}}
    \bigr).
\end{equation}
We will need \(\sigma_{w,\dims}^2\) for \ref{it: step 2 induction step} and
\ref{it: step 3 induction step}, but for now it is not important.
Due to the diagonal structure we have for our new column
\(Z(\rv_{[0:\red{\dimV_\timestep})}; \Param_\timestep)\) by \eqref{eq:
decomposition in distribution}
\begin{align}
    \nonumber
    &Z(\rv_{[0:\red{\dimV_\timestep})}, \Param_\timestep)
    \\
    \nonumber
    &= \E[Z(\rv_{[0:\dimV_\timestep)}, \Param_\timestep) \mid \filt_{\timestep-1}]
    + \sqrt{\Cov[Z(\rv_{[0:\dimV_\timestep)}, \Param_\timestep) \mid \filt_{\timestep-1}]} \begin{pmatrix}
        Y_0\\ \vdots\\ Y_{\dimV_\timestep}
    \end{pmatrix}
    \\
    \nonumber
    &= \E[Z(\hat{v}_{[0:\dimV_\timestep)}, \Param_\timestep) \mid \filt_{\timestep-1}]
    + \frac1{\sqrt{\dims}} \sqrt{
        \green{\rcov^{\black{v},\dims}_\timestep}
        - \blue{\rcov^{\black{v},\dims}_{[0:\timestep),\timestep}}^T
        [\magenta{\rcov^{\black{v},\dims}_{[0:\timestep)}}]^{-1}
        \blue{\rcov^{\black{v},\dims}_{[0:\timestep),\timestep}}
        }\begin{pmatrix}
        Y_0\\ \vdots\\ Y_{\dimV_\timestep}
    \end{pmatrix}
    \\
    &\underset{\dims\to\infty}{\overset{p}\to}
    \begin{pmatrix}
        \limf_\timestep
        \\
        \gamma_\timestep^{([0:\dimV_\timestep))}
    \end{pmatrix}.
    \label{eq: convergence new col}
\end{align}
This is exactly the convergence of the `new column' required for \ref{it: step 1 induction step}.

\paragraph*{\ref{it: step 2 induction step}: The rank of \(V_{\timestep+1}\)}

The last two steps follow fairly quickly. Recall that \(V_{\timestep+1}\) is the
\emph{previsible} running span of evaluation points defined in \eqref{eq: vector
space of evaluation points}. Since it is previsible, it only includes gradients
up to time \(\timestep\) and the Gram-Schmidt candidate \eqref{eq: definition of
v candidate} of its most recent addition is therefore given by
\[
    \tilde{\rv}_\timestep := \nabla\rf_\dims(\Param_\timestep) - P_{V_\timestep} \nabla\rf_\dims(\Param_\timestep)
    = \sum_{k=\dimV_\timestep}^{\dims-1} \langle \rw_i, \nabla\rf_\dims(X_\timestep)\rangle \rw_i
    = \sum_{k=\dimV_\timestep}^{\dims-1} D_{\rw_i}\rf_\dims(X_\timestep) \rw_i,
\]
where \(P_{V_\timestep}\) is the projection onto \(V_\timestep\).

Now we naturally want to analyze the directional derivatives
\(D_{\rw_i}\rf_\dims(X_\timestep)\) which make up \(\tilde{\rv}_\timestep\).
By representation \eqref{eq: decomposition in distribution} and our
formula for the conditional covariance \eqref{eq: conditional covariance}
we have in distribution
\[
    D_{\rw_i}\rf_\dims(\Param_\timestep)
    = \underbrace{\E[D_{\rw_i}\rf_\dims(\Param_\timestep)\mid \filt_{\timestep-1}]}_{\overset{\eqref{eq: conditional expectation in direction w}}=0}
    + \sqrt{\tfrac1\dims \sigma_{w,\dims}} Y_{i+1}.
\]
We now already see where our law of large numbers is going to come from. But we
first need to take a closer look at the residual variance \(\sigma_{w,\dims}^2\) defined
in \eqref{eq: residual variance}. We get convergence in probability of the residual variance
to a value strictly greater zero
\begin{align}
    \label{eq: limiting residual variance}
    \sigma_{w,\dims}^2
    &= \bigl(
        \green{\kernel_3\bigl(\tfrac{\|\Param_\timestep\|^2}2,\tfrac{\|\Param_\timestep\|^2}2, \|\Param_\timestep\|^2\bigr)}
        - \blue{\rcov^{\black{w},\dims}_{[0:\timestep),\timestep}}^T
        [\magenta{\rcov^{\black{w},\dims}_{[0:\timestep)}}]^{-1}
        \blue{\rcov^{\black{w},\dims}_{[0:\timestep),\timestep}}
    \bigr)
    \\
    \nonumber
    \underset{\dims\to\infty}&{\overset{p}\to}
    \bigl(
        \green{\kernel_3\bigl(\tfrac{\|y_\timestep\|^2}2,\tfrac{\|y_\timestep\|^2}2, \|y_\timestep\|^2\bigr)}
        - \blue{\Sigma^{\black{w},\infty}_{[0:\timestep),\timestep}}^T
        [\magenta{\Sigma^{\black{w},\infty}_{[0:\timestep)}}]^{-1}
        \blue{\Sigma^{\black{w},\infty}_{[0:\timestep),\timestep}}
    \bigr)
    =:\sigma_{w,\infty}^2 > 0,
\end{align}
using the convergence of the block matrices and the following lemma.

\begin{lemma}\label{lem: strict postive definite w direction}
    \(\magenta{\Sigma^{\black{w},\infty}_{[0:\timestep)}}\) is strictly positive definite
    and \(\sigma_{w,\infty}^2>0\).
\end{lemma}
\begin{proof}
    We are going to show with a similar argument as in Lemma~\ref{lem: strict
    positive cov matrix} that the matrix
    \[
        \Sigma^{w,\infty}_{[0:\timestep]}
        := \begin{bmatrix}
            \magenta{\Sigma^{\black{w},\infty}_{[0:\timestep)}}     
            & \blue{\Sigma^{\black{w},\infty}_{[0:\timestep),\timestep}}
            \\
            \blue{\Sigma^{\black{w},\infty}_{[0:\timestep),\timestep}}^T
            & \green{\kernel_3\bigl(\tfrac{\|y_\timestep\|^2}2,\tfrac{\|y_\timestep\|^2}2, \|y_\timestep\|^2\bigr)}
        \end{bmatrix}
    \]
    is strictly positive definite. Before we do so, let us quickly argue why this
    finishes the proof. Since \(\Sigma^{w,\infty}_{[0:\timestep]}\) is then strictly positive 
    definite, it has a cholesky decomposition \(L\) such that
    \[
        \Sigma^{w,\infty}_{[0:\timestep]} = LL^T
        = \begin{bmatrix}
            L_\timestep & 0
            \\
            l^T & \sigma
        \end{bmatrix}
        \begin{bmatrix}
            L_\timestep^T & l
            \\
            l &\sigma
        \end{bmatrix},
    \]
    which implies that \(L_\timestep\) is the cholesky decomposition of \(\magenta{\Sigma^{\black{w},\infty}_{[0:\timestep)}}\),
    \(l = L_\timestep^{-1}\blue{\Sigma^{\black{w},\infty}_{[0:\timestep),\timestep}}\) and
    \[
        \sigma = \sqrt{ 
            \green{\kernel_3\bigl(\tfrac{\|y_\timestep\|^2}2,\tfrac{\|y_\timestep\|^2}2, \|y_\timestep\|^2\bigr)}
            - \blue{\Sigma^{\black{w},\infty}_{[0:\timestep),n}}^T
            [\magenta{\Sigma^{\black{w},\infty}_{[0:\timestep)}}]^{-1}
            \blue{\Sigma^{\black{w},\infty}_{[0:\timestep),n}}
        }
        = \sigma_{w,\infty}.
    \]
    Since we have
    \[
        \det(L) = \sigma_{w,\infty}\det(L_\timestep),
    \]
    strict positive definiteness of \(\Sigma^{w,\infty}_{[0:\timestep]}\) and therefore \(0\neq\det(\Sigma^{w,\infty}_{[0:\timestep]}) = \det(L)^2\)
    implies
    \[
        \sigma_{w,\infty}^2 > 0
        \quad \text{and} \quad 
        \det(L_\timestep)\neq 0.
    \]
    But \(\det(L_\timestep)\neq 0\) also implies that
    \(\magenta{\Sigma^{\black{w},\infty}_{[0:\timestep)}}\) has to be strictly
    positive definite, since \(L_\timestep\) is its cholesky decomposition.

    What is therefore left to prove is the strict positive definiteness of
    \(\Sigma^{w,\infty}_{[0:\timestep]}\). For this note that
    \[
       \rcov^{w,\dims}_{[0:\timestep]} 
        := \begin{bmatrix}
            \magenta{\rcov^{\black{w},\dims}_{[0:\timestep)}}     
            & \blue{\rcov^{\black{w},\dims}_{[0:\timestep),n}}
            \\
            \blue{\rcov^{\black{w},\dims}_{[0:\timestep),n}}^T
            & \green{\kernel_3\bigl(\tfrac{\|\Param_\timestep\|^2}2,\tfrac{\|\Param_\timestep\|^2}2, \|\Param_\timestep\|^2\bigr)}
        \end{bmatrix}
    \]
    is the plug-in covariance matrix of \(D_{\rw_i}\rf_\dims(\Param_{[0:\timestep]})\).
    Where we use the term "plug-in" covariance to say that we treat the
    evaluation points \(\Param_{[0:\timestep]}\) and the direction \(\rw_i\) as
    deterministic. Since \(\Param_{[0:\timestep]}\) are contained in \(V_\timestep\),
    we again map them isometrically to \(\real^{\dimV_\timestep}\). But this time
    we view \(\real^{\dimV_\timestep}\) as a subspace of
    \(\real^{\dimV_\timestep+1}\) and map the additional vector \(\rw_i\) to
    \(e_{\dimV_\timestep+1}\) such that \(\rcov^{w,\dims}_{[0:\timestep]}\)
    is a covariance matrix of \(\nabla\rf_\dims\) in
    \(\real^{\dimV_\timestep+1}\).
    
    Now we finish the proof with the same limiting argument as in
    Lemma~\ref{lem: strict positive cov matrix}.  Since
    \(\nabla\rf_\dims\) is strictly positive definite by
    assumption, we get that \(\Sigma_{[0:\timestep]}^{w,\infty}\) is strictly
    positive definite as the covariance matrix of
    \((\partial_{\dimV_\timestep+1}\rf_\dims)\) at the points
    \(y_{[0:\timestep]} = (y_0, \dots, y_\timestep) \subseteq
    \real^{\dimV_\timestep+1}\) of \ref{ind: representation}. These are the
    limiting representations and non of the points \(y_k\) are equal by
    \ref{ind: evaluation points asymptotically different}.
\end{proof}

Using the convergence of the residual variance \(\sigma_{w,\dims}^2\to
\sigma_{w,\infty}^2\) allows us to finally make our law of large
numbers argument
\begin{equation}
    \label{eq: v_n tilde}
    \|\tilde{\rv}_\timestep\|^2
    = \sum_{i=\dimV_\timestep}^{\dims-1} (D_{\rw_i}\rf_\dims(X_\timestep))^2
    = \frac{\sigma_{w,\dims}^2}\dims\sum_{i=\dimV_\timestep+1}^{\dims}  Y_i^2 
    \;\underset{\dims\to\infty}{\overset{p}\to} \;
    \sigma_{w,\infty}^2  > 0.
\end{equation}
This implies that \(V_\timestep\) has full rank asymptotically \ref{ind: full
dimension}.

Assuming the last gradient is always used (and not just in the asymptotic limit),
we can get \(\dimV_{\timestep+1} = \dimV_\timestep+1\) almost surely,
since \(\sigma_{w,\dims}^2> 0\) almost surely is sufficient for
\(\|\tilde{\rv}_\timestep\|^2>0\) almost surely. The use of the
most recent gradient ensures the points \(\Param_k\) are different by an inductive
argument similar to Lemma~\ref{lem: conv info, full rank -> different eval pts}.
This in turn ensures, using the strict positive definiteness of
\((\rf_\dims,\nabla\rf_\dims)\), that \(\sigma_{w,\dims}^2> 0\) almost surely with a similar
argument as in Lemma~\ref{lem: strict postive definite w direction}.

\paragraph*{\ref{it: step 3 induction step}: Convergence of the new corner element}

The last step follows immediately by definition of \(\rv_{\dimV_\timestep} =
\frac{\tilde{\rv}_\timestep}{\|\tilde{\rv}_\timestep\|}\) in
Definition~\ref{def: previsible orthonormal coordinate system}
and \eqref{eq: v_n tilde}
\[
    D_{\rv_{\dimV_\timestep}}\rf_\dims(X_\timestep)
    = \bigl\langle
        \nabla\rf_\dims(X_\timestep),
        \tfrac{\tilde{\rv}_\timestep}{\|\tilde{\rv}_\timestep\|}
    \bigr\rangle
    = \|\tilde{\rv}_\timestep\|
    \underset{\dims\to\infty}{\overset{p}\to}
    \sigma_{w,\infty} =: \gamma_{\timestep}^{(\dimV_\timestep)}
    > 0.
\]

 \section{Outlook}\label{sec: outlook}

In this section we want to discuss possible generalizations and
extensions of our results. Possible generalizations include
\begin{enumerate}[label={G\arabic*.},ref={G\arabic*}]
	\item\label{it: strict positive def removal}
	A \textbf{removal of the strict positive definiteness} assumption.
	The strict positive definiteness was used in the proof of Theorem~\ref{THM:
	ASYMPTOTICALLY DETERMINISTIC BEHAVIOR VARIANT} to show that the limiting
	covariance matrices are invertible. It also ensured the full rank of \(V_\timestep\)
	\ref{ind: full dimension}, which in turn was used with the asymptotic use
	of the last gradient to ensure asymptotically different evaluation points
	\ref{ind: evaluation points asymptotically different}. The asymptotically
	different evaluation were used in turn to keep the covariance matrices
	strictly positive definite in the next iteration.

	A removal of the strict positive definiteness assumption will therefore
	remove the need for the assumption of the asymptotic use of the most
	recent gradient, but it will also lose the implications \ref{ind: full
	dimension} and \ref{ind: evaluation points asymptotically different}.

	There are two avenues to remove this assumption:
	\begin{enumerate}
		\item the use of \emph{generalized matrix inverses}, since the Gaussian
		conditional distribution can be formed with these (cf.
		Theorem~\ref{thm: conditional gaussian distribution}). The key
		would be to show that these generalized inverses would still
		converge with \(\dims\to\infty\).

		\item a \emph{perturbation argument}. One can slightly perturb the
		function \(\rf_\dims\) into \(\rf_{\dims,\epsilon}\) such that
		\((\rf_{\dims,\epsilon},\nabla\rf_{\dims,\epsilon})\) is strictly
		positive definite (by adding an independent random function with this
		property with an \(\epsilon\) weight). We then have \(\lim_{\epsilon\to
		0} \rf_{\dims,\epsilon} = \rf_\dims\).
		The difficulty is now to argue that we can exchange limits, as we are interested
		in 
		\[
			\lim_{\dims\to\infty}\rf_\dims
			= \lim_{\dims\to\infty}\lim_{\epsilon\to 0}\rf_{\dims, \epsilon}
		\]
		but want to consider \(\lim_{\epsilon\to
		0}\lim_{\dims\to\infty}\rf_{\dims, \epsilon}\), where we can apply the
		existing theory for the strict positive definite case to the inner
		limit.
	\end{enumerate}

	\item Removal of the \textbf{Gaussian} assumption. Our argument is
	essentially based on a law of large numbers applied to the norm of
	gradients. Since uncorrelated directional derivatives (obtained from
	isotropy) are sufficient for independence in the Gaussian case, the squared
	directional derivatives are also uncorrelated and the law of large numbers
	applies to the gradient norms. In the non-Gaussian case, the uncorrelated
	directional derivatives no longer guarantee uncorrelated squares and a more
	sophisticated argument is necessary.

	\item Relaxation of the \textbf{isotropy} assumption. In view of the
	general definition of input invariance (Definition~\ref{def: distributional
	input invariance}) it is natural to ask how small the set of invariant
	transformations can be before our results break down. In place of
	rotation invariance, one might consider exchangeable directions
	for example (which is captured by the set of dimension permutations).
	Since the partial derivatives would then be exchangeable random variables,
	their norm should still converge by laws of large numbers for exchangeable
	random variables. Unfortunately, our custom adapted coordinate
	system would not work in this more general framework. Thus a different
	approach for a proof would be necessary.
\end{enumerate}

Perhaps even more interesting than generalizations of our results might be the
pursuit of the following extensions
\begin{enumerate}[label={E\arabic*.},ref={E\arabic*}]
	\item\label{it: explicit derivation} The derivation of more \textbf{explicit representations} for the
	limiting information \(\liminfo_\timestep\) (which includes limiting
	function values \(\limf_\timestep=\limf_\timestep(\gsa, \mu, \kernel,
	\radius)\)).

	\item Ideally, the explicit representation \ref{it: explicit derivation}
	would allow for \textbf{convergence proofs} and \textbf{meta optimization}
	over the optimizer \(\gsa\) such that the best optimization algorithm for
	high dimensional problems can be found.

	\item An analysis of \textbf{stochastic gradient descent} and similar
	stochastic algorithms which do not have access to the true cost but
	rather noisy evaluations of \(\rf_\dims\). Since every evaluation is
	equipped with independent noise, this analysis is most likely similar
	but more difficult than the perturbation argument suggested to
	remove the strict positive definiteness assumption in \ref{it: strict positive def removal}.

	\item A very difficult extension would capture algorithms with
	\textbf{component-wise learning rates} \autocite[such as Adam
	by][]{kingmaAdamMethodStochastic2015}. Since component-wise learning rates
	cannot easily be expressed in a dimension-free manner, it is unclear how to
	argue that the algorithm converges in a sense. That is because the
	convergence of the algorithm follows from the key assumption of continuity
	in the dimensionless information in our proof.  Nevertheless, practical
	experiments confirm that these algorithms exhibit the same behavior (cf.\
	Adam in Figure~\ref{fig: self loss plots}). Perhaps an analysis of the Hessian
	with tools from random matrix theory would be suitable for an attack on this
	problem.

	\item The derivation of \textbf{concentration bounds}. Given that the
	conditional variance is bounded by the unconditional variance
	which scales with \(\frac1\dims\), we would expect concentration
	inequalities to contract exponentially in the dimension \(\dims\). That is
	because both Gaussian random variables and Chi2 random variables (the
	squared gradient norms), contract exponentially in the variance.
\end{enumerate}

\section*{Acknowledgements}
\addcontentsline{toc}{section}{Acknowledgements}

This research was supported by the RTG 1953, funded by the German Research
Foundation (DFG). The authors would like to thank Yan Fyodorov for his
hospitality and helpful discussions at King's College. We thank him along with
Mark Sellke
and Antoine Maillard
for their time and valuable insights on spin glasses. 
\printbibliography[heading=subbibliography]

\begin{subappendices}
\section{Appendix}

Recall that \(\kernel_3>0\) was used in \eqref{eq: LLN induction start} to prove
that the gradient has positive length. This fact follows from strict positive
definiteness, which we are now going to prove.

\begin{lemma}
    \label{lem: kappa_3 positive} 
    Let \(\kernel\) be a (non-stationary) isotropic covariance kernel valid in all
    dimensions (equivalently in \(\ell^2\), cf. Lemma \ref{lem: valid in all
    dimensions}).  Let \(\rf\sim \normal(\mu, \kernel)\) and assume \((\rf,
    \nabla\rf)\) to be strict positive definite, then for any
    \(\param\in \real^\dims\)
    \[
        \kernel_3\bigl(\tfrac{\|\param\|^2}2, \tfrac{\|\param\|^2}2, \|\param\|^2\bigr)
        > 0.
    \]
\end{lemma}
\begin{proof}
    Since the norm \(\|\cdot\|^2\) is rotation invariant, we can assume without
    loss of generality \(x = \lambda e_1\) for the standard basis vector \(e_1
    \in \real^\dims\). Since \(\partial_2 \rf\) is strict positive definite, we know that
    \[
        0< \Var(\partial_2\rf(\param))
        = \Cov(\partial_2 \rf(\lambda e_1), \partial_2\rf(\lambda e_1))
        = \tfrac1\dims \kernel_3\bigl(\tfrac{\|\param\|^2}2, \tfrac{\|\param\|^2}2, \|\param\|^2\bigr),
    \]
    where we have used \eqref{eq: cov df, df} of Lemma~\ref{lem: cov of derivatives, non-stationary isotropic case}
    and \(\|\lambda e_1\| = \|\param\|\) in the last equation.
\end{proof}

In the stationary isotropic case, we can say more. We do not only have \(\kernel_3
= -\ikernel'(0)>0\), but we have \(\ikernel'(r) <0\) for all \(r\ge 0\).

\begin{lemma}
    \label{lem: covariance derivative negative}
    Let \(\ikernel\) be a stationary isotropic covariance kernel
    valid in all dimensions.
    If \(\rf\sim\normal(\mu, \ikernel)\) is \emph{not} almost surely
    constant, then \(\ikernel'(r)<0\) for all \(r\ge 0\).
\end{lemma}
\begin{proof}
    By Corollary~\ref{lem: constant random functions} we have
    \(\schoenbergMeas((0,\infty)) > 0\) for the Schoenberg measure
    \(\schoenbergMeas\) in Schoenbergs's characterization \eqref{eq: schoenberg
    charact}.  Thus by the continuity of measures there exists
    \(a,b>0\) such that \(\schoenbergMeas([a,b]) > 0\). Then by
    \eqref{eq: schoenberg charact} we have
    \begin{align*}
        -\ikernel'(r)
        &= \int_{[0,\infty)} t^2 \exp(-t^2r)\schoenbergMeas(dt)
        \\
        &\ge \int_{[a,b]} t^2 \exp(-t^2r)\schoenbergMeas(dt)
        \\
        &\ge \schoenbergMeas([a,b])\inf_{s\in [a,b]}s^2 \exp(-s^2r) > 0.
        \qedhere
    \end{align*}
\end{proof} 
\end{subappendices}
 		}
	\end{refsection}

	\part{Supervised machine learning}
	\label{part: supervised machine learning}

	\begin{refsection}
		\chapter{Random objective from exchangeable data}
\label{chap: random objective from exchangeable data}

In this chapter we explain how random objective functions
arise from exchangeable data in supervised machine learning.
While this motivation is specific to this application, it offers a more
elegant rationale compared to the assumption of random objectives out of
necessity, given the intractability of worst-case optimization. Additionally,
this derivation can inform distributional
assumptions about the random objective. In particular
we show for a random linear model in Section \ref{sec: random linear model} that
the stationarity assumption is untenable in contrast
to the isotropy assumption.

\begin{quote}
	\textit{``The biggest lesson that can be read from 70 years of AI research is that general
	methods that leverage computation are ultimately the most effective, and by a
	large margin."}\par\raggedleft--- \textup{Richard Sutton}, ``The Bitter Lesson''
\end{quote}

The supervised machine learning task is to find a model\footnote{
	The model \(\model=\model_\param\) is often parametrized by parameters \(\param\),
	such as the weights of an artificial neural network.
} \(\model\)
which best predicts the labels \(Y\) based on input data \(X\). For this purpose the
distance of the prediction \(\model(X)\) from the labels \(Y\) is measured
by a loss function \(\loss\). The goal is then to find a model \(\model\)
which minimizes the accumulated cost over the usage time of the model. Assuming the
model is used indefinitely after training the cost function\footnote{
	If the model is parametrized, then we define
	\(\Cost(\param) := \Cost(\model_\param)\).
} is thus given by
\[
	\Cost(\model) = \lim_{n\to \infty} \sum_{k=1}^n \loss(\model(X_k), Y_k).
\]
The existence of this limit is usually guaranteed by the (strong) law of large numbers, (S)LLN,
and the assumption of independent, identically distributed (iid) data
\((Z_k)_{k\in \nat} = (X_k, Y_k)_{k\in \nat}\). Instead, we will utilize
the more general assumption of \emph{exchangeable} data and explain why the iid assumption
is inconsistent with the typical choices of the loss function in machine learning.

We recall, that a sequence of random variables is called exchangeable if its distribution
is invariant under finite permutation, i.e. the order of the data does not matter.
By de~Finetti's seminal representation theorem \autocite[e.g.][Thm.~1.49]{schervishTheoryStatistics1995}
the sequence \((Z_k)_{k\in \nat}\) is exchangeable if and only if there exists a random
probability measure \(\rP_Z\) such that, conditional on \(\rP_Z\), the data \((Z_k)_{k\in \nat}\)
is iid with distribution \(\rP_Z\), i.e.
\[
	\Pr(Z_k \in A \mid \rP_Z) = \rP_Z(A).
\]
This also implies that the SLLN is applicable conditional on \(\rP_Z\), i.e.
\[
	\Cost(\model) = \E[ \loss(\model(X_k), Y_k) \mid \rP_Z].
\]
In essence, de Finetti's representation theorem states that the exchangeable setting
always corresponds to the two stage model: First, sample the distribution
\(\rP_Z\), and then, in the second stage, sample iid data \((Z_k)_{k\in \nat}\)
from this distribution.\footnote{
	Details in Section \ref{sec: parametric statistical model}.
}
Observe that unless \(\rP_Z\) is deterministic, which
corresponds to the iid case, the cost
function \(\Cost\) is a random function in general.\footnote{
	as a mnemonic device recall that we denote random functions and now
	also random measures in bold.
} The optimization of random
	functions is precisely
the subject of Bayesian optimization but leaves open the important question: Which distributional
assumptions are reasonable for \(\Cost\)? I.e. what prior should be chosen?
 \section{Choice of cost function}

The most common cost functions in machine learning are the crossentropy for
categorical data or generative tasks and the mean squared error (MSE) for
regression tasks. In this section we recapitulate that, not only is the MSE
effectively a crossentropy loss for a Gaussian model, but both cost functions
can be interpreted as maximum likelihood estimators of \(\rP_Z\)
\autocite[e.g.][Section 5.6]{goodfellowDeepLearning2016}. A maximum likelihood
estimation (MLE) is in turn effectively a maximum a posteriori estimation (MAP) with a
uninformative (uniform) prior on \(\rP_Z\). In particular this contradicts dirac
priors on \(\rP_Z\), or in other words deterministic \(\rP_Z\), which
corresponds to iid data. A consistent Bayesian model \emph{must} therefore
generalize the iid assumption on the data to exchangeablility.

By \textcite[Theorem~1.31]{schervishTheoryStatistics1995} the posterior
distribution is given by
\begin{equation}
	\label{eq: posterior}
	\overbrace{\Pr(\rP_Z \in d P_Z \mid Z_1,\dots, Z_n)}^{\text{posterior}}
	\;\propto\; \overbrace{\Pr(\rP_Z \in dP_Z)}^{\text{prior}}
	\overbrace{\prod_{k=1}^n \frac{dP_Z}{d\nu}(Z_k)}^{\text{likelihood}},
\end{equation}
where \(\nu\) is a reference measure
such that \(\Pr(Z_k \in dz \mid \rP_Z) = \rP_Z(dz)\) is almost surely continuous
w.r.t. \(\nu\) and we used that the \(Z_k\) are iid conditional on \(\rP_Z\)
with distribution \(\rP_Z\) to factorize the likelihood. Note that the choice of reference measure \(\nu\)
only changes the multiplicative constant of the likelihood. This implies that
this choice is irrelevant for the MLE. But since it is not possible to maximize
non-discrete distributions without a reference measure, we cannot directly maximize
the posterior to obtain the MAP. And if we
construct densities of the posterior to obtain an MAP, the
choice of reference measure does affect the outcome.

Bayesians would of course argue against reducing the posterior to a single
estimator. And, if you have to choose, Bayesian decision theory asks for a
loss function measuring the loss \(\loss(\hat{P}_Z, P_Z)\) of selecting
\(\hat{P}_Z\) when the actual distribution is given by \(P_Z\). The resulting
minimization problem is then given by 
\[
	\min_{\hat{P}_Z} \int \loss(\hat{P}_Z, P_Z) \; \Pr(\rP_Z \in dP_Z \mid Z_1,\dots, Z_n).
\]
To get back to the intuition of a single most likely \(P_Z\), one can consider the
0-1 loss function \(\loss(\hat{P}_Z, P_Z) = \ind_{d(\hat{P}_Z, P_Z)> \epsilon}\) and let
\(\epsilon\to 0\) forcing the choice of the single most likely \(P_Z\). Since \(\epsilon\)-balls
are intimately related to the Hausdorff (Lebesgue) measure, it is perhaps not surprising that
the limiting choice can often be shown to be the MAP with the Hausdorff measure
as the reference measure \autocite{bassettMaximumPosterioriEstimators2019}.

While the Hausdorff measure does not always exist on the space of \(\rP_Z\), in particular
if \(Z\) has continuous components \autocite[Sec.  5.3, Thm.
4 and following]{gelfandGeneralizedFunctionsVolume1964}, uniform priors also often do not exist as
measures and are taken as improper priors.\footnote{A solution to this problem
may be to replace the `countable additivity' requirement on measures with `finite
additivity' \autocite[p.\@ 21]{schervishTheoryStatistics1995}.}

Regardless, assuming a suitable (Hausdorff) measure exists to convert the distribution
of \(\rP_Z\) in \eqref{eq: posterior} into a density \(p\), we can maximize the
posterior to obtain the MAP. Since the logarithm and the multiplication with
\(\frac1n\) are strictly monotone transformations, the MAP is equivalently given
by
\[
	\text{MAP} = \argmin_{P_Z} -\underbrace{\frac1n\log(p(P_Z))}_{\text{prior}} +  \underbrace{\frac1n\sum_{k=1}^n -\log\Bigl(\frac{dP_Z}{d\nu}(Z_k)\Bigr)}_{=H^\nu(\rP_Z^n, P_Z)},
\]
where \(\rP_Z^n = \frac1n \sum_{k=1}^n \delta_{Z_k}\) is the empirical measure
and the \(\nu\)-crossentropy is given by
\[
	H^\nu(\mu_1, \mu_2) := \int -\log\Bigl(\frac{d\mu_2}{d\nu}(z)\Bigr) \mu_1(dz).
\]
Due to \(H^{\tilde{\nu}}(\mu_1,\mu_2) = H^\nu(\mu_1, \mu_2) + H^{\tilde{\nu}}(\mu_1, \nu)\)
the minimization of the crossentropy in \(\mu_2\) is independent of the
reference measure \(\nu\). Further observe that the influence of the prior decays as more
data is collected (\(n\to \infty\)). But note that the prior must have support on the space \(P_Z\)
is selected from. This rules out dirac priors for \(\rP_Z\) and therefore the iid
setting. Uninformative (uniform) priors imply that the MAP is equivalent to the MLE and
the MLE is always equivalent to the minimization of the crossentropy.

The crossentropy derived above does not yet match the crossentropy
used in supervised machine learning. The reason for this is that in
supervised learning we only wish to learn the conditional distribution
\(\rP_{Y\mid X}\) with 
\[
	\rP_Z(dz) = \rP_{Y\mid X}(dy \mid x) \rP_X(dx),
\]
where \(\rP_X(dx) = \rP_Z(dx \times \range) = \Pr(X_k\in dx \mid \rP_Z)\) is the
marginal distribution of \(X_k\) conditional on \(\rP_Z\). So we only model
the conditional distribution \(P_{Y\mid X}^\model\) parametrized by a model \(\model\)
and turn this into a distribution over \(Z\) by selecting an arbitrary
distribution \(\nu_X\) for \(X\). This can be something reasonable like
\(\rP_X\) or its empirical version \(\rP_X^n\), but we will see that this
selection has no effect on optimization. Assuming \(P^\model_{Y\mid X}\) is
absolutely continuous w.r.t \(\nu_Y\) for every \(x\), we can plug the full
distribution
\[
	P^{\model}_Z(dz) = P_{Y\mid X}^\model(dy \mid x) \nu_X(dx)
\]
into the \(\nu\)-crossentropy with reference measure \(\nu = \nu_Y \times \nu_X\) to obtain
\begin{align*}
	H^\nu(\rP^n_Z, P^\model_Z)
	&= \frac1n \sum_{k=1}^n - \log\Bigl(\frac{P^\model_{Y\mid X}}{d\nu_Y}(Y_k \mid X_k)\Bigr)
	+ \underbrace{\frac1n \sum_{k=1}^n -\log\Bigl(\underbrace{\frac{d\nu_X}{d\nu_X}(X_k)}_{=1}\Bigr)}_{= 0}.
\end{align*}

\subsection{Generative models and classification}

Generative models are supposed to sample (e.g. images) from a distribution
\(\rP_{Y\mid X}(dy \mid x)\) over the domain \(\range\) (of images) that depends
on a prompt \(x\). This is structurally similar to classification where labels \(Y\) must
be predicted based on input data \(x\). But in classification tasks the set of possible labels
\(\range\) is generally a finite set of categories. The model for this conditional
distribution is in both cases typically set up like this: A model \(\model\)
maps inputs \(x\) from \(\domain\) to their `energy' or `logits' in
\(\real^\range\). These are then turned into a conditional distribution via a
`softmax layer' with inverse temperature \(\beta\)
\[
	p^\model_{Y\mid X}(y\mid x)
	= \frac{
		\exp\bigl(\beta\,\model(x)[y]\bigr)
	}{
		\sum_{j\in \range}\exp\bigl(\beta\, \model(x)[j]\bigr)
	},
\]
where \(p^\model_{Y\mid X}\) is the density of \(P^\model_{Y\mid X}\) w.r.t. the
counting measure over the set of possible labels. For generative models with
infinite sample space \(\range\), the sum in the denominator (known as the `partition
function') turns into a general integral and becomes intractable. Methods to
treat this case are known as `energy based models' \autocite[see
e.g.][]{lecunTutorialEnergybasedLearning2007,songHowTrainYour2021}.
Since the model typically includes parameters to determine
the scale, \(\beta\) might seem superfluous. However, the temperature \(\beta\)
is often used after training to tune `creativity'.\footnote{With \(\beta\to 0\) the
softmax distribution approaches to the uniform distribution and thus becomes more `creative'.
Whereas with \(\beta\to\infty\) it concentrates on the most probable classes and
approaches the `argmax', increasing its safety/reproducibility.}
During training however, we can assume without loss of generality \(\beta=1\).

Using the counting measure for \(\nu_Y\), the sample cost induced by the crossentropy is 
then given by
\begin{align*}
	\Cost_n(\model)
	&= H^\nu(\rP_Z^n, P^\model_Z)
	\\
	&= -\frac1n \sum_{k=1}^n \model(X_k)[Y_k]
	+ \frac1n \sum_{k=1}^n \log\Bigl(\sum_{y\in \range} \exp\bigl(\model(X_k)[y]\bigr)\Bigr).
\end{align*}
The long term cost is thus given by
\begin{align*}
	\Cost(\model)
	&= \lim_{n\to\infty}\Cost_n(\model)
	= H^\nu(\rP_Z, P^\model_Z)
	\\
	&= \E\Bigl[-\model(X)[Y]\bigm| \rP_Z\Bigr]
	+ \E\Bigl[\log\Bigl(\sum_{y\in \range} \exp\bigl(\model(X)[y]\bigr)\Bigr) \bigm| \rP_Z \Bigr],
\end{align*}
where we dropped the indices of \((X, Y)\) since they are iid conditional on \(\rP_Z\).

\subsection{Regression}

In regression tasks, the space of labels \(\range\) is assumed to be a (Banach) vector space
and for our purposes given by \(\range = \real^\dims\). While an expectation
(averaging) would have been ill defined in the categorical
case\footnote{(Banach) vector spaces are the most general target space for
integrals to be well defined (cf. Bochner or Pettis integral).}, this setting
allows for the definition of a \emph{target function} representing the best prediction
of the label
given the input \(x\) for a given relation \(\rP_Z\)
\[
	\rf(x) := \E[Y \mid X=x, \rP_Z].
\]
The model \(\model\) attempts to approximate this target function \(\rf\) and it is
reasonable to model the conditional distribution with Gaussian noise around
\(\model\), i.e. 
\[
	P^{\model, \sigma^2}_{Y\mid X}(dy \mid x)
	= \normal\bigl(\model(x), \sigma^2\identity\bigr)(dy).
\]
We emphasize that this merely implies fitting a Gaussian model to the data
without making a Gaussian assumptions about the true distribution.\footnote{
	Of course the quality of the fitted Gaussian model is best when the
	noise \(\noise\) on the target function \(\rf\) is Gaussian. Indeed,
	considering linear models \(\model\),
	\textcite{fengOptimalConvex$M$estimation2024} find that the performance of the
	MSE minimizer deteriorates with non-Gaussian noise. They further device a
	method to automatically adapt the loss function (and thus the implicit model)
	to the noise distribution.
} 
With the Lebesgue measure as the reference measure \(\nu_Y\) it is straightforward
to verify, that minimization of the empirical or limiting MSE 
\[
	\Cost_n(\model) = \frac1n \sum_{k=1}^n \|\model(X_k) - Y_k\|^2
	\quad\overset{n\to \infty}{\longrightarrow}\quad \Cost(\model) = \E\bigl[\|\model(X) - Y\|^2\mid \rP_Z\bigr],
\]
over the model \(\model\) is equivalent to minimizing the crossentropy between
the empirical distribution \(\rP_Z^n = \frac1n \sum_{k=1}^n \delta_{Z_k}\) and
the Gaussian model \(P^{\model, \sigma^2}_{Y\mid X}\) for the empirical MSE
\[
	H^\nu(\rP_Z^n, P^{\model,\sigma^2}_Z)
	= \frac12\log(2\pi\sigma^2)
	+ \frac1n\sum_{k=1}^n \frac{\|\model(X_k) - Y_k\|^2}{2\sigma^2},
\]
and equivalent to minimizing \(H^\nu(\rP_Z, P^{\model,\sigma^2}_Z)\) resp. for the
limiting MSE. This is possible since in both cases the choice of
\(\sigma^2\) has no influence on the optimization over \(\model\). Thus one can
optimize over the model \(\model\) first to obtain \(\model^*\) and then
analytically determine that \(\sigma^2 = \Cost_n(\model^*)\) (or \(\sigma^2 =
\Cost(\model^*)\) resp.) is the minimizing choice for the variance \(\sigma^2\)
in the distributional model. The MSE can therefore be interpreted as the
variance of the fitted Gaussian model.

 \section{Noise distribution}
\label{sec: stochastic loss}

Observe that the stochastic losses \(\loss(\model(X), Y)\) are
exactly of the form as the noisy observations in Definition~\ref{def:
feasible optimization algorithm}. Specifically, the noise is given by
\[
    \noise_k(\model) := \loss(\model(X_k), Y_k) - \Cost(\model)
\]
such that the observations are given by
\[
    \loss_k(\model) := \loss(\model(X_k), Y_k) = \Cost(\model) + \noise_k(\model).
\]
Clearly, the \(\noise_k(\model)\) are centered and identically distributed
conditional on \(\rP_Z\).
\clearmargin
 \section{Parametric statistical model}
\label{sec: parametric statistical model}

We have motivated the random measure \(\rP_Z\) from the axiomatic assumption
of exchangable random variables and de Finetti's theorem
\autocite[e.g.][Thm.~1.49]{schervishTheoryStatistics1995}, which essentially
states that the exchangable case corresponds to the two stage model where we
first sample \(\rP_Z\) and then iid data from \(\rP_Z\).

A different, equivalent way to set this up is to start with a parametric family
of probability distributions, a `statistical model'
\[
    \mathcal P = \{P_\problem : \problem \in \problemSpace\}.
\]
Then one can choose a random variable \(\Problem\) in \(\problemSpace\)
such that \(P_\Problem\) becomes a random measure. To ensure no measurability
problems arise here we may assume the statistical model is of the form \(P_\problem(B) = \kernel(\problem; B)\)
for some probability kernel \(\kernel\). 
One may then select \((Z_k)_{k\in \nat}\) iid conditional on \(\Problem\) such that
\[
    \Pr(Z_k \in B \mid \Problem) =\kernel(\Problem; B) =  P_\Problem(B).
\]
But then the \(Z_k\) are clearly exchangeable and \(P_\Problem = \rP_Z\) and
the non-parametric version follows from the parametric setup.

Similarly, we can turn the non-parametric setting into a parametric setting.
Simply let \(\problemSpace\) be the space of measures, define \(P_\problem :=
\problem\) and let \(\Problem= \rP_Z\).

Observe that the parametric model circumvented the measure theoretical problem
of proving the existence of a kernel \(\kernel\) with \(\kernel(\rP_Z;B) =
\rP_Z(B)\) by taking it as definition. And similarly it is a bit easier to prove
statements that hold for all \(\problem \in \problemSpace\) by proving things
about \(P_\problem\). For this reason this notation is used in Chapter
\ref{chap: optimization landscape of shallow neural networks}.

The disadvantage of this parametric model is the additional notation that is
introduced without a clear interpretation of what this variable \(\problem\)
actually represents.
 \section{Random linear model}
\label{sec: random linear model}

Using a random linear model as an example we will show in this section, that the
cost function typically cannot be assumed to be stationary.
Isotropy however remains a reasonable assumption for cost functions.

We assume that \(X\) is sampled according to some
probability distribution \(\Pr_X\) and the labels \(Y=\rf(X)\) are the output of
a random linear target function \(\rf=f_\Problem\). That is we define
the parametric statistical model
\[
    P_\problem(dx, dy)
    = \Pr_X(dx)\, \delta_{f_\problem(x)}(dy)
\]
and sample \(X_k,Y_k\) iid from \(P_\Problem\), where \(\Problem=(\Tparam, \bb)\)
defines the random linear target function
\[
    \rf(x) = f_\Problem(x) = \Tparam^T x + \bb.
\]
Using a linear model for regression
\[
    \model_\param(x) = \param^T x + \bias,
\]
where \(\param=(\param, \bias)\) with some abuse of notation,
the mean squared error cost function is then of the form
\begin{align*}
    \Cost(\param)
    &= \E[(\model_\param(x) - \rf(x))^2 \mid \Problem]
    \\
    &= \E\Bigl[\bigl( (\param - \Tparam)^T X + (\bias - \bb)\bigr)^2\mid \Problem\Bigr]
    \\
    &= (\param - \Tparam)^T \underbrace{\E[XX^T]}_{=\identity} (\param - \Tparam)
    + (\param - \Tparam)^T \underbrace{\E[X]}_{=0}(\bias - \bb)
    + (\bias - \bb)^2
    \\
    &= \|\param - \Tparam\|^2 + (\bias - \bb)^2
    \\
    &= \|(\param, \bias) - (\Tparam, \bb)\|^2,
\end{align*}
\clearmargin
assuming `whitened'\footnote{
    i.e.\ centered, decorrelated and standardized input \(X\) such that \(\E[X]=0\) and
    \(\E[XX^T] = \identity\).
} input in the third equation.
Via some more abuse of notation, i.e. \(\Tparam := (\Tparam, \bb)\) and
\(\param=(\param, \bias)\) we disregard the bias without loss of generality.
Therefore we have
\[
    \Cost(\param) = \|\param - \Tparam\|^2
    = \underbrace{\overbrace{
        \|\param\|^2 + \E\|\Tparam\|^2
    }^{\text{deterministic}}
    - 2\langle \param, \Tparam\rangle}_{=:\tilde{\Cost}(\param)}
    + \underbrace{
        \|\Tparam\|^2-\E\|\Tparam\|^2
    }_{
        \text{random const.}
    }
\]
For the purpose of optimization, constants are irrelevant and we may therefore
pretend \(\|\Tparam\|^2-\E\|\Tparam\|^2=0\). Observe that by a change
of the origin to \(\E[\Tparam]\) we may assume \(\E[\Tparam] = 0\), since
\(\Cost(\param + \E[\Tparam]) = \|\param - (\Tparam-\E[\Tparam])\|\).
The random function \(\tilde\Cost\) is therefore centered with covariance
function
\[
    \C_{\tilde\Cost}(\param_1, \param_2)
    = 4\E[\langle \param_1, \Tparam\rangle\langle \Tparam, \param_2\rangle]
    = 4 \param_1^T \E[\Tparam \Tparam^T] \param_2.
\]
And for \(\Tparam\sim\normal(0, \sigma^2 \identity)\) we further have
that \(\tilde\Cost\) is a Gaussian random function with
\begin{align}
    \tag{mean}
    \E[\tilde\Cost(\param)] &= c + \|\param\|^2
    \\
    \tag{covariance function}
    \C_{\tilde\Cost}(\param_1, \param_2) &= 4\sigma^2 \langle \param_1, \param_2\rangle.
\end{align}
where \(c=\E\|\Tparam\|^2\). Observe that the random function \(\tilde\Cost\) is
clearly non-stationary but still isotropic.

\begin{prop}[Distribution of the noise]
    In the bias free linear model with \(\Tparam\sim\normal(0,\sigma^2
    \identity_{\dims\times\dims})\) for whitened and Gaussian input
    \(X\sim\normal(0, \identity_{\dims\times\dims})\) the
    noise
    \[
        \noise(\param)
        := \loss(\model_\param(X),Y) - \Cost(\param)
    \]
    is centered, i.e.\ \(\E[\noise(\param)\mid \Tparam]=0\), with covariance
    \begin{align*}
        \E[\noise(\param) \noise(\tilde \param) \mid \Tparam] 
        &= 2 \langle \param - \Tparam, \tilde\param - \Tparam\rangle^2
        \\
        \E[\noise(\param) \noise(\tilde\param)]
        &= 2\bigl(\langle \param, \tilde \param\rangle + \sigma^2 \dims\bigr)^2.
    \end{align*}
\end{prop}
\begin{proof}
The noise is given by
\begin{align*}
    \noise(\param)
    &= \loss(\model_\param(X),Y) - \Cost(\param)
    \\
    &= (\param^T X - \Tparam^T X)^2 - \|\param - \Tparam\|^2
    \\
    &= (\param- \Tparam)(XX^T - \identity)(\param - \Tparam).
\end{align*}
Since the noise is always centered by definition, we skip right to the calculation
of the covariance
\begin{align*}
    \E[\noise(\param)\noise(\tilde{\param})\mid \Tparam]
    &= \sum_{i,j=1}^\dims \sum_{l,k=1}^\dims \begin{aligned}[t]
        \E\Bigl[
        &(\param-\Tparam)_i (X_i X_j - \delta_{ij})(\param-\Tparam)_j
        \\
        &\qquad \cdot (\tilde{\param}-\Tparam)_k (X_k X_l - \delta_{kl})(\tilde{\param}-\Tparam)_l
    \mid \Tparam \Bigr]
    \end{aligned}
    \\
    &= \sum_{i,j,l,k=1}^\dims \begin{aligned}[t]
        &(\param-\Tparam)_i (\param-\Tparam)_j
        (\tilde{\param}-\Tparam)_k (\tilde{\param}-\Tparam)_l
        \\
        &\qquad \cdot \E\bigl[(X_i X_j - \delta_{ij})(X_k X_l - \delta_{kl})\mid \Tparam\bigr]
    \end{aligned}
\end{align*}
Since \(X\sim\normal(0, \identity_{\dims\times\dims})\), independent of
\(\Tparam\), a case by case analysis results in
\begin{align*}
    &\E\bigl[(X_i X_j - \delta_{ij})(X_k X_l - \delta_{kl})\mid \Tparam\bigr]
    = \E[(X_i X_j - \delta_{ij})(X_k X_l - \delta_{kl})]
    \\
    &= \left\{\begin{aligned}
        &\E[(X_i^2-1)^2] &&= \Var(X_i^2) &= 2, \quad && i=j=l=k\\
        &\E[X_i^2 X_j^2] &&= \E[X_i^2]\E[X_j^2] &= 1, \quad && [i=k\neq j=l] \text{ or } [i=l\neq j=k]\\
        &0, && & && \text{ else }
    \end{aligned}
    \right.
    \\
    &= \delta_{ik}\delta_{jl} + \delta_{il}\delta_{jk}.
\end{align*}
With \(\Tparam \sim \normal(0, \sigma^2 \identity)\) this implies
\begin{align*}
    &\E[\noise(\param) \noise(\tilde{\param})\mid \Tparam]
    \\
    &=\sum_{i,j,l,k=1}^\dims
        (\param-\Tparam)_i (\param-\Tparam)_j
        (\tilde{\param}-\Tparam)_k (\tilde{\param}-\Tparam)_l
        [\delta_{ik}\delta_{jl} + \delta_{il}\delta_{jk}]
    \\
    &= 2 \sum_{i,j=1}^\dims
        (\param-\Tparam)_i (\param-\Tparam)_j
        (\tilde{\param}-\Tparam)_i (\tilde{\param}-\Tparam)_j
    \\
    &= 2\langle \param - \Tparam, \tilde\param - \Tparam\rangle^2.
\end{align*}
And the unconditional expectation is given by
\begin{align*}
    &\E[\noise(\param) \noise(\tilde{\param})]
    \\
    &= 2 \sum_{i,j=1}^\dims
        \E[(\param-\Tparam)_i (\param-\Tparam)_j
        (\tilde{\param}-\Tparam)_i (\tilde{\param}-\Tparam)_j]
    \\
    &= 2 \sum_{i,j=1}^\dims  \E\Bigl[
        \Bigl(\param_i \tilde{\param}_i - (\param_i+\tilde{\param}_i)\Tparam_i +\Tparam_i^2\Bigr)
        \Bigl(\param_j \tilde{\param}_j - (\param_j+\tilde{\param}_j)\Tparam_j +\Tparam_j^2\Bigr)
    \Bigr]
    \\
    &= 2 \sum_{i,j=1}^\dims  \E\Bigl[
        \Bigl(\param_i \tilde{\param}_i - (\param_i+\tilde{\param}_i)\Tparam_i +\Tparam_i^2\Bigr)
        \underbrace{\E\Bigl[\param_j \tilde{\param}_j - (\param_j+\tilde{\param}_j)\Tparam_j +\Tparam_j^2\mid \Tparam_i \Bigr]}_{=\param_j\tilde{\param}_j + \sigma^2}
    \Bigr]
    \\
    &= 2 \sum_{i,j=1}^\dims  (\param_i\tilde{\param}_i + \sigma^2)(\param_j\tilde{\param}_j + \sigma^2)
    \\
    &= 2 \left(\sum_{i=1}^\dims  (\param_i\tilde{\param}_i + \sigma^2)\right)^2
    = 2 (\langle \param, \tilde{\param}\rangle + \sigma^2\dims)^2.
    \qedhere
\end{align*}
\end{proof}

\printbibliography[heading=subbibliography]
 	\end{refsection}

	\begin{refsection}
		{
			\renewcommand*{\param}{\theta}
			\chapter{Optimization landscape of shallow neural networks}
\label{chap: optimization landscape of shallow neural networks}

In this chapter we examine the optimization landscape of a broad class of regression 
problems with squared error loss for shallow ANNs using analytic activation 
functions. The main finding is that, in an appropriate sense, almost all such 
problems exhibit a “nice” optimization landscape. More precisely, we show that the cost 
landscape is typically Morse on the domain of \emph{non-degenerate} parameters. 
Conversely, the set of \emph{degenerate} parameters -- those for which the same response function can be accomplished with fewer neurons -- has 
significantly smaller Hausdorff dimension. Such parameters do not fully exploit 
the network's representational capacity, and inherently have redundancies that prevent 
the optimization landscape from being Morse at those points.

In order to state our results we start with a formal definition of shallow neural networks. The definition uses a graph structure that will prove to be useful later.

\begin{definition}[Shallow neural network]\label{def:ANN}
	\label{def: shallow neural network}
    A \emph{(dense) shallow neural network} (briefly called \emph{ANN}) is a tuple \(\network = (\mathbb V, \activation)\)  consisting of
\begin{itemize}
\item a tuple $\mathbb V=(V_0,V_1,V_2)$ of finite disjoint sets $V_0$, $V_1$ and
$V_2$ (the neurons of the input \(V_\tin := V_0\), hidden \(V_1\) and output
layer \(V_\tout := V_2\)) and
\item a measurable function \(\activation:\real\to \real\) (the activation
function).
\end{itemize}
\end{definition}

\begin{definition}Let \(\network=(\mathbb V, \psi)\) be an ANN.
\begin{enumerate}\item We call the directed graph $G=(V,E)$  given by 
\[ V=V_0\cup V_1\cup V_2\text{ \ and \ } E= \bigcup_{k=0}^1 V_k\times V_{k+1}\]
the \emph{ANN-graph} of \(\network\).
\item We call \(\Theta=\Theta_{\network}=\real^E\times \real^{V_1\cup V_2}\) \emph{parameter space of the network \(\network\)} and  
every tuple  \(\param=(\weight, \bias)\in \Theta=\real^E\times \real^{V_1\cup V_2}\) a \emph{parameter of the network~\(\network\)}. We refer to \(\weight\) as the 
 (edge) \emph{weights} and to \(\bias\) as the 
        \emph{biases}.  
        \item      
        For every parameter $\theta=(\weight,\bias)\in\Theta$ we call
        \begin{align}\label{eq:reali}
        \Psi_\param:\real^{V_{\tin}}\to \real^{V_{\tout}}  , \ x\mapsto  \Bigl(\bias_l
        + \sum_{j\in V_1} \activation\Bigl(
            \bias_j + \sum_{i\in V_{\tin}}  x_i\weight_{ij}\Bigr)\weight_{jl} 
        \Bigr)_{l\in V_{\tout}}
        \end{align}
the   \emph{response function} of the parameter $\param$.
\end{enumerate}\end{definition}
Typically, the underlying network is clear from the context and it is therefore omitted in the notation.

We study regression problems where the input data lies in \(\real^{V_\tin}\) and
the labels in \(\real^{V_\tout}\). These can be formally described by a
distribution \(\Pr_X\) on \(\real^{V_\tin}\) (the distribution of the input data)
and a probability kernel \(\kernel\) from \(\real^{V_\tin}\)  to \(\real^{V_\tout}\)
(the conditional distribution of the label given the input data). Our aim is to
show that for a fixed distribution \(\Pr_X\) for ``most'' kernels \(\kernel\) the
respective optimization landscape is Morse on the efficient domain. For this we analyze
random regression problems where the kernel in the regression problem itself is
random.

\begin{definition}\label{def:873}
    Let \(\network\) be an ANN.  A \emph{measurable family of regression
    problems} is a tuple \(\mathfrak R=(\Pr_X,\kernel,\loss)\) consisting of
    \begin{enumerate}
        \item a distribution \(\Pr_X\) on \(\real^{V_\tin}\) (the distribution of the input data \(X\))
        \item a measurable set \((\problemSpace,\mathcal M)\) (the statistical model space)
        \item a probability kernel $\kernel$ mapping model and input data from
        \(\mathbb M\times \real^{V_\tin}\) to a probability distribution over labels in
        \(\real^{V_\tout}\) and
                \item a measurable function \(\loss:\real^{V_\tout}\times\real^{V_\tout}\to [0,\infty]\) (the \emph{loss}).
\end{enumerate}
\end{definition}

When dealing with measurable families of regression problems we will always
associate the setting with a measurable space that is equipped with a family of
distributions (\(\Pr_\problem)_{\problem\in\problemSpace}\) together with a
\(\real^{V_{\tin}}\)-valued random variable \(X\) (the input data) and a
\(\real^{V_{\tout}}\)-valued random variable \(Y\) (the label) such that
under  every distribution \(\Pr_\problem\) with \(\problem\in\problemSpace\),
\(\Pr_X\)  is the distribution of \(X\) and \(\kernel(\problem,\cdot;\cdot)\) is the 
conditional distribution of \(Y\) given \(X\), i.e.,
\[
    \Pr_\problem(Y\in B \mid X) = \kernel(\problem, X; B), \text{ \ a.s.}
\]

Then the loss that we incur when using a shallow ANN for the prediction defines
a cost landscape in the sense of the following definition.

\begin{definition}[Cost function]\label{def:J_and_f}
    Let \(\network\) be an ANN as in Definition \ref{def:ANN}
and \(\mathfrak R\) a measurable family of regression problems as in Definition~\ref{def:873} and a function \(R:\Theta\to\real\) (regularization). The family of functions  $(J_{\mathbf m})_{\mathbf m\in\mathbb M}$ given by
    \[
        J_{\problem}\colon \Theta \to (-\infty,\infty],
        \  \param\mapsto
        \E_\problem\bigl[\loss(\response_\param(X), Y)\bigr]+ \regularizer(\param) 
\]
the \emph{(regularized) cost functions of \((\network, \mathfrak R,R)\)}. 
\end{definition}

A useful concept for the analysis of the MSE cost function is the `target
function' representing the best possible predictor.

\begin{definition}[Family of \(L^p\)-integrable regression problems, target function] Let \(p\in[1,\infty)\). 
A measurable family of
    regression problems  \(\mathfrak R\) is said to be \emph{\(L^p\)-integrable}, if for every   \(\problem\in \problemSpace\) the label   is \(L^p\)-integrable, i.e.,
    \[
        \forall \mathbf m\in\mathbb M:\; \E_\problem\bigl[\|Y\|^p\bigr]<\infty.
    \]
For a family of \(L^1\)-integrable
    regression problems \(\mathfrak R\), for every   \(\problem\in \problemSpace\) the function
    \[
        \target_\problem(x) := \int  y\,\kernel(\problem,x; dy) \overset{\Pr_X\text{-a.s.}}= \E_\problem[Y \mid X=x]
    \]
    is well-defined for \(\Pr_X\)-almost all \(x\in \real^{V_\tin}\).
    We call \(\target_\problem\) the \emph{target function of \(\problem\)}.
\end{definition}

The conditions assumed for our main result are collected in the following
definition.
\begin{definition}\label{def: standard model}
    The \emph{standard setting} is a tuple \((\network,  \mathfrak R, \regularizer,\Problem)\)
    consisting of 
    \begin{itemize}
    \item an ANN $\network$ with one dimensional output $\#V_\tout=1$ and
    analytic activation function \(\activation\),
    \item a  family
    of  \(L^2\)-integrable regression problems \(\mathfrak R\) with squared-error loss
    \(\loss(\hat{y}, y) = (\hat{y}-y)^2\) and compact support 
    \(\Domain:=\mathrm{supp}(\Pr_X)\) of the input distribution \(X\),
    \item  analytic convex (regularization) function \(\regularizer:\Theta\to \real\) and
    \item an $\problemSpace$-valued random variable $\Problem$ such that
    the random target function \(\rf:=\target_\Problem\) is an \emph{weakly universal Gaussian random function} in the sense that for every continuous test function \(\phi:\real^{V_\mathrm{in}}\to\real\) the random variable
    \[ 
    \langle \phi, f_\mathbf{M}\rangle_{\Pr_X}:=\int \phi(x) \,\target_\Problem(x) \,\Pr_X(dx)
    \]
    is Gaussian and has strictly positive variance whenever \(\phi\not\equiv 0\)
    on \(\Domain\).\footnote{
        Note that for every \(\problem\in\problemSpace\),
        \(\target_\problem\) is in \(L^2(\Pr_X)\subseteq L^1(\Pr_X)\) and \(\phi\) is uniformly
        bounded on the compact support of \(\mathbb P_X\) so that the integral
        \(\langle \phi, \target_\Problem\rangle_{\Pr_X}\) is for \emph{all} realizations
        of \(\mathbf M\) well-defined. It  is further  measurable since
        \(f_\cdot(\cdot)\) is product measurable by Fubini's theorem.
    }
    \end{itemize}
\end{definition}

To understand the cost landscape of models $\problem\in\problemSpace$  we need
to divide the space of all parameters $\param\in \paramSpace=\real^E\times \real^{V_1\cup
V_2}$ into two domains. For the activation functions \(\sigmoid\) and \(\tanh\)
we define
\begin{itemize}
\item the \emph{efficient domain} by
\begin{equation}
    \label{eq: defintion efficient 0} 
    \efficient_0
    := \Biggl\{
        \param =(\weight, \bias)\in \real^{E \times (V\setminus V_0)}:
        {
        \begin{aligned}[c]
            \scriptstyle \weight_{j\placeholder}
            &\scriptstyle\neq 0
            &&\scriptstyle \forall j \in V_1,
            \\[-0.5em]
            \scriptstyle\weight_{\placeholder j}
            &\scriptstyle\neq 0
            &&\scriptstyle\forall j\in V_1, 
            \\[-0.5em]
            \scriptstyle(\weight_{\placeholder i}, \bias_i)
            &\scriptstyle\neq \pm(\weight_{\placeholder j}, \bias_j)
            &&\scriptstyle\forall i,j\in V_1\text{ with }i\not=j
        \end{aligned}
        }
    \Biggr\},
\end{equation}
\item the \emph{redundant domain} by \((\real^E\times \real^{V_1\cup V_2})\backslash \efficient_0\).
\end{itemize}

Our main structural result, namely the fact that the cost is typically Morse, is
only true on the efficient domain. The restriction onto the efficient domain is
natural since any redundant parameter lies on a path of constant response
such that no local minimum in the redundant domain can be a strict local
minimum. In particular the cost cannot be Morse on the whole set of parameters.

Our main result states that for ``most'' statistical models $\problem$
the realization of the MSE is Morse on the efficient domain. 

\begin{theorem}[Almost all optimization landscapes are Morse on the efficient domain]\label{thm:main1}
	Let \((\network, \mathfrak R,R, \Problem)\) be a standard setting
    (Definition~\ref{def: standard model}). Assume \(\activation\in \{\sigmoid,
    \tanh\}\) about the activation function and that the support of
    \(\Pr_X\) contains a non-empty open set.
	Almost surely, the regularized cost \(\cost_\Problem\colon \efficient_0\to \real\) with
	\[  
        \cost_\Problem(\param) = \E_\Problem\bigl[\loss(\response_\param(X), Y)\bigr]+ \regularizer(\param) 
	\]
	is a Morse function. Equivalently, it holds that 
	\[
		\Pr\Bigl(
			\exists \param \in \efficient_0:
			\nabla J_{\mathbf M}(\param) = 0, \,\det(\nabla^2 J_{\mathbf M}(\param)) = 0
		\Bigr) = 0.
	\]
\end{theorem}

In Section~\ref{sec:Morse} we actually prove a version of this theorem for general
analytic activation functions (see Theorem~\ref{thm:main1.1}). This requires
a notion of the efficient domain that is an implicitly defined set. In
Section~\ref{sec: characterization of efficient params} we then prove that for the activation functions
\(\mathrm{sigmoid}\) and \(\mathrm{tanh}\) the implicitly defined version of the
efficient domain agrees with the one used in the latter theorem.

\begin{remark}[Generalization of the Gaussian assumption]
    While we assume that the target function is weakly universal Gaussian in the
    standard setting (Definition~\ref{def: standard model}), our main theorem
    (Theorem~\ref{thm:main1}) is a statement about null sets. Since null sets
    remain null sets for measures that are absolutely continuous with respect to
    such Gaussian measures and mixtures thereof, it is straight-forward to
    generalize the statement to significantly more general distributions of the
    random target function \(\target_\Problem\). That is, the statement remains
    true if the distribution of \(\target_\Problem\) can be written as a mixture
    of measures that are absolutely continuous with respect to weakly universal
    Gaussian measures!
\end{remark}
\begin{remark}[Weak universality]
    \label{rem: weak universality}
    For a better understanding of weak universality consider the stronger
    assumption\footnote{
        On the positive probability event in \eqref{eq: positive ball probability} we have
        \[
            \begin{aligned}[t]
            &|\langle \phi, \target_\Problem\rangle_{\Pr_X}-  \|\phi\|_{\Pr_X}^2|
            \\
            &= |\langle \phi, \target_\Problem- \phi \rangle_{\Pr_X}|
            \\
            &\le  \epsilon \|\phi\|_{\Pr_X}.
            \end{aligned}
        \]
        Since \(\phi\) is continuous and non-zero on the support of \(\Pr_X\)
        this implies \(\|\phi\|_{\Pr_X}>0\). Choosing
        \(\epsilon\in (0,\|\phi\|_{\Pr_X})\) we conclude that \(\langle \phi,
        f_\mathbf{M}\rangle_{\Pr_X}>0\) with strictly positive probability. The same
        argument applied to \(-\phi\) gives that \(\langle \phi,
        f_\mathbf{M}\rangle_{\Pr_X}<0\) with strictly positive probability.
        Consequently, \(\langle\phi, f_\mathbf{M}\rangle_{\Pr_X}\) has positive variance.
    } that all continuous functions
    \(\phi\colon \real^{V_\tin}\to\real^{V_\tout}\) lie in the support of
    \(\Pr_{\target_\Problem}\), when \(\target_\Problem\) is a random element
    in \(L^2(\Pr_X)\). I.e. for every continuous
    \(\phi\colon \real^{V_\tin}\to\real^{V_\tout}\) and \(\epsilon>0\), one
    has that 
    \begin{equation}
        \label{eq: positive ball probability}    
        \Pr\bigl(\|f_\mathbf{M}-\phi\|_{\Pr_X}<\epsilon\bigr)>0.
    \end{equation}
    This is a \emph{universality} assumption \autocite[Thm.\@
    3.6.1]{micchelliUniversalKernels2006, carmeliVectorValuedReproducing2010,
    bogachevGaussianMeasures1998} and intuitively means that no continuous
    function \(\phi\) can be ruled out as the target function
    \(\target_\Problem\) ex ante. We believe this is a natural assumption for
    a learning problem.

    In the proof of Theorem~\ref{thm:main1}, we will actually work with
    an even weaker assumption than weak universality: It would suffice to assume
    the non-degeneracy for real-analytic test functions \(\phi\) only.
\end{remark}

A natural question is whether local minima on the efficient domain exist and whether the restriction to the efficient domain in
Theorem \ref{thm:main1} is an artefact of our proof. This will be the content of
Sections \ref{sec: neighborhood of redundant parameters}-\ref{sec: existence of
redundant critical points}.  Intuitively we show, for the standard
unregularized setting with activation \(\activation\in \{\sigmoid, \tanh\}\) and
the additional regularity assumption \eqref{eq: positive ball probability} in
Remark \ref{rem: weak universality}, that we have the following:
\begin{itemize}
    \item For every open set \(U\subset \paramSpace\) containing an efficient
    point, the probability is strictly positive that the  loss has a local
    minimum in \(U\), see Theorem \ref{thm: existence of efficient minima}.
    
    \item With strictly positive probability, there exist critical points in the
    set of redundant parameters (Theorem \ref{thm: redundanct critical points exist}) and all redundant critical points have
    a direction of zero curvature (the determinant of the Hessian is zero), see Theorem \ref{thm: critical points of redundant type}.
\end{itemize}

It is therefore \emph{impossible} to prove the MSE to be a Morse function
on the redundant domain, since critical points may exist and those always
violate the Morse condition.

\paragraph*{Outline}

In Section \ref{sec:Morse} we prove a more general version of Theorem \ref{thm:main1}
which is applicable to all analytic activation functions. However, for
general analytic activation functions the efficient domain has to be defined in
an implicit way. Specifically, we will prove in Theorem~\ref{thm:main1.1} for
the standard setting that the cost landscape is almost surely Morse on the set
of \emph{polynomially efficient parameters} (Definition~\ref{def: polynomial
independence}). In Section~\ref{sec: characterization of efficient params} we
show that the various
definitions of efficient parameter domains coincide for \(\activation\in
\{\sigmoid, \tanh\}\) (Theorem~\ref{thm: characterization of efficient networks}).
With this result Theorem~\ref{thm:main1} becomes a direct corollary of
Theorem~\ref{thm:main1.1}. In Section \ref{sec: neighborhood of redundant parameters}
we prove for any redundant parameter \(\param\) that there exists a straight line 
of parameters \(\param(t)\), starting in \(\param\), where the response, and
therefore the cost, remains unchanged, i.e. \(\response_{\param(0)} =
\response_{\param(t)}\). In Section \ref{sec: existence of
efficient critical points} we prove that efficient local minima exist
with positive probability. We use this fact in Section \ref{sec: existence of
redundant critical points} to prove that redundant critical points exist with positive
probability. To this end we extend an efficient critical parameter of a smaller
network to a redundant critical parameter of a larger network.

\section{MSE is Morse on efficient domain}
\label{sec:Morse}

In this section, we will prove that for a standard model the random
loss-landscape is Morse on the \emph{polynomially}, efficient domain. For
general activation functions~\(\activation\) we have to work with a different
notion of the efficient
domain than \(\efficient_0\) introduced in \eqref{eq: defintion efficient 0}. As we will show in Section~\ref{sec: characterization of efficient
params} the definition coincides with \(\efficient_0\) whenever
\(\activation\in\{\sigmoid, \tanh\}\) and the support of \(\Pr_X\)
contains an open set.

\begin{definition}[Polynomial efficiency]
	\label{def: polynomial independence}
    Let \(\network = (\nodes, \activation)\) be an ANN (Definition~\ref{def: shallow neural network}), \(n\in\nat_0\) and \(m=(m_{\emptyset}, m_0, \dots, m_n)\in\nat_0^{n+2}\).
	\begin{enumerate}[label=(\roman*)]
        \item
        A parameter \(\param\in\Theta\) is called \emph{\(m\)-polynomially
        independent on \(\Domain\)} if \(\activation\) is
        \(n\)-times differentiable and the equation
        \[
            0 = P^{(\emptyset)}(x)
            + \sum_{j\in V_1}
            \sum_{k=0}^n P^{(k)}_j(x) \activation^{(k)}\Bigl(\bias_j + \sum_{i\in V_0}x_i\weight_{ij} \Bigr)
            \quad \forall x\in \Domain
        \]
        considered in all polynomials \(P^{(\emptyset)}\) and \((P^{(k)}_j:j\in V_1,
        k\in\{0,\dots,n\})\) of at most degree \(m_\emptyset\) and  \(m_k\),
        respectively, has only the trivial solution where all polynomials are
        identically zero.
        Here, \(\activation^{(k)}\) denotes the \(k\)-th derivative of the activation function \(\activation\).
        \item
        A parameter \(\param\in\Theta\) is called \emph{\(m\)-polynomially efficient on \(\Domain\)}, if 
        \begin{enumerate}[label={(\alph*)}]
            \item
            all neurons are used meaning that for all  \(k\in V_1\) one has
            \[
                \weight_{k \placeholder} = (\weight_{kl})_{l\in V_\tout} \not\equiv
                0,
            \]
            and
            \item it is \(m\)-polynomially independent.
        \end{enumerate}
	We denote by \(\polyEfficient^m=\polyEfficient^m(\Domain)\) the set of all
	\(m\)-polynomially efficient parameters.
    \end{enumerate}
\end{definition}

\begin{theorem}[MSE is a Morse function on polynomially efficient parameters]
	\label{thm: morse function on efficient networks}\label{thm:main1.1}
	Let \((\network, \mathfrak R, R, \Problem)\) be the standard setting
    (Definition~\ref{def: standard model}).
	Then the MSE cost is almost surely a Morse function on the set
	\(\polyEfficient:=\polyEfficient^{(0,0,1,2)}(\Domain)\) of \((0,0,1,2)\)-polynomially
	efficient parameters on the support \(\Domain\) of \(\Pr_X\), i.e.
	\[
		\Pr\Bigl(
			\exists \param \in \polyEfficient:
			\nabla\Cost(\param) = 0, \det(\nabla^2\Cost(\param)) = 0
		\Bigr) = 0
	\]
\end{theorem}

Before we explain the methodology of our proof we first derive a crucial
representation for the MSE cost \(\cost_\problem\) given by 
\[
    \cost_\problem(\param)
    = \E_\problem\bigl[ \|\response_\param(X) - Y\|^2\bigr] +\regularizer(\param) 
\]
with convex regularizer \(R\). Recall that \(\Psi_\param\) is the realization function of the ANN as introduced in~(\ref{eq:reali}).

\begin{prop}[Decomposition of the MSE cost]
    \label{prop: decomposition}
    For \(\problem\in \problemSpace\) and \(\param\in\paramSpace\)
    one has
    \begin{equation}\label{eq: decomposition}
        \cost_\problem(\param)
        = \regularizer(\param) + \|\response_\param\|_{\Pr_X}^2
        - 2\underbrace{
            \langle \response_\param, \target_\problem\rangle_{\Pr_{X}}
        }_{=:\hat{\cost}_\problem(\param)}
        + \E_\problem[\|Y\|^2].
    \end{equation}
    Here \(\|\cdot\|_{\Pr_X}\) is induced by \(\langle \phi, \varphi\rangle_{\Pr_X} := \int \langle\phi(x), \varphi(x)\rangle \Pr_X(dx)\).
\end{prop}
\begin{proof}
    With the Pythagorean formula we have
    \begin{align*}
        \cost_\problem(\param)
        &= \regularizer(\param) + \E_\problem[\|\response_\param(X) - Y\|^2]
        \\
        &= \regularizer(\param) + \E_\problem[\|\response_\param(X)\|^2] - 2\E_\problem[\langle \response_\param(X), Y\rangle] + \E_\problem[\|Y\|^2].
    \end{align*}
    Since \(\target_\problem(X)=\E_\problem[Y|X]\) we conclude with the tower property
    \[
        \E_\problem[\langle \response_\param(X), Y\rangle] 
        = \E_\problem[\langle \response_\param(X), \target_\problem(X)\rangle]
        \overset{\text{def.}}= \langle \response_\param, \target_\problem \rangle_{\Pr_X}
        =\hat{\cost}_\problem(\param).
        \qedhere
    \]
\end{proof}

In view of  \eqref{eq: decomposition} we observe that the term
\(\E_\problem[\|Y\|^2]\) does not depend on the parameter \(\param\) and it is
thus irrelevant when it comes to deciding whether \(\cost_\problem\) is Morse or
not. In the proof we then argue that the event where the stochastic process
\((R(\param)+\|\response_\param\|_{\rP_X}^2-2\hat{\cost}_\Problem(\param))_{\param\in\Theta}\) is not Morse on the
polynomially efficient parameters is a \emph{``thin set''}.

Unfortunately, our setting is not immediately covered by the arguments of
\citet{adlerRandomFieldsGeometry2007}. Roughly speaking, their approach is as
follows. If there is a parameter \(\param\) that is critical with its Hessian
having a zero eigenvalue, then this satisfies
\begin{align}\label{eq: morse violation}
    \nabla \Cost(\param)=0 \quad\text{and}\quad
    \det(\nabla^2  \Cost(\param))=0.
\end{align}
Note that the latter is a collection of \(\dim(\paramSpace)+1\) real equations in
\(\dim(\paramSpace)\) real variables and intuitively one would expect that, under
appropriate non-degeneracy assumptions, the equation does not have solutions.
The equations in~\eqref{eq: morse violation} depend on the collection of first
order differentials \(\rg_1(\param)\) and of second order differentials
\(\rg_2(\param)\). As shown in Lemma~11.2.10 of
\citet{adlerRandomFieldsGeometry2007} solutions of~\eqref{eq: morse violation} would
not exist, if for every \(\param\) under consideration (in our
case the efficient domain) the combined vector \((\rg_1(\param),\rg_2(\param))\)
has locally uniformly bounded Lebesgue density. 

Unfortunately,  in our situation, many second order differentials are degenerate and the
result is not applicable. To bypass this problem we proceed as follows.  In the
following  subsection, we will first analyze the stochastic process 
\[
    \hat{\Cost}
    =(\hat{\cost}_\Problem(\param))_{\param\in \paramSpace}
    =(\langle \response_\param,\target_\Problem\rangle_{\Pr_X})_{\param \in \paramSpace} .
\]
This process is obtained by applying a \(\param\)-dependent linear functional on
the random target function \(\rf=\target_\Problem\) and thus \(\hat{\Cost}\) is
a Gaussian process since \(\rf\) is Gaussian by assumption. We will
collect in \(\rg_1(\param)\) all first order differentials and in
\(\rg_2(\param)\) the \emph{`centered'}\footnote{
    \label{footnote: centered}
    This is only true if \(\target_\Problem\) is centered, which we are not
    willing to assume. But this provides the right intuition, since we subtract a
    deterministic term (which is not necessarily the mean).
} and \emph{non-degenerate} second order differentials. The
non-degeneracy of the combined collection \(\rg= (\rg_1,\rg_2)\) is shown in
Proposition~\ref{prop: non-degenrate g} and follows from the polynomial
independence that is assumed in the polynomially efficient domain (cf.~Definition~\ref{def:
polynomial independence}).

In Proposition~\ref{prop: graph of g is not in U}
we show that the generalization of the volume argument of
\citet[Lemma~11.2.10]{adlerRandomFieldsGeometry2007}
given in Lemma~\ref{lem: generalization 11.2.10} is
applicable to \(\rg\). The generalization of the volume argument
is necessary since we want to show that the process
\[
    \bigl(\regularizer(\param) + \|\response_\param\|_{\Pr_X}^2-2 \hat{\Cost}(\param)\bigr)_{\param\in\paramSpace}
\]
never satisfies \eqref{eq: morse violation} on \(\polyEfficient\).
While \citet{adlerRandomFieldsGeometry2007} considered level sets, we move
the model indepeendent term \(\regularizer(\param) + \|\response_\param\|_{\Pr_X}^2\) to the other side in \eqref{eq:
morse violation} and therefore need to consider function graph intersections.
Proposition~\ref{prop: graph of g is not in U} would then immediately yield the
Morse property if all second order derivatives would be contained in
\(\rg_2(\param)\).

The subsequent subsection (Section~\ref{sec: proof of morse result}) finishes
the proof of Theorem~\ref{thm: morse function on efficient networks}.
To do so we carefully craft a thin set \(U\) as the zero set of a function \(F\)
which \(\rg\) may not intersect. Although \(\hat {\mathbf J}(\param)\) has
degenerate second order differentials in the last layer,
the additional deterministic term \(\regularizer(\param) + \|\response_\param\|_{\Pr_X}^2\) that
is strictly convex in the last layer helps us out. More explicitly, we
will design a real analytic function \(F\) taking an outcome of \((\param,
\rg_2(\param))\) to a real value in such a way that for all
\(\param\in\polyEfficient\)
\[
    \nabla \mathbf J(\param)=0
    \; \implies \;
    \det(\nabla^2 \mathbf J(\param))= F(\param,\rg_2(\param)).\ 
\]
Consequently, if there is a parameter \(\param\in\polyEfficient\) with 
\[ \nabla \mathbf J (\param)=0 \text{ \ and \ } \det(\nabla^2  \mathbf J (\param))=0,
\]
then we also found a solution to 
\[ \rg_1 (\param)=0 \text{ \ and \ } F(\param , \rg_2 (\param))=0.
\]
This is again a collection of \(\dim(\Theta)+1\) real equations in \(\dim(\Theta)\) variables and we formally conclude with Proposition~\ref{prop: graph of g is not in U} that, almost surely, no solutions exist. Note that we are now able to proceed since the latter equations only make use of the non-degenerate differentials of first and second order of \(\hat{\mathbf J}(\param)\).

\subsection{Analysis of \texorpdfstring{\((\hat{\Cost}(\param))\)}{(ˆJ(θ))}}

In this section we analyze the stochastic process \((\hat{\Cost}(\param))_{\param\in \paramSpace}\). We will
\begin{itemize}
    \item derive representations for differentials of \((\hat{\Cost}(\param))_{\param\in \paramSpace}\) (Lemma~\ref{lem: differentiability})

    \item show non-degeneracy of the combined vector
    \((\rg_1(\param),\rg_2(\param))\) for all \(\param\in\polyEfficient\), where
    \(\rg_1(\param)\) and \(\rg_2(\param)\) are constituted by all first order
    differentials and certain second order differentials  of \(\hat{\mathbf
    J}(\param)\), respectively (Proposition~\ref{prop: non-degenrate g})

    \item generalize the volume argument of \citet{adlerRandomFieldsGeometry2007} to our needs
    (Lemma~\ref{lem: generalization 11.2.10}).
\end{itemize}
The combination of all these results leads to Proposition~\ref{prop: graph of g is not in U} which will allow us to show the Morse property in the subsequent section.

\begin{lemma}[Differentiability]
    \label{lem: differentiability}
    Let \(k\in\nat\), \(\network = (\nodes, \activation)\) be an ANN with a \(C^k\) activation function
    \(\activation\) and \(\# V_\mathrm{out}=1\), let \(\mathfrak R\) be a family of \(L^1\)-integrable regression problems with  \(\Pr_X\) having compact domain and let \(\mathbf m\in \mathbf M\).
    Then \[\hat{\cost}_\problem(\param) = \langle \response_\param,
    \target_\problem \rangle_{\Pr_X}\] is in \(C^k\) and its partial derivatives satisfy
    \[
        \partial_\param^\alpha \hat{\cost}_\problem(\param)
        = \langle \partial_\param^\alpha \response_\param, \target_\problem\rangle_{\Pr_X}
    \]
    for all multi-indices \(|\alpha|\le k\).
\end{lemma}
Note that the lemma implies that in the standard setting all differentials of  \((\hat{\Cost}(\param))_{\param\in \paramSpace}\) define again Gaussian processes. This fact together with the representations for the differentials will be the basic tool in the analysis of  \((\hat{\Cost}(\param))_{\param\in \paramSpace}\) and we mostly will not give reference to the lemma when using it.

\begin{proof}
By assumption, one has that \(\E_\problem[\|Y\|] < \infty\). Recall that \(\target_\problem(x) = \E_\problem[Y \mid X=x]\) so that with the \(L^1\)-contraction property of the conditional expectation
    \begin{align}\label{eq8723}
        \int |f_\mathbf{m}(x)|\, \Pr_X(dx)
        = \E_\problem\bigl[|\E_\problem[Y\mid X]|\bigr]
        \le
        \E_\problem\bigl[\E_\problem[|Y| \mid X]\bigr]
        = \E_\problem[|Y|^2] < \infty.
    \end{align}
     
    Fix \(\param\in\Theta\) and \(v\in \Theta\).  By assumption, \(\Pr_X\) has
    compact support \(\Domain\) and the directional derivative \(D_v^\param
    \response\) in direction \(v\) in the \(\param\) component is a continuous
    mapping on \(\Theta\times \Domain\). In particular, it is  uniformly bounded
    on the compact set \(\overline{B(\param,\|v\|)}\times\Domain\), say by the
    constant \(C\). Consequently, for every \(t\in(0,1]\), one has that
     \begin{align}
        \nonumber
        \frac 1t(\hat{\cost}_\problem(\param + tv) - \hat{\cost}_\problem(\param))
        &= \frac1t\int (\response_{\param + tv}(x)- \response_\param(x)) \target_\problem(x)\, \Pr_X(dx)
        \\
        \nonumber
        \overset{\text{FTC}}&= \int \int_0^1 D_v^\param\response_{\param + stv}(x)  \, \target_\problem(x)\, ds\,\Pr_X(dx)
\end{align}
Now note that \(C |f_\mathbf{m}|\) is an integrable majorant due to~(\ref{eq8723}). Using that for every \(s\in [0,1]\) and \(x\in\domain\), \(\lim _{t\downarrow 0} D_v^\param\response_{\param + stv}(x)= D_v^\param\response_{\param}(x)\) by continuity of the differential it follows with dominated convergence that
\[\lim_{t\downarrow 0}
\frac 1t(\hat{\cost}_\problem(\param + tv) - \hat{\cost}_\problem(\param))= \int  D_v^\param\response_{\param}(x)  \, \target_\problem(x)\,\Pr_X(dx).
\]
Recall that \((\param,v)\mapsto D_v^\param\response_{\param}(x)\) is continuous and we get again with dominated convergence that the latter integral is continuous in the parameters \(\theta\) and \(v\). This proves that \(\hat J_\mathbf{m}\) is \(C^1\) and that the upper identity holds.

By induction, one obtains the general statement. The induction step can be carried out exactly as above by using that the assumptions imply that \(\Psi\) is \(k\)-times continuously differentiable as mapping on \(\Theta\times \domain\).
\end{proof}

As indicated before there are second order differentials that degenerate.
However, for all first order and some second order differentials this is not the
case. In the next step we will show this. 
We consider the stochastic processes
\(\rg_1=(\rg_1(\param))_{\param\in\param}\) and
\(\rg_2=(\rg_2(\param))_{\param\in\param}\) defined by
\begin{align}
        \label{eq: definition of g1}
    \rg_1(\param)&:=\nabla\Cost(\param) \qquad\text{and}\\
    \label{eq: definition of g2}
    \rg_2(\param)
    &:=\Bigl(
        \bigl(\partial_{\bias_j}^2\hat{\Cost}(\param)\bigr)_{j\in V_1},
        \bigl(\partial_{\bias_j}\partial_{\weight_{ij}}\hat{\Cost}(\param)\bigr)_{\substack{j\in V_1\\i\in V_0}},
        \bigl(\partial_{\weight_{ij}}\partial_{\weight_{kj}}\hat{\Cost}(\param)\bigr)_{\substack{j\in V_1\\ i,k\in V_0\\i\le k}}
    \Bigr),
\end{align}
where we assume some total order on the input neurons \(V_0\) such that \(i\le k\)
makes sense for \(i,k\in V_0\).
Note that \(\rg_1\) utilizes the un-centered \(\Cost\) to ensure \(\nabla\Cost(\param) = 0\)
translates to \(\rg_1(\param)=0\) whereas \(\rg_2(\param)\) utilizes the
`centered'\footref{footnote: centered}
\(\hat{\Cost}\). This is because \(\rg_2\) does not contain all second order differentials
and a translation function \(F\) is necessary to get from \(\rg_2\)
to \(F(\param, \rg_2(\param)) = \det(\nabla^2\Cost(\param))\). Constructing
\(F\) in turn is more straightforward with the `mean'
\(\|\response_\param\|_{\Pr_X}^2\) built into \(F\) (cf.~Section~\ref{sec: proof
of morse result}).

\begin{prop}
    \label{prop: non-degenrate g}
    For an ANN \(\network = (\nodes, \activation)\) let \(\# V_\tout=1\).
    We consider \(\rg_i\) as defined in \eqref{eq: definition of g1} and \eqref{eq: definition of g2}
    based on the standard Gaussian setting (Definition~\ref{def: standard model}).
    Then for every parameter \(\param\in\mathcal \polyEfficient\) the  Gaussian
    random vector \((\rg_1(\param), \rg_2(\param))\) is
    non-degenerate meaning that its covariance has full rank. 
\end{prop}

\begin{proof}
    Recall that \(\hat{\Cost}(\param) = \langle \response_\param, \rf\rangle_{\Pr_X}\). 
    Since the variance does not depend on the mean we can assume without loss of
    generality \(\rg_1(\param) = \nabla \hat{\Cost}(\param)\) in this proof. 
    Let \(I_1\) and \(I_2\) be  index sets such that\footnote{
        With slight misuse of the notation, we ignore the ordering of the
        differentials in the representation of \(\rg_2(\param)\).
    }
    \[
        \rg_1(\param)=\nabla \hat{\Cost}(\param)
        = (\partial_{\param_i}\hat{\Cost}(\param))_{i\in I_1}
        \text{ \ and \ }
        \rg_2(\param)
        = (\partial_{\param_i}\partial_{\param_j}\hat{\Cost}(\param))_{(i,j)\in I_2}.
    \] 
    To show that \((\rg_1(\param),\rg_2(\param))\) is non-degenerate it suffices to show that the only vector \((\lambda_\mathbf {i})_{\mathbf i\in I_1\cup I_2}\in \real^{I_1\cup I_2}\) for which the linear combination
    \[
        \sum_{i\in I_1}\lambda_i \partial_{\param_i}\hat{\Cost}(\param)
        + \sum_{(i,j)\in I_2}\lambda_{i,j}\partial_{\param_i}\partial_{\param_j}\hat{\Cost}(\param)
    \]
    has zero variance is \((\lambda_\mathbf {i})_{\mathbf i\in I_1\cup I_2}\equiv 0\). 
    Lemma~\ref{lem: differentiability} ensures that the differentials exist and
    that they can be moved into the inner product defining \(\hat{\Cost}\).
    This implies
    \begin{align*}
        &\Var\Biggl(
            \sum_{i\in I_1}\lambda_i \partial_{\param_i}\hat{\Cost}(\param)
            + \sum_{(i,j)\in I_2}\lambda_{i,j} \partial_{\param_i}\partial_{\param_j}\hat{\Cost}(\param)
        \Biggr)
        \\
        &= \Var\Biggl(
            \Biggl\langle
            \underbrace{
                \sum_{i\in I_1}\lambda_i\partial_{\param_i}
                \response_\param
                + \sum_{(i,j)\in I_2}\lambda_{i,j}\partial_{\param_i}\partial_{\param_j}
                \response_\param
            }_{
                =:\phi
            },\;
            \rf
        \Biggr\rangle_{\Pr_X}\Biggr).
    \end{align*}
    By weak universality (Definition~\ref{def: standard model}) of the Gaussian
    process \(\rf\) on \(\Domain\) it follows that the latter variance is zero
    if and only if \(\phi\equiv 0\) on the support \(\Domain\) of \(\Pr_X\).
    It is therefore sufficient to prove that
    there exists no non-trivial linear combination of
    \[
        \Bigl(  \nabla \response_\param,
            \bigl(\partial_{\bias_j}^2\response_\param\bigr)_{j\in V_1},
            \bigl(\partial_{\bias_j}\partial_{\weight_{ij}}\response_\param\bigr)_{j\in V_1,i\in V_0},
            \bigl(\partial_{\weight_{ij}}\partial_{\weight_{kj}}\response_\param\bigr)_{j\in V_1, i,k\in V_0, i\le k}
        \Bigr),
    \]
    which is zero on \(\Domain\). We need to ensure that in \(\param\) all
    derivatives \(\partial_{\param_i}\) and
    \(\partial_{\param_i}\partial_{\param_j}\) (\(i\in I_1, (i,j)\in I_2\)) of
    the response \(\response_\param\)  are linearly independent as functions
    on~\(\Domain\). We recall that 
    \[
        \response_\param(x)
        = \bias_{\outSgt}
        + \sum_{j\in V_1}
        \activation\bigl(\bias_j + \langle x, \weight_{\placeholder j}\rangle\bigr)
        \weight_{j\outSgt}
    \]
   and note that we need to ensure linear independence of the following
    derivatives of the response \(\response_\param\)
    \begin{align}
        &x\mapsto 1 
        \tag{\(\partial \bias_\outSgt\)}
        \\
        &x \mapsto \activation\bigl(\bias_j + \langle x, \weight_{\placeholder j}\rangle \bigr)
        && j\in V_1
        \tag{\(\partial \weight_{j\outSgt}\)}
        \\
        &x \mapsto
        \activation'\bigl(\bias_j + \langle x, \weight_{\placeholder j} \rangle\bigr)
        \weight_{j\outSgt}
        && j\in V_1
        \tag{\(\partial \bias_j\)}
        \\
        &x \mapsto \activation'\bigl(\bias_j + \langle x, \weight_{\placeholder j} \rangle\bigr)
        \weight_{j\outSgt}x_i
        &&  j\in V_1,\; i\in V_0 
        \tag{\(\partial \weight_{ij}\)}
        \\
        &x \mapsto \activation''\bigl(\bias_j + \langle x, \weight_{\placeholder j} \rangle\bigr)
        \weight_{j\outSgt}
        && j\in V_1
        \tag{\(\partial \bias_j^2\)}
        \\
        &x \mapsto \activation''\bigl(\bias_j + \langle x, \weight_{\placeholder j} \rangle\bigr)
        \weight_{j\outSgt}x_i
        && i\in V_0, j\in V_1
        \tag{\(\partial \bias_j\partial \weight_{ij}\)}
        \\
        &x \mapsto \activation''\bigl(\bias_j + \langle x, \weight_{\placeholder j} \rangle\bigr) 
        \weight_{j\outSgt}x_i x_k
        && j\in V_1,\; i,k\in V_0
        \tag{\(\partial \weight_{ij}\partial\weight_{kj}\)}
    \end{align}
    Recall that \(\phi=0\) on \(\Domain\) is a linear combination of the derivatives above with the
    prefactors \((\lambda_\mathbf{i})\). To distinguish between the
    types of the indices we write \(\lambda_{\bias_*}\),
    \(\lambda_{\bias_j}\), \(\lambda_{\bias_j,\bias_j}\), \(\lambda_{w_{i,j}}\)
    and so on to refer to the coefficients in front of the respective
    differentials.  Note that in the above list all but the first differential
    can all be assigned to a particular neuron \(j\) in the first layer \(V_1\).
    For a fixed neuron \(j\in V_1\) we denote these contributions to \(\phi\)
    by \(\phi_j\) defined as
        \begin{align}\begin{split}\label{eq: polynomial neuron differentials}
        \phi_j(x)
        := 
            &\underbrace{\lambda_{\weight_{j\outSgt}}}_{=:P_j^{(0)}}
            \activation\bigl(\bias_j + \langle x, \weight_{\placeholder j}\rangle\bigr)
           + \underbrace{\Bigl(\lambda_{\bias_j} + \sum_{i\in V_0} \lambda_{\weight_{ij}}x_i\Bigr)}_{=:P_j^{(1)}(x)}
            \activation'\bigl(\bias_j + \langle x, \weight_{\placeholder j}\rangle\bigr)
            \\
            &+ \underbrace{\Bigl(\lambda_{\bias_j^2}
            + \sum_{i\in V_0}\lambda_{\bias_j\weight_{ij}}x_i
            + \sum_{i,k\in V_0:i\le k}\lambda_{\weight_{ij}\weight_{kj}}x_i x_k\Bigr)}_{=:P_j^{(2)}(x)}
            \activation''\bigl(\bias_j + \langle x, \weight_{\placeholder j}\rangle\bigr).
        \end{split}
    \end{align}  
   Together with \(P^{(\emptyset)}(x):=\lambda_{\bias_\outSgt}\) we therefore obtain 
   \[
        \phi(x)
        = \lambda_{\bias_\outSgt} + \sum_{j\in V_1}\phi_j(x)
        = P^{(\emptyset)}(x)+\sum_{j\in V_1}\sum_{k=0}^2 P_j^{(k)}(x)
        \activation^{(k)}\bigl(\bias_j+\langle x,\weight_{\placeholder j}\rangle\bigr).
    \]
    Since \(\phi(x)=0\) for all \(x\in \Domain\) we can conclude with polynomial
    efficiency of \(\param\) that all the polynomials \(P^{(\emptyset)}\) and
    \(P_j^{(k)}\) (\(j\in V_1,k=0,1,2)\) are zero.
    
    Now inspect the definition of the polynomials in \eqref{eq: polynomial neuron differentials} again.
    In the definition of every polynomial no monomial appears twice so that all
    coefficients have to be equal to zero. Since all coefficients appear in the
    polynomials, indeed all coefficients have to be equal to zero. This finishes
    the proof.
\end{proof}

We will combine the latter statement with a  generalization of a result of
\citet{adlerRandomFieldsGeometry2007}.  It will allow us to conclude that the non-degeneracy of the
distributions of \(\rg(\param):=(\rg_2(\param),\rg_1(\param))\) will allow us to
show that certain ``thin'' events (represented in terms of properties of the
function graph of \(\rg\)) have probability zero. In the final step, we will   define
appropriate ``thin'' events and verify the assumptions of the next lemma.

\begin{lemma}[Generalization {\textcite[Lemma~11.2.10]{adlerRandomFieldsGeometry2007}}]
    \label{lem: generalization 11.2.10}
    Let \(\dims,\dims'\in\nat\), \(U\subseteq \real^{d+d'}\) and
    \(W\subseteq\real^\dims\)  be measurable sets. Let \((\rg(w))_{w\in W}\) be
    an \(\real^{\dims'}\)-valued, Lipschitz-continuous stochastic process and suppose
    that for some constants \(C,\rho\in(0,\infty)\) one has for every \((w,v)\in
    U\cap(W\times \real^{\dims'})\) that the distribution of \(\rg(w)\) confined to
    \(B_{\real^{\dims'}}(w,\rho)\) is absolutely continuous w.r.t.\ Lebesgue measure
    with the density being bounded by \(C\). If
    \[
        \mathcal H_{\dims'}(U) =0,
    \]
    then the probability of the graph of \((\rg(w))_{w\in W}\) intersecting \(U\) is zero,
    i.e.
    \[
        \Pr\bigl(\exists w\in W: (w,\rg(w))\in U\bigr) = 0.
    \]
\end{lemma}
\begin{proof}
    Without loss of generality we can assume that \(U\subseteq W\times
    \real^{\dims'}\) (otherwise we replace \(U\) by  \(U\cap( W\times
    \real^{d'})\)).  Let \(L,\epsilon\in(0,\infty)\). Since by assumption
    \(\mathcal H_{d'}(U)=0\) there exist a \(U\)-valued sequence
    \((w_i,v_i)_{i\in\nat}\) and a \((0,\epsilon)\)-valued sequence
    \((r_i)_{i\in\nat}\) such that
    \[
        U\subseteq \bigcup_{i\in\nat} B((w_i,v_i),r_i) \text{ \ and \ } \sum_{i\in\nat} r_i^{d'} \le \epsilon.
    \] 
    Now note that for an arbitrary  Lipschitz function \(g:W\to \real^{\dims'}\)
    with Lipschitz constant \(L\) we have the following: if there exists \(w\in
    W\) with \((w,g(w))\in U\), then there exists an index \(i\in\nat\) with
    \(\|g(w_i)- v_i\|<(1+L)r_i\).  Indeed, in that case there exists an
    index \(i\in\nat\) with \(\|(w,g(w))-(w_i,v_i)\|<r_i\), since balls of
    radius \(r_i\) around \((w_i, v_i) \) cover \(U\), and
    together with the Lipschitz continuity we get that
       \begin{align*}
        \|g(w_i) - v_i\|
        &\le \|g(w_i)- g(w)\| + \|g(w)- v_i\|
        \\
        &\le  L\|w_i - w\| + \|g(w)- v_i\|
        \le (1+L) r_i.
    \end{align*}
    Consequently, we get that for the events
    \[
    \mathbb U=\{\exists w\in W: (w,\rg(w))\in U\}\text{ \ and \ } \mathbb L^{(L)}=\{\rg \text{ is $L$-Lipschitz cont.}\} 
    \]
    one has that
    \[ \mathbb U\cap \mathbb L^{(L)} \subseteq \bigcup_{i\in\nat} \bigl\{ \metric (\rg(w_i), v_i)<(1+L) r_i\bigr\}.
    \]
Now suppose that \(\epsilon\in(0,\infty)\) was chosen sufficiently small to
    guarantee that \((1+L)\epsilon<\rho\) and conclude using the Lebesgue
    measure \(\lebesgue^{\dims'}\) on \(\real^{\dims'}\) that
    \begin{align*}
        \Pr(\mathbb U\cap \mathbb L^{(L)} )
        &\le   \sum_{i=1}^\infty \Pr\bigl(\rg(w_i)\in B_{(1+L)r_i}(v_i) \bigr)
        \\
        &\le \sum_{i=1}^\infty \int_{ B_{(1+L)r_i}(v_i)} \frac{d\Pr_{\rg(x_i)}}{d\lebesgue^{d'}} d\lebesgue^{d'}
        \\
        &\le C\sum_{i=1}^\infty \lebesgue^{d'}(B_{(1+L)r_i}(v_i))
        \le  C \lebesgue^{d'}(B_{1}(0)) (1+L)^{d'} \epsilon.
    \end{align*}
    By letting \(\epsilon\) go to zero we conclude that \(\Pr(\mathbb U\cap
    \mathbb L^{(L)} )=0\). This is true for every \(L\in(0,\infty)\) and a
    union over rational \(L\) finishes the proof.
\end{proof}

\begin{prop}\label{prop: graph of g is not in U}
    Assume the standard setting (Definition~\ref{def: standard
    model}). Let \(\rg=(\rg(\param))_{\param\in\param}\) be the process given by
    \[
        \rg(\param)=(\rg_2(\param),\rg_1(\param)),
    \]
    where \(\rg_1\) and \(\rg_2\) are as defined in~\eqref{eq: definition of g1}
    and~\eqref{eq: definition of g2}. Let \(\dims,\dims'\in\nat\) be the dimensions of
    \(\param\) and  the target space of \(\rg\). Moreover, let \(U\subset
    \real^{\dims+\dims'}\) be a measurable set
    with
    \[
        \mathcal H_{d'}(U)=0,
    \]
    then 
    \[
        \Pr(\{\exists \param\in\polyEfficient: (\param,\rg(\param))\in U\})=0.
    \] 
\end{prop}

\begin{proof}
    Since \(\polyEfficient\) is an open set in \(\paramSpace\) which is
    separable as a finite dimensional real vector space, we can cover it by a
    countable number of compact balls contained in \(\polyEfficient\).
    Specifically, about any rational point in \(\polyEfficient\) we take a
    closed ball contained in \(\paramSpace\). Then it is sufficient to show
    the claim for any such ball as the countable union of zero sets is still
    a zero set. We thus consider such a compact subset \(W\subseteq
    \polyEfficient\) and aim to show
    \[
        \Pr(\{\exists \param\in W: (\param,\rg(\param))\in U\})=0.
    \]

    Since we have that \(\hat{\Cost}\) has continuous differentials up
    to third degree (Lemma \ref{lem: differentiability}),
    the process \(\rg\) is continuous as it only consists of first and second
    order differentials (and a continuous mean in the case of \(\rg_1\)). This
    then implies that the covariance kernel of \(\rg\) is continuous
    \autocite[e.g.][Theorem~3]{talagrandRegularityGaussianProcesses1987}.
    Moreover the third order differentials are continuous and thereby bounded on
    compact sets, which implies \(\rg\) is almost surely Lipschitz continuous on \(V\) and
    since the covariance kernel of \(\rg\) is continuous, the function
    \[
        \gamma(\param) = \det(\Cov(\rg(\param)))
    \]
    is continuous on \(W\) and therefore assumes a minimum in \(W\)
    as \(W\) is compact.  Since the covariance is positive definite, this
    minimum must be greater or equal than zero and by Proposition~\ref{prop:
    non-degenrate g} it cannot be zero since \(W \subseteq
    \polyEfficient\). And since \(\rg(\param)\) is Gaussian by Lemma~\ref{lem:
    differentiability} and assumption on \(\target_\Problem\) its Lebesgue density is bounded by the density at its
    mean given by
    \[
        (2\pi)^{-\dim(\paramSpace)/2} \det(\Cov(\rg(\param)))^{-1/2} \le (2\pi)^{-\dim(\paramSpace)/2} (\inf_{\param \in W}\gamma)^{-1/2}=:C< \infty.
    \]
    We can now finish the proof by application of Lemma~\ref{lem: generalization
    11.2.10}.
\end{proof}
\begin{remark}\label{rem: C3 too strong}
    While we have assumed analytic activation functions in the standard
    setting (Definition~\ref{def: standard model}) we only required
    the activation functions to be in \(C^3\) so far.  As we only work with the
    gradient and Hessian even the assumption of \(C^3\) activation functions
    appears too strong. And indeed in the proof above we only used this fact to
    show that \(\rg\) is almost surely Lipschitz.
    \citet{adlerRandomFieldsGeometry2007} highlights a similar issue after the
    proof of their Lemma~11.2.10, which we generalized in Lemma~\ref{lem:
    generalization 11.2.10}. \citet{adlerRandomFieldsGeometry2007}
    proceed to generalize their result using a growth condition on the modulus of
    continuity in place of Lipschitz continuity (cf. Lemma 11.2.11). A similar
    generalization should be possible for Lemma~\ref{lem: generalization 11.2.10}.
    But since we need analytic activation functions anyway for Lemma~\ref{lem: F
    definition} and therefore Theorem~\ref{thm: morse function on efficient
    networks}, we avoid the complications of this generalization.
\end{remark}

\subsection{Proof of Thm~\ref{thm: morse function on efficient networks}}
\label{sec: proof of morse result}

The main task of this section is to prove existence of the function \(F\) announced in the end of the introduction to Section~\ref{sec:Morse}.
We will show the following.

\begin{lemma}[Definition of \(F\)]
    \label{lem: F definition}
    In the standard setting (Definition~\ref{def: standard model})
    let  \(\rg_2=(\rg_2(\param))_{\param\in\param}\) be the
    \(\real^{I_2}\)-valued process as defined in~\eqref{eq: definition of g2},
    with \(I_2\) being the respective index set. Then there exists a non-zero,
    real-analytic function
    \[ 
        F:
            \polyEfficient \times \real^{I_2}
            \to \real,
    \]
    such that for every \(\param\in\polyEfficient =\polyEfficient(\Domain)\) with \(\nabla\Cost(\param) = 0\) we have that
    \[
        \det(\nabla^2 \Cost(\param))
        = F(\param, \rg_2(\param)).
    \]
\end{lemma}

Before we prove this lemma, we show that it finishes the proof of
Theorem~\ref{thm: morse function on efficient networks}.
We have that
\begin{align*}
    &\Pr\Bigl(\exists \param \in \polyEfficient: \nabla\Cost(\param) = 0,\; \nabla^2\Cost(\param)=0\Bigr)
    \\
    \overset{\text{Lemma~\ref{lem: F definition}}}&\le 
    \Pr\Bigl(\exists \param \in \polyEfficient: \nabla\Cost(\param) = 0,\; F(\param, \rg_2(\param))=0\Bigr)
    \\
    \overset{\rg_1(\param)=\nabla\Cost(\param)}&=
    \Pr\Bigl(\exists \param \in \polyEfficient: (\param, \rg_2(\param)) \in F^{-1}(0),\; \rg_1(\param) \in \{0\} \Bigr)
    \\
    &= \Pr\Bigl(\exists \param \in \polyEfficient: (\param, \rg(\param)) \in \underbrace{F^{-1}(0)\times \{0\}}_{=: U} \Bigr).
\end{align*}
To apply Proposition~\ref{prop: graph of g is not in U} we only need that \(U\) has
sufficiently small Hausdorff dimension (specifically smaller dimension than
the target space of \(\rg\)). And indeed, since \(F\) is a non-zero,
real-analytic function its zero set is one dimension smaller than the dimension \(d':=\#I_1+\# I_2\) of its domain, see
\autocite{mityaginZeroSetReal2020}. This entails that
\[
\mathcal H_{d'}(F^{-1}(0))=0\text{ \ \ and \ \ } \mathcal H_{d'}(U)=0.
\]
By definition, \(d'\) is also the dimension of the target space of \(\rg\)   and Proposition~\ref{prop: graph of g is not in U} entails that
\[
\Pr\bigl(\exists \param \in \polyEfficient: (\param, \rg(\param)) \in U\bigr).
\]
This finishes the proof of Theorem~\ref{thm: morse
function on efficient networks}.

\begin{remark}[Analytic activation]
    Observe that the analytic activation function was only used to ensure
    \(F\) is analytic (cf.~Remark \ref{rem: C3 too strong}). Lemma~\ref{lem: F definition} is
    therefore the appropriate place to search for generalizations.
\end{remark}

\begin{remark}[Efficient is sufficient]
    The function in Lemma~\ref{lem: F definition} can be defined on
    the efficient domain \(\efficient=\efficient(\Domain)\) (cf.~Definition~\ref{def: efficient}),
    i.e.  \(F:\efficient\times\real^{I_2}\to\real\). This is
    a superset of the polynomially efficient domain \(\polyEfficient\)
    in general and coincides with the efficient domain for certain activation
    functions (cf.~Theorem~\ref{thm: characterization of efficient networks}).
    We conduct the proof with \(\efficient\) but readers may replace
    this with \(\polyEfficient\).
\end{remark}
        
\begin{proof}[Proof of Lemma~\ref{lem: F definition}]
    To prove Lemma~\ref{lem: F definition} we need to construct a function \(F\)
    with the following properties
     \begin{enumerate}[label=\text{\normalfont (P\arabic*)}]
        \item\label{item: F is det hessian}
        at any \(\param\in\efficient\) with \(\nabla\Cost(\param) = 0\) we have
        \[
            \det(\nabla^2 \Cost(\param))
            = F(\param, \rg_2(\param)),
        \]
        \item\label{item: F is analytic}
        \(F\) is analytic, and

        \item\label{item: F non-zero}
\(F\) is non-zero, i.e. there exists an input to \(F\) where \(F\) is non-zero.
    \end{enumerate}
    Recall that we have by definition for some total order on \(V_0\)
    \[
        \rg_2(\param)
        =\Bigl(
            \bigl(\partial_{\bias_j}^2\hat{\Cost}(\param)\bigr)_{j\in V_1},
            \bigl(\partial_{\bias_j}\partial_{\weight_{ij}}\hat{\Cost}(\param)\bigr)_{\substack{j\in V_1\\i\in V_0}},
            \bigl(\partial_{\weight_{ij}}\partial_{\weight_{kj}}\hat{\Cost}(\param)\bigr)_{\substack{j\in V_1\\ i,k\in V_0\\i\le k}}
        \Bigr).
    \]
    In order to be able to plug \(\rg_2\) into \(F\) it must thus be of the form
    \[
        F:
        \efficient \times V_1 \times (V_1\times V_0) \times (V_1\times \text{Sym}(V_0^2)) \to \real,
    \]
    where \(\text{Sym}(V_0^2) = \{(i,k) \in V_0^2 : i\le k\}\).
    To satisfy \ref{item: F is det hessian} we need to reconstruct the
    determinant of the Hessian \(\nabla^2\Cost(\param)\) on the basis of the differentials in  \(\rg_2(\param)\) (that originate from the Hessian \(\nabla^2\hat{\Cost}(\param)\)). 
First recall that by Proposition~\ref{prop: decomposition}, one has that
    \begin{align*}
        \det(\nabla^2 \Cost(\param))
        &= \det\Bigl(\nabla^2 \Bigl[\regularizer(\param) + \|\response_\param\|_{\Pr_X}^2 - 2\hat{\Cost}(\param) + \E_\Problem[\|Y\|^2]\Bigr]\Bigr)
        \\
        &= \det\bigl(\nabla^2[\regularizer(\param)+ \|\response_\param\|_{\Pr_X}^2] - 2\nabla^2\hat{\Cost}(\param)\bigr).
    \end{align*}
    Since \(\param\mapsto \nabla^2[\regularizer(\param) +
    \|\response_\param\|_{\Pr_X}^2]\) is a deterministic function, we can absorb
    it into the definition of \(F\) and
    construct a function \(\tilde{F}\) that reproduces
    \(\nabla^2\hat{\Cost}(\param)\) from \(\rg_2(\param)\). That is, assuming we had
    a function \(\tilde{F}\) with \(\tilde{F}(\param, \rg_2(\param)) =
    \nabla^2\hat{\Cost}(\param)\) in all critical points, we can define
    \[
        F(\param, x)
        := \det\bigl(\nabla^2[\regularizer(\param) + \|\response_\param\|_{\Pr_X}^2] - 2\tilde{F}(\param, x)\bigr)
    \]
    If all entries of the matrix valued function \(\tilde{F}\) are analytic
    functions it is therefore sufficient for all entries of \(\param\mapsto\nabla^2 [\regularizer(\param) +
    \|\response_\param\|_{\Pr_X}^2]\) to be analytic for \(F\) to be
    analytic, because the determinant is a sum and product of these analytic
    components.
    And by the assumption that \(X\) has compact support in the
    standard setting (Definition~\ref{def: standard model}) and Theorem~5.1 in
    \citet{dereichConvergenceStochasticGradient2024}, \(\param\mapsto
    \|\response_\param\|_{\Pr_X}^2\) is analytic, and
    so are its differentials.
    \ref{item: F is analytic} thus follows if all entries of
    \(\tilde{F}\) are analytic since we assumed \(\response\) to be analytic in
    the standard setting (Definition~\ref{def: standard model}).

    Our strategy is therefore to show that all entries/differentials of
    \(\nabla^2\hat{\Cost}(\param)\) fall into one of the following categories:
    \begin{enumerate}[label=\text{\normalfont (\alph*)}]
        \item \label{item: regular}
        the partial differential is contained in  \(\rg_2\) in which case  the respective matrix entry in \(\tilde{F}\) is identical to the related input  (analytic)
        \item\label{item: second order is zero}
        the partial differential is zero (analytic), or
        \item\label{item: deterministic functions}
        there is an (analytic) deterministic function of \(\param\) (that we
        still need to construct) such that the partial differential coincides
        with the function value whenever
        \(\nabla\Cost(\param)=0\).
    \end{enumerate}
    
    To enact this strategy, we now consider all the second order derivatives in
    \(\nabla^2\hat{\Cost}(\param)\) that are not contained in \(\rg_2\) and
    categorize them into \ref{item: second order is zero} or \ref{item:
    deterministic functions}. Recall that by Lemma~\ref{lem: differentiability}
    \[
        \hat{\Cost}(\param)
        = \langle \response_\param, \rf\rangle_{\Pr_X},
        \quad
        \nabla \hat{\Cost}(\param)
        = \langle \nabla\response_\param, \rf\rangle_{\Pr_X},
        \quad
        \nabla^2 \hat{\Cost}(\param)
        = \langle \nabla^2\response_\param, \rf\rangle_{\Pr_X}.
    \]
    We therefore need to consider the derivatives of the
    response function. Since the response \(\response_\param\) is essentially a
    weighted sum over the index set \(V_1\) with all parameters appearing only
    in one of the summands, we have that second order partial derivatives that
    belong to two different neurons \(j\neq l\) of the hidden layer \(V_1\)
    vanish. This then directly implies that these differentials also vanish for
    \(\hat{\Cost}\).
    Specifically, we have that
    \begin{align*}
        \partial_{\bias_j}\partial_{\bias_l} \response_\param &= 0
        &&\forall j \neq l\in V_1
        \\
        \partial_{\bias_j}\partial_{\weight_{il}}\response_\param &= 0
        &&\forall j \neq l\in V_1,\; \forall i\in V_0
        \\
        \partial_{\weight_{ij}}\partial_{\weight_{kl}}\response_\param &= 0
        &&\forall j \neq l\in V_1,\; \forall i,k\in V_0 
        \\
        \partial_{\weight_{ij}}\partial_{\weight_{l\outSgt}}\response_\param &= 0
        &&\forall j \neq l\in V_1,\; \forall i\in V_0
        \\
        \partial_{\weight_{ij}}\partial_{\bias_\outSgt}\response_\param &= 0
        &&\forall j\in V_1,\; \forall i\in V_0
        \\
        \partial_{\bias_j}\partial_{\bias_\outSgt} \response_\param &= 0
        &&\forall j \in V_1
\end{align*}
    Let us refer to ``inner parameters'' as the parameters that appear inside the activation
    functions (the connections from the input layer to the hidden layer and
    the biases of neurons in the hidden layer). All second order
    derivatives of these inner parameters are either included in \(\rg_2\)
    or mix derivatives of two different hidden neurons and therefore vanish.
    All second order differentials with respect to inner parameters
    (exclusively) are thus of type \ref{item: regular} or \ref{item: second
    order is zero}.

    Since the response is linear in the outer (remaining)
    parameters \(\weight_{j\outSgt}\) and \(\bias_\outSgt\) taking two
    derivatives in these direction also results in zero, i.e.,
     \begin{align*}
        \partial_{\beta_{\outSgt}}^2 \response_\param = 0,
        \quad \partial_{\beta_{\outSgt}}\partial_{\weight_{j\outSgt}} \response_\param = 0
        \qquad \text{and}\qquad
        \partial_{\weight_{j\outSgt}}^2 \response_\param = 0
        \quad \forall j \in V_1.
       \end{align*}
    Most derivatives are thus in category \ref{item: second order is zero}.

    The only derivatives left are therefore the mixtures of outer derivatives
    \(\weight_{j\outSgt}\) with inner derivatives \(\bias_j\) and
    \(\weight_{ij}\) of the same hidden neuron \(j\in V_1\). These will be in category
    \ref{item: deterministic functions}. For those observe:
    \[
        \partial_{\weight_{ij}} \response_\param(x)
        = \weight_{j\outSgt}
        \activation'\bigl(
            \bias_j
            + \langle x, \weight_{\placeholder j}\rangle
        \bigr) x_i.
        \]
    Using that \(w_{j\outSgt}\) is non-zero as \(\param \in \efficient\) we get that
    \[
        \partial_{\weight_{j\outSgt}} \partial_{\weight_{ij}} \response_\param(x) =   \activation'\bigl(
            \bias_j
            + \langle x, \weight_{\placeholder j}\rangle
        \bigr) x_i
        =\tfrac 1{w_{j\outSgt}} 
        \partial_{\weight_{ij}}
        \response_\param(x).
    \]
    Since we are allowed to move differentiation into the
    inner products by Lemma \ref{lem: differentiability} we get
    \[
        \partial_{\weight_{j\outSgt}}\partial_{\weight_{ij}}
        \hat{\Cost}(\param)
        =
        \bigl\langle
            \partial_{\weight_{j\outSgt}} \partial_{\weight_{ij}} \response_\param(x),
            \rf
        \bigr\rangle
        = \tfrac 1{\weight_{j\outSgt}}
        \bigl\langle
            \partial_{\weight_{ij}} \response_\param(x),
            \rf
        \bigr\rangle= \tfrac 1{\weight_{j\outSgt}}
        \partial_{\weight_{ij}}
        \hat{\Cost}(\param).
    \]
    Now recall by Proposition~\ref{prop: decomposition} we have that
    \begin{equation}
        \label{eq: gradient decomposition} 
        \nabla\Cost(\param)
        = \nabla_\param[\regularizer(\param) +\|\response_\param\|_{\Pr_X}^2] - 2\nabla\hat{\Cost}(\param)
    \end{equation}
    so that in every critical point \(\param\) of \(\Cost\) with \(\nabla\Cost(\param)=0\)
    we get that
    \[
        \frac1{2\weight_{j\outSgt}} \partial_{\weight_{ij}}[
            \regularizer(\param) + \|\response_\param\|_{\Pr_X}^2
        ]
        = \frac1{\weight_{j\outSgt}} \partial_{\weight_{ij}}\hat{\Cost}(\param)
        = \partial_{\weight_{j\outSgt}}\partial_{\weight_{ij}}
        \hat{\Cost}(\param).
    \]
    For the definition of \(\tilde F\) we therefore use the deterministic
    function
    \[
        \param \mapsto \frac1{2\weight_{j\outSgt}}
        \partial_{\weight_{ij}}[\regularizer(\param) +
        \|\response_\param\|_{\Pr_X}^2]
    \]
    for the component where the differential
    \(\partial_{\weight_{j\outSgt}}\partial_{\weight_{ij}}\) should be.
    
      We proceed in complete analogy with the differential \(\partial_{\weight_{j\outSgt}}\partial_{\bias_j}\) and obtain that for every  critical efficient parameter \(\param\)
    \[
        \partial_{\weight_{j\outSgt}}\partial_{\bias_j}
        \hat{\Cost}(\param)
        = \frac{1}{2\weight_{j\outSgt}}\partial_{\bias_j}
        [\regularizer(\param)+ \|\response_\param\|_{\Pr_X}^2]
    \]
    and we use the respective deterministic function to define the corresponding component of \(\tilde F\).
    
	Note that  \(\param \mapsto \|\response_\param\|_{\Pr_X}^2\) is
    analytic, the deterministic functions that we use as substitutes for the
    remaining second order differentials in \(\tilde F\) are therefore analytic
    on \(\efficient\) (where \(\weight_{j\outSgt} \neq 0\)). Thus we constructed
    a function \(F\) satisfying properties
    \ref{item: F is det hessian} and \ref{item: F is analytic}.
    
    It remains to show  that \(F\) is a non-zero function \ref{item: F
    non-zero}.   
    To show this we will arrange the second order differentials appropriately.     
    We put the outer parameters
    \(\gamma := (\bias_\outSgt, (\weight_{j\outSgt})_{j\in V_1})\)
    at the end of the vector. Recall that all these components fall into
    category \ref{item: second order is zero} so that, in particular,
    \[
        \nabla_\gamma^2[\regularizer(\param) + \|\response_\param\|_{\Pr_X}^2] -2 \tilde F_\gamma(\theta,x)
        = \nabla_\gamma^2 [\regularizer(\param) + \|\response_\param\|_{\Pr_X}^2],
    \]
    where \(F_\gamma\) is \(F\) restricted to the output components with pairs from the \(\gamma\) coordinates.

    Let the remaining inner parameters be given by
    \[
        \alpha := \Bigl(
            (\bias_j)_{j\in V_1},
            (\weight_{ij})_{\substack{i\in V_0\\j\in V_1}}
        \Bigr).
    \]
    Recall that the second order differential with respect to two parameters
    from \(\alpha\) belongs to category \ref{item: regular} or \ref{item: second
    order is zero} and that the diagonal belongs to \ref{item: regular}. The
    diagonal can therefore be fully controlled and the other entries are either
    zero naturally or can be set to zero. Thus, for every parameter
    \(\param\in\efficient\) and any given \(\lambda \in \real\) we can find
    $x_\lambda$ with
    \[
        \tilde F_\alpha(\theta,x_\lambda)
        =  -\tfrac\lambda2 \identity,
    \]
    where \(\identity\) is the identity matrix and \(F_\alpha\) is \(F\)
    restricted to the output components with pairs from the \(\alpha\)
    coordinates. Consequently, for this choice of \(x\) we have that
    \[
        \nabla^2 \|\response_\param\|_{\Pr_X}^2-2 \tilde F(\theta,x_\lambda)
        =\begin{pmatrix}
            \nabla_\alpha^2[\regularizer(\param) + \|\response_\param\|_{\Pr_X}^2] + \lambda \identity & B(\param)\\
            B(\param)^\transpose & \nabla_\gamma^2[\regularizer(\param) + \|\response_\param\|_{\Pr_X}^2]
        \end{pmatrix}.
    \]
    Marked with a \(B(\param)\) are the mixed derivatives with parameters from
    \(\alpha\) and \(\gamma\). These are either of type \ref{item: second order
    is zero} or \ref{item: deterministic functions} and in particular they 
    are functions of \(\param\). More than \(F\) being non-zero, we will
    show for any fixed \(\param\) there exists \(\lambda\) and thus
    \(x_\lambda\) such that \(F(\param, x_\lambda) \neq 0\). For this note
    that \(F(\param, x_\lambda)\) is given by the determinant of the equation
    above and the determinant of a block matrix is given by
    \[
        \det\begin{bmatrix}
            A & B\\
            B^\transpose & D
        \end{bmatrix}
        = \det(D)\det(A - BD^{-1}B^\transpose)
    \]
    We will show that \(D := D(\param):=\nabla_\gamma^2 [\regularizer(\param)+ \|\response_\param\|_{\Pr_X}^2]\) is a
    strictly positive definite matrix and by doing so we show that it has full
    rank and therefore non-zero determinant. Since the eigenvalue \(\lambda\) of
    \(\lambda \identity\) can be directly controlled with the selection of
    \(x_\lambda\), selecting a sufficiently large \(\lambda \gg 0\) ensures
    that the eigenvalues of
    \[
        A- BD^{-1}B^\transpose
        = \lambda \identity + \underbrace{\nabla_\alpha^2[\regularizer(\param) + \|\response_\param\|_{\Pr_X}^2]
        - B(\param) D^{-1}B(\param)}_{\text{some matrix}}
    \]
    are all strictly positive. Thus we can ensure the second determinant is
    non-zero.
    
    What is left is thus the proof that \(\nabla_\gamma^2[\regularizer(\param)+
    \|\response_\param\|_{\Pr_X}^2]\)
    strict positive definite. Since \(\regularizer(\param)\) is assumed to
    be convex \(\nabla_\gamma^2\regularizer(\param)\) is positive semi-definite.
    It is thus sufficient to prove strict positive definiteness for
    \(\nabla_\gamma^2\|\response_\param\|_{\Pr_X}^2\).
    
    Let \(c\in \real^{\# V_1+1 }\) and recall \(\gamma=(\bias_\outSgt,
    (\weight_{j\outSgt})_{j\in V_1})\in\real^{\# V_1+ 1}\). Using that the
    iterated differentials from \(\gamma\) belong to category \ref{item: second
    order is zero} we conclude that
    \begin{align*}
       c^\transpose\nabla_{\gamma}^2\|\response_\param\|_{\Pr_X}^2c
       &= \sum_{i,j=0}^{\#V_1} c_i c_j \partial_{\gamma_i}\partial_{\gamma_j}\|\response_\param\|_{\Pr_X}^2
       \\
       &= \sum_{i,j=0}^{\# V_1} c_i c_j \partial_{\gamma_i}\partial_{\gamma_j}
       \int \response_\param^2(x) \,\Pr_X(dx)
       \\
       &= \int\sum_{i,j=0}^{\# V_1} c_i c_j \partial_{\gamma_i}\partial_{\gamma_j} \response_\param(x)^2 \,\Pr_X(dx).
       \\
       \overset{\partial_{\gamma_i}\partial_{\gamma_j}\response_\param \equiv 0}&= 2 \int\sum_{i,j=0}^{\# V_1} c_i c_j (\partial_{\gamma_i}\response_\param(x))(\partial_{\gamma_j} \response_\param(x))\, \Pr_X(dx).
       \\
       &= 2 \Bigl\|\sum_{i=0}^{\# V_1}c_i \partial_{\gamma_i} \response_\param\Bigr\|_{\Pr_X}^2.
    \end{align*}
   Consequently, \(\nabla_\gamma^2\|\response_\param\|_{\Pr_X}^2\) is always 
    positive definite. 
   To see strict positive definiteness we analyze solutions \(c\) for which the latter norm is zero.
    This requires that
    \begin{equation}
        \label{eq: sum of weighted outer parameter derivatives}    
        \sum_{i=0}^{\# V_1} c_i \partial_{\gamma_i} \response_\param = 0,
        \qquad \text{\(\Pr_X\)-almost surely.}
    \end{equation}
    Identifying \(V_1\) with \(\{1,\dots, \# V_1\}\) this implies that
    \[
        \partial_{\gamma_0} \response_\param
        = \partial_{\bias_\outSgt} \response_\param = 1
        \quad\text{and}\quad 
        \partial_{\gamma_i} \response_\param
        = \partial_{\weight_{i\outSgt}} \response_\param
        = \activation\bigl(\bias_i + \langle \weight_{\placeholder i}, x\rangle\bigr),
    \]
    for every \(i\in V_1\), 
   so that (\ref{eq: sum of weighted outer parameter derivatives}) implies that 
    \[
        c_0 +\sum_{i=1}^{\# V_1} c_i \,\activation\bigl(\bias_i + \langle \weight_{\placeholder i}, x\rangle\bigr)= 0
        \qquad \text{\(\Pr_X\)-almost surely.}
   \]
   Since the term above is continuous in \(x\) it is thus zero for all \(x\) in
   the support \(\Domain\) of \(\Pr_X\).
   Efficiency (Definition~\ref{def: efficient}) of \(\param\) then implies that
   \(c\equiv0\) is the unique solution of~(\ref{eq: sum of weighted outer
   parameter derivatives}).
This proves that
    \(\nabla_\gamma^2\|\response_\param\|_{\Pr_X}^2\) is strictly positive
    definite and therefore finishes the proof of \ref{item: F non-zero}.
\end{proof}

\section{Characterization of the efficient domain}
\label{sec: characterization of efficient params}

In the following \(\network=(\mathbb V,\activation)\) is a fixed ANN.
We start with a more natural axiomatic definition of efficient parameters
than the polynomial efficiency we required for our main result.

\begin{definition}[Efficient and redundant parameters] 
	\label{def: efficient}
	A parameter \(\param=(\weight,\bias)\) is called
	\emph{efficient}, if 
	\begin{enumerate}[label={(\alph*)}]
		\item \label{it: disused neuron}
		all neurons are used meaning that for all \(k\in V_1\) one has
        \[
			\weight_{k \placeholder} = (\weight_{kl})_{l\in V_\tout} \not\equiv 0,
		\]

		\item \label{it: linear combination of neurons}
		the equation
	\[
		\lambda_\emptyset + \sum_{j\in V_1} \lambda_j
		\activation\Bigl(\bias_j + \sum_{i\in V_0} x_i\weight_{ij}\Bigr)
		= 0 
        \qquad \forall x\in \Domain
    \]
	has the unique solution \(\lambda_j = 0\) for all \(j \in V_1 \uplus \{\emptyset\}\).
	\end{enumerate}
	We denote by \(\efficient=\efficient(\Domain)\) the \emph{efficient domain}, which is the set of
	all parameters \(\param=(\weight, \bias)\) that are efficient.	A
	parameter that is not efficient is called \emph{redundant}.
\end{definition}

In the remark below we introduce categories of redundant parameters
and hope to convey the intuition of this definition of `efficiency'.

\begin{remark}[Taxonomy of redundant parameters]
    \label{rem: taxonomy of redundant parameters}
    A parameter can be redundant for various reasons:
    If property \ref{it: disused neuron} does not hold we call it
    a \emph{deactivation redundancy} since the output of the hidden neuron
    \(k\) is ignored.
    If the property \ref{it: linear combination of
    neurons} does not hold there exists a neuron which can be
    linearly replicated by other
    neurons meaning that  there exists a neuron \(k\in V_1\) and  \(\lambda_j\in \real\) for all \(j\in V_1\uplus \{\emptyset\}\) with
    \[
        \activation\Bigl(\bias_k + \sum_{i\in V_0}  x_i\weight_{ik} \Bigr)
        = \lambda_\emptyset + \sum_{j\in V_1\setminus \{k\}}\lambda_j
        \activation\Bigl(\bias_j + \sum_{i\in V_0} x_i\weight_{ij}\Bigr)
        \quad \forall x\in \Domain.
    \]
    If  \(\lambda_\emptyset\) is the only \(\lambda_j\) not equal to zero in
    this linear combination, we speak of a \emph{bias redundancy}, as the
    neuron \(k\) is constant like the bias (typically \(\weight_{\placeholder k}
    = 0\), cf. Lemma~\ref{lem: bias redundancy characterization}). If there is another neuron \(j\in V_1\) such that \((\bias_j,
    \weight_{\placeholder j}) = (\bias_k, \weight_{\placeholder k})\)
    the neuron \(k\) can be trivially linearly combined from the others and we speak of a
    \emph{duplication redundancy}. We call all other cases a \emph{generalized duplication redundancy}.
    While deactivation, bias and duplication redundancies occur independently of
    the activation function (cf.~Lemma~\ref{lem: general redundancies}),
    generalized duplication redundancies only occur due to symmetries of
    the activation function (cf. Lemma \ref{lem: sigmoid sign symmetric
    redundancy}, Lemma \ref{lem: tanh sign symmetric redundancy} and Example
    \ref{ex: softplus}).
\end{remark}

Obviously, a configuration is $(0,0)$-polynomially efficient if and only if it
is efficient in the sense of Definition~\ref{def: efficient}. But in general the
assumption of polynomial efficiency is stronger. 

As we will show next both notions generally coincide in the case where the
activation function is either \(\sigmoid\) or \(\tanh\)! Further we will show
that they also coincide with the explicit set representation \eqref{eq:
defintion efficient 0}.

\begin{theorem}[Characterization of efficient networks for \(\sigmoid\) and \(\tanh\)]
    \label{thm: characterization of efficient networks}
    Assume that in the considered ANN \(\network\) the activation function
    \(\activation\) is either  \(\sigmoid\) or \(\tanh\), further assume
    that \(\Domain\) contains an open set. Then for every \(n\in\nat_0\) and
    \(m=(m_{\emptyset}, m_0, \dots, m_n)\in\nat_0^{n+2}\) one
    has
    \[
        \efficient(\Domain)=\polyEfficient^m(\Domain)=\efficient_0,
    \]
   where, as in \eqref{eq: defintion efficient 0},
    \[
        \efficient_0
        := \Biggl\{
            \param =(\weight, \bias)\in \real^{E \times (V\setminus V_0)}:
            {
            \begin{aligned}[c]
                \scriptstyle \weight_{j\placeholder}
                &\scriptstyle\neq 0
                &&\scriptstyle \forall j \in V_1,
                \\[-0.5em]
                \scriptstyle\weight_{\placeholder j}
                &\scriptstyle\neq 0
                &&\scriptstyle\forall j\in V_1, 
                \\[-0.5em]
                \scriptstyle(\weight_{\placeholder i}, \bias_i)
                &\scriptstyle\neq \pm(\weight_{\placeholder j}, \bias_j)
                &&\scriptstyle\forall i,j\in V_1\text{ with }i\not=j
            \end{aligned}
            }
        \Biggr\}.
    \]
\end{theorem}
\begin{proof}[Proofsketch]
    Since the activation function is either \(\sigmoid\)
    or \(\tanh\) and therefore real-analytic, the equations in the definition of
    the efficient set (Definition~\ref{def: efficient}) and the definition of
    polynomial independence (Definition~\ref{def: polynomial independence}) are
    real-analytic. Since analytic functions which are zero on an open
    set are zero anywhere, the open set contained in \(\Domain\) therefore
    ensures that we can assume without loss of generality \(\Domain =
    \real^{V_\tin}\)

    Suppressing \(\Domain\) in the notation of \(\efficient\) and
    \(\polyEfficient^m\) the proof is then established by showing that
    \(\efficient\subseteq \efficient_0\subseteq \polyEfficient^m\subseteq
    \efficient\). Note that \(\polyEfficient^m\subseteq \efficient\) is trivial
    as polynomials can always chosen to be constant.
\begin{itemize}
	\item
	To prove ``\(\efficient \subseteq \efficient_0\)'' we will show that any
	parameter \(\param\notin \efficient_0\) is not in \(\efficient\). More
	explicitly, we will construct a redundancy and show that one of the
	properties \ref{it: disused neuron} or \ref{it:  linear combination of
	neurons} does not hold.
    
    \begin{remark}
    As part of this proof in Subsection~\ref{subsec: explicit open set}, we will
    prove that \(\efficient\) is always contained in
    \begin{equation}
       \label{eq: bound efficient} 
        \bar{\efficient}
        :=\Biggl\{
            \param =(\weight, \bias)\in \real^{E \times (V\setminus V_0)}:
            {
            \begin{aligned}[c]
                \scriptstyle \weight_{j\placeholder}
                &\scriptstyle\neq 0
                &&\scriptstyle \forall j \in V_1,
                \\[-0.5em]
                \scriptstyle\weight_{\placeholder j}
                &\scriptstyle\neq 0
                &&\scriptstyle\forall j\in V_1, 
                \\[-0.5em]
                \scriptstyle(\weight_{\placeholder i}, \bias_i)
                &\scriptstyle\neq (\weight_{\placeholder j}, \bias_j)
                &&\scriptstyle\forall i\neq j\in V_1
            \end{aligned}
            }
        \Biggr\}
    \end{equation}
    regardless of the activation function \(\activation\) (Lemma~\ref{lem:
    general redundancies}). With a counterexample (Example~\ref{ex: softplus})
    we further show that the \emph{sign-symmetric redundancy} \((\weight_{\placeholder
    i}, \bias_i) = -(\weight_{\placeholder j}, \bias_j)\) is specific to the
    activation functions \(\activation\in \{\sigmoid, \tanh\}\). This shows
    that the explicit definition of \(\efficient_0\) does not generalize well.
    \end{remark}
    
    \item
    It will be harder to prove that ``\(\efficient_0\subseteq \polyEfficient^m\)''
    and we will first consider the case with one dimensional input (i.e.
    \(\# V_\tin =1\)) in Subsection~\ref{subsec: polynomial independence, 1-dim}.
    This proof relies on the complex poles of the meromorphic activation
    function \(\activation\in\{\sigmoid, \tanh\}\). We will then show that
    the one-dimensional result also implies the general result in
    Subsection~\ref{subsec: polynomial independence, general case}.
    
    \begin{remark}
        For the transfer from the one-dimensional result to the general result we
        we do not make use of the assumption \(\activation\in\{\sigmoid,
        \tanh\}\). We only use that the result holds for the \(1\)-dimensional case.
        This suggests that this part of the proof should be transferrable to other
        activation functions except for the fact that \(\efficient_0\) may be different
        for other activation functions in general and we use the specific structure
        of \(\efficient_0\).
    \qedhere
    \end{remark}
\end{itemize}
\end{proof}

\subsection{Proof of \texorpdfstring{\(\efficient\subseteq\efficient_0\)}{E ⊆ E0}}
\label{subsec: explicit open set}

Take \(\param\notin \efficient_0\), then it is sufficient to prove this
parameter to be not efficient, i.e. \(\param\notin\efficient\).
For this we are going to consider the possible constraints of \(\efficient_0\)
the parameter \(\param\) can violate and match them with the types of redundancies
we classified in Remark~\ref{rem: taxonomy of redundant parameters}. Recall
\[
    \efficient_0
    \overset{\text{def}}= \Biggl\{
        \param =(\weight, \bias)\in \real^{E \times (V\setminus V_0)}:
        {
        \begin{aligned}[c]
            \scriptstyle \weight_{j\placeholder}
            &\scriptstyle\neq 0
            &&\scriptstyle \forall j \in V_1,
            \\[-0.5em]
            \scriptstyle\weight_{\placeholder j}
            &\scriptstyle\neq 0
            &&\scriptstyle\forall j\in V_1, 
            \\[-0.5em]
            \scriptstyle(\weight_{\placeholder i}, \bias_i)
            &\scriptstyle\neq \pm(\weight_{\placeholder j}, \bias_j)
            &&\scriptstyle\forall i,j\in V_1\text{ with }i\not=j
        \end{aligned}
        }
    \Biggr\}.
\]
For \(\theta\notin \efficient_0\) one of the inequalities must be violated.
We consider all possibilities:

\begin{enumerate}
    \item \emph{Deactivation redundancy:}
    \label{it: deactivation redundancy} 
    If there is a neuron \(j\in V_1\) such that \(\weight_{j \placeholder} = 0\),
    then the parameter \(\param\) has a deactivation redundancy (Remark~\ref{rem: taxonomy of redundant parameters})
    as the output of the neuron is ignored and the parameter violates
    requirement \ref{it: disused neuron} of Definition~\ref{def: efficient}. Thus \(\param\notin \efficient\).

    \item \emph{Bias redundancy:}
    \label{it: bias redundancy}
    There is a neuron \(k\in V_1\) such that \(\weight_{\placeholder k} = 0\).
    Since this implies that the output of neuron \(k\) is constant
    and equal to \(\activation(\bias_k)\) irrespective of the input \(x\) it
    falls into the category of bias redundancies (Remark~\ref{rem: taxonomy of
    redundant parameters}). Specifically, the realization function can be
    replicated by removing the neuron \(k\) and adjusting the bias. This
    allows for a non-trivial linear combination
    \begin{equation}
        \label{eq:lin_indep1}    
    	\lambda_\emptyset + \sum_{j\in V_1} \lambda_j
		\activation\Bigl(\bias_j + \sum_{i\in V_0} x_i\weight_{ij}\Bigr)
		= 0 
    \end{equation}
    with
    \[
     \lambda_j=\begin{cases}
        - \psi(\beta_k), &\text{ if } j=\emptyset,\\
        1, &\text{ if } j=k,\\
        0, &\text{ if } j\in V_1\backslash\{k\}.
    \end{cases}
    \]
    This is a violation of \ref{it: linear combination of neurons} from
    Definition~\ref{def: efficient} and we thus have \(\param\notin \efficient\).

    \item \emph{Duplication redundancy:}
    \label{it: duplication redundancy}
    There are two neurons \(k,\ell\in V_1\) such that their parameters are equal
    \((\weight_{\placeholder k}, \bias_k) = (\weight_{\placeholder \ell},
    \bias_\ell)\). We call this a duplication redundancy since both neurons have
    identical parameters and therefore identical output. Here
    \[
        \lambda_j=\begin{cases}
            1, &\text{ if } j=k,\\
            -1,&\text{ if } j=\ell,\\
            0,&\text{ if } j\in \{0\}\cup(V_1\backslash\{k,\ell\})
        \end{cases}
    \]
    defines a nontrivial solution for \eqref{eq:lin_indep1}, violating \ref{it: linear combination of neurons}
    of Definition~\ref{def: efficient}. Thus \(\param\not\in \efficient\).
    Note that deactivation redundancies fall into the category of `generalized
    deactivation redundancies' in Remark~\ref{rem: taxonomy of redundant
    parameters}. 
\end{enumerate}

Observe that so far, we have not made use of \(\activation\in\{\sigmoid,
\tanh\}\). This leads to the following upper bound on the set of
efficient parameters regardless of the activation function.

\begin{lemma}[General redundancies]
    \label{lem: general redundancies}
    Let \(\network=(G, \activation)\) with \(G=(V,E)\) be a shallow
    neural network. Then the set of efficient parameters \(\efficient\)
    satisfies
    \[
        \efficient
        \subseteq\Biggl\{
            \param =(\weight, \bias)\in \real^{E \times (V\setminus V_0)}:
            {
            \begin{aligned}[c]
                \scriptstyle \weight_{j\placeholder}
                &\scriptstyle\neq 0
                &&\scriptstyle \forall j \in V_1,
                \\[-0.5em]
                \scriptstyle\weight_{\placeholder j}
                &\scriptstyle\neq 0
                &&\scriptstyle\forall j\in V_1, 
                \\[-0.5em]
                \scriptstyle(\weight_{\placeholder i}, \bias_i)
                &\scriptstyle\neq (\weight_{\placeholder j}, \bias_j)
                &&\scriptstyle\forall i\neq j\in V_1
            \end{aligned}
            }
        \Biggr\}
        =: \bar \efficient.
    \]
\end{lemma}
\begin{proof}
    The general arguments \ref{it: deactivation redundancy}, \ref{it: bias redundancy} and \ref{it: duplication redundancy}
    imply \(\bar{\efficient}^\complement \subseteq \efficient^\complement\).
\end{proof}

Continuing with our proof of \(\efficient\subseteq \efficient_0\)
there is one possible constraint violation left for \(\param\notin \efficient_0\):

\begin{enumerate}[resume]
    \item \emph{Sign-symmetric redundancy:}
    \label{it: sign-symmetric redundancy}
    There are two neurons \(i,j\in V_1\) such that \((\weight_{\placeholder i},
    \bias_i) = -(\weight_{\placeholder j}, \bias_j)\). The reason this results in
    a redundancy are symmetries of the activation function \(\activation\).
    Details in Lemma~\ref{lem: sigmoid sign symmetric redundancy} and \ref{lem: tanh sign symmetric redundancy}.
\end{enumerate}

\begin{lemma}[\(\sigmoid\)]
    \label{lem: sigmoid sign symmetric redundancy}
    Let \(\activation=\sigmoid\) and let $\param$ be a parameter
    such that there exist \(k,\ell\in V_1\) with \((\weight_{\placeholder k},
    \bias_k) =- (\weight_{\placeholder \ell}, \bias_\ell)\). Then \(\param\)
    is not efficient.
\end{lemma}

\begin{proof}
    Observe that we have for all \(x\in\real\) that
    \begin{equation}
        \label{eq: sigmoid sign symmetry}
        \sigmoid(-x)
= \frac{e^{-x}}{1+ e^{-x}}
        = \frac{1+e^{-x}}{1+e^{-x}} - \frac{1}{1+e^{-x}}
        = 1 - \sigmoid(x).
    \end{equation}
   This allows for a non-trivial linear combination 
    \[
    		\lambda_0 + \sum_{j\in V_1} \lambda_j
		\activation\Bigl(\bias_j + \sum_{i\in V_0} x_i\weight_{ij}\Bigr)
		= 0 
    \]
    via
    \[
        \lambda_j= \begin{cases}
            1 , &\text{ if } j\in\{k,\ell\},\\
            0,&\text{ if } j\in V_1\backslash\{k,\ell\},\\
            -1,&\text{ if } j=\emptyset.
        \end{cases}
    \]
    This nontrivial solution violates \ref{it: linear combination of neurons}
    of Definition~\ref{def: efficient} and thus implies \(\param\not\in \efficient\).
\end{proof}

\begin{lemma}[\(\tanh\)]
    \label{lem: tanh sign symmetric redundancy}
    Let \(\activation=\tanh\) and let $\param$ be a parameter
    such that there exist \(k,\ell\in V_1\) with \((\weight_{\placeholder k},
    \bias_k) =- (\weight_{\placeholder \ell}, \bias_\ell)\). Then \(\param\)
    is not efficient.
\end{lemma}
\begin{proof}
    Observe that we have for every \(x\in\real\)
        \begin{equation}
        \label{eq: tanh sign symmetry}    
        \tanh(-x) = \frac{e^{-x} - e^{x}}{e^{-x} + e^x} = - \tanh(x).
    \end{equation}
    Similarly, to the proof of Lemma~\ref{lem: sigmoid sign symmetric redundancy}
    we use this symmetry to construct a non-trivial linear combination via
     \[
        \lambda_j=\begin{cases}
            1 , &\text{ if } j\in\{k,\ell\},\\
            0,&\text{ if } j\in \{\emptyset\}\cup V_1\backslash\{k,\ell\}
        \end{cases}
    \]
    violating \ref{it: linear combination of neurons} of Definition~\ref{def: efficient}
    and finishing the proof.
\end{proof}

The redundancy that is treated in the Lemmas \ref{lem: sigmoid sign symmetric redundancy}
and \ref{lem: tanh sign symmetric redundancy} is caused by certain
symmetries in the particular activation function. In general, the particular
structure of the activation function can cause complex additional redundancies.

\begin{example}[Softplus]
    \label{ex: softplus}
Consider the case of the softplus activation
    function \(\activation(x) = \ln(1+\exp(x))\). If we have \((\weight_{\placeholder i}, \bias_i) =
    -(\weight_{\placeholder j}, \bias_j)\), then
    \begin{align*}
        \activation\bigl(\bias_i + \langle x,  \weight_{\placeholder i}\rangle \bigr)
        &= \ln\bigl(1+\exp\bigl[-\bigl(\bias_j + \langle x, \weight_{\placeholder j}\rangle\bigr)\bigr]\bigr)
        \\
        &= -\bigl(\bias_j + \langle x,\weight_{\placeholder j}\rangle \bigr)
        + \ln\bigl(1+\exp\bigl[\bias_j + \langle x, \weight_{\placeholder j}\rangle\bigr]\bigr)
        \\
        &= \bigl(\bias_i + \langle x, \weight_{\placeholder i}\rangle \bigr)
        + \activation\bigl(\bias_j + \langle x,\weight_{\placeholder j}\rangle\bigr).   
    \end{align*}
    So the neurons are equal up to a linear term, which can be cancelled out by another
    sign pair \((\weight_{\placeholder k},\bias_k) = -(\weight_{\placeholder l},
    \bias_l)\) and appropriate \(\weight_{k\placeholder}\) and
    \(\weight_{l\placeholder}\). For this it is important to keep in mind that
    we only need to take care of the first order term, as the zero order term
    can be absorbed by the bias \(\bias_m\) for \(m\in V_2\).

    Pruned networks in the case of the softplus function may therefore
    include one sign symmetry, but not more. The set of efficient parameters is
    therefore slightly larger. On the other hand polynomial efficiency
    results in the set \(\efficient_0\) again if \(P^{(\emptyset)}\) in Definition~\ref{def:
    polynomial independence} may be of degree \(1\) (i.e. \(m_\emptyset \ge 1\))
    as the linear term can be absorbed by a polynomial.
\end{example}

\subsection{Proof of \texorpdfstring{\(\efficient_0\subseteq \polyEfficient^m\)}{E ⊆ EmP} in the case of one dimensional input} \label{subsec: polynomial independence, 1-dim}

Let \(\param\in\efficient_0\). We need to show \(\param\) to be
\(m\)-polynomially independent (Definition~\ref{def: polynomial independence}).
In this first step we assume one dimensional input, i.e. \(V_\tin = \{\inSgt\}\).

For ease of notation, we define \(\alpha_j := \weight_{\inSgt j}\) for all \(j\in V_1\)
and write with slight misuse of notation \(x := x_{\inSgt}\). We thus have to prove that if the representation
\begin{equation}
    \label{eq:lindep2}
    0 = P^{(\emptyset)}(x)
    + \sum_{j\in V_1}
    \sum_{k=0}^n P^{(k)}_j(x) \activation^{(k)}\bigl(\bias_j + \alpha_j x\bigr) 
\end{equation}
is true for all \(x\in\real\),\footnote{
    Recall that we could translate \(x\in \Domain\) to \(x\in \real^{V_\tin}\) without
    loss of generality since the activation functions \(\tanh\) and
    \(\sigmoid\) are analytic and Theorem~\ref{thm: characterization of
    efficient networks} assumes an open set is contained in \(\Domain\).
} then 
all polynomials \(P^{(k)}_j\)  in the equation are zero. We will use complex analysis to show this.
Recall that the activation functions we consider are given by
\begin{align*}
    \sigmoid(x) &= \frac1{1+e^{-x}}\\
    \tanh(x) &= \frac{e^x -e^{-x}}{e^x + e^{-x}} = \frac{1- e^{-2x}}{1+ e^{-2x}}.
\end{align*}
Both activation functions are given as quotients of entire functions and are thus
meromorphic functions on \(\{x\in \complex: e^{-x}\not=-1\}\) and \(\{x\in
\complex: e^{-2x}\not=-1\}\), respectively. The singularities form
an infinite chain on the imaginary axis of the form
\[
    z_m = cm i \qquad (m\in \integer_\mathrm{odd}),
\]
where \(\integer_\mathrm{odd}\) denotes the odd integers and \(c=\pi\) in the
case of \(\sigmoid\) and \(c=\frac12\pi\) in the case of \(\tanh\).  

For every neuron $j\in V_1$ and $k=0,\dots,n$ the singularities of the meromorphic function
\[
    x\mapsto \psi^{(k)}(\alpha_j x+\beta_j)
\]
are of the form 
\[
    z_m^{(j)}
    =\frac{cm i-\beta_j}{\alpha_j}\qquad(m\in\integer_\mathrm{odd}).
\]
In particular, the singularities $S_j:=\{z_m^{(j)}:m\in\mathbb Z_\mathrm{odd}\}$ do not depend on $k$ and we call them the singularities of  neuron $j$. Consequently, the right hand side of (\ref{eq:lindep2}) always defines a meromorphic function on \[\complex\backslash\bigcup_{j\in V_1} S_j.\]
In order to have equality in (\ref{eq:lindep2}), in particular, all singularities have to be liftable.

Now suppose we are given a nontrivial solution to (\ref{eq:lindep2}). 
We will show that this would cause a contradiction as consequence of the following principles that we will state in detail and prove below:

\begin{enumerate}[label={\arabic*)}]
    \item \label{it: active neurons}
    We call a neuron \(j\in V_1\) \emph{active}, if one of the polynomials
    \(P_j^{(k)}\) with $k\in\{0,\dots,n\}$ is nonzero. If the set of
    active neurons is non-empty, there is an active neuron \(j^*\in V_1\) such
    that infinitely many of its singularities in \(S_{j^*}\) are not ``served''
    by another active neuron meaning that the singularity does not lie in one of
    the singularity sets of the other active
    neurons.
    
    \item\label{it: non-liftable singularity}
    If \(j\in V_1\), \(z\in S_j\) and \(k\in\{0,\dots,n\}\) with
    \(P_j^{(k)}(z)\not=0\), then the merophormorphic function
    \[
        \sum_{k=0}^n P_j^{(k)}(x) \psi^{(k)}(\alpha_j x+\beta_j)
    \]
    has a non-liftable singularity in \(z\).
\end{enumerate}
Indeed, the two facts then almost immediately imply a contradiction if a
non-trivial solution to (\ref{eq:lindep2}) would exist: we fix \(j^*\) as in
\ref{it: active neurons} and observe that for infinitely many \(z\in S_{j^*}\),
the singularity \(z\) is not served by other neurons. This entails that for all
these \(z\)

\[
    x\mapsto \sum_{k=0}^n P_{j^*}^{(k)}(x) \psi^{(k)}(\alpha_{j^*} x+\beta_{j^*})
\]
has  a liftable singularity in \(z\). Now \ref{it: non-liftable singularity} implies that for all these \(z\) and \(k=0,\dots,n\), \(P_{j^*}^{(k)}(z)=0\). Since there is no polynomial that satisfies this for an infinite number of points \(z\), we produced a contradiction showing that $j^*$ is not active.

It remains to prove the two statements. 
\begin{lemma}Let \(\theta\in\efficient_0\) and suppose there exists a nontrivial solution to (\ref{eq:lindep2}). There is an active neuron \(j^*\in V_1\) (meaning that there exists $k\in\{0,\dots,n\}$ with $P_k^{(j^*)}\not\equiv0$) such that 
\[ \# \Bigl(S_{j^*}\Big\backslash\bigcup_{j\text{ active}:\ j\not=j^*}S_j\Bigr)=\infty.
\]
\end{lemma}

\begin{proof}First note that if none of the neurons $j\in V_1$ would be active then one has that $0\equiv P^{(\emptyset)}$ which would imply that the solution is trivial.

Now fix an active neuron $j^*\in V_1$ that maximizes $|\alpha_j|$ over all active neurons~$j$ and denote by \(V_1^*\) the set of all active neurons with 
\[
 \frac{\beta_j}{\alpha_j}=\frac{\beta_{j^*}}{\alpha_{j^*}}.
\]
If there were another \(j\in V_1^*\backslash \{v^*\}\) with \(|\alpha_j|=|\alpha_{j^*}|\), then
we get together with 
\(  \frac{\beta_j}{\alpha_j}=\frac{\beta_{j^*}}{\alpha_{j^*}}
\)
 that \((\alpha_j,\beta_j)\) equals either \((\alpha_{j^*},\beta_{j^*})\) or \((-\alpha_{j^*},-\beta_{j^*})\). But this is not possible since \(\param\in\efficient_0\). Hence, for all \(j\in V_1^*\backslash \{v^*\}\) we have that \(|\alpha_j|<|\alpha_{j^*}|\).

Now let \(j\in V_1^*\backslash \{v^*\}\). Then 
 $z_m^{(j^*)}$ with \(m\in\integer_\mathrm{odd}\) is in $S_j$ iff there exists \(m'\in\integer_\mathrm{odd}\) with
\(m'/\alpha_j=m/\alpha_{j^*}\) or, equivalently,
\begin{align}\label{eq123}
\frac{ \alpha_{j^*}}{\alpha_{j}}m'=m.
\end{align} 
This entails that whenever \(\alpha_{j^*}/\alpha_{j}\) is not a rational we have that $S_j\cap S_{j^*}=\emptyset$. Now suppose that \(\alpha_{j^*}/\alpha_{j}\) is a rational
and let $p\in\integer$ and $q\in\nat$ such that the  representation  \(\alpha_{j^*}/\alpha_{j}=p/q\) is minimal (meaning that $|p|$ and $q$ are minimal). Since \(|\alpha_{j^*}|>|\alpha_{j}|\) we have that \(|p|\ge 2\). The integers $p$ and $q$ have no common divisor and equation  (\ref{eq123}) can only be true if $m'$ is a multiple of \(q\). Consequently, for all but at most one odd prime number \(m\) (namely \(|p|\)) we have that
\(z_m^{(j^*)}\not \in S_j\). Since there are only finitely many neurons in $V_1^*$ we conclude that all but finitely many odd primes \(m\) correspond to singularities \(z_m^{(j^*)}\) of neuron \(j^*\) that are not served by other active neurons.  
\end{proof}

\begin{lemma}
    \label{lem: retaining singularities}
    Let \(\alpha \in\real\backslash\{ 0\}\), \(\bias\in \real\) and \(z^*\in\complex\) be a singularity
    of a meromorphic function \(\activation\).  Let \(m\in\nat\) and for every  \(k\in \{0,\dots,m\}\) let
    \(P^{(k)}:\complex\to\complex\) be a polynomial such that either \(P^{(k)}
    \equiv 0\) or \(P^{(k)}\bigl(\frac{z^* - \bias}\alpha\bigr)\neq 0\).
    Then the polynomial
    \[
        \Phi(x) = \sum_{k=0}^m P^{(k)}(x) \activation^{(k)}(\alpha x + \bias),
    \]
     has a non-liftable singularity in  \(\frac{z^* - \alpha}\bias\) or
     \(P^{(k)} \equiv 0\) for all \(k=0,\dots,m\).
\end{lemma}
\begin{proof}
    By defining
    \[
        \tilde{\Phi}(y) := \Phi\Bigl(\frac{y-\bias}\alpha\Bigr)
        = \sum_{k=0}^m P^{(k)}\Bigl(\frac{y-\bias}\alpha\Bigr) \activation^{(k)}(y),
    \]
    and absorbing \(\alpha,\bias\) into the polynomials we can assume without loss of generality
    that \(\alpha=1\) and \(\bias=0\). Since these polynomials are zero if and only if
    the modified polynomials are zero.

    As a meromorphic function, all singularities are poles, i.e., there exists
    a smallest order \(r\) such that
    \[
        (z-z^*)^{r} \activation(z) 
    \]
    can be continuously extended in \(z^*\). Moreover its Laurent series is then
    of the form
    \[
        \activation(z) = \sum_{l=-r}^\infty a_l (z-z^*)^l 
    \]
    with \(a_{-m}\neq 0\). This representation allows us to argue that the order of the
    pole increases with every derivative and therefore we cannot cancel out the singularities
    with a simple linear combination. So either we retain the singularity, or we have to
    set everything to zero.
    More specifically, we have
    \begin{align*}
        \activation^{(k)}(z)
        &= \sum_{l=-r}^\infty a_ll\cdots(l-k+1) (z-z^*)^{l-k}
        \\
        &= \sum_{l=-(r+k)}^\infty a_{l+k}(l+k)\cdots(l+1) (z-z^*)^l
    \end{align*}
    Let us assume that \(P^{(m)}\neq 0\), then we have \(P^{(m)}(z^*)\neq 0\) and thus
    \begin{align*}
        &\Bigl|(z-z^*)^{r+m-1}\Phi(z)\Bigr|
        \\
        &\ge 
            \Biggl|
                \underbrace{|(z-z^*)^{r+m-1}P^{(m)}(z) \activation^{(m)}(z)|}_{\to \infty}
                - \underbrace{\Bigl|(z-z^*)^{r+m-1}\sum_{k=0}^{m-1}P^{(k)}(z) \activation^{(k)}(z)\Bigr|}_{\to c \in \real}
            \Biggr|
        \\
        &\to \infty \qquad (z\to z^*).
    \end{align*}
    Therefore \(\Phi\) has a singularity at \(z^*\) (more specifically a pole of
    order at least \(r+m\)).
    If \(P^{(m)} \equiv 0\), we repeat the same argument with \(m-1\) until we have
    either a pole at \(z^*\) or \(P^{(m)} \equiv \dots \equiv P^{(0)} \equiv 0\).
\end{proof}

\subsection{Proof of \texorpdfstring{\(\efficient_0\subseteq \polyEfficient^m\)}{E ⊆ E}, general case}
\label{subsec: polynomial independence, general case}

In this section we reduce the proof of \(\efficient_0 \subseteq \polyEfficient^m\)
to the one dimensional result we have already proven in Subsection~\ref{subsec:
polynomial independence, 1-dim}. For an element \(\param \in \efficient_0\)
recall that this requires that the equation
\[
    0 = P^{(\emptyset)}(x)
    + \sum_{j\in V_1}
    \sum_{k=0}^n P^{(k)}_j(x) \activation^{(k)}\bigl(\bias_j + \langle x, \weight_{\placeholder j}\rangle \bigr),
    \quad \forall x\in \real^{V_\tin}
\]
has the unique solution of all polynomials being zero (Definition~\ref{def:
polynomial independence}). Since this equation holds for all inputs \(x\in
\real^{V_\tin}\), it holds in particular for \(1\)-dimensional slices
\[
    x_v(\lambda):= \lambda v, \qquad \lambda\in \real
\]
for directions \(v\in \real^{V_\tin}\). The following proof hinges on the fact
that the mappings
\[
    \lambda \mapsto  P_j^{(k)}(x_v(\lambda))
\]
are polynomials in \(t\), which we can prove to be zero with the one dimensional
result. If sufficiently many appropriate directions \(v\) are chosen
and all directional polynomials are zero, then Theorem~\ref{thm: optimal
polynomial slicing} allows us to deduce that the original polynomial is zero.
But to apply the one dimensional result, we need to ensure that we do not
introduce degeneracies with the directions we choose. For example
\[
    \lambda\mapsto \langle x_v(\lambda), \weight_{\placeholder j}\rangle = \lambda \langle v, \weight_{\placeholder j}\rangle
\]
may be constantly zero if the direction \(v\) is orthogonal to
\(\weight_{\placeholder,j}\). This introduces a bias redundancy. For the number
of directions
\(
    N= \binom{\max(m) +\, \#V_\tin -1}{\max(m)}
\)
we therefore want to select \(v^{(1)}, \dots, v^{(N)}\in \real^{V_\tin}\) such that
the following conditions hold at the same time:
\begin{enumerate}
    \item the directions \(v^{(l)}\) characterize polynomials in the sense of
    Theorem~\ref{thm: optimal polynomial slicing}. This is required
    for us to deduce that the original polynomials are zero.
    \item \(\alpha_k^{(l)} := \langle v^{(l)}, \weight_{\placeholder k}\rangle  \neq 0\) for all \(k\in V_1\).

    \item \((\alpha_i^{(l)}, \bias_i) \neq \pm (\alpha_j^{(l)}, \bias_j)\) for all \(i, j\in V_1\) with \(i\neq j\).
\end{enumerate}
The last two requirements ensure non-redundancy for the parameters of
the surrogate neural networks.

\paragraph*{Why is this possible?} Select iid entries
\(v^{(l)}_i\sim\normal(0,1)\) with \(l\in \{1,\dots,N\}\) and \(i\in
V_\mathrm{in}\). Then we satisfy the conditions of Theorem~\ref{thm: optimal
polynomial slicing} and
almost surely have the first condition. Since \(\weight_{\placeholder k} \neq 0\)
for all \(k\in V_1\) by assumption, we also have almost surely
\[
    \alpha_k^{(l)} = \langle v^{(l)}, \weight_{\placeholder k} \rangle \neq 0.
\]
That is we have the second condition almost surely. For the last condition,
we use the assumption \((\weight_{\placeholder i}, \bias_i)\neq
\pm(\weight_{\placeholder j}, \bias_j)\). Let us only consider
the case where ``\(\pm\)'' is ``\(+\)''. The other
case is analogous. Then by assumption, we have
\[
        (\weight_{\placeholder i}, \bias_i) - (\weight_{\placeholder j}, \bias_j) \neq 0
\]
and thus either \(\weight_{\placeholder i}-\weight_{\placeholder j} \neq 0\) (case 1)
or \(\bias_i - \bias_j \neq 0\) (case 2).  This implies
\[
    (\alpha_i^{(l)}, \bias_i) - (\alpha_j^{(l)}, \bias_j) 
= \bigl(
        \langle v^{(l)}, \underbrace{
            \weight_{\placeholder i} - \weight_{\placeholder j}
        }_{\overset{\text{case 1}}{\neq} 0}
        \rangle,
        \underbrace{\bias_i - \bias_j}_{\overset{\text{case 2}}{\neq} 0}
    \bigr)\neq 0.
\]
Note that due to the iid selection of the entries of \(v^{(l)}\) the inequality of the
tuple with zero only holds almost surely in case 1.

In summary, with the selection of random \(v^{(l)}\) we can satisfy all three conditions
almost surely. In particular, there \emph{exist} directions \(v^{(1)}, \dots, v^{(N)}\)
which satisfy all three conditions simultaneously.

\paragraph*{\(1\)-dimensional slices of the \(\# V_\tin\)-dimensional input}

For the network \(\network = (\mathbb V, \activation)\)  we consider the \(1\)-dimensional network
\(\tilde{\network} := (\tilde{\mathbb{V}}, \activation)\) whose input layer is reduced to a
single node
\[
    \tilde{\mathbb{V}} := (\{\inSgt\}, V_1, V_\tout).
\]
From parameters \(\param=(\weight, \bias)\) of \(\network\) and direction \(v^{(l)}\) we
construct parameters \(\param^{(l)}=(\weight^{(l)}, \bias)\) of \(\tilde{\network}\) by retaining
the bias \(\bias\) and all connections from the hidden layer \(V_1\) to the
output \(V_\tout\) that is \(\weight^{(l)}_{\placeholder k} :=
\weight_{\placeholder k}\) for all \(k\in V_\tout\). For the input layer
connections, we set
\[
    \weight^{(l)}_{\inSgt j} := \alpha_j^{(l)} = \langle v^{(l)}, \weight_{\placeholder j}\rangle
    \qquad \forall j\in V_1.
\]
Then \(\param^{(l)} \in \efficient_0\) since we have
\begin{align*}
    \weight^{(l)}_{\inSgt i} &\neq 0 \quad \forall i\in V_1
    &&\text{due to } \alpha_j^{(l)} \neq 0,\\
    \weight^{(l)}_{i\placeholder } &\neq 0 \quad \forall i\in V_1
    &&\text{due to } \weight_{i \placeholder}^{(l)} = \weight_{i \placeholder} \neq 0,\\
    (\weight^{(l)}_{\inSgt i}, \bias_i) &\neq \pm (\weight^{(l)}_{\inSgt j}, \bias_j) \quad \forall i\neq j\in V_1
    &&\text{due to } (\alpha_i^{(l)}, \bias_i) \neq \pm (\alpha_j^{(l)}, \bias_j).
\end{align*}

\begin{remark}[Response function slices]
    The response functions \(\response_{\param^{(l)}}\) of the new parameter
    \(\param^{(l)}\) has the following relation with the response
    \(\response_\param\)
    \[
        \response_{\param^{(l)}}(\lambda)
        = \response_\param(\lambda v^{(l)})
        = \response_\param(x_{v^{(l)}}(\lambda))
        \qquad \forall \lambda \in \real.
    \]
\end{remark}

\paragraph*{Using the slices for the proof of polynomial independence}
We are now finally ready to prove polynomial independence for multi-dimensional
input. To do so, assume we have
\[
    0 = P^{(\emptyset)}(x)
    + \sum_{j\in V_1}
    \sum_{k=0}^m P^{(k)}_j(x) \activation^{(k)}\Bigl(\bias_j + \langle x, \weight_{\placeholder j}\rangle \Bigr) 
    \quad \forall x\in \real^{V_\tin}.
\]
In particular we can select \(x = \lambda v^{(l)}\) to obtain
\[
    0 = P^{(\emptyset)}(\lambda v^{(l)})
    + \sum_{j\in V_1}
    \sum_{k=0}^m P^{(k)}_j(\lambda v^{(l)}) \activation^{(k)}\Bigl(\bias_j + \weight^{(l)}_{\inSgt j} \lambda \Bigr) 
    \quad \forall \lambda\in \real.
\]
By the polynomial independence of the \(1\)-dimensional input network slices
\(\network^{(l)}\), we then have that the \(1\)-dimensional polynomial
slices
\[
    \lambda \mapsto P_j^{(k)}(\lambda v^{(l)})
\]
are all identically zero. As we selected the directions \(v^{(1)},\dots, v^{(N)}\)
to characterize polynomials in the sense of Theorem~\ref{thm: optimal
polynomial slicing}, we thus have \(P_j^{(k)}\equiv 0\) for all \(j\in V_1\) and
all \(k=\bias,0,\dots, n\). That is, polynomial independence for the case of
multivariate input.

\subsection{Polynomial slicing}

The main tool to translate the \(1\)-dimensional input result to the general
case are slices of polynomials that characterize the full polynomial.
This is formalized in the following theorem, which is proven in the remainder
of this section.

\begin{theorem}[Optimal polynomial slicing]
    \label{thm: optimal polynomial slicing}
    Let \(\dims, n\in \nat\) and \(N=\binom{n+\dims -1}{n}\).
    
    \begin{enumerate}[label=\text{\normalfont(\Roman*)}]
        \item \label{item: almost all selections work}
        \textbf{Almost all selections of directions \(v_1,\dots,v_N\) characterize the \(\dims\)-variate polynomials of degree \(n\).} If the
        matrix \((v_1,\dots, v_N)\in \real^{\dims\times N}\) is selected
        randomly with a density with respect to the Lebesgue measure on
        \(\real^{\dims\times N}\) (in particular there exist
        such \(v_i\)), then the following property holds almost surely:
        
        For any \(\dims\)-variate polynomial \(p\in
        \real[x_1,\dots,x_\dims]\) of order \(n\) we have \(p\equiv 0\)
        if and only if all slices in the directions \(v_i\) are zero, i.e.
        \[
            p_{v_i}(\lambda) := p(\lambda v_i) = 0 \quad \forall \lambda \in \real,\quad\forall i=1,\dots,N.
        \]

        \item \label{item: optimality of number of directions}
        \textbf{The number \(N\) of directions is optimal.} That is,
        for any smaller selection of directions \(v_1,\dots, v_M\in\real^\dims\) with \(M<N\)
        there exists a \textbf{non-zero} \(\dims\)-variate polynomial
        \(p\in \real[x_1,\dots,x_\dims]\) of order \(N\) such that all the
        slices \(p_{v_i}\)
        are identically zero.

    \end{enumerate}
\end{theorem}

A crucial object in the proof of this theorem is the vector of all monomials of
degree \(n\). In the multivariate case, the definition of such a vector requires
an ordering of the tuple of powers.
\begin{definition}[Vector of monomials]
    \label{def: vector of monomials}
    We define
    \[
        \mon_{n}(x) = \Biggl(
            \prod_{i=1}^\dims x_i^{r_i}
        \Biggr)_{r\in R} \quad\text{for}\quad x\in \real^\dims
    \]
    where \(R\) is a subset of \(\dims\)-tuples of non-negative integers
    that form monomials of exactly degree \(n\)
    \[
        R = \Bigl\{
            r \in \nat_0^\dims : \sum_{i=1}^\dims r_i = n
        \Bigr\}.
    \]
    With the following injection into the ordered set of
    non-negative integers \(\nat_0\)
    \[
        \phi : \begin{cases}
            R \to \nat_0
            \\
            r \mapsto \sum_{i=1}^{d}r_{i}(n+1)^{i-1},
        \end{cases}
    \]
    we equip \(R\) with the pullback of this order.  That is we define \(r <
    \tilde{r}\) if and only if \(\phi(r) < \phi(\tilde{r})\).
    In other words, we assume \(R\) has reverse lexicographic ordering.

    This ensures \(\mon_n(x)\) to be a vector and not just an unordered set.
    \end{definition}

\begin{prop}[Independent monomials]
    \label{prop: independent monomials}
    Let \(\dims, n\in \nat\) and \(N=\binom{n+\dims -1}{n}\).
    Then there exist \(v_1,\dots, v_N\) such that \(\mon_n(v_i)\) are linearly
    independent.

    Almost all selections of \(v_1,\dots, v_N\) have this property. That is,
    if the directions \((v_1,\dots, v_N)\in \real^{\dims\times N}\) are
    random variables with a density with respect to the Lebesgue measure
    on \(\real^{\dims\times N}\), then \(\mon_n(v_i)\) are almost surely
    linearly independent.
\end{prop}

\begin{proof}
    Using \(0< a_1<\dots < a_N\) with \(a_i\in \real\),
    we define
    \[
        v_k
        = \Bigl(a_k, a_k^{n+1}, \dots, a_k^{(n+1)^{\dims-1}}\Bigr)
    \]
    Then we have by definition of \(v_k\), \(\mon_n\) and \(\phi\)
    \[
        \mon_n(v_k)
        = \Biggl(\prod_{i=1}^\dims (v_k^{(i)})^{r_i}\Biggr)_{r\in R}
        = \Biggl(\prod_{i=1}^\dims a_k^{r_i(n+1)^{i-1}}\Biggr)_{r\in R}
        = (a_k^{\phi(r)})_{r\in R}.
    \]
    Therefore we have
    \begin{equation}
        \label{eq: generalized vandermonde}
        (\mon_n(v_1),\dots, \mon_n(v_N))
        = \begin{pmatrix}
            a_1^{\lambda_1} & \cdots & a_N^{\lambda_1}
            \\
            \vdots & & \vdots
            \\
            a_1^{\lambda_N} & \cdots & a_N^{\lambda_N}
        \end{pmatrix},
    \end{equation}
    with \((\lambda_1, \dots, \lambda_N) = (\phi(r))_{r\in R} \subseteq \nat_0\), where
    the size of the set \(R\) is given by Lemma~\ref{lem: number of monomials} and we
    have \(\lambda_1 < \dots < \lambda_N\) by the ordering defined for \(R\)
    in Definition~\ref{def: vector of monomials}. But the matrix
    \eqref{eq: generalized vandermonde} is a generalized Vandermonde
    matrix as in Lemma~\ref{lem: det generalized vandermonde} and its
    determinant is thus not equal to zero by Lemma~\ref{lem: det generalized vandermonde}. The
    monomial vectors are therefore linearly independent.

    Observe that \(\det(\mon_n(v_1),\dots, \mon_n(v_N))\) is a multivariate
    polynomial in the entries of \(v_i\). In particular it is a (real) analytic
    function. By \citet{mityaginZeroSetReal2020} the zero set of a real analytic
    which is not identically zero is a Lebesgue zero set. Since we found an
    example above where this determinant is non-zero, we ruled out that the
    function is identically zero. Thus almost all selections of \(v_1,\dots,
    v_N\) result in a non-zero determinant.
\end{proof}

\begin{lemma}[Number of monomials]
    \label{lem: number of monomials}
    \(|R| = N = \binom{n+\dims -1}{n}\)
\end{lemma}
\begin{proof}
    The number \(N\) is also sometimes referred to as ``\(\dims\) multichoose \(n\)''
    and denoted by
    \[
        \left(\!\!\binom{\dims}{n}\!\!\right)
        = \binom{n+\dims -1}{n}
    \]
    as it is equal to the number of ways to create a multiset of size \(n\) from
    \(\dims\) elements. In our case, we are picking an \(n\)-sized multiset of
    \(x_i\) to finally multiply all elements together to obtain a monomial.
    Details of the proof can be found in textbooks such as
    \citet[25-26]{stanleyEnumerativeCombinatoricsVolume2011} or \citet[Sec.
    3.2]{riordanIntroductionCombinatorialAnalysis2002}
\end{proof}

\begin{lemma}[Generalized Vandermonde]
    \label{lem: det generalized vandermonde}
    Let \(0< a_1<\dots< a_N\) for \(a_i \in \real\) and \(\lambda_1<\dots<\lambda_N\)
    with \(\lambda_i\in \real\). Then we have that
    the following generalized Vandermonde matrix has non-zero determinant, i.e.
    \[
        \det\begin{pmatrix}
            a_1^{\lambda_1} & \cdots & a_N^{\lambda_1}
            \\
            \vdots & & \vdots
            \\
            a_1^{\lambda_N} & \cdots & a_N^{\lambda_N}
        \end{pmatrix} \neq 0.
    \]
\end{lemma}
\begin{proof}
    The proof is adapted from a stack exchange answer \autocite{szwarcHowProveThat2022}.
    We conduct an induction over \(N\) and note that for the induction start
    \(N=1\) the conclusion obviously holds.

    For the induction step \(N-1 \to N\), assume that there exist
    \(c_1,\dots, c_N\in \real\) such that the rows weighted
    by \(c_i\) of the generalized Vandermonde sum to zero, i.e. we have for all 
    columns \(k\)
    \[
        c_1 a_k^{\lambda_1} + \dots + c_N a_k^{\lambda_N} = 0.
    \]
    In order to prove linear independence of these rows and thus that
    the determinant is zero, we only need to show that these equations
    imply \(c_i =0\) for all \(i=1,\dots, N\).

    Dividing the equations above by \(a_k^{\lambda_1}\), we observe
    that the \(a_k\) are zeros of the function
    \[
        f(x) = c_1 + c_2 x^{\lambda_2 - \lambda_1} + \dots + c_N x^{\lambda_N-\lambda_1}
    \]
    Since \(\lambda_k - \lambda_1>0\) by assumption, \(f\) is a continuously
    differentiable function. Between any two points where \(f\) is zero
    there is therefore a point where its derivative
    \[
        f'(x)
        = c_2(\lambda_2 - \lambda_1)x^{\lambda_2 - \lambda_1-1}
        + \dots
        + c_N(\lambda_N - \lambda_1)x^{\lambda_N - \lambda_1 - 1}
    \]
    is zero by the mean value theorem. In the gaps of \(a_1 <  \dots < a_n\) are thus
    \(N-1\) points \(0< u_1< \dots, u_{N-1}\) such that \(f'\) is zero at all \(u_i\).
    Since by induction hypothesis we have
    \[
        \det\begin{pmatrix}
            u_1^{\tilde{\lambda}_1} & \cdots & u_{N-1}^{\tilde{\lambda}_1}
            \\
            \vdots & & \vdots
            \\
            u_1^{\tilde{\lambda}_{N-1}} & \cdots & u_N^{\tilde{\lambda}_{N-1}}
        \end{pmatrix} \neq 0
    \]
    for \(\tilde{\lambda}_i = \lambda_{i+1} - \lambda_1 -1\), we have that the
    rows of this matrix are linearly independent. And since
    we have that all \(u_k\) are zeros of \(f'\), we have for
    the weighted colum sums
    \[
       0 = c_2(\lambda_2-\lambda_1) u_k^{\tilde{\lambda}_1}+ \dots + c_N(\lambda_N-\lambda_1) u_k^{\tilde{\lambda}_{N-1}}
    \]
    for all \(k\). By linear independence of the rows this implies
    that the coefficients \(c_i(\lambda_i-\lambda_1)\) for \(i\ge 2\) have  to
    be be zero. Since \(\lambda_i -\lambda_1 >0\), this
    implies \(c_2 = \dots = c_N = 0\). We thus have
    \(f \equiv c_1\) and since the \(a_k\) are zeros of \(f\)
    this also implies \(c_1=0\). Thus all \(c_i\) are zero which is what we needed to prove. 
\end{proof}

\subsubsection*{Proof of polynomial slicing (Theorem~\ref{thm: optimal polynomial slicing})}

\begin{enumerate}[wide]
    \item[\ref{item: almost all selections work}. ] 
    For the proof of the first statement of the theorem we intend to use the
    directions \(v_1,\dots, v_N\) of Proposition~\ref{prop: independent
    monomials} which result in linearly independent monomials. Note that these
    are only monomials of \emph{exactly} degree \(n\), while the polynomials of
    degree \(n\) admit all monomials \emph{up to} degree \(n\).

    We address this difference with an induction over the degree and by sorting
    the monomials of the polynomials into buckets with the same degree. The
    induction step is enabled by the following lemma.

    \begin{lemma}
        If the monomials \(\mon_m(v_1), \dots, \mon_m(v_N)\) span the space,
        then the lower degree monomials \(\mon_k(v_1), \dots, \mon_k(v_N)\) also
        span the space for all degrees \(k\leq m\).
    \end{lemma}
    \begin{proof}
        Without loss of generality assume \(k=m-1\).  Let \(K\) be the length of
        \(\mon_k(x)\) and \(M\) be the length of \(\mon_m(x)\) for \(x\in
        \real^\dims\). Choose any \(y\in \real^K\).
        We now have to prove that there is a linear combination of \(\mon_k(v_i)\)
        equal to \(y\).

        Observe that the vector
        \begin{equation}
            \label{eq: monomial subset}
            x_1 \mon_k(x)
            = x_1 \mon_{m-1}(x) \qquad x \in \real^\dims
        \end{equation}
        contains a subset of the entries of \(\mon_m(x)\). We obtain the vector
        \(\tilde{y} \in \real^M\) from \(y\in \real^K\) by setting the positions
        of all other entries to zero. Since the monomials \(\mon_m(v_1), \dots,
        \mon_m(v_N)\) span the space, there exists a linear combination
        \[
            \tilde{y} = \sum_{i=1}^N c_i \mon_m(v_i).
        \]
        By the observation \eqref{eq: monomial subset} this implies
        \[
            y = \sum_{i=1}^N \underbrace{c_i v_i^{(1)}}_{=:\tilde{c_i}}\mon_k(v_i).
        \]
        Thus the \(\mon_k(v_i)\) span the space.
    \end{proof}

    With this lemma we can now finish the proof of the first statement of the theorem.
    As already mentioned we take the directions \(v_1,\dots, v_N\) of
    Proposition~\ref{prop: independent monomials} and note that with this lemma
    we have that \(\mon_m(v_1),\dots, \mon_m(v_N)\) span the space for all
    \(m\leq n\). We proceed by induction over \(m\) up to \(n\). That is,
    we assume that the polynomial \(p\) is of degree \(m\) and assume that
    \(\mon_m(v_1),\dots, \mon_m(v_N)\) span the space but are not necessarily
    linearly independent (this is only the case for \(m=n\)).

    The base case \(m=0\) is trivially true as one direction is enough to figure out
    if a constant polynomial is zero.

    For the induction step \(m-1 \to m\) let \(p\) be a \(\dims\)-variate
    polynomial of degree \(m\). We decomposition the polynomial into
    \[
        p(x) = \sum_{k=0}^n p^{(k)}(x),
    \]
    where the polynomials \(p^{(k)}\) consist of all monomials of exactly degree \(k\).
    For all \(k< m\) we then have for all \(i\)
    \begin{equation}
        \label{eq: lower order monomials disappear}    
        \lim_{\lambda\to \infty} \frac{p^{(k)}(\lambda v_i)}{\lambda^m} = 0.
    \end{equation}
    With the assumption that the slices of \(p\) are zero, i.e. \(p_{v_i}(\lambda) = 0\)
    we thus obtain
    \begin{align}
        \nonumber
        0 &= \lim_{\lambda\to \infty} \frac{p_{v_i}(\lambda)}{\lambda^m}
        \overset{\text{def.}}= \lim_{\lambda\to \infty} \frac{p(\lambda v_i)}{\lambda^m}
        = \lim_{\lambda\to \infty} \sum_{k=0}^n \frac{p^{(k)}(\lambda v_i)}{\lambda^m}
        \overset{\eqref{eq: lower order monomials disappear}}= \lim_{\lambda\to \infty} \frac{p^{(m)}(\lambda v_i)}{\lambda^m}
        \\
        \label{eq: monomials of degree n zero}
        &= p^{(m)}(v_i).
    \end{align}
    For the last equation we used that \(p^{(m)}\) consists only of monomials of
    exactly degree \(m\), for which we have
    \begin{equation}
        \label{eq: monomial scaling} 
        \mon_m(\lambda x) = \lambda^m \mon_m(x).
    \end{equation}
    Since \(p^{(m)}\) consists only of monomials of exactly degree \(m\),
    there exists a vector \(q\in \real^M\) with \(M\) the length of
    \(\mon_m(x)\) for \(x\in \real^\dims\) such that
    \begin{equation}
        \label{eq: p as linear combination of monomials}
        p(x) = q^T \mon_m(x).
    \end{equation}
    With \eqref{eq: monomials of degree n zero} we thus have
    \[
        q^T \underbrace{(\mon_m(v_1), \dots, \mon_m(v_N))}_{\in \real^{M\times N}} = 0.
    \]
    As the monomials \(\mon_m(v_1), \dots, \mon_m(v_N)\) span the space
    \(\real^M\), the matrix has rank \(M\) and thus \(q=0\). By \eqref{eq: p as
    linear combination of monomials} we thus have \(p^{(m)}\equiv 0\). Therefore
    \(p\) is of degree \(m-1\) and we finish the proof of the
    first claim of the theorem using the induction assumption.

    \item[\ref{item: optimality of number of directions}. ]
    Let \(N = \binom{n+\dims -1}{n}\) and let \(v_1,\dots, v_M\) be fewer
    directions \(M<N\). To prove the statement, we construct a non-zero
    \(\dims\)-variate polynomial of degree \(n\) (more specifically it only
    consists of monomials of exactly degree \(n\)), which is zero in all
    directions \(v_i\).
    To do so, consider the matrix
    \[
        A:= (\mon_n(v_1), \dots, \mon_n(v_M)) \in \real^{N\times M}.
    \]
    Since \(M<N\) it is at most of rank \(M\) there exists \(0\neq q\in \real^N\)
    such that \(q^T A = 0\). We then define the polynomial \(p(x) = q^T \mon_n(x)\).
    Then by construction this polynomial is zero at all \(v_i\). Moreover, by the
    scaling property of the monomials \eqref{eq: monomial scaling} we have
    \[
        p_{v_i}(\lambda) = q^T \mon_n(\lambda v_i) = \lambda^n \underbrace{q^T \mon_n(v_i)}_{= p(v_i)} = 0. 
        \quad \forall \lambda \in \real.
    \]
    Thus we have \(p_{v_i} \equiv 0\) for all \(i=1,\dots,M\) but \(p\neq 0\)
    since \(q \neq 0\).
    
    The intuition is in essence, that the scaling property of the monomials
    \eqref{eq: monomial scaling} implies that we only collect a single information
    point for each direction \(v_i\). To ensure the polynomial is zero, we thus
    have to collect enough points \(v_i\) to ensure the polynomial has to be zero.
    As the space of polynomials of exactly degree \(n\) has dimension \(N\), this
    is the number of points required.
\end{enumerate}

\section{The neighborhood of redundant parameters}
\label{sec: neighborhood of redundant parameters}

In this section, we analyze the redundant domain.  We will show that every redundant parameter lies on a line of (redundant) parameters for which the realization function is identical. In the setting with \emph{no regularization}, i.e., \(R\equiv 0\), this entails that redundant parameters always have a degenerate Hessian (in the sense that its determinant is zero) and it can never be a strict local minimum.\smallskip

In a second step,  we will show that for redundancies that are not deactivation redundancies typically either all or no points on the latter line are critical points of the optimization landscape.

\begin{theorem}[Neighborhood of redundant critical points]
	\label{thm: critical points of redundant type} 	Let \(\network\) be an ANN and assume \(\param\) is redundant, i.e. \(\param
	\in \paramSpace\setminus \efficient(\Domain)\). Then there exists a straight line  \(\ell\subset \Theta\) containing \(\param\) such that for all \(\vartheta\in \ell\) 
	\[
		\response_{\vartheta} = \response_{\param}, \text{ \ on }\Domain.
\]    
   \end{theorem}
   
   \begin{proof}
    If \(\param\) has a deactivation redundancy (Definition~\ref{def: efficient}
    \ref{it: disused neuron}) and \(w_{k\bullet}= 0\) for a \(k\in V_1\), then changing the parameters \(w_{i,j}\) and \(\beta_j\) (\(i\in V_0, j\in V_1\)) does have no impact on the response and clearly the respective set contains a line.
    
  Now suppose that there is \(\lambda\not\equiv 0\) such that
      \begin{equation}
        \label{eq: b redundancy assumption} 
  \lambda_\emptyset + \sum_{j\in V_1}\lambda_j
        \activation\Bigl(\bias_j + \sum_{i\in V_0}x_i \weight_{ij} \Bigr)=0
        \quad \forall x\in \Domain.
    \end{equation}
    We define \(\param(t) = (\weight(t),\bias(t))\) for \(t\in
    \real\) as follows: We retain the weights connecting the input to the first
    layer and its biases, i.e.
    \[
        \weight_{ij}(t) := \weight_{ij}
        \text{ \ \ and \ \ }
        \bias_j(t):=\bias_j
        \qquad
        \forall i\in V_0, j\in V_1,
    \]
    and in the second layer we add multiples of \(\lambda\) in an appropriate way:
    \begin{align*}
        \weight_{jk}(t) &:= \weight_{jk} + t\lambda_j\text{ \ \ and \ \ }\bias_k(t):= \bias_k + t\lambda_\emptyset
        &&\forall j\in V_1, k\in V_2.
    \end{align*}
	Since \(\lambda\not=0\), the function \((\theta(t))_{t\in \real}\)
	parametrizes a line \(\ell\) that contains \(\param\) and basic linear
	algebra implies with \eqref{eq: b redundancy assumption} that the response
	does not depend on the choice of \(t\): In terms of
    \[
	\psi _j(x):= \activation\Bigl(\bias_j + \sum_{i\in V_0}x_i \weight_{ij} \Bigr)\qquad \forall j\in V_1, x\in\Domain,
	\]
	one has for  every \(k\in V_2\) and \(x\in\Domain\) that
     \begin{align}
        \label{eq: outer parameters influence on response}
        (\response_{\param(t)}(x))_k
        &= \bias_k(t) + \sum_{j\in V_1} \weight_{jk}(t) \activation_j(x)
        \\
        \nonumber
        &= \bias_k + \sum_{j\in V_1} \weight_{jk} \activation_j(x)
        + t\underbrace{
            \Bigl(
                \lambda_0
                + \sum_{j\in V_1} \lambda_j\activation_j(x)
            \Bigr)
        }_{\overset{\eqref{eq: b redundancy assumption}}=0}
        \\
        \nonumber
        &=(\response_{\param}(x))_k.
    \end{align}
    \end{proof}

\begin{theorem}\label{thm:672345}
  Let \(\network\) be an ANN, \(X\) and \(Y\) be
\(\real^{V_\tin}\)- and \(\real^{V_\tout}\)-valued random variables,
\(\ell: \real^{V_\tout}\times \real^{V_\tout}\to [0,\infty)\) a \(C^1\)-function
such that the cost landscape
\[
  \cost(\param)= \E[\ell(\Psi_\param(X), Y)]
\]
is \(C^1\) on \(\paramSpace\), and differentiation and integration can be interchanged.
Let the parameter \(\param=(w,\beta)\in\Theta\) exhibit a bias or duplication
redundancy (cf.\@ Remark \ref{rem: taxonomy of redundant parameters}) and let \(\param(t)\) be the parametrization of the line as
introduced after (\ref{eq: b redundancy assumption}).
Then either
\begin{itemize} \item for all \(t\in\real\), \(\param(t)\) is a critical parameter or
\item there it at most one  \(t\in\real\), for which \(\theta(t)\) is critical. 
\end{itemize}
If \(\# V_\tout =1\) and there are no deactivation redundancies, then
\(\param\) being critical implies that  \(\param(t)\) is critical for all  \(t\in
\real\).
\end{theorem}
\begin{proof} 
  In the  following, \(i\in V_0, j\in V_1,k,l\in V_2, t\in\real\) and \(x\in\Domain\) are arbitrary. Consider 
  \[
    \psi _j(x):= \activation\Bigl(\bias_j + \sum_{i\in V_0}x_i \weight_{ij} \Bigr)\text{ \ \ and \ \ }\psi_j'(x):=\psi'\Bigl(\bias_j + \sum_{i\in V_0}x_i \weight_{ij} \Bigr)
  \]
  Then we get for the inner differentials of the realization function
  \[
    \partial_{w_{i,j}} \bigl(\Psi_{\theta(t)}(x)\bigr)_l=\psi_j'(x)x_i w_{j,l}(t) \text{ \ \ and \ \ }  \partial_{\beta_{j}} \bigl(\Psi_{\theta(t)}(x)\bigr)_l=\psi_j'(x) w_{j,l}(t)
  \]
  and for the outer differentials
  \[
    \partial_{w_{j,k}} \bigl(\Psi_{\theta(t)}(x)\bigr)_l=\delta_{k,l} \psi_j(x)\text{ \ \ and \ \ }  \partial_{\beta_{k}} \bigl(\Psi_{\theta(t)}(x)\bigr)_l=\delta_{k,l},
  \]
  where \(\delta\) denotes the Kronecker-Delta. By assumption, we have that
  \begin{equation}
    \label{eq: cost representation}
    \nabla  \cost(\param(t))
    = \E\Bigl[
        \sum_{l\in V_\tout}
        \partial_{\hat{y}_l}\loss(\response_{\param(t)}(X), Y) \nabla ( \response_{\param(t)}(x))_l
    \Bigr].
  \end{equation}
  With the above identities we thus get for the inner and outer differentials
  \begin{align*}
    \partial_{w_{i,j}} \cost(\param(t))
    & = \sum_{l\in V_\tout} w_{j,l}(t)  \, \E\Bigl[
      \partial_{\hat{y}_l}\loss(\response_{\param}(X), Y)\, \psi_j'(X)\,X_i \Bigr]
    \\ 
    \partial_{\beta_{j}} \cost(\param(t))
    & = \sum_{l\in V_\tout} w_{j,l}(t)  \, \E\Bigl[
      \partial_{\hat{y}_l}\loss(\response_{\param}(X), Y)\, \psi_j'(X) \Bigr],
    \\
    \label{eq: outer weight derivative}
    \partial_{w_{j,k}} \cost(\param(t))
    & = \E\Bigl[
      \partial_{\hat{y}_k}\loss(\response_{\param}(X), Y)\, \psi_j(X) \Bigr]
    = \partial_{w_{j,k}} \cost(\param),
    \\
\partial_{\beta_{k}} \cost(\param(t))
    &= \E\Bigl[
      \partial_{\hat{y}_k}\loss(\response_{\param}(X), Y) \Bigr]
    = \partial_{\beta_{k}} \cost(\param).
  \end{align*}
  By the latter two identities,
  the derivatives with respect to the outer parameters do not depend on \(t\).
  The inner derivatives can be expressed in terms of
  \[
    a_{i,j,l}
    = \E\Bigl[          
      \partial_{\hat{y}_l}\loss(\response_{\param}(X), Y)\, \psi_j'(X)\,X_i 
    \Bigr]
    \text{ \ and \ }
    b_{j,l}
    =  \E\Bigl[          
      \partial_{\hat{y}_l}\loss(\response_{\param}(X), Y)\, \psi_j'(X)\,X_i
    \Bigr],
  \]
  specifically
  \begin{align*}
    \partial_{w_{i,j}} \cost(\param(t))
    &= \sum_{l\in V_\tout} w_{j,l}(t) \,a_{i,j,l}
    = \sum_{l\in V_\tout} (w_{j,l}+t\lambda _j) a_{i,j,l},
    \\
    &=\sum_{l\in V_\tout} w_{j,l}  a_{i,j,l} + t \lambda_j \sum_{l\in V_\tout}  a_{i,j,l} \\
    \partial_{\beta_{j}} \cost(\param(t))
    &= \sum_{l\in V_\tout} w_{j,l}(t) \, b_{j,l}
    =\sum_{l\in V_\tout} (w_{j,l}+t\lambda _j) b_{j,l}
    \\
    &=\sum_{l\in V_\tout} w_{j,l}  b_{j,l} + t \lambda_j \sum_{l\in V_\tout}  b_{j,l}.
  \end{align*}       
  Note that all differentials \(\partial_{w_{i,j}} \cost(\param(t))\) and
  \(\partial_{\beta_{j}} \cost(\param(t))\) are affine functions in \(t\).
  Hence, each differential is either zero for all \(t\in\real\) or at most one
  point \(t\). Consequently, on the line \((\param(t))_{t\in\real}\) either all
  points are critical or there is at most one point that is critical.
  In the case of \(\# V_\tout =1\), and no deactivation redundancies, i.e.
  \(\weight_{j\placeholder} \neq 0\), a critical point in \(t=0\) implies \(a_{i,j,l} = b_{j,l}=0\)
  for all \(i,j\) and thereby all \(\param(t)\) are critical.
\end{proof}

\section{Existence of efficient critical points}
\label{sec: existence of efficient critical points}

In this section, we analyze the existence of local minima in the standard setting
(Def.~\ref{def: standard model}). More explicitly, we prove that for every open
set \(U\subseteq \paramSpace\) containing a polynomially efficient parameter \(\param\)
one has with strictly positive probability that the random
unregularized squared error loss contains  a local minimum in~\(U\).   The
result illustrates that local minima may exist in the unregularized
setting. In the case of non-trivial regularization the cost typically tends to
infinity when the parameter \(\param\) tends to infinity. In this case the
existence of (local) minima is  trivial.

\begin{restatable}[Efficient minima exist with positive probability]{theorem}{efficientExist}
	\label{thm: existence of efficient minima}
	Assume that we are in the unregularised (i.e., \(R\equiv0\)) standard setting (Definition~\ref{def: standard model}) and that the random target function \(\rf=(f_{\mathbf M}(\param))_{\param\in\Theta}\) additionally satisfies that for all continuous functions \(\phi:\real^{V_\tin}\to \real\) and \(\delta \in(0,\infty)\) one has
	\[
	\Pr(\|\rf-\phi\|_{\Pr_X}<\delta)>0.
	\]

	Then every non-empty, open set \(U\) of \((0,0,1)\)-polynomially efficient parameters
contains a local minimum of
	the MSE cost with positive probability, i.e.
	\[
		\Pr\Bigl(
			\exists \param \in U:
			\param \text{ is a local minimum of } \Cost
		\Bigr) > 0.
	\]
\end{restatable}

\begin{proof}[Proof of Theorem~\ref{thm: existence of efficient minima}]

	Recall that for every \(\problem \in \problemSpace\) the MSE cost function is of the form
	\begin{align*}
		\cost_\problem(\param)
		&= \E_\problem[\|\response_\param(X) - Y\|^2]
		\\
		\overset{(*)}&= \E_\problem[\|\response_\param(X) - \target_\problem(X)\|^2]
		+ \underbrace{\E_\problem[\|\target_\problem(X)-Y\|^2 ]}_{\text{`noise' const. in }\param}.
	\end{align*}
	For \((*)\) we note that \(\E_\problem[Y-\target_\problem(X)\mid X] =0\) by
	definition of the target function \(\target_\problem(x) = \E_\problem[Y\mid X=x]\),
	and the mixed term therefore disappears.
Since we assumed \(\# V_\tout=1\), the norm is simply a square. Consequently, we get by interchanging differentiation and integration (this can be justified in complete analogy to Lemma~\ref{lem: differentiability}) that
	\begin{align}
		\label{eq: cost gradient}
		\nabla \cost_\problem(\param)
		&= 2\int(\response_\param(x) - f_\mathbf{m}(x))\nabla_\param\response_\param(x) \,\Pr_X(dx) \quad \text{and}
		\\
		\label{eq: cost hessian}
		\nabla^2 \cost_\problem(\param)
		&= 
			2\E\bigl[\nabla_\param\response_\param(X)\nabla_\param\response_\param(X)^\transpose\bigr]
			+ 2 \int (\response_\param(x) - f_\mathbf{m}(x))\nabla^2\response_\param(x) \Pr_X(dx).
	\end{align}
	Note that if the realization \(f_\mathbf{m}\) is very close to the response \(\Psi_{\param}\) for an  efficient parameter \(\param\in U\), then we informally have that
	\[
		\nabla \cost_\problem(\param)\approx 0
		\quad\text{and}\quad
		\nabla^2\cost_\problem(\param)\approx 2\E\bigl[
			\nabla_\param\response_\param(X)\nabla_\param\response_{\param}(X)^\transpose
		\bigr]
		=:2 G_{\param}.
	\]
	Our proof strategy is therefore to 
	\begin{itemize}
		\item show that \(G_\param\) is strictly positive definite,

		\item carry out a spectral analysis for  \(\nabla^2\cost_\problem(\param)\)

		\item show that in the case that the target function \( f_\mathbf{m}\) is
		close to \(\Psi_\param\), there exists a local minimum in the neighborhood of
		an efficient parameter \(\param_0\).
	\end{itemize}
\paragraph*{Strict positive definiteness of \(G_\param\).}
Let \(v\in \real^{\dim(\param)}\). We need to show that 
	 \[v^\transpose G_\param v\ge0 \text{ \ \  and \ \ }   [v^\transpose G_\param v=0 \ \Rightarrow \ v=0].\] 	 
	One has 
	\[
		 v^\transpose G_\param v
		= v^\transpose \E\bigl[
			\nabla_\param\response_\param(X)\nabla_\param\response_\param(X)^\transpose
		\bigr] v
		= \Bigl\| \langle v, \nabla_\param\response_\param(\cdot)\rangle \Bigr\|_{\Pr_X}^2\ge 0.
	\]
	Now suppose that \(v^\transpose G_\param v=0\). Then 
	\[
		\langle v, \nabla_\param\response_\param(\cdot)\rangle = 0, \text{ \ \ \(\Pr_X\)-almost surely.}
	\]
	 The latter function is analytic (since
	\(\response_\param\) is analytic and \(\Pr_X\) is compactly supported). Consequently,  it is zero on the entire support \(\Domain\)
	of \(\Pr_X\).

	Recall that \(\# V_\tout=1\) and let \(V_\tout = \{\outSgt\}\).  The
	derivatives of \(\response_\param\) are then given by
	\begin{align}
		&x\mapsto 1 
		\tag{\(\partial \bias_\outSgt\)}
		\\
		&x \mapsto \activation\bigl(\bias_j + \langle x, \weight_{\placeholder j}\rangle \bigr)
		&& j\in V_1
		\tag{\(\partial \weight_{j\outSgt}\)}
		\\
		&x \mapsto
		\activation'\bigl(\bias_j + \langle x, \weight_{\placeholder j} \rangle\bigr) \weight_{j\outSgt}
		&& j\in V_1
		\tag{\(\partial \bias_j\)}
		\\
		&x \mapsto \activation'\bigl(\bias_j + \langle x, \weight_{\placeholder j} \rangle\bigr) \weight_{j\outSgt}x_i
		&&  j\in V_1,\; i\in V_0 
		\tag{\(\partial \weight_{ij}\)}
	\end{align}
	We  thus get the representation 
	\[
		 \langle v, \nabla_\param \response_\param(x)\rangle
		= v_{\bias_\outSgt} + \sum_{j\in V_1} \Phi_j(x)
		\qquad \forall x \in \Domain
	\]
	with
	\[
		\Phi_j(x)
		:= \underbrace{v_{\weight_{j\outSgt}}}_{=:P_0^{(j)}}
		\activation\bigl(\bias_j + \langle x, \weight_{\placeholder j}\rangle\bigr)
		+ \underbrace{\Bigl(v_{\bias_j} w_{j*} + \sum_{i\in V_0} v_{\weight_{ij}}  w_{j*}x_i\Bigr)}_{=:P_1^{(j)}(x)}
		\activation'\bigl(\bias_j + \langle x, \weight_{\placeholder j}\rangle\bigr),
	\]
	where we write \(v_{\param_i}:= v_i\) (to easily refer to the particular
	parameters). Recall that since \(\param\) is \((0,0,1)\)-polynomially
	efficient the function \(\langle v, \nabla_\param
	\response_\param(\cdot)\rangle\) can only be zero on \(\Domain\), if all
	polynomials (`coefficients') are zero. By assumption, every \(w_{j*}\not=0\)
	so that we get that indeed all entries of the vector \(v\) are zero. Thereby we
	proved that \(G_\param\) is strictly positive definite and that the minimal
	eigenvalue \(\underline{\lambda} _\param\) of \(\param\) is strictly
	positive.	
	
	\paragraph*{Spectral analysis  of \(\nabla^2 \cost_\problem(\param)\).}
	Observe that in terms of 
	\[
		\overline \lambda_\param
		:= \sup_{\|v\|=1}\bigl\| v^T \nabla_\param^2\response_\param(\cdot)v \bigr\|_{\Pr_X}
		\le \Bigl(\int \| \nabla^2\Psi_\param (x)\|_{\text{op}}\,\Pr_X(dx)\Bigr)^{1/2}
	\] 
	one has for  \(v\in\paramSpace\) that by the Cauchy-Schwarz inequality
	\begin{align}
		\nonumber
		v^\transpose \nabla^2\cost_\problem(\param) \,v
		\overset{\eqref{eq: cost hessian}}&= 2\Bigl(
			v^\transpose G_\param v
			+ \int (\response_\param(x) - f_\mathbf{m}(x)) v^\transpose \nabla_\param^2\response_\param(x) v \;\Pr_X(dx) 
		\Bigr)
		\\
		\label{eq: lower bound hessian cost}
		&\ge 2\Bigl(
			\underline{\lambda} _\param
			- \|\response_\param - f_\mathbf{m}\|_{\Pr_X} \overline{\lambda}_\param
			\Bigr)\|v\|^2.
	\end{align}

	\paragraph*{Finding a local minimum.}
	We will prove that there exists a local minimum in \(B_\delta(\param_0)\)
	for \(\delta>0\), if there exists a lower bound on the spectrum of the
	Hessian \(\rho>0\) such that
	\begin{equation}
		\label{eq: sufficient for local minimum}
		\|\nabla \cost_\problem(\param_0)\| 
		< \tfrac{\delta}2 \rho
		\qquad \text{and}\qquad
		\bigl[\nabla^2\cost_\problem(\param) \succeq \rho
		\quad \forall \param \in B_\delta(\param_0)\bigr].
	\end{equation}
	Since \(\cost_\problem\) is thereby \(\rho\)-strongly convex on
	\(B_\delta(\param_0)\) \autocite[Thm.\@
	2.1.11]{nesterovLecturesConvexOptimization2018} we have for all \(\param \in
	B_\delta(\param_0)\) \autocite[Def.\@ 2.1.3]{nesterovLecturesConvexOptimization2018}
	\begin{align*}
		\cost_\problem(\param)
		&\ge \cost_\problem(\param_0) + \langle \nabla\cost_\problem(\param_0), \param- \param_0\rangle
		+ \tfrac\rho2 \|\param - \param_0\|^2
		\\
		&\ge \cost_\problem(\param_0)
		- \|\nabla \cost_\problem(\param_0)\| \|\param - \param_0\| + \tfrac\rho2\|\param-\param_0\|^2.
	\end{align*}
	For all parameters \(\param\) on the boundary \(\partial B_\delta(\param_0)\) of the ball
	we therefore have
	\[
		\cost_\problem(\param)
		\ge \cost_\problem(\param_0)
		- \|\nabla \cost_\problem(\param_0)\| \delta + \tfrac\rho2 \delta^2
		\overset{\eqref{eq: sufficient for local minimum}}> \cost_\problem(\param_0).
	\]
	The minimum, which the cost \(\cost_\param\) assumes on the (compact) closed ball
	\(\overline{B_\delta(\param_0)}\), can therefore not be on the boundary. Consequently,
	there must be a local minimum in \(B_\delta(\param_0)\).
	
	\paragraph*{Finishing the proof}
	By reducing the size of the open set \(U\) if necessary, we can
	assume without loss of generality that its closure \(\closure{U}\) is
	compact and also contained in the set of polynomially efficient parameters
	\(\polyEfficient^{(0,0,1)}\).
	Now suppose that \(\theta_0\) is an arbitrary efficient element of \(U\).
	The continuous maps \((\underline {\lambda}_\param)_{\param\in\paramSpace}\) and
	\((\overline {\lambda}_\param)_{\param\in\paramSpace}\) both attain their
	minimum and maximum on the compact set \(\overline U\), where
	all \(\underline\lambda_\param >0\) and all \(\overline\lambda_\param<\infty\) so that
	\[
		\underline {\lambda}
		= \min _{\param\in \overline U}\underline {\lambda}_\param>0 \text{ \ \ and \ \ }\overline {\lambda}=\max _{\param\in \overline U}\overline {\lambda}_\param<\infty.
	\]
	To satisfy \eqref{eq: sufficient for local minimum} with \(\rho:=\underline{\lambda}\),
	let \(\epsilon:= \underline{\lambda}/(2\overline{\lambda})\) and select
	\(\delta\)
	sufficiently small such that
	\begin{itemize}
		\item \(B_\delta(\param_0) \subseteq U\subseteq \polyEfficient^{(0,0,1)}\)
		(this ensures a minimum in \(U\))
		\item  \(\|\response_\param - \response_{\param_0}\|_{\Pr_X}\le
		\epsilon/2\) for all \(\param \in B_\delta(\param_0)\) (using continuity of \(\response_\param\)).
	\end{itemize}
	Then choose \(r\in (0, \epsilon/2)\) such that 
	\(\|\nabla_\param\response_{\param_0}\|_{\Pr_X} < \tfrac{\delta}{2r} \rho\).
	Consequently the inequality \(\|\response_{\param_0}-\target_\problem\|\le r\) implies
	\begin{enumerate}
		\item by \eqref{eq: cost gradient} and Cauchy's inequality
		\[
			\|\nabla\cost_\problem(\param_0)\|
			\le \|\response_{\param_0} - \target_\problem\|_{\Pr_X} \|\nabla\response_{\param_0}\|_{\Pr_X}
			< \tfrac{\delta}{2} \rho,
		\]
		\item
		and for all \(\param\in B_\delta(\param_0)\)
		\[
			\|\response_\param - \target_\problem\|_{\Pr_X}
			\le \|\response_\param - \response_{\param_0}\|_{\Pr_X} + \|\response_{\param_0} - \target_\problem\|_{\Pr_X}
			\le \epsilon .
		\]
		Using \(\epsilon= \underline{\lambda}/(2\overline{\lambda})\) and
		\eqref{eq: lower bound hessian cost} we can lower bound spectrum by \(\rho\),
		i.e.\ for all \(v\) such that \(\|v\|=1\)
		\[
			v^\transpose \nabla^2\cost_\problem(\param) \,v
			\ge 2(\underline{\lambda} - \epsilon
			\overline\lambda)=\underline\lambda\overset{\text{def.}}=\rho.
		\]
	\end{enumerate}
	This means we satisfy \eqref{eq: sufficient for local minimum}
	if \(\|\response_{\param_0}-\target_\problem\|\le r\).
	By assumption on the random statistical model
	\(\mathbf M\) the latter property holds with strict positive probability.
\end{proof}
 \section{Existence of redundant critical points}
\label{sec: existence of redundant critical points}

Since the set of redundant parameters is generally a thin set with respect to
the Lebesgue measure (e.g.\ \(\efficient_0\) in Theorem~\ref{thm: characterization of efficient networks}), one may
reasonably hope that this set does not contain any critical points of
the MSE with probability one. In that case the MSE would be a Morse function
over the entire set of parameters with probability one. Unfortunately, this
hypothesis is wrong in general as the following theorem shows. We further break
down the set of redundant parameters, using the taxonomy introduced in
Remark~\ref{rem: taxonomy of redundant parameters}, to make more precise
statements about the existence of redundancies which are required to be of a
certain type.

\begin{theorem}[Redundancies cannot be ruled out in general]
    \label{thm: redundanct critical points exist}
    Assume the standard setting (Definition~\ref{def: standard model})
    without regularization, i.e. \(\regularizer\equiv 0\).

    Assume \(\polyEfficient^{(0,0,1)}\) contains an open set\footnote{
        \label{footnote: e.g. sigmoid tanh}
        e.g. \(\activation\in \{\sigmoid, \tanh\}\) and an open set in the support of \(\Pr_X\) by Theorem~\ref{thm: characterization of efficient networks}
    } and that there is at least one hidden neuron (\(\#V_1\ge 1\)),
	then, with \textbf{positive probability}, critical points of the MSE
	\(\Cost\) do \textbf{exist} in the sets of
	\begin{enumerate}[label={\normalfont(\Alph*)},series=redundantSets]
		\item\label{it: rule out not possible - all} redundant parameters,
		\item\label{it: rule out not possible - dup}
		redundant parameters that only admit duplication redundancies (assuming \(\# V_1\ge 2\))
		\item\label{it: rule out not possible - bias and deac} 
        redundant parameters that only admit bias \emph{and} deactivation redundancies,
	\end{enumerate}
\end{theorem}

\begin{proof}[Proof (outline)]
    Clearly, the existence of redundant critical points \ref{it: rule out not possible - all}
    follows from the existence of critical points with more specific redundancies, i.e. \ref{it: rule out
    not possible - dup} or \ref{it: rule out not possible - bias and deac}. So we
    only need to prove \ref{it: rule out not possible - dup} and \ref{it: rule out not possible - bias and deac}.
    To do so, we make use of the fact that we have proven efficient
    critical points exist with positive probability (Theorem~\ref{thm: existence
    of efficient minima}). Using \(\#V_1\ge 1\) we can therefore find a critical
    point of a smaller network with \(\#V_1-1\)
    hidden neurons with positive probability. We then carefully \textbf{extend} this
    network and its parameters by a redundancy in a fashion that retains the
    criticality of the parameters. But for a duplication redundancy we obviously
    need at least two hidden neurons. Details follow in Section~\ref{sec: extending}.
\end{proof}

Conversely, pure bias redundancies can be ruled out.

\begin{prop}[Pure bias redundancies can be ruled out]
    \label{prop: pure bias redundancies can be ruled out}
    Assume the standard unregularized setting (Definition \ref{def: standard model}).
    If \(\activation'(x)\neq 0\) for all \(x\in \real\), the support \(\Domain\)
    of \(\Pr_X\) contains an open set and an efficient parameter is automatically \((1,0,1)\)-polynomially efficient,\footref{footnote: e.g. sigmoid tanh}
    then, with \textbf{probability one}, critical points of the MSE \(\Cost\)
    do \textbf{not exist} in the set of
    \begin{enumerate}[resume*=redundantSets]
        \item\label{it: no pure bias redundancies} redundant parameters that
        only admit bias redundancies.
    \end{enumerate}
    
\end{prop}
\begin{proof}[Proof (outline)]
    The proof of \ref{it: no pure bias redundancies} relies on \textbf{pruning} the bias
    redundancies to obtain an efficient parameter for a smaller network. This
    efficient parameter must then also be a critical point of the cost and satisfy
    an additional condition. We then show that there are almost surely no efficient
    critical points which satisfy this additional condition and thereby
    rule out critical bias redundancies. Details follow in Section \ref{sec: critical bias red do not exist} 
    after we outline a general pruning process in Section \ref{sec: pruning}.
\end{proof}

\subsection{Extending (Proof of Theorem \ref{thm: redundanct critical points exist})}
\label{sec: extending}

Using the following lemma to extend an efficient critical point
of a smaller network, \ref{it: rule out not possible - dup} and \ref{it: rule
out not possible - bias and deac} clearly follow from the existence of such
an efficient critical point in the smaller network positive probability.
This follows from Theorem~\ref{thm: existence of efficient
minima}, for which we require the existence of an open set in the set
\(\polyEfficient^{(0,0,1)}\).

\begin{lemma}[Extension]
    Assume the setting of Theorem \ref{thm: redundanct critical points exist}
    and without loss of generality \(V_1 = \{1,\dots, \# V_1\}\). Define
    the reduced ANN to be \(\tilde \network = (\tilde{\mathbb V}, \activation)\)
    with neurons \(\tilde{\mathbb V}:= \bigl(V_0, V_1\setminus\{1\},
    V_2\bigr)\).
    Assume that the parameter \(\tilde \param\) of the network \(\tilde \network\)
    is a critical point of \(\cost_\problem\).
    Then there exists a parameter \(\param\) of the network
    \(\network\) such that it is a critical point of \(\cost_\problem\) and
    either
    \begin{enumerate}
        \item\label{it: single dup} \(\param\) only has a single duplication redundancy and no
        other redundancies (if we further assume \(\#V_1\ge 2\)), or
        \item\label{it: only deac and bias} \(\param\) only has a deactivation and bias redundancy at the
        same neuron and no other redundancies.
    \end{enumerate}
\end{lemma}
\begin{remark}
    We will not make use of the fact that the loss \(\loss\) is the squared
    error. We only require sufficient regularity such that derivatives may be
    moved into the expectation (e.g. Lemma~\ref{lem: differentiability}).
\end{remark}

\paragraph*{Proof of \ref{it: single dup}.}

We are going to construct a parameter \(\param\) with
a single duplication redundancy from the parameter \(\tilde \param\) of the
reduced network. Assume that the parameters we do not mention are kept as is.
Our plan is to duplicate the neuron \(2\) so we define for all \(i\in V_0\)
\begin{equation}
    \tag{duplication}
    \weight_{i1} := \tilde\weight_{i2}
    \qquad\text{and}\qquad 
    \bias_1 := \tilde \bias_2.
\end{equation}
To ensure neither neuron is deactivated pick \(\lambda \in
\real\setminus\{0,1\}\) and
define
\begin{equation}
    \tag{`convex' combination}    
    \weight_{1l} := \lambda \tilde \weight_{2l}
    \qquad\text{and}\qquad
    \weight_{2l} := (1-\lambda)\tilde \weight_{2l}.
\end{equation}
Clearly, \(\param\) is in the set of parameters which only admit duplication
redundancies.

It is straightforward to see, that the response must remain the same, i.e.
\(\response_\param = \response_{\tilde{\param}}\), as we have just split one
identical neuron into a `convex' combination of two identical ones. That is
\begin{align*}
    &(\response_{\param}(x))_l
    \\
    &= \bias_l + \sum_{j=1}^{\# V_1} \activation\bigl(\bias_j + \langle x, \weight_{\placeholder j}\rangle\bigr)\weight_{jl}
    \\
    &= 
    \bias_l
    + \activation\bigl(\bias_{2} + \langle x, \weight_{\placeholder 2}\rangle\bigr)
    \underbrace{\bigl((1-\lambda)\weight_{1 l} + \lambda \weight_{2 l}\bigr)}_{
        = \tilde{\weight}_{2 l} 
    }
    + \smashoperator{\sum_{j=3}}^{\# V_1}\activation\bigl(\bias_j + \langle x, \weight_{\placeholder j}\rangle\bigr) \weight_{jl}
    \\
    &= \tilde{\bias}_l + \sum_{j\in V_1\setminus\{1\}} \activation\bigl(\tilde{\bias}_j + \langle x, \tilde{\weight}_{\placeholder j}\rangle\bigr)\tilde{\weight}_{jl}
    \\
    &=(\response_{\tilde{\param}}(x))_l.
\end{align*}
With this fact under our belt, we can now consider the derivatives
of the cost. Recall that
we denote by \(\partial_{\hat y_l}\loss\) the partial derivative of the loss
\(\loss(\hat{y}, y)\) with respect to the \(l\)-th component of the prediction
\(\hat y\). For \(l\in V_2\), \(j\in V_1\) and \(i\in V_1\) we then have
\begin{align*}
    \partial_{\bias_l} \cost_\problem(\param)
    &=
    \E_\problem\bigl[
        \partial_{\hat{y}_l}\loss(\response_{\param}(X), Y) \underbrace{\partial_{\bias_l}(\response_{\param}(X))_l}_{=1}
    \bigr]
    \overset{\response_\param=\response_{\tilde\param}}=
    \partial_{\tilde\bias_l} \cost_\problem(\tilde\param)
    = 0,
    \\
    \partial_{\weight_{jl}} \cost_\problem(\param)
    &=
    \E_\problem\bigl[\partial_{\hat{y}_l}\loss(\response_{\param}(X), Y)
        \underbrace{
            \partial_{\weight_{jl}}(\response_{\param}(X))_l
        }_{
            =\activation(\tilde{\bias}_j + \langle \tilde{\weight}_{\placeholder j}, X\rangle)
        }
    \bigr]
    \overset{\response_\param=\response_{\tilde\param}}=
    \partial_{\tilde\weight_{jl}} \Cost(\tilde\param) = 0,
    \\
    \partial_{\bias_j} \cost_\problem(\param)
    &=
    \E_\problem\Bigl[
        \sum_{l\in V_2}\partial_{\hat{y}_l}\loss(\response_{\param}(X), Y)
        \underbrace{\weight_{jl}}_{=\mathrlap{\begin{cases}
            (1-\lambda \delta_{j2})\tilde\weight_{jl}
            & j\neq 1 
            \\
            \lambda \tilde\weight_{2l}
            & j = 1 
        \end{cases}
        }\qquad}
        \activation'(\bias_j
        + \langle X, \weight_{\placeholder j}\rangle)
    \Bigr]
    \\
    \overset{\response_\param=\response_{\tilde\param}}&=
    \begin{cases}
        (1-\lambda \delta_{j2})\partial_{\tilde\bias_j} \cost_\problem(\tilde\param) & j\neq 1
        \\
        \lambda \partial_{\tilde\bias_2} \cost_\problem(\tilde\param)  & j= 1
    \end{cases}
    \\
    &= 0,
    \\
    \partial_{\weight_{ij}} \cost_\problem(\param)
    &= \E\Bigl[
        \sum_{l\in V_2}\partial_{\hat{y}_l}\loss(\response_{\param}(X), Y)
        \activation'(\bias_j
        + \langle X, \weight_{\placeholder j}\rangle)
        \weight_{jl}X_i
    \Bigr]
    \\
    \overset{\response_\param=\response_{\tilde\param}}&=
    \begin{cases}
        (1-\lambda \delta_{j2})\partial_{\tilde\weight_{ij}} \cost_\problem(\param) & j\neq 1
        \\
        \lambda \partial_{\tilde\weight_{i2}} \cost_\problem(\param)  & j= 1
    \end{cases}
    \\
    &= 0.
\end{align*}
the parameter \(\param\) is thereby clearly a critical point with no other redundancies
except for a single duplication.

\paragraph*{Proof of \ref{it: only deac and bias}.}
To construct a parameter \(\param\) with a deactivation and bias redundancy
from the reduced network, define
\[
    \weight_{\placeholder 1} = 0,
    \qquad
    \weight_{1 \placeholder} = 0,
\]
select \(\bias_1 \in \real\) arbitrarily and retain all parameters of \(\tilde
\param\) for the other neurons.

Since the additional neuron is deactivated, the response remains the same, i.e.
\(\response_\param = \response_{\tilde\param}\).
And since \(\tilde\param\) is a critical point, it is straightforward
to check that the derivatives with respect to the old parameters remain
the same and are thereby zero. For the derivatives with respect to the new
parameters let us consider the outer derivatives first
\begin{align*}
    \partial_{\weight_{1 l}}\cost_\problem(\param)
    &=
    \E_\problem\bigl[
        \partial_{\hat{y}_l}\loss(\response_\param(X), Y)
        \activation\bigl(
            \bias_1 + \langle X, \weight_{\placeholder 1}\rangle
        \bigr)
    \bigr]
    \\
    \overset{\weight_{\placeholder 1}=0}&=
    \underbrace{\E_\problem\bigl[
        \partial_{\hat{y}_l}\loss(\response_\param(X), Y)
    \bigr]}_{=\partial_{\tilde\bias_l}\cost_\problem(\tilde\param) = 0}
    \activation(\bias_{1}).
\end{align*}
In the last equation we used that the response remains the same.
The derivatives with respect to the inner derivatives on the other hand are all
zero due to the deactivation with the outer parameter
\(\weight_{1\placeholder}=0\).

\subsection{Pruning}
\label{sec: pruning}

The following result shows that for any redundant parameter there exists an
efficient parameter of a smaller network with equal response function on
the support \(\Domain\) of the input \(X\). We also show that criticality is
retained under the standard assumption of \(\# V_\tout=1\) (cf.
Definition~\ref{def: standard model}).

\begin{prop}
    \label{prop: pruning}
    Let \(\network = (\mathbb V, \activation)\) with \(\mathbb V=(V_0, V_1,
    V_2)\) be a shallow ANN as in Definition~\ref{def: shallow neural network}.
    Assume that the parameter \(\param\) of this network is redundant. Then
    there exists a pruned network \(\tilde{\network} = ((V_0, \tilde{V}_1, V_2),
    \activation)\) with \(\tilde V_1 \subseteq V_1\) and an \textbf{efficient}
    parameter \(\tilde \param\) of this pruned network such that the response
    remains the same, i.e.
    \[
        \response_\param(x) = \response_{\tilde \param}(x) \qquad \forall x\in \Domain.
    \]
    Furthermore, if \(\# V_\tout=1\) and \(\tilde \param\) was a critical point
    of some cost function
    \[
        \cost_\problem(\param) = \E_\problem[\loss(\response_{\param}(X), Y)]
    \]
    for some loss function \(\loss\) and expectation \(\E_\problem\) induced by
    some distribution \(\Pr_\problem\), then (assuming
    suitable regularity on \(\loss\) and \(\activation\) such that derivatives
    may be moved into the expectation \(\E_\problem\), cf. Lemma~\ref{lem:
    differentiability}) the pruned parameter \(\tilde \param\) is also a
    critical point of \(\cost_\problem\).
\end{prop}

\begin{proof}
    A parameter is redundant if either of the two criterions \ref{it: disused neuron} 
    or \ref{it: linear combination of neurons} in the Definition of Efficiency
    \ref{def: efficient} are violated. Recall, that we called the violation
    of criterion \ref{it: disused neuron} a deactivation redundancy
    (Remark~\ref{rem: taxonomy of redundant parameters}).
    
    \paragraph*{Deactivation pruning}
    It is straightforward to see that the response of an ANN does not change
    if all the deactivated hidden neurons, i.e. all \(j\in V_1\) where
    \(\weight_{j\placeholder} =0\), are removed from \(V_1\). Similarly,
    it is straightforward to show that critical parameters remain critical
    since the previous gradient contains all the partial derivatives with
    respect to the remaining parameters in the pruned network.

    \paragraph*{Pruning the other redundancies}
    Until the parameter is efficient we will iteratively remove a single neuron,
    while ensuring that the response stays the same and critical points remain
    critical. Since there only a finite number of neurons, this procedure will
    eventually terminate -- if there are no hidden neurons left, then there
    is only a bias on the output which is clearly efficient (Definition
    \ref{def: efficient}). We therefore only describe the procedure of a single
    step.

    Note that if a pruning step reintroduces deactivation redundancies, we
    interject a deactivation pruning step. We can therefore always assume there
    are no deactivation redundancies at the beginning of a pruning step.

    If the parameter \(\param\) is redundant without deactivation redundancies
    \ref{it: disused neuron}, then there must be a hidden neuron \(k\in V_1\) that
    can be linearly combined from the others (cf. Remark~\ref{rem: taxonomy of
    redundant parameters}), i.e.
    \[
        \activation\Bigl(\bias_k + \sum_{i\in V_0}  x_i\weight_{ik} \Bigr)
        = \lambda_\emptyset + \sum_{j\in V_1\setminus \{k\}}\lambda_j
        \activation\Bigl(\bias_j + \sum_{i\in V_0} x_i\weight_{ij}\Bigr)
        \quad \forall x\in \Domain.
    \]
    We define a parameter \(\tilde \param\) for the pruned network
    \(\tilde \network = ((V_0, V_1\setminus\{k\}, V_2), \activation)\) using the
    parameter \(\param\) from the old network. Specifically, we retain all inner
    parameters and define the outer parameters to be
    \[
        \tilde \weight_{jl} := \weight_{jl} +  \weight_{kl}\lambda_j
        \qquad \tilde \bias_l := \bias_l + \weight_{kl}\lambda_\emptyset 
        \qquad \forall l\in V_2,\; j\in V_1\setminus\{k\}.
    \]
    Using this definition it is straightforward to check that the response
    remains the same, i.e. \(\response_\param = \response_{\tilde \param}\).
    
    In the case \(\# V_\tout=1\), i.e. \(V_\tout = \{\outSgt\}\),
    we need to show that this pruning step retains criticality. Recall that
    the derivatives are given by
    \begin{equation}
        \label{eq: previously critical}    
        \partial_{\param_i}\cost_\problem(\param)
        = \E_\problem\bigl[\partial_{\hat y}(\response_\param(X), Y) \partial_{\param_i}\response_\param(X)\bigr].
    \end{equation}
    Since the derivatives of the response
    \(\partial_{\tilde \bias_\outSgt}\response_{\tilde\param}\) and
    \(\partial_{\tilde \weight_{j\outSgt}}\response_{\tilde\param}\) with respect to the outer
    parameters only contain inner parameters (which we have not changed) and the
    response remains the same \(\response_{\tilde \param} = \response_\param\),
    we immediately get the criticality of the partial derivatives with
    respect to the outer parameters. What is left to consider are the derivatives
    with respect to the inner derivatives. Since \(\param\) had no deactivation
    redundancies, we have \(\weight_{j \outSgt} \neq 0\) for all \(j\in V_1\).
    In particular we have for all \(j\in V_1\setminus \{k\}\)
    \begin{align*}
        \partial_{\tilde \bias_j}\response_{\tilde \param}(x)
        &= \activation'(\bias_j + \langle x, \weight_{\placeholder j}\rangle)
        \tilde\weight_{j \outSgt}
        = \partial_{\bias_j}\response_{\param}(x) \frac{\tilde \weight_{j \outSgt}}{\weight_{j\outSgt}}
        \\
        \partial_{\tilde \weight_{ij}}\response_{\tilde \param}(x)
        &= \activation'(\bias_j + \langle x, \weight_{\placeholder j}\rangle)
        x_i\tilde\weight_{j \outSgt}
        = \partial_{\weight_{ij}}\response_{\param}(x) \frac{\tilde \weight_{j \outSgt}}{\weight_{j\outSgt}}
    \end{align*}
    The derivatives of the response therefore only change up to a constant that
    can be moved out of the expectation in \eqref{eq: previously critical}.
    Together with \(\response_{\tilde \param} = \response_\param\) this yields criticality.
\end{proof}

\subsection{Bias redundancies (Proof of Proposition \ref{prop: pure bias redundancies can be ruled out})}
\label{sec: critical bias red do not exist}

Recall that a bias redundancy as defined in Remark \ref{rem: taxonomy of redundant parameters}
implies that \(\activation_k(x) := \activation(\bias_k + \langle x,
\weight_{\placeholder
k}\rangle)\) is constant on \(\Domain\).

\begin{lemma}[Bias redundancy characterization]
    \label{lem: bias redundancy characterization}
    Let the activation function \(\activation\) be injective and assume
    \begin{equation}
        \label{eq: orthogonal}    
        \{x-y: x,y\in \Domain\}^\perp = \{0\}.
    \end{equation}
    Then a bias redundancy at neuron \(k\) implies \(\weight_{\placeholder k} = 0\).
\end{lemma}
Observe that \eqref{eq: orthogonal} is satisfied as soon as \(\Domain\) contains
an open set.
\begin{proof}
    Let there be a bias redundancy at neuron \(k\). If there were
    \(x,y\in \Domain\) such that \(\langle x - y, \weight_{\placeholder k}\rangle \neq 0\),
    then \(\activation_k(x) \neq \activation_k(y)\) due to
    injectivity. Consequently \(\langle x-y, \weight_{\placeholder k}\rangle = 0\)
    for all \(x,y\in \Domain\) and \eqref{eq: orthogonal}
    thereby implies \(\weight_{\placeholder k}  = 0\).
\end{proof}

Recall that we assumed in Proposition \ref{prop: pure bias redundancies can be
ruled out} \(\# V_\tout=1\),  \(\activation'(x) \neq 0\) for all \(x\in \real\),
which also implies \(\activation\) is injective as it is
strictly monotonous, and an open set in \(\Domain\) such that Lemma~\ref{lem:
bias redundancy characterization} is satisfied.

In Section~\ref{sec: pruning} we discussed a general pruning procedure that
proceeds in steps, removing one neuron at a time. Consequently, this procedure
is path dependent. If a different neuron were removed first one might
end up with a different pruned network and parameter. In the case where we
only have bias redundancies we can do better.
Assume the requirements of Lemma~\ref{lem: bias redundancy characterization}
are satisfied, let \(I\subseteq V_1\) be the maximal set
such that \(\weight_{\placeholder j} = 0\) for all \(j\in I\).
We then define the pruned network \(\tilde{\network} := ((V_0, V_1\setminus I, V_2), \activation)\)
in a single step: The parameter \(\tilde \param\) 
retains all the weights and biases from \(\param\)
restricted to the pruned ANN-graph except for
\[
    \tilde \bias_l := \bias_l + \sum_{j\in I}\weight_{jl}\activation(\bias_j)
    \qquad l \in V_2.
\]
It is straightforward to show that the response then remains the same.
In the following lemma we will relate the criticality of \(\tilde \param\)
to that of \(\param\).

\begin{lemma}[Characterization of critical bias redundancies]
        \label{lem: characterization of critical bias redundancies}
        Assume the setting of Proposition \ref{prop: pure bias redundancies can be ruled out}.
        If a critical point \(\param\) of the cost \(\cost_\problem\) only has
        bias redundancies, then the parameter \(\tilde \param\) 
        of the pruned network is also a critical point of \(\cost_\problem\) and
        the following equation is satisfied
        \begin{equation}
            \label{eq: condition for extended is crit. point}
            \E_\problem\bigl[
                \partial_{\hat y}\loss(\response_{\tilde\param}(X), Y)X_i
            \bigr]
            = 0 \qquad \forall i\in V_0.
        \end{equation}
        This is furthermore sufficient, i.e. if \(\tilde \param\) is critical
        and \eqref{eq: condition for extended is crit. point} is satisfied
        then the original parameter \(\param\) is critical.
\end{lemma}
\begin{remark}
    We do not make use of the squared error loss function and only
    require sufficient regularity that derivatives may be moved into the
    expectation.
\end{remark}
\begin{proof}
    ``\(\Rightarrow\)'':
    That \(\nabla\cost_\problem(\param)=0\) implies \(\nabla\cost_\problem(\tilde{\param})=0\)
    is a straightforward exercise since the response remains the same and we retain
    almost all parameters except for the outer bias which does not occur in
    any of the partial derivatives of the response. Since we assume
    \(V_1=\{\outSgt\}\) in this section we furthermore have
    \begin{equation}
        \label{eq: derive bias critical condition}
        0 = \partial_{\weight_{ij}}\cost_\problem(\param)
        = \weight_{j\outSgt} \E_\problem\bigl[
            \partial_{\hat y}\loss(\response_\param(X), Y)X_i
        \bigr]\activation'(\bias_j).
    \end{equation}
    And since we assume \(\activation'(x) \neq 0\) for all \(x\in \real\)
    in this section and \(\weight_{j\outSgt} \neq 0\) since we ruled out
    deactivation redundancies, we obtain \eqref{eq: condition for extended is crit. point}.
    
    ``\(\Leftarrow\)'':
    Using \eqref{eq: condition for extended is crit. point}
    and \(\nabla\cost_\problem(\tilde{\param}) = 0\) we now have to prove
    \(\nabla\cost_\problem(\param)=0\). The directional derivatives of the
    parameters that remained on the pruned ANN-graph are zero because they
    coincide with those of \(\tilde{\param}\). What is left are therefore
    the directional derivatives of the parameters attached to the nodes \(j\in
    I\). Using \eqref{eq: derive bias critical condition}
    in reverse with \eqref{eq: condition for extended is crit. point}
    we obtain that \(\partial_{\weight_{ij}}\cost_\problem(\param) =0\)
    What is left are therefore the biases \(\bias_j\) with \(j\in I\).
    Those are given by
    \[
        \partial_{\bias_j} \cost_\problem(\param)
        = \E\bigl[
            \partial_{\hat y}\loss(\response_\param(X), Y)
            \activation'\bigl(\bias_j\bigr)
        \bigr]\weight_{j\outSgt} 
        = \underbrace{\partial_{\tilde{\bias}_\outSgt}\cost_\problem(\tilde \param)}_{=0}
        \activation'\bigl(\bias_j\bigr)\weight_{j\outSgt}.
        \qedhere
    \]
\end{proof}

\subsubsection{The pruned condition a.s.\ never happens (Proof of \ref{it: no pure bias redundancies})}

With the characterization of critical bias redundancies (Lemma~\ref{lem: characterization
of critical bias redundancies}), proving the non-existence of bias redundancies
is equivalent to proving that a parameter vector
of an efficient network can never be a critical point which also
satisfies \eqref{eq: condition for extended is crit. point}.
To prove \ref{it: no pure bias redundancies} we therefore simply have to show
that \eqref{eq: condition for extended is crit. point} almost
surely never coincides with an efficient critical point.

Since the efficient parameters are automatically
\((1,0,1)\)-polynomially efficient by assumption, we
need to show that the set \(\polyEfficient^{(1,0,1)}\)
almost surely does not contain critical points that
satisfy \eqref{eq: condition for extended is crit. point}.
Recall that we assume the squared error \(\loss(\hat{y}, y) = (\hat y - y)^2\)
in \ref{it: no pure bias redundancies} and therefore
\eqref{eq: condition for extended is crit. point}
reduces to
\[
    0 = \E_\problem\bigl[(\response_\param(X) - Y)X_i\bigr]
    \overset{\text{tower}}= \E_\problem\bigl[(\response_\param(X) - \target_\problem(X))X_i\bigr]
    \qquad \forall i \in V_0.
\]
For the random cost \(\Cost=\cost_\Problem\) this means
\[
    0 = \int (\response_\param(x) - \rf(x))x_i \Pr_X(dx)
    = \langle \response, \pi_i\rangle_{\Pr_X} - \langle \rf, \pi_i\rangle_{\Pr_X}
\]
with random target \(\rf=\target_\Problem\) and projection \(\pi_i:x\mapsto x_i\).
This combination can be captured by	the level set \(\rg^{-1}(0)\) of
\[
    \rg:\begin{cases}
        \polyEfficient^{(1,0,1)} \to \real^{\dim(\param)}\times \real^{V_\tin}
        \\
        \param \mapsto
        \Bigl(
            \nabla\Cost(\param),
            \bigl(
                \langle \response, \pi_i\rangle_{\Pr_X} - \langle \rf, \pi_i\rangle_{\Pr_X}
            \bigr)_{i\in V_\tin}
        \Bigr).
    \end{cases}
\]
Since \(\polyEfficient^{(1,0,1)}\subseteq \real^{\dim(\param)}\), \(\rg\) is a
mapping into a larger dimension. Its level sets are therefore empty with
probability one by Lemma~\ref{lem: generalization 11.2.10},
assuming we can prove it to be non-degenerate for every \(\param\in
\polyEfficient^{(1,0,1)}\). This will therefore be the finial step of
the proof. Since the variance of \(\rg\) is independent of the mean,
we may consider \(\hat{\Cost}(\param)=\langle \response_\param, \rf\rangle_{\Pr_X}\) as introduced in Proposition~\ref{prop:
decomposition} instead and similarly prune \(\langle \response_\param, \pi_i\rangle_{\Pr_X}\), i.e. we may consider
\[
    \hat{\rg}(\param)
    = \Bigl(\nabla_\param\langle \response_\param, \rf\rangle_{\Pr_X}, (\langle \rf, \pi_i\rangle_{\Pr_X})_{i\in V_\tin}\Bigr)
    = \Bigl(\langle \nabla_\param\response_\param, \rf\rangle_{\Pr_X}, (\langle \pi_i, \rf\rangle_{\Pr_X})_{i\in V_\tin}\Bigr)
\]
Similar to our argument in the proof of Proposition~\ref{prop: non-degenrate g}
it is therefore sufficient to prove the linear independence of the following
functions
\begin{align}
    &x\mapsto 1 
    \tag{\(\partial \bias_{\outSgt}\)}
    \\
    &x \mapsto \activation\bigl(\bias_j + \langle x, \weight_{\placeholder j}\rangle\bigr)
    && j\in V_1 
    \tag{\(\partial \weight_{j\placeholder}\)}
    \\
    &x \mapsto \activation'\bigl(\bias_j + \langle x, \weight_{\placeholder j}\rangle\bigr)\weight_{j\placeholder} x_i
    && j\in V_1, i\in V_0
    \tag{\(\partial \weight_{ij}\)}
    \\
    &x \mapsto \activation'\bigl(\bias_j + \langle x, \weight_{\placeholder j}\rangle\bigr)\weight_{j\placeholder} 
    && j\in V_1 
    \tag{\(\partial \bias_j\)}
    \\
    &x\mapsto x_i
    && i\in V_0
    \tag{\(\pi_i\)},
\end{align}
where the last equation is the extra condition that follows from \eqref{eq: condition for extended is crit. point}. But their linear
independence follows from the \((1,0,1)\)-polynomial independence
(Definition~\ref{def: polynomial independence}).

Note, that we did not require second order derivatives, but
first order polynomials in the affine term due to the extra condition
\eqref{eq: condition for extended is crit. point}.
This explains why we required \((1,0,1)\)-polynomial independence
instead of \((0,0,1,2)\)-polynomial independence as in Proposition~\ref{prop:
non-degenrate g}.

  			\printbibliography[heading=subbibliography]
		}
	\end{refsection}

\end{document}